\documentclass[11pt]{amsbook}

\usepackage{latexsym,amssymb}%,times,mathptm}
\usepackage{amsmath}
\usepackage{amscd}
\usepackage{emlines}
\usepackage{psboxit}
\PScommands
\usepackage{color}
\usepackage{xypic}

\makeindex

\newtheorem{Theorem}{Theorem}[section]
\newtheorem{Lemma}[Theorem]{Lemma}
\newtheorem{lemma}[Theorem]{Lemma}
\newtheorem{Corollary}[Theorem]{Corollary}
\newtheorem{corollary}[Theorem]{Corollary}
\newtheorem{Proposition}[Theorem]{Proposition}
\newtheorem{proposition}[Theorem]{Proposition}
\newtheorem{Conjecture}[Theorem]{Conjecture}
\newtheorem{conjecture}[Theorem]{Conjecture}

\newtheorem{Remark}[Theorem]{Remark}
\newtheorem{remark}[Theorem]{Remark}
\newtheorem{Example}[Theorem]{Example}
\newtheorem{example}[Theorem]{Example}
\newtheorem{Definition}[Theorem]{Definition}
\newtheorem{definition}[Theorem]{Definition}
\newtheorem{Question}[Theorem]{Question}
\newtheorem{question}[Theorem]{Question}
\newtheorem{Questions}[Theorem]{Questions}

\newtheorem{Problem}[Theorem]{Problem}

\newtheorem{Notation}[Theorem]{Notation}

\theoremstyle{remark}

\numberwithin{equation}{Theorem}
\def\d{\mbox{\rm d}}
\def\e{\mathrm{e}}

\def\H{{\mathrm H}}
\def\m{\mathfrak{m}}
\def\n{\mathfrak n}
\def\l{\mathfrak l}
\def\p{\mathfrak p}
\def\q{\mathfrak q}

\def\p{{\mathfrak  p}}
\def\Q{{\mathcal{ Q}}}
\def\q{{\mathfrak  q}}

\def\H{\mathrm{H}}

\def\Gr{\mathrm G}

\newcommand{\rme}{\mathrm{e}}

\newcommand{\rmh}{\mathrm{h}}

\newcommand{\rmH}{\mathrm{H}}

\newcommand{\rmK}{\mathrm{K}}

\newfont{\frank}{eufm10 scaled\magstep1}

\def\p{\mbox{\frank p}}
\def\q{\mbox{\frank q}}
\def\m{\mbox{\frank m}}

\def\Q{\mbox{\frank Q}}

\newfont{\frankk}{eufm10 scaled\magstephalf}

\def\bbz{{{\mathbb Z}}}
\def\bbq{{{\mathbb Q}}}

\def\bbk{{{\mathbb K}}}
\def\bbc{{{\mathbb C}}}
\def\bbp{{{\mathbb P}}}
\def\bbn{{{\mathbb N}}}

\def\XX{{\bf X}}
\def\KK{{\bf K}}
\def\LL{{\bf L}}
\def\vv{{\bf v}}

\def\xx{{\bf x}}
\def\yy{{\bf y}}
\def\zz{{\bf z}}
\def\TT{{\bf T}}
\def\tt{{\bf t}}
\def\SS{{\bf S}}
\def\YY{{\bf Y}}

\def\ZZ{{\bf Z}}
\def\ff{{\bf f}}
\def\g2{{\bf g}}
\def\FF{{\bf F}}
\def\hh{{\bf h}}
\def\RR{{\bf R}}
\def\AA{{\bf A}}
\def\BB{{\bf B}}
\def\CC{{\bf C}}
\def\DD{{\bf D}}
\def\GG{{\bf G}}

\def\II{{\bf I}}
\def\JJ{{\bf J}}
\def\hh{{\bf h}}
\def\aa{{\bf a}}
\def\bb{{\bf b}}
\def\cc{{\bf c}}
\def\qq{{\bf q}}
\def\pp{{\bf p}}
\def\ttt{{\bf t}}
\def\HH{{\bf H}}

\newcommand{\mi}[1]{#1\index{#1}}

\def\ad{ \mbox{\rm ad}}
\def\adeg{\mbox{\rm arith-deg}}
\def\ann{\mbox{\rm ann}}
\def\Ass{\mbox{\rm Ass}}

\def\bar#1{{\overline{#1}}}

\def\bdeg{\mbox{\rm bdeg}}

\def\br{\mbox{\rm br}}

\def\cl{\overline}
\def\codim{\mbox{\rm codim }}

\def\coker{\mbox{\rm coker }}

\def\ddet0{\mbox{\rm det}_0 }
\def\Deg{\mbox{\rm Deg}}

\def\depth{\mbox{\rm depth }}

\def\ds{\displaystyle}

\def\ecodim{\mbox{\rm ecodim}}
\def\edeg{\mbox{\rm edeg}}
\def\emb{\mbox{\rm emb}}

\def\End{\mathrm{End}}
\def\Ext{\mbox{\rm Ext}}

\def\Fitt{\mbox{\rm Fitt}}

\def\gdeg{\mbox{\rm gdeg}}
\def\Gr{\mbox{\rm  $G$}}
\def\gr{\mbox{\rm gr}}
\def\grade{\mbox{\rm grade }}

\def\hdeg{\mbox{\rm hdeg}}
\def\height{\mbox{\rm height }}
\def\Hom{\mbox{\rm Hom}}

\def\image{\mbox{\rm image }}

\def\jdeg{\mbox{\rm jdeg}}

\def\l{\lambda}
\def\lar{\longrightarrow}

\def\length{\mbox{\rm length}}

\def\mult{\mbox{\rm mult}}

\def\pp{{\mathbf p}}

\def\QED{\hfill$\Box$}

\def\rad{\mbox{\rm rad}}
\def\rank{\mbox{\rm rank}}
\def\rar{\rightarrow}
\def\red{\mbox{\rm r}}

\def\Rees{\mbox{${\mathcal R}$}}
\def\reg{\mbox{\rm reg}}
\def\reltype{\mbox{\rm reltype}}

\def\Spec{\mbox{\rm Spec}}

\def\sup{\mbox{\rm sup}}
\def\supp{\mbox{\rm supp}}
\def\surj{\twoheadrightarrow}
\def\surjects{\twoheadrightarrow}
\def\Symi{\mbox{$\mathcal S$}}
\def\Sym{\mbox{\rm Sym}}

\def\tn{\mbox{\rm tn}}
\def\Tor{\mbox{\rm Tor}}

\numberwithin{section}{chapter}
\newcommand\HRule{\noindent\rule{\linewidth}{2.5pt}}

\begin{document}

\begin{titlepage}
%\vspace*{\stretch{1}}\HRule

%\end{document}
%\begin{flushright}
\begin{center}
\fontsize{30}{40}\selectfont{\bf Complexity Degrees of \\ Algebraic Structures}
\end{center}
%\end{flushright}

\HRule
\begin{center}
\vspace{2cm}
\textsc{Wolmer V Vasconcelos \\ Department of Mathematics\\ Rutgers University \\ 110 Frelinghuysen Rd, Piscataway\\  New Jersey 08854-8019} \\  \textrm{vasconce@math.rutgers.edu}
\end{center}
\vspace*{\stretch{2}}
\begin{center}
\textsc{December 2013}
\end{center}
\end{titlepage}

\author{
{Wolmer V Vasconcelos}}%\\
%\thanks{Partially
% supported by the NSF Grants 0097093 and 0500379}
%\vspace{-0.75mm}\\
%\vspace{-1mm}}
%\footnotetext{AMS 2000 {\it Mathematics Subject
%Classification}. Primary 13H10; Secondary 13F20, 05C90.}
%\address{
%Department of Mathematics,\\
%Rutgers University \\
%110 Frelinghuysen Rd, Piscataway\\
% New Jersey 08854-8019}
%\email{vasconce@math.rutgers.edu}
%\end{center}

%\date=\today

%\maketitle

\tableofcontents

%\chapter*{Preface}

%{\bf Warning:} These notes are in a  preliminary condition.
% The material is often the dump
%of manuscripts and whole sections will be re-written. Only forgiving
%friends, from whom we are soliciting suggestions and comments, should
%look at them.

%\medskip

%In response to comments, suggestions and corrections,
%the content undergoes frequent upd
%ates.

\chapter*{Prologue}
%\addcontentsline{toc}{chapter}{Introduction}

\noindent
{\bf A Word  to Readers:} These notes are in a  preliminary condition.
We gratefully solicit suggestions and comments
 from kind readers and colleagues with regard to:

 \begin{itemize}

 \item[$\bullet$] missed or misplaced references

 \item[$\bullet$] repeated typos or awkward notation

 \item[$\bullet$] missing related topics

 \item[$\bullet$] inconsistent use of fonts

 \end{itemize}

\medskip

\noindent
We are grateful to numerous  colleagues [which will eventually be listed] who endured conversations about the topics of these notes. To
 Jooyoun [J] Hong we give a special word of thanks for all the support, mathematical and technical.

%\begin{flushleft}
%Wolmer V Vasconcelos\\
%Department of Mathematics\\
%Rutgers University\\
%New Brunswick, NJ 08854 \\
%vasconceATmath.rutgers.edu
%\end{flushleft}

\chapter*{Abstract}

\noindent
We aim at
studying
collections $\mathcal{C}(\RR)$ of algebraic structures defined over a
commutative ring $\RR$ and investigating the
complexity of significant constructions carried out on these objects.
Noteworthy
 are the category of finitely generated modules (over local rings, or
graded, particularly vector bundles), finitely generated algebras or objects
suitable of more fundamental decompositions and construction of
various closure operations---such as the determination of global sections and
the (theoretical or machine) computation
 of integral closures of algebras. Typically, the
operations  express
 smoothing processes that enhance the structure of the algebras,
enable them to support new constructions  (including
analytic ones),  divisors acquire a group
structure, and the cohomology tends to slim down.
The assignment of measures of size, via a multiplicity theory, to the
algebras and to  the  {\em construction} itself
is a novel aspect to the subject.
One of its specific  goals is to develop a comprehensive theory
of normalization in commutative algebra  and  a broad set of multiplicity functions to be
used as complexity benchmarks.

\bigskip

{\bf Short Abstract}

\bigskip

\noindent
We aim at
studying
collections  of algebraic structures defined over a
commutative ring  and investigating the
complexity of significant constructions carried out on these objects.
The assignment of measures of size, via a multiplicity theory, to the
algebras and to  the   construction itself
is a novel aspect to the subject.

\chapter*{Introduction}

\noindent
We seek to define numerical functions, $\Deg$, which
assigns to an algebraic structure $\AA$ the value $\Deg(\AA)$, and to
a smoothing operation
\[ \AA \leadsto \BB,\]
a value $\Deg(\leadsto)$ to the transformation itself. The  naive
expectations  set $\Deg(\leadsto)=|\Deg(\BB)-\Deg(\AA)|$, or
$\Deg(\leadsto)=\Deg(\BB/\AA)$,  but a
closer accounting may often be given. In the case of an algebra or
module $\AA$ the functions used $\Deg(\AA)$ are enhanced versions of the ordinary
multiplicities $\deg(\AA)$ where the difference $\Deg(\AA)-\deg(\AA)$
should represent deeper properties of $\AA$ not read off $\deg(\AA)$.

There are two interlocked grand themes arising from this general
endeavor. First, the
development of  multiplicities, that is numerical (degree)
functions $\Deg: \mathcal{C}(\RR)\rar \bbn$ with suitable properties,
such as
with regard to hyperplane sections and particular short exact
sequences. The functions treated here  come in two families, the
{\em arithmetic degrees} and the {\em  cohomological
degrees}. They can be understood as providing a {\em geometric volume} in
the classical manner, and a {\em cohomological volume}.
The first are
directly related to the standard multiplicity and include the
classical arithmetic and geometric degrees, the Buchsbaum-Rim
multiplicity and the more recent $j$-multiplicity and $\jdeg$.
Among the  cohomological degrees, usually referred to as big Degs, are
$\hdeg$ and $\bdeg$. These two  functions have different
 purposes, the first
addressing  issues of genericity, while the other
has Cohen-Macaulay structures as benchmarks and seeks to measure
deviations from corresponding  baselines. A fruitful aspect in these
developments are the numerous questions seeking to relate them
 to the Castelnuovo-Mumford regularity in the case of
graded structures.

The other theme is the character of  the specific problems
the multiplicities are being used
on, focused on the search for  estimates for the number of steps that
general algorithms must use to calculate the closure of some operation, to
determine the number of generators of an algebra or module, or just
to bound the number of Hilbert functions of strata of elements in
$\mathcal{C}(\RR)$.

   The novel perspective of
its subject, linking commutative algebra and computational complexity,
 has matured into numerous
explicit  challenges of great technical difficulty. Rooted in the
current worldview that takes Cohen-Macaulay structures as the basic
blocks of the field, we  seek  to track how those
objects undergo smoothing operations and the new structures that
arise.

%\medskip

 This monograph deals with several themes in
 commutative algebra
and computational algebra centered on the construction and analysis
of integral closure of algebraic structures of great interest along with  the development of
degree functions.
It ranges from the practical construction of algorithms for integral
closures and
the development of settings to analyze their complexity to the
examination of several classes of blowup rings.
We hold the view that Cohen--Macaulay structures are ubiquitous in
algebra, are central to commutative algebra and algebraic geometry
  and
represent unique computational efficiencies and that without an
understanding of the ways it appears in the structure of non
Cohen-Macaulay structures   the prospects for successful large scale
symbolic computations are rather dim.
The project
  is focused on studying graded algebras, blowup algebras,
 graded rings associated to ideals of combinatorial interest,
the
 Cohen--Macaulayness of certain toric varieties and the role they may
  play in the theory of
Gr\"obner bases.

The notes are organized into 5 chapters and one appendix, each split into sections:

\begin{enumerate}
\item   {Chapter 1:}  Degree/Multiplicity Functions
\item   {Chapter 2:} Normalization of Algebras, Ideals and Modules
\item   {Chapter 3:} The Hilbert Coefficients and the Character of Local Rings
\item   {Chapter 4:} Multiplicity-Based Complexity of Derived Functors
\item   {Chapter 5:} The Equations of Ideals
\item   {Chapter 6:  Appendix}
\end{enumerate}

From all of these topics we select some high points.

\medskip

\subsubsection*{Cohomological degrees}

The multiplicity functions of algebraic geometry and commutative
algebra have wide usage in sizing up a module or an algebra. Its
drawbacks are that ignores the lower dimensional components of the
structure and does not behave too well with respect to hyperplane
sections. To account for the components, Bayer and Mumford
 introduced the {\em arithmetic degree}. To
account for a fuller behavior of hyperplane sections, {\em extended
degrees} were introduced. They all coincide on Cohen-Macaulay rings
and algebras. Unlike multiplicities, that have a combinatorial
foundation, extended degrees are both combinatorial and very
homological. The need to express local behavior in global terms
presents many technical challenges for their computation.

Let $(\RR, \mathfrak{m})$ be a Noetherian local ring (or a Noetherian
graded algebra) and let $\mathcal{M}(\RR)$ be the category of finitely
generated $\RR$--modules (or the appropriate category of graded
modules). A {\em degree function}\index{degree function} is simply a
numerical function
\[\mathbf{d}: \mathcal{M}(\RR) \mapsto \mathbb{N}\]
The more interesting of them initialize on modules on finite length,
are additive on {\em certain} short exact sequences
and have mechanisms that control how the functions behave under
generic hyperplane sections. Some of these functions are:
 Classical degree (Multiplicity) and
Castelnuovo--Mumford regularity.

Most degree functions $\mathbf{d}$ are derived from the leading
coefficients of Hilbert
polynomials of certain filtrations. Refinements involves assembling
$\mathbf{d}$ by adding several of these  coefficients, so that a
value such as $\mathbf{d}(A)$ may capture several elements of the
structure of $A$.

The most demanding requirement on these functions
are those  regarding their behavior with respect to generic hyperplane sections.
When $A$ has positive depth and  $h\in R$ is a generic hyperplane
section, the comparison
\[ \mathbf{d}(A) \leftrightarrow \mathbf{d}(A/hA) \]
is the principal divider among the degrees. For those directly
derived from the classical multiplicity, one always has
\[ \mathbf{d}(A) \leq \mathbf{d}(A/hA). \]
For the other family of degrees, the cohomological multiplicities,
the relationship is reversed,
\[ \mathbf{d}(A) \geq \mathbf{d}(A/hA). \]
It is this feature that makes them appealing as complexity
benchmarks, motivating
 our interest  on
 cohomological (or extended) degrees (\cite{DGV}), particularly
 on the homological degree and bdeg (the so-called Big
 Degs)\index{extended degree}
 (\cite{Gunston}):

\medskip

\noindent
{\bf Definition.}
{\rm
  A {\em cohomological
degree}, or {\em extended
multiplicity}
 is a function
\[\Deg(\cdot) : {\mathcal M}(\RR) \mapsto {\mathbb N},\]
that satisfies
\begin{itemize}
\item[\rm {(i)}] If $L = \Gamma_{\mathfrak m}(A)=\H_{\m}^0(A)$ is the submodule of
elements of $A$ which are annihilated by a
power of the maximal ideal and $\overline{A} = A/L$, then
\begin{eqnarray*}
\Deg(A) = \Deg(\overline{A}) + \lambda(L),
\end{eqnarray*}
where
$ \lambda(\cdot)$ is the ordinary length function.
\item[\rm {(ii)}] (Bertini's rule)
 If $A$ has positive depth and $h$ is a generic hyperplane section,
  then
\begin{eqnarray*}
\Deg(A) \geq \Deg(A/hA).
\end{eqnarray*}
\item[{\rm (iii)}] (The calibration rule) If $A$ is a Cohen--Macaulay
module, then
\begin{eqnarray*}
\Deg(A) = \deg(A),
\end{eqnarray*}
where $\deg(A)$ is the ordinary multiplicity of the module $A$.
\end{itemize}
}

\medskip

It was not obvious that such functions even existed.
An infinite
 number of such notions however arose from the construction of
 $\hdeg(\cdot)$ (\cite{hdeg}) simply by using Samuel's
multiplicities.
The philosophy of introducing these functions is the following.
Taking Cohen-Macaulay structures (e.g. algebras over a local ring)
 as a departure, with its rich variety
of techniques, it becomes clear that  the next class of objects
to be examined are those that are Cohen-Macaulay on the punctured
spectrum. To broaden this to more general algebras is accomplished by
the technique of the extended degrees.

A main characteristic  of a big Deg is that it satisfies a form of
Bertini's property ($\Deg(E/hE)$ $\leq \Deg(E)$, for a hyperplane
section) and therefore can be used as a measure of complexity.
The most significant aspect of these particular degrees lies in the fact that
they satisfy a number of rules of calculation with regard to
certain exact  sequences.
In certain situations, they are often more flexible than those followed by
Castelnuovo--Mumford's {\em regularity}
$\reg(\cdot)$ (\cite{BM}, \cite{EG}). One of the comparisons
  shows that
for any standard graded algebra $A$ (\cite{DGV}, and its extension
in \cite{Nagel03}):
\[
 \Deg(A)> \reg(A),
\]
  for any $\Deg(\cdot)$ function, in particular for $\hdeg(\cdot)$.
The fact that it comes coded by an explicit formula will permit
making many {\em a priori} estimates for the number of generators in
terms of
multiplicities.
In the case of graded modules one of properties of $\Deg(\cdot)$ is
that it roughly mediates between the Castelnuovo-Mumford regularity
and the Betti numbers.

\medskip

These degrees have a greater strength than those defined using primarily the primary
decomposition data. For instance, whenever a cohomological degree function is defined, we have 
a decreasing chain
\[ \Deg(A) \geq \adeg(A) \geq \gdeg(A) \geq \jdeg(A).\]

\subsubsection*{Normalization of algebras, ideals and modules}
A general theory of the integral closure of algebraic
 structures--as a full fledged integration
of algorithms, cost estimates, along
 with tight
 predictions at termination--is highly desirable but it is an as yet
 unrealized goal.
Several
 stretches of
 this road
have been built and we will begin to describe the role of
multiplicities in the undertaking.

We will now describe the setting in which our discussion takes place.
Let $\RR$ be a geometric/arithmetic  integral domain, and
let $\AA$ be a graded algebra (more generally with a
$\bbz^r$-grading)
\[ \AA= \RR[x_1, \ldots, x_n], \quad \mbox{\rm  $x_i$ homogeneous},\]
of integral closure
\[ \bar{\AA}= \RR[y_1, \ldots, y_m], \quad \mbox{\rm $y_i$
homogeneous}.\]

The construction of $\bar{\AA}$ usually starts out from an
observation of Emmy Noether in her last paper \cite{Noether}: If
$d\in \AA$ is a regular element of $\AA$ contained in the Jacobian
ideal, then
\[ \AA \subseteq  \bar{\AA} \subseteq \frac{1}{d} \AA.\]

One approach to the  construction
 of $\bar{\AA}$ (\cite{deJong}) consists in finding processes $\mathcal{P}$ that
create extensions
\[ \AA \mapsto \mathcal{P}(\AA) \subset \bar{\AA},\]
\[ \mathcal{P}^n(\AA) = \bar{\AA}, \quad n\gg 0.\]
typically $\mathcal{P}(\AA)= \Hom_{\AA}(I(\AA), I(\AA))$, where
$I(\AA)$ is
some ideal connected to the conductor/Jacobian of $\AA$.
Another approach (\cite{SiSw}), in characteristic $p$, is distinct.
Starting with $\frac{1}{d}\AA$, and using the Frobenius map, it
constructs a descending sequence
\[ V_0=\frac{1}{d} \AA \supset V_1 \supset \cdots \supset V_m=
 \bar{\AA}.
\]

 Predicting properties of $\bar{\AA}$ in $\AA$, possibly by
adding up the changes  $	\mathcal{P}^i(\AA) \rar
\mathcal{P}^{i+1}(\AA)$ (or the similar terms  $V_i$ in the
descending sequences)
 is a natural question in need of
clarification.

\medskip To realize this objective, we select the following specific
goals:
\begin{enumerate}
\item[{\rm (a)}]  We
develop numerical indices for
$\bar{\AA}$, that is for $\AA_i=\mathcal{P}^{i}$,
find $r$ such that \[\bar{\AA}_{n+r}= \AA_n \cdot \bar{\AA}_r, \quad n\geq
0.\]

\item[{\rm (b)}]
How many ``steps'' are there between $\AA$ and
$\bar{\AA}$,
\[\AA= \AA_0 \subset \AA_1 \subset \cdots \subset \AA_{s-1}\subset
\AA_s =
\bar{\AA},\]
where the $\AA_i$ are constructed by an effective
process?

\item[{\rm (c)}]
Express $r$ and $s$ in
terms of invariants of $\AA$.

\item[{\rm (d)}]
Generators of $\bar{\AA}$: Number  and
distribution of their degrees.

\end{enumerate}

Let $\RR$ be a commutative Noetherian ring, and let $\AA$ be a
finitely generated graded algebra (of integral closure $\bar{\AA}$).
   We want to gauge the `distance' from $\AA$ to $\bar{\AA}$.
   To this end, we introduce two degree functions whose properties
   help to determine the  position    of the algebras in the chains.

\begin{enumerate}

\item[{\rm (a)}]  (tracking number) $\tn: \AA\rar \bbn$,    for a certain class of algebras
(includes $\RR$ a field).
We shall use for distance
\[ \tn(\AA)-\tn(\bar{\AA})\]

\item[{\rm (b)}] $\jdeg:$ a multiplicity on all finitely generated graded
$\AA$-modules.    We shall use for distance
\[ \jdeg(\bar{\AA}/\AA)\]

\end{enumerate}

The first of these will be a {\em determinant}.
Consider a graded $\RR=k[z_1, \ldots,
z_r]$--module $E$:    Suppose $E$ has torsion free rank $r$.
    The {\em determinant} of $E$ is the graded module
\[ {\det}_{\RR}(E)=  (\wedge^r E)^{**}.\]

\[ {\det}_{\RR}(E) \simeq \RR[-c(E)],    \quad \tn(E):= c(E)\]
The integer $c(E)$ is the {\bf tracking number} of $E$.    If $E$
has no associated primes of codimension $1$, $\tn(E)$ can be read off
the Hilbert series $H_E(\ttt)$ of $E$:
\[ H_E(t) = [\![E]\!]=
 \frac{\hh(\ttt)}{(1-\ttt)^{d+1}},\quad  \tn(E) = \hh'(1)\]
   This shows that $\tn(E)$ is independent of $\RR$.
It comes with the following elementary property.

\medskip

\noindent
{\bf Proposition.}
Let $\RR=k[z_1, \ldots, z_d]$ and $E\subset F$ two graded
$\RR$--modules of the same rank. If $E$ and $F$ have the condition
$S_2$ of Serre,
\[ \tn(E)\geq \tn(F),\]
with equality if and only if $E=F$.
Therefore, if
\[ E_0 \subset E_1 \subset \cdots \subset E_n\] be a sequence of
distinct
f.g. graded $\RR$-modules of the same rank.
Then
\[ n \leq \tn(E_0)-\tn(E_n).\]

%In Section~\ref{trackingnumber} we will treat this function in a
%broader setting.

\medskip

The other function, $\jdeg$, has different properties.
 One of the properties $\jdeg$ addresses  our main issue.
 Let $\RR$ be a Noetherian domain and $\AA$ a semistandard
graded $\RR$-algebra, that is an algebra finite over a standard
graded $\RR$-algebra.\index{semistandard graded algebra}
Consider a sequence of integral  graded extensions
\[A \subseteq A_0 \rar A_1 \rar A_2 \rar \cdots \rar A_n=\overline{A},\]
where the $A_i$ satisfy the $S_2$ condition of Serre.
Then $n \leq \jdeg(\overline{A}/A)$
(\cite[Theorem 2]{jdeg1}).

The main goal
becomes to express $\jdeg(\overline{A}/A)$ in terms of {\it priori}
data, that is, one wants to view this number as an
  invariant of $A$. Relating it to other, more accessible,
  invariants of $A$ would bring
 considerable predictive power to the function $\jdeg$. This is the
 point of view of \cite{ni1} and \cite{HUV}, when $\jdeg$ is the
 Hilbert-Samuel multiplicity or the Buchsbaum-Rim multiplicity, respectively.

\medskip

\subsubsection*{The Hilbert coefficients and the character of local rings}

A parameter ideal  of a Noetherian local ring $A$ of dimension $d$ is
an ideal $Q$ of dimension $0 $ generated by $d$ elements. Despite
the transparency, their properties code for numerous properties of
the ring $A$ itself. It is worthwhile
to make use of its properties as a barometer of the Cohen-Macaulay
property or of the lack of it.

We treat the results of \cite{chern7} dealing with the
comparative study of two integers associated to a parameter ideal in
a local ring, one derived from a Hilbert function, the other from a
Koszul complex.
 They  play  similar  roles as predictors of the
Cohen-Macaulay property, a fact that led to  conjectural
relationships between them.
 To state the questions and our
results on them, first of all let us fix some notation and recall
terminology.

 Let $(\RR, \mathfrak{m})$ be a Noetherian local ring of dimension
$d>0$, and let $I$ be an $\mathfrak{m}$-primary ideal. There is a
great deal of interest on the set of $I$-good filtrations of $\RR$.
More concretely, on the set of multiplicative, decreasing filtrations
of $\RR$ ideals, \[\mathcal{A}=\{ I_n \; \mid \; I_0=\RR, \;
I_{n+1}=I I_n, n\gg 0 \},\] integral over the $I$-adic filtration,
conveniently coded in the corresponding Rees algebra and its
associated graded ring \[ \Rees(\mathcal{A}) = \sum_{n\geq 0} I_nt^n,
\quad \gr_{\mathcal{A}}(\RR) = \sum_{n\geq 0} I_n/I_{n+1}. \]

Our focus here is on a set of filtrations both broader and more narrowly
defined.   Let
$M$ be a finitely generated $\RR$-module and $I$ an $\m$-primary
ideal. The Hilbert polynomial of the associated graded module \[
\gr_I(M)= \bigoplus_{n\geq 0} I^nM/I^{n+1}M,\] more precisely the
values of $\sum_{k\leq n}\lambda(I^kM/I^{k+1}M)$ for large $n$, can be assembled as \[
P_{M}(n)=\sum_{i=0}^{d} (-1)^i\rme_i(I,M) {{n+d-i}\choose{d-i}}, \] where
$d=\dim M$.  In most of our discussion, locally, either $I$ or $M$ are fixed,
and by simplicity we set $\rme_i(I,M)=\rme_i(M)$ or $\rme_i(I,M)=\rme_i(I)$,
accordingly.  These integers are the Hilbert coefficients
\index{Hilbert coefficients} of $M$ relative to $I$.
Occasionally the first Hilbert coefficient $\rme_1(I, M)$
is referred to as the {\em Chern coefficient} of $I$ relative to $M$
(\cite{chern}).

%\medskip

A great deal of progress has been achieved on
  how the values of $\rme_1(Q,\RR)$
 codes for  structural information about the ring
$\RR$ itself. More explicitly one defines the set
 \[\Lambda (M) = \{\rme_1(Q,M) \; \mid \; Q ~\operatorname{is ~a ~parameter
~ideal ~for} M \},
\] and examines what its structure expresses about $M$.
In case $M=\RR$, this set was analyzed in \cite{chern3} and related papers for the following extremal  properties:
\begin{enumerate}
\item[{\rm (a)}] $0\in \Lambda(\RR)$;
\item[{\rm (b)}] $\Lambda(\RR)$ contains a single element;
\item[{\rm (c)}] $\Lambda(\RR)$ is bounded.
\end{enumerate}

%\medskip

\noindent The following describes some of the main results:

\smallskip

\noindent
{\bf Theorem.} { Let $\RR$ be an  unmixed Noetherian local ring.}
\begin{enumerate}
\item[{\rm (a)}] (\cite[Theorem 2.1]{chern3})
{\em
 $\RR$ is  Cohen-Macaulay if and only
if $\rme_1(Q)=0$ for some parameter ideal $Q$. }

\item[{\rm (b)}] (\cite[Theorem 4.1]{GO1}) {\em  $\RR$ is  a Buchsbaum
 ring if and only
if $\rme_1(Q)$ is constant for all parameter ideals $Q$. }

\item[{\rm (c)}]
 (\cite[Proposition~4.2]{chern3}) {\em $\RR$ is a
Cohen-Macaulay ring on the punctured spectrum if and only if the set
$\Lambda(\RR)$ is finite.}

\end{enumerate}

\medskip

The other integer  arising in our investigation and
 which we want to contrast to $\rme_1(Q)$ is the
following.
Let $Q = (a_1, a_2, \ldots, a_d)$ be a parameter ideal in $A$. We
denote by $\rmH_i(Q)$~($i \in \mathbb Z)$ the $i$--th homology module of
the Koszul complex $\rmK_{\bullet}(Q; A)$  generated
by the system $a_1, a_2, \ldots, a_d$ of parameters in $A$. We put
\[\chi_1(Q) = - \sum_{i \ge 1}(-1)^i\ell_A(\rmH_i(Q))\] and call it
the first {\em Euler characteristic}\index{Euler number} of $Q$; hence \[\chi_1(Q) = \ell_A(A/Q) -
\rme_0(Q)\] by a classical result of Serre (see \cite{AB},
\cite{Serrebook}).

%\medskip

%\newpage

\begin{small}
%\begin{center}
\begin{table}
\centering
\caption{Comparative properties of the first Chern and Euler  numbers
of a parameter ideal $Q$ of a local ring $A$.
 We call attention to the more stringent  requirements
for the $\rme_1(Q)$.
}
\label{tab:3}
\begin{tabular}{|c|c|}
%\caption{hello}
\hline%\noalign{\smallskip}
% &  \\
%\hline
 &  \\
 $-\e_1(Q)\geq 0$ (\cite{MV})
&          $\chi_1(Q)\geq 0$ (\cite[Appendice II]{Serrebook})\\
 &  \\
\hline
 &  \\
 $A$ unmixed \& $\e_1(Q)= 0\Rightarrow $ $A$ C.-M. (\cite{chern3},
 \cite{chern5})
&          $\chi_1(Q)=0 \Rightarrow $ $A$ C.-M. (\cite[Appendice
II]{Serrebook})  \\
  & \\
\hline
  & \\
$A$ unmixed \& $\e_1(Q)$ constant $\Rightarrow$ $A$ Buchsbaum (\cite{GO1})
 & $\chi_1(Q)$ constant $\Rightarrow$ $A$ Buchsbaum (\cite{SV1})\\
 & \\
\hline
\end{tabular}
\end{table}
%\end{center}
\end{small}

%(These comparisons make plausible some form of the following,  which
%is the overall goal of the present research. The actual thrust is on the
%identification of natural settings where this question can be
%effectively addressed.

\medskip

\subsubsection*{Multiplicity-based complexity of derived functors}

 Let $(\RR, \mathfrak{m})$ be a Noetherian local ring  and let $A$ and $B$
 be finitely generated $\RR$-modules.
Motivated by the occurrence of the derived functors of $\Hom_R(A,
\cdot)$ and $A\otimes_R\cdot$ in several constructions based on
$A$ we seek to develop gauges for the sizes for these modules. In the
case of graded modules, a rich {\bf degree} based theory   has been
developed centered on the notion of Castelnuovo regularity. It is
particularly well-suited to handle complexity properties of tensor
products.

Partly driven by its use in normalization processes, we seek to
develop a {\em length} based theory ground on the notion of extended
multiplicity. Starting with a fixed Deg function--and the choice will
be $\hdeg$--we examine various scenarios of the following
question: Express in terms of $\hdeg(A)$ and $\hdeg(B)$ quantities
such as
\begin{itemize}
\item[{$\bullet$}] $\nu(\Hom_R(A,B))$;
\item[{$\bullet$}] $\nu(\Hom_R(\Hom_R(A,R),R))$;
\item[{$\bullet$}] $\nu(\Ext_R^i(A,B)), \ i\geq 1$;
\item[{$\bullet$}] $h_0(A\otimes_RB)=
\lambda(\H_{\mathfrak{m}}^0(A\otimes_RB))$.
\end{itemize}

The sought-after estimations are polynomial functions on
$\hdeg(A)$ and $\hdeg(B)$,  whose
coefficients are given in terms of invariants of $\RR$. The first of
these questions was treated in \cite{Dalili} and \cite{DV2}, who
refer to it as the HomAB question.

\medskip

The {\em HomAB} question
  asks for uniform estimates for the number of generators of
 $\Hom_R(A,B)$ in terms of invariants of $\RR$, $A$ and $B$. The
 extended question asks for the estimates of the number of generators
 of $\Ext_R^{i}(A,B)$ (or of other functors).
%A variation asks the similar question with $A$ and $B$ replaced by
%functors of the category $\mathcal{C}$ of $\RR$-modules, noteworthy
%being $\Hom_{\mathcal{C}}(\Ext_R^1(A, \cdot), \Ext_R^1(B, \cdot))
%$.

In addition to the appeal of the question in basic homological
algebra, such modules of endomorphisms appear frequently in several
constructions, particularly in the algorithms that seek the integral
closure of algebras (see \cite[Chapter 6]{icbook}).
The algorithms used tend to use rounds of operations of the form:
\begin{itemize}
\item[$\bullet$] $\Hom_R(E,E)$: ring extension;
\item[$\bullet$] $A \rar \tilde{A} $: $S_2$-ification of $A$;
\item[$\bullet$] $I:J$: ideal quotient.
\end{itemize}

 Another natural
application is for primary decomposition of ideals, given the
prevalence of computation of ideal quotients.

\medskip

\subsubsection*{The equations of an ideal}

Let $\RR$ be a Noetherian ring  and let $I$ be an
 ideal. By the {\em equations} of $I$ is meant a free presentation of
 the Rees algebra $\RR[It]$ of $I$,
\[ 0 \rar \LL \lar \SS = \RR[\TT_1, \ldots, \TT_m] \stackrel{\psi}{\lar}
\RR[It] \rar 0,  \quad \TT_i \mapsto f_it .\]
More precisely, we refer to the ideal $\LL$ as the {\em ideal of
equations} of $I$.
The ideal $\LL$ depends on the chosen set of generators of $I$, but all
of its significant cohomological properties, such as the integers that bound the
degrees of minimal generating sets of $L$, are independent of the presentation
$\psi$.
The examination of $\LL$ is one pathway to the
unveiling of
the properties of $\RR[It]$.
It codes the syzygies of all powers of $I$, and therefore
 is a carrier of not just algebraic properties of $I$, but of  analytic
ones
as well.
This is a notion that can be extended to Rees algebras $\Rees= \sum_{n\geq
0}I_nt^n$ of more general
filtrations. Among these  we will be interested in the integral
closure of $\RR[It]$.

\medskip

We are going to make a hit list of properties that would be desirable
to have about $\LL$, and discuss its significance.
Foremost among them  is to estimate the degrees of the generators of
$L$,
the so-called {\em relation type} of the ideal $I$: $\reltype(I)$.
If $\LL$ is generated by forms of degree $1$, $I$ is said to be of linear type
 (this is independent of the set of generators). The Rees algebra
 $\RR[It]$ is then the symmetric algebra $\Symi=\Sym(I)$ of $I$.
Such is the case when
 the $f_i$ form a regular sequence, $L$ is then generated by the
 Koszul forms $f_i\TT_j-f_j\TT_i$, $i<j$.
There are other kinds of kinds of sequences with this property (see
\cite{HSV1} for a discussion).

\medskip

The general theme we are going to pursue is the following: How to
bound the relation type of an ideal in terms of its numerical
invariants known {\it a priori}? Given the prevalence of various
multiplicities among these invariants most of our questions depend on
them.

\begin{enumerate}
\item[{\rm (a)}] Let $(\RR, \mathfrak{m})$ be a Noetherian local ring of dimension
$d$ and let $I$
be an $\mathfrak{m}$-primary ideal. Is there a low polynomial $f(e_0,e_1,
\ldots, e_d)$ on the Hilbert coefficients of $I$ such that
\[ \reltype(I)\leq f(e_0,\rme_1, \ldots, e_d)?\]

\item[{\rm (b)}] Among the equations in $\LL$, an important segment arises from
the {\em reductions} of $I$. That is, subideals $J\subset I$ with the
property that $I^{r+1}= JI^r$, for some $r$, the so-called {\em
reduction number} of $I$ relative to $J$; in notation, $\red_J(I)$.
 Minimal reductions are
important choices here and there are various bounds for the
corresponding values of $\red_J(I)$ ($\red(I)$ will be the minimum
value of $\red_J(I)$).  Among these formulas, for Cohen-Macaulay
local rings $(\RR, \mathfrak{m})$ of dimension $d$, and
$\mathfrak{m}$-primary ideal $I$, one has
\[ \red(I) \leq \frac{d}{o(I)} e_0(I)-2d+1,
\] where $o(I)$ is the largest integer $s$ such that $I\subset
\mathfrak{m}^s$. These formulas are treated in \cite[Section
2.2]{icbook}. One task is to derive such formulas for arbitrary
Noetherian local
rings. We plan to carry this out using {\em extended multiplicities}.
It will require a broadening of our understanding of related issues:
Sally modules of reductions.

\item[{\rm (c)}] If the ideal $I$ is not $\mathfrak{m}$-primary, which {\em
multiplicities} play the roles of the Hilbert coefficients?

\item[{\rm (d)}] A related question, originally raised by C. Huneke, attracted
considerable attention: whether for ideals $I$ generated by systems
of parameters, there is a uniform bound for $\reltype(I)$.
Although the answer was ultimately negative,
significant classes of rings do have this property. Wang obtained
several beautiful bounds for rings of low dimensions
(\cite{Wang97},\cite{Wang97a}), while Linh and Trung
(\cite{LinhTrung}) derived uniform bounds for all generalized
Cohen-Macaulay rings.

\item[{\rm (e)}]
Let $\mathcal{F}$ be the special fiber of $\Rees$,
that is $\mathcal{F}= \Rees\otimes \RR/\mathfrak{m}$. Since
\[\reltype(\mathcal{F})\leq \reltype(\Rees),\] it may be useful to
 constrain
 $\reltype(\mathcal{F})$ in terms of the Hilbert coefficients of the
 affine algebra $\mathcal{F}$. A rich vein to examine is that of
 ideals which are near complete intersections, more precisely in case
 of ideals of deviation at most two.
Among other questions, one should examine the problem of {\em
variation} of the Hilbert polynomials of ideals $I\subset I' \subset
\overline{I}$, where the latter is the integral closure of $I$. For
ideals of finite co-length, one of the questions (given the equality
$\rme_0(I)=e_0(I')$), is to study the change $\rme_1(I')-e_1(I)$ in terms
of the degrees of $\mathcal{F}$ and $\mathcal{F}'$.

%\item[{\rm (f)}] Almost all the previous questions can be asked for Rees
%algebras of modules or even of filtrations whose components are
%modules. The requisite multiplicities have already been built.

\end{enumerate}

We will treat mainly ideals of small deviations--particularly
  ideals of codimension $r$ generated
by $r+1$ or $r+2$ elements--because of their ubiquity in the
construction
of birational maps.
Often, but not
  exclusively, $I$ will be an ideal of finite
co-length in a local ring, or in a ring of polynomials over a field.

%Our focus on $\Rees$ is shaped by the following fact. The class of
%ideals $I$ to be considered will have the property that both  its
%symmetric algebra $\Sym(I)$ and the normalization $\Rees'$ of $\Rees$
%have amenable properties, for example, one of them (when not both) is
%Cohen-Macaulay. In such case,  the diagram
%\[ \Symi \twoheadrightarrow \Rees \subset \Rees'\]
%gives a convenient platform from which to examine $\Rees$.

\chapter{Degree/Multiplicity Functions}

\section*{Introduction}
\noindent
Let $\RR$ be a Noetherian ring.
Two of the basic structures defined over $\RR$
to be considered here are a module $M$ defined
by a free presentation
\[ \RR^m \stackrel{\varphi}{\lar} \RR^n \lar M \rar 0,\]
and an $\RR$-algebra $\AA$ given by generators and relations,
\[ \AA = \RR[\TT_1, \ldots, \TT_n]/(f_1, \ldots, f_m).\]

One of our goals is to extract numerical information from $\varphi$
and $(f_1, \ldots, f_m)$, and track it as $M$ and $\AA$ undergo
transformations of various kinds.
By a {\em degree} of a module $M$ or an algebra $\AA$ we mean a numerical measure of
information carried by $M$ or $\AA$. It must serve the purposes of allowing
comparisons between modules/algebras  and to exhibit flexible calculus rules
that track the degree
under some basic constructions, such as hyperplane sections,
normalization (in the case of algebras), number of generators of
cohomology modules, degrees of modules of syzygies, etc.

The premier example of a degree (vector space dimension excluded)
is that of the {\em multiplicity} of a module.
Let $(\SS, {\mathfrak m})$ be a local ring and let $M$ be a finitely
generated $\SS$--module. The Hilbert--Samuel\index{Hilbert-Samuel
function} function of $M$ is
\begin{eqnarray*}\label{hilbertf}
H_{M}: n \mapsto \lambda(M/{\mathfrak m}^{n+1}M),
\end{eqnarray*}  which for $n \gg 0$ is given
by a polynomial:
\[H_{M}(n) = \frac{\deg(M)}{d!}n^d + \mbox{\rm lower order terms}.\]

 The integer
 $d=\dim M$ is the  {\em dimension}\index{dimension
of a module} of $M$, and the integer
 $\deg(M)$ is its {\em
multiplicity}\index{multiplicity of a module}.
Under certain conditions, $\deg(M)$ can be interpreted as the {\em
volume} of a manifold, or of a polytope.

This {\em degree} arises in the setting of finitely generated
graded modules over graded $\RR$-algebra $\AA$, for an Artinian ring
$\RR$: From $\AA=\RR[z_1, \ldots, z_d]$ and $M = \bigoplus_{n\in \bbz} M_n $
(more general grading will be considered later), under fairly broad
conditions the function, referred as its {\em Hilbert function},
\[ H_M: n \mapsto \lambda(M_n)\]
provides a wealth of information, with a {\em degree} coding a great deal
of it.

In reality the Hilbert functions mentioned above are distinct, but
are closely related, that is, we will call  either $H_M(n)$ or
\[ \sum_{k\leq n} H_M(k)\]
as the Hilbert function of $M$, or refer to the second kind as the
Hilbert--Samuel function.

This however is a very narrow view of what is desirable. It is
important to obtain {\em degree like} information from $M$ even when
$\RR$ is not Artinian, by which one means shared properties relative to
certain exact or hyperplane sections. Another aim is to define
functions like $H_M$ when $M$ is not necessarily a Noetherian module.
For instance, if $(\RR,\mathfrak{m})$ is a Noetherian local ring and
$J\subset I$ are $\mathfrak{m}$-primary ideals, the module whose
components are
\[ M_n = I^n/J^n\]
has a nice Hilbert function, since
\[ \lambda(M_n)= \lambda(\RR/J^n)-\lambda(\RR/I^n),\]
a difference of ordinary Hilbert functions.

To give an example of the challenge we will undertake, consider a
case inspired by the {\em Buchsbaum-Rim multiplicity}. Let $\RR$ be a
local Noetherian ring,
 let $M$ be a
submodule of the free $\RR$-module $F=\RR^r$ such that
$\lambda(F/M)<\infty$. Setting $\SS=\RR[x_1, \ldots, x_d]=\bigoplus_n S_n
$, let $M^n\subset S_n$ denote the $\RR$-submodule generated by the
products of $n$ elements of $M$. The function
$ \lambda(S_n/M^n)$
has properties similar to (\ref{hilbertf}) despite the fact that
unlike the former it lacks a {\em direct} interpretation as a Hilbert
function of a Noetherian graded module. (It should be said though that
there are {\em indirect} ones.) The challenge will consist in
attaching a degree when the $M$ is an arbitrary submodule of $F$.

%The classical case, when $\RR$ is a field and
%$\AA$
%is a standard graded algebra, say $\AA=k[x_1, \ldots, x_n]/I$,  has an
% important property of this degree regarding its computability. If $<$
%is a term ordering for the ring of polynomials, then by a theorem of
%Macaulay
%\[\deg(\AA)= \deg(k[x_1, \ldots, x_n]/in(I)),\]
%where $in(I)$ is the corresponding initial ideal of $I$.
%There is a similar observation with regards to $\AA$-modules.

One of the paths we employ to obtain new degrees from old ones is the
following. Let $\RR$ be a Noetherian ring and $\AA$ a finitely
generated graded $\RR$-algebra. Usually we will require that $\AA$ be a
standard graded algebra\index{standard graded algebra}, but later we
shall discuss more general  gradings.
Suppose we are equipped with a numerical function $\lambda$ on  the category
of finitely generated $\RR$-modules. Let $M=\bigoplus_{n\in \bbz} M_n $ be a finitely generated
graded $\AA$-module, to each we associate the Laurent  power series
\[ [\![M]\!]= \bbp_M(t) = \sum_{n\in \bbz} \lambda(M_n)t^n.\]
Under some conditions, this function is a rational function,
\[ [\![M]\!] = \frac{h_M(t)}{(1-t)^{d}},\]
$h_M(t) \in \bbz[t,t^{-1}]$. From this representation is extracted
the Hilbert polynomial
\[P_M(n) = \sum_{i=0}^{d-1} (-1)^i\rme_i(M)
{{n+d-i-1}\choose{d-i-1}}, \]
with the property that $\lambda(M_n)= P_M(n)$ for $n\gg 0$.

There is an obvious degradation in the quality of information about
$M$ when we move from $\bbp_M(\ttt)$ to one of the individual
coefficients of $P_M(\ttt)$
\[ \bbp_M(\ttt) \rar P_M(\ttt)\rar \{ \rme_0(M),\rme_1(M), \ldots, \rme_{d-1}(M)\}
\rar\rme_i(M).\]
Nevertheless it will be   the $\rme_i(M)$, particularly on $\rme_0(M)$
and $\rme_1(M)$, that we look for carriers of information on $M$. The
reasons are both for the practicality it involves but also justified
by the surprising amount of information these coefficients may pack. 
Motivated by the celebrated theorem of Serre on multiplicities expressed by
Euler characteristics (cf. \cite[Theorems 4.7.4, 4.7.6]{BH}), there is a great deal
of attention given to the relationships between $\e_1(M)$ and the Cohen-Macaulayness
of $M$.

\medskip

In the study of Hilbert functions associated to a filtration $\mathcal{M}=\{
M_n, n\in \bbz\}$ of a module $M$, it is
common to consider the (integral) function
\[ H_{\mathcal{M}}(n) = \lambda(M/M_{n+1}).\]
There will be (under appropriate finiteness conditions) a rational
function
\[ \bbp_{\mathcal{M}}(\ttt) = \frac{h_{\mathcal{M}}(\ttt)}{(1-\ttt)^{d+1}},\]
$h_{\mathcal{M}}(\ttt) \in \bbz[\ttt,\ttt^{-1}]$. From this representation is extracted
the Hilbert polynomial
\[P_{\mathcal{M}}(n) = \sum_{i=0}^{d} (-1)^i\rme_i(M)
{{n+d-i}\choose{d-i}}, \]o
with the property that $\lambda(M/M_{n+1})= P_{\mathcal{M}}(n)$ for $n\gg 0$.

\medskip

We introduce a word of caution in the definitions of Hilbert functions associated to a filtration 
$\mathcal{M} = \{M_n, n \geq 0\}$.  If $\lambda$ is a numerical function with properties that mimic
those of a length function on Artinian modules, the core definition of the Hilbert function
\[ n \mapsto \lambda(M_n/M_{n+1})\]
can be supplemented by the Hilbert functions
\[ n \mapsto \sum_{k\leq n} \lambda(M_k/M_{k+1})\]
and 
\[ n \mapsto \lambda(M/M_{n+1}).\]
The latter may not be equal as the function may not be additive on short exact sequences. We will meet such 
functions when we discuss cohomological degrees. Under appropriate conditions--obviously when
the quotients $M/M_k$ are Artinian but  also when $\mathcal{M}$ is an adic-filtration and
$\lambda(\cdot)$ is the classical multiplicity--we have equivalencies in these definitions.

\medskip

\subsubsection*{Primary decomposition}
The function $\deg(\cdot)$ may fail to capture features which
are significant for $M$. For instance, if $M$ is filtered by a
chain of submodules
\begin{eqnarray*}
M = M_1\supset M_2 \supset \cdots \supset M_r\supset 0,\quad
M_i/M_{i+1}\simeq \SS/{\mathfrak p}_i,
\end{eqnarray*}
 where
${\mathfrak p}_i $ is a prime ideal,
then \[\deg(M) = \sum\limits_{\dim \SS/{\mathfrak p}_i =\dim M}\deg(\SS/{\mathfrak p}_i ).\]
Alternatively, if ${\mathfrak p} $ is an associated prime
of $M$ of dimension $\dim \SS/{\mathfrak p}= \dim M$, the localization $M_{\mathfrak p} $ is an
Artinian module whose length $\lambda(M_{\mathfrak p})$ is the number
of times ${\mathfrak p} $ occurs in (\ref{filt1}), then
\[\deg(M) = \sum\limits_{\dim \SS/{\mathfrak p} =\dim M}
\lambda(M_{\mathfrak p})\deg(\SS/{\mathfrak p} ).\]

This formulation still ignores  the contributions of the lower
dimensional components. Partly to address this, other  degree
functions have been
defined that collects information about the associated primes of $M$
in all codimensions.
What kind of information is this? One of these indices
 is the following: If $\mathfrak{p}\in
 \Spec(\SS)$, the module $H^0_{\mathfrak{p}}(M_{\mathfrak{p}})$ is an
 Artinian $\RR_{\mathfrak{p}}$-module. Its length,
\[ \mult_M(\mathfrak{p}):=
\lambda(\H^0_{\mathfrak{p}}(M_{\mathfrak{p}})),
\] is called the {\em multiplicity length} $M$ at $\mathfrak{p}$.
Since $\mult_M(\mathfrak{p})\neq 0$ if and only if $\mathfrak{p}\in
\Ass(M)$, its value is a good measure of the weight of $\mathfrak{p}$
as an associated prime of $M$.
It can be used directly as follows.

\medskip

%\begin{Definition}[\cite{BM}]{\rm
 The {\em arithmetic degree}\index{arithmetic
degree}
 of $M$ is the integer
\[\adeg(M) = \sum\limits_{{\mathfrak p}{ \in \Ass(M)}} \mult_M({\mathfrak p})\cdot \deg
\SS/{\mathfrak p}=\sum\limits_{{\mathfrak p}} \mult_M({\mathfrak p})\cdot \deg
\SS/{\mathfrak p}.\]
%}\end{Definition}

\medskip

A related degree is $\gdeg(M)$, the {\em geometric
degree}\index{geometric degree} of $M$,
 where in
this sum one adds only the terms  corresponding to minimal
associated primes of $M$.

One feature of these degrees is that they can be expressed in a way
which does not require the actual knowledge of the associated primes
of $M$. Both $\adeg(M)$ and $\gdeg(M)$ are put together from a sum of
multiplicities of certain modules derived from $M$ in a way that does
not require the values of  $\mult_M(\mathfrak{p})$.

Another type of arithmetic degree will be defined if $\SS$ is a
standard graded algebra  over a Noetherian ring $\RR$, not
necessarily local. Using the notion of {\em $j$-multiplicity} of
Flenner, O'Carroll and Vogel (\cite{FOV}), one attaches a new degree,
$\jdeg$, to any finitely generated graded $\SS$-module $M$ (\cite{jdeg1}).
It arises as follows. For any prime $\mathfrak{p}\in \Spec (\RR)$, the
module
 $H=\H^0_{\mathfrak{p}}(M_{\mathfrak{p}})$ is a finitely generated
 graded module over a standard graded ring over an Artin ring. It
 makes sense to consider the multiplicity
\[ \deg(\H^0_{\mathfrak{p}}(M_{\mathfrak{p}})).\]
If the dimension of $H$ is $\dim M_{\mathfrak{p}}$, its multiplicity
is called the $j$-multiplicity of $M$ at $\mathfrak{p}$. It is
denoted by $j_{\mathfrak{p}}(M)$. For completeness, one sets
$j_{\mathfrak{p}}(M)=0$
 if $\dim H< \dim (M_{\mathfrak{p}})$.
The $\jdeg$ of $M$ is defined by the (finite) sum
\[ \jdeg_{\RR}(M)
:=\sum_{\mathfrak{p}\in \Spec \RR} j_\mathfrak{p}(M).\]
It is defined globally and will have sensitivities not shared by the
other
 degrees derived from the
classical multiplicity.

\medskip

Let $(\RR, \mathfrak{m})$ be a Noetherian local ring (or a Noetherian
graded algebra) and let $\mathcal{M}(\RR)$ be the category of finitely
generated $\RR$-modules (or the appropriate category of graded
modules). A {\em degree function}\index{degree function} is a
numerical function $\mathbf{d}: \mathcal{M}(\RR) \mapsto \mathbb{N}$.
The more interesting of them initialize on modules of finite length
and have mechanisms that control how the function behaves under
generic hyperplane sections. Thus, for example, if $L$ is a given
module of finite length, one function may  require that
\[ \mathbf{d}(L) := \lambda(L),\]
or when $L$  a graded module
\[ L = \bigoplus_{i\in \mathbb{Z}} L_i,\] the requirement is that
\[ \mathbf{d}(L) := \sup\{ i \mid L_i\neq 0\}.\]

In this Chapter we study some of these functions,
particularly those that are not treated extensively in the
literature. One of our aims is to provide an even
 broader approach to Hilbert functions that does not
depend so much on the base rings being Artinian. We shall be looking at
the associated multiplicities and seek ways to relate them to
the reduction numbers of the algebras. Among the {\em multiplicities}
treated here and in the next section are:

\begin{itemize}
\item[{$\bullet$}] Classical multiplicity\index{classical multiplicity}
\item[{$\bullet$}] Arithmetic degree and geometric degrees
\item[{$\bullet$}]  $\jdeg$
\item[{$\bullet$}] Castelnuovo-Mumford
regularity
\item[{$\bullet$}] Cohomological degree
\item[{$\bullet$}] Homological degree
\item[{$\bullet$}] bdeg
\end{itemize}

\[
\diagram
 & \mbox{\rm Degree\ families} \dlline \dline \drline &  \\
multiplicity\dline & Castelnuovo \ regularity & extended \ degrees\dline \\
arithmetic   \ degrees \dline & & homological \ degree \& bdeg  \\
 adeg,\  gdeg, \ jdeg & & \\
\enddiagram
\]

\section{Hilbert Functions}

\subsection{Introduction}

Let $\RR$ be a Noetherian ring, $\SS=\RR[x_1, \ldots, x_n]$ a standard
graded algebra and $M$ a finitely generated graded $\SS$-module. We
will pursue  approaches to extract numerical invariants of the
module $M$ from various Hilbert functions attached to the module.

 \medskip

There are two prototypes here.
Let $(\RR,\mathfrak{m})$ be an Artinian local ring, $\SS=\RR[x_1, \ldots,
x_n]$ a ring of polynomials, graded by setting $\deg x_i=1$, and let
$M$ be a finitely generated $\SS$-module. The {\em Hilbert function}
 of $M$
is the assignment
\[ H_M(n):= \lambda(M_n).\]
It is coded as a Laurent series with a rational representation called
the {\em Hilbert-Poincar\'e series} of $M$
\[ \bbp_M(\ttt)= \sum_{n\in \bbz} H_M(n)\ttt^n = \frac{h_M(\ttt)}{(1-\ttt)^d}\in
\bbz[[\ttt,\ttt^{-1}]],\]
where $h_M(\ttt)$ is a Laurent polynomial called the {\em
$h$-polynomial} of $M$. If $d$ is chosen so that $h_M(1)\neq 0$, then
$d$ is the Krull dimension of the module $M$, $d=\dim M$.

\medskip

For values of $n> \deg h_M(\ttt)-d$, the Hilbert function is given by a
polynomial $P_M$ of degree $d-1$
\[ H_M(n) =P_M(n) = \sum_{i=0}^{d-1} (-1)^i\rme_i(M)
{{n+d-i-1}\choose{d-i-1}}. \]

\begin{Definition} \index{a-invariant of a module}{\rm
 The integer $a(M)=\deg h_M(\ttt) - \dim M$ is called
the $a$--invariant of $M$.
}\end{Definition}

The first iterated Hilbert function of $M$,
\[ H^1_M(n) := \sum_{i\leq n} \lambda(M_i), \]
carries some additional information in the corresponding Hilbert
polynomial
\[ H_M^1(n) =P_M^1(n) = \sum_{i=0}^{d} (-1)^i\rme_i(M)
{{n+d-i}\choose{d-i}}. \]

The coefficients $\rme_i(M)$ are labeled the {\em Hilbert coefficients}
of $M$. The leading
 coefficient, $\rme_0(M),$ is the {\em multiplicity} of $M$, and it is one of
its most important invariants.
 The next coefficient, $\rme_1(M)$, under certain conditions, is called the {\em
Chern number} of $M$\index{Chern coefficient}. It plays an ubiquitous role in the complexity
of $M$. All these coefficients can be obtained from the
$h$-polynomial $h_M(\ttt)$,
\[\rme_i(M) = \frac{h_M^{(i)}(1)}{i!}, \quad i = 0,\ldots, d.\]

The other Hilbert function is that associated to  {\em good
filtrations} over a local ring $(\RR, \mathfrak{m})$. Let $M$ be a
finitely generated $\RR$-module and let $\mathfrak{q}$ be an
$\mathfrak{m}$-primary ideal The collection of subsets of $M$, $\mathcal{M}=\{ \mathfrak{q}^nM,
\ n\geq 0 \}$,  is the $\mathfrak{q}$-{\em adic filtration}  of $M$.
More generally, say that a decreasing collection of subsets of $M$,
$\mathcal{F}= \{ M_n, \ n\geq 0\}$, is a $\mathfrak{q}$-good
filtration if $M_0=M$, $\mathfrak{q}M_{n}\subset M_{n+1}$,  and
$M_{r+1}= \mathfrak{q}M_r$ for  $r\gg 0$.
There is a graded module
attached to $\{ M, \mathcal{M}\}$,
\[ \gr_{\mathcal{F}}(M):= \bigoplus_{n\geq 0} M_n/M_{n+1},\]
called the {\em associated graded module}.

The module $\gr_{\mathcal{F}}(M)$ has a Hilbert function in the
manner defined above. These filtrations occur often in ordinary
operations with modules. For example, given the
$\mathfrak{q}$-filtration of $M$, and a submodule $N\subset M$,
the submodules
$\{ N\cap \mathfrak{q}^nM\}$ define a good filtration on $N$.

\begin{Definition}{\rm Let $(\RR, \m)$ be a Noetherian local ring
and.
$M$ a finitely generated $\RR$-module. The Hilbert coefficients
associated to the $\m$-adic filtration will be denoted by
$\rme_0(M)=\deg(M)$, $\rme_1(M)$, etc.
}\end{Definition}

\subsubsection*{Comparison of Hilbert functions}

Let $\AA$ be a standard graded algebra over an Artin ring   and let
$M$ a finitely generated graded $\AA$-module. Suppose $I$ is an ideal
of finite co-length,
generated by $1$-forms. We seek comparisons between the Hilbert
functions
\[ \lambda(M_n), \quad \lambda(M/I^{n+1}M),\]
particularly about the information they share. A case in mind is when
$\AA$ is a polynomial ring over a field and $I$ is its maximal
homogeneous ideal.

%Another
%way it occurs is when on a Rees algebra of an ideal $I$, one
%carries out operations to determine its integral closure,
%\[ R[It] \lar \bar{R[It]}=\AA = \sum_{n\geq 0} A_nt^n.\]
%The process

%\begin{itemize}

%\item[{$\square$}] Macaulay, Green, Stanley
%\item[{$\square$}] Hilbert Coefficients
%\item[{$\square$}]
%\end{itemize}
\medskip

\subsubsection*{Change of rings} Let $(\RR, \m)\rar (\SS,
\n)$\index{change of rings and multiplicity}
be a local homomorphism of rings of the same Krull dimension and such
that $\lambda(\SS/\m \SS)< \infty$. For any $\m$-primary ideal $I$,
$I\SS$ is $\n$-primary.

\begin{Proposition} \label{degandcr}
If $M$ is a finitely generated $\RR$-module,
the Hilbert polynomials associated to the filtrations  $I^nM$ of $M$
and $I^n M\otimes_{\RR}\SS$ have the same degree, and their
multiplicities satisfy
\[ \deg_{I\SS}(M\otimes_{\RR}\SS) \leq \deg_I(M)\cdot \lambda(\SS/\m
\SS).\]
Equality holds in case $\SS$ is $\RR$--flat, such as the case when
$\SS$ is the $\m$-adic  completion of
$\RR$, or $\SS=\RR[x]_{\m \RR[x]}$. In both of these cases the
Hilbert functions coincide.
\end{Proposition}

\subsubsection*{Syzygies}
Two general methods to obtain the Hilbert
functions of modules or algebras work by indirection: Syzygies and
Gr\"{o}bner bases.

\medskip

One of the most used methods to obtain Hilbert functions of graded
modules  and their
associated polynomials is via resolutions.
Let $k$ be a field, $\SS= k[x_1,\ldots,x_d]$ a polynomial ring over $k$,
 and $M$ a finitely generated
graded $\SS$-module.
Then $M$, as an $\SS$-module, admits a finite graded free resolution:
\[
0\rar\bigoplus_j\SS[-j]^{b_{pj}}\lar\cdots\lar
\bigoplus_j\SS[-j]^{b_{0j}}\lar M\rar 0.
\]
 The Hilbert-Poincar\'e series of $M$ is
\[ \bbp_M(t)= \frac{\sum_{jp} (-1)^j t^{b_{pj}}}{(1-t)^d}.\]

\subsubsection*{Macaulay Theorem}
The classical case, when $\RR$ is a field and
$\AA$
is a standard graded algebra, say $\AA=k[x_1, \ldots, x_n]/I$,  has an
 important property of this degree regarding its computability.

\begin{Theorem}\index{Macaulay theorem}
If $<$
is a term ordering for the ring of polynomials, then  the algebras $\AA$ and
$k[x_1, \ldots, x_n]/in(I))$,
where $in(I)$ is the corresponding initial ideal of $I$, have the
same Hilbert functions.
There is a similar assertion with regard to $\AA$-modules.
\end{Theorem}

\subsection{General gradings}

\section{The Classical Multiplicity for Local Rings}

\subsection{Introduction}
Let $(\RR, {\mathfrak m})$ be a Noetherian local ring and let $M$ be a finitely
generated $\RR$--module. The basic multiplicity function arises as
follows.

\begin{definition}{\rm  Consider
the Hilbert function of $M$ is
\[H_{M}: n \mapsto \lambda(M/{\mathfrak m}^{n+1}M),\]  which for $n \gg 0$ is given
by a polynomial:
\[H_{M}(n) = \frac{\deg(M)}{d!}n^d + \mbox{\rm lower order terms}.\]
 The integer $\deg(M)$ is the {\em
multiplicity}\index{multiplicity of a module} of the module $M$.
$($The integer $d=\dim M$ is its {\em dimension}$)$.\index{dimension
of a module}
}\end{definition}

\subsection{Rules of computation} Let us derive some rules
that govern the computation of multiplicities.

\bigskip

\subsubsection*{Exact sequences}
 The behavior
of $\deg(M)$ and $\dim M$ with regard to submodules and primary
decomposition  follows from:

\begin{Proposition} Let $\RR$ be a local ring and
let
\[ 0 \rar L \lar M \lar  N \rar 0\]
be an exact sequence of $\RR$--modules. Then
\[ \deg(M)= \left\{  \begin{array}{ll}
 \deg(L)+ \deg(N) & \mbox{\rm if \ } \dim L = \dim N = \dim M,\\
 \deg(L) & \mbox{\rm if \ }  \dim N < \dim M,\\
 \deg(N) & \mbox{\rm if \ } \dim L < \dim M.
\end{array} \right. \]
\end{Proposition}

\begin{proof}
Consider the induced exact sequence
\[ 0 \rar L/L\cap \mathfrak{m}^{n+1} M \lar M/\mathfrak{m}^{n+1} M \lar
N/\mathfrak{m}^{n+1} N \rar 0.
\]
The assertions will follow from direct comparison of the Hilbert
functions of $N/\mathfrak{m}^{n+1} N$ and $N/N\cap\mathfrak{m}^{n+1} M$: It
will suffice to show that the
 module whose components are   on the left of the exact sequence
\[ 0 \rar N\cap\mathfrak{m}^{n+1} M/ \mathfrak{m}^{n+1} N
\lar   N/\mathfrak{m}^{n+1} N \lar N/N\cap\mathfrak{m}^{n+1} M\rar 0,
\]

\[G=\bigoplus_{n\geq 0}N\cap\mathfrak{m}^{n+1} M/ \mathfrak{m}^{n+1} N,\]
has dimension $<\dim N$. By the Artin-Rees Lemma, there is an integer
such that for $n\geq r$,
\[ N\cap\mathfrak{m}^{n+1} M= \mathfrak{m}^{n+1-r}(\mathfrak{m}^rM \cap N)
\subset \mathfrak{m}^{n+1-r}N.
\]
In particular, the Hilbert polynomial of $G$ has degree at most that
of  the module $G'$ whose components is a sum
\[ \mathfrak{m}^{n}N/
\mathfrak{m}^{n+2}N \oplus \cdots \oplus \mathfrak{m}^{n-r}N/
\mathfrak{m}^{n+2-r}N.
\] The Hilbert polynomial of $G'$ has degree $\dim N-1$ (and
multiplicity $r\cdot \deg(N)$). 
\end{proof}

\bigskip

Similar rules apply to more general filtrations. Let us consider one
of these.

\medskip

Let $(\RR, \mathfrak{m})$ be a Noetherian local
ring, let $I$ be an $\mathfrak{m}$-primary ideal
 and   let $\mathcal{A}=\{I_n, n\geq 0\}$ be an $I$-good filtration.
For a finitely generated $\RR$-module $M$, denote
  by $\rme_i(M)$ the Hilbert coefficients of $M$
for the  filtration $\mathcal{A}M= \{ I_nM, n\geq 0 \}$.
The function
\[H_{M}: n \mapsto \lambda(M/I_{n+1}M),\]   for $n \gg 0$ is given
by a polynomial:
\[H_{M}(n) = \frac{\deg_{\mathcal{A}}(M)}{d!}n^d + \mbox{\rm lower order terms}.\]

\begin{definition}{\rm The integer $\deg_{\mathcal{A}}(M)$ is
called the {\em Hilbert-Samuel multiplicity} of $M$ relative to
$I$.\index{$\deg_I$, the Samuel-Hilbert multiplicity}
It will be denoted also by $\deg_I(M)$, and at times even by
$\deg(M)$ if no confusion arises.
}\end{definition}

\begin{Proposition} \label{abc}
Let
\[ 0 \rar A \lar B \lar C \rar 0\]
be an exact sequence of finitely generated $\RR$-modules.
If $r=\dim A< s= \dim B$, then $\rme_i(B)=\rme_i(C)$ for $i< s-r$.
\end{Proposition}

\begin{proof} From the diagram
\[
\diagram
0 \rto  & A \cap I_nB \rto\dto & I_nB \rto \dto & I_nC
\rto\dto  &0 \\
0 \rto          &  A \rto                 & B \rto  & C \rto  &0
\enddiagram
\]
we have the equality
\[ \lambda(B/I_nB)= \lambda(A/I_nB\cap A) + \lambda(C/I_nC).
\] Since $\lambda(A/I_nA)\geq \lambda(A/I_nB\cap A)$, the Hilbert
polynomials of $B$ and $C$ have the same coefficients in the
specified range. 
\end{proof}

%hello
\bigskip

Exact sequences
\[ 0 \rar A \lar B \lar C \rar 0\]
where $\dim C < \dim B$ give rise to slightly different equalities.
We leave the proof for the reader.

\begin{example}\label{degabc}{\rm Let $\RR$ be a Cohen-Macaulay local
ring and let $\{ x_1, \ldots, x_n\}$, $n\geq 2$,  be a set of elements such that
any pair forms a regular sequence. Set
\begin{eqnarray*}
 \xx & = & x_1\cdots x_n,\\
z_i & = & \xx/x_i, \quad i=1\ldots n.
\end{eqnarray*}
We claim that
\begin{eqnarray} \label{degxx}
 \deg(\RR/(z_1, \ldots, z_n)) \leq \frac{1}{2} ((\deg(\RR/(\xx)))^2-n).
\end{eqnarray}
We argue by induction on $n$, the formula being clear for $n=2$.

Consider the exact sequence
\[ 0 \rar (x_1, z_2, \ldots, z_n )/(z_1, z_2, \ldots, z_n) \lar
\RR/(z_1, z_1) \lar \RR/(x_1,z_1) \rar 0.\]
Since
\[ (x_1,z_1)/(z_1, \ldots, z_n)\simeq \RR/(z_2', \ldots, z_n'),
\] where $z_i'$, $i\geq 2$, denotes the products from elements in the set $\{x_2,
\ldots, x_n\}$ using the formation rule of the $z_i$.

Adding the multiplicities of the modules of the same dimension, we
have
\[ \deg(\RR/(z_1, \ldots, z_n)) = \deg(\RR/(z_1,\ldots, z_n)) + \deg(\RR/(z_2',
\ldots, z_n')). \]
As
 \[\deg(\RR/(x_1, z_1))= \deg(\RR/(x_1))\cdot \deg(\RR/(z_1))=
\deg(\RR/(x_1))\cdot \sum_{j\geq 2}\deg(\RR/(x_j)) ,\] and by
induction
\[ \deg(\RR/(z_2', \ldots, z_n'))= \sum_{2\leq i< j\leq
n}\deg(\RR/(x_i))\cdot \deg(\RR/(x_j)),\]
we have
\[ \deg(\RR/(z_1, \ldots, z_n)) = \sum_{1\leq i< j \leq n}
\deg(\RR/(x_i))\cdot \deg(\RR/(x_j)).
\]

The rest of the calculation is clear.
There are similar formulas in case every subset of $k$ elements of
$\{x_1, \ldots, x_n\}$ forms a regular sequence.

}\end{example}

\subsubsection*{Hyperplane section}
Let $\AA$ be a standard graded algebra over the Artin local ring
$(\RR,\mathfrak{m})$,
\[ \AA = \RR[\AA_1] = \RR + \AA_1 + \cdots + \AA_n + \cdots.\]
We denote by $P= (\mathfrak{m}, \AA_1)$. This ideal play many roles in
the study of the Hilbert functions of graded $\AA$-modules.
We examine one aspect grounded on the classical notion of {\em
generic hyperplane section}.\index{generic hyperplane section}

\medskip

Let $M = \bigoplus_{n}M_n$ be a finitely generated graded
$\AA$-module, and denote by
\[ \mathcal{P}=\{P_1, \ldots, P_m\} = \Ass_{\AA}(M)\]
the set of associated primes. Let $h$ be a homogeneous element, $h\in
A_s$ that is not contained in any $P_i\neq P$. Consider the submodule
\[ \ann_M(h) = \{x\in M \mid hx = 0\}. \]
If $\ann_M(h)\neq 0$, its only associated prime ideal is $P$, and
therefore $P^r\cdot \ann_M(h)=0$, for $r\gg 0$, which means that
$\ann_M(h)$ has only a finite number of graded components, and is
therefore a module of finite length.

\medskip

Converted into Hilbert functions, this property gives rise to
an exact sequence of graded modules,
\[ 0 \rar \ann_M(h)[-s] \lar M[-s] \lar M \lar M/hM \rar 0,\]
with
\[ H_M(n)-H_M(n-s)= H_{M/hM}(n), \quad n\gg 0. \]
In particular, $\dim M/hM = \dim M-1$.

\medskip

There are variations of this behavior, we recall two of them. Let
$\RR$ be a Noetherian ring $\q$ an $\RR$-ideal. Let $M$ be a finitely
generated $\RR$-module equipped with a decreasing  $\q$-good
\index{good filtration} filtration
of submodules, $\mathcal{F}=\{M_i\}$.

\begin{definition}\index{superficial element}{\rm  The element $x\in
\q$ is {\em superficial} for $\mathcal{F}$--or by abuse of
terminology for $M$--if there is an integer $c$ such that
\[ (M_{n+1}:_M x)\cap M_c=M_n
\] for all $n\geq c$.

}\end{definition}

This implies that  multiplication by $\bar{x}\in \q/\q^2$ on
$\gr(\mathcal{F})$ is injective in degrees $n\geq c$. The existence
of such elements usually require  some cardinality condition on
the units of $\RR/\q$. Let us give the fuller characterization of
these elements from \cite{RV10}.

\medskip

The setup is the following.

\begin{itemize}

\item Set $G=\gr_{\mathcal{F}}(M)$,
 $G_0= \bigoplus_{j\geq 0} \q^j/\q^{j+1}$, $Q=\bigoplus_{j\geq 1}
 \q^j/\q^{j+1}$.
\item Let $\Ass(M)= \{\p_1, \ldots, \p_s\}$.
\item Let $(0_M)= P_1\cap \cdots \cap P_s$ be an irreducible
decomposition of the trivial submodule of $M$,
$\p_i = \sqrt{0:
M/P_i}$. Suppose $\q\not \subset \p_i$, for $i=1, \ldots, r$,
$\q \subset \p_i $, $i=r+1, \ldots, s$.

\item Setting $N=P_1\cap \cdots\cap P_r$, $\Ass(M/N)= \{\p_1, \ldots,
\p_r\}$, $N\cap M_j=0$ for $j\gg 0$.

\item Then $N=P_1\cap \cdots \cap P_r=\{ z\in M\mid \q^n z=0, n\gg
0\}$.

\item Let $(0_G)= T_1\cap \cdots \cap T_u\cap T_{u+1}\cap \cdots \cap
T_v$ be an irreducible decomposition of the trivial submodule of $G$
as a $G_0$-module.

\item Set $\Ass_{G_0}(G)= \{\Q_1, \ldots, \Q_v\}$,
$\Q_i = \sqrt{0:
G/T_i}$. Suppose $Q\not \subset \Q_i$, for $i=1, \ldots, u$,
$Q \subset \Q_i $, $i=u+1, \ldots, v$.

\item Then $H=T_1\cap \cdots \cap T_u=\{ z\in G\mid Q^n z=0, n\gg
0\}$.

\end{itemize}

\begin{Theorem} \label{supelt} Let $x\in \q$, and set $x^*= x +
\q/\q^2$.
The following conditions are equivalent:

\begin{enumerate}
\item[{\rm (a)}] $x$ is $M$-superficial for $\q$;

\item[{\rm (b)}] $x^*\notin \bigcup_{i=1}^m \mathfrak{B}_i$;

\item[{\rm (c)}] $H:_{\gr(\mathcal{F})}x^*=H$;

\item[{\rm (d)}] $N:x=N$ and $M_{n+1}\cap x M= M_n$ for all large $n$;

\item[{\rm (e)}] $(0:_{\gr(\mathcal{F})} x^*)_n=0$ for all large $n$;
\item $M_{n+1}:x=M_n + (0:_M x)$ and
$M_n \cap (0:_M x)=0$ for all large $n$.

\end{enumerate}
\end{Theorem}

We shall refer to $x$ as an $M$-{\em filter regular} element of degree
$1$.

\begin{Proposition} \label{filtreg} Let $\AA$ and $M$ be as above. If
$M\neq 0$, for some $s>0$ there exists an $M$-filter regular element
of degree $s$. Moreover, if $\RR$ has infinite residue field, $s$ can
be arbitrarily chosen.
\end{Proposition}

When the residue field of $\RR$ is infinite, being filter regular is a
generic property of the elements in $A_s\setminus \
\mathfrak{m}A_s$.

\medskip

One of the applications of this notion is to associated graded rings
and modules with good filtrations.
Let $(\RR, \mathfrak{m})$ be a Noetherian local ring, and  let $M$ be a
finitely generated $\RR$-module and let $\mathfrak{q}$ be an
$\mathfrak{m}$-primary ideal. Let
 $\mathcal{M}=\{ \mathfrak{q}^{n+1}M,
\ n\geq 0 \}$ be a   $\mathfrak{q}$-good   filtration  of the finitely
generated $\RR$-module $M$.
Let $\AA = \gr_\mathfrak{q}(\RR)$ and set
\[ \BB=\gr_{\mathcal{M}}(M):= \bigoplus_{n\geq 0} M_n/M_{n+1},\]
for the  { associated graded module} of $\mathcal{M}$.

Let $h\in \mathfrak{q}$. Let $\mathcal{M}'$ be the filtration induced
$M_n+hM/hM$, and denote by $\CC$ its associated graded module. We are
going to compare $\BB$ to $\CC$.

\begin{Proposition} Suppose $(\RR, \mathfrak{m})$ be a Noetherian local
ring, and $\BB$ and $\CC$ are as above. Assume that $h\in \mathfrak{q}$ is
filter regular relative to $\BB$. Then there is a surjection
\[ 0\rar L \lar \bar{\BB} = \BB\otimes \RR/(h) \lar \CC \rar 0\]
such that $L_n=0$ for $n\gg 0$.
\end{Proposition}

\begin{proof} The mapping $\bar{\BB}\rar \CC$, in degree $n$, is given by the
natural homomorphism
\[     M_n/ (M_{n+1} + hM_{n-1})\lar (M_n+hM)/(M_{n+1}+ hM) \rar 0. \]

Denote by $h^*$ the image of $h\in \AA_1=
\mathfrak{q}/\mathfrak{q}^2$.
By the choice of $h$, there exists an integer $c$
 such that $\ann_{\BB}(h^*)\cap \BB_n=
0$ for $n>c$.  It is easy to verify that this implies that for all
$n>c$, $(M_n:_M h)\cap M_c = M_{n-1}$, which will show that $L_n=0$,
as desired. 
\end{proof}

\begin{Corollary} Let $M$ and $h$ be as above.
\begin{enumerate}
\item[{\rm (a)}]  If $h^*\in \AA_1=\mathfrak{q}/\mathfrak{q}^2$ is regular on
$\BB$, then $\CC=\bar{\BB}$.

\item[{\rm (b)}] If $\mathfrak{q}$ is $\mathfrak{m}$-primary, the Hilbert
coefficients of $\BB$ and $\CC$ satisfy
\[\rme_i(\BB)=\rme_i(\CC), \quad 0\leq i< \dim M.\]
\end{enumerate}
\end{Corollary}

We shall now treat a remarkable property of $\BB$ discovered by J.
Sally and finally clarified in \cite[Lemma 2.2]{HucMar97}.

\begin{Proposition}\label{SallyMachine} Let $M$ and $h$ be as above.
If $\depth \CC>0$, then $h^*\in \AA_1$ is a regular element on $\BB$.

\end{Proposition}

\begin{proof} Since $\ann(h^*)\cdot \AA_{+}^n=0$, for $n\gg 0$, it suffices
by the previous Corollary to prove that $\depth \AA>0$.
Let $y$ be an element in $\mathfrak{q}^t$
such that the image of $y$ in $\AA_t$ is  a regular element on $\CC$.
  Then \[ M_{n+tj}:y^t= M_n+hM, \quad \forall n,j.\]
 Since $h$ is
filter regular,  there exists an integer c such that
$(M_{n+j} : h^j)\cap  M_c\subset M_n$ for all
$j\geq 1$ and $n>c$. Let $n$ and $j$ be arbitrary and $p$ any integer
greater than $c/t$.
 Then
\[ y^p(M_{n+j} : h^j)\subset (M_{n+pt+j} : h^j)\cap M_c\subset
M_{n+pt}.\]

Therefore
\[(M_{n+j} : h^j)\subset (M_{n+pt} : y^p)\subset (M_n, hM).\]
Thus $(M_{n+j} : h^j)
=M_n+h(M_{n+j} : h^{j+1})$
 for all $n$ and $j$. Iterating this formula $n$ times, we get
that
\[ (M_{n+j} : h^j)
\subset
M_n+ hM_{n-1} + h^2M_{n-2} + \cdots + h^n(M_{n+j} : h^{j+n})=M_n.\]
Hence $h^*$ is a regular element on $\BB$.
\end{proof}

\subsubsection*{Rules of computation} Because of our usage of certain multiplicity functions, particularly
$\hdeg_I(\cdot)$, we recall some basic observations about the
calculation of Hilbert coefficients.

\medskip

%{\bf The reader be warned:}
The terminology of {\em generic hyperplane section} should be
replaced by that of
{\em filter regular element} whenever we
are dealing with Samuel's multiplicity and its Hilbert polynomials.

\medskip

For $s=1$, the $h$ are called {\em generic
hyperplane section}, although there are other usages for this
expression (see \cite[(22.6)]{Nagata}).

\begin{Proposition} \label{genhs} Let $\RR$ and $M$ be as above.
\begin{enumerate}
\item[{\rm (a)}] If
$M\neq 0$ and $h$ is a filter regular element of degree $1$, then
the Hilbert coefficients of $M$ and $M/hM$ satisfy
\begin{eqnarray*}\rme_i(M)& =&\rme_i(M/hM), \quad i< d-1.\\
\rme_{d-1}(M) &=& e_{d-1}(M/hM) + (-1)^{d} \lambda(0:_Mh).
\end{eqnarray*}

\item[{\rm (b)}]
Let
$ 0 \rar A \lar B \lar C \rar 0$
be an exact sequence of finitely generated $\RR$-modules.
If $r=\dim A< s= \dim B$, then $\rme_i(B)=\rme_i(C)$ for $i< s-r$.

\item[{\rm (c)}]
  If $r=0$ and $s\geq 3$, then
\begin{eqnarray*}\rme_i(B)&=&\rme_i(C), \quad i<s \\
\rme_s(B) &= & \rme_s(C) + (-1)^s \lambda(A).
\end{eqnarray*}

\item[{\rm (d)}] \label{e1dim1} If $M$ is a module of dimension $1$, then
$\rme_1(M)=-\lambda(\H_{\m}^0(M))$.

\item[{\rm (e)}] \label{e1dim2} If $M$ is a module of dimension $2$, then
\begin{eqnarray*}
\rme_1(M)&=&\rme_1(M/hM) + \lambda(0:_Mh)\\
&=&- \lambda(\H_{\m}^0(M/hM))+ \lambda(0:_Mh).
\end{eqnarray*}

\item[{\rm (f)}] Furthermore, if $M_0=\H_{\m}^0(M)$ and $M'=M/M_0$, then
\[\rme_1(M)=\rme_1(M').\]

\item[{\rm (g)}]
Let
$ 0 \rar A \lar B \lar C \rar 0$
be an exact sequence of finitely generated $\RR$-modules of the same
dimension. Then
\begin{eqnarray*}
\rme_1(B)\geq\rme_1(A)  +\rme_1(C).\end{eqnarray*}
\end{enumerate}

\end{Proposition}
\begin{proof}
We just observe (6): Since $h$ is regular on $M'$, we have the exact
 sequence $0\rar M_0/hM_0 \rar M/hM \rar M'/hM'\rar 0$,
 from which we get the exact sequence
\[0\rar \H_{\m}^0(M_0/hM_0) \lar
\H_{\m}^0( M/hM)) \lar \H_{\m}^0(M'/hM')\rar 0,\]
which along with the equalities $0:_Mh=0:_{M'}h$ and
$\lambda(0:_{M'}h)=\lambda(M'/hM')$ since $M'$ is a finite length
module, give the desired assertion. 
\end{proof}

Let us add the following observation:

\begin{Proposition}\label{genhs2}
Let $(\RR,\m)$ be a Noetherian local ring, $M$ a finitely generated
$\RR$-module and $h\in \m$ be such that $0:_Mh$ has finite length.
Then
\[ \lambda(0:_Mh)\leq \lambda(\H_{\m}^0(M/hM)).\] 
Furthermore,
\begin{itemize}
\item
If $\dim M>0$ then $h$
is a parameter;
\item
 if  $\dim M>1$ and $M/hM$ is Cohen-Macaulay then $M$ is
Cohen-Macaulay.
\end{itemize}
\end{Proposition}

\begin{proof} Consider the exact sequences induced by multiplication by $h$ on
$M$,
$0\rar 0:_Mh \rar M\rar L \rar 0$ and $0\rar L \rar M \rar
M/hM\rar 0$.
Taking local cohomology we have the exact sequences
\[ 0\rar 0:_Mh \lar \H_{\m}^0(M)\lar  \H_{\m}^0(L)\rar
\H_{\m}^1(0:_Mh)=0,
\]
\[ 0\rar  \H_{\m}^0(L)\lar  \H_{\m}^0(M)\lar  \H_{\m}^0(M/hM).
\]
Counting lengths, we have the desired inequality.

\medskip

To prove that $h$ is a parameter if $\dim M>0$, consider the exact
sequence
\[ 0 \rar \H_{\m}^0(L)=L_0 \lar L \lar L'\rar 0.\]
Reducing mod $h$, we have
the exact sequence
\[ 0\rar 0:_{L_0} h = 0:_Lh \rar 0:_{L'}h \rar L_0/hL_0\rar L/hL \rar
L'/hL'\rar 0.\]
Thus if $0:_{L'}h\neq 0$, it has finite length, which is not possible
since $\H_{\m}^{0}(L')=0$. This proves that
$h$ is a parameter for $L'$ and thus for $L$ as well.
The final assertion follows easily.
 \end{proof}

\begin{corollary}[\cite{MSV10}] \label{e1notpos}
If $M$ is a module of positive dimension  then $\rme_1(M)\leq 0$.
\end{corollary}

\begin{Remark}{\rm A major technical difficulty arise because of the
lack of explicit methods to get hold of $\rme_i(M)$, for $i\geq 2$. For
instance, for $i=2$, the observation above says that it suffices to
know how handle $\dim M=2$: When $\dim M=d>2$, the passage from $M$
to $M/(\xx)M$, where $\xx $ is a superficial sequence of $d-2$, can e
well controlled, that is we can estimate $\H_{\m}^0(M/(\xx)M)$ is
terms of invariants of $M$ (see Theorem~\ref{torsionhdeg}).
}\end{Remark}

\subsubsection*{Multiplicity length of a module}\index{multiplicity length of a module}
We shall now begin to develop mechanisms to allow to discern the role
of associated primes on the multiplicity of a module.

%Let $(\RR,\mathfrak{m})$ be a Noetherian local ring and let $M$ be a
%finitely generated $\RR$-module.

%Used in tandem with the
%following definition it gives some flexible rules of computation.
\medskip

Let $\RR$ be a Noetherian ring and let $M$ be a finitely generated
$\RR$-module. For each $\mathfrak{p}\in \Spec(\RR)$, we denote the
$\RR_{\mathfrak{p}}$-module
\[ \H^0_{\mathfrak{p}}(M_{\mathfrak{p}}) =
 \H^0_{\mathfrak{p}\RR_{\mathfrak{p}}}(M_{\mathfrak{p}})
 \] by $ \Gamma_{\mathfrak{p}}(M_{\mathfrak{p}}) $. (Note that this  is
a $\RR_{\mathfrak{p}}$-module of finite length.)
It is rarely accessible to direct computation, but as we shall see it
occurs in several summation formulas that can be computed.

\begin{Definition}{\rm Let $\RR$ be a Noetherian ring and
 $M$  a finitely $\RR$--module.
For a prime ideal ${\mathfrak p}\subset R$, the integer\index{length multiplicity}
\[ \mult_M({\mathfrak p})= \lambda(\Gamma_{\mathfrak p}(M_{\mathfrak p}))\] is the {\em length
multiplicity} of  ${\mathfrak p}$ with respect to $M$.
}\end{Definition}

This number $\mult_M({\mathfrak p})$, which vanishes if ${\mathfrak p}$ is not
an associated prime of $M$, is a measure of the contribution of
${\mathfrak p}$ to the primary decomposition of the null submodule of
$M$. It is not usually accessible through direct computation except
in exceptional cases.
\medskip

\subsubsection*{Associativity  formulas}\index{associativity formulas for multiplicities} They assemble the
classical multiplicity from local multiplicities. 

\begin{Proposition}[First associativity formula] \label{assocformula}
\index{first associativity formula} Let $(\RR, \mathfrak{m})$ be a Noetherian local
ring and let $M$ be a finitely generated $\RR$-module. Then
\[ \deg(M) = \sum_{\dim \RR/\mathfrak{p}=\dim M}
\mult_M(\mathfrak{p})\deg \RR/\mathfrak{p}.\]
\end{Proposition}

This formula may fail to capture features which
are significant for $M$. For instance, if $M$ is filtered by a
chain of submodules
\begin{eqnarray}\label{filt1}
M = M_1\supset M_2 \supset \cdots \supset M_r\supset 0,\quad
M_i/M_{i+1}\simeq \SS/{\mathfrak p}_i,
\end{eqnarray}
 where
${\mathfrak p}_i $ is a prime ideal,
then \[\deg(M) = \sum\limits_{\dim \SS/{\mathfrak p}_i =\dim M}\deg(\SS/{\mathfrak p}_i ).\]
Alternatively, if ${\mathfrak p} $ is an associated prime
of $M$ of dimension $\dim \SS/{\mathfrak p}= \dim M$, the localization $M_{\mathfrak p} $ is an
Artinian module whose length $\lambda(M_{\mathfrak p})$ is the number
of times ${\mathfrak p} $ occurs in (\ref{filt1}).
This means that all modules  sharing a filtration with the same
$\mathfrak{p}_i$ occurring with the same multiplicities, have the
same value for $\deg(M)$.

\bigskip

The following (\cite[Theorem 24.7]{Nagata})
is the general version of Theorem~\ref{assocformula}.

\begin{Theorem}[General associativity formula]\index{general associativity formula for multiplicities}\label{assoformula2}
Let $(\RR, \m)$ be a Noetherian local ring and $M$ a finitely generated $\RR$--module of dimension
$d$. Let $I = (x_1, \ldots,x_r, \ldots, x_d)$ be a system of parameters for $M$, and set
$J = (x_1, \ldots, x_r)$. Let $\pp$ run over the set of  minimal  prime ideals of $J$. Then
\[ \deg(I;M) = \sum_{\pp} \deg (I+\pp/\pp; \RR/\pp)\cdot \deg(J_{\pp}; M_{\pp}).\]
\end{Theorem}

\subsubsection*{Fitting ideal} The multiplicity of a module may be coded
into its Fitting ideal in several cases of interest.

\begin{Theorem}\label{Fittdeg}
Let  $\RR$ be a Noetherian local ring and
let $E$ be a torsion module with a free presentation
\[ \RR^m \stackrel{\varphi}{\lar} \RR^n \lar E \rar 0.\]
Its {\em Fitting ideal}\index{Fitting ideal of a module} is
\[ \Fitt(E) = I_n(\varphi).\]
Then
\[ \deg(E) = \deg(\RR/\Fitt(E))\]
in the following cases:
\begin{enumerate}
\item[{\rm (a)}] For every prime ideal $\Fitt(E)\subset \mathfrak{p}$,
$E_{\mathfrak{p}}$ has projective dimension $1$ over
$\RR_{\mathfrak{p}}$;

\item[{\rm (b)}] $m=n+g-1$, where $g= \grade \Fitt(E) $.
\end{enumerate}

\end{Theorem}

\begin{proof}
 This follows by applying Proposition~\ref{assocformula} to
$\RR/\Fitt(E)$, according to \cite[Theorem~4.5(2)]{BR}. 
\end{proof}

\begin{Question}{\rm Does the formula of
Theorem~\ref{Fittdeg} apply to other interesting cases of modules, such as Cohen-Macaulay or  of
finite projective dimension?
}\end{Question}

\subsubsection*{Estimation of Hilbert polynomials}\index{Hilbert polynomial: estimation}

Let $\RR$  be a Noetherian local ring or
 a standard graded algebra over the Artinian local ring $\RR_0$.  If $M$ is a finitely generated graded module of dimension $d$, [or an associated graded module]   and Hilbert \[ P_M(n) = \sum_{i=0}^{d-1} (-1)^i \rme_i(M) {{n+d-i-1}\choose{d-i-1}},\]
finding the Hilbert coefficients can be rather demanding. Making use of the Hilbert function $\lambda(M_n)$ would
require knowledge of the $a$-invariant $a(M) $ of $M$. For instance, if $a(M)$ is known, at least some upper bound of
it is at hand, a simple interpolation calculation available, in say, {\it Macaulay2}, would provide the
values.

\begin{Proposition} \label{calcei} \index{Hilbert coefficients: calculation} 
 Suppose $b\geq a(M)$.
Consider the $d\times d$ matrix $\AA$ with binomial entries
\[ \AA = \left[ {{j+d-i-1}\choose {d-i-1}} \right], \quad 0\leq i \leq d-1, b+1\leq j \leq b+d,
\]
and set the column vector
\[ \mathbf{c}^t = \left[ \lambda(M_j) \right], \quad b+1\leq j \leq b+d.\]

For $\mathbf{x}^t = [x_0, x_1, \ldots, x_{d-1}]$, solving for
\[ \AA \cdot \mathbf{x}  = \mathbf{c},\]
would give
\[ \rme_i(M) = (-1)^i x_i, \quad 0\leq i\leq d-1.\]
\end{Proposition}

\begin{Example}{\rm 
Let $\RR = k[x,y,z]$ and $I = (x^2, y^2, z^2, xy-xz, xz - yz)$. 
For $J = (x^2,y^2, z^2)$ we have $I^3 = JI^2$, so we guess $b=1$: The matrix
\[ \AA  = \left[
%\begin{array}{ccc}
%{{2+3-1}\choose{2}}&{{3+3-1}\choose{1}}&{{4+3-1}\choose{0}}\\
%{{2+3-2}\choose{2}}&{{3+3-2}\choose{1}}&{{4+3-2}\choose{0}}\\
%{{2+3-3}\choose{2}}&{{3+3-3}\choose{1}}&{{4+3-3}\choose{0}}
%\end{array} \right]
% = \left[
 \begin{array}{ccc}
6 & 3 & 1 \\
10 & 4 & 1 \\
15 & 5 & 1
\end{array} \right] \]
and
\[ \mathbf{c} = [\lambda(I^2/I^3), \lambda(I^3/I^4), \lambda(I^4/I^5) ] = [36, 64, 100].  \]
Solving 
\[ [\e_0, \e_1, \e_2] = [8, 4, 0],\]
{  by Theorem~\ref{Huckaba}(a) $\e_1 \leq \lambda(I/J) + \lambda(I^2/JI) = 3+3$}. 

\medskip
Now we should compare to taking $b = 2$: We get the same result.
}\end{Example}

\subsection{Reductions and Koszul homology}
The following notion introduced by Northcott and Rees (\cite{NR})
plays a critical role in the whole theory of multiplicities.

\begin{definition}{\rm Suppose that $\RR$ is a commutative ring and $J,
I$  ideals of $\RR$ with $J\subset I$.
Then $J$ is a {\em {reduction}}\index{reduction of an
ideal} of
$I$ if $I^{r+1}= JI^r $
for some integer $r$; the least such integer
 is the {\em \mi{reduction number}} of
 $I$ with respect to $J$. It is denoted by $r_J(I)$.
}\end{definition} \label{relrednum}

If $(\RR, \mathfrak{m})$ is a Noetherian local of dimension $d$, with an
 infinite residue field,
and $I$ is an $\mathfrak{m}$-primary ideal, there a reduction $J$
generated by $d$ elements, $J= (x_1, \ldots, x_d)$.

\medskip

These sequences have the following properties. The second result are
classical theorems of Auslander-Buchsbaum and Serre.

\begin{Theorem}[{\rm \cite[Corollary 4.6.10]{BH}}] Let $\RR$ and $I$ be as above, and let $M$ be a
finitely generated $\RR$-module of dimension $d$. If $J$ is a reduction
of $I$, $\deg_J(M)=\deg_I(M)$.
\end{Theorem}

\begin{Theorem}[{\cite[Theorems 4.7.4, 4.7.6]{BH}}] Let $\bbk(\xx;M)$
be the Koszul complex of $\xx$ with coefficients in $M$ and denote by
$\H_i(\xx;M)$, $i\geq 0$ the homology modules. Then
\begin{enumerate}
\item[{\rm (a)}] $\deg_{\xx}(M)= \sum_{i\geq 0}(-1)^i \lambda(\H_i(\xx;M))$.

\item[{\rm (b)}] $\deg_{\xx}(M) \leq \lambda(M/\xx M) $, with equality if
and only if $\xx$ is a regular sequence on $M$.

\item[{\rm (c)}] More precisely, $\chi_1(\xx;M)= \sum_{i\geq 1}(-1)^{i-1}
\lambda(\H_i(\xx;M))\geq 0$, with equality if and only $\xx$ is a
regular sequence on $M$.
\end{enumerate}
\end{Theorem}

\subsection{Monomial ideals}

(For other
materials see  \cite[Section 7.3]{icbook}.)
Let $f_1 = {\mathbf x}^{v_1}, \ldots , f_r = \mathbf{x}^{v_r}$
be a set of monomials generating the ideal $I$. The integral closure
$\overline{I}$
 is generated by the monomials whose exponents have the form
\[\sum_{i=1}^r r_iv_i + \epsilon \subset \mathbb{N}^n\]
such that $r_i \geq 0$  and $\sum r_i=1$  and $\epsilon$
 is a positive vector with
entries in $[0,1)$.

Consider the integral closure of $I^m$. According to the valuative
criterion, $\overline{I^m}$  is  the integral closure of the ideal
generated by the $m$th powers of the $f_i$.
 This means that the generators of $\overline{I^m}$ are defined by the
exponent vectors of the form

\[\sum r_imv_i + \epsilon,\]
   with $r_i$ and $\epsilon$ as above. We rewrite  this as
\[         m(\sum r_iv_i + \frac{\epsilon}{m}),\]
so the vectors enclosed   must have denominators dividing $m$.

This is the first resemblance to Ehrhart polynomials.
 Suppose that $I$ is of finite co-length, so that $\lambda(\RR/\overline{I})$
 is  the number   of lattice points `outside' the ideal. The
 convex hull $C(V) $ of the $v_i$'s partitions the positive octant
 into $3$ regions: an unbounded connected region, $C(V)$ itself, and
the complement   $\mathcal{P}$ of the union of other two. It is this
 bounded region that is most pertinent to our calculation (see also
 \cite[p. 235]{Stanley}, \cite{Teissier2}).

 To deal with $\overline{I^m}$ we are going to use the set
 of vectors $v_i$, but change the scale by $1/m$.  This means that
 each $I^m$ determines the same $\mathcal{P}$.
 We get that the number of lattice points is the length $\ell_m$ of
 $
\RR/\overline{I^m}$.
 Note that  $\ell_m (1/m)^n$
is the volume of all these hypercubes, and therefore it is a Riemann
sum of the volume of $\mathcal{P}$.
 Thus the limit is just the
volume of the region defined by the convex hull of the $v_i$.
 This number, multiplied by $n!$,   is
also the multiplicity of the ideal. This is the second connection with
the Ehrhart polynomials. We want to emphasize that the volume of
$\mathcal{P}$ may differ from $\lambda(\RR/I)$,
$\lambda(\RR/\overline{I})$, or even from $\lambda(\RR/I_c)$, for the ideal
generated by the lattice points  in the convex hull of $V$.

\bigskip

Let us sum up some of these relationships between multiplicities and
volumes of polyhedra.\index{multiplicity and volume of polyhedra}

\begin{proposition} \label{multasvol} Let $I$ be a
 monomial ideal of the ring $\RR=k[x_1, \ldots,
x_n]$ generated by $\mathbf{x}^{v_1}, \ldots,
\mathbf{x}^{v_m}$. Suppose that $\lambda (\RR/I)< \infty$. If
$\mathcal{P}$ is the region of $\mathbb{N}^n$ defined  as above
by $I$, then

\begin{eqnarray}
{\fbox{$
\begin{array}{c}
\ \\
\  \rme(I) = n! \cdot  \textrm{\rm Vol}(\mathcal{P}) .\\
\ \\
\end{array}  $}}
\end{eqnarray}
\end{proposition}

\begin{example}{\rm
We consider an example of an ideal  of $k[x,y,z]$. Suppose that
\[ I = (x^a, y^b, z^c, x^{\alpha}y^{\beta}z^{\gamma}), \quad
\frac{\alpha}{a}+ \frac{\beta}{b} + \frac{\gamma}{c}< 1.\]
The inequality ensures that the fourth monomial does not lie in the
integral closure of the other three. A direct calculation shows that
\[ \rme(I)= ab\gamma + bc\alpha + ac\beta.\]
It is easy to see that if
 $a\geq 3\alpha$, $b\geq
3\beta$ and $c\geq 3\gamma$, the ideal $J= (x^a-z^c,y^b-z^c,
x^{\alpha}y^{\beta}z^{\gamma})  $
is a minimal reduction and that the reduction number satisfies the
inequality
$\red_J(I)\leq 2$.

}\end{example}

\subsection{Localization and multiplicity loci}
We will now quote two important theorems of M. Nagata (\cite{Nagata})
on the behavior of
the multiplicity under localization. These results have a different
character from those others in this section as they require
structure theory.

\begin{Theorem}[{\rm \cite[Theorem~40.1]{Nagata}}] \label{Nagata40.1}
 Let $\RR$ be a
Noetherian local ring and let $\mathfrak{p}$ be an analytically
unramified prime ideal. If
\[ \dim \RR/\mathfrak{p} + \height \mathfrak{p} = \dim \RR,\] then
$\deg \RR_{\mathfrak{p}}\leq \deg \RR$.
\end{Theorem}

The next result asserts that the multiplicity  function $\deg: \Spec(\RR)\rar \bbz$ is
{\em upper semi-continuous}:

\begin{Theorem}[{\rm \cite[Theorem~40.3]{Nagata}}] \label{Nagata40.3}
 Assume that $\RR$ is a locally
analytically unramified Noetherian ring. For any integer $r>0$,
\[ M_r(\RR) = \{ \mathfrak{p}\in \Spec(\RR) \mid \deg \RR_{\mathfrak{p}}\geq
r\} \]
is a closed set.
\end{Theorem}

\section{Arithmetical Degrees: arith-deg, gdeg, jdeg}

\subsection{Introduction} Let $\RR$ be a Noetherian ring, $\AA$ a
standard graded $\RR$-algebra and $M$ a finitely generated graded
$\AA$-module. (Included is the case when $\RR=\AA$.) The module $M$ has
two sets of associated primes, $\Ass_{\AA}(M)$ and $\Ass_{R}(M)$.
(Both sets are finite.) Along with these sets of primes one has
several functions on $\Spec(\AA)$ and $\Spec(\RR)$:

\begin{enumerate}
\item[{\rm (a)}] If $P\in \Spec(\AA)$, the {\em multiplicity length}
$\mult_M(P)$;

\item[{\rm (b)}] If $\mathfrak{p}\in \Spec(\RR)$, the dimension and multiplicity
of $\H_{\mathfrak{p}}(M_{\mathfrak{p}})$;

\item[{\rm (c)}] If $\RR$ is a local ring, and $\mathfrak{p}\in \Spec(\RR)$,
$\deg(\RR/\mathfrak{p})$;

\item[{\rm (d)}]  $\mathfrak{p}\in \Spec(\RR)$, the local multiplicity
$\deg(\RR_{\mathfrak{p}})$.

%\item
\end{enumerate}

They assemble into several {\em degree} functions. We have already
observed in (\ref{assocformula}) that
\[ \deg(M) = \sum_{\dim \RR/\mathfrak{p}=\dim M}
\mult_M(\mathfrak{p})\deg \RR/\mathfrak{p}.\]
In this section we introduce three other degrees that have this
overall character: They are functions of the form
\[ \ff(M) = \sum_{\mathfrak{p}\in \Spec(\RR)}
a_M(\mathfrak{p}) \cdot b(\mathfrak{p}),
\] where $a_M(\cdot)$ and $b(\cdot)$ mimic the functions listed above.

\subsection{Arithmetic and geometric degrees}

We now define two new degrees. They are extensions of the classical
multiplicity but seek to incorporate data contained in the primary
decomposition of the modules.

\medskip

 Throughout the ring $ \AA$  will
be as before, either a standard graded algebra  or a local ring.

\begin{Definition}{\rm Let $M$ be a finitely generated $\AA$--module.
The {\em arithmetic degree}\index{arithmetic
degree}
 of $M$ is the integer
\begin{eqnarray}
&&
\adeg(M) = \sum_{{\mathfrak p}{ \in \Ass(M)}} \mult_M({\mathfrak p})\cdot \deg
\AA/{\mathfrak p}.\label{defarithdeg}
\end{eqnarray}
}\end{Definition}

When applied to a module $M=\AA/I$,
it is clear that this number gives a numerical signature for the
primary decomposition of $I$.

The definition  applies to general  modules, not just graded
modules, although its main use is for graded modules over
a standard algebra $\AA$ over a local ring $\RR$.

\medskip

 If all the associated primes
of $M$ have the same dimension, then $\adeg(M)$ is just the
multiplicity  $\deg(M)$ of $M$, which is obtained from its Hilbert
polynomial. % (If no confusion arises, we shall also call $\deg(M)$ the
%{\em geometric degree} $\deg(M)$ of $M$; see however \cite{BM} for a
%finer distinction between $\deg(M)$ and $\deg(M)$.)

\medskip

In general, $\adeg(M)$ can be put together as follows.  Collecting  the
associated primes of $M$ by their dimensions:
\[\dim M = d_1 > d_2 > \cdots > d_n,\]
\[ \adeg(M) = a_{d_1}(M) + a_{d_2}(M) + \cdots + a_{d_n}(M),\]
where $a_{d_i}(M)$ is contribution of all primes in $\Ass(M)$ of dimension
$d_i$.

The integer $a_{d_1}(M)$, in the case of a graded module $M$, is its
multiplicity $\deg(M)$. Bayer and Mumford (\cite{BM}) have refined
$\deg(M)$ into the following integer:

\begin{Definition}{\rm The {\em geometric degree}\index{geometric degree} of $M$ is the integer
\begin{eqnarray} \label{gdeg}
&&
\gdeg(M) = \sum_{{\mathfrak p}\ { \mbox{\rm  \ minimal } \in
 \Ass(M)}} \mult_M({\mathfrak p})\cdot \deg \AA/{\mathfrak p}.
\end{eqnarray}
}\end{Definition}

If $\AA = \SS/I$, we want to view $\reg(\AA)$ and $\adeg(\AA)$ as  two basic measures of
complexity of $\AA$.   \cite{BM} has a
far-flung survey of $\reg (\AA)$ and $\adeg (\AA)$ in terms of the degree
data of a presentation $\AA = k[x_1, \ldots, x_n]/I$.

\bigskip

\subsubsection*{Stanley--Reisner rings}\index{Stanley-Reisner rings}
An important class of rings arises from monomial ideals generated by
squarefree elements. They are elegantly coded in the following notion.

\begin{Definition}{\rm
A {\em simplicial complex}\index{simplicial complex}
 $\triangle$ on vertex set $V$
is a family of subsets of $V$ satisfying
\begin{enumerate}
\item[{\rm (a)}] If $x \in V$ then $ \{ x \} \in \triangle$;

\item[{\rm (b)}] If $F \in \triangle $ and $G \subset F$ then $G \in \triangle$.
\end{enumerate}
The elements of $\triangle$ are called {\em faces} or
{\em simplices}. If $|F|=p+1$ then $F$ has {dimension} $p$.
 We define the dimension of
the complex as dim$(
\triangle) = \max_{F \in \triangle}(\mbox{dim} F)$.
}\end{Definition}

\begin{Definition}{\rm
Given any field $k$ and any simplicial complex $\triangle$ on the
finite vertex set $V=\{x_{1},\ldots,$ $x_{n}\}$ define the
{\em face ring} or {\em Stanley-Reisner ring}
 $k[\triangle]$ by
\[
k[\triangle ]= k[x_1, \ldots, x_n]/I_{\triangle},
\]
where
\[
 I_{\triangle}=( x_{i_{1}}x_{i_{2}} \cdots x_{i_{r}}\ |\ i_{1}\!<\!i_{2}\!<\!
\ldots <\!i_{r} \mbox{ and } \{x_{i_{1}},\ldots,x_{i_{r}} \}
\not\in\triangle).
\]
We need only consider the minimal non-faces of $\triangle$ to arrive
at a set of generators for $I_{\triangle}$.}\end{Definition}

Let $f_{p}$ be the number of $p$-simplices. Since
$\emptyset\in\triangle$ and $\dim \emptyset=-1$, $f_{-1}=1$. Also
$f_{0}$ is the number of vertices, thus $f_{0}=n$. We call the
$d$-tuple ${\bf f}(\triangle)=(f_{0},\ldots,f_{d-1})$, consisting of the number of faces
in each dimension, the {\bf f}--{\em vector} of $\triangle$.
The Hilbert series of $k[\triangle]$ is determined entirely by
${\bf f}(\triangle)$,
\[H_{k[\triangle]}({\bf t}) = \frac{h_0+h_1{\bf t}+ \cdots + h_{d}{\bf
t}^d}{(1-{\bf t})^{d}}\]
where the $h_i$ are certain linear combinations of the $f_j$, in
particular \begin{eqnarray} \label{Eulerchar}
h_d = (-1)^{d-1} (\chi(\triangle)-1), \ \chi(\triangle)=
\sum_{i=0}^{d-1}(-1)^if_i. \end{eqnarray}
We quote two elementary facts about these rings (see \cite{BH}).

\begin{Proposition} Let $\triangle$ be a simplicial complex on the
vertex set $V = \{x_1, \dots, $ $ x_n\}$. Then:
\begin{itemize}
\item[{\rm (a)}] $\dim k[\triangle] =1+\dim \triangle=d$.
\item[{\rm (b)}] $I_{\triangle}= \bigcap_F(x_{i_1}, \ldots, x_{i_r},
x_{i_j}\notin F)$.
\end{itemize}
\end{Proposition}

As a consequence we have that
\begin{eqnarray*}
\gdeg(k[\triangle])= \adeg(k[\triangle])&=& \mbox{\rm number of maximal faces of
$\triangle$};\\
\deg(k[\triangle])&=& \mbox{\rm number of faces of $\triangle$ of
maximal dimension}.
\end{eqnarray*}

One elementary use of the arithmetic degree of a module is: a Cayley-Hamilton 
type result:

\begin{Theorem} \label{arithdegCH}\index{Cayley-Hamilton theorem} 
Let $\RR$ be a finitely generated algebra over a field $K$
 and $M$ a finitely generated $\RR$-module. If $\varphi : M\rar M$ is
an endomorphism of $M$ there is a monic polynomial $P(\xx) \in \RR[\xx]$ with
$\deg (P(\xx) )\leq \adeg  ( M)$ such that $P(\varphi)=0$.
\end{Theorem}

\begin{proof} Let $\pp_1,\ldots, \pp_n$ be the associated primes of $M$. Suppose $\pp = \pp_1$ is not
contained in any of the other primes. Set $M_1 = \H_{\pp}^0(M)$ and $M' = M/M_1$. $\varphi$ restricts to
an endomorphism $\varphi_1: M_1 \rar M_1$ which induces an endomorphism $\varphi': M' \rar M'$.

\medskip

It is easy to verify that $\Ass (M) = \Ass (M_1) \cup \Ass (M')$ and 
\[\adeg(M) = \adeg(M_1) + \adeg(M'),\]
indeed if $\qq$ is an associated prime of $M'$ which does not contain $\pp$,  then $M_{\qq}= M'_{\qq}$,
so $\qq \in \Ass(M)$. On the other hand, if $\pp\subset \qq$ and $x'\in M'$ with $\qq x'=0$, for $x\in M$ mapping
to $x'$, $\pp x \in M_1$ and therefore $x\in M_1$.
Furthermore, if $P_1(\xx)$ and $P'(\xx)$ are monic polynomials such that $P_1(\varphi_1)=0$ and 
$P'(\varphi') = 0$.
Consider the endomorphism $\psi = P'(\varphi)P_1(\varphi)$ of $M$. On $M_1$ induces the trivial
endomorphism $P'(\varphi_1)P_1(\varphi_1)$, so it induces an endomorphism on $M'$, $P_1(\varphi')P'(\varphi')$,
which is also trivial. 

\medskip

Thus by induction on the number of associated primes, it suffices to consider the case $\Ass(M)$ has
a unique element $\pp$. If $\pp$ is a maximal ideal, $M$ is an artinian module and
$\adeg(M) = \lambda(M) \deg \RR/\p = \lambda(M)$. Since $\lambda(M) \geq \nu(M)$, it suffices
to find $P(\xx)$ of degree $\leq \nu(M)$. But this is always possible using the standard determinant trick.

\medskip

We can replace $\RR$ by $\RR/\ann(M)$, that is we may assume that $\pp$ is a minimal prime ideal of $\RR$. 
Let $\SS$ be a Noether normalization of $\RR$ and view $M$ and $\varphi$ as defined over $\SS$. $M$ is torsionfree 
over $\SS$ of rank given by $\adeg(M)$. By \cite[Proposition 1.119]{icbook}, there is a monic 
polynomial $P(\xx) \in \SS[\xx] $ such that $P(\varphi) =0$. We lift it to a monic polynomial of
the same degree with the desired property. 
\end{proof}

\begin{Question}{\rm Does the assertion requires that $\RR$ be an affine domain?
}\end{Question}

\subsection{Computation of the arithmetic degree of a module}

We show how a program with the capabilities of  {\em CoCoA}
(\cite{CoCoA}) or
{\em Macaulay 2}
(\cite{Macaulay2}) can be used to compute the arithmetic degree of a graded
module $M$ without availing itself of any primary decomposition. More
details about the explicit computation of $\adeg(M)$ and $\gdeg(M)$
can be found in \cite[p. 234]{compu}.

\medskip

 Let $\SS=k[x_1, \ldots, x_n]$ and suppose $\dim M= d\leq n$. It
suffices to construct graded modules $M_i, \  i=1\ldots n, $ such
that
\[a_i(M) = \deg(M_i).\]

For each integer $i\geq 0$, denote $L_i = \mbox{\rm Ext}_{\SS}^i(M,\SS)$. By local
duality (\cite[Section 3.5]{BH}), a prime ideal ${\mathfrak p}\subset
\SS$ of height $i$
is associated to $M$ if and only if $(L_i)_{\mathfrak p}\neq 0$;
furthermore $\lambda((L_i)_{\mathfrak p})=\mult_M({\mathfrak p})$.

We are set to find a path to $\adeg(M)$:
 Compute for each $L_i$ its degree
$\rme_1(L_i)$  and codimension $c_i$. Then choose $M_i$ according to
\[ M_i = \left\{ \begin{array}{ll}
0 & \mbox{\rm if $c_i> i$}\\
L_i & \mbox{\rm otherwise.}
\end{array} \right. \]

\begin{Proposition}\label{aritgdegformula} For a graded $\SS$--module $M$ and for each integer
$i$ denote by $c_i$ the codimension of $\mbox{\rm Ext}^i_{\SS}(M,\SS)$. Then
\begin{eqnarray}
&&
 \adeg(M) = \sum\limits_{i=0}^n \ \lfloor {\frac{i}{c_i}}
\rfloor \  \deg(\mbox{\rm
Ext}^i_{\SS}(M,\SS)).
  \label{aritdegformula}
\end{eqnarray}
This can also be expressed as

\begin{eqnarray}
&&
 \adeg(M) = \sum\limits_{i=0}^n   \deg(\mbox{\rm
Ext}^i_{\SS}(\mbox{\rm Ext}^i_{\SS}(M,\SS),\SS)).
  \label{aritdegformula2}
\end{eqnarray}
\end{Proposition}

Note that
 in the first  formula (one sets $\lfloor{\frac{0} {0}}\rfloor =
1$)  gives a sum of sums of terms some of which are not always
available.

%\subsubsection*{Hyperplane sections}

%\begin{Definition}{\rm Let $k$ be a field and let
% $A$ be a finitely generated ${\mathbb N}$-graded $k$-algebra,
%\[A = k + A_1 + \cdots + A_i + \cdots .\] $A$ is a
% {\em standard algebra}\index{standard algebra} if it is generated by
% its elements of degree $1$, that is $A_i= (A_1)^i$. A {\em
% hyperplane section}\index{hyperplane section} $h$ of $A$ is a
%form in $A_1$ that is not contained in any minimal prime of $A$
%(sometimes stricter conditions are imposed).
%}\end{Definition}

%If $k$ is an infinite field, most elements of $A_1$ are hyperplane
%sections. A fruitful method to probe $A$ is to compare the properties
%of $A$ with those of $A/(h)$ for some hyperplane section $h$.

%\begin{Theorem} \label{hughs} Let $(\RR,{\mathfrak m})$ be a Noetherian
%local ring of dimension
%$d>0$ and let $M$ be a finitely generated $\RR$--module. Let $x\in
%{\mathfrak m}$ and consider the short exact sequence induced by
%multiplication by $x$,
%\[ 0 \rar L \lar M \stackrel{x}{\lar} M \lar G \rar 0.\]
%If $L$ is a module of finite length, then $\lambda(\H^0_{\mathfrak m}(G))\geq
%\lambda(L)$. Moreover, if $d=1$ then
%\[\lambda(\H_{\mathfrak m}^0(G)) = \lambda(L) +
%\lambda(\overline{M}/x\overline{M}),\] where
%$\overline{M}=M/{\rm torsion}$.
%\end{Theorem}

%Note that the last formula
%follows immediately if $M$ decomposes into a direct sum of
%a torsionfree plus torsion summands. In case $\RR$ is a standard graded
%algebra over a field of characteristic zero, the first assertion is
%contained in the refined statement of \cite[Proposition 3.5]{HuL}.

\medskip

\subsubsection*{Arithmetic/geometric degrees and hyperplane sections}

To  examine the critical behavior of $\adeg(M)$ and $\gdeg(M)$ under  hyperplane
sections we  focus on the following (see \cite{BM}, \cite{MVY},
 \cite{hdeg}):
\begin{Theorem} \label{adegxhp}
Let $\AA$ be a standard graded algebra, let
 $M$ be a finitely generated graded $A$-module  and let $h\in \AA_1$ be a regular element on
$M$.
Then
\begin{eqnarray}
\adeg(M/hM) \geq \adeg(M).
\end{eqnarray}
\end{Theorem}

\begin{proof} We first show that if
 $M$ is a finitely generated $\RR$--module and
  ${\mathfrak p}$ is an associated prime of $M$, then for any $h\in
  \RR$ that
is regular on $M$,
any minimal prime $P$ of $({\mathfrak p},h)$ is an associated prime of
$M/hM$.

 We may assume that
$P$ is the unique maximal ideal of $\RR$. Let $L = \Gamma_{\mathfrak p}(M)$;
note that we cannot have $L  \subset
 hM$ as otherwise since  $h$ is regular on $M$ we
would have $L = hL$. This means that there is a mapping
\[ \RR/{\mathfrak p} \lar M,\]
which on reduction modulo $h$ does not have a trivial image. Since
$\RR/({\mathfrak p},h)$ is annihilated by a power of $P$, its image will
also be of finite length and non-trivial, and therefore $P \in
\Ass(M/hM)$.

\medskip

We take stock of the relationship
 between the associated primes of a
module $M$ and the associated primes of $M/hM$, where $h$ is a regular
element on $M$. According to the above, for each
associated prime ${\mathfrak p}$ of $M$, for which $({\mathfrak p}, h)$ is not
the unity ideal (which will be the case when $M$ is a graded module and
$h$ is a homogeneous form of positive degree) there is  at least
one associated prime of $M/hM$ containing ${\mathfrak p}$: any minimal prime of $({\mathfrak p},h)$
will do.

  There may be  associated primes
of $M/hM$, such as ${\mathfrak n}$, which do not arise in this manner. We indicate this in the
diagram below:
%\begin{verbatim}
$$
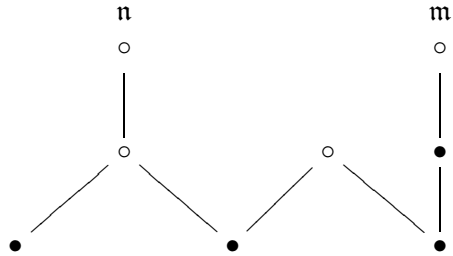
\diagram
& \mbox{\begin{tabular}{c} ${\mathfrak n}$ \\ $\circ$\end{tabular}}
\dline &       & & \mbox{\begin{tabular}{c} ${\mathfrak m}$ \\
$\circ$\end{tabular}}
\dline \\
& \circ\dlline\drline & &\circ \drline\dlline & \bullet\dline \\
\bullet & & \bullet & & \bullet
\enddiagram
$$
%\end{verbatim}
In the expression for $\adeg(M/hM)$ we are going to keep apart primes
such as ${\mathfrak m}$ (which we call associated primes
of the first kind) and primes as ${\mathfrak n}$, which are not minimal
over  $({\mathfrak p},h)$, for any ${\mathfrak p}\in \Ass(M)$.

\medskip

Finally we consider the proof proper.
For  ${\mathfrak m}~\in \Ass(M/hM)$ of the first kind, let \[A({\mathfrak
m})=\{
{\mathfrak
p}_1, \ldots, {\mathfrak p}_r\} \] be the set of associated primes of $M$ such
that ${\mathfrak m}$ is minimal over each ideal $({\mathfrak p}_i,h)$. Note
that these prime ideals have the same dimension, $\dim (A/{\mathfrak m})+1$. Denote
by $I$ the product of the ${\mathfrak p}_i$ and let $L({\mathfrak m})=L = \Gamma_I(M)$. We
have an embedding
\[ L/hL \hookrightarrow M/hM\]
since $L \cap hM = hL$. As a consequence, $\mult_{L/hL}({\mathfrak m})
\leq \mult_{M/hM}({\mathfrak m})$.
On the other hand,
\begin{eqnarray} \label{samemult}
\mult_L({\mathfrak p}_i) & = & \mult_M({\mathfrak p}_i)
\end{eqnarray}
for any prime in the set $A({\mathfrak m})$.

Let  now
\begin{equation}
\label{filtration1} L = L_0 \supset L_1 \supset L_2 \supset \cdots
\end{equation}
be a filtration of $L$ whose factors are of the form $A/{\mathfrak q}$ for
some  homogeneous prime ideal ${\mathfrak q}$. The number of times in which
the  prime ${\mathfrak p}= {\mathfrak p}_i$ occurs in the filtration is $\mult_L({\mathfrak p})$.
Consider the effect of  multiplication by $h$ on the terms of this
filtration: By the snake
lemma, we have exact sequences
\[ 0 \rar {}_hL_{i+1} \lar {}_h L_i \lar {}_hA/{\mathfrak q} \lar
L_{i+1}/hL_{i+1} \lar  L_{i}/hL_i \lar A/({\mathfrak q}, h) \rar 0.\]

Localizing at ${\mathfrak m}$ we get sequences of modules of finite length
over $A_{\mathfrak m}$. Adding these lengths, and taking into account the
collapsing that may occur, we get

\begin{equation} \label{filtration2}
  \lambda(L_{\mathfrak m}/hL_{\mathfrak m}) =
\sum  \lambda((A/({\mathfrak
q},h))_{\mathfrak m})-  \lambda(({}_{h}A/{\mathfrak q})_{\mathfrak m}),
\end{equation}
where  ${}_{h}A/{\mathfrak q}$ is $ A/{\mathfrak q}$ if $h \in {\mathfrak q}$ and
zero otherwise. Note that some of the ${\mathfrak q}$ may occur repeatedly
according to
 the previous observation.
After the localization at ${\mathfrak m}$ only the ${\mathfrak
q}$ corresponding to the ${\mathfrak p}_i$'s survive:
\begin{eqnarray}
 \mult_{L/hL}({\mathfrak m})&=&
 \lambda(L_{\mathfrak m}/hL_{\mathfrak m})
= \sum_{{\mathfrak p}_i\subset {\mathfrak m}} \mult_M({\mathfrak
p}_i)  \lambda((A/({\mathfrak
p}_i,h))_{\mathfrak m}) \label{filtration3}.
\end{eqnarray}

We thus have

\begin{eqnarray}
\adeg(M/hM)  & = & \sum_{{\mathfrak m}} \mult_{M/hM}({\mathfrak
m})\deg(A/{\mathfrak m})\nonumber\\
&+&   \sum_{{\mathfrak n}} \mult_{M/hM}({\mathfrak n})\deg(A/{\mathfrak n})\nonumber \\
& \geq  & \sum_{{\mathfrak m}} \mult_{M/hM}({\mathfrak m})\deg(A/{\mathfrak
m})\nonumber \\
& \geq  & \sum_{{\mathfrak m}} \mult_{L({\mathfrak m})/hL({\mathfrak m})}({\mathfrak
m})\deg(A/{\mathfrak m})\nonumber \\
&=& \sum_{\mathfrak m}(\sum_{{\mathfrak p}_i \in A({\mathfrak m})} \mult_M({\mathfrak
p}_i)\lambda((A/({\mathfrak p}_i,h))_{\mathfrak m}) )\deg (A/{\mathfrak
m}) \label{eq11a} \\
&=&   \sum_{{\mathfrak p}_i}    \mult_M({\mathfrak p}_i) (\sum_{{\mathfrak
p}_i\in A({\mathfrak m})}\lambda((A/({\mathfrak p}_i,h))_{\mathfrak m}) \deg
(A/{\mathfrak m})) \nonumber \\
&=&   \sum_{{\mathfrak p}_i}    \mult_M({\mathfrak p}_i) \deg(A/{\mathfrak
p}_i)\nonumber \\
&=& \adeg(M)\nonumber,
\end{eqnarray}
where in equation (\ref{eq11a}) we first used the equality
(\ref{samemult}), and then derived from each $L({\mathfrak m})$ the
equality provided by the computation of multiplicities in
(\ref{filtration2}) and (\ref{filtration3}), while
taking into account that
\begin{eqnarray*}
\deg(A/{\mathfrak p}_i)&=&\deg(A/({\mathfrak p}_i,h))\\
&=& \sum_{ {\mathfrak p}_i\in A({\mathfrak m})} \lambda((A/({\mathfrak p}_i,
h))_{\mathfrak m})\deg
( A/{\mathfrak m}),
\end{eqnarray*}
the last equality by \cite[Corollary 4.6.8]{BH}. 
\end{proof}

For $\gdeg$ the effect of hyperplane sections cuts the other way:

\begin{Theorem} \label{gdegxhp}
Let $\AA$ be a standard graded algebra, let
 $M$ be a finitely generated graded $\AA$-module  and let $h\in \AA_1$ be a regular element on
$M$.
Then
\begin{eqnarray}
\gdeg(M/hM) \leq \gdeg(M).
\end{eqnarray}
\end{Theorem}

\subsubsection*{Associated graded rings/Initial ideals} There are
comparisons between the arithmetic and geometric degrees of a ring
$\RR$ (or of an $\RR$-module $M$) and the degrees of associated graded
rings $\gr_{\mathcal{A}}(\RR)$ for several filtrations. Some (consider
the case of the  $I$-adic filtration) occur
simply by considering the passage $\RR \rar \RR[It, t^{-1}]  $ followed
by the hyperplane section $t^{-1}$. Other detailed information is
found by a direct analysis (see \cite{Vinai}).

\medskip
In the case of polynomial rings over fields,
a related comparison is that between the degrees of $k[x_1, \ldots,
x_n]/I$, for some homogeneous ideal $I$, and of $k[x_1, \ldots,
x_n]/{\rm in}_{>}(I)$, where $>$ is a term ordering. (See some
bibliography discussed in \cite[Chapter 9]{compu}.)

\subsection{$\jdeg$: a variation of the arithmetic degree}

To construct and develop this  degree,
we will use extensively the notion of {\bf j}--multiplicity
introduced and developed by  Achilles and Manaresi (\cite{AMa93}), Flenner, O'Carroll and Vogel (\cite{FOV}).
It is employed in \cite{Pham} to develop the notion of $\jdeg$.
Here we will  recall the construction of $\jdeg$ and some of its
properties that are required in our examination of normalization.

\bigskip

First though we briefly discuss a family of {\em transforms} from the category of modules
to the category of graded modules. They are simply the composition of two functors: Let $\RR$ be a 
Noetherian ring and $I\subset \m$ be two $\RR$-ideals. For an $\RR$-module $M$ consider the
composition 
\[ M \mapsto \gr_I(M) = \bigoplus_{n\geq 0} I^nM /I^{n+1}M \mapsto \H_{\m}^0(\gr_I(M)).\]
Its properties are being for various purposes as we will have opportunity to examine. A very
important case is that of a local ring, with $\m$
its maximal ideal. Some of its properties are rather intriguing, beginning with the
determination of the Krull dimension of $\H^0_{\m}(\gr_I(M))$. 

\medskip

\subsubsection*{The $j$-multiplicity of a graded
module.}\index{$j$-multiplicity of a module}
Let $(\RR,\mathfrak{m})$ be a Noetherian local ring and
 $\AA$ be a finitely generated graded $\RR$-algebra.
 For a finitely generated graded
$\AA$-module $M$,
 ${\H}^0_\mathfrak{m}(M)$ is a graded
submodule of $M$
which is annihilated by a sufficiently large power of $\mathfrak{m}$.
 Therefore, ${\H}^0_\mathfrak{m}(M)$ can be considered as a graded module over
the Artinian ring  $A/\mathfrak{m}^kA$ for some $k \gg 0$. We will make use of
its Hilbert polynomial to define a new multiplicity function of $M$.

\begin{definition}{\rm
Consider an integer $d$ such that $d \geq \dim M$. We set
\begin{displaymath}
j_d(M)=\left\{\begin{array}{ll}
0&\textrm{if $\dim {\H}^0_\mathfrak{m}(M)< d$},\\
\deg ({\H}^0_\mathfrak{m}(M)) &\textrm{if $\dim {\H}^0_\mathfrak{m}(M)=
d$}.
\end{array}\right.
\end{displaymath}
Moreover, we set $j(M)=j_{\dim M}(M)$, and call it the
$j$-multiplicity of $M$.
}\end{definition}

\begin{remark}{\rm
If $\RR$ is an Artinian local ring, then $j(M)$ is the ordinary
multiplicity $\deg(M)$.
}\end{remark}

Let $(\RR,\mathfrak{m})$ be a Noetherian local ring and $\AA$
be a  standard graded $\RR$-algebra
Let $M$ be a finitely generated graded $\AA$-module.
Recall that the {\em analytic spread} $\ell(M)$ of $M$
is the Krull dimension of $M/\mathfrak{m}M$ as an  $\AA$-module, in particular 
 $\ell ( M ) \leq \ell(\AA)= \dim\AA/\mathfrak{m}\AA$-module.

\begin{lemma} Let $M$ be a finitely generated graded $\AA$-module. Then
$\ell(M)=\dim(M/\mathfrak{m}^nM)$  for any $n\geq 1$.
\end{lemma}

\begin{proof}
It is enough to show that $\dim M/\mathfrak{m}M =\dim
M/\mathfrak{m}^nM$  for any $n \ge 1$.
Suppose that $\mathfrak{m}^n=(a_1,\dots,a_r)$, the surjection
$M^{\oplus{r}}/\mathfrak{m}M^{\oplus{r}} \rar
\mathfrak{m}^nM/\mathfrak{m}^{n+1}M$
 given by sending
$(m_1,\dots,m_r) \longmapsto \sum_{k=1}^r a_kr_k$
shows
$\dim \mathfrak{m}^nM/\mathfrak{m}^{n+1}M \le \dim M^{\oplus{r}}/\mathfrak{m}M^{\oplus{r}}= \dim M/\mathfrak{m}M$.
Observe that the exact sequence
\[0\rar  \mathfrak{m}^{n-1}M/\mathfrak{m}^nM
 \rar M/\mathfrak{m}^nM \rar M/\mathfrak{m}^{n-1}M \rar 0\]
implies $\dim M/\mathfrak{m}^nM = {\max}\{\dim
\mathfrak{m}^{n-1}M/\mathfrak{m}^nM, \dim M/\mathfrak{m}^{n-1}M \}$.
The conclusion then follows immediately by induction.
\end{proof}

\begin{proposition} \label{elljdeg}
The equality $\dim {\H}^0_\mathfrak{m}(M)=\dim M$ holds iff $\ell(M)=\dim(M)$.
\end{proposition}

\begin{proof} Consider the short exact sequence
\[0 \rar {\H}^0_\mathfrak{m}(M)=H \rar M \rar C \rar 0.\]
We note that if $C\neq 0$, it has positive depth.
We know that $\mathfrak{m}^nH=0$ for $n \gg 0$, so we have the
following  exact sequence:

\begin{eqnarray*} \label{elldim}
0 \rar H \cap \mathfrak{m}^nM \rar H \rar M/\mathfrak{m}^nM=M'
 \rar C/\mathfrak{m}^nC=C' \rar 0. \end{eqnarray*}
It follows from the Artin-Rees lemma that
$H \cap \mathfrak{m}^nM=\mathfrak{m}^{n-r}(H \cap \mathfrak{m}^rM)=0$
 for some $r\gg 0$ and $n>r$.
The sequence %(\ref{elldim})
above  becomes
\[0 \rar H \rar M/\mathfrak{m}^nM=M' \rar C/\mathfrak{m}^nC=C' \rar 0\]
Note that $\dim M'=\ell(M)$ by the previous lemma and $\dim H \le
\dim M'\le \dim M$,  so if $\dim H =\dim M$ then $\ell(M)=\dim M$.
\medskip

Conversely, if $\ell(M)=\dim M$ then
$\dim M'=\ell(M)=\dim M \ge \dim C > \dim C'$, the inequality because $\depth C\geq 1$.
It then shows $\dim H=\dim M'$.
\end{proof}

This notion has several properties that contrast to $\deg(\cdot)$, while
others are similar.
First we recall that this {\bf j}-multiplicity behaves well
with respect to short exact sequences (see \cite[Proposition 6.1.2]{FOV}).

\begin{proposition} \label{jmultses}
Let $(\RR,\mathfrak{m})$ be a Noetherian local ring and $\AA$
 be a standard  graded $\RR$-algebra.
Assume that
\begin{eqnarray} \label{jes}
0 \rar M_1 \rar M_2 \rar M_3 \rar 0
\end{eqnarray}
 is an exact sequence of finitely generated graded $\AA$-modules. Then
 for $d = \dim M_2$
\[ j_d(M_2)=j_d(M_1)+j_d(M_3).\]
In particular, if $\dim M_1=\dim M_2=\dim M_3$ then
\[ j(M_2)=j(M_1)+j(M_3).\]
\end{proposition}

\begin{proof}
Taking local cohomology  we obtain the exact sequence
\begin{eqnarray}\label{jes2}
 0 \rar \H^0_{\mathfrak{m}}(M_1) \lar \H^0_{\mathfrak{m}}(M_2)
\lar \H^0_{\mathfrak{m}}(M_3) \lar H \rar 0,
\end{eqnarray} where $H \subset  \H^1_{\mathfrak{m}}(M_1)$.
 If $\dim M_3< d$, $j_d(M_2) = j_d(M_1)$ since $\dim
 \H^0_{\mathfrak{m}}(M_3)< d$. Assume $\dim M_3=d$. We claim that
 $\dim H< d$ and therefore it does not contribute to the multiplicity
 count.
Let $\mathfrak{m}\subset P\subset A$ be an associated prime of $H$, with $\dim
A/P=d$. Localizing (\ref{jes2}) at $P$ we get a short exact sequence of
Artinian $A_P$-modules. Noting that  for $i=1,2,3$,
$\H^0_{\mathfrak{m}}(M_i)_P = (M_i)_P  $, we get $H_P=0$, which is a
contradiction.
\end{proof}

\subsubsection*{Hilbert coefficients of Achilles-Manaresi polynomials}
Let $(\RR, \m)$ be a Noetherian local and $M$ a finitely generated module over the algebra $\AA[A_1]$. 
We can attach to $\H=\H_{\m}^{0}(M)$ a Hilbert function
\[ n \mapsto \lambda(\H_n).\]
The corresponding Hilbert series and Hilbert polynomial will be still written as
$P(M; \ttt)$ and $H(M; \ttt)$.
We use a different notation for the coefficients of these functions:

\[ H(M;\ttt) = \sum_{i = 0}^{r-1} (-1)^i j_i(M) {{\tt+ r-i-1}\choose{r-i-1}}, \quad r=\ell(M). \]
 If $r=1$ this polynomial does not provide for $j_1(M)$, so we use instead
 the function
\[ n \mapsto \sum_{k\leq n} \lambda(\H_k)=
 \sum_{i = 0}^{r} (-1)^i j_i(M) {{n+ r-i}\choose{r-i}}. \]

If
$\H \neq 0$, for effect of comparisons if need be, we regrade $\H$ so that $\H_n=0$ for $n< 0$ and $\H_0 \neq 0$. 
The coefficients $j_i(M)$ are integers but unlike the usual case of an Artinian local $\RR$ it is
very hard to calculate being  related to $M$ via  a mediating construction.
%In addition, some general relationships that are known to exist between the standard coefficients
%$\e_0, \e_1, \e_2$, for instance, are not known. 

\medskip

\bigskip

\subsubsection*{The $\jdeg$ of a module}\index{$\jdeg$ of a graded module}
\label{const}
Next we introduce the {\bf j}-multiplicity of a module $M$ at a given prime $\mathfrak{p}$ of $\RR$. These numbers will 
then be assembled into the $\jdeg$ of $M$.

\begin{definition}{\rm
Let $\RR$ be a Noetherian  ring and $\AA$
 be a standard  graded $\RR$-algebra. For a
 finitely generated graded $\AA$-module $M$, a prime ideal
 $\mathfrak{p}$ of $\RR$, the $j$-multiplicity of $M$ at
 $\mathfrak{p}$ is the integer
\[
j_{\mathfrak{p}}(M)=\left\{ \begin{array}{ll}
\deg ( {\H}^0_{\mathfrak{p}}(M_{\mathfrak{p}}))
 & \textrm{if $\dim  {\H}^0_{\mathfrak{p}}(M_{\mathfrak{p}})
 = \dim M_{\mathfrak{p}}$} \\
0 &\textrm{otherwise}.\\
\end{array}\right.
\]
}\end{definition}

\begin{definition}{\rm
Let $\RR$ be a Noetherian  ring and
 $\AA$ be a standard graded
$\RR$-algebra. For a finitely generated graded
$\AA$-module $M$, the $\jdeg$ of $M$ is the integer
\[\jdeg(M):=\sum_{\mathfrak{p}\in \Spec (\RR)} j_\mathfrak{p}(M).\]
}\end{definition}

\begin{remark}\label{nonzerojdeg}{\rm
Since the summands depend on $\RR$, the  strict notation should be
$\jdeg_{\RR}(M)$. It is easy to see that this is a
 finite sum. It is also a straightforward verification to show that
 $\jdeg(M)=0$ iff $M=0$.
}\end{remark}

If $\RR$ is an integral domain, then $\jdeg(\RR[x])=1$, which shows that
$\jdeg$ does not always code lots of information. This means that if
$\AA\subset \BB$ are very different algebras, one may have
$\jdeg(\AA)=\jdeg(\BB)$. We are going to see that $\jdeg(\BB/\AA)$ is
very discriminating.

We note that
 if $\RR$ is an Artinian local ring, then $\jdeg(M)$
just gives the ordinary multiplicity $\deg(M)$.

\medskip

\subsubsection*{Bounds} Some properties of $\jdeg$ follow directly from
the definition: If $\SS$ is a multiplicative subset of $\RR$,
$\jdeg(\SS^{-1}M)\leq \jdeg(M)$. Others require additional hypotheses
(\cite{Pham}):

\begin{Theorem}\label{compare-deg}Let $(\RR,\mathfrak{m})$ be an analytically unramified
Cohen-Macaulay local ring. Let $\AA$ be a finitely generated standard graded
$\RR$-algebra. Then for every finitely generated graded $\AA$-module $M$
\[\jdeg(M)\leq \gdeg(M).\]
\end{Theorem}

\subsubsection*{Symmetric algebras} Let $\RR$ be an integral
domain and let $M$ be an $\RR$-module of rank $r\geq 1$. Let $\SS$ be the
symmetric algebra of $M$, $\SS=\Sym_{\RR}(M)$. If $M$ is a free (or
projective) $\RR$-module, $\jdeg(\SS)=\jdeg (\RR)=1$. For more general
modules though the associated primes of $\SS$ are difficult to
determine (see \cite[lookup]{alt} for references). Let us consider
one special case.

\begin{proposition} \label{jdegvs} Let $(\RR, \mathfrak{m})$ be a
Noetherian local domain of dimension $d>0$ and let $M$ be a non-free module of
rank $r$ that is free on the punctured spectrum, and let
$\SS=\Sym_{\RR}(M)$.
Then
\[ \jdeg(\SS) = \left\{ \begin{array}{ll}
1 + j_{\mathfrak{m}} (\SS), & \mbox{\rm if $\nu(M) \geq  \dim \RR+\rank(M)$} \\
1, & \mbox{\rm otherwise }
\end{array}
\right.
\]

\end{proposition}

\begin{proof} The only possible primes in $\Ass_{\RR}(\SS)$ are $(0)$, and
$\mathfrak{m}$ if $\SS$ is not an integral domain. According to
Proposition~\ref{elljdeg}, we must check whether  $\ell(M)=\dim \SS$,
 and which by
\cite[lookup]{alt} is equivalent to the assertion on the dimension of
$\SS$.
\end{proof}

To obtain more explicit bounds, one needs other kind of data. Let us
consider one of these that will be used later.

\begin{Theorem} \label{jdegvb1} Let $(\RR, \mathfrak{m})$ be a
Noetherian local domain of dimension $1$, with a finite integral
closure. Let $M$ be a module of rank $r$ that is minimally generated
by $n=r+1$ elements. Let
\[ 0 \rar L \lar \RR^n \lar M \rar 0\]
be a minimal presentation of $M$. Denote by $c(L)$ the $\RR$-ideal
generated by  the entries of the elements of $L$, that is the Fitting
ideal $\Fitt_{r}(M)$.
 If $\SS=\Sym_{\RR}(M)$, then
\[ \jdeg(\SS) \leq 1 + e_0(c(L))+ \nu(L) \cdot  \lambda(\bar{\RR}/\RR).\]
\end{Theorem}

\begin{proof} $L$ is a module of rank $1$, whose generators define the
symmetric algebra of $M$,
\[ 0 \rar (L)=L\BB \lar \BB = \RR[x_1, \ldots, x_n] \lar \SS =
\Sym_{\RR}(M) \rar
0.\]
If $\RR$ is integrally closed, that is, a discrete valuation domain (or
more generally a PID), $(L)$ would be generated by a $1$-form,
$(L)= (c\ff)$, where the entries of $\ff$ generate the unit ideal. In
this case $M \simeq \RR/(c)\oplus \RR^r $, and therefore
$j_{\mathfrak{m}}(\SS) = \lambda(\RR/(c))=\deg(\RR/(c))$. To take advantage
of this elementary fact, let us change make a change of rings,
$\RR\rar \bar{\RR}$.
Consider the exact sequence
\[ 0 \rar \RR \rar \bar{\RR} \rar \bar{\RR}/\RR\rar 0.\]
Tensoring by $\SS$, we obtain the exact complex
\[  \Tor_1^{\RR}(\bar{\RR}/\RR, \SS) \lar \SS \lar \bar{\SS} =
\bar{\RR}\otimes_{\RR}\SS \lar \bar{\RR}/\RR \otimes_{\RR} \SS \rar 0.
\] It gives the inequality
\begin{eqnarray*} j_{\mathfrak{m}}(\SS)
& \leq &
j_{\mathfrak{m}}(\bar{\SS}) +  \deg(\Tor_1^{\RR}(\bar{\RR}/\RR,\SS) \\
& \leq & \lambda(\bar{\SS}/c(L)\bar{\SS}) +
\deg(\Tor_1^{\RR}(\bar{\RR}/\RR,\SS) \\
& \leq & e_0(c(L)) +  \nu(L)\lambda(\bar{\RR}/\RR),
\end{eqnarray*}
where we used   that
$\Tor_1^{\RR}(\RR/\mathfrak{m},\SS) = \RR/\mathfrak{m} \otimes_{\RR}
L\BB  $.  As for
the equality $\lambda(\bar{\SS}/c(L)\bar{\SS}) = e_0(c(L))$, this is a
general fact.
\end{proof}

Observe that $L$ is isomorphic to an ideal of $\RR$, and therefore
$\nu(L)\leq \deg \RR$. As for $\lambda(\bar{\RR}/\RR)$, need to recall a
bound in terms of $\deg(\RR)$.

\bigskip

If $\dim \RR>1$, but $M$ is still free on the punctured spectrum  and
$\dim \RR + \rank(M)=\nu(M)$, it is harder to estimate $\jdeg(\SS)$.
One reason is that the ideal $L\BB$ now has codimension
$\nu(M)-\rank(M)=\dim \RR> 1$, although it is still generated by $1$-forms.
One
approach, is the following. Let $H = \H^0_{\mathfrak{m}}(\SS)\subset \SS$.
Denote by $\mathcal{L}$ the inverse image of $H$ in $\BB = \RR[x_1,
\ldots, x_n]$. Since $\SS/H$ is an integral domain, we actually have
the primary decomposition
\[ L\BB = \mathcal{L} \cap Q,
\] where $Q$ is $\mathfrak{m}\BB$-primary. We can express $\mathcal{L}$
as
\[ \mathcal{L} =\bigcup_{r\geq 1} (L\BB :_{\BB} \mathfrak{m}^r).\]

The ideal $\mathcal{L}$ is preferably obtained by elimination
(\cite{syl1}).
The index $r$ where the operation stabilizes will be called the
of $M$ will be  denoted by $\mbox{\rm ie}(M)$.

\begin{proposition} \label{jmsym}
Let $\RR$, $M$ and $r$ be as above. Then
\[ j_{\mathfrak{m}} (\SS)\leq \lambda(\RR/\mathfrak{m}^r).
\]
\end{proposition}

\begin{proof} We note that $\mathfrak{m}^r H=0$. Consider the exact sequence
\[ 0 \rar H \lar \SS \lar \SS/H \rar 0.\] Tensoring it by
$k=\RR/\mathfrak{m}$, we have the exact sequence
\begin{eqnarray} \label{hssh}
 \Tor_1^{\RR}(k, {\SS}/H) \lar k\otimes_{\RR} H \lar k \otimes_{\RR}
 {\SS} \lar
k\otimes_{\RR} \SS/H \rar 0.
\end{eqnarray}

Let us examine the dimensions of the terms at the end. Since $\dim \SS
= \dim \RR+\rank(M)$ and  $\mathfrak{m}$ is not an associated prime
of $\SS/H$,
$\dim k\otimes_{\RR} \SS/H<  \dim \SS$.

At the other end, $\Tor_1^{\RR}(k,\SS/H)$ is a finitely generated
$\SS/H$-module that is annihilated by $\mathfrak{m}$. Thus its
dimension is also bounded by $\dim \SS-1.$ Counting the multiplicities
in (\ref{hssh}) we only have to account for the modules of dimension
$\dim \SS$. This gives $\deg (k\otimes_{\RR} H)= \deg( k\otimes_{\RR}
\SS)= 1,
$ and therefore using a composition series for $\RR/\mathfrak{m}^r$, we
get
$ \deg(H) = \deg(H\otimes_{\RR} \RR/\mathfrak{m}^r)\leq
\lambda(\RR/\mathfrak{m}^r),$
as desired. 
\end{proof}

\begin{remark}{\rm
A computational approach to this calculation is the following. Let
$I=(a_1, \ldots, a_n)\subset \RR$ be an ideal, let
\[ 0 \rar \mathcal{L} \lar \BB = \RR[\TT_1, \ldots, \TT_n] \lar \RR[It] \rar 0 \]
be a presentation of its Rees algebra. Setting
 \begin{eqnarray*}
L &=& (Q, I), \\
L' & = & L : \mathfrak{m}^{\infty},
\end{eqnarray*}
one has $H = L'/L$. Consider the following example computed by R.
Villarreal,
\[ I = (x,y,z)^2+w(x,y) \subset \RR= k[x,y,z,w].\]
Using {\it Macaulay 2} (\cite{Macaulay2}), one verifies that
\begin{eqnarray*}
 L : \mathfrak{m}^{\infty} & = &  L : \mathfrak{m}^{3}, \\
\deg(H) &=& 12, \\
\lambda(\RR/\mathfrak{m}^3)&= &15.
\end{eqnarray*}
}\end{remark}

\subsection{Associated graded rings and invariance}

Let $\RR$ be a Noetherian ring and let $I$ be an ideal. We are going to
examine some properties of $\jdeg(\cdot)$ as it applies to modules
over the Rees algebra $\RR[It] $ of the ideal $I$. We set $\AA=\RR[It]$.

\begin{proposition} Let $\BB=\bigoplus_{i\geq 0}B_i$ be a finitely generated
graded $\AA$-subalgebra of $\overline{\AA}$. If
 $\BB$ satisfies the condition $S_2$ of Serre then
the filtration of $\{B_i, i\geq 0\}$ is decreasing, i.e. $\gr (\BB)$ is well-defined.
\end{proposition}
\begin{proof} This follows from \cite[Proposition 4.6]{icbook}.
\end{proof}

 An important property of $\jdeg$ is
  the following invariance. It will provide the backbone for many
  calculations.

\begin{Theorem}\label{jdeggr}
Let $\RR$ be a reduced
 universally catenary  ring of dimension
$ d$ and let $I$ be
an ideal of $\RR$.
 Let $\BB$ be a finitely generated  graded $\AA$-subalgebra
 with $\AA=\RR[It]\subsetneq \BB\subset
 \overline{\AA}$, and assume
 that $\BB$
 satisfies the condition $S_2$ of Serre.
Then $\jdeg (\gr(\AA))=\jdeg (\gr(\BB))$.
\end{Theorem}

\begin{proposition} Let $\RR$ be a Noetherian ring, let $A$
be an $\RR$-module and let $\Phi:A \rar A$ be a nilpotent
 endomorphism of $A$. Then $\dim \ker \Phi =\dim A=\dim\coker \Phi$.
\end{proposition}

\begin{proof} This follows from  \cite[Exercise 12.8]{Eisenbudbook}.
Indeed
\[ \rad(\ann(\ker \Phi))=\rad(\ann(A))
=\rad(\ann(\coker \Phi))\]
 so the $\RR$-modules
$\ker \Phi$, $A$, $\coker \Phi$ have the same support. Hence, they all
have the same dimension. 
\end{proof}

Theorem~\ref{jdeggr} is a consequence of the following calculation.

\begin{proposition} \label{prejdeggr}
Let $\RR$ be a reduced
 universally catenary  ring of dimension
$ d$ and let $I$ be
an ideal of $\RR$ containing regular elements.
 Let $B$ be a finitely generated  graded $\AA$-subalgebra defined by a
 decreasing filtration,
 with $\AA=\RR[It]\subsetneq \BB\subset
 \overline{\AA}$.
Then
\[j_\mathfrak{p}(\gr(\AA))=j_\mathfrak{p}(\gr(\BB))\] for
every $\mathfrak{p}
\in \Spec (\RR)$.
\end{proposition}

\begin{proof}
 Let $C=\BB/\AA$, consider the following
 diagram induced by multiplication by $t^{-1}$:
\[
\diagram
& & & 0\dto & \\
& & & C_0\dto & \\
0 \rto & \AA[+1] \rto\dto^{t^{-1}}  & \BB[+1]
\rto\dto^{t^{-1}}  & C[+1] \rto\dto^{t^{-1}} & 0 \\
0 \rto & \AA \rto  & \BB \rto & C \dto\rto & 0\\
& & & C'\dto & \\
& & &  0  &
\enddiagram
\]
which induces the exact sequence
\[0 \rar C_0 \rar \gr(\AA) \rar \gr(\BB) \rar C' \rar 0.\]
Note that $\dim (C_0)_\mathfrak{p}=\dim C'_\mathfrak{p}=\dim C_\mathfrak{p}$
for every $\mathfrak{p} \in \Spec(\RR)$ because the multiplication by
$t^{-1}$ is nilpotent.
 Also, since the {\bf j}-multiplicity is additive with respect to
 exact sequences of modules of the same dimension, using the vertical exact sequence,
one has $j_\mathfrak{p}(C_0)=j_\mathfrak{p}(C') $ for every $
\mathfrak{p}$ since $j_\mathfrak{p}(C)= j_\mathfrak{p}(C[+1])$.
From the other  exact sequence we obtain
 $j_\mathfrak{p}(\gr(\AA))=j_\mathfrak{p}(\gr(\BB))$ for
every $\mathfrak{p}
\in \Spec (\RR)$. 
\end{proof}

Obviously, the most significant part of this proof is the fact that
if $I$ is integral over $J$, then
 $j_\mathfrak{p}(\gr_J(\RR))=j_\mathfrak{p}(\gr_I(\RR))$ for
every $\mathfrak{p}\in \Spec (\RR)$. Later, we shall discuss a theorem of
Flenner and
Manaresi (\cite{FM}) (see Theorem~\ref{FlennerM}) with the proof of a
converse. It is the hard part of a characterization of integrality.

\medskip

\subsubsection*{Approximation complexes} Let $(\RR, \mathfrak{m})$ be a
Noetherian local ring of dimension $d$, and let $I$ be an ideal. If
$\mathcal{A}= { I_n, n\geq 0}$ is an $I$-good filtration, for
\[ \AA = \sum_{n\geq 0} I_nt^n,\]
asserts that
\[ \jdeg(\gr_I(\RR)) = \jdeg(\gr(\AA)).\]
We can further profit by selecting a minimal reduction $J$ of $I$,
and seeking to determine $\jdeg(\gr_J(\RR))$. We are going to discuss
two cases: $\dim \RR/I=1, 2$.

\medskip

Let $(\RR,\mathfrak{m})$ be a Gorenstein local ring of dimension $d$
and let $\mathfrak{p}$
be a prime ideal of dimension $1$, that is $\codim \mathfrak{p}=d-1$.
Assume that $\mathfrak{p}$ is generically a complete intersection. In
order to approach the calculation of $\jdeg(\gr_{\mathfrak{p}}(\RR))$,
we select a minimal reduction $J $ of $\mathfrak{p}$. If
$\ell(\mathfrak{p})=d-1$, $J=\mathfrak{p}$ and
$\jdeg(\gr_{\mathfrak{p}}(\RR))= \deg(\RR/\mathfrak{p})$.

\medskip

We are going to assume that
$\ell(\mathfrak{p})=d$, and that $J$ is one of its minimal reductions
with $\nu(J)=d$. Under these conditions, from the theory of
approximation complexes (\cite[lookup]{alt}), there is  an exact
complex ($\BB = \RR[\TT_1, \ldots, \TT_d]$)
\[ 0 \rar \H_1(J) \otimes_{\RR} \BB[-1] \lar \BB_0=
\H_0(J)\otimes_{\RR}\BB \lar G=\gr_J(\RR)
\rar 0,
\] where $\H_i(J)$ are the Koszul homology modules on a minimal set of
generators of $J$.

\begin{Theorem} \label{jdegdim1} Let $\mathfrak{p}$ be a one-dimensional prime
ideal as above. Set $\SS=\RR/\mathfrak{p}$.
 If $Z_1$ is the module of syzygies of $J$,
\[ \jdeg(\gr_{\mathfrak{p}}(\RR)) \leq 1 + \lambda(\mathfrak{p}/J) +
\lambda(\overline{\SS}/\SS) \cdot \nu(L) +
\lambda(\overline{\SS}/c(Z_1)\overline{\SS}).
\]
\end{Theorem}

\begin{proof} $\H_1(J)$ is the canonical module of $\H_0(J)=\RR/J$. According to
\cite[lookup]{alt}, $\H_1(J)$ is isomorphic to the canonical module of
$\SS=\RR/\mathfrak{p}$. It is thus a torsion free $\SS$-module, while
$\mathfrak{p}/J$ is an $\RR$-module of finite length.
This implies that the submodule $\mathfrak{p}\BB_0$
maps one-one onto $\mathfrak{p}G$.

\medskip

The exact sequence
\[ 0 \rar \mathfrak{p}G =  \mathfrak{p}B_0  \lar G \lar G_0 =
\RR/\mathfrak{p}
\otimes_{\RR} G \rar 0
\] gives the cohomology exact sequence
\[ 0 \rar
\mathfrak{p}B_0  \lar \H^0_{\mathfrak{m}}(G)) \lar
\H^0_{\mathfrak{m}}(G_0)
 \rar 0,
\]
and therefore
\[
j_{\mathfrak{m}}(G) = \lambda(\mathfrak{p}/J) +
j_{\mathfrak{m}}(G_0).
\]

We are now in the setting of Theorem~\ref{jdegvb1}:
\[ \RR/\mathfrak{p}\otimes_{\RR} \gr_J(\RR) =
\Sym_{\RR/\mathfrak{p}}(J/\mathfrak{p}J),
\] where $J/\mathfrak{p}J$ is an
$\RR/\mathfrak{p}$-module of rank $d-1$, free on the punctured spectrum
of $\RR/\mathfrak{p}$, to which we apply the estimation of
Theorem~\ref{jdegvb1}. 
\end{proof}

\begin{example}{\rm
Let $k$ be a field of characteristic $\neq 2$. In the normal
hypersurface ring $\RR=k[x,y,z]/(x^2+yz)$ let $J=\mathfrak{p}= (x,y)$.
Set $\SS= \RR/\mathfrak{p}$; then $\bar{\SS}=\SS=k[z]$ and $c(Z_1) = (z)$.
It follows that
\[ \jdeg(\gr_{\mathfrak{p}}(\RR))=2.\]
}
\end{example}

Let $\mathfrak{p}$ be a prime ideal of dimension $2$ that is a
complete intersection on the punctured spectrum (e.g. $\RR$ is a
regular local ring and $\RR/\mathfrak{p}$ is normal).
We shall assume that $\mathfrak{p}$ is Cohen-Macaulay. Suppose
$\ell(I)=d$ (as otherwise the analysis is simple), and let $J$ be a
minimal reduction (generated by $d$ elements). The approximation
complex associated to $J$,
\[ 0 \rar \H_2(J)\otimes_{\RR} \BB[-2]
\lar \H_1(J)\otimes_{\RR} \BB[-1]
\rar \H_0(J)\otimes_{\RR} \BB
\lar \gr_J(\RR) \rar 0,
\] where $\H_i(J)$ are the Koszul homology modules of $J$, and
$\BB=\RR[\TT_1, \ldots, \TT_d]$, can be quickly examined. We set ourselves in
the context of \cite[Theorem 5.3.4]{alt}.

\begin{proposition} Let $(\RR,\mathfrak{m})$ be a Gorenstein local ring
of dimension $d$, and let $I$ be a prime ideal of codimension $d-2$.
Suppose $I$ is Cohen-Macaulay, and a complete intersection in
codimension $d-1$, $\ell(I)=d$, and let $J$ be a minimal reduction.
The the Koszul homology modules of $J$ satisfy: $\H_2(J) $ is the
canonical module of $\RR/I$ and $\depth \H_1(J)\geq 1$.
\end{proposition}

The proof is analogous to that of \cite[Theorem 5.3.4]{alt}. A
consequence is that both $\H_2(J)$ and $\H_1(J)$ are $\RR/I$-modules,
since $I \H_1(J)$ vanishes at each localization in codimension $d-1$.
As was the case in Theorem~\ref{jdegvb1}, these facts lead to the
exact sequences: (set $B_0 = \BB/J\BB$)

\[ 0 \rar I G =  I B_0  \lar G \lar G_0 =
\RR/I
\otimes_{\RR} G \rar 0
\]
 and the exact sequence of $\RR/I$-modules
\[ 0 \rar \H_2(J)\otimes_{\RR} \BB[-2]
\lar \H_1(J)\otimes_{\RR} \BB[-1]
\rar \RR/I\otimes_{\RR} \BB
\lar G_0 \rar 0. \]

\begin{Theorem}\label{jdegdim2} Let $I$ be as above and $J$ a minimal
reduction of $I$. Then
\[ \jdeg(G) = \lambda(I/J) + \jdeg(G_0)=  \lambda(I/J) + 1+
j_{\mathfrak{m}}(G_0)
.\]
\end{Theorem}

We note that the second sequence says that $G_0$ is the symmetric
algebra of the $\RR/I$-module $J/IJ$, and by the acyclicity lemma $G_0$
is Cohen-Macaulay. The issue is how to make use of the approximation
complex beyond what is already done in Proposition~\ref{jmsym}.

%\subsection{Buchsbaum-Rim multiplicity}

\subsection{Determinants and arithmetic multiplicities}
We give a discussion of the expression of several properties of
modules by associated ideals and their arithmetic multiplicities.

\begin{definition}{\rm
Let $\RR$ be a
Noetherian domain and $E$ a finitely generated $\RR$-module.  If $r$ is
the torsion free rank of $E$, its {\em determinant}\index{determinant
of a module} is
\[{\det}_{\RR}(E)= (\wedge^r E)^{**}= \Hom_{\RR}(\Hom_{\RR}(\wedge^rE, \RR), \RR).\]
For simplicity, we will often just write $\det(E)$ for ${\det}_{\RR}(E)$.
}\end{definition}

\begin{remark}{\rm We note several properties of the notion.

\begin{enumerate}
\item[{\rm (a)}] If $\RR$ is a normal domain, then ${\det}_{\RR}(E)$ is a reflexive
ideal. It is clear that it is isomorphic to a divisorial ideal of
$\RR$. In case $\RR$ is a graded ring and $E$ is homogeneous, ${\det}_{\RR}(E)$
is a homogeneous module with a homogeneous isomorphism with a
divisorial ideal of $\RR$.

\item[{\rm (b)}] If $E^{**}$ is the bidual of $E$, then
${\det}_{\RR}(E)={\det}_{\RR}(E^{**})$. Thus we may assume that $E$ is
reflexive.

\item[{\rm (c)}] The set of (homogeneous) isomorphism classes of the ${\det}_{\RR}(E)$
can be made into a group by usual definition
\[ I\cdot J = (IJ)^{**}.\]

\item[{\rm (d)}] If $E\subset F$ are modules of the same rank $r$, the induced
homomorphism $\wedge^r E \rar \wedge^r F$ leads to a natural
inclusion \[ {\det}_{\RR}(E)\subset {\det}_{\RR}(F).\]

\item[{\rm (e)}] Under the condition above, the ideal
\[ I = \ann({\det}_{\RR}(F)/{\det}_{\RR}(E)) \] is divisorial (and graded when
the modules are homogeneous). Indeed,
\[ I = \Hom_{\RR}({\det}_{\RR}(F), {\det}_{\RR}(E))\subset
\Hom_{\RR}({\det}_{\RR}(F), {\det}_{\RR}(F))\subset \RR.
\]

\end{enumerate}
}\end{remark}

\begin{Theorem} \label{divchains} Let $\RR$ be a Noetherian normal domain, and let
$E\subset F$ be finitely generated $\RR$-modules of the same rank. Set
$I=\ann(\det(F)/\det(E))$. Then $I$ is a divisorial ideal and  denote by
\[ I = \bigcap_1^s \mathfrak{p}_i^{(e_i)}\]
its primary decomposition. Then
\begin{itemize}
\item[{\rm (a)}] Moreover, if $E$ and $F$ are reflexive modules the
maximum length of the chains of distinct reflexive modules between $E$ and
$F$, \[ E_0=E \subset E_1 \subset \cdots \subset E_n=F
\] is given by
\[ \jdeg(\RR/I) = \sum_{i=1}^s\rme_i\]

\item[{\rm (b)}] If $\RR$ is a local ring, or a standard graded
algebra and the modules are homogeneous, then
\[ \jdeg(\RR/I) \leq \deg(\RR/I).\]
\end{itemize}
\end{Theorem}

\begin{proof} The definition of $\det(\cdot)$ localizes properly. In
particular, for every prime ideal $\mathfrak{p}$ of height $1$, we
have that $\det(E)_{\mathfrak{p} }= {\det}(E_{\mathfrak{p}})$. Using
the structure  theorem for finitely generated modules over a PID, it
follows that if $E\subset F$ are reflexive modules of the same rank and
$\det(E)=\det(F)$ then $E=F$. Indeed, in the exact sequence
\[ 0 \rar E \lar F \lar F/E \rar 0,\] by the Depth Lemma (\cite[Lemma
3.52]{icbook}) the associated primes of $F/E$
have height one since $F$ is torsion free and $E$ has the  condition
$S_2$ of Serre.

%y\mathfrak{p}

Note that as we move along the elements of the chain
\[ E_0=E \subset E_1 \subset \cdots \subset E_n=F, \]
the corresponding ideals $I_i= \ann(\det(F)/\det(E_i))$ give rise to
an increasing
chain of {\em distinct} divisorial ideals. This already shows hat the
length of such chains is bounded by $\sum_{i=1}^s\rme_i$.

To build a chain of this length, we induct on $\sum_{i=1}^s\rme_i$. As
already observed, the prime ideals that occur in the decomposition of
$I=\ann(\det(F)/\det(E))$ are the associated primes of $F/E$.

Let $\mathfrak{p} $ be one of the $\mathfrak{p}_i $, say
$\mathfrak{p} =\mathfrak{p}P_1$ and let
$\bar{x}\in F/E$ be an element annihilated by $\mathfrak{p} $. Define
$\rme_1$ as the bidual of the submodule $(E,x)\subset F$. It is an easy
verification that
\[ I_1 = \ann(\det(F)/\det(E_1)) =
\mathfrak{p}_1^{(e_1-1)}\bigcap_{i=2}^s
\mathfrak{p}_i^{(e_i)}.
\]

The second assertion, that $\jdeg(\RR/I)\leq \deg(\RR/I)$ (in the graded
case) is a property of $\jdeg(\cdot)$ that has been discussed already.
\end{proof}

The actual significance of $I$ is the following:

\begin{corollary} \label{dete1} Suppose that $\RR$ is a normal standard graded
algebra (or a normal Noetherian local domain)
and that $E\subset F$ are homogeneous modules of the same rank. If
$E$ and $F$ are reflexive modules, then
\[ \deg(\RR/I)=\rme_1(E)-e_1(F).\]
\end{corollary}

\begin{proof} The difference $\rme_1(E)-e_1(F)$ is the multiplicity of $F/E$. On
the other hand, by Proposition~\ref{assocformula},
\[ \deg(F/E) = \sum_{i=1}^s \lambda((F/E)_{\mathfrak{p}_i})
\deg(\RR/\mathfrak{p}_i).\] After noting that $\rme_i=
\lambda((F/E)_{\mathfrak{p}_i})$, we have the desired formula. 
\end{proof}

\section{Castelnuovo--Mumford Regularity}

\subsection{Introduction}

This is one of most useful of the degree functions and has excellent
treatments in \cite{Eisenbudbook} and \cite{EG}. It has several
equivalent formulations, one of which is the following: Let $\RR=
k[x_1, \ldots, x_d]$ be a ring of polynomials over the field $k$
with the standard grading,  and
let $\AA$ be a finitely generated graded $\RR$-module with a minimal
graded resolution
\[ 0 \rar F_n \lar F_{n-1} \lar \cdots \lar F_1 \lar F_0 \rar 0,\]
\[ F_j = \bigoplus_{j} \RR[-a_{ij}].\] Then $\reg(\cdot)$ is defined by
(see below for a more precise formulation)
\[ \reg(\AA) = \sup\{ a_{i,j}-i\}.\]

Before we point out
 a more abstract formulation of  the Castelnuovo-Mumford regularity,
 or simply {\em regularity},
one of its important properties shows up: The
Hilbert-Poincar\'e-series of $\AA$,
\[ \bbp_{\AA}(\ttt)= \frac{\sum_{ij} (-1)^j \ttt^{a_{ij}}}{(1-\ttt)^d},\]
 encodes the Hilbert function $H_{\AA}(\ttt)$ and Hilbert polynomial
 $\mathcal{P}_{\AA}(\ttt)$, and from the expression for $\bbp_{\AA}(t)$, one has
\[ H_{\AA}(n)=\mathcal{P}_{\AA}(n), \quad n\geq \reg(\AA).\]

\bigskip

Another use of {\em regularity} is for the estimation of various
indices attached to algebraic structures. Consider for instance a
standard graded algebra $\AA$ over the infinite field $k$. If $\dim
\AA=d$, by Noether Normalization, there are subalgebras, $\BB=k[z_1,
\ldots, z_d]$, $z_i\in A_1$, $d=\dim \AA$, such that $\AA$ is finite over
$\BB$. This implies that
\[ A_{n+1} = (z_1, \ldots, z_d)A_n.\]
The smallest such degree $n$ is called the {\em reduction number of
$\AA$ relative to $\BB$}, \index{reduction number of an algebra relative
to a subalgebra} $\red_{\BB}(\AA)$. Its value may vary with the choice of
$\BB$, but the inequality
\[ \red_{\BB}(\AA) \leq \reg(\AA)\]
will always hold.

\bigskip

\subsubsection*{Cohomological formulation}
Let $\RR= \bigoplus_{n\geq 0}R_n=R_0[R_1]$ be a finitely generated
graded algebra over the Noetherian ring $\RR_0$. For any graded $\RR$-module
$F$, define
\[ \alpha(F) = \left\{ \begin{array}{ll}
\sup\{ n\mid F_n \neq 0\} & \textrm{\rm if  } F \neq 0,\\
-\infty & \textrm{\rm if } F = 0.
\end{array}\right.
\]

 Let $\RR= \bigoplus_{n\geq 0} R_n$ be a standard graded
 ring of Krull dimension $d$ with  irrelevant maximal ideal $M = ({\mathfrak m},
\RR_{+})$ (note that $(R_0, {\mathfrak m})$ is  a local ring),
 and let $E$ be a
finitely generated graded $\RR$-module.
Some local cohomology modules $\H^i_J(E)$
of $E$
often give rise to graded modules with $\alpha(\H^i_J(E))< \infty$.
This occurs, for instance, when $J = \RR_{+}$ or $J= M$. This fact
gives rise to several numerical measures of the cohomology of $E$.

\begin{definition}
{\rm
 For any finitely generated graded $\RR$-module $F$, and
for each integer $i\geq 0$, the integer
\[a_i(F) = \alpha(\H^i_M(F))\]
is the $i$th {\em $a$-invariant}\index{$a$-invariant} of $F$.
}\end{definition}

We shall also make use of the following complementary notion.

\begin{definition}
{\rm.
 For any finitely generated graded $\RR$-module $F$, and
for each integer $i\geq 0$, the integer
\[\underline{a}_i(F) = \alpha(\H^i_{R_{+}}(F))\]
is the $i$th {\em
$\underline{a}$-invariant}\index{$\underline{a}$-invariant}  of $F$.
}\end{definition}

By abuse of terminology, if $F$ has Krull dimension $d$, we shall
refer to $a_d(F)$ as simply the $a$-invariant of $F$.
In the case $F=\RR$, if
 $\omega_{\RR}$ is the canonical module of $\RR$,
 by local duality it follows that
\begin{eqnarray}
a(\RR)&=& -\inf \  \{ \ i \ | \ (\omega_{\RR})_i
\neq 0 \ \}. \label{ainvariant1}
\end{eqnarray}

The $\underline{a}_i$-invariants are usually assembled into the
{\em Castelnuovo-Mumford regularity} of $F$\index{Castelnuovo--Mumford
regularity}
\[ \reg(F) = \sup\{ \underline{a}_i(F)+i\mid i\geq 0\}.\]

%This definition can be extended in three ways. If $N$ is an graded
%ideal of $\RR$ and $F$ is a finitely generated graded $\RR$-module the
%local cohomology modules $H_N^i(M)$ are graded and still have the
%property of vanishing in all high degrees. This means that the
%integers below are all well-defined.

\begin{example}{\rm
For several general examples, we direct the reader to \cite[Section
3.6]{BH}. If $\RR$ is Cohen-Macaulay and $I$ is an ideal of positive
height such that $\Rees(I)$ is Cohen-Macaulay, it is easy to see that
$a(\Rees(I))=-1$; ultimately this is a consequence of the fact that
for a field $k$ the canonical ideal of $k[t]$ is $tk[t]$.
}\end{example}

\begin{example}{\rm
The $a$- and $\underline{a}$-invariants are closely related in
several cases of interest.
Suppose that $(\RR, \mathfrak{m})$ is a Cohen-Macaulay local ring of
dimension $d>0$, and $I$  an ideal of dimension one that is
generically a complete intersection. Set $\Gr =\gr_I(\RR)$, choose
$x\in \mathfrak{m}\setminus I$ so that $(I,x)$ is
$\mathfrak{m}$-primary (in particular $I_x$ is a complete
intersection by hypothesis). We let $\mathcal{M}$ be the maximal
homogeneous ideal of $\Gr$ and set $\Gr' = \Gr_x$. Note that $\Gr'$ is a ring
of polynomials so all of its invariants are known.
Let us apply \cite[Proposition 3.15]{icbook} with $J = \Gr_{+}$. We
have exact sequences
\[ \H_{\Gr_{+}}^i(\Gr) \cong \H_{\mathcal{M}}^i(\Gr), \quad i
<d-1\]
\[ 0 \rar \H_{\mathcal{M}}^{d-1}(\Gr)\lar \H_{\Gr_{+}}^{d-1}(\Gr)\lar
\H_{\Gr'_{+}}^{d-1}(\Gr')\lar
\H_{\mathcal{M}}^{d}(\Gr)\lar   \H_{\Gr_{+}}^{d}(\Gr) \rar 0.
\]
Therefore
\begin{eqnarray*}
 a_i(\Gr)& =& \underline{a}_i(\Gr), \quad i< d-1, \\
a_{d-1}(\Gr) & \leq &\underline{a}_{d-1}(\Gr) \leq \max
\{a_{d-1}(\Gr), d-1, \} \\
\max\{ d-1, \underline{a}_d(\Gr)\} \geq a_{d}(\Gr) & \geq &
\underline{a}_d(\Gr).
\end{eqnarray*}
}\end{example}

\subsubsection*{The Castelnuovo-Mumford regularity of Rees algebras}
The following calculation has been carried out by several authors
 (\cite[Proposition 4.1]{JU96},
\cite[Lemma 4.8]{Oo88}, \cite[Theorem
1.1]{Tr98}).

\begin{proposition} \label{regRregG}
Let $\RR$ be a Noetherian local ring and $I$ a proper ideal. Then
$\reg(\RR[It])=\reg(\gr_{I}(\RR))$.
\end{proposition}

The proofs use  the paired short exact sequences (see \cite{Hu0}):
\begin{eqnarray}
0 \rar I\cdot \RR[It]=\mathcal{R}_{+}[+1] \lar \RR[It]=\mathcal{R} \lar {\rm
gr}_I(\RR)=\Gr
\rar 0 \label{hunekeeq10}
\end{eqnarray}
\begin{eqnarray}
0 \rar \mathcal{R}_{+} \lar \mathcal{R} \lar \RR \rar 0,   & &
\label{hunekeeq20}
\end{eqnarray}
with the tautological  isomorphism
\begin{eqnarray*}
&& \mathcal{R}_+=It\cdot \RR[It] \cong \mathcal{R}_{+}[+1]=I\cdot \RR[It]
\end{eqnarray*} playing a pivotal role.

\subsection{Basic properties}
The role of the regularity in the connection between the Hilbert
series and the Hilbert polynomial of modules over polynomial rings
was already observed. The following is more precise and general
(\cite[Theorem 4.4.3]{BH}).

\begin{Theorem} Let $\RR$ be a local Artinian ring, $\AA$ a standard
graded algebra over $\RR$, and $M$ a finitely generated graded
$\AA$-module. Then
\[ H_M(n)-\mathcal{P}_M(n) = \sum_{j\geq
0}(-1)^j\lambda(\H^j_{\AA_{+}}(M)_n),\] in particular
$H_M(n)=\mathcal{P}_M(n)$ for $n> \reg(M)$.
\end{Theorem}

We briefly describe the behavior of $\reg(\cdot)$ with regard to some
exact sequences.

\begin{Proposition}\label{addiofreg} Let $\RR$ be a standard graded algebra, and let
\[ 0 \rar A \lar B \lar C \rar 0\] be an exact sequence of finitely
generated $\RR$-modules $($and homogeneous homomorphisms$)$, then
\[ \reg(B)\leq \reg(A)+ \reg(C),
\]
Similarly,
\[ \reg(A)\leq \reg(B)+ \reg(C), \quad \mbox{\rm and  }\quad
\reg(C)\leq \reg(A)+ \reg(B).\]
\end{Proposition}

\begin{proof} The functions
$\underline{a}_i(\cdot) = \alpha(\H^i_{R_{+}}(\cdot))$ clearly
satisfy
$\underline{a}_i(B)\leq
\underline{a}_i(A)+\underline{a}_i(C)$.
Thus $a_i(B)+i \leq \reg(A)+\reg(C)$ for each $i$.
\medskip

The other assertions have similar proofs.
\end{proof}

A more precise statement (\cite[Proposition 4]{Ooishi82}) asserts:

\begin{Proposition}\label{addiofreg2} Let $\RR$ be a standard graded algebra, and let
\[ 0 \rar A \lar B \lar C \rar 0\] be an exact sequence of finitely
generated $\RR$-modules $($and homogeneous homomorphisms$)$. Then 
\begin{itemize}
\item[{\rm (i)}] 
If $\reg(B) > \reg(C)$, then $\reg(A) = \reg(C)$.
\item[{\rm (ii)}]
If $\reg(B) < \reg(C)$, then $\reg(A) = \reg(C) +1$.
\item[{\rm (iii)}]
If $\reg(B) = \reg(C)$, then $\reg(A) \leq \reg(B)+1$.
\end{itemize}
\end{Proposition}

The function $\reg(\cdot)$ can also be characterized by its initialization  on
modules of finite length  (see \cite[Proposition
20.20]{Eisenbudbook}):

\begin{proposition}\label{charCMreg} If $h$ is a linear form of $\RR$ whose annihilator
$0:_Ah$ has finite length, then
\[ \reg(A) = \max\{\reg(0:_Ah), \reg(A/hA)\}.\]
\end{proposition}

The proof goes through an analysis of the exact sequence
\[ 0 \rar (0:_Ah)[-1] \lar A[-1] \lar A \lar A/hA \rar 0.\]
As a long exercise, the reader can verify that  $\reg(\cdot)$ is the
unique degree function satisfying these rules.

%Let $\Deg$ be an extended degree function

%\subsection{Local cohomology}

\subsection{Bounds on syzygies}

Let $k$ be a field, $\SS= k[x_1,\ldots,x_m]$ a polynomial ring over $k$,
$\RR=\SS/I$ a homogeneous $k$-algebra, and $M$ a finitely generated
graded $\RR$-module.
Then $M$, as an $\SS$-module, admits a finite graded free resolution:
\[
0\rar\oplus_j\SS[-j]^{b_{pj}}\lar\cdots\lar
\oplus_j\SS[-j]^{b_{0j}}\lar M\rar 0.
\]

\begin{Definition}{\rm
The  {\em Castelnuovo-Mumford\index{Castelnuovo--Mumford regularity}}
regularity,
or simply the {\em Castel\-nuovo}
 {\em regularity}, of $M$  is
the integer
\[
\mbox{\rm reg } M=\max\{j-i\colon b_{ij}\neq 0\}.
\]
}\end{Definition}

In other words, $\mbox{\rm reg } M=\max\{\alpha_{+}(\mbox{\rm
Tor}^{\SS}_i(M,k))-i\colon i\in{\mathbb Z}\}$
where for a graded module $N$ with $N_j=0$ for large $j$, we set
$\alpha_{+}(N)=\max\{j\colon N_j\neq 0\}$. Let $q$ be an integer.
The module $M$ is called {\em $q$-regular} if $q\geq \mbox{\rm reg }(M)$,
equivalently, if $\mbox{\rm Tor}^{\SS}_i(M,k)_{j+i}=0$ for all $i$ and all $j>q$.
\bigskip

Eisenbud and Goto \cite{EG} give an interesting interpretation of
regularity. Denoting by $M_{\geq q}$ the truncated graded $\RR$-module
$\bigoplus_{j\geq q}M_j$, one has (see also \cite[Appendix B]{compu}):
\begin{Theorem}
\label{Ego}
The following conditions are equivalent:

\noindent {\rm (a)} $M$ is $q$-regular;

\noindent {\rm (b)} $\H^i_{\mathfrak m}(M)_{j-i}=0$ for all $i$ and all $j>q$;

\noindent {\rm (c)} $M_{\geq q}$ admits a linear $\SS$-resolution, i.e.,
a graded resolution of the form
\[
0\rar \SS[-q-l]^{c_l} \lar
\cdots\lar \SS[-q-1]^{c_1}\lar \SS[-q]^{c_0}\lar
M_{\geq q}\rar 0.
\]
\end{Theorem}

\begin{proof}
(a) $\Leftrightarrow$ (c): By definition, the module
$M_{\geq q}$ has a linear resolution if
and only if for all $i$
\[
\mbox{\rm Tor}^{\SS}_i(M_{\geq q},k)_r=\H_i({\mathbf x};M_{\geq q})_r=0
\]
for $r\neq i+q$. Here $\H({\mathbf x};M)$ denotes the Koszul homology of
$M$ with respect to the sequence ${\mathbf x}=x_1,\ldots,x_m$.

Since $(M_{\geq q})_j=0$ for $j< q$, we always have
$\H_i({\mathbf x};M_{\geq q})_r=0$ for $r<i+q$, while for $r>i+q$
\[
\H_i({\mathbf x};M_{\geq q})_r=\H_i({\mathbf x};M)_r = \mbox{\rm
Tor}^{\SS}_i(M,k)_r.
\]
Thus the desired result follows.

(b) $\Rightarrow$ (c): We may assume $q=0$, and $M=M_{\geq 0}$.
Then it is immediate that $\H^0_{\mathfrak m}(M)$
is concentrated in degree 0. This implies that
$M=\H^0_{\mathfrak m}(M)\oplus M/\H^0_{\mathfrak m}(M)$. The first summand is a
direct sum of copies of $k$. Hence $M$ is $0$-regular if and
only if $M/\H^0_{\mathfrak m}(M)$ is $0$-regular. In other words we may
assume that $\mbox{\rm depth } M>0$. Without any problems we may further assume
that $k$ is infinite. Then there exists an element $y\in S$ of
degree $1$ which is $M$-regular. From the cohomology exact sequence
associated with
\[
0\rar M[-1]\stackrel{y}{\lar} M \lar M/yM\rar 0
\]
we see that $M/yM$ is $0$-regular. By induction on the dimension
on $M$, we may suppose that $M/yM$ has linear $S/yS$-resolution.
But if $F$ is a minimal graded free $\SS$-resolution, the $F/yF$
is a minimal graded $\SS/y\SS$-resolution of $M/yM$. This implies
that $F$ is a linear $\SS$-resolution of $M$.

(c) $\Rightarrow$ (b): Again we may assume $q=0$, and $M=M_{\geq 0}$.
Then $M$ has a linear
resolution
\[
\cdots\lar \SS[-2]^{c_2}\lar \SS[-1]^{c_1}\lar \SS^{c_0}\lar M\rar 0.
\]
Computing $\mbox{\rm Ext}^i_{\SS}(M,\SS)$ with this resolution we see at once that
$\mbox{\rm Ext}^i_{\SS}(M,\SS)_j=0$ for $j<-i$. By duality (see \cite[Section
3.6]{BH}) there exists an isomorphism of graded $\RR$-modules
\[
\H^{i}_{\mathfrak m}(M)\simeq \Hom_k(\mbox{\rm Ext}^{m-i}_{\SS}(M,\SS[-m]),k).
\]
Therefore, $\H^i_{\mathfrak m}(M)_{j-i}=0$ for all $j>0$, as desired.
\end{proof}

\subsection{Regularity of some derived functors}\index{regularity:
Ext and Tor}

Let $\RR=k[x_1, \ldots, x_n]$ be a polynomial ring in $n$
indeterminates over the field $k$. We briefly describe some results
of Chardin, Ha and Hoa (\cite{CHH9}) on regularity bounds for
derived functors of $\Hom(\XX,\YY)$ and $\XX\otimes \YY$.
The first of their results (\cite[Lemma 2.1]{CHH9}) is a beautiful
cohomological calculation:

\begin{Theorem}\label{CHH2.1} Let $\FF_{\bullet}$ be a graded complex
of free $\RR$-modules with \[\FF_i= \bigoplus_{f_i\leq j \leq
b_i}\RR[-j]^{\beta_{ij}}.\] Set $T_i= \sum_j\beta_{ij}$. Then
\[ \reg(\H_i(\FF_{\bullet})) \leq \max\{ b_i, b_{i+1},
[T_{i+1}(b_i-f_{i+1})]^{2^{n-2}} +f_{i+1}+2,
[T_{i}(b_{i-1}-f_{i})]^{2^{n-2}} +f_{i}\}. \]
\end{Theorem}

For two finitely generated graded $M$ and $N$,
it leads directly to estimates for $\reg(\Tor_i^{\RR}(M,N))$. Let us
state their
result on $\reg(\Ext_{\RR}^i(M,N))$ (\cite[Theorem 2.3(2)]{CHH9}).
For a finitely generated graded $\RR$-module $P$, set  the following
notation for the Betti number, initial degree ($\mathrm{indeg}$)
and regularity of $\Tor_i^{\RR}(P,k)$:
\begin{eqnarray*}
\beta_i(P) &= & \dim_k (\Tor_i^{\RR}(P,k))\\
f_i(P)  &= & \mathrm{indeg}( \Tor_i^{\RR}(P,k))\\
\reg_i(P)  &= & \reg( \Tor_i^{\RR}(P,k)).
\end{eqnarray*}

\begin{Theorem}\label{CHH2.5} Let $M$ and $N$ be finitely generated
graded modules over the polynomial ring $\RR$. With the notation
above, set $T_i=\sum_{p-q=i} \beta_p(M)\beta_q(N) $,
$r_M= \reg(M)-\mathrm{indeg}(M)$,
$r_N= \reg(N)-\mathrm{indeg}(N)$, and
$\delta=\mathrm{indeg}(N)-\mathrm{indeg}(M)$. Then
\[ \reg(\Ext_{\RR}^i(M,N)) + i \leq (r_M+r_N+1)^{2^{n-2}} \max \{T_i,
T_{i+1}\}^{2^{n-2}} +1-\delta.
\]
\end{Theorem}

This establishes the existence of polynomials of $\reg(M)$ and
$\reg(N)$ and of the various Betti numbers that bound the regularity
of the Exts and Tors. Kia Dalili reformulated these assertions in a
convenient format in terms of extended degree functions.

\subsubsection*{Questions} There are several questions about the
expected value of the regularity. One of the boldest is (\cite{EG}):

\begin{Conjecture}[Eisenbud-Goto] Let $\RR$ be a standard graded
domain\index{Eisenbud-Goto conjecture}
over an algebraically closed field $k$. Then
\[\reg(\RR)\leq \deg(\RR) - \codim (\RR), \]
where $\codim \RR= \dim_k(\RR_1) - \dim \RR.$
\end{Conjecture}

\section{Cohomological Degrees}

\subsection{Introduction}

Taking Cohen-Macaulay structures,    with the  rich variety
of techniques they support,
(e.g. algebras over a local ring) as
baselines,
 it becomes clear that  the next class of objects
to be examined are those that are Cohen-Macaulay on the punctured
spectrum. Multiplicity theory is still a handy technique here.
It is
to broaden this to very general algebras that
the technique of the extended degrees was introduced.

\medskip

\subsection{Definition of cohomological degrees}\index{cohomological
degree}
Let $(\RR, \mathfrak{m})$ be a Noetherian local ring (or a standard
graded algebra over an Artinian local ring) of infinite residue field. We denote by
$\mathcal{M}(\RR)$ the category of finitely generated $\RR$-module (or
the corresponding category of graded $\RR$-modules).

A general class of these functions was introduced in \cite{DGV} , while
a prototype was defined earlier in \cite{hdeg}. In his thesis
(\cite{Gunston}), T. Gunston carried out a more formal examination
of such functions in order to introduce
his own construction of a new cohomological degree. One of the points
that must be taken care of is that of an appropriate {\em generic
hyperplane} section. Let us recall the setting.

\medskip

Throughout we suppose that the residue field $k$ of $\RR$ is infinite.
Moreover, in some of the discussions we will assume that $\RR$ is a
quotient of a Gorenstein ring. Both conditions are realized by
considering two changes of rings: $\RR\rar \SS=\RR[X]_{\mathfrak{m}[X]}\rar
\widehat{\SS}$, the latter being the completion relative to the maximal ideal.

\begin{definition}{\rm If $(\RR,\mathfrak{m})$ is a local ring, a {\em notion of
genericity}\index{genericity, notion of} on $\mathcal{M}(\RR)$ is a
function
\[ U: \{\textrm{\rm isomorphism classes of $\mathcal{M}(\RR)$} \}
\lar \{ \textrm{\rm non-empty subsets of $\mathfrak{m}/
\mathfrak{m}^2$}\}
\] subject to the following conditions for each $A\in
\mathcal{M}(\RR)$:
 \begin{itemize}
\item[{\rm (a)}]If $f-g\in \mathfrak{m}^2$ then $f\in U(A)$ if and only if
$g\in U(A)$.

\item[{\rm (b)}] The set $\overline{U(A)}\subset
\mathfrak{m}/\mathfrak{m}^2$ contains a non-empty Zariski-open subset.

\item[{\rm (c)}] If $\depth A> 0$ and $f\in U(A)$, then $f$ is regular on
$A$.
\end{itemize}
}\end{definition}

There is a similar definition for graded modules. We shall usually
switch notation, denoting the algebra by $\SS$.

\medskip

Another extension  is that associated to an
$\mathfrak{m}$-primary ideal $I$ (\cite{Linh}): A notion of
genericity on $\mathcal{R}$ with respect to $I$ is a function
\[ U: \{\textrm{\rm isomorphism classes of $\mathcal{M}(\RR)$} \}
\lar \{ \textrm{\rm non-empty subsets of $I/
\mathfrak{m}I$}\}
\] subject to the following conditions for each $A\in
\mathcal{M}(\RR)$:
 \begin{itemize}
\item[{\rm (a)}]If $f-g\in \mathfrak{m}I$ then $f\in U(A)$ if and only if
$g\in U(A)$.

\item[{\rm (b)}] The set $\overline{U(A)}\subset
I/\mathfrak{m}I$ contains a non-empty Zariski-open subset.

\item[{\rm (c)}] If $\depth A> 0$ and $f\in U(A)$, then $f$ is regular on
$A$.
\end{itemize}

%\medskip
%Another extension  is that associated to an
%$\mathfrak{m}$-primary ideal $I$ (\cite{Linh}): A notion of
%genericity on $\mathcal{R}$ with respect to $I$ is a function
%\[ U: \{\textrm{\rm isomorphism classes of $\mathcal{M}(\RR)$} \}
%\lar \{ \textrm{\rm non-empty subsets of $I\setminus
%\mathfrak{m}I$}\}
%\] subject to the following conditions for each $A\in
%\mathcal{M}(\RR)$:
% \begin{itemize}
%\item[{\rm (a)}]If $f-g\in \mathfrak{m}I$ then $f\in U(A)$ if and only if
%$g\in U(A)$.
%\item[{\rm (b)}] The set $\overline{U(A)}\subset
%I/\mathfrak{m}I$ contains a non-empty Zariski-open subset.
%\item[{\rm (c)}] If $\depth A> 0$ and $f\in U(A)$, then $f$ is regular on
%$A$.
%\end{itemize}

\medskip

 Fixing a notion of
genericity $U(\cdot)$ one has the following extension of the
classical multiplicity.

\begin{definition}{\rm A {\em cohomological
degree}, or {\em extended
multiplicity function},
 is a function \label{Degnu}
\[\Deg(\cdot) : {\mathcal M}(\RR) \mapsto {\mathbb N},\]
that satisfies the following conditions.
\begin{itemize}
\item[{\rm (a)}]  If $L = H_{\mathfrak m}^0(M)$ is the submodule of
elements of $M$ that are annihilated by a
power of the maximal ideal and $\overline{M} = M/L$, then
\begin{eqnarray}\label{hs0}
\Deg(M) = \Deg(\overline{M}) + \lambda(L),
\end{eqnarray}
where
$ \lambda(\cdot)$ is the ordinary length function.

\item[{\rm (b)}] (Bertini's rule)
 If $M$ has positive depth, there is $h\in \mathfrak{m} \setminus
 \mathfrak{m}^2$,
  such that
\begin{eqnarray} \label{Bertini}
\Deg(M) \geq \Deg(M/hM).
\end{eqnarray}

\item[{\rm (c)}] (The calibration rule) If $M$ is a Cohen-Macaulay
module, then
\begin{eqnarray} \label{calibration}
\Deg(M) = \deg(M),
\end{eqnarray}
where $\deg(M)$ is the ordinary multiplicity of  $M$.
\end{itemize}
}\end{definition}

These functions will be referred to as {\bf big Degs}. If $\dim \RR =0$,
$\lambda(\cdot)$ is the unique Deg function. For $\dim \RR=1$, the
function $\Deg(M)=
\lambda(L) + \deg(M/L)$ is the unique extended degree. When $d\geq 2$, there are several  big Degs.
An explicit $\Deg$, for all dimensions,  was introduced in \cite{hdeg}.

\begin{definition}{\rm
For an extended degree $\Deg$ the quantity
\[ I(M)=\Deg(M)-\deg(M)\]
is called the {\em Cohen-Macaulay defficiency} of $M$.
\index{Cohen-Macaulay defficiency}
}\end{definition}

We refer to $h$ as a {\em superficial} element relative to $M$ and
$\Deg$. We will abuse the terminology when $h$ is chosen to be
superficial for $M/L$.

\medskip

To define such functions on general local rings one makes use to the
standard changes of rings. For example, if $(\RR, \mathfrak{m})$ is a
Noetherian local ring and $\XX$ is a set of indeterminates,
$\SS=\RR[\XX]_{\mathfrak{m}[\XX]}$ is a Noetherian local ring of
dimension $d$ and the flat change of rings $\RR\rar \SS$ preserves
multiplicity. Its large residue field permits the development of
notions of genericity and the definition of extended degrees.

\medskip

In one special case this notion becomes  familiar.

\begin{proposition} \label{Deg1}
 Let $\RR$ be a local Noetherian ring and  $M$
a finitely generated $\RR$-module of Krull dimension one. For any
extended multiplicity function $\Deg(\cdot)$,
\[ \Deg(M)= \deg(M)+\lambda(\Gamma_{\mathfrak{m}}(M))= \adeg(M).\]
\end{proposition}

If $I$ is an $\mathfrak{m}$-primary ideal, in the proposition above,
replacing $\deg(M)$ by the Samuel's multiplicity $e(I;M)$ would
result in the extended multiplicity $\Deg_I(M)$.
This usually means that $M$ is a finitely generated graded module
over a ring such as $\gr_I(\RR)$.

\medskip

\subsubsection*{Functions with the Bertini rule} \index{Bertini's rule}
Both the notion of genericity and the Bertini property above occur in
limited ways. We will encounter functions
\[ \TT: \mathcal{M }(\RR)\rar \bbz,\]
such that for certain classes of modules $M$ satisfies
\[ \TT(M/hM)\leq \TT(M),\]
for generic elements $h$.

\subsection{General properties of extended degrees}
Let $(\RR, \mathfrak{m})$ be a Noetherian local ring and let $\AA$ be
standard graded algebra over $\RR$. We are going to derive basic
properties of extended degree functions defined over $\RR$ or $\AA$.
Throughout $\Deg$ is an unspecified extended degree function.
We examine some
relationships between the $\Deg(M)$ of a module and some data
expressed in its projective resolution, typically the degrees and
the ranks of the higher modules of syzygies.

\begin{remark} \label{axiom2a} {\rm
The axiom (ii) above admits variations. Let $M$ be a finitely
generated $\RR$-module and set $H=\Gamma_{\mathfrak{m}}(M)$ and $N=M/H$.
Suppose $h$ is superficial
with regard $N$ and $\Deg$. Then we have a short exact sequence
\[ 0 \rar H/hH \lar M/hM \lar N/hN \rar 0,\]
and therefore \[
\Deg(M/hM) =\Deg(H/hH)+ \Deg(N/hN) \leq \Deg(H)+\Deg(N)= \Deg(M).\]
}\end{remark}

\subsubsection*{Cohen-Macaulay defficiency} \index{Cohen-Macaulay
defficiency} The following assertion is a justification for the
terminology.

\begin{Proposition} For any extended degree $\Deg$, the
Cohen-Macaulay defficiency
\[I(M) = \Deg(M)- \deg(M)\]
vanishes if and only if $M$ is Cohen-Macaulay.
\end{Proposition}
\begin{proof}
We begin by observing that the axioms guarantee that $\Deg(M)\geq
\deg(M)$.
Let us use induction on the dimension of $M$. If $\dim M\leq 1$, the
assertion holds since $\Deg(M)=\deg(M) +
\lambda(\Gamma_{\mathfrak{m}}(M))$.
\medskip

If $\dim M\geq 2$ and $L=\lambda(\Gamma_{\mathfrak{m}}(M))$, and $h$
is a superficial element for $\Deg(M')$, $M'=M/L$,
then
\[ \Deg(M)= \lambda(L) + \Deg(M') \geq \lambda(L) + \Deg(M'/hM').
\]
The assumption implies that $L=0$ and $\Deg(M'/hM')=
\deg(M')=\deg(M)$.
By induction $M'$ is Cohen-Macaulay and therefore $M$ is
Cohen-Macaulay as well

\end{proof}

\subsubsection{Betti numbers}
We first want to emphasize the control that extended degrees have
over the Betti numbers of the modules.

\begin{Theorem}\label{lowerDeg} Let $M$ be a module of dimension $d$
and let $\xx=\{x_1, \ldots, x_d\}$ be a superficial sequence relative
to $M$ and the extended degree $\Deg$. Then
\[ \lambda(M/(\xx)M) \leq \Deg(M).\]
\end{Theorem}

\begin{proof} This is a straightforward calculation. If $M_0=
\H^0_{\mathfrak{m}}(M)$ and $M'=M/H$, the exact sequence
\[ 0 \rar H/x_1 H \lar M/x_1M \lar M'/x_1M' \rar 0\]
gives
\[\Deg(M')= \Deg(M)-\lambda(H)\leq \Deg(M) - \lambda(H/x_1H),\]
and allows an easy induction on the dimension of the module. 
\end{proof}

 We assume
that $(\RR,\mathfrak{m}, k=\RR/\mathfrak{m})$
is a Cohen-Macaulay local ring.  For any finitely generated $\RR$-module $A$,
we denote by $\beta^{\RR}_i(A)$ its $i$-th Betti number\index{Betti numbers
of a module!big deg}, and by $\mu^i_{\RR}(A)$
its $i$-th {\em Bass number}\index{Bass numbers of a module!big deg}:
\begin{eqnarray*}
\beta_i^{\RR}(A) & = & \dim_k \Tor^R_i(k,A) \\
\mu^i_{\RR}(A) & = & \dim_k \Ext_R^i(k,A).
\end{eqnarray*}

We recall a classical result on the Betti numbers of the residue field
of $\RR$ under the change of rings $\RR\rar \RR/(\xx)$ where $\xx $ is
a regular element of $\RR$. This is best expressed in terms of the {\em
Poincar\'e series} of $\RR$\index{Poincar\'e series of a local ring}:
\[ P(\RR) = \sum_{i\geq 0} \beta^{\RR}_i(k)\ttt^i.\]

\begin{Theorem}\label{GL342}{\rm (\cite[Corollary 3.4.2]{GL})} Let
$\xx$ be a regular element of $\RR$. Put $\RR'=\RR/(\xx)$. Then
\begin{enumerate}
\item[{\rm (a)}] {\rm ([Tate])} If $\xx\in \m^2$ then
\[ P(\RR') = \frac{P(\RR)}{1-\ttt^2}.\]

\item[{\rm (b)}] If $\xx \in \m \setminus \m^2$ then
\[ P(\RR') =   \frac{P(\RR)}{1+\ttt}.
\]
\end{enumerate}
\end{Theorem}

\begin{Theorem}\label{Degandbetti} Let $\AA$ be a finitely generated
$\RR$-module. For any $\Deg(\cdot)$ function and any integer $i \geq
0$,
\begin{eqnarray*}
 \beta_i^{\RR}(A) &\leq & \beta_i^{\RR}(k)\cdot \Deg(A), \\
 \mu^i_{\RR}(A) & \leq & \mu_{\RR}^i(k)\cdot \Deg(A),
\end{eqnarray*} in particular,
$ \nu(A)\leq \Deg(A)$.
\end{Theorem}

\begin{proof} If $L$ is the submodule of $A$ of finite support, the exact
sequence
\[ 0\rar L \lar A \lar A'\rar 0\]
gives $\beta_i(A) \leq \beta_i(L) + \beta_i(A')$ (for simplicity we
set $\beta_i(A)=\beta_i^{\RR}(A)$). For the summand
$\beta_i(L)$,
by induction on the length of $L$, one has that
 $\beta_i(L)\leq \lambda(L)\beta_i(k)$. For the other summand, one
 chooses a sufficiently generic hyperplane section (good for both $A'$
 and $\RR$ if need be, in particular an element $\xx \in \m \setminus
 \m^2$) and setting $\RR'=\RR/(\xx)$ gives
\[ \beta_i(A') = \beta_i^{\RR'}(A'/\xx A')).\]
Now  by an induction on the dimension of $A'$, one has
\[ \beta_i^{\RR'}(A'/\xx A') \leq \beta_i^{\RR'}(k)\Deg(A'/\xx A').
\]
Finally, since $\Deg(A'/\xx A')\leq \Deg(A')$, we appeal to
Theorem~\ref{GL342}(b) (saying  that $\beta_i^{\RR'}(k)\leq
\beta_i^{\RR}(k)$) to
prove the assertion.

 One can use
a similar argument for $\mu_{\RR}^i(A)$. 
\end{proof}

\medskip

\subsubsection*{Degs and arithmetic degree}\index{extended degree versus arithmetic degree}

\begin{Theorem} Let $(\RR,\m)$ be a Noetherian
 local ring and $\Deg$ an extended degree function on $\RR$-modules. If $M$ is a finitely generated $\RR$-module
 then 
 \[ \adeg(M) \leq \Deg(M).\]
 \end{Theorem}
 
 \begin{proof}
Since both functions are additive on exact sequences $0 \rar \H^0_{\m}(M) \rar M \rar M'\rar 0$
and    have the same value on modules of finite support, we may assume that $\H_{\m}^0(M)=0$. Let $h$
be a suitable hyperplane section, then
\[ \adeg(M) \leq \adeg(M/hM) \leq \Deg(M/hM) \leq \Deg(M),\]
the first inequality by Theorem~\ref{adegxhp}, the second by induction on the dimension, and the
third by the Bertini property of extended degrees.
 \end{proof}

\subsection{Degs and Castelnuovo-Mumford regularity}\index{Degs and
regularity}
Let $\RR$ be a standard graded algebra and $M$ a finitely generated
graded $\RR$-module. Suppose $\alpha(M)$ is the maximum of the degrees
of a minimal generating set of $M$. If $M$ is a module of finite
length, it is clear that
\[ \reg(M) \leq \alpha(M) + \lambda(M)-1.\]
We will prove an extension of this inequality where $\lambda(M)$ is
replaced by $\Deg(M)$ (\cite[Theorem 2.5]{Nagel03}).

\begin{Theorem} \label{Nagel} Let $\RR=k[x_1, \ldots, x_n]$ be a ring of polynomials
over the infinite field $k$, and
 let $M$ be a finitely generated graded $\RR$-module of positive dimension.
Then we have for every sufficiently general linear form $h$ and for
any extended degree $\Deg$,
\[
\Deg(M/hM) -\reg(M/hM)\leq   \Deg(M) - \reg(M).\]
\end{Theorem}

\begin{proof} By Proposition~\ref{charCMreg}, we have
\[ \reg(M) = \max\{ \reg(0:_Mh),\reg (M/hM)\}, \] so that if
$\reg(M)=\reg(M/hM)$ we can make use of Remark~\ref{axiom2a}
asserting  that $\Deg(M/hM)\leq \Deg(M)$ for all sufficiently generic $h's$.

We may thus assume $\reg(M/hM)< \reg(M)$.
%
%We are going to consider several exact sequences and process the
%two functions $\reg(\cdot)$ and $\Deg(\cdot)$ on them.
For any module $A$, denote $A^*= A/\Gamma_{\mathfrak{m}}(A)$.
Consider the exact sequence
 \[ 0 \rar H \lar M \lar N \rar 0,\] where
$H=\Gamma_{\mathfrak{m}}(M)$.  Note that $\dim N>0$.
Since $h$ is regular on $N$, we have the exact sequence
\[ 0 \rar H/hH \lar M/hM \lar N/hN \rar 0.\]
Taking the local cohomology of this sequence we obtain the exact
sequences
\[ 0 \rar H/hH \lar \H_{\mathfrak{m}}^0(M/hM) \lar
\H_{\mathfrak{m}}^0(N/hN)\rar 0.\]
and the isomorphism
\[  \H_{\mathfrak{m}}^i(M/hM) \simeq
 \H_{\mathfrak{m}}^i(N/hN), \quad i\geq 1, \]
$(M/hM)^*\simeq (N/hN)^*$. In particular we have
 that $\reg(N/hN)\leq \reg(M/hM)$.
\medskip

We now process $\Deg(M/hM)$:
\begin{eqnarray*}
\Deg(M/hM) &=& \lambda(\H^0_{\mathfrak{m}}(M/hM) +
 \Deg((M/hM)^*) \\
&=& \lambda(\H^0_{\mathfrak{m}}(M/hM) +
 \Deg((N/hN)^*) \\
&\leq &
\sum_{j\leq \reg(M/hM)} \lambda(\H_{\mathfrak{m}}^0(M)_j) +
\lambda(\H_{\mathfrak{m}}^0(N/hN)_j) +
 \Deg((N/hN)^*) \\
&=& \sum_{j\leq \reg(M/hM)} \lambda(\H_{\mathfrak{m}}^0(M)_j) +
\Deg(N/hN)\\
&\leq & \sum_{j\leq \reg(M/hM)} \lambda(\H_{\mathfrak{m}}^0(M)_j) +
\Deg(N).\\
\end{eqnarray*}
To complete the proof it suffices to add to this expression the
observation
\[ \reg(M) - \reg(M/hM)\leq \sum_{j= \reg(M/hM)+1}^{\reg(M)}
\lambda(\H_{\mathfrak{m}}^0(M)_j).\]
Summing up
 we obtain the desired inequality. 
 \end{proof}

\begin{Theorem}\label{DegvsReg} Let $\AA$ be a finitely generated graded $\RR$-module,
generated by elements of degree at most $\alpha(\AA)$. Then for any extended
degree function $\Deg$,
\[ \reg(\AA) \leq  \Deg(\AA)+ \alpha(\AA).\]
\end{Theorem}

\subsection{The homological degree}

To  establish the existence of cohomological degrees in arbitrary
dimensions, we describe in
some detail one such function introduced in \cite{hdeg}. If $a$ and
$b$ are integers, we set ${a\choose b}=0$ if $a<b$, and ${a\choose b}=1$
if $a=b$ and possibly negative.

\begin{definition}\label{hdegdef}{\rm Let  $M$ be a finitely generated graded module
over the graded
 algebra $\AA$ and  $\SS$  a
 Gorenstein graded algebra mapping onto $\AA$, with maximal graded ideal
${\mathfrak m}$.  Set $\dim \SS=r$,
 $\dim M= d$.
The {\em homological degree}\index{homological degree}\index{hdeg,
 homological degree} of $M$\label{hdegnu}
is the integer
\begin{eqnarray}
 \hdeg(M) &=& \deg(M) +
%+  \sum\limits_{i=r-d+1}^{r} {{d-1}\choose{i-r+d-1}}\cdot
% \hdeg(\mbox{\rm Ext}^i_{\SS}(M,\SS))
  \label{homologicaldegree} \\
&& \sum\limits_{i=r-d+1}^{r} {{d-1}\choose{i-r+d-1}}\cdot
 \hdeg(\mbox{\rm Ext}^i_{\SS}(M,\SS)).\nonumber %\\
%&&+
% \adeg(\mbox{\rm Ext}^r_S(M,S)).\nonumber
\end{eqnarray}
This expression becomes more compact when $\dim M=\dim \SS=d>0$:
\begin{eqnarray}
 \hdeg(M) & = &\deg (M) + \label{homologicaldegree2}\\
&&
 \sum\limits_{i=1}^{d} {{d-1}\choose{i-1}}\cdot
 \hdeg(\mbox{\rm Ext}^i_{\SS}(M,\SS)). \nonumber%\\
%&&+
% \adeg(\mbox{\rm Ext}^d_S(M,S)). \nonumber
\end{eqnarray}
}\end{definition}

The definition of $\hdeg $ can be extended to any Noetherian local ring $\SS$ by
setting $\hdeg(M)= \hdeg(\widehat{\SS}\otimes_{\SS}M)$. On other occasions,
we may also assume that the residue field of $\SS$ is infinite, an
assumption that can be realized by replacing $(\SS, \mathfrak{m})$ by
the local ring $\SS[X]_{\mathfrak{m}\SS[X]}$. In fact, if $X$ is any set
of indeterminates, the localization is still a Noetherian ring, so
the residue field can be assumed to have any needed cardinality, as we shall
assume in the proof.

\begin{remark} \label{Deg2} {\rm
Consider the case when $\dim M = 2$; we assume that $\dim \SS= 2$ also.
The expression for $\hdeg(M)$ is
now
\[ \hdeg(M) = \deg(M) + \hdeg(\Ext_{\SS}^1(M,\SS)) +
\hdeg(\Ext_{\SS}^2(M,\SS)).
\] The last summand, by duality, is the length of the submodule
$\Gamma_{\mathfrak{m}}(M)$. The middle term is a module of dimension at most
one, so can be described according to Proposition~\ref{Deg1} by the
equality
\[ \hdeg(\Ext_{\SS}^1(M,\SS)) =
\deg(\Ext_{\SS}^1(\Ext_{\SS}^1(M,\SS),\SS)) +
\deg(\Ext_{\SS}^2(\Ext_{\SS}^1(M,\SS),\SS)).
\]
}\end{remark}

There are alternative ways to define the  $\hdeg$ of an $\RR$-module
$M$ in a manner that does not require the direct presence of the
Gorenstein ring $\SS$. For simplicity we assume $\dim \SS=\dim M$.
 If we take for $\SS$ the completion of $\RR$
and denote by $E=E_{\SS}(\SS/\m)$ the injective envelope of its
residue field. For each integer $i$ define
\[ M_i = \Hom_{\SS}(\H_{\m}^i(M), E).\]
Note that $M_i$ is a finitely generated $\SS$-module.

\medskip

We leave to the reader to show that:

\begin{proposition} \label{homologicaldegree2a} For $\RR$ and $M$ as above,
\begin{eqnarray*}
 \hdeg(M) = \left\{ \begin{array}{ll}
\lambda(M), & \mbox{\rm if $\dim M=0$} \\
\deg (M) +
 \sum\limits_{i=0}^{d-1} {{d-1}\choose{i}}\cdot
 \hdeg(M_i), & \mbox{\rm if $\dim M>0$}.
\end{array}\right.
%&&+
% \adeg(\mbox{\rm Ext}^d_S(M,S)). \nonumber
\end{eqnarray*}
\end{proposition}

\begin{Theorem} \label{defdheg} The function
$\hdeg(\cdot)$ is a cohomological degree.
\end{Theorem}

The proof requires a special notion of generic hyperplane
sections\index{generic hyperplane section} that fits the concept of
{\em genericity} defined earlier.
Let  $\SS$ be a Gorenstein standard graded ring with infinite residue
field and $M$  a finitely generated  graded module over
$\SS$.  We recall that a {\em superficial} element of order $r$ for $M$
 is an element
$z\in S_r$ such that $0\colon_{M} z$ is a submodule of $M$ of finite
length.

\begin{definition} {\rm A {\em special hyperplane section} of $M$ is an element
$h\in \SS_1$ that is superficial for all the iterated Exts
\[M_{i_1,i_2, \ldots, i_p}= \Ext_{\SS}^{i_1}(\Ext_{\SS}^{i_2}(\cdots
(\Ext_{\SS}^{i_{p-1}}(\Ext_{\SS}^{i_p}(M,\SS),\SS),\cdots, \SS))),
\] and all sequences of integers $i_1\geq i_2\geq\cdots \geq i_p\geq
0$.
}\end{definition}

By local duality it follows that, up to shifts in grading,  there are
only finitely many such modules. Actually, it is enough to consider
those sequences in which $i_1\leq \dim S$ and $p\leq  2 \cdot \dim S$,
which ensures the existence of such $1$-forms as $h$. It is clear
that this property holds for generic hyperplane sections.

\medskip

The following result establishes $\hdeg(\cdot)$ as a {\it bona fide}
cohomological degree.

\begin{Theorem}\label{hdegvsghs}
Let $\SS$ be a standard Gorenstein graded algebra and  $M$  a
finitely generated graded module of depth at least $1$.
 If $h\in \SS_1$ is  a generic
hyperplane section on $M$, then
\[\hdeg(M)\geq \hdeg(M/hM).\]
\end{Theorem}

\begin{proof} This will require
several technical reductions. We assume that $h$ is a {regular},
generic hyperplane section for  $M$ that is regular on $\SS$.
We also assume  that $\dim M=\dim \SS = d$, and derive several exact sequences
from
\begin{eqnarray} \label{reghs}
0 \rar M\stackrel{h}{\lar} M\lar N\rar 0.
\end{eqnarray}
\end{proof}

For simplicity, we write $M_i=\Ext_{\SS}^i(M,\SS)$, and
$N_i=\Ext_{\SS}^{i+1}(N,\SS)$ in the case of $N$. (The latter because $N$
 is a module of
dimension $\dim \SS-1$ and $N_i=\Ext_{\SS/(h)}^i(N,\SS/(h))$.)
\medskip

Using this notation, in view of the binomial coefficients in the
definition of $\hdeg(\cdot)$, it will be enough to prove:

\begin{lemma} \label{hdegmmn} For $M$ and $h$ as above,
\begin{eqnarray*}
\hdeg(N_i)&\leq & \hdeg(M_i)+ \hdeg(M_{i+1}), \ \mbox{\rm for \ }i\geq
1.
\end{eqnarray*}
\end{lemma}

\begin{proof}
The sequence (\ref{reghs})
gives rise to the cohomology long sequence
\begin{eqnarray*}
0 \rar M_0\lar M_0\lar N_0 \lar M_1\lar M_1 \lar N_1\lar M_2 \lar
\cdots \\
\cdots \lar M_{d-2} \lar M_{d-2} \lar N_{d-2} \lar M_{d-1} \lar M_{d-1}\lar
N_{d-1}\rar 0,
\end{eqnarray*}
which are broken up into shorter exact sequences as follows:
\begin{eqnarray}
0 \rar L_i \lar M_i \lar \widetilde{M}_i \rar 0{\ } \label{hs1} \\
0 \rar \widetilde{M}_i \lar M_i \lar G_i \rar 0{\ } \label{hs2}\\
0 \rar G_i \lar N_i \lar L_{i+1} \rar 0. \label{hs3}
\end{eqnarray}

We note that all $L_i$ have finite length, because of the condition on $h$.
For $i=0$, we have the usual relation
$ \deg(M)= \deg(N)$.
When $\widetilde{M}_i$ has finite length, then $M_i, G_i$ and $N_i$ have
finite length, and
\begin{eqnarray*}
\hdeg(N_i)=\lambda(N_i)&=& \lambda(G_i) + \lambda(L_{i+1})\\
&\leq & \hdeg(M_i)+
\hdeg(M_{i+1}).
\end{eqnarray*}
\end{proof}

It is a similar relation that we want to establish for all other cases.

\begin{proposition} \label{hs4} Let $\SS$ be a Gorenstein graded algebra and let
\[0 \rar A \lar B \lar C \rar 0\]
be an exact sequence of graded modules. Then
\begin{enumerate}%[{\rm (ii)}]
\item[{\rm (a)}] If $A$ is a module of finite length, then
\[\hdeg(B) = \hdeg(A)+ \hdeg(C).\]
\item[{\rm (b)}] If $C$ is a module of finite length, then
\[\hdeg(B) \leq \hdeg(A)+ \hdeg(C).\]
\item[{\rm (c)}] Moreover, in the previous case,  if $\dim B=d$, then
\[ \hdeg(A) \leq \hdeg(B)+ (d-1)\hdeg(C).\]
\item[{\rm (d)}] If $C$ is a module of finite length and $\depth
B\geq 2$, then
\[\hdeg(A) = \hdeg(B)+ \hdeg(C).\]
\end{enumerate}
\end{proposition}

\begin{proof} They are all  clear if $B$ is a module of finite length so we
assume that $\dim B=d \geq 1$.

\medskip

(a) This is immediate since $\deg (B)=\deg(C)$ and the cohomology sequence
gives
\begin{eqnarray*}
\Ext_{\SS}^i(B,\SS) &=& \Ext_{\SS}^i(C,\SS), \ 1\leq i\leq  d-1, \ \mbox{\rm and}\\
\lambda(\Ext_{\SS}^d(B,\SS)) & = & \lambda(\Ext_{\SS}^d(A,\SS))+ \lambda(\Ext_{\SS}^d(C,\SS)).
\end{eqnarray*}

(b) Similarly we have
\begin{eqnarray*}
\Ext_{\SS}^i(B,\SS) &=& \Ext_{\SS}^i(A,\SS), \ 1\leq i< d-1,
\end{eqnarray*}
and the exact sequence
\begin{small}
\begin{eqnarray} \label{hs41}
0 \rar \Ext_{\SS}^{d-1}(B,\SS)\rar \Ext_{\SS}^{d-1}(A,\SS) \rar
\Ext_{\SS}^{d}(C,\SS) \\
 \rar \Ext_{\SS}^{d}(B,\SS) \rar \Ext_{\SS}^{d}(A,\SS) \rar 0.\nonumber
\end{eqnarray}
\end{small}
If $\Ext_{\SS}^{d-1}(A,\SS)$ has finite length, then
\begin{eqnarray*}
\hdeg(\Ext_{\SS}^{d-1}(B,\SS)) &\leq & \hdeg(\Ext_{\SS}^{d-1}(A,\SS)) \\
\hdeg(\Ext_{\SS}^{d}(B,\SS)) & \leq & \hdeg(\Ext_{\SS}^{d}(A,\SS)) +
\hdeg(\Ext_{\SS}^{d}(C,\SS)). \end{eqnarray*}
Otherwise, $\dim \Ext_{\SS}^{d-1}(A,\SS)= 1 $, and
\begin{eqnarray*}
\hdeg (\Ext_{\SS}^{d-1}(A,\SS)) &=& \deg (\Ext_{\SS}^{d-1}(A,\SS)) +
\lambda(\Gamma_{\mathfrak m}(\Ext_{\SS}^{d-1}(A,\SS))).
\end{eqnarray*}
Since we also have
\begin{eqnarray*}
\deg(\Ext_{\SS}^{d-1}(B,\SS))& = & \deg(\Ext_{\SS}^{d-1}(A,\SS)), \\
\lambda(\Gamma_{\mathfrak m}(\Ext_{\SS}^{d-1}(B,\SS)))& \leq  & \lambda(\Gamma_{\mathfrak
m}(\Ext_{\SS}^{d-1}(A,\SS))), \end{eqnarray*}
we again obtain the stated bound.

\medskip

(c) In the sequence (\ref{hs41}), if $\dim \Ext_{\SS}^{d-1}(B,\SS)=0$,
then
\begin{eqnarray}\label{hs42}
 \lambda(\Ext_{\SS}^{d-1}(A ,\SS))\leq \lambda(\Ext_{\SS}^{d-1}(B ,\SS)) +
\lambda(C)
\end{eqnarray}
and also
\[ \lambda(\Ext_{\SS}^{d}(A ,\SS))\leq \lambda(\Ext_{\SS}^{d}(B ,\SS)).
\] When taken into the formula for $\hdeg(A)$, the binomial
coefficient
${d-1}\choose{d-2}$ gives the desired factor for $\lambda(C)$.

On the other hand, if $\dim \Ext_{\SS}^{d-1}(B,\SS)=1$, we also have
\[ \hdeg(\Ext_{\SS}^{d-1}(A ,\SS))\leq \hdeg(\Ext_{\SS}^{d-1}(B ,\SS)) +
\lambda(C),\]
the dimension one case of (\ref{hs42}).

\medskip

(d) This follows by applying the definition of $\hdeg$ to the exact
sequence.
\end{proof}

Suppose that $\dim \widetilde{M_i}\geq 1$.
% From (\ref{hs1}) and
%Proposition~\ref{hs4}(a),
%\[\hdeg(M_i)= \hdeg(\widetilde{M_i}) + \lambda(L_i),\]
%while from (\ref{hs3}) and
From Proposition~\ref{hs4}(b) we have
\begin{eqnarray}\label{hs5}
 \hdeg(N_i) \leq \hdeg(G_i) + \lambda(L_{i+1}).
\end{eqnarray}
We must now relate $\hdeg(G_i)$ to $\deg({M}_i)$.
Apply the functor $\Gamma_{\mathfrak m}(\cdot)$ to the sequence
(\ref{hs2}) and consider the commutative diagram
\[
\begin{array}{ccccccccc}
0 &\rar& \widetilde{M}_i &\lar& M_i &\lar& G_i& \rar& 0\\
&& \uparrow && \uparrow && \uparrow && \\
0 &\rar& \Gamma_{\mathfrak m}(\widetilde{M}_i) &\lar& \Gamma_{\mathfrak m}(M_i)
 &\lar& \Gamma_{\mathfrak m}(G_i)&
\end{array},
\] in which we denote by $H_i$ the image of the natural map
\[\Gamma_{\mathfrak m}(M_i)\lar \Gamma_{\mathfrak m}(G_i).\]
Through the snake lemma, we obtain
the exact sequence
\begin{eqnarray} \label{hs6}
0 \rar \widetilde{M_i}/\Gamma_{\mathfrak m}(\widetilde{M}_i) \stackrel{\alpha}{\lar}
{M_i}/\Gamma_{\mathfrak m}(M_i) \lar
G_i/H_i \rar 0.
\end{eqnarray}

Furthermore,  from (\ref{hs1}) there is  a  natural isomorphism
\begin{eqnarray*}
\beta: {M_i}/\Gamma_{\mathfrak m}(M_i)\cong
  \widetilde{M_i}/\Gamma_{\mathfrak m}(\widetilde{M}_i),
\end{eqnarray*}
while from (\ref{hs2}) there is a natural injection
\begin{eqnarray*}
 \widetilde{M_i}/\Gamma_{\mathfrak m}(\widetilde{M}_i)
\hookrightarrow{M_i}/\Gamma_{\mathfrak m}(M_i),
\end{eqnarray*} whose composite with $\beta$
 is induced by multiplication by $h$
on $M_i/\Gamma_{\mathfrak m}(M_i)$. We may thus replace
$  \widetilde{M_i}/\Gamma_{\mathfrak m}(\widetilde{M}_i)$ by
${M_i}/\Gamma_{\mathfrak m}(M_i)$ in (\ref{hs6}) and take $\alpha$ as
multiplication by $h$:
\[0 \rar M_i/\Gamma_{\mathfrak m}(M_i) \stackrel{h}{\lar}
M_i/\Gamma_{\mathfrak m}(M_i) \lar G_i/H_i \rar 0.\]

\medskip
Observe that
since \[\Ext_{\SS}^j(M_i/\Gamma_{\mathfrak m}(M_i),\SS) = \Ext_{\SS}^j(M_i,\SS), \
j<\dim S,\]
 $h$ is still a regular, generic
 hyperplane section for $M_i/\Gamma_{\mathfrak
m}(M_i)$.
By induction on the dimension of the module,
we have
\[\hdeg(M_i/\Gamma_{\mathfrak
m}(M_i)) \geq \hdeg(G_i/H_i).\]

Now from Proposition~\ref{hs4}(a), we have
\[\hdeg(G_i)= \hdeg(G_i/H_i) + \lambda(H_i).\] Since these summands are
bounded, by
$\hdeg(M_i/\Gamma_{\mathfrak
m}(M_i)) $ and $\lambda(\Gamma_{\mathfrak m}(M_i))$ respectively,
 (in fact, $\lambda(H_i)=
\lambda(L_i)$),  we have
\[\hdeg(G_i)\leq \hdeg(M_i/\Gamma_{\mathfrak
m}(M_i)) +\lambda(\Gamma_{\mathfrak m}(M_i)) = \hdeg(M_i),\]
the last equality by Proposition~\ref{hs4}(a) again.
Finally, taking this estimate into (\ref{hs5}) we get
\begin{eqnarray} \label{hs7}
\hdeg(N_i)&\leq &\hdeg(G_i) + \lambda(L_{i+1})\\
& \leq & \hdeg(M_i) +
\hdeg(M_{i+1}), \nonumber \end{eqnarray}
to establish the claim.

\begin{remark}{\rm That equality does not always hold is shown by the
following example. Suppose that $\RR=k[x,y]$ and $M= (x,y)^2$. Then $\hdeg(M) =
4$,  but $\hdeg(M/hM)=3$ for any hyperplane section $h$.
 To get an
example of a ring one takes the idealization of $M$.
}\end{remark}

\begin{remark}{\rm
It should be emphasized that the {\em weight} binomial coefficients in the
definition of $\hdeg$ were chosen to enable the Bertini property at
the expense of its behavior on other short	 exact sequences. Suppose
$\RR=k[x,y,z]$, and $M = R \oplus R/(x,y)$. Then
\[ \hdeg(M) = \deg(\RR) + {3-1\choose 2-1} \deg(\RR/(x,y)) = 1+ 2=3.\]
}\end{remark}

\begin{example}{\rm This example shows how the function $\hdeg$
captures important aspects of the structure of the module.
We
  recall the notion of a
{\em sequentially} Cohen-Macaulay module\index{sequentially
Cohen-Macaulay module}. This is a module $M$ having a filtration
\[ 0 = M_0 \subset M_1 \subset \cdots \subset M_r = M,\]
with the property that each factor  $M_i/M_{i-1}$ is Cohen-Macaulay
and
\[
\dim M_i/M_{i-1} < \dim M_{i+1}/M_i, \ i=1, \ldots, r-1 .\]
If $L$ is the leftmost nonzero submodule in such a chain, it follows
easily that $N = M/L$ is also sequentially Cohen-Macaulay ($L$ is
Cohen-Macaulay by hypothesis) and we have
\begin{eqnarray*}
\Ext_{\SS}^p(M,\SS) & = & \Ext_{\SS}^p(L,\SS), \ p = \dim L \\
\Ext_{\SS}^i(M,\SS) & = & \Ext_{\SS}^i(N,\SS), \ i < p
\end{eqnarray*}
with the other $\Ext$'s vanishing. This gives the expression
\begin{eqnarray*}
	\hdeg(M)& = &\deg(M)   +
 \sum\limits_{i=1}^{d} {{d-1}\choose{i-1}}\cdot
 \deg(\Ext_{\SS}^i(M,\SS))  \\
& = & \deg(M) +
 \sum\limits_{i=1}^{d} {{d-1}\choose{i-1}}\cdot
 \deg(\Ext_{\SS}^i(M,\SS)  \\
\end{eqnarray*}
in multiplicities of all the $\Ext$ modules. The multiplicity
$\deg(\Ext_{\SS}^i(M,\SS))$ can be written as the ordinary multiplicity of
one of the factors of the filtration.
}\end{example}

\subsubsection*{Generalized Cohen-Macaulay modules}

The expression for the function $\hdeg(\cdot)$ arrives in known
territory if $M$ is a generalized Cohen-Macaulay module.

\begin{Definition}\index{generalized
Cohen-Macaulay module}
 Let $(\RR,\m)$ be a Noetherian local  ring and $M$ a
finitely generated $\RR$-module. $M$ is a {\em generalized
Cohen-Macaulay module} if $M_{\p}$ is Cohen-Macaulay for all prime
ideals $\p\neq \m$.

\end{Definition}

A typical characterization is:

\begin{Theorem}
 Let $(\RR,\m)$ be a Noetherian local  ring and $M$ a
finitely generated $\RR$-module. $M$ is a {\em generalized
Cohen-Macaulay module} if and only if the modules
$\rmH_{\m}^i(M)$ are finitely generated for all $i<\dim M$.
\end{Theorem}

For these modules the expression
 $\hdeg(M)$, by local duality, converts into:

\begin{proposition} \label{hdeggcm} If $M$ is a generalized
Cohen-Macaulay module of dimension $d$, then
\begin{eqnarray*}
 \hdeg(M) & = &\deg (M)
+ \sum\limits_{i=0}^{d-1} {{d-1}\choose{i}}\cdot
 \lambda(\H^{i}_{\mathfrak{m}}(M)).
\end{eqnarray*}
\end{proposition}

\subsubsection*{Buchsbaum modules} This is an important class of
rings introduced by St\"uckrad and Vogel.
They are studied in great detail in \cite{StuckradVogel}.

\begin{Definition}\index{Buchsbaum module}
 Let $(\RR,\m)$ be a Noetherian local  ring and $M$ a
finitely generated $\RR$-module. $M$ is a {\em Buchsbaum
 module} if  the partial Euler characteristic
 $\chi_1(\xx;M)$\index{Euler number} is
 independent of the  system of parameters $\xx$ for $M$. In other
 words, $\lambda(M/(\xx)M)-\deg_{\xx}(M)$ is independent of $\xx$.
This integer is called the {\em Buchsbaum
invariant}\index{Buchsbaum invariant} of $M$.
\end{Definition}

Two
of its properties are:

\begin{Theorem}
 Let $(\RR,\m)$ be a Noetherian local  ring and $M$ a
finitely generated $\RR$-module. Then
\begin{enumerate}
\item[{\rm (a)}] $M$ is a Buchsbaum module if and only if every system of
parameters for $M$ is a $d$-sequence relative to $M$.

\item[{\rm (b)}] If $M$ is a Buchsbaum module then $\H_{\m}^i(M)$ are
$\RR/\m$-vector spaces for $i< \dim M$. $($In particular they are
generalized Cohen-Macaulay modules.$)$ The converse does not hold
true.
\end{enumerate}
\end{Theorem}

\begin{Corollary}
 If $M$ is a Buchsbaum module, then
\[ \hdeg(M) = \deg(M) + I(M),\] where $I(M)$ is the Buchsbaum
invariant of $M$.
\end{Corollary}

An effective test for the property is given by the following theorem of
K. Yamagishi (\cite{Ya91}).\footnote{We thank T. Kawasaki for the
reference.}

\begin{Theorem} Let $(\RR, \m)$ be a regular local ring of dimension $n$.
A finitely generated $\RR$-module $M$ is a Buchsbaum module
if and only if its  Bass numbers are given by
\begin{eqnarray*}
\lambda(\Ext_{\RR}^i(\RR/\m, M)) &=&
\sum_{j=0}^i\binom{n}{n-j}\lambda(\H_{\m}^j(M)) \\
&=& \sum_{j=0}^i\binom{n}{i-j}\lambda(\Ext_{\RR}^{n-j}( M, \RR))
\end{eqnarray*} for all $i<\dim M$.
\end{Theorem}

\begin{Remark}{\rm
Note that this requires that $\Ext_\RR^i(M, \RR)$, for $i>0$, be a
module of finite support. In the case of Buchsbaum modules, it
suffices to check that $\m \cdot \Ext_\RR^i(M, R)=0$, $i>0$.
\medskip

A related theme is how to decide whether an $\RR$-module $M$ is
Cohen-Macaulay on the punctured spectrum. If $\RR$ is a regular local
ring and the annihilator of $M$ is equidimensional a criterion can be
cast as follows. Suppose $\codim M=r$ and
\[ \cdots \lar F_{r+1} \stackrel{\varphi}{\lar} F_r \lar \cdots \]
is a projective resolution of $M$.
For $\p\neq \m$, $M_{\p}$ is Cohen-Macaulay if and only if $\mbox{\rm
proj dim }_\RR M_{\p}=r$, a condition equivalent to saying that the
image of $\varphi_{\p}$ splits off $(F_r)_{\p}$. In other words that the
ideal $I(\varphi)$ of maximal minors of $\varphi$ is $\m$-primary.
\medskip

It would be less cumbersome to find a more amenable module and make
use of the rigidity of Tor since $\RR$ is a regular local ring. Let
us sketch this with a non amenable module! The following is an
equivalence:
\[ \mbox{\rm $M$ is free on the punctured spectrum} \Leftrightarrow
\mbox{\rm $\Tor_{r+1}^\RR(M,M)$ has finite support.}
\]

}\end{Remark}

\subsection{Homological torsion}\index{homological torsion}
There are
other combinatorial expressions of the terms $\hdeg (\Ext_{\SS}^i(M,\SS))$ that
behave well under hyperplane section. It will turn
out to be useful in certain characterizations of generalized
Cohen-Macaulay modules and in the study of Hilbert coefficients.

\begin{definition} {\rm
Let $\SS$ be a ring of dimension $r$ and $M$ an $\SS$-module of dimension $d\geq 2$. Its {\em homological
torsion} is the integer
\[ \TT(M)
=
  \sum\limits_{i=r-d+1}^{r-1} {{r-2}\choose{i-d+r-1}}\cdot
 \hdeg(\mbox{\rm Ext}_{\SS}^i(M,\SS)).\]
If $I$ is an $\m$-primary, the integer obtained using
$\hdeg_I$ will be denoted $\TT_I(M)$.
\medskip

If $r=d$, this formula becomes
\[ \TT(M)
=
  \sum\limits_{i=1}^{d-1} {{d-2}\choose{i-1}}\cdot
 \hdeg(\mbox{\rm Ext}^i_{\SS}(M,\SS)).\]

}\end{definition}

If $d=2$, $\TT(M)= \hdeg(\Ext^1_{\SS}(M,\SS))= \deg(\Ext_{\SS}^1(M,\SS))
+ \deg(\Ext_{\SS}^2(\Ext_{\SS}^1(M,\SS), \SS))
$.

\medskip

The restriction to $d\geq 2$ explanation:
$\lambda(M)$ if $d=0$, or $\lambda(\H_{\m}^0(M))$ if $d=1$
are
quantities we liken to the ordinary torsion.

\begin{example}{\rm If $M$ is a generalized Cohen-Macaulay module of
dimension $d\geq 2$,
\[ \TT_I(M)\leq \sum_{i=1}^{d-1}{d-2\choose
i-1}\lambda(\H^i_{\mathfrak{m}}(M)),\]
with equality if $M$ is Buchsbaum.
}\end{example}

\begin{Theorem} \label{torsionhdeg} Let $M$ be a module of dimension
$d\geq 3$ and let
$h$ be a generic hyperplane section. Then
 $\TT(M/hM)\leq \TT(M)$.
\end{Theorem}

 \begin{proof} We use induction on $d$. If $d\geq 3$,
using the notation of Lemma~\ref{hdegmmn}, in particular setting
$N=M/hM$, we have
\[ \TT(N) =  \sum\limits_{i=1}^{d-1} {{d-2}\choose{i-1}}\cdot
 \hdeg(\mbox{\rm Ext}^i_{\SS}(N,\SS)).
\]
By Lemma~\ref{hdegmmn},
\[ \hdeg(\Ext^i_{\SS}(N,\SS)) \leq \hdeg(\Ext^i_{\SS}(M,\SS))+
\hdeg(\Ext^{i+1}_{\SS}(N,\SS)),
 \]
which by rearranging gives $\TT(N)\leq \TT(M)$, as desired.
\end{proof}

Higher order {\em homological torsion}\index{higher order homological
torsion} functions can be similarly defined. For this reason it is
appropriate to write $\TT^{(1)}(M):= \TT(M).$

\begin{definition} {\rm Let $M$ be a $\SS$-module of dimension $d\geq
3$. Its {\em second order  homological
torsion} is the integer
\[ \TT^{(2)}(M)
=
  \sum\limits_{i=1}^{d-2} {{d-3}\choose{i-1}}\cdot
 \hdeg(\mbox{\rm Ext}^i_{\SS}(M,\SS)).\]
If $I$ is an $\m$-primary, the integer obtained using
$\hdeg_I$ will be denoted $\TT_I^{(2)}(M)$.

\medskip

More generally, let $M$ be a module of dimension $d$. For each $j\leq
d$,
 set
\[ \TT^{(j)}(M)
=
  \sum\limits_{i=0}^{d-j} {{d-j-1}\choose{i-1}}\cdot
 \hdeg(\mbox{\rm Ext}^i_{\SS}(M,\SS)).\]

}\end{definition}

\begin{Theorem} \label{torsion2hdeg} Let $M$ be a module of dimension
$d\geq 4$ and let
$h$ be a generic hyperplane section. Then
 $\TT^{(2)}(M/hM)\leq \TT^{(2)}(M)$.
A similar assertion holds for all $\TT^{(j)}.$
\end{Theorem}

\begin{corollary} \label{allthdegs}
Let $M$ be a module of dimension $r=d+1\geq 3$. Then
\[ \hdeg(M)> \TT^{(1)}(M) \geq \TT^{(2)}(M) \geq \cdots \geq \TT^{(d)}(M).\]
\end{corollary}

\begin{remark} \label{rigidity1}
{\rm
There are a number of rigidity questions about the values of the
Hilbert coefficients $\rme_i(I)$. Typically they have the form
\[ |e_1(I)| \geq |e_2(I)| \geq \cdots \geq |e_d(I)|.\]
Two of these cases are (i) parameter ideals and (ii) normal ideals,
or more generally the case of the normalized filtration. Because
values of $\TT^{(i)}_I(\RR)$ have been used, in a few cases,  to bound the $\rme_i(I)$,
the descending chain in Corollary~\ref{allthdegs} argue for an
underlying rigidity.
This runs counter to the known formulas for the values of the
$\rme_i(I)$ for general ideals (e.g. (\cite[Theorem 4.1]{RTV})).
}
\end{remark}

\subsubsection*{Generalized Cohen-Macaulay modules} These modules can be
characterized in terms of their homological torsions.

\begin{proposition} \label{hdegvsgcm} Let $(\RR, \m)$ be a Gorenstein local ring of
dimension $d\geq 2$ and $M$ a finitely generated $\RR$-module. Then $M$ is
a generalized Cohen-Macaulay module if and only if
$T_I(M)$ is bounded
 for every $\m$-primary ideal $I$.
More precisely, if $I=(x_1, \ldots, x_d)$ is a parameter
ideal\index{parameter ideal}, and
$I_n=(x_1^n, \ldots, x_d^n)$, it suffices that all $T_{I_n}(M)$ be
bounded.
\end{proposition}

\begin{proof} If $d=2$, $T(M)=\hdeg (\Ext_{\RR}^1(M,\RR))$. This is a module of
dimension at most one. If the dimension is $1$,
$\deg_{I_n}(\Ext_{\RR}^1(M,\RR) $ is a polynomial of degree $1$ for
large $n$.

\medskip

If $d\geq 3$, let $h$ be a generic hyperplane section. It will
suffice to show that $M/hM$ is a generalized Cohen-Macaulay module.
For that
we can use Theorem~\ref{torsionhdeg} and induction. 
\end{proof}

\subsubsection*{Castelnuovo-Mumford regularity and homological degree}
The following result of Chardin, Ha and Hoa (\cite{CHH9} provides a
polynomial bound for the homological degree in terms of the
Castelnuovo-Mumford regularity.
(For $d=2$, a similar bound had been found by Gunston in his thesis
(\cite{Gunston}).)

\begin{Theorem} \label{CHH5.1}
 Let $M$ be a non-zero
finitely generated graded $\RR$-module of dimension $d > 0$. Denote by
$n$ the minimal number of generators of $M$ and by $\alpha(M)$ the
maximal degree of the a minimal set of homogeneous generators. Then
\begin{eqnarray}
\hdeg(M) \leq \left[n{{\reg (M)-\alpha(M) + n}\choose{n}}\right]^{2^{(d-1)^2}}.
\end{eqnarray}
\end{Theorem}

\subsubsection*{Specialization and torsion} \label{torsion and
specialization}  One of the uses of extended
degrees is the following. Let $M$ be a module and $\xx=\{ x_1,
\ldots, x_r\}$ be a superficial sequence for the module $M$ relative
to  an extended degree
$\Deg$. How to estimate the length of $\H^0_{\mathfrak{m}}(M)$ in
terms of the initial data of $M$?

Let us consider the case of $r=1$. Let $H=\H^0_{\mathfrak{m}}(M)$
and write
\begin{eqnarray} \label{hmm} 0 \rar H \lar M \lar M' \rar 0.
\end{eqnarray}
Reduction modulo $x_1$ gives
the exact sequence
\begin{eqnarray} \label{hhmm}
 0 \rar H/x_1H \lar M/x_1 M\lar M'/x_1 M'\rar 0.
\end{eqnarray}
From the first sequence we have $\Deg(M)=\Deg(H)+\Deg(M')$, and from
the second
\[
 \Deg(M/x_1M) - \Deg(H/x_1H) = \Deg(M'/x_1M') \leq \Deg(M').\]

Taking local cohomology of the second exact sequence yields the short exact
sequence
\[ 0 \rar H/x_1H \lar \H^0_{\mathfrak{m}}(M/x_1 M)
\lar \H^0_{\mathfrak{m}}(M'/x_1 M')\rar 0,\]
from we have the estimation
\begin{eqnarray*} \Deg( \H^0_{\mathfrak{m}}(M/x_1 M) )& =&
 \Deg(H/x_1H) + \Deg( \H^0_{\mathfrak{m}}(M'/x_1 M')) \\
& \leq &  \Deg(H/x_1H) + \Deg( M'/x_1 M') \\
& \leq &  \Deg(H) + \Deg( M') = \Deg(M). \\
\end{eqnarray*}

We resume these observations as:

\begin{proposition} \label{specialtor} Let $M$ be a module and let $\{x_1, \ldots, x_r\}$
be a superficial sequence relative to $M$ and $\Deg$. Then
\[\lambda(\H^0_{\mathfrak{m}}(M/(x_1, \ldots, x_r)M))
\leq \Deg(M).
\]
\end{proposition}

Now we derive a more precise formula using $\hdeg$. It will be of use
later.

\begin{Theorem} \label{Degreddim1} Let $M$ be a module  of
dimension  $d\geq 2$
and let $\xx=\{x_1, \ldots, x_{d-1}\}$ be a superficial sequence for
$M$ and $\hdeg$. Then
\[
\lambda(\H^0_{\mathfrak{m}}(M/(\xx)M)) \leq
\lambda(\H^0_{\mathfrak{m}}(M))+ \TT(M).
% \sum\limits_{i=1}^{d-1} {{d-2}\choose{i-1}}\cdot
% \hdeg(\mbox{\rm Ext}^i_{\SS}(M,\SS))\leq \hdeg(M)-\deg(M). \nonumber%\\
\]
\end{Theorem}

\begin{proof}
Consider the exact sequence
\[ 0 \rar H = \H^0_{\mathfrak{m}}(M) \lar M \lar M' \rar 0.\]
We have $\Ext_{\SS}^i(M,\SS) = \Ext_{\SS}^i(M',\SS)$ for $d>i\geq 0$, and
therefore $T(M) = T(M')$. On the other hand, reduction mod $\xx$ gives
\begin{eqnarray*}
\lambda(\H^0_{\mathfrak{m}}(M/(\xx)M)) &\leq &
\lambda(\H^0_{\mathfrak{m}}(M'/(\xx)M')) +
\lambda(H/(\xx)H) \\
&\leq &
\lambda(\H^0_{\mathfrak{m}}(M'/(\xx)M')) +
\lambda(H),
\end{eqnarray*}
which shows that it is enough to prove the assertion for $M'$.

\medskip

If $d>2$, we apply Theorem~\ref{torsionhdeg}, to pass to $M'/x_1M'$.
This reduces all the way to the case $d=2$. Let $M$ be a module of
positive depth.
Write  $h=\xx$. The assertion requires that
$\lambda(\H^0_{\mathfrak{m}}(M/hM)) \leq \hdeg(\Ext_{\SS}^1(M,\SS)).$
We have the cohomology exact sequence
\[ \Ext^1_{\SS}(M,\SS) \stackrel{h}{\lar} \Ext^1_{\SS}(M,\SS)
\lar \Ext^2_{\SS}(M/hM,\SS)\lar \Ext^2_{\SS}(M,\SS)=0,
\] where
\[ \lambda(\H^0_{\mathfrak{m}}(M/hM)) =\hdeg(\Ext_{\SS}^2(M/hM,\SS)).\]

If $\Ext_{\SS}^1(M,\SS)$ has finite length the assertion is clear.
Otherwise $L=\Ext_{\SS}^1(M,\SS)$ is a module of dimension $1$ over a
discrete valuation domain $V$ with $h$ for its parameter. By the
fundamental theorem for such modules,
\[ V = V^r \oplus (\bigoplus_{j=1}^s V/h^{e_j}V), \]
so that multiplication by $h$ gives
\[ \lambda(L/hL) = r+s \leq r+ \sum_{j=1}^s e_j = \hdeg(L).\]

An alternative argument at this point is to consider the exact
sequence (we may assume $\dim
\SS=1$)
\[ 0 \rar L_0 \lar L \stackrel{h}{\lar} L \lar L/hL \rar 0, \]
where both $L_0$ and $L/hL$ have finite length. If
$F$ denotes the image of the multiplication by $h$ on $L$, we have
the exact sequences
$ 0 \rar L_0 \lar L \lar F \rar 0 $ and $0 \rar F \lar L \lar L/hL
\rar 0$.
 Dualizing
we have $\Hom_{\SS}(L, \SS)= \Hom_{\SS}(F, \SS) $  and
the  exact sequence
\[ 0 \rar \Hom_{\SS}(L,\SS) \stackrel{h}{\lar} \Hom_{\SS}(L,\SS) \lar
\Ext_{\SS}^1(L/hL, \SS) \lar \Ext_{\SS}^1(F,\SS),  \]
which shows that
\[ \lambda(L/hL) \leq \deg(L) + \lambda(\H_{\m}^0(F))
\leq \deg(L) + \lambda(\H_{\m}^0(L))
 = \hdeg(L), \]
as desired.
\end{proof}

\subsection{Cohomological degrees and Samuel multiplicities}
There are  variations of cohomological degree functions   obtained by using Samuel's
notion of multiplicity.
\index{cohomological degrees and Samuel multiplicities}

\begin{definition}{\rm
Let $(\RR, \mathfrak{m})$ be a Noetherian
local ring and $I$ is an $\mathfrak{m}$-primary ideal.
We denote by $\hdeg_I(\cdot)$ the function obtained
 by replacing,
 in the
definition of $\hdeg(\cdot)$, $\deg(M)$
by Samuel's $e(I;M)$. \index{\textrm{$\hdeg_I$}}
}
\end{definition}

Let us make a rough comparison between $\hdeg(M)$ and $\hdeg_I(M)$.

\begin{proposition}
Let $(\RR,\mathfrak{m})$ be a Noetherian local ring and let $I$ be an
$\mathfrak{m}$-primary ideal. Suppose $\mathfrak{m}^r\subset I$. If
$M$ is an $\RR$-module of dimension $d$, then
\[ \hdeg_I(M) \leq r^d\cdot \deg(M) +  r^{d-1}\cdot(\hdeg(M)-\deg(M)).
\]
\end{proposition}

\begin{proof} If $r$ is the index of nilpotency of $\RR/I$, for any $\RR$-module
$L$ of dimension $s$,
\[ \lambda(L/(\mathfrak{m}^r)^n L)\geq \lambda(L/I^nL).\]
The Hilbert polynomial of $L$ gives
\[  \lambda(L/(\mathfrak{m}^r)^n L)= \deg(M) {\frac{r^s}{s!}}n^s +
\mbox{\rm lower terms}.\]
We now apply this estimate to the definition of $\hdeg(M)$, taking
into account that its terms are evaluated  at modules of decreasing
dimension.
\end{proof}

For later reference we rephrase Proposition~\ref{specialtor} for use
with Samuel's multiplicities.

\begin{Theorem} \label{boundtorsion} Let $(\RR,\mathfrak{m})$ be a
Noetherian local ring and let $I$ be an $\mathfrak{m}$-primary ideal,
and let $M$ be a finitely generated
$\RR$-module of dimension $d\geq 1$.
Let $\xx=\{x_1, \ldots, x_r\}$ be a superficial sequence in $I$
relative to $M$ and $\hdeg_I$. Then
\[ \hdeg_I(M/(\xx)M) \leq \hdeg_I(M).\] Moreover, if $r<d$ then
\[ \lambda(\H^0(M/(\xx)M))\leq \hdeg_I(M)-e(I;M).\]

\end{Theorem}

\begin{corollary}\label{boundRmodJ} Let $(\RR, \mathfrak{m})$ be a
Noetherian local ring of dimension $d>0$ and infinite residue field. For an
$\mathfrak{m}$-primary $I$ there is a minimal reduction $J$ such that
\[ \lambda(\RR/J)\leq \hdeg_I(\RR).\]
In particular, the set of all $\lambda(\RR/J)$ for all parameter
ideals with the same integral closure is bounded.
\end{corollary}

\subsection{The extreme cohomological  degree: bdeg}

In his thesis (\cite{Gunston}), Tor Gunston introduced the following
cohomological degree:

\begin{definition} \label{bdegdef}{\rm Let $\RR$ be a Noetherian local
ring
with infinite residue field or a standard  graded algebra over an
infinite field. For a finitely generated $\RR$-module (graded if
required) $M$,  $\bdeg(M)$ is the integer
\[ \bdeg(M) = \min \{\Deg(M) \mid  \Deg \ \mbox{\rm is a
cohomological degree function}\}.\]
}\end{definition}

This function
 is well-defined since cohomological degrees exist (e.g. $\hdeg$).
  It is
 obviously a cohomological degree itself. Gunston proved that
\cite[Theorem 3.1.3]{Gunston}):

\begin{Theorem}[Gunston]\label{Gunstonth}
Suppose that $M$ is a finitely generated $\RR$-module of positive depth.
Then for a generic hyperplane section $h$,
\[
\bdeg(M) = \bdeg(M/hM).\]
\end{Theorem}

\begin{Corollary} \label{bdegvslength} If $M$ is a finitely generated $\RR$-module of dimension
$d$ there are there is a generic superficial sequence $\xx=\{x_1,
\ldots, x_d\}$ for $M$ such that
\[ \bdeg(M) \leq \lambda(\H_{\m}^0(M))+
\sum_{i=1}^d \lambda(\H_{\m}^0(M/(x_1, \ldots, x_i)M))\leq
(d+1)\hdeg(M)-d \deg(M).\]

\end{Corollary}

\begin{proof} The first inequality is a direct consequence of
Theorem~\ref{Gunstonth}, the second from Theorem~\ref{boundtorsion}.
\end{proof}

%There are few effective formulas for $\bdeg$. Here is one that holds
%in dimension two. It takes into account
% the variability of $\rme_1(\cdot)$ with the parameter ideal
%(see
%Chapter~\ref{chapchern}).

%\begin{Theorem}
%Let $(\RR, \m\}$ be a Noetherian local ring and let $M$ be a module
%of dimension $2$. For parameter ideal $I$,
%\[\bdeg_I(M)= e_0(I;M) -\rme_1(I;M) + e_2(I;M).
% \]
%\end{Theorem}

%\begin{proof}
%Let $M'=M/\H_{m}^0(M)$. By Proposition~\ref{genhs} it suffices to
%prove the equality for $M'$. (Recall that
%$\lambda(\H_{\m}^0(M))=e_2(I;M)$.)
%Using the notation of
%Proposition~\ref{genhs},
%\begin{eqnarray*} \bdeg_I(M)= \bdeg_I(M'/hM')&=& \deg_I(M'/hM')+
%\lambda(\H_{\m}^0(M'/hM'))\\
%&=& e_0(I;M')-e_1(I;M').
%\end{eqnarray*}
%\end{proof}
\bigskip

\subsubsection*{Rules of computation} One of the difficulties of computing cohomological degrees lies on
their behavior on short exact sequences. The degree $\bdeg$ has
slightly more amenable properties, according to
the following rules.

\begin{proposition}\label{Gu2} Let $\RR$ be a standard graded algebra
and $A$ a finitely generated graded $\RR$-module with Hilbert
function $H_A(\tt)$. Then
\[ \bdeg(A)\leq \sum_{j=\alpha(A)}^{\reg(A)} H_A(j).\]

\end{proposition}

\begin{proof}
Set $L = \H_{\m}^0(A)$ and consider the exact sequence
\[ 0 \rar L \lar A \lar A'= A/L \rar 0.\]
We may assume $\dim A\geq 1$. Let $h$ be a generic hyperplane
section. The proof follows by induction from   the following inequalities:
\begin{eqnarray*}
\reg(A') & = & \reg(A'/hA'), \\
\alpha(A) & \geq & \alpha(A'/hA'),\\
H_A(j) &=& H_L(j) + H_{A'}(j) \geq H_L(j) + H_{A'/hA'}(j),  \\
\bdeg(A) &=& \lambda(L) + \bdeg(A') = \lambda(L) + \bdeg(A'/hA')\\
\reg(A)  &=& \max\{\reg(L),  \reg(A')\}= \max\{\reg(L),
\reg(A'/hA')\}.
\end{eqnarray*}

\end{proof}

\begin{Proposition}$($\cite[Proposition
3.2.2]{Gunston}$)$\label{Gunston322}
 Suppose $A,B,C$ are finitely generated
$\RR$-modules and
\[ 0 \rar A \lar B \lar C \rar 0\]
is an exact sequence. If $\lambda(C)<\infty$,  then
\[ \bdeg(B)\leq \bdeg(A)+
\lambda(\Gamma_{\m}(B))-\lambda(\Gamma_{\m}(A)).\]
In particular, we have
\[\bdeg(B)\leq \bdeg(A)+\lambda(C), \quad \mbox{and if }\ \depth
B>0, \ \mbox{then  } \bdeg(B)\leq \bdeg(A).  \]
\end{Proposition}

A calculation in \cite[Proposition 3.2]{Dalili} then leads to:

\begin{Proposition} \label{Dalili3.2} Suppose $A,B,C$ are finitely generated
$\RR$-modules and
\[ 0 \rar A \lar B \lar C \rar 0\]
is an exact sequence. Then
\begin{enumerate}
\item[{\rm (a)}] $\bdeg(B)\leq \bdeg(A)+\bdeg(C)$;
\item[{\rm (b)}] $\bdeg(A)\leq \bdeg(B) + (\dim A-1)\bdeg(C)$.
\end{enumerate}
\end{Proposition}

\begin{proof}
We prove only (1).
The case where $A$
 has finite length is a consequence of $\bdeg$ being an extended
 degree. For the other cases,
one make use of indiction on the dimension, repeated use of snake
lemma and the following special case.
\end{proof}

\begin{Lemma}
Suppose $A,B,C$ are finitely generated
$\RR$-modules and
\[ 0 \rar A \lar B \lar C \rar 0\]
is an exact sequence. If $C$ has finite length then
\[
\bdeg(B)\leq \bdeg(A)+\bdeg(C).\]
\end{Lemma}

\begin{proof} For a finitely generated $\RR$-module, denote by $A_0$
its submodule of finite support; set $A'=A/A_0$.
If $A=A_0$, the assertion follows from the additivity of the length
 function. Clearly we may replace $A$ by $A/A_0$ and $B$ by $B/A_0$;
 thus changing notation  we may assume that $A$ has positive depth.
If follows that $B_0$ embeds in $C$. We may now replace $B$ by
$B/B_0$ and $C$ by $C/B_0$.

Changing notation again, we may assume that both $A$ and $B$ have
positive depth. Let $h$ be an appropriate hyperplane section for $A$
and $B$, that is $\bdeg(A)= \bdeg(A/hA)$ and $\bdeg(B)=\bdeg(B/hB)$.
Reduction mod $h$ gives the exact sequence
\[ 0 \rar E \lar A/hA \lar B/hB \lar C/hC\rar 0,\]
and thus by induction ($E$ has finite length and
$\bdeg(E)=\bdeg(C/hC)$),
\begin{eqnarray*} \bdeg(B)=\bdeg(B/hB) &\leq& \bdeg((A/hA)/E) +
\bdeg(C/hC)\\
&\leq & \bdeg(A) + \bdeg(C),
\end{eqnarray*}
as desired.
\end{proof}

We continue with the proof of Proposition~\ref{Dalili3.2}. If $C$ has
positive depth, $A_0=B_0$ and as in the Lemma we replace $A$ and $B$
by $A/A_0$ and $B/B_0$, respectively. Picking an appropriate
hyperplane section $t$, reduction gives an exact sequence with the
same $\bdeg$s but in lower dimension.
\medskip
Consider the diagram

\[
\diagram
& & & 0\dto & \\
& 0\dto & & C_0\dto & \\
0 \rto & A \rto\dto  & B
\rto\dto  & C \rto\dto & 0 \\
0 \rto & A_1 \rto\dto  & B \rto & C/C_0 \dto\rto & 0\\
& C_0\dto & & 0 & \\
& 0 & &    &
\enddiagram .
\]
Since $C/C_0$ has positive depth (or vanishes),
\[ \bdeg(B) \leq \hdeg(A_1) + \bdeg(C/C_0)=
\hdeg(A_1) + \bdeg(C)- \bdeg(C_0),
\] while the Lemma gives
\[ \bdeg(A_1) \leq \bdeg(A) + \bdeg(C_0),\]
which proves the assertion. \QED

Unfortunately the  delicate task of bounding $\bdeg(C)$ in terms of
$\bdeg(B)$ and $\bdeg(C)$ is
not possible: Just consider a module with a free resolution
\[ 0 \rar F \stackrel{\varphi}{\lar} F \lar C \rar 0.\] Then
$\hdeg(C) = \deg(\RR/(\det(\varphi)))$.

\begin{Example}{\rm Let $G$ be a graph on a vertex $V$ set indexed by the set of
variables and edge set $E$ and set $I=I(G)$ its edge ideal.
\begin{itemize}
\item[{\rm (a)}]
 Let $I$ be
monomial ideal defined by the complete graph $K_n$: $I=(x_ix_j, i<
j)$. Consider the exact sequence (a so-called {\em deconstruction
sequence} in the terminoloy of R. Villarreal)
\[ 0 \rar \RR/I:x_{n} \lar \RR/I \lar \RR/(I,x_n) \rar 0.
 \]
$I:x_{n}= (x_1, \ldots, x_{n-1})$ and $\RR/(x_n, I)= \RR'/I'$, where
$\RR'=k[x_1, \ldots, x_{n-1}]$ and $I'$ is the ideal corresponding to
the graph $K_{n-1}$.
 Using (\ref{Dalili3.2}) and induction we get
\[ \bdeg(\RR/I)= \deg(\RR/I)=n, \quad \bdeg(I) \leq 1 +n(n-1).\]

\item[{\rm (b)}] For a graph with a minimal vertex cover $\{x_1, \ldots,
x_c\}$--in other words $(x_1, \ldots, x_c)$ is a minimal prime of $I$
of maximal dimension--one has to repeat the deconstruction step
above $c-1$ times: From
\[ I  = x_1I_1 + x_2I_2 +  \cdots + x_c I_c, \]
\[ 0 \rar \RR/J_1 \lar \RR/I \lar \RR/(I,x_1) \rar 0
 \] that is
\[ 0 \rar \RR/J_1 \lar \RR/I \lar \RR/(I,x_1) \rar 0\]
$J_1 = (I_1 + x_2I_2 +  \cdots + x_c I_c).$

%\[\bdeg(\RR/I)\leq {{n+e-2}\choose {e-1}}, \quad
\[\bdeg(\RR/I)\leq {c+1}, \quad
\bdeg(I) \leq 1 + (n-1)(c+1).\]

\end{itemize}

}\end{Example}

\begin{Proposition} Let $I$ be a monomial ideal of $\RR$ generated by
the set of monomials $\{m_1, \ldots, m_r\}$ of degree $\leq d$. Then
\[ \bdeg(\RR/I)\leq {{n+d-1}\choose{d}}.\]
% \sum_1^r \deg(m_i) -r+1.\]
\end{Proposition}
\begin{proof}
Denote $f(p,q)$ a bound for $\bdeg(\RR/I)$ valid for ideals generated
by monomials of degree $p<d$ or in $p<n$ variables. The
reconstruction sequence and  Proposition~\ref{Dalili3.2} assert that
we can
take
\[ f(d,n)= {{n+d-1}\choose{d}}.\]

As the case of graph ideals indicate this often overstates the value
of $\bdeg(\RR/I)$.
 We use induction on the number $r$ of monomials. If $x_n$ is a
variable present in one of the monomials $m_i$, say $m_i=x_nm_i'$,
$i\leq s$ and $ m_i\notin (x_n)$ for $i> s$, $I:x_n=(m_1', \ldots,
m_s')$, we have the exact sequence
\[ 0 \rar \RR/I:x_{n}=\RR/(m_1', \ldots, m_s') \lar \RR/I \lar
\RR/(I,x_n)= \RR'/(m_{s+1}, \ldots, m_r) \rar 0.
 \]
\end{proof}

This gives, using Theorem~\ref{Degandbetti}, the following uniform
bound for the Betti numbers of monomial ideals.

\begin{Corollary} Let $I=(m_1, \ldots, m_r)$ be a monomial ideal of
$\RR=k[x_1, \ldots, x_n]$, $\deg m_\leq d$. Then the Betti numbers of $\RR/I$ are
bounded by
\[ \beta_i(\RR/I)
\leq \bdeg(\RR/I)\cdot {{n}\choose{i}}\leq
 {{n+d-1}\choose{d}}\cdot  {{n}\choose{i}}.\]
 \end{Corollary}

A more interesting calculation is:

\begin{Theorem} Let $\RR$ be a Buchsbaum local ring with infinite
residue field  and $\xx=\{x_1,
\ldots, x_d \}$ be a system of parameters. Then
\[ \bdeg_{(\xx)}(\RR)= \lambda(\RR/(\xx)).\]
\end{Theorem}

\begin{proof} We use induction with basic properties of $d$-sequences
(Proposition~\ref{dseqprop}):

\begin{itemize}
\item[{\rm (a)}] Set $I=(\xx)$. Since $\H_{\m}^0(\RR)=0:x_1$, consider the exact
sequence
\[ 0 \rar (0:x_1) \lar \RR \lar \RR'\rar 0.\]
From (\ref{dseqprop})  and  basic properties of Buchsbaum rings (see
\cite{StuckradVogel}), we have that  $(0:x_1)\cap I = 0$ and $\RR'$ is a Buchsbaum ring
of positive depth.

\item[{\rm (b)}] In the equality $\bdeg(\RR) = \lambda(\H_{\m}^0(\RR))+
\bdeg(\RR')$, we make use of another property $\bdeg(\cdot)$:
There
exists a generic element $x\in I'=I\RR'$ such that
\[ \bdeg(\RR')= \bdeg(\RR'/(x)).
\]

\item[{\rm (c)}] If $d>1$, $\RR'$ is a Buchsbaum ring of dimension $d-1$. Thus
by induction $\bdeg(\RR'/(x))= \lambda(\RR'/I')$.

\item[{\rm (d)}] Let us return to the exact sequence above (and the case $d=1$).
Tensoring it by $\RR/I$, we get the complex

\[  (0:x_1)/I(0:x_1) \lar \RR/I \lar \RR'/I'\rar 0.\]
Since $(0:x_1)\cap I$, we obtain the exact sequence
\[ 0 \rar (0:x_1) \lar \RR/I \lar \RR'/I'\rar 0,\]
\end{itemize}

\noindent which gives the formula $\bdeg_I(\RR)=\lambda(\RR/I)$. 
\end{proof}

\begin{Corollary} Let $\RR$ be a Buchsbaum local ring. For any system
of parameters $\xx=\{x_, \ldots, x_d\}$, $\lambda(\RR/(\xx))$ depends
only on the integral closure of $(\xx)$.
\end{Corollary}

The following follows from a general property of $\Deg$:

\begin{Corollary} Let $\SS$ be a regular local ring of dimension $n$
and $\RR=\SS/L$ a Buchsbaum ring of dimension $d$. If $\xx=\{x_1,
\ldots, x_d\}$ is a system of parameters for $\RR$, then the Betti
numbers of $\RR$ are bounded by
\[ \beta_i(\RR)\leq \lambda(\SS/(L,\xx))\cdot {n\choose i}.\]
\end{Corollary}

\subsubsection*{Semiadditive degree functions} If $\RR$ is a Noetherian
local ring a {\em semiadditive degree
function}\index{semiadditive degree function} is a mapping
\[ {\bf d}:\mathcal{M}(\RR)\rar \mathbb{Q}\]
such that if $0\rar A\rar B \rar C \rar 0$ is a short exact sequence
of modules in $\mathcal{M}(\RR)$ then
\[ {\bf d}(B)\leq {\bf d}(A)+ {\bf d}(C).\] 
Such functions are susceptible, as the example above indicates, to the
derivation of binomial bounds.
\medskip

We have seen several degrees with this property, starting with
$\nu(\cdot)$ (and more generally Bass and Betti numbers) and including $-e_1(\xx;\cdot)$, $\chi_1(\xx;\cdot)$,
$\reg(\cdot)$, $\bdeg(\cdot)$.

\begin{Question}{\rm Is $\hdeg(\cdot)$ semiadditive?}
\end{Question}

\begin{Question} {\rm Let $\RR$ be a Cohen-Macaulay local ring of
dimension $d$. Consider the set of rational numbers
\[ \frac{\hdeg(M)- \hdeg(M/hM)}{\deg(M)}\]
over all $M\in \mathcal{M}(\RR)$ and all generic hyperplane sections. Is
this set finite,  or, more generally,  bounded? Can it be
expressed  as an invariant of $\RR$?
}

\end{Question}

\begin{Question}{\rm
Let $\xx$ be a system of parameters of the Noetherian local ring
$\RR$. If $\xx$ is a d-sequence, find an estimation for
$\hdeg_{(\xx)}(\RR)$.
}
\end{Question}

\begin{Question}{\rm
Let $\RR=k[\Delta]$ be the Stanley-Reisner ring of the simplicial
complex $\Delta$. Find estimations for $\hdeg(\RR)$ and/or
$\bdeg(\RR)$.
}\end{Question}

%\chapter{Applications of Degree Functions}

%\section*{Introduction}

%\section{Counting Hilbert Functions with Multiplicities}

%\subsection{Introduction}

%\subsection{Hilbert polynomials and Castelnuovo-Mumford regularity}

\chapter{Normalization of Algebras, Ideals and Modules}\index{normalization}

\section*{Introduction}

In chapter we treat various metrics to gauge  the passage from an algebraic structure $\AA$ to its integral closure
$\bar{\AA}$. The main topics are:
\begin{itemize}
\item[{$\bullet$}] Tracking numbers
\item[{$\bullet$}] Complexity of normalization of algebras
\item[{$\bullet$}] Normalization of ideals
\item[{$\bullet$}] Sally modules of filtrations
\item[{$\bullet$}] Normalization of modules
\end{itemize}

\section{The Tracking Number of an Algebra}\index{tracking number!graded
module or algebra}\label{trackingnumber}

\subsection{Introduction}
The terminology {\em tracking number} (or {\em twist}) refers to the
use of integers as locators, or tags, for modules and algebras in
partially ordered sets. A forerunner of this use was made in the
previous section,
 when divisorial ideals were employed to bound chains
of algebras with Serre's property $S_2$.

There are immediate generalizations of this notion to more general
gradings. Here we have in mind just applications to graded algebras
that admit a standard Noether normalization, that is, the so-called
{\em semistandard}\label{semistandard graded algebra} graded algebras.
In a brief discussion we treat abstract tracking numbers.

\subsection{The determinant of a module or an algebra: field case}

We now give a more conceptual explanation of the boundedness of the
chains of the previous section  based on another family of divisors attached to the
extensions (\cite{DV}). Its added usefulness will include
extensions to all characteristics.

Let $E$ be a finitely generated graded module
over the ring of polynomial over the field $k,$
 $\RR=k[x_1, \ldots, x_d]$. If $\dim E=d$,
denote by $\det_\RR(E)$ the determinantal divisor of $E$: if $E$ has
multiplicity $e$,
\[ \det(E) = (\wedge^rE)^{**}\cong  \RR[-\delta].\]

\begin{definition}{\rm
 The integer $\delta$ will be called the {\em
tracking number} of $E$: $\delta= \mbox{\rm tn}(E)$. \label{trnu}
}\end{definition}

\subsubsection*{Chern coefficients}\index{Chern coefficient}

We first develop the basic properties of the tracking number. It may help to
begin with these examples:

\begin{example} \label{tnex1} {\rm
 If $\AA$ is a homogeneous domain over a field $k$,  $\RR$  a
homogeneous Noether normalization and $\SS$  a hypersurface ring
over which $\AA$ is birational,
\[ \RR \subset \SS= \RR[t]/(f(t)) \subset \AA\]
($f(t)$ is a homogeneous polynomial of degree $e$), we have
$(\wedge^r\SS)^{**}= \RR[-{{e}\choose{2}}]$. As a consequence $\tn(\AA)\leq
{{e}\choose{2}}$.

Another illustrative example is that of the fractionary ideal $I =
(x^2/y,y)$. This $I$ is positively generated and $I^{**}= \RR(1/y)$:
$\tn(I)=-1$.
}\end{example}

The following observation
shows the use of tracking numbers   to locate
the members of certain chains of modules.

\begin{proposition}
 If
 $E\subset F$ are modules with the same multiplicity that satisfy
 condition $S_2$, then $\tn(E)\geq \tn(F)$ with equality  only if
 $E=F$.
\end{proposition}

\begin{proof} Consider the  exact sequence
\[ 0 \rar E \stackrel{\varphi}{\lar} F \lar C \rar 0,\]
and we shall key on the annihilator of $C$ (the {\em conductor} of $F$
into $E$). Localization at a height $1$ prime $\mathfrak{p}$ gives a
free resolution of $C_{\mathfrak{p}}$ and the equality
\[ \det(F) = \det(\varphi_{\mathfrak{p}})\cdot \det(E).\]
Therefore $\det(E)=\det(F)$ if and only if $C$ has codimension at
least two, which is not possible if $C\neq 0$ since it has depth at
least $1$. 
\end{proof}

\begin{corollary}\label{boundlength} If the  modules in the strictly
increasing  chain
\[ E_0 \subset E_1 \subset \cdots \subset E_n\]
have the same multiplicity and satisfy  condition $S_2$,
then $n\leq \tn(E_0)-\tn(E_n)$.
\end{corollary}

\begin{remark}{\rm
The inequality $\tn(E_0)-\tn(E_n)\geq n$ is rarely an equality, for the
obvious reasons. Let $E\subset F$ be distinct graded modules of the
same multiplicity satisfying $S_2$. Then $\det(E) = f\cdot \det(F)$
where $f$ is a homogeneous polynomial of degree given by $\tn(E)-\tn(F)$. If
\[ f = f_1^{a_1}\cdots f_r^{a_r}\] is a primary decomposition of $f$,
a  tighter measure of the ``spread'' between $E$ and $F$ is the
integer
\[ a_1 + \cdots + a_r,\]
rather than
\[ \tn (E)-\tn(F)=a_1\cdot \deg(f_1)\cdots a_r\cdot \deg(f_r).\]
}\end{remark}

\begin{proposition}\label{div3} If the complex
\[ 0 \rar A \stackrel{\varphi}{\lar} B \stackrel{\psi}{\lar} C \rar 0\]
of finitely generated graded
$\RR$-modules
  is an exact sequence
of free modules in every localization $\RR_{\mathfrak{p}}$ at height
one primes, then
$ \tn(B)= \tn(A) + \tn(C) $.
% In particular, if $B$ is a free
%$\RR$-module, then $\det(A)=-\det(C)$.
\end{proposition}

\begin{proof} We break up the complex into  simpler exact  complexes:
\[ 0 \rar \ker(\varphi) \lar A \lar A'=\textrm{\rm image}(\varphi)  \rar 0 \]

\[ 0 \rar A' \lar \ker(\psi) \lar \ker(\psi)/A' \rar 0 \]

\[0 \rar B'=\textrm{\rm image}(\psi) \lar C \lar C/B' \rar 0 \]
and
\[ 0 \rar \ker(\psi) \lar B \lar B' \rar 0.\]

We note that  $\codim \ker(\varphi)\geq 1$, $\codim
C/B'\geq 2$, $\codim \ker(\psi)/A' \geq 2$ by hypothesis,
 so that we have the
equality of determinantal divisors:
\[\det(A)=\det(A')=
\det(\ker(\psi)), \  \textrm{\rm and} \  \det(C)= \det(B').\] What this all means is
that we may assume that the given complex is exact.

\smallskip

Suppose that $r= \rank(A)$ and $\rank(C)=s$ and set $n=r+s$. Consider the
pair $\wedge^rA$, $\wedge^sC$. For $v_1, \ldots, v_r\in A$ and  $u_1,
\ldots, u_s \in C$, choose $w_i$ in $B$ such that $\psi(w_i)=u_i$ and
consider
\[ v_1 \wedge \cdots \wedge v_r \wedge w_1 \wedge \cdots \wedge
w_s\in \wedge^{n}B. \]
Different choices for $w_i$ would produce elements in $\wedge^{n}B$
that differ from the above by terms that  contain at least $r+1$
factors of the form
\[ v_1 \wedge \cdots \wedge v_r \wedge v_{r+1} \wedge \cdots,\]
with $v_i \in A$. Such products are torsion elements in $\wedge^nB$.
This implies that modulo torsion we have a well-defined pairing:
\[ [\wedge^r A/\textrm{\rm torsion}] \otimes_R [\wedge^s C/\textrm{\rm
torsion}] \lar [\wedge^n B/\textrm{\rm torsion}].
\] When localized at primes $\mathfrak{p}$ of codimension at most $1$,
the complex becomes an exact complex of projective
$\RR_{\mathfrak{p}}$-modules and the pairing is an isomorphism. Upon
taking biduals and the $\circ $ divisorial composition, we obtain the asserted
isomorphism.
\end{proof}

\begin{corollary}\label{tneuler} Let
\[ 0 \rar A_1 \lar A_2 \lar \cdots \lar A_n \rar 0\]
be a complex of graded  $\RR$-modules and
homogeneous homomorphisms that is an exact complex of free modules
in codimension $1$. Then
\[ \sum_{i=1}^n (-1)^i \tn(A_i)=0.\]
\end{corollary}

\begin{proposition} Suppose that  $\RR=k[x_1, \ldots, x_d]$ and let
\[ 0 \rar A \lar B \lar C \lar D \rar 0,\] be an exact sequence of
graded $\RR$-modules and homogeneous homomorphisms. If $\dim B=\dim
C=d$,
$\codim A\geq
1$ and $\codim D\geq 2$, then $\tn(B)=\tn(C)$.
\end{proposition}

\begin{corollary}\label{corsup} If $E$ is a graded $\RR$-module of dimension $d$,
then
 \[\tn(E)= \tn(E/\textrm{\rm mod torsion})=\tn(E^{**}).\]
\end{corollary}

Let $\AA$ be a homogeneous algebra defined over a field $k$ that admits
a Noether normalization $\RR=k[x_1, \ldots, x_d]$; then clearly
$\tn_\RR(\AA)= \tn_{\RR'}(\AA')$, where $K$ is a field extension of $k$, $\RR'=
K\otimes_k\RR$ and $\AA' = K\otimes_k\AA $.
Partly for this reason, we can always define the tracking number of an algebra
by first enlarging the ground field.  Having done that and chosen a
Noether normalization $\RR$ that is a standard graded algebra, it will
follow that $\tn_\RR(\AA)$ is independent of $\RR$: the $\RR$-torsion
submodule $\AA_0$ of $\AA$ is actually an ideal of $\AA$ whose definition
is independent of $\RR$.

\medskip

\subsubsection*{Calculation rules}

We shall now derive several rules to facilitate the computation of
tracking numbers.

\begin{proposition} \label{tne1}  Let $E$ be a finitely generated graded module
over the polynomial ring $\RR=k[x_1, \ldots, x_d]$. If $E$ is
torsion free over  $\RR$, then
$\tn(E)=\e_1(E)$, the first Chern number
 of $E$.
\end{proposition}

\begin{proof}
 Let
\[ 0 \rar \oplus_{j}\RR[-\beta_{d,j}] \rar \cdots \rar \oplus_j
\RR[-\beta_{1,j}] \rar \oplus_j \RR[-\beta_{0,j}]\rar E \rar 0\]
be a (graded) free resolution of $E$. The integer
\[ \e_1(E) = \sum_{i,j}(-1)^{i}\beta_{i,j}\]
is (see \cite[Proposition 4.1.9]{BH})
 the next to the leading
Hilbert coefficient of $E$. It is also the integer that one gets by
taking the alternating product of the determinants in
   the free graded resolution (see Corollary~\ref{tneuler}).
\end{proof}

\begin{example}{\rm Suppose that $\AA = k[(x,y,z)^2]$  a subalgebra of the ring
of polynomials $k[x,y,z]$. Setting the weight of the indeterminates to
$\frac{1}{2}$, $\AA$ becomes  a standard graded algebra with $\RR=k[x^2,y^2,
z^2]$ as a Noether normalization. A calculation easily shows that the
Hilbert series of $\AA$ is
\[ H_\AA(t) =\frac{h_\AA(t)}{(1-t)^3}= \frac{1+3t}{(1-t)^3}.\]
Thus $\deg(\AA) = h_\AA(1) =4$, $\tn(\AA)= \e_1(\AA) = h'_\AA(1) = 3$.
}\end{example}

In general, the connection between the tracking number and the first
Hilbert coefficient has to be `adjusted' in the following manner.

\begin{proposition} \label{e1adj} Let $E$ be a finitely generated graded module
over $\RR=k[x_1, \ldots, x_d]$. If $\dim E=d$ and $\rme_0$ is its torsion
submodule, \[ 0 \rar E_0 \lar E \lar E' \rar 0,\]
then \[ \tn(E)= \tn(E') = \e_1(E')=e_1(E) + \hat{\e}_0(E_0),\]
where $\hat{\e}_0(E)$ is the multiplicity of $\rme_0$ if $\dim E_0=d-1$,
and $0$ otherwise.
\end{proposition}

\begin{proof} Denote by $H_\AA(t)$ the Hilbert series of an $\RR$-module $\AA$
(see \cite[Chap. 4]{BH}) and write
\[ H_\AA(t) = \frac{h_\AA(t)}{(1-t)^d},\]
if $\dim \AA=d$. For the exact sequence defining $E'$, we have
\[ h_E(t) = h_{E'}(t) + (1-t)^r h_{E_0}(t),\]
where $r=1$ if $\dim E_0=d-1$, and  $r\geq 2$ otherwise. Since
\[\rme_1(E) = h_{E}'(1) = h_{E'}'(1) +
r(1-t)^{r-1}|_{t=1}h_{E_0}(1),\]
the assertion follows. 
\end{proof}

\begin{remark}{\rm This suggests a reformulation of the notion of
tracking number. By using exclusively the Hilbert function, the
definition could be extended to all finite modules over a graded
algebra.
}\end{remark}

\begin{corollary} \label{tntimes} Let $E$ and $F$ be graded $\RR$-modules of dimension
$d$.
Then
\[ \tn(E\otimes_\RR F)= \deg(E)\cdot \tn(F)+\deg(F)\cdot \tn(E).\]
\end{corollary}

\begin{proof} By Corollary~\ref{corsup}, we may assume that $E$ and $F$ are
torsion free modules.
 Let $\mathbb{P}$ and $\mathbb{Q}$ be minimal projective
resolutions of $E$ and $F$ respectively. The complex
$\mathbb{P}\otimes_\RR\mathbb{Q}$ is acyclic in codimension $1$, by the
assumption  on $E$ and $F$. We can then use Corollary~\ref{tneuler},
\[ \tn(E\otimes_\RR F) = \sum_{k\geq 0} (-1)^k \tn(\oplus_{i+j=k}
\mathbb{P}_i\otimes_\RR\mathbb{Q}_j).
\] Expanding gives the desired formula.
\end{proof}

Very similar to Proposition~\ref{e1adj} is:

\begin{proposition}\label{e1adj2}
Let $E\subset F$ be graded  torsion free $\RR$-modules of the
same multiplicity. If $E$ is reflexive, then
\[ \tn(E)= \tn(F) + \deg(F/E).\]
\end{proposition}

\subsubsection*{Simplicial complexes} Let $\SS=k[x_1, \ldots, x_n]$, for a
field $k$, and $I$ a monomial ideal. It should be a straightforward
calculation to determine $\tn(\SS/I)$. Here is a special case:

\begin{Theorem} Let $\Delta$ be a simplicial complex on the vertex
set $V=\{x_1, \ldots, x_n\}$, and denote by $k[\Delta]$ the
corresponding Stanley-Reisner ring. If $\dim k[\Delta]=d$, then
\[ \tn(k[\Delta])= df_{d-1}- f_{d-2} + f'_{d-2},\]
where $f_i$ denotes the number of faces of dimension $i$, and
$f'_{d-2}$ denotes the number of maximal faces of dimension $d-2$.
\end{Theorem}

\begin{proof} Set $k[\Delta]= \SS/I_{\Delta}$, and consider the decomposition
 $I_{\Delta}=
I_1\cap I_2$, where $I_1$ is the intersection of the primary
components of dimension $d$
and $I_2$ of the remaining components. The exact sequence
\[ 0 \rar I_1/I_{\Delta} \lar \SS/I_{\Delta}\lar \SS/I_1 \rar 0\]
gives, according to Proposition~\ref{e1adj},
\[ \tn(k[\Delta]) =\rme_1(k[\Delta]) + \hat{e}_0(I_1/I_{\Delta}).\]
From the Hilbert function of $k[\Delta]$ (\cite[Lemma 5.1.8]{BH}), we
have  that $\rme_1= df_{d-1}-f_{d-2}$, while if $I_1/I_{\Delta}$ is a
module of dimension $d-1$, its multiplicity is the number of maximal
faces of dimension $d-2$. 
\end{proof}

\begin{Theorem} \label{initheorem} Let $\SS=k[x_1, \ldots, x_n]$
 be a ring of polynomials
and $\AA = \SS/I$ a graded algebra. For a monomial ordering $>$, denote
by $I' = \textrm{in}_{>}(I)$ the initial ideal associated with $I$ and
set $\BB = \SS/I'$. Then $\tn(\BB)\geq \tn(\AA)$.
\end{Theorem}

\begin{proof} Let $J$ be the component of $I$ of maximal dimension and
consider the exact sequence
\[ 0 \rar J/I \lar \SS/I \lar \SS/J \rar 0.\] Now $\dim J/I < \dim \AA$, and
therefore $\tn(\AA) = \tn(\SS/J)=\rme_1(\SS/J)$. Denote by $J'$ the corresponding
initial ideal of $J$, and consider the sequence
\[ 0 \rar J'/I' \lar \SS/I' \lar \SS/J' \rar 0.\]
Noting that $\SS/I$ and $\SS/J$ have the same multiplicity, as do
$\SS/I'$ and $\SS/J'$  by Macaulay's theorem, we find that $\dim J'/I' < \dim
\AA$. This means that
\[ \tn(\SS/I') = \tn(\SS/J') =\rme_1(\SS/J') + \hat{e}_0(J'/I')=\rme_1(\SS/J) +
\hat{e}_0(J'/I')= \tn(\AA) + \hat{e}_0(J'/I').
\]
\end{proof}

\begin{example}{\rm Set $\AA=k[x,y,z,w]/(x^3-yzw, x^2y-zw^2)$. The
Hilbert series of this (Cohen-Macaulay) algebra is
\[ H_\AA(t)= \frac{h_\AA(t)}{(1-t)^2} = \frac{(1+t+t^2)^2}{(1-t)^2},\]
so that
\[ \tn(\AA) =\rme_1(\AA) = h'_\AA(1) = 18.\]

Consider now the algebra $\BB = k[x,y,z,w]/J$, where $J$ is the initial
ideal of $I$ for the Deglex order. A calculation with {\em Macaulay 2} gives
\[ J = (x^2y, x^3, xzw^2, xy^3zw, y^5zw).\] By Macaulay's Theorem,
$\BB$ has the same Hilbert function as $\AA$. An examination of the
components of  $\BB$ gives an exact sequence
\[ 0 \rar \BB_0 \lar \BB \lar \BB' \rar 0,\] where $\BB_0$ is the ideal of
elements with support in codimension $1$. By Corollary~\ref{corsup},
\[ \tn(\BB) = \tn(\BB')=\rme_1(\BB').\] At  same time one has the following equality of
$h$-polynomials,
\[ h_\BB(t) = h_{\BB'}(t) + (1-t)h_{\BB_0}(t),\]
and therefore
\[\rme_1(\BB') =\rme_1(\BB) + e_0(\BB_0).\] A final calculation of
multiplicities gives
$\rme_0(\BB_0)=5$, and
\[ \tn(\BB)= 18+ 5 = 23.\]
The example shows that $\tn(\AA)$ is independent of the Hilbert
function of the algebra.
}\end{example}

\subsection{Bounding tracking numbers}

We now describe how the technique of generic hyperplane sections
leads to bounds of various kinds. We are going to assume that the
algebras are defined over infinite fields.

One of the important properties of the tracking number is that it
will not change
under hyperplane sections as long as the dimension of the ring is at least
3. So one can answer questions about the tracking number just by
studying the 2-dimensional case. The idea here is that tracking number is more or less the
same material as $\rme_1$ and hence cutting by a superficial element
will not change
it unless the dimension  drops below 2.

\begin{proposition} \label{hypersec} Let $E$ be a finitely generated
graded module of dimension $d$ over $\RR=k[x_1,\dots,x_d]$ with $d>2$.
  Then for a  general  element $h$ of degree one,   $\RR'=\RR/(h)$ is
  also a polynomial ring, and
  $\tn_\RR(E)=\tn_{\RR'}(E')$, where $E'=E/hE$.
\end{proposition}

\begin{proof} First we shall prove the statement for a torsion free module
$E$. Consider the exact sequence
\[0 \rar E \lar E^{**} \lar C \rar 0.\] Note that
$C$ has codimension at least $2$ since after localization at any height $1$
prime $E$ and $E^{**}$ are equal. Now for a linear form $h$ in $\RR$ that
is a   superficial element for $C$
we can tensor the above exact sequence with $\RR/(h)$ to get the
complex
 \[ \Tor_1(C,\RR/(h)) \lar
E/hE \lar  E^{**}/hE^{**} \lar C/hC \rar  0.\] Now as an $\RR$-module,
 $C/hC$ has
codimension at least $3$, so as an $\RR'=\RR/(h)$ module it has codimension at
least $2$. Also as $\Tor_1(C,\RR/(h))$ has codimension at least $2$ as
an  $\RR$-module,  it is a torsion $\RR/(h)$-module. Hence we have
$\tn_{\RR'}(E/hE)=\tn_{\RR'}(E^{**}/hE^{**})$. But $E^{**}$ is a torsion free
$\RR/(h)$-module, so \[ \tn_{\RR'}(E/hE)=\rme_1(E^{**}/hE^{**})
=e_1(E^{**})=\tn_\RR(E^{**}) = \tn_\RR(E).\]

To prove the statement for a general $\RR$-module $E$, we consider the short
exact sequence
\[0 \rar E_0 \lar E \lar  E' \rar 0, \] where $\rme_0$ is the torsion
submodule of $E$. But ${E'}$ is torsion free, so by the first case we know that
$\tn_{\RR/(h)}({E'}/h{E'})=\tn_\RR({E'})$ for a general linear element
$h$ of $\RR$.

Now if in addition we restrict
ourselves to those $h$ that are superficial for $E$ and $\rme_0$, we can tensor the
above exact sequence with $\RR/(h)$ and get
\[0=\Tor_1({E'},\RR/(h)) \lar  E_0/hE_0
\lar E/hE \lar {E'}/h{E'} \rar  0,\] but since $\rme_0/hE_0$
 is a torsion $\RR'$-module, where $\RR'=\RR/(h)$, we have  $\tn_{\RR'}(E/hE) =
\tn_{\RR'}({E'}/h{E'})=\tn_\RR({E'})=\tn_\RR(E)$. 
\end{proof}

We shall now derive the first of our general bounds for $\tn(E)$ in
terms of the Castelnuovo-Mumford regularity $\reg(E)$ of the
module. For terminology and basic properties of the $\reg(\cdot)$
function, we shall use
\cite[Section 20.5]{Eisenbudbook}.

\begin{Theorem}\label{tnvsreg} Suppose that  $\RR=k[x_1, \ldots, x_d]$
and
that $E$ a finitely
 generated graded $\RR$-module of dimension $d$. Then
\[ \tn(E)\leq \deg(E)\cdot \reg(E).\]
\end{Theorem}

\begin{proof} The assertion is clear if $d=0$. For $d\geq 1$, if $\rme_0$
denotes the submodule of $E$ consisting  of the elements with finite
support, then
$\deg(E)=\deg(E/E_0)$, $\tn(E)=\tn(E/E_0)$ and $\reg(E/E_0)\leq
\reg(E)$, the latter according to \cite[Corollary
20.19(d)]{Eisenbudbook}. From this reduction, the assertion is also
clear if $d=1$.

If $d\geq 3$, we use a hyperplane section $h$ so that
$\tn_\RR(E)=\tn_{\RR/(h)}(E/hE)$ according to
Proposition~\ref{hypersec}, and $\reg(E/hE)\leq \reg(E)$ according
to \cite[Proposition 20.20]{Eisenbudbook}. (Of course,
$\deg(E)=\deg(E/hE)$.)

With these  reductions, we may assume that $d=2$ and that $\depth
E>0$. Denote by $\rme_0$  the torsion submodule of $E$ and consider
the exact sequence
\[ 0 \rar E_0 \lar E \lar E' \rar 0.\]
Noting that either $\rme_0$ is zero or $\depth E_0>0$,
on taking local cohomology with
respect to the maximal ideal $\mathfrak{m}=(x_1,x_2)$, we have the
exact sequence
\[0 \rar \H^1_{\mathfrak{m}}(E_0) \rar  \H^1_{\mathfrak{m}}(E) \rar
\H^1_{\mathfrak{m}}(E') \rar  \H^2_{\mathfrak{m}}(E_0) =0
\rar  \H^2_{\mathfrak{m}}(E) \rar \H^2_{\mathfrak{m}}(E') \rar 0,
\] from which we get
\[ \reg(E)=\max\{\reg(E_0), \reg(E')\}.\]
This provides the final reduction to $d=2$ and $E$ torsion free. Let
\[ 0 \rar \bigoplus_{j=1}^s \RR[-b_j] \lar \bigoplus_{i=1}^r \RR[-a_i]
\lar E \lar 0 \]
be a minimal projective resolution of $E$. From
Corollary~\ref{tneuler}, we have
\[ \tn(E)= \sum_{i=1}^ra_i-\sum_{j=1}^s b_j.\]
Reducing this complex modulo a hyperplane section $h$, we get a
minimal free resolution for the graded module $E/hE$ over the PID
$\RR/(h)$. By the basic theorem for modules over such rings, after
basis change we may assume that
\[ b_j=a_j+c_j,\quad  c_j>0,\quad j=1\ldots s.\]
Noting that $\alpha=\reg(E) = \max\{a_i, b_j-1\mid i=1\ldots r, j=1\ldots s
\}\geq 0$ and $\deg(E)=r-s$, we have
\begin{eqnarray*}
 \deg(E)\reg(E)-\tn(E)& =& (r-s)\alpha  - \sum_{i=1}^ra_i+
\sum_{j=1}^s(a_j+c_j) \\
&=& \sum_{i=s+1}^{r}(\alpha-a_i) + \sum_{j=1}^s c_j \\
&\geq & 0,
\end{eqnarray*}
as desired. 
\end{proof}

\subsection{Positivity of tracking numbers}\index{tracking number!positivity}

We shall
now prove our main result, a somewhat surprising positivity result
  for a
reduced homogeneous algebra $\AA$.
Since such algebras
already admit a general upper bound for $\tn(\AA)$ in terms of its
multiplicity,  together these statements are
useful in the  construction of integral closures by all algorithms
that use intermediate extensions that satisfy  condition $S_2$.

\begin{Theorem}\label{plustheorem} Let $\AA$ be a reduced non-negatively
graded
algebra that is finite over a standard  graded Noether normalization
$\RR$.
Then $\tn(\AA)\geq 0$. Moreover, if $\AA$ is an integral domain and $k$ is
algebraically closed, then $\tn(\AA)\geq \deg(\AA)-1$.
\end{Theorem}

\begin{proof} Let $\AA=\SS/I$, $\SS=k[x_1, \ldots, x_n]$, be a graded presentation of $\AA$.
From our earlier discussion, we may assume that $I$ is height-unmixed
(as otherwise the lower dimensional components gives rise to the
torsion part of $\AA$, which is dropped in the calculation of $\tn(E)$
anyway).

Let $I=P_1\cap \cdots \cap P_r$ be the primary decomposition of $I$,
and define the natural exact sequence
\[ 0 \rar  \SS/I \lar \SS/P_1 \times \cdots \times \SS/P_r \lar C \rar 0,\]
from which a calculation with Hilbert coefficients gives
\[ \tn(\AA) = \sum_{i=1}^r \tn(\SS/P_i) + \hat{\rme}_0(C).\]
This shows that it suffices to assume that $\AA$ is a domain.

Let $\overline{\AA}$ denote the integral closure of $\AA$. Note that
$\overline{\AA}$ is also a non-negatively graded  algebra and that the same
Noether normalization $\RR$ can be used. Since $\tn(\AA)\geq
\tn(\overline{\AA})$, we may assume that $\AA$ is integrally closed.

Since the cases $\dim \AA\leq 1$ are trivial, we may assume $\dim \AA=d\geq
2$. The case $d=2$ is also clear since $\AA$ is then Cohen-Macaulay.
Assume then that $d> 2$. We are going to change the base field using
rational extensions of the form $k(t)$, which do not affect the
integral closure condition. (Of course we may assume
 that  the base field
is infinite.)

If  $h_1$ and $h_2$ are linearly independent hyperplane
sections in $\RR$ they define a regular sequence in $\AA$, since the algebra
is normal and therefore satisfies condition $S_2$. Effecting a
change of ring of the type $k\rar k(t)$ gives a hyperplane section
$h_1-t\cdot h_2\in \RR(t)$, which is a prime element in $\AA$, according to
Nagata's trick (\cite[Lemma 14.1]{Fossum}).
Clearly we can choose $h_1$ and $h_2$ so that $h_1-t\cdot h_2$ is a generic
hyperplane section, for the purpose of applying
Proposition~\ref{hypersec} to $\AA$. This completes the reduction to
domains in dimension $d-1$.

The last assertion follows in the reduction to the case where  $\dim
\RR$ is $2$, as
$\tn(\AA)\geq \tn(\overline{\AA})$ and $\deg(\AA)=\deg(\overline{\AA})$.
The algebra $\overline{\AA}$ is Cohen--Macaulay so that its
$h$--polynomial $h(t)$ has only nonnegative coefficients and it
follows easily that $h'(1)\geq h(1)-1$, since $h(0)=1$ as $k$ is
algebraically closed.
\end{proof}

We single out:

\begin{corollary} Let $\AA$ be  a semistandard graded domain over
an algebraically closed field. If $\AA$ has $S_2$ and $\tn(\AA)=
\deg(\AA)-1$, then $\AA$ is normal.
\end{corollary}

\begin{example}{\rm
Let $\AA$ be a ring of polynomials over an algebraically closed field,
$\dim \AA\geq 4$, and let $I$ be a codimension two prime ideal
generated by $3$ quadrics. If $I$ is Cohen-Macaulay then $A/I$ is
normal. Indeed, the resolution of $\AA/I$ is given by a Hilbert-Burch
complex
\[ 0 \rar \AA[-3]^2 \lar \AA[-2]^3\lar \AA \lar \AA/I \rar 0.\]

A straightforward computation of the first Hilbert coefficients (see
\cite[Proposition 4.1.9]{BH}), shows that $\deg(\AA/I)=3$, and $\tn(\AA/I)=2$.
}\end{example}

One application is to the study of  constructions of the integral closure
of an affine domain in arbitrary characteristics.

\begin{Theorem} Let $\AA$ be a semistandard graded domain over a field
$k$ and let $\overline{\AA}$ be its integral closure. Then any chain
\[ \AA \subset \AA_1 \subset \cdots \subset \AA_n = \overline{\AA}\]
of distinct subalgebras satisfying Serre's condition $S_2$
has length at most ${{\rme}\choose{2}}$, where $\rme=\deg(\AA)$. If $k$ is
algebraically closed, such a  chain  has length at most
${{\rme-1}\choose{2}}$.
\end{Theorem}

\begin{proof} It will suffice, according to Corollary~\ref{boundlength}, to
show that $0\leq \tn(\AA)\leq {{\rme}\choose{2}}$.
The non-negativity having been  established in Theorem~\ref{plustheorem},
we now prove the upper bound.

If $k$ is a field of characteristic zero, then  by the theorem of the
primitive element $\AA$ contains a hypersurface ring $\SS=\RR[t]/(f(t))$,
where $\RR$ is a ring of polynomials $\RR=k[z_1, \ldots, z_d]$,
$\deg(z_i)=1$, and $f(t)$ is a homogeneous polynomial in $t$ of
degree $\rme$. As $\tn(\AA)\leq \tn(\SS)={{\rme}\choose{2}}$, the assertion
 holds in this case. (We also observe that when $\dim A\geq 3$,
any hyperplane section, say $h$, used to reduce the dimension that
were employed in the proof of Theorem~\ref{plustheorem}, could be
chosen so that the image of $\SS$ in $\overline{\AA}/h\overline{\AA}$ would
be $\SS/(h)$, and therefore we would maintain the same upper bound.)

\medskip

To complete the proof in
 other characteristics we resort to the following construction
 developed in \cite[Proposition 6.4]{icbook}.
If $\AA = \RR[y_1, \ldots, y_n]$, let
\[ E = \sum \RR y_1^{j_1}\cdots y_n^{j_n}, \quad 0\leq j_i <
r_i=[F_i:F_{i-1}] \]
be the  $\RR$-module of rank $e$ constructed there.
As the rank satisfies the equality  $\rme= \prod_{i=1}^n r_i$,
 $e$ is the
number of `monomials' $y_1^{j_1}\cdots y_{n}^{j_n}$. Their linear
independence over $\RR$ is a simple verification. Note also that there
are monomials of all degrees between $0$ and $(r_1-1, \ldots, r_n-1)$.
Thus according to Proposition~\ref{tne1},
 the bound for
$\tn(E)$ is obvious, with equality holding only when $E$ is a
hypersurface ring over $\RR$. (The precise value for $\tn(E)$ could be
derived from Corollary~\ref{tntimes}.)

\medskip

 For the case of  an algebraically closed field $k$ the last assertion
follows from Theorem~\ref{plustheorem}: $\tn(\AA)\geq e-1$, so that the
length of the chains is at most $\tn(E)-\tn(\overline{\AA})\leq
{{\rme}\choose{2}}-(\rme-1)= {{\rme-1}\choose{2}} $.
\end{proof}

\begin{remark}{\rm When the operations used to create the chain of
divisorial extensions are those described by
\cite[Proposition 6.66]{icbook}, the last bound  further reduces to
\[\displaystyle  \left \lceil \frac{(\rme-1)(\rme-2)}{4} \right \rceil .\]

There is another way in which the length of the divisorial chains may
shorten. Let $\AA \subset \BB\subset \overline{\AA}$ be divisorial
extensions of a graded algebra of dimension $d$ as above and let
$\RR=k[z_1, \ldots, z_d]$ be a Noether normalization. The canonical
module of $\BB$ is the graded module $\omega_\BB=\Hom_\RR(\BB, \RR[-d])$. We
recall that $B$ is {\em quasi-Gorenstein}\index{quasi-Gorenstein
algebra} if $\omega_\BB\simeq \BB[-a]$. In terms of their tracking
numbers this means that
\begin{eqnarray*} \tn(\omega_\BB)& = & -\tn(\BB) + d\cdot \deg(\BB) \\
& = & \tn(\BB) + a\cdot \deg(\BB),
\end{eqnarray*}
and therefore
\[ \tn(\BB) = \deg(\BB) \cdot \textrm{\rm half-integer}.\] A consequence
is that if more than one quasi-Gorenstein extension occur in a same
divisorial chain, they must lie fairly far apart.
}\end{remark}

\subsection{Algebras with Noether normalization}
\label{noether}

We now treat a very different venue of applications of
Theorem~\ref{jdegintclos}.
Let $\RR$ be a Noetherian domain and let $\BB=\RR[x_1, \ldots, x_n]$ be a
finitely generated $\RR$--algebra. We will assume that $\BB$ is
$\RR$-torsion free.
By abuse of terminology, by  {\em
arithmetic Noether normalizations} we shall mean statements asserting
the existence of extensions
\[ \RR \subset \SS= \RR[y_1, \ldots, y_r] \subset \BB, \]
$\BB$ finite over $\SS$, where $r$ depends uniformly on the Krull
dimensions of $\BB$ and $\RR$, or on some of the fibers. We shall refer
to $\SS$
 as a  Noether
normalization of $\BB$.

 Part of the significance of the existence of $\SS$ lies with the fact
 that it gives rise to a fibration
  $\Spec(\BB) \rar \Spec(\SS)$, that may add to the understanding
 of the geometry of $\BB$. The  `simpler' $\SS$ would bring greater
 benefit. Another motivation lies in that $\SS$ may serve as a
 platform for certain constructions on $\BB$, such as of its integral
 closure.

What are expected values for $r$? The answers tend to come in two
groups, depending on $\BB$ being a standard graded algebra or not, or
whether it is an integral domain or not. (The values for $r$ varies by
one accordingly.) For instance, if  $\BB = \RR  + \BB_{+}$, where $\RR$ is
an integral domain of field of fractions $K$,
by a
standard dimension formula (see \cite[Lemma 1.1.1]{alt}),
\[ \dim \BB = \dim \RR + \height Q= \dim \RR + \dim K\otimes_\RR \BB, \]
$\dim K\otimes_\RR \BB= K[x_1, \ldots, x_n]=n$ is the lower bound for $r$.

In one of his first papers, Shimura (\cite{Sh54}) proves such a result
when $\RR$ is a ring of algebraic integers.

\begin{Theorem}[\cite{Sh54}] Let $\RR$ be a ring of algebraic
integers of field of fractions $K$ and
let $B$ be a homogeneous $\RR$--algebra that is a domain. Then there is
a graded subalgebra $\SS= \RR[y_1, \ldots, y_r]$, $r= \dim
K\otimes_\RR \BB$,
over which $\BB$ is integral.
\end{Theorem}

Note that $\SS$ is actually a ring of polynomials over $\RR$, since $\dim
K\otimes_\RR \BB= \dim \BB-\dim \RR$. Shimura also remarked that this bound is not
achieved if  $\RR = \mathbb{C}[t]$.

\subsubsection*{Tracking numbers and $\jdeg$}

Let $\RR$ be a normal domain of field of fractions $K$, and
let $\AA$ be a semistandard graded algebra finite over the polynomial
graded subring $\SS=\RR[y_1, \ldots, y_r]$.
One can define the {\em determinant} of $\AA$ as an
$\SS$-module
\[ \det{_{\SS}}(\AA)=  (\wedge^r \AA)^{**},\]
where $r$ is the rank of $\AA$ as an $\SS$-module.

In case $\RR$ is a field, the graded module $\det{_{\SS}}(\AA)$ isomorphic
to a unique module $\SS[-n]$, so the integer $n$ can be employed as a
marker for $\AA$. It was used in \cite{DV} to define the {\em tracking
number} of $\AA$, $\tn_\RR(\AA)=n$. If $\RR$ is not a field, the integer $n$ is
still well-defined (as $\tn_K(K\otimes_\RR\AA)$), but it will not suffice
to fix $\det{_{\SS}}(\AA)$ entirely, particularly when given two such
$\SS$-algebras
$\AA\subset \BB$ we want to compare $\det{_{\SS}}(\AA)$ and
$\det{_{\SS}}(\BB)$.

To address this issue, we proceed as follows.
To the  rank $1$
reflexive module $\det{_{\SS}}(\AA)$,
 one can
attach several homogeneous divisorial ideals of $\SS$, all in the same divisor
class group. Suppose $I\subset \SS$ is one of these, and
\[ I = (\bigcap \mathfrak{p}_i^{(r_i)}) \cap (\bigcap
\mathfrak{q_j}^{(s_j)}),\]
is its primary decomposition, where we denote by $\mathfrak{p}_i$ the
primes that are extended from $\RR$, and $\mathfrak{q}_i$ those that
are not. This means  that  $\mathfrak{q}_j\cap \RR= (0)$ and thus
$\mathfrak{q}_jK\SS= (f_j)K\SS$, $\deg(f_j)> 0$.
 We associate a {\em
degree} to $I$ by setting
\[ \deg(I) = \sum_i r_i + \sum_j s_j \cdot  \deg(f_j).\]

\begin{definition} {\rm Let $\RR $ and $ \SS$ as above and let $\AA
\subset \BB
\subset \overline{\AA}$ be finitely generated $\RR$-algebras. Set
${\det}_\SS(\BB)= (\wedge^r \BB)^{**}$ and consider the natural image of
$(\wedge^r \AA)^{**}$ in ${\det}_\SS(\BB)$. Let
\[ I = \ann({\det}_\SS(\BB)/{\det}_\SS(\AA))
={\det}_\SS(A):_\SS{\det}_\SS(\BB).
\]
 The {\em relative tracking number} of $\AA$ in $\BB$ is the integer
\[ \tn(\AA,\BB) = \deg(I).\]
}
\end{definition}

\begin{example}{\rm
If
$ \AA = \mathbb{Z}[x,y,z]/(z^3+xz^2+x^2y)$, $\SS= \mathbb{Z}[x,y]$,
then $\overline{\AA}= \AA[z^2/x]$. It follows that
$ \tn(\AA,\overline{A})=1$.
}\end{example}

It is clear that $\tn(\AA,\BB)$ is well-defined. We note that in
\cite{DV}, in case $\RR$ is a field, one defines an
absolute
tracking number simply by declaring $\tn(\AA)$ as the degree of the
free rank one $\SS$-module $(\wedge^r \AA)^{**}$. It turns out that
$\tn(\AA,\BB)=
\tn(\AA)-\tn(\BB)$.

\begin{Theorem} \label{tnAB} Let $\RR$ be a normal domain and let $\AA$ be a
semistandard graded $\RR$-algebra with finite integral closure
$\overline{\AA} $ that admits a Noether normalization
$\SS$. If $\BB$ is a graded subalgebra, $\AA\subset \BB \subset
\overline{\AA}$ and $\AA$ satisfies the $S_2$ condition of Serre,
then
\[ \jdeg(\BB/\AA) = \tn(\AA,\BB).\]
\end{Theorem}

\begin{proof} We first clarify $\Ass_\RR(\BB/\AA)$. If $\mathfrak{p}$ is one of
these primes and $\mathfrak{p}\neq (0)$, the ideal
 $P=\mathfrak{p}\SS$ has the same height as  $\mathfrak{p}$ and
 therefore cannot consist of zerodivisors of $\BB/\AA$ if $\height  \
 \mathfrak{p}\geq 2$ since $\BB/\AA$ has the condition
$S_1$ of Serre (as an $\AA$--module or as an $\SS$-module). Thus
 $\mathfrak{p}$ has height $1$.

We now argue the asserted equality. Let $K$ be the field of fractions
of $\RR$.
If $\mathfrak{p}= (0)$, we note
that $\jdeg(K\BB/K\AA)$ is just the difference $\rme_1(K\AA)-\rme_1(K\BB)$ of
Hilbert coefficients of the graded modules $K\AA$ and $K\BB$. According
to \cite[Proposition 3.1]{DV}, for modules such as $K\AA$ and $K\BB$,
these values are the tracking numbers $\tn(K\AA)$ and $\tn(K\BB)$,
respectively. The difference $\tn(\AA)-\tn(\BB)$ clearly accounts for the
second summand in the definition of $\deg(I)$, for $I=
\ann({\det}_\SS(\BB)/{\det}_\SS(\AA))$.

Now we must account for the $j$-multiplicity of $\BB/\AA$ for a prime
$\mathfrak{p}$ of height $1$. Localizing at $\mathfrak{p}$, we may
assume that $\RR$ is a discrete valuation domain
and $\mathfrak{p}$ is its maximal ideal. Let $H$ be the submodule of
$\BB/\AA$ of support in $\mathfrak{p}$.
We claim that the multiplicity of $H$ is the corresponding integer
$r_i$  for $\mathfrak{p}_i= \mathfrak{p}$ in the expression of
$\deg(I)$. Consider the exact sequences of
$\SS$-modules
\[ 0 \rar \AA \lar \BB \lar \BB/\AA \rar 0,\]
\[ 0 \rar H \lar \BB/\AA \lar C \rar 0 ,\]
and the
prime ideal $P = \mathfrak{p}\SS$. Note that $P$ is the only associated
prime ideal of $H$ as an $\SS$-module and no associated prime of $C$ is
contained in it, and in  particular  $(\BB/\AA)_{P}= H_P$.
Applying the additivity formula for multiplicities to $H$ and
noting that $\deg(\SS/P)= 1$, we have $\deg(H) = \lambda(H_P)$.
Since $\SS_P$ is a discrete
valuation domain  the length of the torsion module $(\BB/\AA)_P$
is given  by $r_i$.

We have therefore matched each summand that occurs in $\deg(I)$ to
one  in definition of $\jdeg(\BB/\AA)$. 
\end{proof}

\begin{remark}{\rm It is not difficult to see that for the context
here, the assertion of
 Theorem~\ref{jdegintclos} could be improved by saying that the
 number of extensions between $\AA$ and $\overline{\AA}$ is bounded by
\[ \sum_i r_i + \sum_j s_j,\] rather than the whole of $\deg(I)$. (A
similar observation is given in \cite{ni1}.)
}\end{remark}

\subsection{Abstract tracking numbers}\index{tracking number!abstract}

If $\RR$ is an integrally closed  local ring, it is unclear  how to
construct
tracking numbers for its $\RR$-modules. While the determinant of an
$\RR$--module $E$ (here $\dim E = \dim \RR$) can be formed,
\[ \det(E) = (\wedge^rE)^{**},\]
there does not seem to be a  natural way to attach a degree to it.

In the case of an $\RR$--algebra $\AA$ of finite integral closure
$\overline{\AA}$, there is an {\it ad hoc} solution for the set of
submodules of $\overline{\AA}$ of rank $r=\rank_\RR(\AA)$.
The construction proceeds as follows.
Let $F$ be a normalizing free $\RR$--module,
\[ F = \RR^r = \RR e_1 \oplus \cdots \oplus \RR e_r \subset \AA
\subset \overline{\AA}.\]
There exists $f$ with  $0\neq f \in \RR$  such that $f\cdot
\overline{\AA} \subset   F$. For the $\RR$--submodule $E$ of
$\overline{\AA}$,

\[ \det(fE) \subset \RR(e_1\wedge \cdots \wedge e_r)= \RR\epsilon.\]
This means that
\[ \det(fE) = I.\RR\epsilon,\]
where $I$ is a divisorial ideal of $\RR$, so $If^{-r}$ is also a
divisorial ideal with a primary decomposition
\[ If^{-r} = \bigcap \p_i^{(r_i)}.\]

\begin{definition}{\rm  The tracking number of $E$ (offset by  $F$)
is the integer
\[\tn(E) = \sum_i r_i\cdot  \deg(\RR/\p_i) + \tn(\RR\epsilon).\]
}\end{definition}

The value
$\tn(E)$ is defined up to an offset but it is independent of $f$. It
will have several of the properties of the tracking number defined for
graded modules and can play the same role in the comparison of the
lengths of chains of subalgebras lying between $\AA$ and $\overline{\AA}$.
When $\RR$ is more general one can define $\tn(E)$ as the supremum of
its local values.

%\subsection{The determinant of an algebra: reductions}

\medskip

\subsubsection*{Questions} In addition to the issues already
raised---particularly  with regard on how to determine the tracking number of
an $\RR$-algebra without a Noether normalization---some of a more general nature are:
\begin{enumerate}
\item[{\rm (a)}]  What are the properties of the $k[x_1, \ldots, x_n]$--modules with positive tracking numbers?
For example, from Corollary~\ref{corsup} the set of such modules is
closed under tensor products.

\item[{\rm (b)}] If $E$ is a  module, $\tn(E^*)= -\tn(E)$, that has the
following consequence. Let
  $E$
be a torsion free $\RR$-module. From Proposition~\ref{Adual}, the
natural homomorphism
\[ E^*\otimes_\RR E \rar \Hom_\RR(E,E) \]
is an isomorphism in codimension $1$, so that
\[ \tn(\Hom_\RR(E,E)) = \tn(E^*\otimes_\RR E) =
\deg(E)\tn(E^*)+\deg(E^*)\tn(E) = 0.\]
\end{enumerate}
\section{Complexity of Normalization}

\subsection{Introduction}

 A general theory of the integral closure of algebraic
 structures--as a full fledged integration
of algorithms, cost estimates, along
 with tight
 predictions at termination--is highly desirable but it is an as yet
 unrealized goal.
Several
 stretches of
 this road
have been built and we will begin to describe the role of
multiplicities in the undertaking.

\medskip

We will now describe the setting in which our discussion takes place.
Let $\RR$ be a geometric/arithmetic  integral domain, and
let $\AA$ be a graded algebra (more generally with a
$\bbz^r$-grading)
\[ \AA= \RR[x_1, \ldots, x_n], \quad \mbox{\rm  $x_i$ homogeneous},\]
of integral closure
\[ \bar{\AA}= \RR[y_1, \ldots, y_m], \quad \mbox{\rm $y_i$
homogeneous}.\]

The construction of $\bar{\AA}$ usually starts out from an
observation of Emmy Noether in her last paper \cite{Noether}: If
$d\in \AA$ is a regular element of $\AA$ contained in the Jacobian
ideal, then
\[ \AA \subseteq  \bar{\AA} \subseteq \frac{1}{d} \AA.\]

One approach to the  construction
 of $\bar{\AA}$ (\cite{deJong}) consists in finding processes $\mathcal{P}$ that
create extensions
\[ \AA \mapsto \mathcal{P}(\AA) \subset \bar{\AA},\]
\[ \mathcal{P}^n(\AA) = \bar{\AA}, \quad n\gg 0.\]
typically $\mathcal{P}(\AA)= \Hom_{\AA}(I(\AA), I(\AA))$, where
$I(\AA)$ is
some ideal connected to the conductor/Jacobian of $\AA$.
Another approach (\cite{SiSw}), in characteristic $p$, is distinct.
Starting with $\frac{1}{d}\AA$, and using the Frobenius map, it
constructs a descending sequence
\[ V_0=\frac{1}{d} \AA \supset V_1 \supset \cdots \supset V_m=
 \bar{\AA}.
\]

 Predicting properties of $\bar{\AA}$ in $\AA$, possibly by
adding up the changes  $	\mathcal{P}^i(\AA) \rar
\mathcal{P}^{i+1}(\AA)$ (or the similar terms  $V_i$ in the
descending sequences)
 is a natural question in need of
clarification.

\medskip To realize this objective, we select the following specific
goals:
\begin{enumerate}
\item[{\rm (a)}]  We
develop numerical indices for
$\bar{\AA}$, that is for $\AA_i=\mathcal{P}^{i}$,
find $r$ such that \[\bar{\AA}_{n+r}= \AA_n \cdot \bar{\AA}_r, \quad n\geq
0.\]

\item[{\rm (b)}]
How many ``steps'' are there between $\AA$ and
$\bar{\AA}$,
\[\AA= \AA_0 \subset \AA_1 \subset \cdots \subset \AA_{s-1}\subset
\AA_s =
\bar{\AA},\]
where the $\AA_i$ are constructed by an effective
process?

\item[{\rm (c)}]
Express $r$ and $s$ in
terms of invariants of $\AA$.

\item[{\rm (d)}]
Generators of $\bar{\AA}$: Number  and
distribution of their degrees.

\end{enumerate}

Let $\RR$ be a commutative Noetherian ring, and let $\AA$ be a
finitely generated graded algebra (of integral closure $\bar{\AA}$).
   We want to gauge the `distance' from $\AA$ to $\bar{\AA}$.
   To this end, we introduce two degree functions whose properties
   help to determine the  position    of the algebras in the chains.

\begin{enumerate}

\item[{\rm (a)}]  (tracking number) $\tn: \AA\rar \bbn$,    for a certain class of algebras
(includes $\RR$ a field).
We shall use for distance
\[ \tn(\AA)-\tn(\bar{\AA})\]

\item[{\rm (b)}] $\jdeg:$ a multiplicity on all finitely generated graded
$\AA$-modules.    We shall use for distance
\[ \jdeg(\bar{\AA}/\AA)\]

\end{enumerate}

The first of these will be a {\em determinant}.
Consider a graded $\RR=k[z_1, \ldots,
z_r]$--module $E$:    Suppose $E$ has torsion free rank $r$.
    The {\em determinant} of $E$ is the graded module
\[ {\det}_{\RR}(E)=  (\wedge^r E)^{**}.\]

\[ {\det}_{\RR}(E) \simeq \RR[-c(E)],    \quad \tn(E):= c(E)\]
The integer $c(E)$ is the {\bf tracking number} of $E$.    If $E$
has no associated primes of codimension $1$, $\tn(E)$ can be read off
the Hilbert series $H_E(t)$ of $E$:
\[ H_E(t) = \frac{\hh(t)}{(1-t)^d},\quad  \tn(E) = \hh'(1)\]
   This shows that $\tn(E)$ is independent of $\RR$.
It comes with the following property.

\begin{Proposition}
Let $\RR=k[z_1, \ldots, z_d]$ and $E\subset F$ two graded
$\RR$--modules of the same rank. If $E$ and $F$ have the condition
$S_2$ of Serre,
\[ \tn(E)\geq \tn(F),\]
with equality if and only if $E=F$.
Therefore, if
\[ E_0 \subset E_1 \subset \cdots \subset E_n\] be a sequence of
distinct
f.g. graded $\RR$-modules of the same rank.
Then
\[ n \leq \tn(E_0)-\tn(E_n).\]
\end{Proposition}

In Section~\ref{trackingnumber} we will treat this function in a
broader setting.

\medskip

The other function, $\jdeg$, has different properties.
 One of the properties $\jdeg$ addresses  our main issue.
 Let $\RR$ be a Noetherian domain and $\AA$ a semistandard
graded $\RR$-algebra.
Consider a sequence of integral  graded extensions
\[\AA \subseteq \AA_0 \rar \AA_1 \rar \AA_2 \rar \cdots \rar
\AA_n=\overline{\AA},\]
where the $\AA_i$ satisfy the $S_2$ condition of Serre.
Then $n \leq \jdeg(\overline{\AA}/\AA)$
(\cite[Theorem 2]{jdeg1}).

The main goal
becomes to express $\jdeg(\overline{\AA}/\AA)$ in terms of {\it priori}
data, that is, one wants to view this number as an
  invariant of $\AA$. Relating it to other, more accessible,
  invariants of $\AA$ would bring
 considerable predictive power to the function $\jdeg$. This is the
 point of view of \cite{ni1} and \cite{HUV}, when $\jdeg$ is the
 Hilbert-Samuel multiplicity or the Buchsbaum-Rim multiplicity, respectively.

\subsection{Normalization and $\jdeg$}
\label{normjdeg}
In this subsection we make use of the  $\jdeg$ of appropriate modules in
order to
give estimates for the complexity of  normalization processes.

Let $\RR$ be a universally catenary ring and $\AA$ a standard graded integral
domain over $\RR$. We examine the character of the chains of graded
subalgebras between $\AA$ and its integral closure $\overline{\AA}$.
We will show that the lengths of certain chains are bounded by
$\jdeg(\overline{\AA}/\AA)$.

\medskip

To prepare the way we first make the following observation.

\begin{proposition} Let $\AA$ be an equidimensional
 Noetherian local ring of dimension
$d>0$ and $M$
a finitely generated module of dimension $d$. Assume that $N$ is a
submodule of $M$ such that $\dim (M/N) =d-1$. If $M$ satisfies the
condition $S_1$ of Serre and $N$ the condition $S_2$, then $M/N$
satisfies $S_1$ and in particular it is equidimensional.
\end{proposition}

\begin{proof} Let $\mathfrak{p}$ be an associated prime of $M/N$.
It is enough to prove that $\mathfrak{p}$ has codimension $1$.
 Applying
$\Hom_\AA(\AA/\mathfrak{p}, \cdot)$ to the exact sequence
\[ 0 \rar N \lar M \lar M/N \rar 0,\]
gives the exact sequence
\[ \Hom_\AA(\AA/\mathfrak{p}, M) \lar  \Hom_\AA(\AA/\mathfrak{p}, M/N)
\lar  \Ext^1_\AA(\AA/\mathfrak{p}, N).
\]
If  $\height \mathfrak{p}>1$, the modules at the ends vanish, and
 $\Hom_\AA(\AA/\mathfrak{p}, M/N)$ would be $0$. 
 \end{proof}

\begin{proposition} \label{quots2} Let $\AA\subset \BB \subset \BB'\subset
\overline{\AA}$ be finitely
generated graded   $\RR$-algebras all satisfying the $S_2$ condition of
Serre.
Then \[\jdeg (\BB'/\AA)=\jdeg(\BB/\AA)+\jdeg(\BB'/\BB).\]
\end{proposition}

\begin{proof} Because $\BB$, $\BB'$ satisfy the condition $S_2$ and have the same
total ring of fractions
\[\dim \BB/\AA=\dim \BB'/\AA=\dim \BB'/\BB=\dim \AA-1.\]

Noting that the quotients are equidimensional,
for every $\mathfrak{p}\in \Spec \RR$,
consider the short exact sequence
\[0 \rar (\BB/\AA)_\mathfrak{p} \rar (\BB'/\AA)_\mathfrak{p}
\rar(\BB'/\BB)_\mathfrak{p} \rar
0.\]
By Proposition~\ref{jmultses},
 \[j_\mathfrak{p}(\BB'/\AA)=j_\mathfrak{p}(\BB/\AA)+j_\mathfrak{p}(\BB'/\BB)\] for
  every ${ \mathfrak{p}}$, therefore
\[\jdeg(\BB'/\AA)=\jdeg(\BB/\AA)+\jdeg(\BB'/\BB).\] 
\end{proof}

\begin{Theorem} \label{jdegintclos} Let $\RR$ be a Noetherian
domain and $\AA$ a semistandard
graded $\RR$-algebra of finite integral closure $\overline{\AA}$.
Consider a sequence of distinct integral  graded extensions
\[\AA \subseteq \AA_0  \rar \AA_1 \rar \AA_2 \rar \cdots \rar
\AA_n=\overline{\AA},\]
where the $\AA_i$ satisfy the $S_2$ condition of Serre.
Then $n \leq \jdeg(\overline{\AA}/\AA)$.
\end{Theorem}

\begin{proof}
By Proposition~\ref{quots2}
and  Remark~\ref{nonzerojdeg},  we have
\[\jdeg(\overline{\AA}/\AA)>
\jdeg(\AA_{n-1}/\AA)>\cdots>\jdeg(\AA_1/\AA_0)>0.\]
Therefore, $n \leq \jdeg(\overline{\AA}/\AA)$.
\end{proof}

To illustrate  the values of  $\jdeg(\overline{\AA}/\AA)$ may take we
consider two classes of applications. First,  let $\AA$ be a
semistandard
 graded domain over a field
$k$ and let $\overline{\AA}$ be its integral closure.
If $\AA$ satisfies the condition $S_2$ and $\AA\neq \overline{\AA}$, the
graded module $
\overline{\AA}/\AA$ has dimension $\dim \AA-1$. Its Hilbert polynomial
satisfies
\[ H_{\overline{\AA}/\AA}(n) = H_{\overline{\AA}}(n)-H_\AA(n).\]
Therefore
\[
\jdeg(\overline{\AA}/\AA)=
\deg(\overline{\AA}/\AA) =\rme_1(\AA) -\rme_1(\overline{\AA}),
\] where $\rme_1(\cdot)$ is the next to the leading coefficient of the
Hilbert polynomial of the graded algebra. These coefficients are
positive and
\[\rme_1(\AA) -\rme_1(\overline{\AA})\leq {{\rme_0(\AA)}\choose{2}},\]
according to \cite{DV}, \cite[Theorem 6.90]{icbook}.

Let now
 $(\RR, \mathfrak{m})$
be  an analytically unramified  Cohen--Macaulay local ring of
dimension $d$ and $I$
an $\mathfrak{m}$--primary ideal. If $\AA= \RR[It]$,
\[\jdeg(\overline{\AA}/\AA)=
\deg(\overline{\AA}/\AA) = \rme_1(\overline{\AA})-\rme_1(\AA),\]
the switch in signs explained by the character of the Hilbert-Samuel
function of $\AA$.
A focus is to bound $ \rme_1(\overline{\AA})$, the main control of
normalization algorithms for $\RR[It]$. Several bounds are treated in
\cite{ni1}, in particular
\[  \rme_1(\overline{\AA}) \leq {\frac{d-1}{2}}\rme_0(\AA),\]
if $\RR$ is a regular local ring.

\bigskip

As we are going to see, these cases can be vastly generalized, to
arbitrary Rees algebras of ideals or modules, and to algebras
which are finite over a polynomial subring.

\subsection{Divisorial generation} Several of the modules that occur
in the construction of normalizations have the property $S_2$ of
Serre. It suggests the following notion of {\em generation}.

\begin{definition}{\rm Let $\RR$ be a Noetherian ring and $M$ a
finitely generated $\RR$-module. $M$ is {\em divisorially
generated}\index{divisorial generation} by
$r$ elements if there exists  a submodule $E\subset M$ generated by $r$ elements
such that $M$ is the  $S_2$-closure of $E$.
}\end{definition}

\begin{Theorem} \label{divgen} Let $\RR$ be a Noetherian
domain and $\AA$ a semistandard
graded $\RR$-algebra of finite integral closure $\overline{\AA}$.
Then $\bar{\AA}$ is the $S_2$-closure of a subalgebra generated by
homogeneous elements
$\AA[z_1, \ldots, z_r]$, for $r\leq 1+\jdeg(\overline{\AA}/\AA)$.
\end{Theorem}

\begin{proof} This is a direct consequence of Theorem~\ref{jdegintclos}. Let
$\AA_0$ be the $S_2$-closure of $\AA$.
If
$\bar{\AA}\neq \AA_0$, choose a homogeneous element $z_1\in \bar{\AA}\setminus
\AA_0$. In the same manner, if needed choose $z_2,\ldots$ and form the chain
\[ \AA_0 \subset \AA_1 \subset \cdots \subset \AA_{r}\subset \bar{\AA}.\]
Let $r= \jdeg(\overline{\AA}/\AA)$ and set
$\BB=\AA[z_1, \ldots, z_r]$. The $S_2$-closure of $\BB$ contains all $\AA_i$. 
\end{proof}

In several cases, $\bar{\AA}$ may actually have higher depth while
$\AA_0$ is Cohen-Macaulay. For example, if $(\RR,\mathfrak{m})$ is a two-dimensional
local normal domain and $I$ is $\mathfrak{m}$-primary with a minimal
reduction $J=(a,b)$, then in the sequence
\[ 0 \rar \AA= \RR[Jt] \lar \bar{\AA} \lar C \rar 0,\]
the module $C$ is Cohen-Macaulay, which can be used to bound its
number of generators. Let us phrase one of these situations in the
context of Corollary~\ref{dete1}. A different and more comprehensive
approach that emphasizes exclusively the multiplicity can be found in
\cite[Section 3]{cxintclos}.

\begin{Theorem} \label{divdete1} Suppose that $\RR$ is a normal standard graded
algebra $($or a normal Noetherian local domain$)$
and that $E\subset F$ are homogeneous modules of the same rank. If
$E$ is Cohen-Macaulay and $F$ is a  reflexive module with $\depth
F\geq \dim R-1$, then
\[ \nu(F)\leq \rme_0(F)+\rme_1(E)-\rme_1(F).\]
\end{Theorem}

\begin{proof} The module $F/E$ is Cohen-Macaulay of multiplicity
$\rme_1(E)-\rme_1(F)$, which will bound $\nu(F/E)$. 
\end{proof}

\begin{remark}{\rm
In the case of ideal above, $\deg(C) =
\bar{\rme}_1(J)-\rme_1(J)=\bar{\rme}_1(J)$, so that $\bar{\RR[Jt]}$ will be
generated over $\RR[Jt]$  by
\[ 1 + \bar{\rme}_1(J)\]
elements.

\medskip

If $I$ is an ideal of codimension one, what is the bound? Suppose $I$ is
unmixed. Then there is a minimal reduction (may assume the residue
field of $\RR$ is infinite) $J=(a,b)$. It is easy to see that $\RR[Jt]$
is Cohen-Macaulay. The module $C=\bar{\RR[Jt]}/\RR[Jt]$ is Cohen-Macaulay
and since at every prime ideal $\mathfrak{p}$ of $\RR$ of codimension
one $\RR[Jt]$ is normal, the support of $C$ is $\{\mathfrak{m}\}$. Thus
$\jdeg(C)=\deg(C)$, therefore $\bar{\RR[Jt]}$ is generated by
\[ 1+ \jdeg(\bar{\RR[Jt]}/\RR[Jt])\]
elements over $\RR[Jt]$.
}\end{remark}

For contrast, wed briefly describe some of the related results of
\cite[Section 3]{cxintclos}. Note that their emphasis in exclusively
on the multiplicity.
Let $\AA = k[x_1, \ldots, x_n]/I$ be a reduced equidimensional affine
algebra over a field $k$ of characteristic zero, $n= \emb(\AA)$, of
dimension $\dim \AA=d$. We call
$\ecodim(\AA)=n-d $ the {\em embedding codimension} of $\AA$.

\begin{lemma}[{\cite[Lemma 6.95]{icbook}}] \label{deglemma}
Let $k$ be a field of characteristic zero,  $\SS$  a reduced and
equidimensional  standard graded $k$-algebra of dimension
$d$ and multiplicity $e$, and assume that $\emb(\SS)\leq d+1$.
Let $\SS \subset \BB$ be a finite and birational extension of graded
rings.
\begin{itemize}
\item[{\rm (a)}] The $\SS$-module $\BB/\SS$ satisfies
\[ \deg(\BB/\SS)\leq (\rme-1)^2.\]
If $\rme\geq 3$ then
\[ \deg(\BB/\SS)\leq (\rme-1)^2-1. \]
If $\rme\geq 4$ or if $\BB$ is $S_2$ but not Cohen-Macaulay, then
\[ \deg(\BB/\SS)\leq (\rme-1)^2-2. \]

\item[{\rm (b)}]
If $\BB$ is Cohen-Macaulay, then
\[ \deg(\BB/\SS ) \leq {{\rme}\choose{2}}.\]
If moreover $k$ is algebraically closed and $\SS$ is a
domain, then
\[ \deg(\BB/\SS)\leq {{\rme-1}\choose {2}} .\]
\end{itemize}
\end{lemma}

\begin{Theorem}[{\cite[Theorem 6.96]{icbook}}] \label{main1}
Let $k$ be a field of characteristic zero and  $\AA$  a reduced
and equidimensional standard graded $k$--algebra of dimension $d$ and
multiplicity $e$. If $\AA\subset \BB$ is a finite and birational
extension of graded rings
with $\depth_\AA \BB \geq d-1 $
then \[\nu_\AA(\BB)\leq (\rme-1)^2 +1\]
and
\[ \emb(\BB) \leq  (\rme-1)^2+d +1.\]
\end{Theorem}

\begin{remark}{\rm  \cite[Theorem 6.93(a)]{icbook} and
Lemma~\ref{deglemma}  show that the estimates
of Theorem~\ref{main1} can be sharpened under suitable additional
assumptions.  Indeed, if
$\rme\geq 3$ then
\[\nu_\AA(\BB)\leq (\rme-1)^2  \ \ \   {\rm and} \ \ \   \emb(\BB) \leq (\rme-1)^2+d. \]
If $\rme\geq 4$ or if $\BB$ satisfies $S_2$ and $\rme\geq 3$, then
\[\nu_\AA(\BB)\leq (\rme-1)^2-1  \ \ \  {\rm and}  \ \ \    \emb(\BB) \leq (\rme-1)^2+d-1. \]
}\end{remark}

\begin{corollary} \label{corfor3}
Let $k$ be a field of characteristic zero and  $\AA$  a reduced
and equidimensional standard graded $k$--algebra of dimension $3$ and
multiplicity $\rme\geq 3$. The integral closure $\BB=\overline{\AA}$
satisfies
the inequality
\[ \emb(\BB) \leq (\rme-1)^2+2.\]
\end{corollary}

\begin{remark}\label{3.5}{\rm  If in the setting of Theorem~\ref{main1} and its proof,
$\TT$ is a homogeneous Noether normalization of $\AA$, then
 $\nu_\TT(\BB)\leq \rme(\rme-1)^2+\rme$,
as can be seen by applying the theorem with $\AA=\SS$.
This bound is not strictly module-theoretic: the
algebra structure of $\BB$ really matters. The assertion fails
if $B$ is merely a finite $\TT$-module, even a graded reflexive $T$-module,
with $\depth_{\TT}\BB\geq d-1$.
For example, let $\TT$ be a polynomial ring in $d$ variables
over an infinite field and $I$ a homogeneous perfect ideal of $\TT$
of height 2 that is generically a complete
intersection but has a large number of generators
(such ideals exist whenever $d\geq 3$).
There is an exact sequence of the form
\[ 0 \rar \TT[-a] \lar E \lar I \rar 0 \,,\]
with $E$ a graded reflexive $\TT$-module. Now indeed $\depth \, E =d-1$ and
$\deg(E)=2$, whereas $\nu_\TT(E)\gg 0$.
%According to Theorem 3.3, $E$ could not be the underlying
%module of a reduced and equidimensional standard graded $k$-algebra
%The other point is that even the quadratic bound may be too large and
%perhaps a bound of type $C\deg(A)^{3/2}$ suffices.
%In examples it is all we have been able to achieve.
}\end{remark}

\section{Normalization of Ideals}

\subsection{Introduction}
Let $\RR$ be a Noetherian domain and $I$ an ideal. By the {\em
normalization} \index{normalization of an ideal} of $I$ we mean the
integral closure of the Rees algebra of $I$,
\[ \AA = \RR[It]= \sum_{n\geq 0}I^nt^n \subset \bar{\AA} = \sum_{n\geq
0} \bar{I^n}t^n. \]

The study of $\bar{\AA}$ has  features that open opportunities beyond
those examined in the previous sections. The ultimate goal is to
describe properties of $\bar{\AA}$ that can be predicted from $\AA$.
To that end, following \cite{ni2}, we will introduce several
multiplicity based measures.

\subsubsection*{$S_2$-ification} For the purpose of this section to
assume that $\RR$ is a Noetherian domain that is a homomorphic image of
a Gorenstein ring. For an ideal $I$, we will assume that $\BB=
\bar{\RR[It]}$ is a finitely generated module over $\AA= \RR[It]$. For
any graded subalgebra $\AA\subset \CC \subset \BB$.

\medskip

 We are going to
collect some properties of $\CC$ (\cite[Theorems 6.17,
6.18, 6.19]{compu}; see also \cite[Proposition 4.4]{icbook}).

\begin{Proposition} Let $\RR$ be a normal domain of dimension $d$,
 and let $\SS\rar \CC$ be
a surjection with $\SS$ a locally Gorenstein ring. If $\dim \SS=r$, then
\begin{enumerate}
\item[{\rm (a)}] $\Ext_\SS^{r-d-1}(\CC, \SS)$ is a canonical module $\omega_{\CC}$ of $\CC$.

\item[{\rm (b)}] The $S_2$-ification of $\CC$ is a graded subalgebra
$\AA\subset \DD \subset \BB$, where
\[\DD = \Hom_{\CC} (\omega_{\CC}, \omega_{\CC}) \simeq
\Ext_\SS^{r-d-1}(\Ext_\SS^{r-d-1}(\CC, \SS),\SS).\]
\end{enumerate}
\end{Proposition}

The other property permits the definition of an associated graded
ring to $\DD$ (\cite[Proposition 4.6]{icbook}).

\begin{Proposition} For any $\AA$-algebra $\DD= \sum_{n\geq 0}D_nt^n\subset
\BB$ satisfying the condition $S_2$ of Serre, the filtration $\{D_n
\}$ is decreasing, and the associated graded ring
\[ \gr_{\DD}(\RR) = \bigoplus_{n\geq 0} D_n/D_{n+1}\] satisfies the
condition $S_1$ of Serre.
\end{Proposition}

\subsection{Indices of normalization}

Suppose $\RR$ is a commutative ring and $J,I$ are ideals of $\RR$ with
$J\subset I$. $J$ is a {\em reduction} of $I$ if $I^{r+1}=JI^r$ for
some integer $r$; the least such integer is the {\em reduction
number} of $I$ relative to $J$. It is denoted
$\red_J(I)$. $I$ is {\em equimultiple} if there is a reduction $J$
generated by $\height I$ elements.

\begin{Definition}{\rm Let $\RR$ be a locally analytically unramified normal domain and
let $I$ be an ideal.
\begin{itemize}
\item[{\rm (a)}]
 The {\em  normalization
index}\index{normalization of an ideal!normalization index} of $I$ is the
smallest integer $s=s(I)$
such that
\[\overline{I^{n+1}} = I\cdot \overline{I^n} \quad  n\geq s.\]

\item[{\rm (b)}]
 The {\em generation
index}\index{normalization of an ideal!generation index} of $I$ is the
smallest integer $s_0=s_0(I)$ such that
\[ \sum_{n\geq 0}\overline{I^n}t^n=
\RR[\overline{I}t, \ldots,
\overline{I^{s_0}}t^{s_0}].\]

\item[{\rm (c)}] The {\em normal relation
 type}\index{normalization of an ideal!normal relation type} of $I$ is
 the  maximum degree of a minimal generating set of  the presentation ideal
\[ 0 \rar M \lar \RR[\TT_1, \ldots, \TT_m] \lar
\RR[\overline{I}t, \ldots,
\overline{I^{s_0}}t^{s_0}] \rar 0.\]
\end{itemize}
}\end{Definition}

\medskip

For example, if $\RR = k[x_1, \ldots, x_d]$ is a polynomial ring over
a field and $I = (x_1^d, \ldots, x_d^d)$, then $I_1 = \overline{I} =
(x_1, \ldots, x_d)^d$. It follows that $s_0(I) = 1$, while $s(I) =
\red_I(I_1) = d-1$.

\medskip

If $(\RR, \m)$ is a local ring, these indices have an expression in
term of the special fiber ring $F$ of the normalization map $\AA
\rar \overline{\AA}$.

\begin{Proposition} \label{indicesnorm} With the above assumptions let \[ F =
\overline{\AA}/(\mathfrak{m}, It)\overline{\AA}= \sum_{n\geq 0} F_n.\]
We have
\begin{eqnarray*}
s(I) & = & \sup\{ n \mid F_n \neq 0\}, \\
s_0(I) & = & \inf \{n \mid F = F_0[F_1, \ldots, F_n]\}.
\end{eqnarray*}
Furthermore, if the index of nilpotency of $F_i$ is $r_i$, then
\[ s(I) \leq \sum_{i=1}^{s_0(I)} (r_i-1).\]
\end{Proposition}

Although these integers areell defined---since $\overline{\AA}$ is
finite over $\AA$---it is not clear, even in case $\RR$ is a regular
local ring, which invariants of $\RR$ and of $I$ have a bearing on the
determination of $s(I)$. An affirmative case is that of a monomial
ideal $I$ of a ring of polynomials in $d$ indeterminates over a
field---when $s\leq d-1$.

\begin{Remark}{\rm For normal ideals the information about the Hilbert function is very detailed--see, for
instance a modern survey in \cite{MMV12}. It remains a goal to have for general normalizations that level
of specificity.}
\end{Remark}

\subsection{Cohen-Macaulay normalization}

Expectably, normalization indices are easier to obtain when the
normalization of the ideal is Cohen-Macaulay. The following is
directly derived from the known  characterizations of
Cohen-Macaulayness of Rees algebras of ideals  in terms of
associated graded rings and reduction numbers (\cite{AHT},
\cite{JK95}, \cite{SUV2}).

\begin{Theorem} \label{redintclos} Let $(\RR,\mathfrak{m})$ be a
Cohen--Macaulay local ring and let $\{ I_n, \ n\geq 0\}$ be a
decreasing multiplicative filtration of ideals, with $I_0=R$,
$I_1=I$, and the property that the corresponding Rees algebra
$\BB=\sum_{n\geq 0}I_nt^n$ is finite over $\AA$. Suppose that $\height
I\geq 1$
 and let $J$ be a minimal reduction of $I$. If $\BB$ is Cohen-Macaulay,
 then
\[{I_{n+1}}= J{I_n}=
I_1{I_n} \  \  \ \  \mbox{{\rm for every \,}} n \geq \ell(I_1)-1,\]
and in particular, $\BB$ is generated over $\RR[It]$ by forms  of
 degrees at most $\ell-1=\ell(I_1)-1$,
\[ \sum_{n\geq 0} {I_n}t^n = R[{I_1}t, \ldots,
{I_{\ell-1}}t^{\ell-1}].\]
\end{Theorem}

The proof of Theorem~\ref{redintclos} relies on substituting in any
of the proofs mentioned above (\cite[Theorem 5.1]{AHT}, \cite[Theorem
2.3]{JK95},
\cite[Theorem 3.5]{SUV2}) the $I$-adic filtration $\{I^n\}$ by the filtration
$\{I_n\}$.

\begin{Corollary} \label{normalpowers}
Let $(\RR, \mathfrak{m})$ be a local analytically unramified normal
Cohen--Macaulay ring and let $I$ be an ideal. If $\overline{\AA}$ is
Cohen-Macaulay then both indices of normalization $s(I)$ and
$s_0(I)$ are at most $\ell(I)-1$. Moreover, if ${I^n}$ is integrally
closed for $n< \ell(I)$, then $I$ is normal.
\end{Corollary}

A case this applies to is that of monomial ideals in a polynomial
ring, since
 the ring $\overline{\AA}$ is Cohen--Macaulay by Hochster's theorem
(\cite[Theorem 6.3.5]{BH}) (see  \cite{RRV}).

%\newpage

\begin{Example}\rm Let $I=I(\mathcal{C})=(x_1x_2x_5,x_1x_3x_4,x_2x_3x_6,x_4x_5x_6)$
be the edge ideal associated to the
clutter

\bigskip

\bigskip

%\bigskip

\bigskip

%\begin{small}
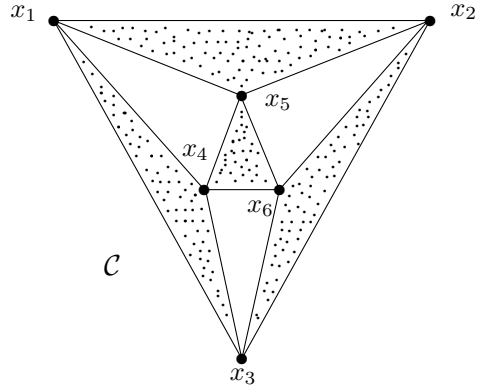
\special{em:linewidth 0.4pt} \unitlength 0.5mm \linethickness{0.4pt}
\begin{picture}(30,100)(-10,30)
\emline{30.00}{130.00}{1}{130.00}{130.00}{2}
\put(80,40){\circle*{2.5}} \put(30,130){\circle*{2.5}}
\put(70,85){\circle*{2.5}} \put(90,85){\circle*{2.5}}
\put(80,110){\circle*{2.5}} \put(130,130){\circle*{2.5}}
\emline{130.00}{130.00}{3}{80.00}{40.00}{4}
\emline{80.00}{40.00}{5}{30.00}{130.00}{6}
\emline{70.00}{85.00}{7}{90.00}{85.00}{8}
\emline{90.00}{85.00}{9}{79.67}{110.00}{10}
\emline{79.67}{110.00}{11}{70.33}{85.00}{12}
\emline{70.33}{85.00}{13}{80.00}{40.00}{14}
\emline{80.00}{40.00}{15}{90.00}{85.00}{16}
\emline{90.00}{85.00}{17}{130.00}{130.50}{18}
\emline{130.00}{130.00}{19}{79.50}{110.00}{20}
\emline{80.00}{110.00}{21}{30.00}{130.00}{22}
\emline{30.00}{130.00}{23}{70.00}{85.00}{24}
\put(67.33,124.00){\circle*{.1}} \put(73.67,124.33){\circle*{.1}}
\put(75.83,124.33){\circle*{.1}} \put(84.67,123.33){\circle*{.1}}
\put(91.33,125.00){\circle*{.1}} \put(90.00,127.67){\circle*{.1}}
\put(89.67,127.67){\circle*{.1}} \put(80.33,119.33){\circle*{.1}}
\put(80.33,119.33){\circle*{.1}} \put(92.33,121.33){\circle*{.1}}
\put(106.33,125.33){\circle*{.1}}
\put(106.33,125.33){\circle*{.1}} \put(99.67,124.67){\circle*{.1}}
\put(99.67,122.33){\circle*{.1}} \put(99.67,122.33){\circle*{.1}}
\put(102.00,127.67){\circle*{.1}}
\put(111.33,127.00){\circle*{.1}}
\put(101.33,127.33){\circle*{.1}} \put(85.67,116.00){\circle*{.1}}
\put(57.33,127.00){\circle*{.1}} \put(46.00,127.33){\circle*{.1}}
\put(50.33,125.00){\circle*{.1}} \put(56.33,122.00){\circle*{.1}}
\put(59.67,128.00){\circle*{.1}} \put(62.00,121.33){\circle*{.1}}
\put(69.00,118.00){\circle*{.1}} \put(78.00,114.67){\circle*{.1}}
\put(71.33,127.33){\circle*{.1}} \put(82.00,127.00){\circle*{.1}}
\put(84.67,119.33){\circle*{.1}} \put(74.00,116.33){\circle*{.1}}
\put(86.33,119.67){\circle*{.1}} \put(62.00,124.67){\circle*{.1}}
\put(116.33,127.33){\circle*{.1}}
\put(103.00,123.00){\circle*{.1}} \put(95.67,123.00){\circle*{.1}}
\put(94.67,127.67){\circle*{.1}} \put(95.33,120.00){\circle*{.1}}
\put(86.67,119.33){\circle*{.1}} \put(90.67,116.67){\circle*{.1}}
\put(81.33,114.33){\circle*{.1}} \put(82.67,126.67){\circle*{.1}}
\put(62.00,125.67){\circle*{.1}} \put(65.67,127.67){\circle*{.1}}
\put(51.00,128.00){\circle*{.1}} \put(40.00,128.33){\circle*{.1}}
\put(38.00,119.00){\circle*{.1}} \put(40.33,115.67){\circle*{.1}}
\put(42.67,112.33){\circle*{.1}} \put(42.67,110.67){\circle*{.1}}
\put(42.67,110.67){\circle*{.1}} \put(46.00,110.00){\circle*{.1}}
\put(43.67,108.67){\circle*{.1}} \put(48.67,105.33){\circle*{.1}}
\put(48.83,105.33){\circle*{.1}} \put(49.00,105.33){\circle*{.1}}
\put(46.00,105.00){\circle*{.1}} \put(47.33,101.67){\circle*{.1}}
\put(47.33,101.67){\circle*{.1}} \put(51.00,101.00){\circle*{.1}}
\put(50.33,97.33){\circle*{.1}} \put(53.00,99.00){\circle*{.1}}
\put(55.00,98.00){\circle*{.1}} \put(52.67,96.33){\circle*{.1}}
\put(52.00,94.00){\circle*{.1}} \put(56.00,94.00){\circle*{.1}}
\put(56.00,94.00){\circle*{.1}} \put(60.33,91.00){\circle*{.1}}
\put(57.67,89.33){\circle*{.1}} \put(58.33,87.00){\circle*{.1}}
\put(60.67,87.00){\circle*{.1}} \put(63.33,87.00){\circle*{.1}}
\put(66.67,86.67){\circle*{.1}} \put(67.67,83.67){\circle*{.1}}
\put(67.67,83.67){\circle*{.1}} \put(63.67,83.00){\circle*{.1}}
\put(59.33,82.33){\circle*{.1}} \put(54.67,90.33){\circle*{.1}}
\put(60.33,84.33){\circle*{.1}} \put(60.33,79.33){\circle*{.1}}
\put(62.67,79.67){\circle*{.1}} \put(65.00,79.67){\circle*{.1}}
\put(67.00,79.67){\circle*{.1}} \put(69.33,80.00){\circle*{.1}}
\put(69.67,76.67){\circle*{.1}} \put(67.00,76.67){\circle*{.1}}
\put(65.00,76.67){\circle*{.1}} \put(65.00,76.67){\circle*{.1}}
\put(62.67,76.67){\circle*{.1}} \put(65.00,72.67){\circle*{.1}}
\put(67.33,70.67){\circle*{.1}} \put(69.67,71.00){\circle*{.1}}
\put(69.00,68.67){\circle*{.1}} \put(69.00,66.67){\circle*{.1}}
\put(71.00,66.00){\circle*{.1}} \put(70.00,64.33){\circle*{.1}}
\put(72.33,63.00){\circle*{.1}} \put(71.67,61.00){\circle*{.1}}
\put(73.67,58.00){\circle*{.1}} \put(73.67,56.67){\circle*{.1}}
\put(75.00,54.33){\circle*{.1}} \put(76.00,51.67){\circle*{.1}}
\put(56.33,87.33){\circle*{.1}} \put(58.67,94.00){\circle*{.1}}
\put(84.33,51.00){\circle*{.1}} \put(84.00,51.67){\circle*{.1}}
\put(87.00,55.67){\circle*{.1}} \put(86.00,58.00){\circle*{.1}}
\put(88.33,57.67){\circle*{.1}} \put(86.67,60.33){\circle*{.1}}
\put(66.67,120.00){\circle*{.1}} \put(54.67,125.00){\circle*{.1}}
\put(106.67,127.67){\circle*{.1}}
\put(108.67,124.00){\circle*{.1}} \put(81.33,123.00){\circle*{.1}}
\put(87.33,125.67){\circle*{.1}} \put(89.00,122.00){\circle*{.1}}
\put(75.00,127.33){\circle*{.1}} \put(82.67,117.00){\circle*{.1}}
\put(92.67,83.33){\circle*{.1}} \put(92.67,79.33){\circle*{.1}}
\put(98.00,83.00){\circle*{.1}} \put(88.00,65.00){\circle*{.1}}
\put(90.00,64.67){\circle*{.1}} \put(92.33,65.67){\circle*{.1}}
\put(94.00,68.33){\circle*{.1}} \put(91.00,69.00){\circle*{.1}}
\put(89.00,70.33){\circle*{.1}} \put(92.00,71.00){\circle*{.1}}
\put(95.00,72.00){\circle*{.1}} \put(92.67,74.00){\circle*{.1}}
\put(90.67,74.33){\circle*{.1}} \put(90.67,76.33){\circle*{.1}}
\put(93.33,76.33){\circle*{.1}} \put(96.33,76.33){\circle*{.1}}
\put(98.33,76.33){\circle*{.1}} \put(99.33,78.33){\circle*{.1}}
\put(96.67,79.00){\circle*{.1}} \put(95.00,80.33){\circle*{.1}}
\put(91.67,81.00){\circle*{.1}} \put(96.00,84.00){\circle*{.1}}
\put(102.00,84.67){\circle*{.1}} \put(98.33,86.00){\circle*{.1}}
\put(95.33,86.67){\circle*{.1}} \put(98.67,88.33){\circle*{.1}}
\put(102.67,88.33){\circle*{.1}} \put(104.67,88.67){\circle*{.1}}
\put(106.00,92.33){\circle*{.1}} \put(102.33,93.67){\circle*{.1}}
\put(100.33,93.67){\circle*{.1}} \put(104.67,95.33){\circle*{.1}}
\put(107.33,96.33){\circle*{.1}} \put(111.00,97.33){\circle*{.1}}
\put(107.67,101.67){\circle*{.1}}
\put(110.67,102.00){\circle*{.1}}
\put(112.33,102.33){\circle*{.1}}
\put(113.67,104.33){\circle*{.1}}
\put(113.33,106.00){\circle*{.1}}
\put(111.33,106.33){\circle*{.1}}
\put(114.00,108.33){\circle*{.1}}
\put(116.33,109.67){\circle*{.1}}
\put(117.00,111.67){\circle*{.1}}
\put(119.00,113.00){\circle*{.1}}
\put(119.67,115.33){\circle*{.1}}
\put(121.67,117.33){\circle*{.1}}
\put(121.67,118.67){\circle*{.1}}
\put(123.67,121.00){\circle*{.1}}
\put(120.33,128.00){\circle*{.1}}
\put(112.67,126.00){\circle*{.1}} \put(86.00,127.33){\circle*{.1}}
\put(70.33,121.33){\circle*{.1}} \put(72.67,120.00){\circle*{.1}}
\put(77.67,121.00){\circle*{.1}} \put(68.67,126.67){\circle*{.1}}
\put(59.33,124.00){\circle*{.1}} \put(68.33,74.00){\circle*{.1}}
\put(62.67,90.00){\circle*{.1}} \put(73.33,87.00){\circle*{.1}}
\put(77.00,89.33){\circle*{.1}} \put(80.00,91.00){\circle*{.1}}
\put(80.33,89.00){\circle*{.1}} \put(84.00,89.00){\circle*{.1}}
\put(84.00,91.33){\circle*{.1}} \put(86.33,87.33){\circle*{.1}}
\put(82.33,95.33){\circle*{.1}} \put(80.33,95.33){\circle*{.1}}
\put(77.33,95.33){\circle*{.1}} \put(79.33,98.33){\circle*{.1}}
\put(80.67,98.67){\circle*{.1}} \put(80.67,100.33){\circle*{.1}}
\put(78.67,100.67){\circle*{.1}} \put(80.00,102.33){\circle*{.1}}
\put(80.00,104.33){\circle*{.1}} \put(80.00,105.67){\circle*{.1}}
\put(95.67,126.00){\circle*{.1}} \put(89.33,61.00){\circle*{.1}}
\put(89.00,67.33){\circle*{.1}} \put(90.33,72.00){\circle*{.1}}
\put(95.67,74.33){\circle*{.1}} \put(100.33,81.33){\circle*{.1}}
\put(93.67,86.00){\circle*{.1}} \put(97.33,90.67){\circle*{.1}}
\put(100.67,90.33){\circle*{.1}} \put(104.33,91.00){\circle*{.1}}
\put(110.00,99.67){\circle*{.1}} \put(105.67,99.00){\circle*{.1}}
\put(103.67,97.00){\circle*{.1}} \put(110.67,104.33){\circle*{.1}}
\put(66.33,67.67){\circle*{.1}} \put(69.67,62.00){\circle*{.1}}
\put(71.33,68.67){\circle*{.1}} \put(71.00,74.00){\circle*{.1}}
\put(65.00,85.00){\circle*{.1}} \put(57.67,84.33){\circle*{.1}}
\put(66.67,82.00){\circle*{.1}} \put(49.33,99.67){\circle*{.1}}
\put(50.67,103.33){\circle*{.1}} \put(46.00,107.67){\circle*{.1}}
\put(60.00,121.00){\circle*{.1}} \put(64.67,119.33){\circle*{.1}}
\put(72.33,115.67){\circle*{.1}} \put(76.33,114.33){\circle*{.1}}
\put(79.33,117.00){\circle*{.1}} \put(76.33,119.33){\circle*{.1}}
\put(74.33,121.33){\circle*{.1}} \put(70.67,124.33){\circle*{.1}}
\put(64.67,123.00){\circle*{.1}} \put(47.33,125.00){\circle*{.1}}
\put(48.33,128.00){\circle*{.1}} \put(43.33,127.33){\circle*{.1}}
\put(54.33,127.33){\circle*{.1}} \put(98.00,127.67){\circle*{.1}}
\put(90.00,119.33){\circle*{.1}} \put(89.00,116.33){\circle*{.1}}
\put(93.00,118.33){\circle*{.1}} \put(98.33,120.67){\circle*{.1}}
\put(106.00,122.33){\circle*{.1}}
\put(103.33,125.67){\circle*{.1}}
\put(108.67,125.67){\circle*{.1}} \put(78.33,87.33){\circle*{.1}}
\put(82.67,87.33){\circle*{.1}} \put(82.00,90.33){\circle*{.1}}
\put(84.33,94.67){\circle*{.1}} \put(81.33,93.00){\circle*{.1}}
\put(81.33,93.00){\circle*{.1}} \put(75.67,91.67){\circle*{.1}}
\put(74.33,89.33){\circle*{.1}} \put(77.67,97.67){\circle*{.1}}
\put(82.67,97.33){\circle*{.1}} \put(81.33,101.67){\circle*{.1}}
\put(79.33,127.33){\circle*{.1}} \put(79.00,124.33){\circle*{.1}}
\put(81.67,125.00){\circle*{.1}} \put(79.00,122.67){\circle*{.1}}
\put(83.33,120.33){\circle*{.1}} \put(53.33,123.33){\circle*{.1}}
\put(80.00,112.67){\circle*{.1}} \put(78.33,92.33){\circle*{.1}}
\put(74.67,90.33){\circle*{.1}} \put(77.67,98.00){\circle*{.1}}
\put(75.67,93.00){\circle*{.1}} \put(78.00,94.67){\circle*{.1}}
\put(75.67,94.33){\circle*{.1}}
\put(22.00,132.00){\makebox(0,0)[cc]{\small $x_1$}}
\put(139.00,132.33){\makebox(0,0)[cc]{\small $x_2$}}
\put(80.33,35.33){\makebox(0,0)[cc]{\small $x_3$}}
\put(89.67,108.33){\makebox(0,0)[cc]{\small $x_5$}}
\put(85.00,80.00){\makebox(0,0)[cc]{\small $x_6$}}
\put(67.67,95.00){\makebox(0,0)[cc]{\small $x_4$}}
\put(45.67,65.00){\makebox(0,0)[cc]{$\mathcal C$}}
%\put(195.67,125.00){\makebox(0,0)[cc]
%{\small $\RR[It]\subsetneq\overline{R[It]}=R_s(I)$}}
\end{picture}
%\end{small}

\bigskip

\medskip

Consider the incidence matrix $A$ of this clutter, i.e., the columns
of $A$ are the exponent vectors of the monomials of $I$. Since
the polyhedron $Q(A)=\{x\vert xA\geq 1;\, x\geq 0\}$ is integral, we
have the equality $\overline{\RR[It]}=\RR_s(I)$, the symbolic Rees
algebra of $I$ (see
\cite[Proposition~3.13]{reesclu}). The ideal $I$ is not normal
because the monomial $x_1x_2\cdots x_6$ is in
$\overline{I^2}\setminus I^2$.

The minimal primes of $I$ are:
$$
\begin{array}{llll}
\mathfrak{p}_1=(x_1,x_6),&\mathfrak{p}_2=(x_2,x_4),&\mathfrak{p}_3=(x_3,x_5),&
\\
\mathfrak{p}_4=(x_1,x_2,x_5),&\mathfrak{p}_5=(x_1,x_3,x_4),&
\mathfrak{p}_6=(x_2,x_3,x_6),&\mathfrak{p}_7=(x_4,x_5,x_6).
\end{array}
$$
For any $n$,
\[I^{(n)}=\bigcap_{i=1}^7\mathfrak{p}_i^n.\] A computation
with {\it Macaulay\/} $2$ (\cite{Macaulay2}) gives that
$\overline{I^2}=(I^2,x_1x_2\cdots x_6)$ and that
$\overline{I^3}=I\ \overline{I^2}$. By
Theorem~\ref{redintclos},
$I^{(n)}=II^{(n-1)}$ for $n\geq \ell(I)=4$, where $\ell(I)$ is the
analytic spread of $I$. As a consequence,
\[\bar{\RR[It]}= \RR[It, x_1x_2x_3x_4x_5x_6t^2].\]

\end{Example}

\subsection{Integrality and multiplicity}

Let $(\RR, \mathfrak{m})$ be a Noetherian local ring and let $J\subset
I$ be ideals. If $J$ is a reduction of $I$, say $I^{r+1}=JI^r$, for
some $r$, this relationship has several expressions in terms of
multiplicities. The simplest of these is when the ideals are
$\mathfrak{m}$-primaries. It follows from the definition of the
Hilbert-Samuel function that $\rme_0(J)=e_0(I)$. The converse holds if
$\RR$ is a universally catenary local ring, according to a
ground-braking
 theorem of Rees (\cite{Reesa}). Several extensions, to more general
 ideals have been proved since. The current state of affairs is that
 provided by \cite{FM} for ideals, and pushed to  fuller generality
 for
 modules in \cite{UV}.

\begin{Theorem}[{\cite[Theorem 3.3]{FM}}] \label{FlennerM}
Let $\RR$ be   an equidimensional complete
 Noetherian local
ring and let $J\subset I$ be $\RR$-ideals.
 Then the following are equivalent:
\begin{enumerate}
\item[{\rm (a)}]  $J$ is a reduction of $I$;

\item[{\rm (b)}] If $\AA$ and $\BB$ are the Rees algebras of $I$ and $J$, then
 $j_\mathfrak{p}(\gr(\AA))=j_\mathfrak{p}(\gr(\BB))$ for
every $\mathfrak{p}\in \Spec (\RR)$.
\end{enumerate}
\end{Theorem}

\subsection{Brian\c{c}on-Skoda theorems}

Let us recall the notion of the Brian\c{c}on-Skoda number of an
ideal. If
$I$ is an ideal of a Noetherian ring $\RR$, the {\it Brian\c{c}on-Skoda
number} $c(I)$ of $I$ is the smallest integer $c$ such that
$\overline{I^{n+c}}\subset J^n$ for every $n$ and every reduction
$J$ of $I$. The motivation for this definition is a result of
\cite{BrSkoda} asserting that for  the rings of convergent power series over
${\mathbb{C}}^n$ (later extended to all regular local rings)
$c(I)< \dim \RR$.
The existence of uniform values for $c(I)$ has been established for
much broader classes of rings (and modules) by Huneke (\cite{Hu00}).
Examination of their
proofs should yield effective bounds.
Here
we will derive a related invariant.

\bigskip

The following extends  \cite[Theorems 2.2 and 4.2]{ni1}, which dealt
with
equimultiple ideals. Let us recall:

\begin{definition}{\rm Let $(\RR, \m)$ be a Noetherian local ring of
dimension $d$.
We say that $\RR$ is unmixed,
if $\operatorname{dim} \widehat{\RR}/\p = d$ for every $\p \in
\operatorname{Ass} \widehat{\RR}$, where $\widehat{\RR}$ denotes the
$\m$-adic completion of $\RR$.\index{unmixed local ring} This is a
notion that can be extended to modules. \index{unmixed module over a
local ring} (Rings and modules with this property are
equidimensional.)
}\end{definition}

There are numerous questions in local algebra that permit changes of
rings of the form $\RR\rar \widehat{\RR}$, noteworthy being
consideration of Hilbert polynomials and their coefficients. In the
case of unmixed local rings the following formulation is useful.

\begin{Proposition} \label{unmixedrep} Let $(\RR, \m)$ be an unmixed
complete local ring. There exists a Gorenstein local ring $\SS$,
$\dim \SS=\dim \RR$, an embedding of $\SS$-modules
\[ \RR \hookrightarrow \SS^n.\]
\end{Proposition}

Let us state this as:

\begin{proposition} \label{embedinfree} If $\SS$ is a Gorenstein local ring and
$\RR=\SS/L$ is local ring such that $\dim \SS=\dim \RR$ and $\RR$ is
equi-dimensional and without embedded primes, then $\RR$ can be
embedded in a finitely generated free $\SS$-module.
\end{proposition}

\begin{proof} We are going to argue that $\RR$ embeds into its
$S_2$-ification $\AA$, which in turn embeds into some $\SS^n$:
\begin{itemize}
\item[{\rm (a)}] $\RR$ has a canonical module: $\Hom_{\SS}(\SS/L, \SS)$.

\item[{\rm (b)}] $\AA= \Hom_{\SS}(\Hom_{\SS}(\SS/L,\SS),
\Hom_{\SS}(\SS/L,\SS))= \Hom_{\SS}(\Hom_{\SS}(\SS/L,\SS), \SS).$

\item[{\rm (c)}] The natural mapping $\RR\rar \AA$ is an isomorphism at each
associated prime of $\RR$, therefore it is an embedding.

\item[{\rm (d)}] Let $\SS^n \lar \Hom_{\SS}(\SS/L,\SS)\rar 0$ be a surjection;
dualizing with $\Hom_{\SS}(\cdot,\SS)$ gives an embedding of $\AA$
into $\SS^n$.
\end{itemize}
Finally, the composition $\RR\rar \AA \rar \SS^n$ gives the assertion.
\end{proof}

\begin{Corollary}  Let $\SS$ is be a local ring and let  $\RR=S/L$
be a
ring
of dimension $\dim \SS$.
If $\RR$
 can embedded into
a MCM $\SS$-module then $\RR$ is unmixed.
\end{Corollary}

\begin{proof} By assumption $\widehat{\RR}$ embeds into a maximal
$\widehat{\SS}$-module, in particular $\widehat{\RR}$ is
equidimensional and has no embedded primes.
\end{proof}

\begin{Theorem} \label{bsjdegI} Let $\RR$ be a reduced, unmixed
ring and
let $I$ be an ideal of Brian\c{c}on-Skoda number $c(I)$. Then
\[ \jdeg(\overline{\RR[It]}/\RR[It]) \leq c(I)\cdot  \jdeg(\gr_I(\RR)).\]
\end{Theorem}

\begin{proof} Let $J$ be a minimal reduction of $I$.
We consider the $\RR[Jt]$-module $\CC=\overline {\RR[It]}/\RR[Jt]$.
Then $\jdeg(\CC) \geq \jdeg(\overline{\RR[It]}/\RR[It])$.
By definition of $c=c(I)$, $\CC$ is a submodule of the graded
$\RR[Jt]$-module
\[\DD = \bigoplus_{n\geq 0} J^{n-c}/J^n.\]
 The inclusions $J^{n-c}\supset J^{n-c+1} \supset
\cdots \supset J^n$
 induce a filtration of $\DD$.
Because local cohomology commutes with  direct sums,
from the Hilbert functions of the factors in this filtration it follows that
\[\jdeg(\DD)=c \cdot \jdeg(\gr_J(\RR))= c\cdot \jdeg(\gr_I(\RR)).\]
Moreover, $\jdeg(\CC) \leq \jdeg(\DD)$, so
$\jdeg(\overline{\RR[It]}/\RR[It])
\leq c(I)\cdot  \jdeg(\gr_I(\RR))$. 
\end{proof}

There are several extensions of this result
if we broaden the notion of Brian\c{c}on-Skoda number of an ideal.

\begin{Theorem} \label{newbs}
Let $k$ be a perfect field, let $(\RR, \mathfrak{m})$ be
  a reduced
 $k$-algebra essentially of finite type.  Suppose that $\RR$ has
 isolated singularities. If $L$ is the Jacobian ideal of $\RR$, then
for any ideal $I$ of analytic spread  $\ell = \dim R$,
\[ \jdeg(\overline{\RR[It]}/\RR[It])\leq (\ell+ \lambda(\RR/L) -1) \cdot
 \jdeg(\gr_I(\RR)).\]
\end{Theorem}

\begin{proof} Let $\DD$ be a $S_2$-ification of $\RR[It]$. From
Theorem~\ref{GBS},  setting
$\CC= \overline{\RR[It]}$ and  $c=\ell -1$, we have
$L \overline{I^{n+c}}\subset D_n$.
Consider  the diagram
\begin{small}
\[
\diagram
 & D_n/D_{n+c} & & C_{n+c}/LC_{n+c} \dto|>>\tip & \\
0 \rto & (LC_{n+c}+ D_{n+c})/D_{n+c} \rto \uto|<\hole|<<\ahook
& C_{n+c}/D_{n+c} \rto
 & C_{n+c}/(LC_{n+c} + D_{n+c}) \rto &0 .
\enddiagram
\] \end{small}
Considering  that the modules in the  short exact sequence have dimension
$\dim \RR$ and that the module on the right is supported in
$\mathfrak{m}$ alone, it follows that
\[ \jdeg(\CC/\DD)= \jdeg((L\CC+ \DD)/\DD) + \jdeg(\CC/L\CC+\DD) =
\jdeg((L\CC+ \DD)/\DD) + \deg(\CC/L\CC+\DD) .
\] Given the embedding on the left and the surjection on the right, we
have
that
\[ \jdeg((L\CC+\DD)/\DD)\leq \jdeg(\DD/\DD[-c])= c\cdot \jdeg(\gr(\DD)),  \]
and
\[ \deg(\CC/L\CC+\DD) \leq \deg(\CC/L\CC) \leq \lambda(\RR/L) \deg
(\CC/\mathfrak{m}\CC).
\] Since $\jdeg(\gr(\DD))= \jdeg(\gr_I(\RR))$ by Proposition~\ref{jdeggr},
 we obtain the  estimate of multiplicities
\[ \jdeg(\CC/\DD) \leq c\cdot \jdeg(\gr_I(\RR)) + \lambda(\RR/L)\cdot
f_0(\CC),
\] where $f_0(C)= \deg(\CC/\mathfrak{m}\CC)$.

Finally, we consider the exact sequence
\[ 0 \rar \mathfrak{m}C_n/C_{n+1}  \lar  C_n/C_{n+1} \lar
C_n/\mathfrak{m}C_n \rar 0.
 \]
Taking into account that
   $\gr(\CC)$ is  a ring  which has the condition $S_1$,
$ \bigoplus_n\mathfrak{m}C_n/C_{n+1}  $ either vanishes or has
 the same dimension as the
ring. Arguing as above, we have
\[ \jdeg(\gr(\CC))=
\jdeg(\bigoplus_n\mathfrak{m}C_n/C_{n+1}) + f_0(\CC).
\]
A final application of Proposition~\ref{jdeggr}, yields either
\[ \jdeg(\overline{\RR[It]}/\RR[It])\leq (\ell+ \lambda(\RR/L) -1) \cdot \jdeg(\gr_I(\RR))\]
as in the assertion of the theorem, or
\[ \jdeg(\overline{\RR[It]}/\RR[It])\leq (\ell+ \lambda(\RR/L) -2) \cdot
\jdeg(\gr_I(\RR)),\]
if $\overline{I}\neq \mathfrak{m}$. 
\end{proof}

The argument above are relatively indifferent to the ideal $I$,
suggesting the convenience of defining the {\em Brian\c{c}on-Skoda
number of the ring $\RR$} as the integer
\[ c(\RR)= \sup\{ c(I) \mid I \subset \RR.\}\]

\begin{Conjecture}{\em Let $\RR$ be a Noetherian ring of finite Krull
dimension. Then
\[ c(\RR[x]) \leq c(\RR) + 1.\]
}\end{Conjecture}

\subsection{Equimultiple ideals}
Let $(\RR, \mathfrak{m})$ be an unramified local domain and let $I$ be
an ideal. Suppose $I$ is {\em equimultiple}\index{equimultiple
ideal}, that is it has a reduction $J$ generated by $\height I$
elements.

\begin{Proposition} Suppose further that $\RR$ is a Cohen-Macaulay
ring, and let $\{\mathfrak{p}_1, \ldots, \mathfrak{p}_m\}$ be the
set of minimal primes of $I$.  Then
\begin{itemize}
\item[{\rm (a)}] $c(I)=\max\{ c(I_{\mathfrak{p}_i})\mid i=1,\ldots,
m\}$, in particular if $\RR$ is a regular local ring then $c(I)\leq
\height I-1$.

\item[{\rm (b)}] $\jdeg(\bar{\RR[It]}/\RR[Jt])= \sum_{i=1}^m
\rme_1(\bar{I_{\mathfrak{p}_i}})$.

\item[{\rm (c)}] ?

\end{itemize}

\end{Proposition}

\begin{Theorem}\label{maintheorem} Let $(\RR, \m)$ be a
reduced local Cohen--Macaulay ring of dimension $d
> 0$ and let $I$ be an $\mathfrak{m}$--primary ideal.
\begin{itemize}
\item[{\rm (a)}] If in addition $\RR$ is an algebra essentially of
finite type over a perfect field $k$ with type $t$, and $\delta
\in {\rm Jac}_k(\RR)$ is a non zerodivisor, then
 \[\overline{\rme}_1(I)
\leq \frac{t}{t + 1} \bigl[ (d-1)\rme_0(I) + \rme_0(I +\delta \RR/\delta
\RR) \bigr].\]

 \item[{\rm (b)}] If the assumptions of $($a$)$ hold,
then
\[
\overline{\rme}_1(I) \leq (d-1)\bigl[\rme_0(I)
-\lambda(\RR/\overline{I})\bigr] + \rme_0(I +\delta \RR/\delta \RR).
\]

\item[{\rm (c)}] If $\RR$ is analytically unramified and $\RR/\m$ is
infinite, then
\[
\overline{\rme}_1(I) \leq c(I) \ {\rm min} \, \{\frac{t}{t + 1} \ \rme_0(I),
\rme_0(I) -\lambda(\RR/\overline{I})\}.
\]
\end{itemize}
\end{Theorem}
\begin{proof} We may assume that $\RR/\m$ is infinite. Then, passing to a
minimal reduction we may suppose that $I$ is generated by a
regular sequence $f_1, \ldots, f_d$. Notice this can only decrease
$c(I)$. Let $\SS$ be a local ring obtained from $\RR$ by a purely
transcendental residue field extension and by factoring out $d-1$
generic elements $a_1, \ldots, a_{d-1}$ of $I$. To be more
precise, $\SS=\RR(\{X_{ij}\})/(a_1, \dots, a_{d-1})$ with $\{X_{ij}\}$
a set of $(d-1)d$ indeterminates and
$a_i=\sum_{j=1}^{d}X_{ij}f_j$. Notice that $\SS$ is also a
birational extension of a localization of a polynomial ring over
$\RR$, and hence is analytically unramified according to
\cite[Theorem 36.8]{Nagata} and \cite[1.6]{R}. Furthermore $\SS$ is a
one--dimensional local ring and the $\SS$--ideal $IS$ is generated
by a single non zerodivisor, say $I\SS=f\SS$. From \cite[Theorem
1]{Itoh} one has
\begin{equation}\label{eq2}
\overline{I}\SS = \overline{I\SS},
\end{equation}
\begin{equation}\label{eq3}
\overline{I^{n}}S = \overline{(IS)^{n}} \quad {\rm for \ } n \gg
0.
\end{equation}
The last fact combined with the genericity of $a_1, \ldots,
a_{d-1}$ yields $ \overline{\rme}_1(I) =\overline{\rme}_1(I\SS)$. Moreover
$\overline{\rme}_1(I\SS)= \lambda(\overline{\SS}/\SS)$ as $\SS$ is a
one--dimensional analytically unramified local ring. Thus
\begin{equation}\label{eq4}
\overline{\rme}_1(I) =\lambda(\overline{\SS}/\SS).
\end{equation}

In the setting of $($a$)$ and $($b$)$ the element $\delta$ is a
non zerodivisor on $\SS$. Furthermore Theorem~\ref{GBS} shows that
\begin{equation}\label{eq1}
\delta I^{d-1} \subset \delta \ \overline{I^{d-1}}\subset I^n
\colon  \overline{I^{n}}  \qquad {\rm for \ every \ } n.
\end{equation}
For $n \gg 0$, by $($\ref{eq3}$)$, $($\ref{eq1}$)$ and since $f^n
\SS$ is contained in the conductor $\SS \colon \overline{\SS}$, we
obtain
\[
\delta f^{d-1} \SS \subset f^n \SS \colon \overline{f^n \SS} = f^n \SS
\colon f^n \overline{ \SS}=\SS \colon \overline{\SS}.
\]
Hence
\begin{equation}\label{eq5}
\delta f^{d-1} \overline{\SS} \subset \SS \colon \overline{\SS}.
\end{equation}

We prove $($a$)$ by computing lengths along the inclusions
\begin{equation}\label{eq6}
\delta f^{d-1} \SS \subset \delta f^{d-1} \overline{\SS} \subset \SS
\colon \overline{\SS} \subset \SS.
\end{equation}
Also recall that
\begin{equation}\label{eq7}
\lambda(\overline{\SS}/\SS) \leq t \ \lambda(\SS/\SS \colon
\overline{\SS})
\end{equation}
by \cite[the proof of 3.6]{HK} (see also \cite[Theorem 1]{BrH}
 and \cite[2.1]{D}). We obtain
\begin{eqnarray*}
\frac{t + 1}{t} \  \overline{\rme}_1(I) & = & \lambda(\overline{\SS}/\SS)
+ \frac{1}{t} \ \lambda(\overline{\SS}/\SS) \qquad \ \ \ \ \qquad
\qquad \qquad  {\rm by}
\ (\ref{eq4})\\
& \leq &  \lambda(\delta f^{d-1}\overline{\SS}/\delta f^{d-1} \SS) +
\lambda(\SS/\SS \colon \overline{\SS}) \qquad \qquad {\rm by}
\ (\ref{eq7}) \\
& \leq  & \lambda(\SS/\delta f^{d-1} \SS)  \qquad \qquad \ \qquad
\qquad \qquad \qquad {\rm by} \
(\ref{eq6}) \\
& =  & (d-1) \ \lambda(\SS/f \SS) + \lambda(\SS/\delta \SS) \\
& =  & (d-1)\rme_0(I) + \rme_0(I +\delta \RR/\delta \RR)  \qquad \qquad {\rm
by \ the \ genericity \ of \ } a_1, \ldots, a_{d-1}.
\end{eqnarray*}

Next we prove part $($b$)$. The inclusion (\ref{eq5}) yields the
filtration
\[
%\begin{equation}
%\label{eq8}
\SS= \delta f^{d-1}\overline{\SS} + \SS \subset f^{d-1}\overline{\SS}
+ \SS
\subset \ldots \subset f^2\overline{\SS} + \SS \subset f \overline{\SS}
+ \SS \subset \overline{\SS}, \]
%\end{equation}
which shows
\begin{equation}\label{eq9}
 \overline{\rme}_1(I) =\lambda(\overline{\SS}/\SS)=  \sum_{i=1}^{d-1}
 \lambda(f^{i-1}\overline{\SS} + \SS/ f^i\overline{\SS} + \SS) +
\lambda(f^{d-1}\overline{\SS} + \SS/ \delta f^{d-1} \overline{\SS} +
\SS).
\end{equation}
Multiplication by $f$ induces epimorphisms of $\SS$--modules
\begin{equation}\label{eq10}
\frac{f^{i-1}\overline{\SS} + \SS}{ f^i\overline{\SS} + \SS}
\twoheadrightarrow \frac{f^{i}\overline{\SS} + \SS}{
f^{i+1}\overline{\SS} + \SS}\ .
\end{equation}
Now (\ref{eq9}) and (\ref{eq10}) show
\begin{equation}\label{eq11}
 \overline{\rme}_1(I) \leq (d-1) \
 \lambda(\overline{\SS}/f\overline{\SS} +\SS) +
\lambda(f^{d-1}\overline{\SS} + \SS/ \delta f^{d-1} \overline{\SS} +
\SS
).
\end{equation}
Next we claim that
\begin{equation}\label{eq12}
\lambda(\overline{\SS}/f\overline{\SS} +\SS) = \lambda(\overline{I}/I).
\end{equation}
Indeed,
\begin{eqnarray*}
\lambda(\overline{\SS}/f\overline{\SS} +\SS)  & = &
\lambda(\overline{\SS}/f\overline{\SS}) -
\lambda(f\overline{\SS} +\SS/f\overline{\SS}) \\
& =&  \lambda(\SS/f \SS) -  \lambda(\SS/\SS \cap f \overline{\SS}) \\
& =  & \lambda(\SS \cap f \overline{\SS}/f \SS)  \\
& =  & \lambda(\overline{f \SS}/f \SS)\\
& =  &  \lambda(\overline{I}/I) \qquad \qquad {\rm by} \
(\ref{eq2}) .
\end{eqnarray*}
On the other hand,
\begin{eqnarray*}
\lambda(f^{d-1}\overline{\SS} + \SS/ \delta f^{d-1} \overline{\SS} +
\SS )
&\leq&
\lambda(f^{d-1}\overline{\SS} / \delta f^{d-1} \overline{\SS}) \\
& =  &   \lambda(\overline{\SS} / \delta \overline{\SS}) \\
& =  &   \lambda (\SS / \delta \SS)\\
& =  &   \rme_0(I +\delta \RR/\delta \RR) \qquad \qquad {\rm by \ the \
genericity \ of \ } a_1, \ldots, a_{d-1}.
\end{eqnarray*}
Therefore
\begin{equation}\label{eq13}
\lambda(f^{d-1}\overline{\SS} + \SS/ \delta f^{d-1} \overline{\SS} +
\SS )
\leq
 \rme_0(I +\delta \RR/\delta \RR).
\end{equation}
 Combining (\ref{eq11}),  (\ref{eq12}) and (\ref{eq13}) we deduce
\begin{eqnarray*}
 \overline{\rme}_1(I) &\leq&  (d-1) \ \lambda(\overline{I}/I) + \rme_0(I
 +\delta \RR/\delta \RR)
 \\
& =  &   (d-1)\bigl[\rme_0(I) -\lambda(\RR/\overline{I})\bigr] + \rme_0(I
+\delta \RR/\delta \RR).
\end{eqnarray*}

Finally we prove part $($c$)$.  Write $b=c(I)$. We first claim
that
\begin{equation}\label{eq14}
 f^{b} \overline{\SS} \subset \SS \colon \overline{\SS}.
\end{equation}
Indeed, for $n \gg 0$
\begin{eqnarray*}
f^n \SS & \supset &
\overline{I^{n + b}}\SS   \\
& = & \overline{(I\SS)^{n + b}}  \qquad \qquad {\rm by} \
(\ref{eq3}) \\
& =  &  \overline{f^{n +b} \SS}  \\
& =  & f^{n +b} \overline{\SS}   \qquad \qquad {\rm since \ } n \gg
0.
\end{eqnarray*}
Therefore $f^{b} \overline{\SS} \subset \SS $, proving (\ref{eq14}).
Now (\ref{eq14}) yields the filtrations
\begin{equation}\label{eq15}
 f^{b} \SS \subset f^b  \overline{\SS} \subset \SS
\colon \overline{\SS} \subset \SS,
\end{equation}
\begin{equation}\label{eq16}
\SS= f^b\overline{\SS} + \SS \subset \ldots \subset f^2\overline{\SS}
+ \SS
\subset f \overline{\SS} + \SS \subset \overline{\SS}.
\end{equation}
Filtration (\ref{eq15}) implies
\begin{eqnarray*}
\frac{t + 1}{t} \  \overline{\rme}_1(I) & = & \lambda(\overline{\SS}/\SS)
+ \frac{1}{t} \ \lambda(\overline{\SS}/\SS)  \ \ \qquad \qquad {\rm
by}
\ (\ref{eq4})\\
& \leq &  \lambda(f^{b}\overline{\SS}/f^{b} \SS) + \lambda(\SS/\SS \colon
\overline{\SS})   \ \qquad {\rm by}
\ (\ref{eq7}) \\
& \leq  & \lambda(\SS/f^{b} \SS)  \ \ \  \qquad \qquad \qquad \qquad
{\rm by} \
(\ref{eq15}) \\
& =  & b \ \lambda(\SS/f \SS) \\
& =  & b \ e_0(I) \qquad \qquad {\rm by \ the \ genericity \ of \
} a_1, \ldots, a_{d-1}.
\end{eqnarray*}
On the other hand filtration (\ref{eq16}) yields
\begin{eqnarray*}
\overline{\rme}_1(I)  & = & \lambda(\overline{\SS}/\SS) \\
& =&   \sum_{i=1}^{b}
 \lambda(f^{i-1}\overline{\SS} + \SS/ f^i\overline{\SS} + \SS) \\
& \leq & b \ \lambda(\overline{\SS}/f\overline{\SS} +\SS)   \qquad
\qquad {\rm by} \
(\ref{eq10})\\
& =  & b \ \lambda(\overline{I}/I) \qquad \qquad \qquad \ \  {\rm
by} \
(\ref{eq12}) \\
& =  &  b  \ (\rme_0(I)
 -\lambda(\RR/\overline{I})).
\end{eqnarray*}
\end{proof}

\begin{Remark}
{\rm The multiplicity $\rme_0(I+\delta \RR/\delta \RR)$ occurring in
Theorem~\ref{maintheorem} can be bounded by
\[
\rme_0(I+\delta \RR/\delta \RR)\leq (d-1)! \  \rme_0(I) \, \deg(\RR/\delta
\RR),
\]
where $\deg(\RR/\delta \RR)$ is the multiplicity of the local ring
$\RR/\delta \RR$. Indeed, \cite[Theorem 3]{Lech} gives
\[
\rme_0(I+\delta \RR/\delta \RR)\leq (d-1)! \  \lambda(\RR/I + \delta \RR) \,
\deg(\RR/\delta \RR) .\] }
\end{Remark}

\begin{Corollary}\label{RLR} Let $(\RR, \mathfrak{m})$
be a regular local ring of dimension $d>0$ and let $I$ be an
$\mathfrak{m}$--primary ideal. Then
\[
\rme_1(I) \leq \overline{\rme}_1(I) \leq (d-1) \ {\rm
min}\{\frac{\rme_0(I)}{2}, \rme_0(I) -\lambda(\RR/\overline{I})\}.
\]
\end{Corollary}
\begin{proof} We may assume that $\RR/\m$ is infinite. The classical
Brian\c{c}on--Skoda theorem gives that $c(I) \leq d-1$, see
\cite[Theorem 1]{LipmanSathaye}. The assertions now follow from
Theorems~\ref{bsjdegI} and \ref{maintheorem}(c). 
\end{proof}

We are now going to use Corollary~\ref{RLR} to bound the length of
divisorial chains for  classes of Rees algebras.

\begin{Corollary} Let $(\RR, \mathfrak{m})$ be a regular local ring of
dimension $d > 0$ and let $I$ be an $\mathfrak{m}$--primary ideal.
Then $(d-1) \ {\rm min}\{\frac{\rme_0(I)}{2}, \rme_0(I)
-\lambda(\RR/\overline{I})\} + 1$ bounds the length of any chain of
graded $\RR$--subalgebras satisfying the condition $S_2$ of Serre
lying between $\RR[It]$ and $\overline{\RR[It]}$.
\end{Corollary}

\begin{proof}  The assertion follows from Theorem~\ref{jdegintclos} and
Corollary~\ref{RLR}. 
\end{proof}

\begin{Remark}{\rm
 Let $(\RR, \mathfrak{m})$ be a regular local ring of
dimension $2$ and let $I$ be an $\mathfrak{m}$--primary integrally
closed ideal. Then $\rme_1(I)  = \overline{\rme}_1(I) = \rme_0(I)
-\lambda(\RR/I)\leq \frac{\rme_0(I)}{2}$. This follows, for instance,
from Corollary~\ref{RLR} combined with the inequality $\rme_1(I) \geq
\rme_0(I) - \lambda(\RR/I)$, see \cite[Theorem 1]{Northcott}. Furthermore by
\cite[2.1]{Hu2} or \cite[3.3]{O}, the equality $\rme_1(I) = \rme_0(I) -
\lambda(\RR/I)$ implies that $I$ has reduction number at most one if
$\RR/\mathfrak{m}$ is infinite, a fact proved in \cite[5.4]{LiTe}.
Thus $\RR[It]$ is Cohen--Macaulay according to \cite[3.1]{valval}
and \cite[3.10]{GS82}.
Now by
\cite[Theorem $2^{'}$,
p. 385]{Z}, $I$ is normal.
}\end{Remark}

\begin{Remark}{\rm  Let $k$ be an infinite field, write $\RR=k[X_1,
\ldots, X_d]_{(X_1, \ldots, X_d)}$, let $\mathfrak{m}$ denote the
maximal ideal of $\RR$, and let $I$ be an $\mathfrak{m}$--primary
$\RR$--ideal generated by homogeneous polynomials in $k[X_1, \ldots,
X_d]$ of degree $s$. Then $ \overline{\rme}_1(I) =
\rme_1(\mathfrak{m}^s)= \frac{d-1}{2} \rme_0(I) (1- \frac{1}{s})\approx
\frac{d-1}{2} \rme_0(I)$. This shows that the estimate of
Corollary~\ref{RLR} is essentially sharp.
}\end{Remark}

\begin{Proposition}\label{boundonBS} Let $k$ be a perfect
field, let $(\RR, \mathfrak{m})$  be a
reduced local Cohen--Macaulay $k$--algebra essentially of finite
type of dimension $d > 0$, and let $\delta \in {\rm Jac}_k(\RR)$ be
a non zerodivisor. Then for any $\mathfrak{m}$--primary ideal $I$,
\[
c(I) \leq  d-1 + \rme_0(I +\delta \RR / \delta \RR).
\]
\end{Proposition}
\begin{proof} We may assume that $\RR/\m$ is infinite. Then, replacing $I$
by a minimal reduction with the same Brian\c{c}on--Skoda number we
may suppose that $I$ is generated by a regular sequence of length
$d$. As in the proof of Theorem~\ref{maintheorem} let $\SS$ be a
local ring obtained from $\RR$ by a purely transcendental residue
field extension and by factoring out $d-1$ generic elements $a_1,
\ldots, a_{d-1}$ of $I$. Write $I\SS=f\SS$ and let $b$ be the smallest
non negative integer with $f^b \overline{\SS} \subset \SS \colon
\overline{\SS}$.
%\[
%b={\rm min} \, \{ \, n \geq 0 \ | \ f^n \overline{S} \subset S
%\colon \overline{S} \ \}.
%\]

We first claim that
\begin{equation}\label{eq17}
c(I) \leq b.
\end{equation}
Indeed, for any integer $n \geq 0$ we have
\[
\overline{I^{n + b }} \SS \subset \overline{I^{n + b} \SS} \subset
f^{n + b }\overline{\SS} \subset f^n \SS,
\]
hence $\overline{I^{n+b}}\SS \subset I^n\SS$.
%(I^n, a_1, \ldots, a_{d-1}).$
As ${\rm gr}_I(\RR)$ is a polynomial ring in $d$ variables over
$\RR/I$, the generic choice of $a_1, \ldots, a_{d-1}$ gives that
%the images of these variables in ${\rm gr}_{IS}(S)$ are still
%algebraically independent over $\RR/I$. Thus
${\rm gr}_I(\RR)$ embeds into ${\rm gr}_{I\SS}(\SS)$. Therefore
$\overline{I^{n+b}} \subset I^{n}$, proving (\ref{eq17}).

By (\ref{eq17}) it suffices to show that $b \leq d-1 + \rme_0(I +
\delta \RR/\delta \RR)$. To this end we may assume $b \geq d-1$. The
definition of $b$ yields the filtration

\begin{equation}\label{eq18}
\SS= f^b\overline{\SS} + \SS \subset \ldots \subset
f^{d-1}\overline{\SS}
+ \SS.
\end{equation}
On the other hand (\ref{eq5}) implies
\begin{equation}\label{eq19}
\SS= \delta f^{d-1}\overline{\SS} + \SS \subset f^{d-1}\overline{\SS}
+ \SS.
\end{equation}
If $f^{i}\overline{\SS} + \SS = f^{i-1}\overline{\SS} + \SS$ for some $b
\geq i \geq d$, then multiplication by $f^{b-i}$ yields
$\SS=f^b\overline{\SS} + \SS=f^{b-1}\overline{\SS} + \SS$, contradicting the
minimality of $b$. Thus (\ref{eq18}) gives
\[
\lambda(f^{d-1}\overline{\SS} + \SS/\SS) \geq b-d+1.
\]
On the other hand from (\ref{eq13}) and (\ref{eq19}) we deduce
\[
\lambda(f^{d-1}\overline{\SS} + \SS/\SS) \leq \rme_0(I + \delta \RR/\delta
\RR).
\]
Thus $b -d +1 \leq \rme_0(I + \delta \RR/\delta \RR).$ 
\end{proof}

\begin{Remark}
{\rm In the setting of the proof of Proposition~\ref{boundonBS},
$c(I)=c(I\SS)=b$. Indeed,  (\ref{eq3}) implies that for $n \gg 0$,
\[
f^{n + c(I)} \overline{\SS} = \overline{I^{n + c(I)} \SS} =
\overline{I^{n + c(I)}} \SS \subset I^n \SS = f^n \SS,
\]
showing that $f^{c(I)} \overline{\SS} \subset \SS.$ Hence $b \leq
c(I)$ and then $b =c(I)$ by (\ref{eq17}). Clearly $c(I\SS) \leq b$.
For $n \gg 0$, $f^{n + c(I\SS)}\overline{\SS} = \overline{f^{n +
c(IS)}\SS} \subset f^n \SS$ and therefore $f^{c(I\SS)} \overline{\SS}
\subset \SS$, showing that $b \leq c(I\SS)$.}
\end{Remark}

\subsection{Ideals of positive dimension}

\section{Sally Modules and Normalization of Ideals}

\subsection{Introduction} Let $(\RR, \mathfrak{m})$
be an analytically unramified local ring of
dimension $d$ and $I$ an $\mathfrak{m}$-primary ideal. Let
$\mathcal{F}= \{ F_n, \ n\geq 0\}$ be a decreasing, multiplicative
filtration of ideals, with $F_0=R$, $F_1= I$, with the property that
the corresponding Rees algebra
\[\Rees = R(\mathcal{F})  = \sum_{n\geq 0} F_nt^n \]
is a Noetherian ring. We will examine in detail the case when
$\mathcal{F}$ is a subfiltration of the integral closure filtration of
the powers of $I$, $F_n\subset \overline{I^n}$.

There are several algebraic structures attached to $\mathcal{F}$,
among which we single out  the associated graded ring of
$\mathcal{F}$ and its Sally modules. The first is
\[ \gr_{\mathcal{F}}(\RR)
= \sum_{n\geq 0} F_n/F_{n+1},\]
whose properties are closely linked to $\RR(\mathcal{F})$. We note that
if $J$ is a minimal reduction of $F_1$, then $\gr_{\mathcal{F}}(\RR)$ is
a finite generated module over $\gr_J(\RR)$, so that it a semi-standard
graded algebra.

To define the Sally module,
 we choose a minimal reduction $J$ of $I$ (if need be, we may
assume that the residue field of $\RR$ is infinite).
Note that $\Rees$ is a finite extension of the Rees algebra $\Rees_0=
\RR[Jt]$ of the ideal $J$. The corresponding Sally module $\SS$ defined by
the
 exact sequence of finitely generated modules over
$\Rees_0$,
\begin{eqnarray}
 0 \rar F_1\Rees_0 \lar \Rees_{+}[+1] \lar S = \bigoplus_{n=
 1}^{\infty}
F_{n+1}/J^nF_1 \rar 0.
\label{sallymod}
\end{eqnarray}

As an $\Rees_0$-module, $\SS$ is annihilated by an
$\mathfrak{m}$--primary ideal. If $S\neq 0$, $\dim S \leq d$, with
equality if $\RR$ is Cohen-Macaulay. The
Artinian module
\[S/JtS = \bigoplus_{k\geq 1} F_{k+1}/JF_k \]
gives some control over the number of generators of $\Rees$ as an
$\Rees_0$-module.
 If $S$ is
Cohen-Macaulay, this number is also its multiplicity.
It would, of course, be more useful to obtain bounds for  the
length of $\Rees/(\mathfrak{m}, \Rees_{+})\Rees$, but this requires
lots more.

\subsection{General properties of Sally modules}
The cohomological properties of $\Rees$, $\gr_{\mathcal{F}}(\RR)$ and
$S$ become more entwined when $\RR$ is Cohen-Macaulay. Indeed, under this
condition, the exact sequence (\ref{sallymod}) and
 the
 exact sequences (originally paired in \cite{Hu0}):

\begin{eqnarray}
0 \rar  \Rees_{+}[+1] \lar \Rees \lar {\rm gr}_{\mathcal{F}}(\RR)
\rar 0\  & & \label{hunekeeq1} \\
& & \nonumber \\
0 \rar \Rees_{+} \lar \Rees \lar R \rar 0,   & & \label{hunekeeq2}
\end{eqnarray}
with the tautological  isomorphism
\begin{eqnarray*}
&&  \Rees_{+}[+1] \cong \Rees_{+}
\end{eqnarray*}
gives a fluid mechanism to pass cohomological information around.

\begin{Proposition} Let $(\RR, \mathfrak{m})$ be a Cohen-Macaulay local
ring of dimension $d$ and $\mathcal{F}$ a filtration as above. Then
\begin{itemize}
\item[{\rm (a)}]
 $\depth \Rees \leq \depth \gr_{\mathcal{F}}(\RR)+1  $,
 with equality if
 $\gr_{\mathcal{F}}(\RR)$ is not Cohen-Macaulay.

\item[{\rm (b)}]  $\depth S\leq \depth \gr_{\mathcal{F}}(\RR) +1 $,
 with equality if
 $\gr_{\mathcal{F}}(\RR)$ is not Cohen-Macaulay.
\end{itemize}
\end{Proposition}

\begin{proof} For (a), see \cite[Theorem 3.25]{icbook}. For (b), it follows
simply because $F_1\Rees_0$ is a maximal Cohen--Macaulay
$\Rees_0$-module. \end{proof}

\subsection{Hilbert functions}

Another connection between $\mathcal{F}$ and $S$ is realized via their
Hilbert functions. Set
\[ H_{\mathcal{F}}(n) = \lambda(\RR/F_n), \quad H_{S}(n-1)=
\lambda(F_{n}/IJ^{n-1}).\]
The associated Poincar\'{e}-series
 \[\bbp_{\mathcal{F}}(t) = \frac{f(t)}{(1-t)^{d+1}},\]
 \[\bbp_{S}(t) = \frac{g(t)}{(1-t)^{d}}\]
are related by
\begin{eqnarray*}
 \bbp_{\mathcal{F}}(t) &=&  {\displaystyle\frac{\lambda(\RR/J)\cdot
 t}{(1-t)^{d+1}}  +
\frac{\lambda(\RR/I)(1-t)}{(1-t)^{d+1}}
 -H_S(t)}\\
&=& {\displaystyle \frac{\lambda(\RR/I)
 + \lambda(I/J)\cdot t}{(1-t)^{d+1}} - H_S(t)}.
\end{eqnarray*}

\begin{Proposition} The $h$-polynomials $f(t)$ and $g(t)$ are related
by
\begin{eqnarray}\label{relh}
f(t) &=&  \lambda(\RR/I)+ \lambda(I/J)\cdot t - (1-t)g(t).
\end{eqnarray}
In particular, if $f(t)= \sum_{n\geq 0}a_it^i$ and $g(t) = \sum_{n\geq
1}b_it^i$, then for $i\geq 2$
\[ a_i= b_{i-1}-b_i.\]
\end{Proposition}

\begin{Corollary} \label{decreasing} If moreover $\gr_{\mathcal{F}}(\RR)$ is
Cohen-Macaulay, then the $h$-vector of $S$ is positive and  non-increasing,
\[ b_i\geq 0, \quad b_1\geq b_2 \geq \cdots \geq 0.\]
In particular, if $b_{k+1}= 0$ for some $k$, then $\Rees$ is generated by its
elements of degree at most $ k$.
\end{Corollary}

\begin{proof} That $b_i\geq 0$ follows because $S$ is Cohen--Macaulay, while
the positivity of the $a_i$'s for the same reason and the difference
relation shows that $b_{i-1}\geq b_i$.
For the other assertion, since $S$ is Cohen-Macaulay, $b_k=
\lambda(F_{k+1}/JF_k)$ so that if it vanishes no fresh generators for
$F_{k+1}$ are needed. 
\end{proof}

\begin{Remark}\label{jun18-05}{\rm The equality (\ref{relh}) has several useful general properties,
of which we remark the following. For $k\geq 2$, one has
\[ f^{(k)}(1) = k g^{(k-1)}(1),\] that is the coefficients $\rme_i$ of the
Hilbert polynomials of $\gr_{F}(\RR)$ and $\SS$ are identical, more
precisely
\[ \rme_{i+1}(\mathcal{F}) =\rme_i(S), \quad i\geq 1.\]
Observe that when $\depth \gr_{\mathcal{F}}(\RR)\geq d-1$, $S$ is
Cohen--Macaulay so its $h$-vector is positive, and therefore all the
$\rme_i$ along with it (see \cite[Corollary 2]{Marley}). }
\end{Remark}

\begin{Corollary}\label{deg4}p
If  $\gr_{\mathcal{F}}(\RR)$ is
Cohen-Macaulay and  $g(t)$ is a polynomial of degree at most $4$,
then
\[ \rme_{2}(\mathcal{F})
 \geq \rme_3(\mathcal{F}) \geq \rme_4(\mathcal{F}) \geq \rme_5(\mathcal{F}). \]
\end{Corollary}

\begin{proof} By our assumption  \[ g(t) = b_1t + b_2t^2 + b_3t^3 +
b_4t^4.\] As
\[ \rme_{k+1}(\mathcal{F}) = \rme_{k}(S) = \frac{g^{(k)}(1)}{k!},
\] we have the following equations:
\begin{eqnarray*}
  \rme_{2}(\mathcal{F})&=& b_1 + 2 b_2 + 3 b_3 + 4 b_4 \\
  \rme_{3}(\mathcal{F}) &=& b_2 + 3 b_3 + 6 b_4 \\
  \rme_{4}(\mathcal{F}) &=& b_3 +4 b_4 \\
  \rme_{5}(\mathcal{F}) &=& b_4
\end{eqnarray*}
Now the assertion follows because $ b_1\geq b_2 \geq b_3 \geq  b_4
\geq 0$ according to Corollary~\ref{decreasing}. 
\end{proof}

This considerably lowers the possible number of distinct Hilbert
functions for such algebras.

\medskip

\begin{Remark}
{\rm The assumptions of Corollary~\ref{deg4} are satisfied for
instance if $\dim \RR\leq 6$ and $\Rees$ is Cohen-Macaulay.}
\end{Remark}

\subsection{Sally modules and Hilbert coefficients}
Depth estimates, the Sally module provides a quick way to establish the positivity of certain Hilbert
coefficients.

\begin{Theorem}\label{DepthSally} Let $(\RR,\mathfrak{m})$ be a Cohen--Macaulay
local ring of dimension $d$ and infinite residue field. Let $I$ be an
ideal generated by a system of parameters,
 $\AA= R[It]$ and $\BB= \sum_{n\ge
0}I_nt^n$ a finite $\AA$--subalgebra of $\RR[t]$.  Define the
Sally module $S_{\BB/\AA}$  by the exact sequence
\[ 0 \rar I_1\AA \lar \BB_{+} \lar S_{\BB/\AA} \rar 0.\]
If $I$ is a primary ideal and $J$ is a minimal reduction, we set
$S_{J}(I)$ for the corresponding Sally module.
Then: 
\begin{enumerate}
\item[{\rm (a)}] $S_{\BB/\AA}$ is an $\AA$--module that is either $0$ or has Krull
dimension $d$, more precisely, has the condition $S_1$ of Serre.
\item[{\rm (b)}] If $J$ is a minimal reduction of $I_1$ and $\AA=\RR[Jt]$,
the multiplicity of $S_{\BB/\AA}$ is $s_0(S_{\BB/\AA})=
\rme_1(\BB)-\rme_0(\BB)+
\lambda(\RR/I_1)$.
\item[{\rm (c)}] The first $3$ Hilbert coefficients,
 $\rme_0(\BB),\rme_1(\BB), \rme_2(\BB)$, of the function $\lambda(\RR/I_n)$ are
non-negative. Furthermore, if $\dim \RR=2$, $\rme_2(I)=0$ if only if the Sally module of $I^q$ vanishes for $q\gg 0$.

\end{enumerate}
\end{Theorem}

\begin{proof}
We only prove the assertions regarding $\rme_2(I)$ (Narita's Theorem \cite{Narita}).

We first consider the case for $d=2$ and without great loss of
generality take $I_n=I^n.$ Let
\[ \lambda(\RR/I^n)= \rme_0{{n+1}\choose{2}} -\rme_1{{n}\choose{1}} + \rme_2,
\quad n \gg 0\]
be the Hilbert polynomial of $I$. Let $q$ be an integer greater than
the postulation number for this function. We have
\begin{eqnarray*} \lambda(\RR/I^{qn}) = \lambda(\RR/(I^q)^n)&=& f_0{{n+1}\choose{2}} -f_1{{n}\choose{1}} +
f_2 \\ &=& \rme_0{{qn+1}\choose{2}} -\rme_1{{qn}\choose{1}} + \rme_2,
\quad n \gg 0.\end{eqnarray*}
Comparing these two polynomials in $n$, gives
\begin{eqnarray*}
f_0 & = & q^2\rme_0 \\
f_1 & = & q\rme_1 + \frac{1}{2}q^2\rme_0-\frac{1}{2}q\rme_0 \\
f_2 &=& \rme_2.
\end{eqnarray*}
By \cite[Theorem 2.11]{icbook}, the multiplicity $s_0(I^q)$ of the
Sally module of $I^q$ is given by
\[ s_0(I^q) = f_1-f_0+ \lambda(\RR/I^q),\] which by the equalities
above yields
$\rme_2(I) = s_0(I^q)$, the multiplicity of the Sally module of $I^q$.
Thus $\rme_2(I)$ is nonnegative and vanishes if
 and only if the reduction number of $I^{q}$ is less than or equal to $1$.

In dimension $d> 2$, we reduce $\RR$ and $I$ modulo a superficial
element.
\end{proof}

\subsection{Number and distribution of generators of normalizations}

Another application of Sally modules is to obtain a bound for the
number of generators (and the distribution of their degrees)
of $\Rees$ as an $\Rees_0$-algebra. (Sometimes the notation is used to
denote the number of module generators.)
A distinguished feature is the front loading of the new generators in
the Cohen-Macaulay case.

\begin{Theorem} Let $\mathcal{F}$ be a filtration as above.
\begin{itemize}
\item[{\rm (a)}]
 If $\depth \gr_{\mathcal{F}}(\RR)\geq d-1$,
 the $\Rees_0$-algebra $\Rees$ can be
generated by $\rme_1(\mathcal{F})$ elements.
\item[{\rm (b)}]
 If $\depth \gr_{\mathcal{F}}(\RR)= d$,
 the $\Rees_0$-algebra $\Rees$ can be
generated by $\rme_0(\mathcal{F})$ elements.
\item[{\rm (c)}]
 If $\Rees$ is Cohen-Macaulay, it can be
 generated by $\lambda(\mathcal{F}_1/J) + (d-2)
\lambda(\mathcal{F}_2/\mathcal{F}_1J)$ elements.
\end{itemize}
%\[ \nu(\Rees) \leq \lambda(F_1/J) +  (e_1(\mathcal{F})-e_0 +
%\lambda(\RR/F_1))\lambda(F_2/JF_1)
%\]
%  Moreover, if $\Rees$ is Cohen-Macaulay, then
%\[ \nu(\Rees) \leq \lambda(F_1/J) +  (d-1)\lambda(F_2/JF_1).\]
\end{Theorem}

\begin{proof}
%Since the $h$-polynomial $g(t)$ of $\SS$ has no intermediate
%coefficient that vanishes, the degree of $g(t)$ is bounded by $g(1)$,
%that is the multiplicity of $\SS$.
 (a) From the relation (\ref{relh}), we
have
\[\rme_1(\mathcal{F}) = f'(1) = \lambda(F_1/J) + g(1).\]
From the sequence (\ref{sallymod}) that defines $\SS$, one has
\[ \nu(\Rees) \leq \nu(F_1/J) + \nu(S)\leq \lambda(F_1/J) + \nu(S).\]
Since $\SS$ is Cohen-Macaulay, $\nu(S)= \rme_0(S)=g(1)$, which combine to
give the
 promised assertion.

\medskip

(b) A generating set for $\Rees$ can be obtained from a lift of a
minimal set of generators for $\gr_{\mathcal{F}}(\RR)$, which is given
by its multiplicity since it is Cohen-Macaulay.

\medskip

(c) Since $\Rees $ is Cohen-Macaulay, its reduction number is $\leq
d-1$. Thus the $h$-polynomial of $\gr_{\mathcal{F}}(\RR)$ has degree
$\leq d-1$, and consequently the $h$-polynomial of the Sally module
has degree at most $d-2$. As the $h$-vector of
$S$ is decreasing, its multiplicity is at most
$(d-2)\lambda(\mathcal{F}_2/\mathcal{F}_1J)$, and we conclude as in
(a).
 \end{proof}

\begin{Remark}{\rm
A typical application is to the case $d=2$ with $\mathcal{F}$ being
the filtration $F_n= \overline{I^n}$.
}\end{Remark}

\begin{Remark}{\rm There are several relevant issues here. The first, to get bounds for
$\rme_1(\mathcal{F})$. %, then to estimate $\lambda(F_2/JF_1)$.
 This
is addressed in \cite{ni1}. For instance, when $\RR$ is a regular local
ring of characteristic zero, $\rme_1(\mathcal{F})\leq
\displaystyle{\frac{(d-1)\rme_0}{2}}$.

%When $\gr_{F}(\RR)$ is Cohen-Macaulay

%To bound $\lambda(F_2/JF_1)$ there
%are several choices, a simple one being
%\begin{eqnarray*} \lambda(F_2/JF_1)&=&\lambda(\RR/JF_1)-\lambda(\RR/F_2) =
%\lambda(\RR/J) + \lambda(J/JF_1)-\lambda(\RR/F_2) \\ &\leq & e_0 +
%d\lambda(\RR/F_1)-\lambda(\RR/F_1) = e_0 + (d-1) \lambda(\RR/F_1)
%.\end{eqnarray*}
}\end{Remark}

%\bigskip

To further these calculations, natural questions are:
\begin{itemize}
%\item Refine the estimations for $\lambda(F_2/JF_1)$. Experiment

\item[$\bullet$] How far away from $\depth \gr_{\mathcal{F}}(\RR)\geq d-1$  can we
go? In particular find classes of filtrations whose $h$-polynomials
are positive.

\item[$\bullet$] Extend estimations to equimultiple ideals. How to approach
more general reductions?

\item[$\bullet$] Look for connections in the spirit of \cite[Theorem 2]{Itoh}.

\item[$\bullet$] How is the issue of no-gaps in the $h$-vector of $S$ related to
the old issue of no-gaps in the canonical module of the Rees algebra?

\item[$\bullet$] In order to extend the bounds for $\nu(\Rees)$, it would be
helpful to prove a version of the Eakin-Sathaye theorem for non-adic
filtrations. One guess is that there is already a version of this in
\cite{icbook}.

\end{itemize}

%(\subsection{$\rme_2(I)$: what is its meaning?}

\section{Normalization  of Modules}

\subsection{Introduction}

Let us recall the notion of the Rees algebra of a module. Let $\RR$ be
a Noetherian ring, let $E$ be a finitely generated torsion free
$\RR$--module having a rank, and choose an embedding $\varphi: E
\hookrightarrow R^r$. The {\em Rees algebra} $\Rees(E)$ of $E$ is
the subalgebra of the polynomial ring $\RR[t_1, \ldots, t_r]$
generated by all linear forms $a_1t_1 + \cdots + a_rt_r$, where
$(a_1 , \ldots , a_r)$ is the image of an element of $E$ in $\RR^r$
under the embedding $\varphi$. The Rees algebra $\Rees(E)$ is a
standard graded algebra whose $n$th component is denoted by $E^n$
and is independent of the embedding $\varphi$ since $E$ is
torsion free and has a rank.

There is a more general notion of Rees algebra of a module (see
\cite{EHU} and the discussion in \cite[pp. 10--11]{icbook}). One of
the reasons for taking the approach here is to  benefit from
the properties of   a fixed embedding
$E\subset R^r$, but also to define additional  structures.

The algebra $\Rees(E)$ is a subring of the polynomial ring $\SS=\RR[t_1,
\ldots, t_r]$. Following \cite{HUV}, we consider the ideal $(E)$ of
$\SS$ generated by the forms $a_1t_1 + \cdots + a_rt_r$.
 Its Rees algebra $\Rees((E))$
(defined over $\SS$) is much larger than $\Rees(E)$, but has an
associated graded  ring $G=\gr_{(E)}(\SS)$, something that $\Rees(E)$
lacks--a fact that  makes the  study of $\Rees(E)$
more difficult than the Rees algebras of
ideals.

\subsection{Reductions and integral closure of modules}

There are several measures of size attached to a Rees algebra
$\Rees(E)$, all derived from ordinary Rees algebras. We survey them
briefly.

\begin{proposition} \label{2.2}
Let $\RR$ be a Noetherian ring of
dimension $d$ and $E$  a finitely generated $\RR$-module of rank
$r$. Then
$$ \dim \Rees(E) = d+r = d + {\rm height} ~~ \Rees(E)_+.
%\label{Krulldim}
$$
\end{proposition}

\begin{proof} See \cite[Proposition 8.1]{icbook}.
\end{proof}

By analogy with the case of ideals, one defines:

\begin{definition}{\rm
 Let $U $ be a submodule of $E$.  We
say that $U$ is a {\em reduction}\index{reduction of a module}
 of $E$, or equivalently $E$ is {\em
integral} over $U$, if $\Rees(E)$ is integral over the $\RR$-subalgebra generated
by $U$.
}
\end{definition}

Alternatively, the integrality condition is expressed by the equations
$\Rees(E)_{s+1} = U \cdot \Rees(E)_s$, $ s \gg 0$.  \ The least
integer $s \geq
0$ for which this equality holds is called the {\em reduction
number}\index{reduction number of a module} of $E$
{\em with respect to} $U$ and denoted by $r_U(E)$.  For any reduction $U$ of $E$
the module $E/U$ is torsion, hence $U$ has the same rank as $E$.
  This follows from the fact that a module of linear type, such as a
free module, admits no proper reductions.

\medskip

Let $E $ be a submodule of $ \RR^r$.
The {\em integral closure}\index{integral closure of a module} of
$E$ in $\RR^r$ is the largest submodule $\overline{E}\subset \RR^r$
having $E$ as a reduction. (A broader definition exists, but we shall
limit ourselves to this case.)

\medskip

If $\RR$ is  a local ring with residue field $k$
then the {\em special fiber} of
 $\Rees(E)$ is the ring ${\mathcal F}(E)= k \otimes_\RR \Rees(E)$;
its Krull dimension
 is called the {\em analytic spread\/}\index{analytic spread of a
module}
 of $E$ and is denoted by $\ell(E)$.

Now assume in addition that $k$ is infinite.  A reduction of $E$ is said to be
{\em minimal} if it is minimal with respect to inclusion.
  For any reduction $U$
of $E$ one has $\nu (U) \geq \ell(E)$ (here $\nu(\cdot)$ denotes the minimal number
of generators function), and equality holds if and only if $U$ is minimal.
Minimal reductions arise from the following construction.
The algebra ${\mathcal
F}(E)$ is a standard graded algebra of dimension $\ell = \ell(E)$ over the
infinite field $k$.  Thus it admits a Noether normalization $k[y_1, \ldots,
y_{\ell}]$ generated by linear forms; lift these linear forms to elements $x_1,
\ldots, x_{\ell}$ in $\Rees(E)_1= E$, and denote by $U$ the submodule generated
by $x_1, \ldots, x_{\ell}$.  By Nakayama's Lemma, for all large $r$ we have
$\Rees(E)_{r+1} = U \cdot \Rees(E)_r$, making $U$ a minimal reduction of $E$.

Having established the existence of minimal reductions, we can define the {\em
reduction number} $r(E)$ of $E$ to be the minimum of $r_U(E)$, where $U$ ranges
over all minimal reductions of $E$.

\begin{proposition} \label{2.3a} Let $(\RR, \mathfrak{m})$ be a
Noetherian local ring with infinite residue field,
of dimension
$d \geq 1$ and  $E$  a finitely generated $\RR$-module having
 rank $r$.
Then
\begin{eqnarray}
\begin{array}{c}
\  r\leq \ell(E) \leq d+r-1.  \\
\end{array}
\end{eqnarray}
 Moreover,   if  $E\hookrightarrow \RR^r$ and $0\neq
 \lambda(\RR^r/E)< \infty$, then
$ \ell(E) = d+r-1  $.
 \end{proposition}

\begin{proof} See \cite[Propositions 8.3, 8.4]{icbook}.
\end{proof}

%\bigskip

\subsubsection*{Integral Closure and Normality} Let $E$ be a module of rank
$r$. Having considered the notion of reductions of $E$ we make the
following definition.
  This
portion of our exposition is highly dependent on \cite{HUV} and
\cite[Chapter 8]{icbook}.

\begin{definition}{\rm Let $E\hookrightarrow \RR^r$ be an embedding of
modules of the same rank. Then
$E$ is {\em integrally closed in $\RR^r$}\index{integrally closed
module} if it is not a proper reduction of any submodule of $\RR^r$.
Further, $E$ is {\em normal in $\RR^r$}\index{normal module} if all components
$\Rees(E)_n$ are integrally closed.
}\end{definition}

\begin{remark}{\rm As defined thus far, the integral closure of a
module has a relative character. Thus if $E$ allows two embeddings,
$\varphi_1: E \rar \RR^r$ and $\varphi_2: E \rar \RR^s$, the integral
closures of $E$ may be different. An {\em absolute  integral
closure}\index{integral closure of a module!absolute closure} for $E$
could be defined as follows: let $K$ be the field of fractions of
$\RR$, and set $T= K\otimes_\RR E$. For any valuation overring $V$ of $\RR$, we
can consider the image of $V\otimes_\RR E$ in $K\otimes_\RR E=T$,
 $VE \subset T$. The intersection $\overline{E} $ of these groups,
 where
\[ E\subset \overline{E}=\bigcap_{V}VE \subset T,\]
is the desired closure. Note that when $\RR$ is integrally closed it
is easy to see that the integral closure of $E$ in any embedding
 coincides with the absolute integral closure.
}\end{remark}

\begin{proposition}\label{intclosrees} Let $\RR$ be a normal domain and
 $E$  a torsion free finitely generated  $\RR$-module. The Rees algebra
\[ \Rees(E) = \RR \oplus E \oplus  E_2 \oplus \cdots\]
is integrally closed if and only if each $\rme_n$ is an  integrally
closed module.
\end{proposition}

\begin{proof} If each component $\rme_n$ is integrally  closed, we have
\[ \bigcap_{V}V\Rees(E) = \sum_{n\geq 0}\bigcap_{V}VE_n = \Rees(E),\]
$V$ running over all the valuation overrings of $\RR$. For each $V$ we
have
$V\Rees(E) = \Rees(VE)$, where $VE$ is a  finitely generated
torsion free $V$-module and therefore $\Rees(VE)$ is a ring of
polynomials over $V$ since $VE$ is a free $V$-module. This gives a
representation of $\Rees(E)$ as an intersection of polynomial rings
and it is thus normal. The converse is similar. 
\end{proof}

Let $\RR$ be an integral domain, and
let $E\subset \RR^r$ be an $\RR$-module of rank $r$. Suppose $U\subset E$
is a submodule. Denote  $\SS= \RR[\TT_1, \ldots, \TT_r]$.

\begin{proposition} \label{UredE} $U$ is a reduction of $E$ if and only if the
$\SS$-ideal $(U)$ is a reduction of the ideal $\SS$-ideal $(E)$.
\end{proposition}

\begin{proof} Another way to express that $U$ is a reduction of $E$ is to say
that $U$ and $E$ have the same integral closure in $\RR^r$, in other
words,
for every valuation $V$ of $\RR$, $VU=VE\subset V^r$. A similar
observation applies to the ideals $(U)$ and $(E)$ of $\SS$. In this
case, instead of using all the valuations of $\SS$, it suffices to
verify that for every valuation $V$ of $\RR$, $(U)V[\TT_1, \ldots, \TT_r] =
(E)V[\TT_1, \ldots, \TT_r]$, which is equivalent to $UV=EV$. 
\end{proof}

The invariants of  the Rees algebras of $E$ and $(E)$ may be
different. It is easy to see that $\ell(E)= \ell((E))$, the latter
the special fiber relative to the maximal graded ideal. On the other
hand, while clearly $\red_U(E)\leq \red_{(U)}((E))$, it is not clear
that equality will always hold. The differences will be sharpened
when we discuss the Brian\c{c}on-Skoda numbers of $E$ and $(E)$.

\subsection{Buchsbaum-Rim multiplicity}\index{Buchsbaum-Rim
multiplicity}

Let us recall  the definition of Buchsbaum--Rim polynomials.
Let $\RR$ be a Noetherian local ring of dimension $d$ and let $E
\subsetneq \RR^r$ be a submodule such that the length $\lambda(\RR^r/E)$
is finite. Let $E^n$ denote the image of the $\RR$--linear map
$S_n(E) \rightarrow \SS=\bigoplus_{n \geq 0}S_n=S(\RR^r)=\RR[\TT_1,
\ldots, \TT_e]$. Buchsbaum and Rim proved that the length
$\lambda(S_n/E^n)$ is a polynomial in $n$ of degree $d+r-1$ for
sufficiently large $n$ (\cite[3.1 and 3.4]{BR}). This polynomial is
called the {\em Buchsbaum--Rim polynomial} of $E$, which is of the
form

\[P(n)= \br(E){{n+d+r-2}\choose{d+r-1}} -
\br_1(E) {{n+d+r-3}\choose{d+r-2}}+ \textrm{\rm lower terms.} \] The
positive integer $\br(E)$ in the Buchsbaum--Rim polynomial is called
the {\em Buchsbaum--Rim multiplicity} of $E$. For any reduction $U$
of $E$, it is known that $\br(U)=\br(E)$ (\cite[3.1]{BUV}).

Let $\RR$ be an analytically unramified normal local ring. Then since
the integral closure $\overline{\Rees(E)}$ is finite over
$\Rees(E)$, one obtains the following polynomial expression for $n
\gg 0$,

\[\lambda(S_n/\overline{E^n})= \overline{P}(n) =
\overline{\br}(E){{n+d+r-2}\choose{d+r-1}}  - \overline{\br}_1(E)
{{n+d+r-3}\choose{d+r-2}}+ \textrm{\rm lower terms.}\] Notice that
$\br(E)=\overline{\br}(E)$ because $\lambda(\overline{E^n}/E^n)$ is
a polynomial of degree at most $d+r-2$ for $n \gg 0$ and
$\br(E)-\overline{\br}(E) \geq 0$. More generally if
$\mathcal{A}=\bigoplus_{n \geq 0} A_n$ is a graded subalgebra such
that $\Rees(E)\subset \mathcal{A} \subset \overline{\Rees(E})$, then
for $n \gg 0$,

\[
\lambda(S_n/A_n) =P_{\mathcal{A}}(n) =
\br(\mathcal{A}){{n+d+r-2}\choose{d+r-1}} - \br_1(\mathcal{A})
{{n+d+r-3}\choose{d+r-2}}+ \textrm{\rm lower terms.}
\]

Let $\RR$ be a Noetherian local ring and let $E \neq 0$ be a finitely
generated torsion free $\RR$--module. The module $E$ is called an {\em
ideal module} if $E^{**}$ is free, where ${}^*$ denotes dualizing
into the ring $\RR$. For example, if $\RR$ is Cohen--Macaulay ring of
dimension at least $2$, then a submodule $E \subsetneq \RR^r$ having
finite length $\lambda(\RR^r/E)$ is an ideal module (\cite[5.1]{ram1}).
An ideal module $E$ has a rank(\cite[5.1]{ram1}), say $e$, and
affords a {\it natural} embedding $E \hookrightarrow E^{**}=\RR^r$.
Composing an epimorphism $\RR^n \twoheadrightarrow E$ with this
embedding, we obtain the $\RR$--linear map $\varphi: \RR^n \rightarrow
\RR^r$ with $\image(\varphi)=E$. Notice that the $i$th Fitting ideal
$\Fitt_i(\RR^r/E)$ of $\RR^r/E$ is the $\RR$--ideal generated by $(e-i)
\times (e-i)$--minors of $\varphi$, and that these ideals only
depend on $E$. Hence we may write $\ddet0(E)$ for $\Fitt_0(\RR^r/E)$.
We define the {\em codimension} of $E$ as the height of $\ddet0(E)$,
the {\em deviation} of $E$ as $\d(E)= \nu(E) - e+1
-\grade(\ddet0(E))$, and the {\em analytic deviation} as $\ad(E)=
\ell(E) -e +1 -\grade(\ddet0(E))$. An ideal module $E$ is said to be
a {\em complete intersection}, an {\em almost complete
intersection}, or {\em equimultiple} if $\d(E) \leq 0$, $\d(E) \leq
1$, or $\ad(E) \leq 0$ respectively. At this point we need the
notion of $\mathfrak{m}$--full modules. A submodule $E \subset \RR^r$
over a local ring $(\RR, \mathfrak{m})$ is called {\em
$\mathfrak{m}$--full} if there is an element $x \in \mathfrak{m}$
such that $\mathfrak{m}E:_{\RR^r} x =E$. Now we are ready to find
bounds for the Buchsbaum--Rim multiplicity.

\begin{Proposition}\label{nubar} Let $(\RR, \mathfrak{m})$ be an analytically unramified Cohen--Macaulay
normal local ring of dimension $d \geq 2$ and let $E \subsetneq \RR^r$
be a submodule such that $\lambda(\RR^r/E) < \infty$. Let $s$ be the
integer such that $\Fitt_{e-1}(\RR^r/E)$ is contained in
$\mathfrak{m}^s$ but not in $\mathfrak{m}^{s+1}$. %Suppose that the integral closure $\overline{\Rees(\mathfrak{m})}$ is finite over the Rees algebra $\Rees(\mathfrak{m})$ of the maximal ideal $\mathfrak{m}$.
Then
\[  e(\mathfrak{m}) s^{d-1}\cdot { {d+r-2}\choose
{e-1}} \leq \br(E) \leq e(\ddet0(E)) \cdot {{d+r-1}\choose{e-1}},
\]where $e(\mathfrak{m})$ and $e(\ddet0(E))$ are the Hilbert--Samuel multiplicities of the ideals
$\mathfrak{m}$ and $\ddet0(E)$ respectively.
\end{Proposition}

\begin{proof} Let $\SS=\bigoplus_{n \geq 0} S_n =S(\RR^r)=\RR[\TT_1, \ldots,
\TT_e]$
and let $\RR^n \stackrel{\varphi}{\longrightarrow} \RR^r \longrightarrow
\RR^r/E \longrightarrow 0$ be a presentation of $\RR^r/E$. Since
$\Fitt_{e-1}(\RR^r/E)=I_1(\varphi) \subseteq \mathfrak{m}^s$ and
$E=\image(\varphi)$, we have
%\[ E^n \subseteq \mathfrak{m}^{sn} S_n \quad \mbox{\rm for all}\; n
%\geq 0.\] Hence
\[ \overline{E^n} \subseteq \overline{\mathfrak{m}^{ns}}S_n  \quad \mbox{\rm for all}\; n
\geq 0.\] Since integrally closed modules are $\mathfrak{m}$--full
(\cite[2.6 and 2.7]{BV}), we obtain
\[ \nu(\overline{\mathfrak{m}^{ns}} S_n) \leq \nu(\overline{E^n}).\]
Now for sufficiently large $n$, we have a polynomial expression
\[\lambda \left(\overline{\mathfrak{m}^{ns}}/\overline{\mathfrak{m}^{ns+1}} \right)=
\overline{e_0}{{ns+d-1}\choose{d-1}}-\overline{e_1}{{ns+d-2}\choose{d-2}}
+ \cdots + (-1)^{d-1}\overline{e_{d-1}}  \leq
\nu(\overline{\mathfrak{m}^{ns}}).\] Let $U$ be a minimal reduction
of $E$. Then $U$ is a complete intersection and using \cite[3.1 and
4.4]{BUV}, we obtain

\begin{eqnarray*}
\nu(\overline{\mathfrak{m}^{ns}}) \cdot {{n+r-1}\choose{e-1}}
&=& \nu(\overline{\mathfrak{m}^{ns}} S_n)
\leq \nu(\overline{E^n})\\
& = & \nu(\overline{U^n}) \leq \br(U) \cdot
{{n+d+r-3}\choose{d+r-2}} + {{n+d+r-3}\choose{d+r-3}}.
\end{eqnarray*}
%Notice that
%\[\begin{array}{lll}
%{{n+d+r-3}\choose{d+r-2}} &=&
%\frac{n(n+r+d-3)(n+r+d-4)\cdots(n+r)}{(e+d-3)\cdots
%e}{{n+r-1}\choose{n}}\\ &&\\{{n+d+r-3}\choose{d+r-3}} &=&
%\frac{(n+r+d-3)(n+r+d-4)\cdots(n+r)}{(e+d-3)\cdots
%e}{{n+r-1}\choose{n}}\\
%\end{array}\]
Therefore we have
%%%%%%%%%%%%%%%%%%%%%%%%%%%
%\[e_0{{ns+d-1}\choose{d-1}}-e_1{{ns+d-2}\choose{d-2}} + \cdots +(-1)^{d-1}e_{d-1} \leq
%\deg(\RR/\ddet0(F))\frac{n(n+r+d-3)(n+r+d-4)\cdots(n+r)}{(e+d-3)\cdots
%e}{{n+r-1}\choose{n}}
%+\frac{(n+r+d-3)(n+r+d-4)\cdots(n+r)}{(e+d-3)\cdots e}\]
%%%%%%%%%%%%%%%%%%%%%%%%%%%

\[ \frac{\overline{e_0} s^{d-1}}{(d-1)!(e-1)!} n^{d+r-2} + \mbox{\rm lower terms} \leq
\frac{\br(U)}{(d+r-2)!}n^{d+r-2} + \mbox{\rm lower terms}.
\]
%\[ \deg(\RR/\ddet0(F)) \geq \frac{1}{(d-1)!} e_0 s^{d-1} (e+d-3)(e+d-4)\cdots e. \]
%\[ \br(E) \geq \frac{1}{(d-1)!} e_0 s^{d-1} (e+d-3)(e+d-4)\cdots e. \]
Since $\br(E)=\br(U)$ and $\rme_0=\overline{e_0}$, we have
\[ \br(E) \geq {e_0 s^{d-1}} \cdot { {d+r-2}\choose {e-1}}. \]
The other inequality follows from the fact that
\[ \lambda(S_n/E^n ) \leq \lambda(S_n/\ddet0(E)^n S_n) \quad \mbox{\rm for all}\; n \geq 0.
\] \end{proof}

\subsubsection{Buchsbaum-Rim multiplicity and $\jdeg$} Let
$E\subset \RR^r$ be a submodule of finite co-length.
Let
us show how to obtain how to obtain the Buchsbaum-Rim multiplicity
$\br(E)$ as an expression of an appropriate $\jdeg$. As a side
effect, it will permit the definition of a global Buchsbaum-Rim
multiplicity.

\begin{Proposition} Let $(\RR, \mathfrak{m})$ be a Noetherian integral
domain of dimension $d$ and let $E$ be a torsion free $\RR$--module of
rank $r$ with a fixed embedding $E \hookrightarrow \RR^r$. Then
\begin{itemize}
\item[{\rm (a)}] The components of $\GG = \bigoplus_{n\geq 0}
E^nS/E^{n+1}$ have a natural grading
\[ G_n = E^n + E^nS_1/E^{n+1} + E^nS_2/E^{n+1}S_1
 + \cdots . \]

\item[{\rm (b)}] There is a decomposition
\[ \GG= \Rees(E) + H,\]
where $\Rees(E)$ is
the Rees algebra
 of $E$ and $H$ is the $\RR$-torsion submodule of $G$.

\item[{\rm (c)}] If  $E\subset \mathfrak{m}\RR^r$ and $\lambda(\RR^r/E)<
\infty$, $H = \H^0_{\mathfrak{m}}(G)$ has dimension $d+r$ and
multiplicity equal to the Buchsbaum-Rim multiplicity of $E$.

\end{itemize}
\end{Proposition}

\begin{proof} Only Part (c) requires attention. The module
$H$ has a natural bigraded
structure as an $\Rees(E)[y_1, \ldots, y_d]$-module which we redeploy
in terms of total degrees as
follows
\[ H = \bigoplus_{n\geq 0} H_n, \quad H_n = \bigoplus_{i=0}^n
E^iS_{n-i}/E^{i+1}S_{n-i-1}.\]
Observe that $\lambda(H_n) = \lambda(S_n/E^n)$, the Buchsbaum-Rim
Hilbert function of the module $E$ (see
\cite{BR}). For large $n$,
\[ \lambda(S_n/E^n) = \textrm{\rm br}(E){{n+d+r-2}\choose{d+r-1}} +
\textrm{\rm lower terms}.
\] (The nonzero integer
$\textrm{\rm br}(E)$
is the Buchsbaum-Rim  multiplicity
of $E$ (see \cite[Section 8.3]{icbook} for details).)
The degree of this polynomial shows that $\dim H =d+r-1$. 
\end{proof}

We could show directly, using properties of Rees algebras, that $\dim
H=d+r-1$ and thereby prove the existence of the polynomial above.

\begin{Corollary} \label{jdegGE}
 Let $E$ be a module as above. Then $\jdeg (\GG) =
\textrm{\rm br}(E)+1$.
\end{Corollary}

\begin{proof} It suffices to note that $\Ass_R(\GG)= \{0, \mathfrak{m}\}$, and
use the calculation above for $j_{\mathfrak{m}}(\GG)$ and take into
account  that
$j_{(0)}(\GG) = \deg (\Rees(E)_{(0)})= \deg (K[\TT_1, \ldots, \TT_r])=1 $.
\end{proof}

We can extract additional degrees information from the decomposition
$\GG= H + \Rees(E)$.

\begin{definition} {\rm Suppose $E\hookrightarrow \RR^r$ is a submodule of
rank $r$. For each $ \mathfrak{p}\in \Spec(\RR)\setminus (0)$, the
Buchsbaum-Rim multiplicity of $E$ is the $j$-multiplicity of $H$,
\[\br_{\mathfrak{p}}(E)= j_{\mathfrak{p}}(H).\]
}
\end{definition}

When $\mathfrak{p}$ is a minimal prime of $\RR^r/E$, this definition
coincides with the Buchsbaum-Rim multiplicity of $E$ (or, of the
embedding $E\hookrightarrow \RR^r$).

\begin{Theorem}[{\cite[Theorem 1.1]{UV}}] Let $\RR$ be a universally catenary integral domain and
let $E$ be a torsion free $\RR$-module of rank $r$ with an embedding
$E\hookrightarrow \RR^r$. Suppose $U$ is a submodule of $E$ of rank
$r$, and fix the induced embedding $U\hookrightarrow \RR^r$. The following conditions are equivalent:
\begin{enumerate}
\item[{\rm (a)}] $\br(U_{\mathfrak{p}})= \br(E_{\mathfrak{p}})$ for each
$\mathfrak{p}\in \Spec(\RR)$.

\item[{\rm (b)}] $U$ is a reduction of $E$.
\end{enumerate}

\end{Theorem}

\begin{proof} Under the conditions, it suffices to apply
Theorem~\ref{FlennerM} together with Proposition~\ref{UredE}. (A more
general proof is that of \cite{UV}.) 
\end{proof}

\subsubsection{Buchsbaum-Rim coefficients}

Several questions could be asked about the Hilbert coefficients of
modules (see \cite{chern5}), that is about the range of values  of $\rme_1(J;M)$. We are going to raise some questions
however about
 modules and their Buchsbaum-Rim
coefficients.

\medskip

Let $(\RR,\m)$ be a
 Noetherian local ring of dimension $d$.
The Buchsbaum-Rim multiplicity (\cite{BR}) arises in the context of an embedding
\[ 0 \rar E \stackrel{\varphi}{\lar} F=\RR^r \lar C \rar 0,\]
where $C$ has finite length.
 Denote by
\[ \varphi: \RR^m \lar \RR^r\]
a matrix with entries in $\m$.
     There is a  homomorphism
\[ {\mathcal S}(\varphi): {\mathcal S}(\RR^m) \lar {\mathcal S}(\RR^r)\] of symmetric algebras
whose image is $\Rees(E)$, and whose cokernel we denote by
$C(\varphi)$,\begin{eqnarray}
0 \rar \Rees(E) \lar {\mathcal S}(\RR^r)=\RR[x_1, \ldots, x_r] \lar C(\varphi)\rar 0. \label{br1}
\end{eqnarray}
This exact sequence (with a different notation) is studied in
\cite{BR} in great detail. Of significance for us is the fact that
$C(\varphi)$, with the grading induced by the homogeneous homomorphism
${ S}(\varphi)$, has components of finite length which for $n \gg 0$:

\begin{Theorem}  \label{degofBR}
$\lambda(S_n(\RR^r)/E^n)$ is a polynomial in $n$ of degree $d+r-1$ for
$n\gg 0$,
\[ \lambda(S_n(\RR^r)/E^n) = \br(E) {{n+d+r-2}\choose{d+r-1}}-\br_1(E)
{{n+d+r-3}\choose{d+r-2}}
  + \textrm{\rm lower
terms}. \]
\end{Theorem}

This polynomial is
 called the {\em Buchsbaum-Rim polynomial} of
$E$\index{Buchsbaum-Rim coefficients of a module}.
The leading coefficient $\br(E)$ is the
 {\em Buchsbaum-Rim multiplicity} of $\varphi$; if the embedding
$\varphi$ is understood, we shall simply denote it by $\br(E)$. This
number  is determined by
an Euler characteristic of the Buchsbaum-Rim complex (\cite{BR}).

\medskip

Assume the residue field is infinite and $ \RR$ is unmixed.
The minimal reductions $U$ of $E$ are minimally generated
by $r+d-1$ elements. If $E=\RR^n$, we refer to $U$ as a {\em
parameter module}.
The corresponding coefficient $\br(U)=\br(E)$,
but $\br_1(U)\leq \br_1(E)$. It is not clear what
the possible values of $\br_1(U)$ are. In similarity to the case of
ideals, there are the following questions:

\begin{enumerate}

\item[{\rm (a)}] $\br_1(U)\leq 0$ according to \cite{HaHy10}
 with equality if $\RR$ is Cohen--Macaulay. When does a converse hold?

\item[{\rm (b)}] When are the values of $\br_1(E)$
bounded, or even constant?

\item[{\rm (c)}] What happens in low dimensions?

\end{enumerate}

Let us begin an exploration, beginning with $\dim \RR=1$. In this
case $U\subset F=\RR^r$ is a module generated by $r $ elements.
Let $F_0=\H_{\m}^0(F)$ and consider the exact complex

\[
\diagram
0 \rto & U_0 \rto \dto &U \rto \dto & U' \rto \dto & 0  \\
0 \rto & F_0 \rto  &F \rto  & F' \rto  & 0
\enddiagram
\]
with $U_0=F_0\cap U$. Note that $U'$ is a free $\RR'$-module, for
$\RR'=\RR/ \RR_0$, $\RR_0=\H_{\m}^0(\RR)$.

The Buchsbaum-Rim polynomial of $U$ is easily seen to be
\[ \lambda(\RR'/(\det \varphi'))
 {{n+r-1}\choose{r}}-\br_1(U)
{{n+r-2}\choose {r-1}}
  + \textrm{\rm lower
terms}. \]

The coefficient $\br_1(U)$  comes (with the sign reversed) from the
multiplicity of $\RR_0(x_1, \ldots, x_r)^n$ (Recall that by
Artin-Rees $  (\RR_0\cdot  F^n)\cap U^n=0 $, for $n\gg 0$),  that is
\[ \br_1(U) = -r\cdot \lambda(\H_{\m}^0(\RR)).\]

To advance, we need a technique like that of Proposition~\ref{genhs},
but applied to Buchsbaum-Rim coefficients.

\begin{conjecture}\label{negbr1}{\rm Let $(\RR, \m)$ be a Noetherian   local
ring and let $U\subset \m\RR^n$ be a parameter module. According to
Hayasaka and Hyry (\cite{HaHy10}),
\begin{eqnarray*}
%&&\psboxit{box 1.0 setgray fill}{\fbox{$
%\begin{array}{c}
%\ \\
 \br_1(U)\leq 0.
%\ \\
%\end{array}  $}}.
\end{eqnarray*}
In analogy with the case of ideals in \cite{chern3},
one asks whether $\RR$ is Cohen-Macaulay if and only if $\RR$ is
unmixed and
\begin{eqnarray*}
&&\psboxit{box 1.0 setgray fill}{\fbox{$\begin{array}{c}
\ \\
 \br_1(U)= 0.\\
\ \\
\end{array}  $}}.
\end{eqnarray*}

}\end{conjecture}

Interestingly enough, \cite[Theorem 1.1]{HaHy10} characterizes the
Cohen-Macaulayness of $\RR$ in terms of a single value of the
Buchsbaum-Rim function of $U$.

\subsection{Brian\c{c}on-Skoda theorem for
modules}\index{Brian\c{c}on-Skoda number of a module}

One application of these estimations is to provide a module version of
Theorem~\ref{bsjdegI}. It will have a surprising development. For that
we
 recall the notion of the Brian\c{c}on-Skoda number of a
module following \cite{HUV}. If
$E$ is a submodule of rank $r$ of   $\RR^r$, the {\it Brian\c{c}on-Skoda
number}
 $c(E)$ of $E$ is the smallest integer $c$ such that
$\overline{E^{n+c}}\subset F^nS_c$ for every $n$ and every reduction
$F$ of $E$. In contrast, the Brian\c{c}on-Skoda number of the ideal
$(E)$ of $\SS$
 is the smallest integer $c=c((E))$ such that
$\overline{(E)^{n+c}}\SS\subset F^n\SS$ for every $n$ and every reduction
$F$ of $E$. In particular this implies that
$\overline{E^{n+c}}\subset F^nS_c$, so that
$c(E)
\leq c((E))$. Note that   one has equality of analytic spreads $\ell(E) =
\ell((E))$.

\medskip

We shall refer to both $c(E)$ and $c((E))$ as the
Brian\c{c}on-Skoda numbers of $E$.
When $\RR$ is a regular local ring, applying the Brian\c{c}on-Skoda
theorem to $\SS$, gives $c((E))\leq \ell(E)-1$.

%\begin{example}{\rm  Let $\RR=K[x,y]$ and $E$ the submodule of the free
%$\RR$--module $\RRe_1\oplus Re_2$ generated by $x^2e_1, y^2e_2$. $E$ is a
%free $\RR$-module so $c(E)=0$. As for the ideal $(E)$ of $S=R[e_1,
%e_2]$, an application of {\sc Normaliz} (\cite{BrunsKoch}) shows that
%its normalization is $S[Et, xye_1e_2t]$; it follows that $c((E))=1$.
%}\end{example}

\medskip

Combining Theorem~\ref{bsjdegI} and
Corollary~\ref{jdegGE} one has:

\begin{Theorem} \label{bsjdegE} Let $\RR$ be a reduced quasi-unmixed
ring and
let $E$ be a module of Brian\c{c}on-Skoda number $c((E))$
and
  Buchsbaum-Rim  multiplicity $\textrm{\rm br}(E)$. Then
\[ \jdeg(\overline{\SS[(E)t]}/\SS[(E)t]) \leq c((E))\cdot  (\textrm{\rm
br}(E)+1).\]
\end{Theorem}

We now give the module version of Theorem~\ref{newbs}.

\begin{Theorem} \label{newbsmod}
Let $k$ be a perfect field, and let $(\RR, \mathfrak{m})$ be
  an integral domain which is a
 $k$-algebra essentially of finite type and dimension $d$.
  Suppose that $\RR$ has
 isolated singularities and  $L$ is the Jacobian ideal of $\RR$.
If $E\subset \mathfrak{m} \RR^r$ be a module such that $\lambda(\RR^r/E)<
\infty$,  then
\[ \jdeg(\overline{\SS[(E)t]}/\SS[(E)t])\leq (d+r+ \lambda(\RR/L) -2) \cdot
 (\br(E)+1).\]
\end{Theorem}

\begin{proof} The argument is nearly identical to that of
Theorem~\ref{newbs}. First, note that $LS$ is the Jacobian ideal of
$\SS=\RR[\TT_1, \ldots, \TT_r]$.
 Let $\DD$ be a $S_2$-ification of $\SS[(E)t]$. From
Theorem~\ref{GBS},  setting
$\CC= \overline{\SS[(E)t]}$ and  $c=\ell(E) -1=d+r-2$, we have
$L \overline{(E)^{n+c}}\subset D_n$. Observe that for all large $n$,
$\CC_{n+1}= (E)\CC_n\subset \mathfrak{m}\CC_n$.
Consider  the diagram
\begin{small}
\[
\diagram
 & D_n/D_{n+c} & & C_{n+c}/LC_{n+c} \dto|>>\tip & \\
0 \rto & (LC_{n+c}+ D_{n+c})/D_{n+c} \rto \uto|<\hole|<<\ahook
& C_{n+c}/D_{n+c} \rto
 & C_{n+c}/(LC_{n+c} + D_{n+c}) \rto &0 .
\enddiagram
\] \end{small}
Considering  that the modules in the  short exact sequence have dimension
$d+r-1$ and that the module on the right is supported in
$\mathfrak{m}$ alone, it follows that
\[ \jdeg(\CC/\DD)= \jdeg((L\CC+ \DD)/\DD) + \jdeg(\CC/L\CC+\DD) =
\jdeg((L\CC+ \DD)/\DD) + \deg(\CC/L\CC+\DD) .
\] Given the embedding on the left and the surjection on the right, we
have
that
\[ \jdeg((L\CC+\DD)/\DD)\leq \jdeg(\DD/\DD[-c])= c\cdot \jdeg(\gr(\DD)),  \]
and
\[ \deg(\CC/L\CC+\DD) \leq \deg(\CC/L\CC) \leq \lambda(\RR/L) \deg
(\CC/\mathfrak{m}\CC).
\] Since $\jdeg(\gr(\DD))= \jdeg(\gr_I(\RR))$ by Proposition~\ref{jdeggr},
 we obtain the  estimate of multiplicities
\[ \jdeg(\CC/\DD) \leq c\cdot \jdeg(\gr_I(\RR)) + \lambda(\RR/L)\cdot
f_0(\CC),
\] where $f_0(C)= \deg(\CC/\mathfrak{m}\CC)$.

Finally, we consider the exact sequence
\[ 0 \rar \mathfrak{m}C_n/C_{n+1}  \lar  C_n/C_{n+1} \lar
C_n/\mathfrak{m}C_n \rar 0.
 \]
Taking into account that
   $\gr(\CC)$ is  a ring  which has the condition $S_1$,
$ \bigoplus_n\mathfrak{m}C_n/C_{n+1}  $ either vanishes or has
 the same dimension as the
ring. Arguing as above, we have
\[ \jdeg(\gr(\CC))=
\jdeg(\bigoplus_n\mathfrak{m}C_n/C_{n+1}) + f_0(\CC).
\]
A final application of Proposition~\ref{jdeggr}, yields either
\[ \jdeg(\overline{\SS[(E)t]}/\SS(E)t])\leq (d+r+ \lambda(\RR/L) -2)
 \cdot (\br(E)+1),\] since
$ \jdeg(\gr_{(E)}(S))= \br(E)+1$. 
\end{proof}

\subsection{Complexity of the normalization of modules}

Let us apply this result to the complexity of the normalization of the
Rees algebra of a module $E\subset \RR^r$, with $\lambda(\RR^r/E) <
\infty$.
Let $\AA$ be the Rees algebra of the $\SS$-ideal $(E)$. By
Theorem~\ref{jdegintclos}, the chain of intermediate graded algebras
$\DD$,
\[\AA\subseteq \DD_0 \subset \DD_1\subset \cdots \subset \DD_s
\subset  \overline{\AA},\]
 satisfying the condition $S_2$, has  length
 bounded by $\jdeg(\overline{\AA}/\AA)$. Each one of  these algebras
 $\DD_i$
contains the subalgebra $\DD' = \DD \cap \overline{\Rees(E)}$, inducing
the chain of subalgebras
\[\Rees(E)\subset \DD_1'\subset \cdots \subset \DD_s' \subset
\overline{\Rees(E)}.\]

Therefore, if $\overline{\SS[(E)t]}$ can be built in $s$ such steps,
$\overline{\Rees(E)}$ will also be built in $s$ steps, { although
we cannot  guarantee the algebras $\DD_i'$ will satisfy the condition
$S_2$.}

The significance here comes when we compare the estimates for $s$
which come from Theorem~\ref{bsjdegE},
\[ s \leq  c((E)) \cdot \jdeg(\gr_{(E)}(\SS))=
 c((E))\cdot (\textrm{\rm br}(E)+1),\]
with
\[ s \leq {{r+c(E)-1}\choose{r}} \cdot  \textrm{\rm br}(E),\]
given by \cite[Theorem 2.3]{HUV}, using chains of subalgebras of
$\overline{\Rees(E)}$
 satisfying the condition $S_2$. Despite the inequality
 $c(E)\leq c((E))$, one is likely to obtain in most cases
   significant enhancements.

\subsection{A remark on the Zariski-Lipman conjecture}

Let $k$ be a field of characteristic zero and let $\RR$ be a finitely
generated $k$--algebra, that is a homomorphic image of a ring of
polynomials $\RR = k[x_1, \ldots, x_n]/I$.
A $k$--derivation of $\RR$ is a $k$-linear mapping $\delta: \RR
\rightarrow \RR$ that satisfies Leibniz rule,
\[\delta(ab)=a\delta(b)+b\delta(a)\] for all pairs of elements of $\RR$.
The set of all such maps is a (often non--commutative) algebra that
is a finitely generated $\RR$--module
$\mathfrak{D}= Der_k(\RR)$. The algebra and
 module structures of $\mathfrak{D}$ often code aspects of
the singularities of $\RR$.

 A more primitive object attached to $\RR$ is its module of K\"ahler
differentials, $\Omega_k(\RR)$, of which $\mathfrak{D}$ is its
$\RR$--dual, $\mathfrak{D}= \textrm{Hom}_\RR(\Omega_k(\RR),\RR)$.
More directly,
 the
structure of $\Omega_k(\RR)$ reflects many properties of $\RR$. Thus the
classical Jacobian criterion asserts that $\RR$ is a smooth algebra
over $k$ exactly when $\Omega_k(\RR)$ is a projective $\RR$--module.
For an algebra $\RR$, without nontrivial nilpotent elements, local
complete intersections are also characterized by saying that the
projective dimension of $\Omega_k(\RR)$ is at most one.
The technical issues are the comparison between the set of
polynomials that define $\RR$, represented by the ideal $I$, and the
syzygies of either $\Omega_k(\RR)$ or $\mathfrak{D}$.

The {\em Zariski--Lipman} conjecture (ZLC) makes
predictions\index{Zariski--Lipman conjecture}
about $\mathfrak{D}$, similar to those properties of $\Omega_k(\RR)$.
The most important of these questions originated with
Zariski--Lipman: {\em If $\mathfrak{D}$ is $RR$--projective then $\RR$ is
a regular ring.}
In \cite{Lipman}, the question is settled affirmatively for rings of
Krull dimension $1$, and in all dimensions the rings are shown to be
normal. Subsequently, Scheja and Storch (\cite{SchSto}) established
the conjecture  for hypersurface
rings, that is when  $\RR$ is defined by a single equation, $I = (f)$.
The last major progress on the question was the proof by Hochster
(\cite{Hochster}) of the
 graded case.

A related set of questions is collected in \cite{Vasconcelos}:
whether the finite projective dimension of either $\Omega_k(\RR)$ or
$\mathfrak{D}$ necessarily forces $\RR$ to be a local complete
intersection. (It is not known whether this is true if $\mathfrak{D}$
is projective.)
 Several lower dimension cases are known, but the most
significant progress was made by
Avramov and Herzog when they  solved the graded case (\cite{AvrHer}).

We offer a comment on a reformulation of (ZLC). For terminology, we
refer to \cite{icbook}.

\begin{Proposition} \mbox{\rm (ZLC)} holds if $\Omega_k(\RR)$ is
integrally closed.
\end{Proposition}

\begin{proof}
 Let $(\RR, \mathfrak{m})$ be the localization of an affine
$k$-algebra $\AA$ of Krull dimension $r$; $r= \dim \RR + \textrm{\rm tr.
deg. }_k(\RR/\mathfrak{m})$.
The hypohesis actually means that $\Omega_k(\RR)$,
modulo torsion, call it $E$,  is integrally closed. It will suffice to
show that $E$ is $\RR$-projective. We will argue by induction on the
Krull dimension of $R$. Consider the natural embedding
\[ E \hookrightarrow \mathfrak{D}^*= E^{**}.\] We may assume that $(\RR,
\mathfrak{m})$ is a local ring and that on the punctured spectrum of
$\RR$ the embedding is an isomorphism.

To make the argument less technical, we assume that $\RR$ is the
localization of $\AA$ at one of its maximal ideals and that $\RR$ is
regular on the punctured spectrum.

If $E$ is not contained in $\mathfrak{m}E^{**}$, $E$ will have a free
summand and therefore by an argument of Lipman-Zariski (cf.
\cite{Lipman}) the completion of $\RR$ has the form $\RR_0[[x]]$. In that
case we could reduce the dimension.

In the other case, $E$ is contained in  $\mathfrak{m}E^{**}$ and the
cokernel is a module of finite length. Since $E$ is integrally closed,
it has the property that it is $\mathfrak{m}$-full, and in particular
for any torsionfree module $F$ with $E\subset F$ and $\lambda(F/E)<
\infty$, the minimal numbers of generators satisfy $\nu(F)\leq
\nu(E)$. For the case of $\mathfrak{m}E^{**}$, this would mean
\[ \nu(\mathfrak{m})\geq  \nu(E) \geq \nu(\mathfrak{m}E^{**})=
\nu(\mathfrak{m})\cdot \textrm{\rm rank} (E), \]
an inequality that  implies $\dim \RR=1$, when (ZLC) is routine.
\end{proof}

\chapter{The Hilbert Coefficients and the Character of Local
Rings}\index{Hilbert coefficients: structure of rings}\label{chapchern}

\section*{Introduction}

Let $\RR$ be a Noetherian local ring.
Our goal here to study
what the values of the Hilbert coefficients of filtrations defined on
$\RR$ reveal about the nature of $\RR$. It is an attempt to express
in a single reading, or in a small number of samplings of these
coefficients, the overall character of $\RR$, that is whether
properties such as  being
Cohen-Macaulay, Buchsbaum, Gorenstein, regular or have finite local
cohomology are present.

\begin{itemize}
\item[{$\bullet$}] Chern coefficients of local rings
\item[{$\bullet$}] Bounding Hilbert coefficients
\item[{$\bullet$}] Euler coefficients
\end{itemize}

\section{The Chern Coefficients of Local Rings}\index{Chern
coefficient}

\subsection{Introduction}
Let $(\RR, \mathfrak{m})$  be a Noetherian
local ring of dimension $d>0$, and let $I$ be an
$\mathfrak{m}$-primary ideal. One of our goals is to study the set of
$I$-good filtrations of $\RR$. More concretely,  we will consider the
set of multiplicative, decreasing  filtrations of $\RR$
 ideals, $\mathcal{A}=\{ I_n,
 I_0=R,  I_{n+1}=I I_n, n\gg 0 \}$, integral over the $I$-adic
 filtration,
conveniently coded in
  the corresponding
Rees algebra and its associated graded ring
\[ \Rees(\mathcal{A}) = \sum_{n\geq 0} I_nt^n, \quad
\gr_{\mathcal{A}}(\RR) = \sum_{n\geq 0} I_n/I_{n+1}.
 \]

We will study certain strata of these algebras. For that we will
 focus on the role of the Hilbert polynomial of the Hilbert--Samuel
function $\lambda(\RR/I_{n+1})$,
\[ H_{\mathcal{A}}^1(n) =P_{\mathcal{A}}^1(n) = \sum_{i=0}^{d} (-1)^i
\rme_i(\mathcal{A})
{{n+d-i}\choose{d-i}}, \quad n\gg 0, \] particularly of its coefficients
$\rme_0(\mathcal{A})$ and $\rme_1(\mathcal{A})$.
These coefficients have an asymptotic character but have also a remarkable
sensitivity to the structural properties of the ring $\RR$ (see
\cite{chern3}).

\medskip

Two of our main goals are to establish relationships between the
 coefficients $\rme_i(\mathcal{A})$, for $i=0,1$, and marginally
 $\rme_2(\mathcal{A})$.  If $\RR$ is a
 Cohen-Macaulay ring, there are numerous related developments, noteworthy
 ones being given and discussed in \cite{Elias8} and \cite{RV05}.
  The situation is very distinct in the non-Cohen-Macaulay
 case. Just to illustrate the issue, suppose $d>2$ and consider a
 comparison between $\rme_0(I)$ and $\rme_1(I)$, a subject that has
 received considerable attention. It is often possible to
 pass to a reduction $\RR\rar \SS$, with    $\dim \SS=2$, or even
 $\dim \SS=1$, so
 that $\rme_0(I)= \rme_0(I\SS)$ and $\rme_1(I)=\rme_1(I\SS)$. If $\RR$ is
 Cohen-Macaulay, this is straightforward. However, in general
 the relationship between $\rme_0(I\SS)$ and $\rme_1(I\SS)$ may involve other
 invariants of $\SS$, some of which may not be easily traceable all the way to
 $\RR$.

Our perspective is partly influenced by the interpretation of the
coefficient $\rme_1$ as a {\em tracking number} introduced in
\cite{DV},  that is, as a numerical tag of the algebra
$\Rees(\mathcal{A})$ in the set of all such algebras with the same
multiplicity.
The coefficient $\rme_1$ under various circumstances is
also called the {\em Chern number or Chern coefficient} of the algebra.

\subsection{Tracking number of a filtration}

We discuss two interpretations for the coefficient
$\rme_1(\mathcal{A})$. They both require special settings.
The first, involves the {\em determinant} of the
associated graded rings with fields of representatives.

\medskip

Let $(\RR, \mathfrak{m})$ be a local ring of dimension $d>0$, with a
coefficient (infinite)  field
$k\subset \RR$, $k\simeq \RR/\mathfrak{m}$.
 For a filtration $\mathcal{A}$ as above, the associated graded ring
\[ G=\gr_{\mathcal{A}}(\RR) = \sum_{n\geq 0} I_n/I_{n+1}\]
admits a Noether normalization
\[ \SS= k[z_1, \ldots, z_d]\rar G,\]
$\deg z_i=1$. The multiplicity $\rme_0(\mathcal{A})$ is the torsionfree
rank of $G$ as an $\SS$-module. We may then define the {\em determinant}
of $G$,
\[ \det(G) := (\wedge^r(G))^{**} \simeq \SS[-\delta], \quad r=
\rme_0(\mathcal{A}).\] $\delta$ is the
tracking number of $G$ as an $\SS$-module and is denoted by $\tn_\SS(G)=\delta$.
\medskip

Applying Proposition~\ref{e1adj} gives:

\begin{proposition} \label{e1asdet} Let $G_0$ be the
torsion $\SS$-submodule of $G$. If $\dim G_0<d-1$, then
\[\rme_1(\mathcal{A}) = \tn_\SS(G).\] Otherwise,
\[\rme_1(\mathcal{A}) = \tn_\SS(G)+\deg(G_0).\]
\end{proposition}

 These equalities suggest that relationships between $\rme_0(\mathcal{A})$
and $\rme_1(\mathcal{A})$ are to be expected, a fact that has been
repeatedly borne out.
As for the  coefficient $\rme_2(\mathcal{A})$ it  can often be
understood as the degree of a {\em determinant}.

\begin{proposition}
Let $\RR$ be a Cohen-Macaulay local ring of dimension $d>0$ with a field
of coefficients as above, and $J$  a minimal reduction
for the filtration $\mathcal{A}$. Let $T=S_{\mathcal{A}/J}$ be the
corresponding Sally module \cite[Theorem 2.11]{icbook}.
 If $T\neq
0$, it is a module of dimension $d$ with the condition $S_1$ of
Serre, and therefore torsionfree as a module over  $\SS=k[z_1, \ldots,
z_d]$ of determinant
\[ \det(S_{\mathcal{A}/J}) := (\wedge^r T)^{**} \simeq \SS[-\delta], \quad r=
\rme_1(\mathcal{A})- \rme_0(\mathcal{A}) + \lambda(\RR/I_1),\] where  the degree
$\delta=\rme_2(\mathcal{A})$.

\end{proposition}

\begin{proof}
This is  another application of
Proposition~\ref{e1adj}. The interpretation of $\delta $ comes  from \cite[Theorem
2.11]{icbook}.
\end{proof}

A second interpretation of $\rme_1(\mathcal{A})$ uses the setting of normalization of ideals.
Let $(\RR, \mathfrak{m})$ be a  Noetherian local domain of
dimension $d$, which is a
quotient of a Gorenstein ring.
For an $\mathfrak{m}$-primary ideal $I$, we are going to consider the
set of all graded subalgebras $\AA$ of the integral closure of
$\RR[It]$,
\[ \RR[It] \subset \AA \subset \bar{\AA} = \bar{\RR[It]}.\]
We will assume that $\bar{\AA}$ is a finite $\RR[It]$-algebra.
We denote the set of these algebras by $\mathfrak{S}(I)$. Each of
these algebras
\[ \AA = \sum I_nt^n,\]
comes with a filtration that  is decreasing
 $I_1\supseteq I_2 \supseteq I_3 \supseteq \cdots $,  it has
an associated graded ring
\[ \gr(\AA) = \sum_{n=0}^{\infty} I_n/I_{n+1}.\]

\begin{Proposition} If $\AA$ has the condition $S_2$ of Serre, then
$\{I_n, n\geq 0\}$ is a decreasing filtration.

\end{Proposition}

\begin{proof} See \cite[Proposition 4.6]{icbook}.
\end{proof}

We are going to describe the role of the Hilbert coefficient
$\rme_1(\cdot)$ in the study of normalization of $\RR[It]$ (\cite{ni1},
\cite{jdeg1}). For each $\AA=\sum A_nt^n\in \mathfrak{S}(I)$, we may consider
the Hilbert polynomial (for $n\gg 0)$
\[ \lambda(\RR/A_{n+1}) = \sum_{i=0}^d (-1)^i\rme_i(\AA)
{{n+d-i}\choose{d-i}}.
\] The multiplicity $\rme_0(\AA)$ is constant across the whole set
$\mathfrak{S}(I)$, $\rme_0(\AA)= e_0(I)$. The next proposition shows the
role of $\rme_1(\AA)$ in tracking $\AA$ in the set $\mathfrak{S}(I)$.

\begin{Proposition} \label{e1tracking}
Let $(\RR, \mathfrak{m})$ be a normal, Noetherian local domain  which is a
quotient of a Gorenstein ring, and let $I$
be an $\mathfrak{m}$-primary ideal. For algebras $\AA$, $\BB$ of
$\mathfrak{S}(I)$:
\begin{enumerate}
\item[{\rm (a)}] If the algebras $\AA$ and $\BB$ satisfy
   $\AA\subset \BB$, then $\rme_1(\AA)\leq\rme_1(\BB)$.

\item[{\rm (b)}] If $\BB$
is the $S_2$-ification of $\AA$, then $\rme_1(\AA)=\rme_1(\BB)$.

\item[{\rm (c)}]  If the algebras $\AA$ and $\BB$ satisfy the condition
$S_2$ of Serre and
   $\AA\subset \BB$, then $\rme_1(\AA)=\rme_1(\BB)$ if and only if
   $\AA=\BB$.
   
\item[{\rm (d)}] $\rme_1(\AA)=e_1(\bar{\AA})$ if and only if $\AA$ satisfies the
condition $R_1$ of Serre.

\end{enumerate}

\end{Proposition}

\begin{proof}	The first two assertions  follow directly from the relationship between Krull
dimension and the degree of Hilbert polynomials.
 The exact sequence of graded $\RR[It]$-modules
\[ 0 \rar \AA \lar \BB \lar \BB/\AA \rar 0 \]
gives that the dimension of $\BB/\AA$ is at most $d-1$. Moreover,
$\dim \BB/\AA= d-1$ if and only if  its multiplicity is \[\deg(\BB/\AA)=
\rme_1(\BB)-\rme_1(\AA)> 0.\]

The last assertion follows because with $\AA$ and $\BB$ satisfying
$S_2$, the quotient $\BB/\AA$ is nonzero, will satisfy $S_1$ and therefore has
Krull dimension $d-1$. 
\end{proof}

\begin{Corollary} Given  a sequence of distinct algebras in
$\mathfrak{S}(I)$,
\[ \AA_0 \subset \AA_1 \subset \cdots \subset \AA_n= \bar{\RR[It]},\]
 that satisfy
the condition $S_2$ of Serre, then
\[ n \leq\rme_1(\bar{\RR[It]})-\rme_1(\AA_0)
 \leq\rme_1(\bar{\RR[It]})-\rme_1(I).
\]
\end{Corollary}

\section{Bounding Hilbert Coefficients}

\subsection{Introduction}
There are two classes of filtrations that will play an important role
in this presentation. If $(\RR, \m)$ is a Noetherian local ring of
dimension $d>0$, we focus on:
\begin{itemize}
\item[$\bullet$] Parameter ideals $I=(x_1, \ldots, x_d)$ and their $I$-adic
filtrations, and
\item[$\bullet$] Integral closure filtrations $\{\bar{I^n}\}$.
\end{itemize}

The first type will be used to examine the structure of $\RR$ itself,
while the second group is an essential element in studying the
structure of the normalization of blowups.

\medskip

The observations above
highlights the importance of having lower bounds for $\rme_1(I)$
and upper bounds for $\rme_1(\bar{\RR[It]})$. For simplicity we denote the
last coefficient as $\bar{\rme}_1(I)$. In the Cohen-Macaulay case, for
any parameter ideal $J$, $\rme_1(J)=0$. Upper bounds for $\bar{e}_1(I)$
were given in \cite{ni1}. For instance, \cite[Theorem 3.2(a),(b)]{ni1}
 show that if $\RR$ is a Cohen-Macaulay  algebra of type $t$,
essentially of finite type over a perfect field $k$
and $\delta$ is a non zerodivisor in ${\rm Jac}_k(\RR)$, then
\[
\bar{\rme}_1(I) \leq \frac{t}{t + 1}\bigl[(d-1)\rme_0(I) + \rme_0(I
+\delta \RR/\delta \RR)\bigr]\] \,\, and
\[ \bar{\rme}_1(I) \leq (d-1)\bigl[\rme_0(I)
-\lambda(\RR/\overline{I})\bigr] + \rme_0(I +\delta \RR/\delta \RR).\]

For the $I$-adic filtration, the sharpest bounds were established by
Elias (\cite{Elias8}) and Rossi-Valla (\cite{RV05}): Let $\RR$ be a Cohen-Macaulay local ring of dimension $d$, then
 \[\rme_1(I) \leq {{\rme_0(I)}\choose{2}}- {{\nu(I)-d}\choose{2}}-
 \lambda(\RR/I) +1.
\]

In \cite{chern} there were listed some questions
and conjectural
statements about the values of $\rme_1$ for very general filtrations
associated to the $\mathfrak{m}$--primary ideals of a local
Noetherian ring $(\RR, \mathfrak{m})$:

%Here is a short list of
%questions/conjectural statements that we shall address, some of higher significance
% than the others.

\begin{enumerate}
\item[{\rm (a)}] (Conjecture 1: the negativity conjecture) Let $\RR$ be an
unmixed local ring. For an ideal $J$, generated by a
system of parameters, $\rme_1(J)<0$ if and only if $\RR$ is not
Cohen-Macaulay. \index{negativity conjecture}

\item[{\rm (b)}] (Conjecture 2: the positivity conjecture) For every  $\mathfrak{m}$-primary
ideal $I$, for its integral \index{positivity conjecture}
closure filtration $\mathcal{A}$
\[ {\rme}_1(\mathcal{A}) \geq 0.\]

\item[{\rm (c)}] (Conjecture 3: the uniformity conjecture)
 For each Noetherian local ring $\RR$, there exist two functions
$\ff_l(\cdot)$, $\ff_u(\cdot)$ defined with some extended
multiplicity degree over $\RR$,
 such that for each  $\mathfrak{m}$-primary ideal
$I$ and any $I$-good filtration $\mathcal{A}$
\[ \ff_l(I) \leq {\rme}_1(\mathcal{A}) \leq \ff_u(I).\]

\item[{\rm (d)}] (Question 4: the independence question) For any two minimal reductions $J_1$, $J_2$ of  an
$\mathfrak{m}$-primary ideal $I$,
\[\rme_1(J_1)=\rme_1(J_2).\]

%\item What role the higher $\rme_i(\mathcal{A})$ play here? Likely, if
%$\gr(\mathcal{A})$ is near Cohen-Macaulay we should start seeing their
%influence.

\end{enumerate}

All of these questions were subsequently settled. Thus
 Conjecture 1
was partly settled in \cite{chern2}, and fully in \cite[Theorem
2.1]{chern3} (For another proof, see Theorem~\ref{negativityconj}.):

\begin{Theorem}
  Let $\RR$ be a Noetherian local ring with $d = \dim \RR > 0$
 and let $Q$ be a parameter ideal in $\RR$. Then following are
 equivalent:
\begin{enumerate}
\item[{\rm (a)}] $\RR$ is Cohen--Macaulay{\rm ;}
\item[{\rm (b)}] $\RR$ is {unmixed} and $\rme_1(Q)=0${\rm ;}
\item[{\rm (c)}] $\RR$ is {unmixed} and $\rme_1(Q) \geq 0$.
\end{enumerate}
\end{Theorem}

\begin{remark}{\rm.
This raises the issue  for what additional classes of rings
is $\rme_1(Q)=0$ equivalent to the Cohen-Macaulayness of $\RR$?
 \cite{chern3, chern5, chern6, chern7} contain several new questions and
developments which we will be reporting on in this Chapter.
}\end{remark}

Conjecture 2 was settled by S. Goto, J. Hong and M. Mandal (\cite[Theorem 1.1]{e1bar}).
Conjecture 3 is in a state of flux regarding the nature of its
setting.  Conjecture 4 met a prompt counterexample, but nevertheless
it has been proved the there is only a finite number of values for
$\e_1(I)$ for parameter ideals with the same integral closure
(\cite{chern7}).

\medskip

There are several other questions that will not be treated here, such
as the likely shape of the bounds for $\rme_i(I)$, $i\geq 1$, in the
non-Cohen-Macaulay case.

\bigskip

\subsubsection*{Reduction to dimension one}
Let $(\RR, \mathfrak{m})$ be a Noetherian local ring of dimension
$d>0$, $I$ an $\mathfrak{m}$-primary ideal and let $\mathcal{A}$ be a
filtration as above. We consider morphisms $\RR\rar \SS$, with $\dim
\SS=1$
and $\rme_i(\mathcal{A})=\rme_i(\mathcal{A}\SS)$, for $i=0,1$.

\medskip

If $\dim \RR =1$, let $T= \H^0_{\mathfrak{m}}(\RR)$ and $\SS=\RR/T$. For large
$n$, we have the exact sequence
\[ 0 \rar T \lar \RR/A_n \lar \SS/\SS A_n\rar 0,\]
and have the equality of Hilbert polynomials (set $\mathcal{A}'=
\SS\mathcal{A}$)
\[ \rme_0(\mathcal{A})n - \rme_1(\mathcal{A}) =
\rme_0(\mathcal{A}')n - \rme_1(\mathcal{A}') + \lambda(T),
\] which gives the equalities:
\begin{eqnarray*}
 \rme_0(\mathcal{A}) &=&  \rme_0(\mathcal{A}') \\
\rme_1(\mathcal{A}) &=& \rme_1(\mathcal{A}')-\lambda(T).
\end{eqnarray*}

In particular, if $\mathcal{A}$ is the filtration defined by the
parameter $x$, $\rme_1(x\SS)=0$, $\rme_1(x)= -\lambda(T)$.

\bigskip

If $\dim R\geq 2$, if we pass to $\SS$ as above, by
Proposition~\ref{abc},
\begin{eqnarray*}
 \rme_0(\mathcal{A}) &=&  \rme_0(\mathcal{A}') \\
\rme_1(\mathcal{A}) &=& \rme_1(\mathcal{A}').
\end{eqnarray*}

To reduce the dimension, for instance from $2$ to $1$, we need a
morphism $\RR \rar \SS$ so that the data $\H^0_{\mathfrak{m}}(\SS)$, needed
to express $\rme_1(\mathcal{A}')$, can be tracked in $\RR$. The morphism
of choice is to take $\SS= \RR/(x)$ where $x$ is a superficial element
for the filtration $\mathcal{A}$. This will guarantee that
$\rme_i(\mathcal{A})=\rme_i(\mathcal{A}')$, for $i=0,1$.

\medskip

To illustrate, consider the case of an ideal $J\subset \RR$, generated
by a system of parameters. Assume $\dim \RR \geq 2$. We choose a set of
generators $\{ x_1, \ldots, x_d\}$ so that for $\SS= \RR/(x_1, \ldots,
x_{d-1})$, $\rme_i(J)=\rme_i(J\SS)$ for $i=0,1$. Since
\[\rme_1(J\SS) = -\lambda(\H^0_{\mathfrak{m}}(\SS)),\]
we need a mechanism to estimate this value in terms of data of the
original $\RR$.
This value occurs in the expression of the big Degs of $\SS$. Thus if
we use the Samuel's multiplicity defined by $J$,
\[ \Deg_J(\SS) = \rme_0(J\SS) + \lambda(\H^0_{\mathfrak{m}}(\SS)).\]
On the other hand,
\[ \Deg_J(\SS)\leq \Deg_J(\RR),\]
so that
\[  0 \leq -\rme_1(J) = -\rme_1(J\SS) =\lambda(\H^0_{\mathfrak{m}}(\SS))\leq \Deg_J(\RR)-\rme_0(J).\]
This last expression is the Cohen-Macaulay defficiency of $\RR$
given in the {\em scale} $\Deg_I$.

\medskip

\subsection{Cohen-Macaulayness versus the negativity of $\rme_1$}

Given the role of the Hilbert coefficient $\rme_1$ as a tracking number
in the normalization of blowup algebras, it is of interest to know
its signature.

Let $(\RR ,\mathfrak{m})$ be a Noetherian local ring of dimension $d$.
If $\RR$ is Cohen-Macaulay, then $\rme_1(J)=0$ for ideal any $J$ generated by a  system of
parameter $x_1, \ldots, x_d$. As a consequence,
for any $\mathfrak{m}$-primary ideal $I$, $\rme_1(I)
\geq 0$.
\medskip

If $d=1$, the property $\rme_1(J)=0$ is characteristic of
Cohen-Macaulayness. For $d\geq 2$, the situation is somewhat
different. Consider the ring $\RR= k[x,y,z]/(z(x,y,z))$. Then for $T=
\H^0_{\mathfrak{m}}(\RR)$ and $\SS=k[x,y]= \RR/T$,  $\rme_1(\RR)=e_1(\SS)=0$.

\medskip

We are going to argue that the negativity of $\rme_1(J)$ is an
expression of the lack of Cohen-Macaulayness of $\RR$ in
numerous classes of rings. To provide a
framework, we state a conjecture posed in \cite{chern}:

\begin{Conjecture}{\rm Let $\RR$ be a Noetherian local ring that admits
an embedding  into a big Cohen-Macaulay module. Then for a parameter ideal $J$,
$\rme_1(J)< 0$ if and only if $\RR$ is not Cohen-Macaulay.
}\end{Conjecture}

Completions of unmixed local rings will satisfy the embedding condition.

\medskip

We next establish the {\em small} version of the conjecture.

\begin{Theorem} \label{e1sMCM} Let $(\RR,\mathfrak{m})$ be a Noetherian local ring of
dimension $d\geq 2$. Suppose there is an embedding
\[ 0 \rar \RR \lar E \lar C \rar 0,\] where $E$ is a finitely generated
maximal Cohen-Macaulay $\RR$-module. If $\RR$ is not Cohen-Macaulay
then for  any parameter ideal $J$, $\rme_1(J)<0$.
\end{Theorem}

\begin{proof} We may assume that the residue field of $\RR$ is infinite. We are
going to argue by induction on $d$. For $d=2$, let $J$ be a parameter
ideal. If $\RR$ is not Cohen-Macaulay, $\depth C=0$.

Let $J=(x,y)$; we may assume that $x$ is a superficial element for
the purpose of computing $\rme_1(J)$ and is also superficial relative
to $C$, that is, $x$ is not contained in any associated prime of $C$
distinct from $\mathfrak{m}$.

Tensoring the exact sequence above by $\RR/(x)$, we get the exact
complex
\[ 0 \rar T = \Tor_1^R(\RR/(x), C) \lar \RR/(x) \lar E/xE \lar C/xC \rar
0,\] where $T$ is a nonzero module of finite support. Denote
by $\SS$ the image of $\RR'=\RR/(x)$ in $E/xE $. $\SS$ is a Cohen-Macaulay ring
of dimension $1$. By the Artin-Rees Theorem, for $n\gg 0$,
$T\cap (y^n)\RR'=0$, and therefore from the diagram

\[
\diagram
0 \rto  & T \cap (y^n)\RR' \rto\dto & (y^n)\RR' \rto \dto &(y^n)\SS
\rto\dto  &0 \\
0 \rto          &  T \rto                 & R' \rto  & S \rto  &0
\enddiagram
\]
the Hilbert polynomial of the ideal $y\RR'$ is
\[ e_0 n -\rme_1 = e_0(y\SS)n + \lambda(T).\]
Thus
\begin{eqnarray}
\rme_1(J) = - \lambda(T)< 0,
\end{eqnarray}
as claimed.

\medskip

Assume now that $d\geq 3$, and let $x$ be a superficial element
for
$J$ and for the modules $E$ and $C$. In the exact sequence
\[ 0 \rar T = \Tor_1^\RR(\RR/(x), C) \lar \RR'= \RR/(x) \lar E/xE \lar C/xC \rar
0\] $T$ is either zero, and we would go on with the induction
procedure, or $T$ is a nonzero module of finite support.

If $T\neq 0$, by Proposition~\ref{abc}, we have that $\rme_1(J\RR')=\rme_1(J(\RR'/T))$, and
the embedding $\RR'/T \hookrightarrow E/xE$. By the induction
hypothesis, it suffices to prove that if $\RR'/T$ is Cohen-Macaulay
then $\RR'$, and therefore $\RR$, will be Cohen-Macaulay. This is the
content of \cite[Proposition 2.1]{HuL}. For convenience we sketch the
proof.

We may assume that $\RR$ is a complete local ring which is a quotient
$\RR= \widehat{\RR}/L$, where $\widehat{\RR}$ is a Gorenstein ring of
dimension $\dim \RR$.
Since $\RR$ is
embedded in a maximal Cohen-Macaulay module, it satisfies the
condition $S_1$ of Serre and therefore $\Ext_{\widehat{\RR}}^{d-1}(\RR,
\widehat{\RR})$ is a module of finite length. By local duality
(\cite[Theorem 3.5.8]{BH}), $\H^1_{\mathfrak{m}}(\RR)$ is a finitely
generated $\RR$-module.

 Consider the exact sequences
\[ 0 \rar T \rar \RR' \rar \SS = \RR'/T \rar 0, \]

\[ 0 \rar \RR \stackrel{x}{\lar} \RR \lar \RR' \rar 0.\]
From the first sequence,
taking local cohomology,
\[ 0 \rar \H_{\mathfrak{m}}^0(T)=T \lar \H_{\mathfrak{m}}^0(\RR')
\lar \H_{\mathfrak{m}}^0(\SS) = 0,\]
 since $\H_{\mathfrak{m}}^i(T)=0$ for $i>0$ and $\SS$ is Cohen-Macaulay of dimension
 $\geq 2$; one also has
$\H_{\mathfrak{m}}^1(\SS)=\H_{\mathfrak{m}}^1(\RR')=0$.

From the second sequence we obtain the exact sequence
\[ 0 \rar \H^1_{\mathfrak{m}}(\RR)\stackrel{x}{\lar}
\H^1_{\mathfrak{m}}(\RR)\lar \H^1_{\mathfrak{m}}(\RR')=0,
\] and therefore since $\H^1_{\mathfrak{m}}(\RR)$ is finitely generated,
%From the second
%sequence,
%since $\RR$ satisfies the condition $S_1$ of Serre, $$
%the associated primes of
%$\RR$ have dimension $d$, $H_{\mathfrak{m}}^1(\RR)$ is a finitely generated $\RR$-module,
 by Nakayama Lemma, $\H_{\mathfrak{m}}^1(\RR)=0$, and
 $T=\H_{\mathfrak{m}}^0(\RR')=0$. 
 \end{proof}

We now analyze what is required to extend the proof to big
Cohen-Macaulay cases. We are going to assume that $\RR$ is an integral
domain and that $E$ is a big balanced Cohen-Macaulay module (see
\cite[Chapter 8]{BH}, \cite{Sharp81}).
Embed $\RR$ into $E$,
\[ 0 \rar \RR \lar E \lar C \rar 0.\]
The argument above ($d\geq 3)$ will work if in the induction argument we can pick
$x\in J$  superficial for the Hilbert polynomial of $J$ and avoids
the finite set of associated primes of $E$ and all associated primes
of $C$ different from $\mathfrak{m}$. It is this last condition that
is the most troublesome.

There is one case when this can be overcome, to wit, when $\RR$ is a
complete local ring and $E$ is  countably generated. Indeed, $C$ will
be countably generated and $\Ass(C)$ will be a countable set. The
prime avoidance result of
\cite[Lemma 3]{Burch} allows for the choice of $x$. Let us apply
these ideas in an important case.

\begin{Theorem} \label{e1geoMCM} Let $(\RR,\mathfrak{m})$ be a
Noetherian local integral domain essentially of finite type over a
field.
 If $\RR$ is not Cohen-Macaulay
then for  any parameter ideal $J$, $\rme_1(J)<0$.
\end{Theorem}

\begin{proof} Let $\AA$ be the integral closure of $\RR$ and $\widehat{\RR}$ its
completion. Tensor the embedding $\RR\subset \AA$ to obtain
\[ 0 \rar \widehat{\RR} \lar \widehat{\RR} \otimes_\RR \AA=\widehat{\AA}. \]
From the properties of pseudo-geometric local rings (\cite[Section
37]{Nagata}), $\widehat{\AA}$ is a reduced semi-local ring with a
decomposition
\[ \widehat{\AA} = \AA_1 \times \cdots \times \AA_r, \]
where each $\AA_i$ is a complete local domain, of dimension $\dim \RR$
and finite over $\widehat{\AA}$.

For each $\AA_i$ we make use of \cite{Griffith76} and
\cite{Griffith78}
and pick a countably generated big balanced Cohen-Macaulay
$\AA_i$-module and therefore $\widehat{\RR}$-module. Collecting the $\rme_i$
we have an
embedding
\[ 0 \rar \widehat{\RR} \lar \AA_1 \times \cdots \times \AA_r\lar  E= E_1
\oplus \cdots \oplus E_r.\]
As $E$ is a countably generated big balanced Cohen-Macaulay
$\widehat{\RR}$-module, the argument above shows that if $\widehat{\RR}$
is not Cohen-Macaulay then $\rme_1(J\widehat{\RR})<0$. This suffices to
prove the assertion about $\RR$. 
\end{proof}

\begin{Remark}{\em
There are other classes of local rings admitting big balanced Cohen-Macaulay
modules. A crude  to handle it would be: Let $E$
be such a module and assume it is generated by $s$ elements. Let $X$
be a set of indeterminates of cardinality larger than $s$.
The local ring $\SS=\RR[X]_{\mathfrak{m}[X]}$ is $\RR$-flat and has the
same depth as $\RR$. If $\SS \otimes_\RR E$ is a big balanced Cohen-Macaulay
$\SS$-module, with its residue field of cardinality larger than that of
a set of generators
of $C$, prime avoidance would again work. It is
not known however where $\SS\otimes_\RR E$ is balanced.
}\end{Remark}

\begin{Example}{\rm We will consider some classes of examples.

\medskip

(i) Let $(\RR, \mathfrak{m})$ be a regular local ring and let $F$ be a
nonzero (finitely generated) free $\RR$-module. For any non-free submodule of $F$,
the idealization (trivial extension) of $\RR$ by $N$, $\SS=\RR\oplus N$ is
a non Cohen-Macaulay local ring. Picking $E=\RR\oplus F$,
Theorem~\ref{e1geoMCM}
 implies that for any parameter ideal $J\subset \SS$, $\rme_1(J)< 0$.
It is not difficult to give an explicit formula for $\rme_1(J)$ in this
case.

\medskip

(ii) Let $\RR= \mathbb{R}+ (x,y)\mathbb{C}[x,y]\subset
\mathbb{C}[x,y]$, for $x,y$ distinct indeterminates. $\RR$ is not
Cohen-Macaulay but its localization at the maximal irrelevant ideal
is a Buchsbaum ring. It is easy to verify that
$\rme_1(x,y)=-1$ and that $\rme_1(\SS)=0$ for the $\mathfrak{m}$-adic
filtration of $\SS$.

Note that $\RR$ has an isolated singularity. For these rings,
\cite[Theorem 5]{jdeg1} can be extended (does not require the
Cohen-Macaulay condition), and therefore describes bounds  for
$\rme_1(\bar{\AA})$ of integral closures. Thus, if $\AA$ is the Rees of
the parameter ideal $J$, one has
\[\rme_1(\bar{\AA})-e_1(J) \leq (d-1+ \lambda(\RR/L))e_0(J),  \]
where $L$ is the Jacobian of $\RR$.

In this example,  one has $d=2$, $\lambda(\RR/L)=1$,
\[\rme_1(\bar{\AA})-e_1(J) \leq 2 e_0(J).  \]

\medskip

(iii) Let $k$ be a field of characteristic zero and let
$f=x^3+y^3+z^3$ be a polynomial of $k[x,y,z]$. Set $\AA=k[x,y,z]/(f)$
and let $\RR$ be the Rees algebra of the maximal irrelevant ideal
$\mathfrak{m}$ of
$\AA$. Using the Jacobian criterion,  $\RR$ is
normal. Because the reduction number of $\mathfrak{2}$ is $2$, $\RR$ is
not Cohen-Macaulay.
 Furthermore, it is easy to verify that $\RR$ is not contained in any
 Cohen-Macaulay domain that is finite over $\RR$.
Let
 $\SS=\RR_{\mathcal{M}}$, where $\mathcal{M}$ is the irrelevant maximal
 ideal of $\RR$.
The first superficial element (in the reduction to dimension two)
can be chosen to be prime. Now one takes the integral which will be a
maximal Cohen-Macaulay module.

The argument extends to geometric  domains in any characteristic
 if $\depth \RR=d-1$.

}\end{Example}

\begin{Remark}{\rm Uniform lower bounds for $\rme_1$ are rare but still
exist in special cases. For example, if $\RR$ is a generalized Cohen-Macaulay
ring, then according to \cite[Theorem 5.4]{GoNi03},
\begin{eqnarray} \label{e1boundgcm}
\rme_1(J) &\geq &- \sum_{i=1}^{d-1}{d-2\choose
i-1}\lambda(\H^i_{\mathfrak{m}}(\RR)),
\end{eqnarray}
with equality if $\RR$ is Buchsbaum.

It should be observed that uniform lower bounds may not always exist.
For instance, if $\AA=k[x,y,z]$, and $\RR  $ is the idealization of $(x,y)$,
then for the ideal $J=(x,y,z^n)$, $\rme_1(J)=-n$.

The Koszul homology modules $\H_i(J)$ of $J$ is a first place where to look for
bounds for $\rme_1(J)$. We recall \cite[Theorem 4.6.10]{BH}, that the
multiplicity of $J$ is given by the formula
\[ \e_0(J)= \lambda(\RR/J) - \sum_{i=1}^d(-1)^{i-1} h_i(J),
\] where $h_i(J)$ is the length of $\H_i(J)$. The summation term is
non-negative and only vanishes if $\RR$ is Cohen-Macaulay. Unfortunately
it does not gives  bounds for $\rme_1(J)$. There is a formula involving
these terms in the special case when $J$ is generated by a
$d$-sequence. Then the corresponding {\em approximation complex} is
acyclic, and the Hilbert-Poincar\'e series of $J$ (\cite[Corollary
4.6]{HSV3}) is
\[ \frac{\sum_{i=0}^d (-1)^i h_i(J)t^i}{(1-t)^d},
\] and therefore
\[\rme_1(J) =  \sum_{i=1}^d (-1)^i i h_i(J).
\]

Later in the notes, we shall prove the existence of lower bounds more
generally, by making use of extended degree functions.

}\end{Remark}

%\subsection{Bounds}

\subsubsection{Cohen-Macaulayness versus the Euler number}

We develop a more abstract approach to the
relationship between the signature of $\rme_1$ and the
Cohen-Macaulayness of a local ring. We also give a more
self-contained proof of the main result that avoids the use of big
Cohen-Macaulay modules.

\medskip

Let $(\RR, \mathfrak{m})$ be a Noetherian local ring of dimension
$d\geq 2$. The enabling idea above was the embedding of $\RR$ into a
Cohen-Macaulay (possibly big Cohen-Macaulay) module $E$,
\[ 0 \rar \RR \lar E \lar C \rar 0.\]
Unfortunately, it may not be always possible to find the appropriate
$\RR$-module $E$. Instead we will seek embed $\RR$ into a Cohen-Macaulay
module $E$ over a ring $\SS$ closely related to $\RR$ for the purpose of
computing associated graded rings of adic-filtrations.

\begin{Theorem}$($\cite[Theorem 3.3]{chern2}$)$ \label{e1dMCM}
Let $(\RR, \mathfrak{m})$ be a Noetherian local domain  of dimension
$d\geq 2$, which is a homomorphic image of a Cohen-Macaulay local
domain $\SS$, let  $\RR = \SS/\mathfrak{p}$. Let $J$ be an $\RR$-ideal generated
by a system of parameters. If $\RR$ is not Cohen-Macaulay then $\rme_1(J)<
0$.
\end{Theorem}

\begin{proof} If $\height \mathfrak{p}\geq 1$, we first replace $\SS$ by
$\SS'=\SS/L$, where $L$ is $\SS$-ideal generated by
a maximal regular sequence in $\mathfrak{p}$.
This means that we may assume that $\dim \RR=\dim \SS$, and that
$\mathfrak{p}$ is a minimal prime of $\SS$. In particular, we have an
embedding of $\SS$-modules
\[ 0 \rar \RR=\SS/\pp \lar E= \SS \lar C \rar 0.\]

At this point we may assume that $\RR$ has an infinite residue field.
Let $J= (x_1, \ldots, x_d)$. We claim that there exists an $\SS$-ideal
$I= (y_1, \ldots, y_d)$ generated by a system of
parameters that lifts $J$, more precisely with $x_i= y_i+
\mathfrak{p}$.
 To choose $y_1$, let $\mathfrak{p}_1, \ldots,
\mathfrak{p}_n$ be the minimal primes of $\SS$ that do not contain
$\mathfrak{p}$. Let $z_1$ be any lift of $x_1$ and let $\mathfrak{p}=
(a_1, \ldots, a_s)$.
For each $\lambda\in R\setminus \mathfrak{m}$, consider the element
\[ z_1 + \lambda a_1 + \cdots + \lambda^s a_s.\]
If all of these elements belong in $\bigcup_{k=1}^n \mathfrak{p}_k$,
as
  $\RR/\mathfrak{m}$ is infinite, there would be $\lambda_1, \ldots,
  \lambda_{s+1}$ with different images in $\RR/\mathfrak{m}$ such that
\[ z_1 + \lambda_i a_1 + \cdots + {\lambda_i}^s a_s\in
\mathfrak{p}_k, \quad i=1, \ldots, s+1,\] contained in a same
$\mathfrak{p}_k$. Since the Vandermonde matrix, determined by
$\lambda_i, 1\leq i \leq s+1$,  is invertible,  all the $a_i$
would be contained in the same prime $\mathfrak{p}_k$, against the
assumption.

\medskip

Now reduce $\RR$ and $\SS$ modulo $(y_1)$. Note that
\[\dim \RR/(y_1)=\dim
\RR/(x_1)=\dim \SS/(y_1)=d-1,\] but that $\RR'=\RR/(x_1)$ may not be embedded in
$\SS'=\SS/(y_1)$. To continue the lifting process, let $\mathfrak{p}_1,
\ldots, \mathfrak{p}_n $ now denote the minimal primes of $(y_1)$
which do not contain $\mathfrak{p}$; these primes have codimension
$1$. We proceed as above to select
$y_2$. In this fashion we arrive at the sequence $y_1, \ldots, y_d$, that
clearly forms a system of parameters of $\SS$.

If we denote $I=(y_1, \ldots, y_d)\SS$, we have $I\RR =J$ and
\[ \gr_I(\RR)= \gr_J(\RR).\]
Therefore, for the purpose of determining   $\rme_1(J)$, we may treat $\RR$ as
an $\SS$-module and use the $I$-adic filtration. Now proceed as in the
proof of Theorem~\ref{e1sMCM}. Let $x$ be a superficial element for
the purpose of determining the Hilbert function of $\gr_I(\RR)$ and
which is not contained in any associated prime of $C$ distinct from
$\mathfrak{m}$. Reduction modulo $(x)$ gives rise to the exact
sequence
\[ 0 \rar T = \Tor_1^\SS(\SS/(x), C) \lar  \RR/(x\RR) \lar \SS/(x) \lar C/xC \rar
 0,\] where $T$ is a  module of finite support. Denote
by $\RR'$ the image of $\RR/(x)$ in $\SS'=\SS/(x) $. $\SS'$ is a Cohen-Macaulay ring
of dimension $d-1$.

 If $d=2$,  $T\neq 0$, as in the proof of Theorem~\ref{e1sMCM}, we have
that
\[\rme_1(J)= -\lambda(T),\] and therefore
$\rme_1(J)$ is always non-positive.
If $d>  2$, setting $J'= J\RR'$, we have $\rme_1(J')=e_1(J)$.
If $\RR'$ is not Cohen-Macaulay, by induction $\rme_1(J')< 0$. To complete
the proof we argue that if $\RR'$ is Cohen-Macaulay, then $\RR$ is also
Cohen-Macaulay.
 We are going to use the argument of
\cite[Proposition 2.1]{HuL}.

Once $\RR$ is embedded in $\SS$, we may assume that $\SS$ is a complete
local ring.
Let us consider the exact sequences
\[ 0 \rar \RR \stackrel{x}{\lar} \RR \lar \RR/(x) \rar 0 \]
and
\[ 0 \rar T \lar \RR/(x) \lar \RR' \rar 0.\]
Taking local cohomology, we have the exact sequence
\[ 0 \rar \H^0_{\mathfrak{m}}(\RR/(x)) \lar
\H^1_{\mathfrak{m}}(\RR) \stackrel{x}{\lar}
\H^1_{\mathfrak{m}}(\RR) \lar \H^1_{\mathfrak{m}}(\RR/(x)),
\] where $\H^1_{\mathfrak{m}}(\RR/(x))$ vanishes since $\RR'$ is
Cohen-Macaulay of dimension $d-1\geq 2$, which gives
$\H_{\mathfrak{m}}^1(\RR)=x\H_{\mathfrak{m}}^1(\RR)$.

\begin{lemma} \label{finiteH1} Let $(\SS, \mathfrak{n})$ be a Cohen-Macaulay complete
local ring of dimension $d\geq 2$, and let $M$ a finitely generated
$\SS$-module with the condition $S_1$ of Serre. Then
$\H_{\mathfrak{n}}^1(M)$ is a module of finite length.
\end{lemma}

\begin{proof}  If $\omega$ is the canonical module of $\SS$,
$\Ext_\SS^{d-1}(M,\omega)$ vanishes at each localization
$\SS_{\mathfrak{p}}$, $\mathfrak{p}\neq \mathfrak{n}$, since $\depth
M_{\mathfrak{p}}\geq 1$.
  By local duality (see \cite[Theorem 3.5.8]{BH}),
 $\H^1_{\mathfrak{m}}(\RR)$ is a finitely generated $\RR$-module.
 \end{proof}

To complete the proof of  Theorem\ref{e1dMCM},
by
  Nakayama Lemma $\H^1_{\mathfrak{m}}(\RR)=0$, and $T$ vanishes as well.
Thus $\RR/(x)=\RR'$, and $\RR$ will be Cohen-Macaulay. 
\end{proof}

%since \H^1_{\mathfrak{m}}(\RR') = \H^1_{\mathfrak{m}}(S/()

\bigskip

To extend the result from integral domains to more general rings,
$\RR=\SS/L$, $\dim \RR=\dim \SS$, where $\SS$ is Cohen-Macaulay, we
discuss how to  embed $\RR$
into a Cohen-Macaulay  $\SS$-module. At a minimum, we need to require that $\RR$
be unmixed.

\medskip

We make  an elementary observation of what it takes to embed
  $\RR$ into a free $\SS$-module.

\begin{proposition} Let $\SS$ be a Noetherian  ring and $L$ an
ideal of codimension zero without embedded components. If $\RR=\SS/L$,
there is an embedding $\RR\rar F$ into a free
$\SS$--module $F$ if and only if $L=0:(0:L)$. In particular, this
condition always holds if the total ring of fractions of  $\SS$ is a Gorenstein ring.
\end{proposition}

\begin{proof} Let $\{a_1, \ldots, a_n\}$ be a generating set of  $0:L $, and
consider the   mapping $\varphi: \SS
\rar \SS^n$, $\varphi(1)=(a_1, \ldots, a_n)$; its kernel
is isomorphic to $0:(0:L)$. This shows that the equality $L=0:(0:L)$
is required for the asserted embedding.

Conversely, given an embedding $\varphi:\SS/L \rar \SS^n$, let $(a_1,
\ldots, a_n)\in \SS^n$ be the imagine of a generator of $\SS/L$. The
ideal $I$ these entries generate  is annihilated by $L$, $I\subset
0:L$. Since $ 0:I=L$, we have $0:(0:L)\subset 0:I= L$.

If the total ring of fractions of $\SS$ is Gorenstein, to prove that
 $0:(0:L)\subset  L$, it suffices to localize at the associated
 primes of $L$, all of which have codimension zero and with a
 localization which is Gorenstein. But the double annihilator
 property is characteristic of such rings. 
 \end{proof}

We are now in a position to give another proof of \cite[Theorem
2.1]{chern3}:

\begin{Theorem}\label{negativityconj} Let $(\RR, \m)$ be an unmixed
 Noetherian local
ring. Then $\RR$ is  Cohen-Macaulay if and only if for some parameter
ideal $Q$, $\rme_1(Q)=0$.
\end{Theorem}

\begin{proof}
Pass to the completion of $\RR$, which leaves unchanged the
assumptions and assertions. We can now write $\RR=\SS/L$, where $\SS$
is a Gorenstein local ring of the same dimension as $\RR$. By
Proposition
~\ref{embedinfree}, $\RR\hookrightarrow \SS^n$, for some integer $n$.
Now the argument of Theorem~\ref{e1sMCM}, combined
with the lifting argument of
Theorem~\ref{e1dMCM}, gives the desired assertions. 
\end{proof}

\begin{remark}{\rm Mandal and Verma (\cite{MV}) clarify this further:
For all Noetherian local rings,
$\rme_1(Q)$ is never positive.
Similar assertions follows for modules, according to \cite{MSV10}.
See also Corollary~\ref{e1notpos}.
}\end{remark}

Another of the questions raised about $\rme_1$ concerned whether the
equality $\rme_1(J_1)=e_1(J_2)$ always hold for minimal reductions of
the same ideal. According to \cite{chern3} this is not always the
case. We give now a family of
examples based on a method of
\cite{GoNi03}.

\begin{example}{\rm Let $(\SS,\mathfrak{m})$ be a regular local ring of
dimension four, with an infinite residue field. Let $P_1, \ldots,
P_r$ be a family of codimension two Cohen-Macaulay such that
$P_i+P_j$ is $\mathfrak{m}$-primary for $i\neq j$.
    Define $\RR=\SS/\bigcap_iP_i$.

 Consider the exact
sequence of $\SS$-modules
\[ 0 \rar \RR \lar \bigoplus_i \SS/P_i \lar L \rar 0.\] Note that $L$ is a
module of finite support; it may be identified to
$\H_{\mathfrak{m}}^1(\RR)$.
Let $J=(a,b)$ be an ideal of $\RR$  forming a system
of parameters, contained in the annihilator of $L$.\footnote{We thank
Jugal Verma for this observation.} We can assume that $a,b\in \SS$ form a regular
sequence in each $\SS/P_i$.
We are going to determine $\rme_1(J).$
For each integer $n$,
tensoring by $\SS/(a,b)^n$ we get the exact sequence
\begin{small}
\[ 0\rar \Tor_1^\SS(L, \SS/(a,b)^n) \rar \RR/(a,b)^n \rar \bigoplus_i
\SS/(P_i, (a,b)^n)
 \rar L \otimes_\SS \SS/(a,b)^n \rar 0.\]
\end{small}
For $n\gg 0, $ $(a,b)^nL=0$, so we have
 that  $L \otimes_\SS \SS/(a,b)^n = L$ and
$ \Tor_1^\SS(L, \SS/(a,b)^n)= L^{n+1}$, from
 the Burch-Hilbert $(n+1)\times n$ resolution of the ideal
$(a,b)^n$. Since the  $\RR_i=\SS/P_i$  are
Cohen-Macaulay, we obtain the following Hilbert-Samuel polynomial:
\[\rme_0(J){{n+2}\choose {2}} -\rme_1(J) {{n+1}\choose{1}} + \rme_2(J) =
(\sum_{i=1}^r \rme_0(J\RR_i))
 {{n+2}\choose{2}}  +
(n+1)\lambda(L)-\lambda(L).\]
 It gives
\begin{eqnarray*}
\rme_0(J\RR) &=& \sum_{i=1}^r \rme_0(J\RR_i) , \\
\rme_1(J\RR) & = & -\lambda(L),\\
\rme_2(J\RR) &=& -\lambda(L).
\end{eqnarray*}

}
\end{example}

\subsubsection*{Filtered modules}

The same relationship discussed above between the signature $\rme_1(J)$
and the Cohen-Macaulayness of $\RR$ holds true for  modules.
More precisely, one can extend the proof of Theorem~\ref{e1dMCM} in
the following manner.

\begin{Theorem}[{\cite[Theorem 5.1]{chern2}}] \label{e1formodules}
Let $(\RR, \mathfrak{m})$ be a Noetherian local ring of dimension
$d\geq 1$
 and let $M$ be a finitely generated module embedded in a
MCM module $E$. If $J$ is an ideal generated by a system of
parameters, then $M$ is Cohen-Macaulay if and only if
$\rme_1(\gr_J(M))\geq 0$.
\end{Theorem}

Let $\RR=k[x_1, \ldots, x_d]$ be a ring of polynomials over the field
$k$, and let $M$ be a finitely generated graded $\RR$-module. Suppose
$\dim M=d$. For $J=(x_1, \ldots, x_d)$ we can apply
Theorem~\ref{e1formodules} to $M$ in a manner that uses the Hilbert
function information of the natural grading of $M$.

\begin{Theorem} \label{e1vsCM0}
Let $\RR=k[\xx]=k[x_1, \ldots, x_d]$, $d\geq 2$,  be a ring of polynomials over the field
$k$, and let $M$ be a finitely generated graded $\RR$-module generated
in degree $0$. If $M$ is torsion free, then $M$
 is a free $\RR$-module if and only if $\rme_1(M)=0$.
\end{Theorem}

\begin{proof} Since $M$ is generated in degree $0$, $M \simeq \gr_{\xx}(M)$.
By assumption, $M$ can be embedded in a free $\RR$-module $E$ (not
necessarily by a homogeneous homomorphism).  Now we apply
Theorem~\ref{e1formodules}. 
\end{proof}

If $M$ is generated in, say, degree $a>0$, we have the equality,
\[ \lambda(M/(\xx)^{n+1}M) = \sum_{k=0}^n \lambda(M_{a+k}),\]
so the Hilbert coefficients satisfy
\begin{eqnarray*}
\rme_0(\gr_{\xx}(M)) &=& \rme_0(M[-a])= \rme_0(M),\\
\rme_1(\gr_{\xx}(M)) &=&\rme_1(M[-a])=\rme_1(M)-a\rme_0(M).
\end{eqnarray*}

This leads  to:

\begin{corollary} \label{e1vsCMa} Let
$\RR=k[\xx]=k[x_1, \ldots, x_d]$, $d\geq 2$,  be a ring of polynomials over the field
$k$, and let $M$ be a finitely generated graded $\RR$-module generated
in degree $a\geq 0$. If $M$ is torsion free, then $\rme_1(M)\leq
a\rme_0(M)$, with equality if and only if
$M$
 is a free $\RR$-module.
\end{corollary}

\begin{question}{\rm
Let $\RR$ be a complete Noetherian local ring of dimension $d\geq 1$ and let
$J$ be an ideal generated by a system of parameters. If all the
associated prime ideals of $M$ have dimension $d$, does $\rme_1(\gr_J(M))\geq
0$ implies that $M$ is Cohen-Macaulay?
}\end{question}

\subsubsection*{Embedding into big balanced Cohen-Macaulay
modules}\index{big balanced CM module}

Let $(\RR, \mathfrak{m})$ be a Noetherian local domain. If $\RR$ has a
big Cohen-Macaulay module $E$, any $0\neq e\in E$ allow for an
embedding $\RR\hookrightarrow E$. In fact, one may assume that $E$ is a
 balanced big Cohen-Macaulay (bbCM for short)
module (see \cite[Section 8.5]{BH} for a discussion).
According to the results of Hochster, if $\RR$
 contains a field, then there is a bbCM
$\RR$-module $E$.

To use the argument in \cite[Theorem 3.2]{chern}, in the
exact sequence
\[ 0 \rar \RR \lar E \lar C\rar 0,\]
we should, given any parameter ideal $J$ of $\RR$ pick an element
superficial for the purpose of building $\gr_J(\RR)$ (if $\dim \RR>2$)
and not contained in any associated prime of $C$ different from
$\mathfrak{m}$.
This is possible if the cardinality of the residue field is larger
than the cardinality of $\Ass (C)$.

\begin{Theorem} \label{e1withBBCM} Let $(\RR,\mathfrak{m})$ be a
Noetherian local integral domain that is not Cohen-Macaulay and let $E$ be a balanced big
Cohen-Macaulay module. If the residue field of $\RR$ has cardinality
larger than the cardinality of a generating set for $E$, then
  $\rme_1(J)<0$
 for  any parameter ideal $J$.
\end{Theorem}

Let
 $\mathbf{X}$ be a
set of indeterminates of larger cardinality than $\Ass(C)$, and
consider $\RR(\XX) = \RR[\XX]_{\mathfrak{m}\RR[\XX]}$. This is a Noetherian
ring (\cite{GH79}), and we are going to argue that if $E$ is a bbCM
$\RR$-module, then
  $\RR(\XX)\otimes_\RR E$ is a bbCM module over $\RR(\XX)$.
S. Zarzuela has kindly pointed out to us the following result:

\begin{Theorem}[{\cite[Theorem 2.3]{Zarzuela}}]
Let $\AA \rightarrow \BB$ a flat morphism of local rings $(\AA,
\mathfrak{m}),  (\BB, \mathfrak{n})$
and $M$ a balanced big Cohen-Macaulay $A$-module. Then, $M
\otimes_\AA \BB$
is a balanced big Cohen-Macaulay $\BB$-module if and only if the following
two conditions hold:
\begin{itemize}

\item[{\rm (a)}]  $\mathfrak{n}(M\otimes_\AA \BB) \neq  M \otimes_\AA
\BB$ and

\item[{\rm (b)}]  For any prime ideal $\mathfrak{q} \in
\supp_\BB(M\otimes_\AA \BB)$,  $(1)$
$\height(\mathfrak{q} / \mathfrak{p}\BB) = \depth
(C_{\overline{\mathfrak{q}}})$  and $(2)$ $\height(\mathfrak{q}) +
\dim(\BB/\mathfrak{q}) = \dim(\BB)$.

\end{itemize}
Here, we denote by $\supp_\AA(M)$ $($small support$)$ the set of prime ideals
in $A$ with at least one non-zero Bass number in the $\AA$-minimal
injective resolution of $M$, $\mathfrak{p} =
\mathfrak{q}\cap \AA,$  $C = \BB/\mathfrak{p}\BB$ and
$\overline{\mathfrak{q}}=\mathfrak{q}C$.

Moreover, if $\mathfrak{q} \in \supp_\BB(M\otimes _\AA \BB)$ then
$\mathfrak{p} \in \supp_\AA(M)$, and
$\height(\AA/\mathfrak{p}) + \dim(\AA/\mathfrak{p}) = \dim(\AA)$.
\end{Theorem}

\begin{corollary} Let $(\RR, \mathfrak{m})$ be a universally catenary
integral domain and let $E$ be a bbCM $\RR$-module. For any set $\XX$
of indeterminates and $\BB=\AA(\XX)= \RR[\XX]_{\mathfrak{m}[\XX]}$,
$\BB\otimes_\AA E$ is a bbCM $\BB$-module.
\end{corollary}

\begin{Theorem} \label{e1Hochster} Let $(\RR, \mathfrak{m})$ be a universally catenary
integral domain containing a field. If $\RR$ is not Cohen-Macaulay,
then $\rme_1(J)<0$ for any parameter ideal $J$.
\end{Theorem}

\subsection{Bounds on $\overline{e}_1(I)$ via the Brian\c{c}on--Skoda number}

\medskip

This is a slight generalization of \cite[Theorem 3.1]{ni1} and
\cite[Theorem 5]{jdeg1}.

\bigskip

We discuss the role of Brian\c{c}on--Skoda type theorems (see
\cite{AHu}, \cite{LipmanSathaye}) in determining some
relationships between the coefficients $\rme_0(I)$ and
$\overline{\rme}_1(I)$. We follow the treatment given in \cite[Theorem
3.1]{ni1}, but formulated for the non Cohen-Macaulay case.
We will use a Brian\c{c}on--Skoda theorem
that works in non--regular rings. We are going to provide a short
proof along the lines of \cite{LipmanSathaye} for the special case
we need: $\mathfrak{m}$--primary ideals in a local Cohen--Macaulay
ring. The general case is treated by Hochster and Huneke in
\cite[1.5.5 and 4.1.5]{HH}. Let $k$ be a perfect field, let $\RR$ be
a reduced and equidimensional  $k$--algebra essentially of finite
type, and assume that $\RR$ is affine with $d={\rm dim} \, \RR$ or
$(\RR, \mathfrak{m})$ is local with $d={\rm dim}\, \RR + {\rm trdeg}_k
R/{\m}$. Recall that the {\it Jacobian ideal}   ${\rm Jac}_k(\RR)$
of $\RR$ is defined as the $d$--th Fitting ideal of the module of
differentials $\Omega_k(\RR)$ -- it can be computed explicitly from
a presentation of the algebra. By varying Noether normalizations
one deduces from \cite[Theorem 2]{LipmanSathaye} that the Jacobian
ideal ${\rm Jac}_k(\RR)$ is contained in the conductor $\RR \colon
\overline{\RR}$ of $\RR$ (see also \cite{Noether}, \cite[3.1]{AB} and
\cite[2.1]{H}); here $\overline{\RR}$ denotes the integral closure
of $\RR$ in its total ring of fractions.

\begin{Theorem}\label{GBS}\index{Brian\c{c}on-Skoda theorem}
Let $k$ be a perfect field, let $\RR$ be a reduced local
 $k$--algebra essentially of finite type, and let
$I$ be an  ideal with a minimal reduction generated by $g$ elements.
Denote by $\DD = \sum_{n\geq 0} D_nt^n$ the $S_2$-ification of $\RR[
It]$.
 Then for every
integer $n$,
\[
{\rm Jac}_k(\RR) \, \overline{I^{n+g -1}} \subset D_n.
\]
\end{Theorem}
\begin{proof} We may assume that $k$ is infinite. Then, passing to a
minimal reduction, we may suppose that $I$ is generated by  $g$
generators. Let $\SS$ be a finitely generated
$k$--subalgebra of $\RR$ so that $\RR=\SS_{\p}$ for some $\p \in {\rm
Spec}(\SS)$, and write $\SS=k[x_1, \ldots, x_e] = k[X_1, \ldots,
X_e]/{\mathfrak a}$ with ${\mathfrak a}=(h_1, \ldots, h_t)$ an ideal of
height $c$. Notice that $\SS$ is reduced and equidimensional. Let
$K= (f_1, \ldots, f_g)$ be an $\SS$--ideal with $K_{\p}=I$, and
consider the extended Rees ring $\BB=\SS[Kt, t^{-1}]$. Now $\BB$ is a
reduced and equidimensional affine $k$--algebra of dimension
$e-c+1$.

Let $\varphi \colon k[X_1, \ldots, X_e, \TT_1, \ldots, \TT_g, U]
\twoheadrightarrow \BB$ be the $k$--epimorphism mapping $X_i$ to
$x_i$, $\TT_i$ to $f_it$ and $U$ to $t^{-1}$. Its kernel has height
$c+g$ and contains the ideal ${\mathfrak b}$ generated by $\{ h_i,
\TT_jU-f_j | 1 \leq i \leq t, 1 \leq j \leq g \}$. Consider the
Jacobian matrix of these generators,
\[
\Theta = \left(
\begin{array}{c|cccc}
\displaystyle\frac{\partial h_i}{\partial X_j} & &   \!\!\!\! 0 & & \\
\hline
& U & & & \TT_1 \\
* & & \ddots & & \vdots \\
& & & U & \TT_g
\end{array}
\right).
\]
Notice that $ I_{c+g}(\Theta)\supset I_c(\left(\frac{\partial
h_i}{\partial X_j}\right))U^{g-1}(\TT_1, \ldots, \TT_g)$. Applying
$\varphi$ we obtain $ {\rm Jac}_k(\BB) \supset I_{c+g}(\Theta) \BB
\supset {\rm Jac}_k(\SS) K t^{-g+2}$. Thus ${\rm Jac}_k(\SS) K
t^{-g+2}$ is contained in the conductor of $\BB$. Localizing at $\p$
we see that ${\rm Jac}_k(\RR) I t^{-g+2}$ is in the conductor of the
extended Rees ring $\RR[It, t^{-1}]$. Hence for every $n$, ${\rm
Jac}_k(\RR) \, I \, \overline{I^{n+g-1}} \subset I^{n+1}$, which
yields
\[
{\rm Jac}_k(\RR) \, \overline{I^{n+g-1}} \subset I^{n+1} \colon I
 \subset D_{n+1}\colon I = D_n,
\]
as $\gr(\DD)_{+}$ has positive grade. 
\end{proof}

This result, with an early application of Theorem~\ref{newbs}, gives
the following estimation.

\begin{Corollary}
Let $k$ be a perfect field, let $(\RR, \mathfrak{m})$ be a reduced local
 $k$--algebra essentially of finite type of dimension $d$ and let
$I$ be an  $\mathfrak{m}$-primary ideal. If the Jacobian ideal $L$ of
$\RR$ is  $\mathfrak{m}$-primary, then for any minimal reduction $J$ of
$I$,
\[ \bar{\rme}_1(I) -\rme_1(J) \leq (d+ \lambda(\RR/L)-1) \rme_0(I).\]
Moreover, if $L\neq \RR$ and $\bar{I}\neq \mathfrak{m}$, one replaces $-1$ by $-2$.

\end{Corollary}

\subsubsection*{Lower bounds for $\bar{\rme}_1$}
Let $(\RR, \mathfrak{m})$ be a Noetherian local ring of dimension
$d\geq 1$ and let $I$ be an $\mathfrak{m}$-primary ideal. If $\RR$ is
Cohen-Macaulay, the original lower bound for $\rme_1(I)$ was provided by
Narita (\cite{Narita} and Northcott (\cite{Northcott}),
\[\rme_1(I) \geq \rme_0(I) - \lambda(\RR/I).\]
It has been improved in several ways (see a detailed discussion in
\cite{RV10}). For non Cohen-Macaulay rings, estimates for $\rme_1(I)$
are in a state of flux.

\medskip

We are going to experiment with a special class of non Cohen-Macaulay
rings and methods in seeking lower bounds for $\bar{\rme}_1(I)$. We are
going to assume that $\RR$ is a normal domain and the minimal reduction
$J$ of $I$ is generated by a $d$-sequence. This is the case of normal
Buchsbaum rings, examples of which can be constructed by a machinery
developed in \cite{EvGr} (see \cite[lookup]{StuckradVogel}).

\medskip

Let $\AA= \RR[Jt]$ and $\BB=\bar{\RR[Jt]}$. The corresponding Sally
module $S = S_{\BB/\AA}$ is defined by the exact sequence
\[ 0 \rar B_1t \AA \lar \BB_{+} \lar S = \bigoplus_{n\geq
2}B_n/B_1J^{n-1} \rar 0. \]

\begin{Lemma} $(\RR, \mathfrak{m})$ be an analytically unramified
local ring and let $J$
be an ideal generated by a system of
parameters $x_1, \ldots, x_d$.
Then
\[ \RR/\overline{J}\otimes_\RR \gr_J(\RR) \simeq \RR/\overline{J}[\TT_1,
\ldots, \TT_d].  \]
\end{Lemma}

\begin{proof} Let $\RR[\TT_1, \ldots, \TT_d]\lar \AA$ be a minimal presentation of
 $\AA=\RR[Jt]$.  According to
  \cite[lookup]{Reesred}, the presentation ideal $\mathcal{L}$ has all
  coefficients in
  $\overline{J}$, that is $\mathcal{L}\subset \bar{J}\RR[\TT_1, \ldots,
  \TT_d]$. 
  \end{proof}

\begin{Proposition}
 Let $(\RR, \mathfrak{m})$ be a normal, analytically unramified
local domain and let $J$
be an ideal generated by a system of
parameters $=(x_1, \ldots, x_d)$ of linear type.
The Sally module $\SS$ defined  above is either $0$ or a
module of dimension $d$ and multiplicity
\[ \rme_0(S)=\deg(S)= \bar{\rme}_1(I)-\rme_0(I) -\rme_1(J)+\lambda(\RR/\bar{I}).\]
\end{Proposition}

\begin{proof} By \cite[Corollary 4.6]{HSV3}, since $\RR$ is normal of dimension
$d\geq 2$, $\depth \gr_J(\RR)\geq 2$. Now we make use of \cite[Lemma
1.1]{Hu0} (see also \cite[Proposition 3.11]{icbook}), $\depth
\RR[Jt]\geq 2$ as well.
Now, by \cite[Theorem 3.53]{icbook},
  $\RR[Jt]$ has the condition $S_2$ of Serre.

Consider the exact sequence
\[ 0 \rar B_1 \RR[Jt] \lar \RR[Jt] \lar  \RR/\overline{J}\otimes_\RR
\gr_J(\RR) \rar 0. \]
By the Lemma, the last algebra is a polynomial ring in $d$ variables.
Therefore $B_1\RR[Jt]$ satisfies the condition $S_2$ of Serre.
Thus in
 the defining
sequence of $\SS$, either $\SS$ vanishes (and $B_n= B_1J^n$ for $n\geq
2$) or $\dim S = d$. In the last case, the calculation of the Hilbert
function (see \cite[Remark 2.17]{icbook})
 of $\SS$ gives the value for its multiplicity.
\end{proof}

%It is also possible to derive the first Hilbert coefficient of $\SS$ in
%terms of $\rme_2(I)$.

\subsubsection*{Local rings} We shall now treat
bounds for $\bar{\rme}_1$ for certain non Cohen-Macaulay local rings. Let
 $(\RR, \mathfrak{m})$ be a local ring  which is
Cohen-Macaulay on the punctured spectrum. Without further ado we
assume that the residue field of $\RR$ is infinite.

Suppose $\dim \RR = d > 0$. By local duality, the condition is
equivalent to $\lambda(\H^i_{\mathfrak{m}}(\RR))< \infty$ for $0\leq i <
d$. Later we shall assume additional properties for $\RR$. Let
$I$ be an $\mathfrak{m}$-primary ideal and
let $\mathcal{A}=\{A_n, n\geq 0\}$ be a $I$-good filtration. We may assume that
$I=J = (x_1, \ldots, x_d)$ is a parameter ideal.

We will consider some reductions on $\RR\rar \RR'$ that preserve the
multiplicities and
$\rme_1(\mathcal{A})=\rme_1(\mathcal{A}')$, for
$\mathcal{A}'=\{A_n\RR', n\geq 0\}$.

\begin{enumerate}

\item[{\rm (a)}]
  If $d\geq 2$ and
$\RR'=\RR/\H^0_{\mathfrak{m}}(\RR)$, that will be the case.
In addition $\H^i_{\mathfrak{m}}(\RR)= \H^i_{\mathfrak{m}}(\RR')
$ for $i\geq 1$.

\item[{\rm (b)}] Another property
is that if $d\geq 2$ and $x_1$ is a superficial element for
$\mathcal{A}$ (that is $\RR$-regular if $d=2$),  then preservation will
hold passing to $\RR_1=\RR/(x_1)$.
As for the lengths of the local cohomology modules, in case $x_1$
is $\RR$-regular, from
\[ 0 \rar \RR \stackrel{x}{\lar} \RR \lar \RR_1\rar 0\]
we have the exact sequence
\[ 0 \rar  \H^0_{\mathfrak{m}}(\RR_1) \rar  \H^1_{\mathfrak{m}}(\RR)
\rar  \H^1_{\mathfrak{m}}(\RR)  \rar  \H^1_{\mathfrak{m}}(\RR_1)  \rar
\H^2_{\mathfrak{m}}(\RR)  \rar \cdots,
\] that gives
\begin{eqnarray*}
\lambda(\H^0_{\mathfrak{m}}(\RR_1)) & \leq  & \lambda(
\H^1_{\mathfrak{m}}(\RR)) \\
\lambda(\H^1_{\mathfrak{m}}(\RR_1)) & \leq & \lambda(
\H^1_{\mathfrak{m}}(\RR)) +
 \lambda( \H^2_{\mathfrak{m}}(\RR)) \\
 & \vdots & \\
 \lambda(\H^{d-2}_{\mathfrak{m}}(\RR_1)) & \leq &
\lambda( \H^{d-2}_{\mathfrak{m}}(\RR))    +
\lambda(\H^{d-1}_{\mathfrak{m}}(\RR)).\\
\end{eqnarray*}

\item[{\rm (c)}] Let us combine the two transformations. Let $T_0 =
\H^0_{\mathfrak{m}}(\RR)$, set $\RR'= \RR/T_0$, let $x\in I$ be a superficial
element for $\mathcal{A}R'$ and set $\RR_1= \RR'/(x)$. As $x$ is regular
on $\RR'$,  we have the exact
sequence
\[ 0 \rar T_0/xT_0 \lar \RR/(x) \lar \RR'/x\RR'\rar 0,\]
and the associated exact sequence
\[ 0 \rar T_0/xT_0 \lar \H^0_{\mathfrak{m}}(\RR/(x)) \lar
\H^0_{\mathfrak{m}}(\RR'/x\RR')
 \rar 0,\]
since $\H^1(T_0/xT_0)=0$. Note that this gives
\[ \RR/(x)/\H^0_{\mathfrak{m}}(\RR/(x)) \simeq
\RR'/x\RR'/\H^0_{\mathfrak{m}}(\RR'/x\RR'). \]

There are two consequences to this calculation:
\begin{eqnarray*}
\lambda(\H^0_{\mathfrak{m}}(\RR/(x))) & \leq &
\lambda(\H^0_{\mathfrak{m}}(\RR)) +
\lambda(\H^0_{\mathfrak{m}}(\RR'/x\RR'))\\
& \leq & \lambda(\H^0_{\mathfrak{m}}(\RR)) +
\lambda(\H^1_{\mathfrak{m}}(\RR))\\
\lambda(\H^i_{\mathfrak{m}}(\RR/(x))) & \leq &
\lambda(\H^i_{\mathfrak{m}}(\RR)) +
\lambda(\H^{i+1}_{\mathfrak{m}}(\RR)), \quad 1\leq i \leq d-2\\
\end{eqnarray*}
\end{enumerate}

\begin{Proposition} \label{reddim1} Let $(\RR, \mathfrak{m})$ be a local ring of
positive depth $I$ and $\mathcal{A}$ as
above. If $d\geq 2$, consider a sequence of  $d-1$ reductions of the
type $\RR\rar
\RR/(x)$,
 and denote by $\SS$ the ring $\RR/(x_1, \ldots, x_{d-1})$. Then $\dim S = 1$ and
\[ \lambda(\H^0_{\mathfrak{m}}(S)) \leq T(\RR)= \sum_{i=1}^{d-1}
{{d-2}\choose{i-1}}\lambda(\H^i_{\mathfrak{m}}(\RR)).
\] Moreover, if $\RR$ is a Buchsbaum ring, equality holds.
\end{Proposition}

Let us illustrate an elementary but useful kind of reduction. Let
$(\RR, \mathfrak{m})$ be a Noetherian  local domain of dimension $d\geq
2$ and let $I$ be
an $\mathfrak{m}$-ideal. Suppose $\RR$ has a finite extension $\SS$
with the condition  $S_2$ of
Serre,
\[ 0 \rar \RR \lar \SS \lar C \rar 0.\]

 Consider the polynomial ring $\RR[x,y]$ and
tensor the sequence by  $\RR'=\RR[x,y]_{\mathfrak{m}[x,y]}$.
Then there are $a,b\in I$
such that
the ideal generated by
 polynomial $\ff=ax+by$ has the  following type of primary
 decomposition:
\[ (\ff) = P \cap Q,\]
where $P$ is a minimal prime of $\ff$ and $Q$ is
$\mathfrak{m}R'$-primary. In the sequence
\[ 0 \rar \RR' \lar S' = \RR'\otimes_\RR \SS \lar C'= \RR'\otimes_\RR C \rar 0 \]
we can find $a\in I$ which is superficial for $C$. Pick now $b\in I$
so that $(a,b)$ has codimension $2$. Now reduce the second sequence
modulo $\ff=ax+by$. Noting that $\ff$ will be superficial for $C'$,
we will have an exact sequence
\[0 \rar T \lar \RR'/(\ff) \lar \SS'/\ff \SS' \lar C'/\ff C' \rar 0,\]
in which $T$ has finite length and $\SS'/\ff \SS' $ is an integral
domain since $a,b $ is a regular sequence in $\SS$. This suffices to
establish the assertion.

\medskip

We are now going to consider generic reductions on $\RR$. Let $\XX$
be a $d\times (d-1)$ matrix $\XX=(x_{ij})$ in $d(d-1)$ indeterminates
and let $\CC$ be the local ring
$\RR[\XX]_{\mathfrak{m}[\XX]}$. The filtration $\mathcal{A}\CC$
has the same Hilbert polynomial as $\mathcal{A}$.
If $I =(x_1, \ldots, x_d)$ we now define the ideal
\[ (f_1, \ldots, f_{d-1})=(x_1, \ldots, x_d) \cdot \XX. \]

\begin{Proposition} Let $\RR$ be a local  integral
domain that is Cohen-Macaulay on the punctured spectrum
and $I$ and $\mathcal{A}$ defined as above. Then
 $\SS= \CC/(f_1, \ldots, f_{d-1})$ is a local ring of dimension one
 such that $\lambda(\H^0_{\mathfrak{m}}(\SS)) \leq T(\RR)$ and
$\SS/\H^0_{\mathfrak{m}}(\SS)$ is an integral domain.

\end{Proposition}

The next result shows the existence of bounds for $\bar{\rme}_1(I)$ as
in \cite{ni1}.

\bigskip

We outline the strategy to find bounds for $\rme_1(\mathcal{A})$.
Suppose $(\RR, \mathfrak{m})$ is a local domain essentially of finite
type over a field and $I$ is an $\mathfrak{m}$-primary ideal. We
denote by $\mathcal{A}$ an $I$-good filtration of
$\mathfrak{m}$-primary ideals. We will seek a reduction
to a one dimensional ring $\RR\rar \RR'$, where $\rme_0(I) = \rme_0(I\RR')$,
$\rme_1(\mathcal{
A}) =\rme_1(\mathcal{A}\RR')$, and in the sequence
\[ 0 \rar T = \H^0_{\mathfrak{m}}(\RR') \lar \RR' \lar \SS \rar 0, \]
$T$ is a prime ideal.

\medskip

Since $\SS$ is a one-dimensional integral domain, we have the following
result.

\begin{Proposition} \label{reddim1a} In these conditions,
\begin{eqnarray*}
\rme_1(\mathcal{A})=\rme_1(\mathcal{A}\RR') & = &\rme_1(\mathcal{A}\SS) -
\lambda(T) \\
\rme_1(\mathcal{A}\SS) &\leq &
\bar{\rme}_1(I) = \lambda(\bar{\SS}/\SS),
\end{eqnarray*}
where $\bar{\SS}$ is the integral closure of $\SS$.
\end{Proposition}

This shows that to find bounds for $\rme_1(\mathcal{A})$ one needs to
trace back to the original ring $\RR$ the properties of $T$ and $\SS$. As
to $T$, this is realized if $\RR$
 is  Cohen-Macaulay on the punctured spectrum  for the reductions described in
 Proposition~\ref{reddim1}.

\begin{Theorem}[Existence of Bounds] \label{ebounds} Let $(\RR, \mathfrak{m})$ be a
local integral domain of dimension $d$,
essentially of finite type over a perfect field, and let  $I$ be an
$\mathfrak{m}$-primary ideal. If $\RR$ is  Cohen-Macaulay on the punctured spectrum
 and $\delta$ is a nonzero element of the Jacobian ideal of $\RR$
then for any $I$-good filtration $\mathcal{A}$ as above,
\[\rme_1(\mathcal{A}) < (d-1) \rme_0(I) +\rme_0((I,\delta)/(\delta))- T.\]
\end{Theorem}

\begin{proof} It is a consequence of the proof of \cite[Theorem 3.2(a)]{ni1}.
The assertion there is that
 \[\bar{\rme}_1(I)
\leq \frac{t}{t + 1} \bigl[ (d-1)\rme_0(I) + \rme_0((I +\delta \RR)/\delta
\RR) \bigr],\]
where $t$ is the Cohen-Macaulay type of $\RR$. Here we apply it to the
reduction in Proposition~\ref{reddim1a}
 \[\bar{\rme}_1(I\SS)
<   (d-1)\rme_0(I) + \rme_0((I +\delta \RR)/\delta
\RR) ,\]
dropping the term involving $t$, over which we lose control in the
reduction process. It should be pointed out that the  element
$\delta$ survives the reduction.
 \end{proof}

\begin{Question}{\rm
There are many kinds of $I$-good filtrations over which it would be
of interest to obtain explicit bounds--e.g. for the $I$-adic
filtration itself. Their study may throw light into the invariants of
$\RR$ that should be looked at in the comparison between
$\rme_0(\mathcal{A})$ and $\rme_1(\mathcal{A})$.
}\end{Question}

\subsubsection*{Extended degrees and $\rme_1$} The derivation of bounds for
$\rme_1(\mathcal{A})$
above required that $\RR$ be a generalized Cohen-Macaulay ring. Let us
do away with this requirement by working with an extended degree
function $\Deg$.
We select an extended degree but measure multiplicity by using the
Hilbert-Samuel multiplicity associated to the ideal $I$.  This seems
natural for the setting.

\medskip

The aim of the game is to get a sense of
$\lambda(\H^0_{\mathfrak{m}}(\RR'))$, at the final step of the
reduction, in terms of invariants of the original ring $\RR$. This is
carried out as follows: Suppose we start out at the integral domain $\RR$ and
select $\ff$ as above. It is clear that in addition to being
superficial of $\RR$, it can also be chosen to be generic for $\Deg$.
This would imply that
\[ \Deg(\RR/(\ff))\leq \Deg(\RR).\]
If we write $\RR'= \RR/(\ff) $,
  $H= \H^0_{\mathfrak{m}}(\RR')$ and $\SS=\RR/H$, then
\[ \Deg(\SS)= \Deg(\RR')- \lambda(H) \leq \Deg(\RR). \]
We can repeat this $d-1$ times and at final reduction $\SS$ still have
that, as in Theorem~\ref{specialtor}, that length of the
$\H^0_{\mathfrak{m}}(\SS)$ is bounded by $\hdeg_I(\RR)-\rme_0(I)$.

\medskip

If we apply Theorem~\ref{boundtorsion} to $\RR$, passing to $\RR'=\RR/(x_1, \ldots, x_{d-1})$, we
have the estimate for $\H^0_{\mathfrak{m}}(\RR')$ so that the formula
(\ref{e1dim1}) can be used:

\begin{Theorem}[Lower Bound for $\rme_1(I)$] \label{lowere1} Let $(\RR, \mathfrak{m})$ be a
Noetherian local ring of dimension $d\geq 1$. If $I$ is an
$\mathfrak{m}$-primary ideal, then
\[\rme_1(I) \geq -\hdeg_I(\RR) + \rme_0(I).\]$($Note that
$\hdeg_I(\RR)-\rme_0(I)$ is the Cohen-Macaulay defficiency of $\RR$ relative
to the degree function $\hdeg_I$.$)$
\end{Theorem}

\begin{Corollary}
 \label{bounde1bbar} Let $(\RR, \mathfrak{m})$ be a
Noetherian local ring of dimension $d\geq 1$. For all parameter
ideals $I$ with the same integral closure $\bar{I}$
\[\rme_1(I) \geq -\hdeg_I(\RR) + \rme_0(I).\]
\end{Corollary}

Note that the right-hand side depends only on $\bar{I}$.

\subsection{Modules with finite local cohomology} The following is a
characterization of modules where the values of $\rme_1(\cdot)$ are
bounded.

\begin{Theorem} \label{e1vsgcm} Let $(\RR, \m)$ be a Gorenstein local ring and $M$ a
finitely generated $\RR$-module. Then $M$ is a generalized
Cohen-Macaulay module if and only the Hilbert coefficients $\rme_1(J;M)$
are bounded for all parameter ideals $J$.
\end{Theorem}

\begin{proof}
We have already remarked that if $M$ is a generalized Cohen-Macaulay
module,
\[\rme_1(J;M) \geq - \sum_{i=1}^{d-1} {{d-2}\choose{i-1}}
\lambda(\H_{\m}^i(M)).
\]

Suppose now that  $\rme_1(J;M)$ is bounded.
We may assume $\dim \RR=\dim M=d\geq 2$.
First we are going to consider the case $d=2$.
We may replace $M$ by $M/\H_{\m}^0(M)$, a change  which leaves both conditions
unaffected. (Regarding  the boundedness of $\rme_1(J;M)$, we appeal to
Proposition~\ref{genhs}.)

Let $h\in J$ be a superficial element. Consider the exact sequence
\[ 0 \rar M \stackrel{h}{\lar} M \lar \bar{M} \rar 0.
\]
Applying $\Hom_{\RR}(\cdot, \RR)$, we obtain the exact complex
\[
 \Ext_{\RR}^1({M}, \RR) \stackrel{h}{\rar}
\Ext_{\RR}^1({M}, \RR) \rar
\Ext_{\RR}^2(\bar{M}, \RR) \rar 0\]
since
$\Ext_{\RR}^2({M}, \RR)=0$.

 This gives the embedding
\[ \Ext_{\RR}^{1}(M,\RR)\otimes \RR/(h)\hookrightarrow
\Ext_{\RR}^2(\bar{M}, \RR).
\]
This last module is the dual of $\H_{\m}^0(\bar{M})$, whose length
according to (\ref{e1dim1}) gives the value of $-\rme_1(J;M)$.
If this value is bounded for all $J$, $\Ext_{\RR}^1(M, \RR)$ cannot
be a module of dimension $1$ since replacing $h$ by a one of its powers
would give a contradiction.

\medskip

If $d\geq 3$, as $\rme_1(J;M)=\rme_1(J;\bar{M})$, $\bar{M}$ is a
generalized Cohen-Macaulay module. Thus for each $1\geq i$,
the exact sequence
\[
\Ext_{\RR}^{i}(\bar{M}, \RR) \rar
 \Ext_{\RR}^i({M}, \RR) \stackrel{h}{\rar}
\Ext_{\RR}^i({M}, \RR) \rar
\Ext_{\RR}^{i+1}(\bar{M}, \RR) \]
has the ends
$\Ext_{\RR}^i(\bar{M}, \RR) $ of finite length. Therefore the modules
$\Ext_{\RR}^i({M}, \RR) $ have dimension at most one. In the case of
equality,  from the
argument above, it would follow that the modules $\Ext_{\RR}^{i+1}(\bar{M},
\RR)$ have unbounded length, which is a contradiction. 
\end{proof}

\subsection{Euler characteristics}

Let $(\RR,\m)$
be a Noetherian local ring of dimension $d>0$, with infinite residue
field. Let $J=(a_1, \ldots, a_d)$ be a parameter ideal. We seek to
make comparisons of the Hilbert coefficient $\rme_1(J)$ to another
invariant of $J$. In \cite{chern}, \cite{chern2} and \cite{chern3},
$\rme_1(J)$ is employed to study the Cohen-Macaulayness of $\RR$. There
is another invariant of $J$ that has this usage, the Euler
characteristic $\chi_1(J)$ of Serre.
We will argue that $\chi_1(J)$ and $\rme_1(J)$ have  similar roles
as predictors of Cohen-Macaulayness, and thus it is natural to
compare one invariant to the other.

\medskip

We begin our discussion with a challenge.
Let $\mathbb{K}$ be the Koszul complex associated to $J$, and denote
by $\lambda(H_i(\mathbb{K}))=h_i(J)$. The first two Euler characteristics
of $\H_{\bullet}(\mathbb({K}))$ have important properties according to a
result of Serre (see \cite[Appendice II]{Serrebook},
\cite[Theorem 4.6.10]{BH}):
\begin{eqnarray*}
\rme_0(J)=\chi_0(J) &=& \sum_{i\geq 0} (-1)^i h_i(J),\\
\chi_1(J) &=& \sum_{i\geq 1} (-1)^{i-1} h_i(J).
\end{eqnarray*}

The coefficient $\chi_1(J)\geq 0$, and vanishes if and only if $\RR$
is Cohen-Macaulay. It does have properties similar to those of
$\rme_1(J)$--more precisely, of its negative--and therefore it is of
interest to attempt a comparison between these coefficients.
We consider special cases, but still include some large classes of
rings and ideals.

\medskip

Suppose $J$ is generated by a $d$-sequence. The corresponding
approximation complex $\mathcal{M}(J)$ ($\SS=\RR[\TT_1, \ldots,
\TT_d]$)
\[ 0\rar \H_d(\mathbb{K})\otimes \SS[-d] \rar \cdots \rar
\H_1(\mathbb{K})\otimes \SS[-1] \rar  \H_0(\mathbb{K})\otimes \SS \rar
\gr_J(\RR)\rar 0  \]
 is acyclic,
 which yields the
Hilbert series of $J$ (see \cite[Theorem 5.6]{HSV3}, \cite[Corollary
12.4]{HSV1}),
\[ H_J(\tt)= \frac{\sum_{i=0}^d (-1)^i h_i(J) \tt^i}{(1-\tt)^d}.\]
This gives for $-\rme_1(J)$ the value
\[ -\rme_1(J) = \sum_{i=1}^d (-1)^{i-1}\cdot i \cdot h_i(J). \]

In the special case of $d=2$, we have
\begin{eqnarray*}
-\rme_1(J) &=& h_1(J)-2\cdot h_2(J)\\
\chi_1(J) &=& h_1(J)-h_2(J),
\end{eqnarray*}
and therefore $\chi_1(J)\geq -\rme_1(J)\geq 0$. In particular
this explains why $\chi_1(J)=0$ implies that $\RR$ is Cohen-Macaulay
(at least when $\RR$ is unmixed).

\begin{Question}\label{e1vschi1}{\rm For any parameter ideal $J$ of an unmixed local
ring, when is
\begin{eqnarray*}
&&\psboxit{box 1.0 setgray fill}{\fbox{$
\begin{array}{c}
\ \\
 \chi_1(J)\geq -\rme_1(J).\\
\ \\
\end{array}  $}}
\end{eqnarray*}
In the case of ideals generated by $d$-sequences, this asks whether
\[ \sum_{i\geq 2} (-1)^i (i-1)h_i(J) \leq 0.
\]
}\end{Question}

One affirmative case follows from a straightforward calculation:

\begin{proposition} If $\RR$ is a Buchsbaum ring the inequality
holds.
\end{proposition}

\begin{proof} We make use of basic properties of Buchsbaum rings as developed
in \cite[Chap. 1, Section 1]{StuckradVogel}. Two of these we use
repeatedly: (i) If $\xx=\{x_1, \ldots, x_n\}$ is part of a system of
parameters, then $(x_1, \ldots, x_{n-1}): x_n=
(x_1, \ldots, x_{n-1}): \m$ (in particular, $\xx$ is a $d$-sequence), and (ii) the Koszul homology modules
$\H_i(\mathbb{K}(\xx))$ are $\RR/\m$-vector spaces for $i\geq 1$.

\medskip

We use the following well-known technique
to relate the Koszul homology of a partial set of parameters
$\xx$ to that of the enlarged
sequence $\xx,y$ obtained by adding one more parameter.

\begin{proposition} \label{growthseq}
Let ${\mathbb C}$ be a chain complex and let ${\mathbb F}= \{F_1,F_0\}$ be a chain
complex of free modules concentrated in degrees $1$ and $0$. Then
for each integer $q\geq 0$ there is an exact sequence
\[0\rar \H_0(\H_q({\mathbb C})\otimes {\mathbb F})\lar \H_q({\mathbb C}\otimes{\mathbb
F})\lar \H_1(\H_{q-1}({\mathbb C})\otimes{\mathbb F})\rar 0.\]
\end{proposition}

Setting ${\mathbb C}=\mathbb{K}(\xx)$,
${\mathbb F}=\mathbb{K}(y)$, and
${\mathbb K}=\mathbb{K}(\xx,y)$, we have equalities of dimensions of
vector spaces
\begin{eqnarray*}
h_1(\mathbb{K}) &=& h_1(\mathbb{K}_0)+ \lambda(O:_{\RR/(\xx)} y),\\
h_i(\mathbb{K}) &=& h_i(\mathbb{K}_0)+
h_{i-1}(\mathbb{K}_0), \quad i\geq 2,\\
h_{n+1}(\mathbb{K}) &=& h_n(\mathbb{K}_0).  \\
\end{eqnarray*}

In case $n=d-1$, for $d=\dim \RR$, $J=(\xx,y)$,
these relations would give
\begin{eqnarray*}
 \chi_1(J) &=& \lambda(O:_{\RR/(\xx)} y ), \\
-\rme_1(J) &=& \lambda(O:_{\RR/(\xx)} y)
+\sum_{i=1}^{n-1}(-1)^i h_i(\mathbb{K}_0).\\
\end{eqnarray*}

For $d=2$, the assertion is clear. For $d>2$, it is enough to show
that $\sum_{i=1}^{n-1}(-1)^i h_i(\mathbb{K}_0)\leq 0$.
In this case, write $\xx_0=\{x_1, \ldots, x_{n-2}\}$ and apply
Proposition~\ref{growthseq}  to get (with some abuse of notation),
\[ \chi_1(\xx) = \sum_{i=1}^{n-1} h_i(\xx) =
 \lambda(O:_{\RR/(\xx_0)} x_{n-1} ),
\] as desired. 
\end{proof}

We also observe that Buchsbaum rings are those for which $\rme_1(J)$ is
constant (\cite{GO1}), which makes one wonder whether $\chi_1(J)$ will
be constant as well.

\medskip

This makes use of the fact that the Koszul homology modules are easy
to manipulate. A more careful, but still direct analysis will work for
generalized Cohen-Macaulay ring.

\begin{remark}{\rm
Let $\RR$ be a generalized Cohen-Macaulay ring of dimension $d\geq 1$
and let $\xx=\{x_1, \ldots, x_{d-1}\}$, which together with $x_d$
give a parameter ideal.
By assumption, the homology modules $\H_s(\mathbb{K}(\xx))$, $s\geq
1$, have finite length. Therefore the two homology modules
$\H_j(\mathbb{K}(\H_s(\mathbb{K}(\xx);x_d)))$, for $j=0,1$, have the
same length. These facts lead again to the expression of
$\chi_1(\xx,x_d)$,
\[
 \chi_1(\xx,x_d) = \lambda(O:_{\RR/(\xx)} x_d ). \]

Observe the following consequence: If $\RR$ is unmixed, then $\RR$ is Cohen-Macaulay if and
only  if $\RR/(\xx)$ is Cohen-Macaulay .

}\end{remark}

\subsection{Chern coefficient and homological torsion}\index{Chern
coefficient and homological torsion}

We now this into a uniform bound for the first Hilbert coefficient of
a module $M$ relative to an ideal $I$ generated by a system of
parameters $\xx =\{x_1, \ldots, x_d\}$.

\begin{Theorem} \label{Degreddim2} Let $M$ be a module  of
dimension  $d\geq 2$
and let $\xx=\{x_1, \ldots, x_{d-1}, x_d\}$ be a superficial sequence for
$M$ and $\hdeg$. Then
\[ -\rme_1(I;M)\leq \TT(M).\]
\end{Theorem}

\begin{proof} We first argue a reduction to the case $d=2$. If $d\geq 3$, we
select for $x_1$ a superficial element for the Hilbert function and
$\hdeg$. This follows from  Proposition~\ref{genhs} and Theorem~\ref{torsionhdeg}
\begin{eqnarray*}
\rme_1(I;M) &=&\rme_1(I/(x_1);M/x_1M),\\
\TT(M) &\geq & \TT(M/x_1M).
\end{eqnarray*}

In the exact  sequence
\[ 0 \rar H = \H^0_{\mathfrak{m}}(M) \lar M \lar M' \rar 0\]
consider reduction mod  the superficial element $h=x_2$:
\begin{eqnarray*}
\lambda(\H^0_{\mathfrak{m}}(M/hM)) &\leq &
\lambda(\H^0_{\mathfrak{m}}(M'/hM')) +
\lambda(H/hH) \\
&\leq &
\lambda(\H^0_{\mathfrak{m}}(M'/hM')) +
\lambda(0:_H h),
\end{eqnarray*}
where $\lambda(H/hH)=\lambda(0:_Hh)$ since $H$ is a module of finite
length (the equality results from the exact sequence $0\rar 0:_Hh
\rar H \rar H \rar h/hH \rar 0$ induced by multiplication by $h$).
Because $h$ is superficial, $0:_Mh = 0:_H h$, which yields the
inequality
\begin{eqnarray*}
\lambda(\H^0_{\mathfrak{m}}(M/hM)) -\lambda(0:_Mh)
&\leq & \lambda(\H^0_{\mathfrak{m}}(M'/hM')),
\end{eqnarray*}
in other words, by Proposition~\ref{genhs}
\begin{eqnarray*}
-\rme_1(I;M) &\leq & -\rme_1(I;M').
\end{eqnarray*}

\medskip

The proof of the inequality $-\rme_1(I;M')\leq \TT(M')$ is as before:
Write  $h=\xx$. The assertion requires that
$\lambda(\H^0_{\mathfrak{m}}(M'/hM')) \leq \hdeg(\Ext_{\SS}^1(M',\SS)).$
We have the cohomology exact sequence
\[ \Ext^1_{\SS}(M',\SS) \stackrel{h}{\lar} \Ext^1_{\SS}(M',\SS)
\lar \Ext^2_{\SS}(M'/hM',\SS)\lar \Ext^2_{\SS}(M',\SS)=0,
\] where
\[ \lambda(\H^0_{\mathfrak{m}}(M'/hM'))
=\hdeg(\Ext_{\SS}^2(M'/hM',\SS)),\]
and we argue as above.
\end{proof}

\section{Euler Characteristics}\label{eulerchar}\index{Euler
characteristic}

\subsection{Introduction}

 Let $Q=(a_1, a_2, \ldots, a_d)$ be a
parameter ideal in a given Noetherian local ring $A$ of dimension $d
>0$. We seek to make comparisons of the first Hilbert coefficient
$\rme_1(Q,A)$ to another invariant of $Q$. For simplicity, let us
denote by $\rme_i(Q,A) = \rme_i(Q)=\rme_i(\aa)$ for $0 \le i \le d$,
where $\aa$ denotes the sequence $\{a_1, a_2, \ldots, a_d\}$. In
\cite{chern3}, \cite{chern5}, \cite{chern6}, \cite{GO1}, \cite{GO2}
and \cite{chern}, $\rme_1(Q)$ is employed to study the
Cohen-Macaulayness of $A$. There is another invariant of $Q$ that has
this usage, the Euler characteristic $\chi_1(Q)$ of Serre.  We will
argue that $\chi_1(Q)$ and $\rme_1(Q)$ have similar roles as
predictors of Cohen-Macaulayness, and thus it is natural to compare
one invariant to the other.

\subsection{Partial Euler characteristics}\index{partial Euler
characteristic}
%\medskip

Let $\mathbb{K}=\mathbb{K}_{\bullet}(\aa;A)$ be the Koszul complex
associated to the sequence $\aa$ and denote by
$\lambda_A(\H_i(\mathbb{K}))=\rmh_i(Q)$ (or, by $\rmh_i(\aa)$). The
first two Euler characteristics of $\H_{\bullet}(\mathbb{K})$ have
important properties according to a result of Serre (see
\cite[Appendice II]{Serrebook}, \cite[Theorem 4.6.10]{BH}):
\begin{eqnarray*}
\rme_0(Q)=\chi_0(Q) = \chi_0(\aa) &=& \sum_{i\geq 0} (-1)^i \rmh_i(Q),\\
\chi_1(Q) = \chi_1(\aa)&=& \sum_{i\geq 1} (-1)^{i-1} \rmh_i(Q).
\end{eqnarray*}

The coefficient $\chi_1(Q)\geq 0$, and vanishes if and only if $A$
is Cohen-Macaulay. It does have properties similar to those of
$\rme_1(Q)$ -- more precisely, of its negative -- and therefore it is of
interest to attempt a comparison between these coefficients.
We consider special cases, but still include some large classes of
rings and ideals. Another extension is that of the following notion
(see
 Corollary~\ref{subsys} for  examples).

%\medskip

\begin{Definition} Let $(\RR, \m)$ be a Noetherian local ring of dimension
$d\geq 1$. A system of elements $\xx=\{x_1, \ldots, x_n\}$ is
suitable for defining {\em  partial Euler characteristics}
if $\H_j(x_1, \ldots, x_i)$ have finite length for all $i$ and $j\geq
1$.  We set
\[ \chi_1(\xx)= \sum_{j\geq 1}^n (-1)^{j-1}\rmh_j(\xx).\]
\end{Definition}

A special type of systems of parameters that play an important role
in this paper is that of {\em d-sequences} and one of its
generalizations.
 They are
  extensions of regular sequences, the first the notion of  {\em
  $d$-sequence} \index{d-sequence} invented by Huneke (\cite{Hu1}),
  and of a {\em proper sequence} (\cite{HSV3}). They play natural
  roles in the theory of the {\em approximation complexes}
  (\cite{HSV3}).

We now collect appropriate properties of these
sequences. (The notion extends {\it ipsis literis} to modules.)
 Let $\RR$ be a commutative ring.

  \begin{Definition} {\rm Let $\xx= \{x_1, x_2, \ldots, x_n\}$ be a sequence of
elements in $\RR$. $\xx$ is a  $d$-sequence if
\[ (x_1, x_2, \ldots, x_i): x_{i+1}x_k = (x_1, x_2, \ldots, x_i): x_k, \quad
\mbox{for } i=0,\ldots, n-1, \ k\geq i+1.\]
}\end{Definition}

Basic properties of d-sequences, useful in induction arguments, are
 (\cite{Hu1}):

\begin{Proposition} \label{dseqprop}
Suppose that ${\xx} = \{x_1, x_2, \ldots, x_n\}$ is a d-sequence in
$\RR$. Then we have the following.
\begin{enumerate}
\item[{\rm (a)}] The images of $x_2, x_3, \ldots, x_n$ in $\RR/(x_1)$ form a
d-sequence.
\item[{\rm (b)}] $[(0)\colon x_1]\cap (x_1, x_2, \ldots, x_n)=(0). $
\item[{\rm (c)}] The images of $x_1, x_2, \ldots, x_n$ in $\RR/[(0)\colon x_1]$ form a
d-sequence.
\item[{\rm (d)}]
$\xx$ is  a regular
sequence if and only if
$ (x_1, \ldots, x_{n-1}): x_n = (x_1, \ldots, x_{n-1})$.
\item[{\rm (e)}] If $\xx$ form a d-sequence then the sequence
$\xx'=\{x_1^*, \ldots, x_n^* \}$ of $\gr_{\xx}(\RR)$ is a d-sequence.
\end{enumerate}
\end{Proposition}

\begin{proof} In (d), the argument is a straightforward calculation.
We may assume $n>1$.
 If $(x_1,\ldots, x_{n-1}):x_n=(x_1,\ldots, x_{n-1})$, we claim that
 $(x_1, \ldots, x_{n-2}):x_{n-1}=(x_1, \ldots, x_{n-2})$.
From $ax_{n-1}\in (a_1, \ldots, x_{n-2})$, we have
 $x_nx_{n-1} a\in (x_1, \ldots, x_{n-2})$ so \[a\in (x_1, \ldots,
 x_{n-2}):x_{n-1}x_n = (x_1, \ldots, x_{n-2}): x_n
\subset  (x_1, \ldots, x_{n-1}): x_n = (x_1, \ldots, x_{n-1}).\]
 Thus $a=a_1x_1 +\cdots + a_{n-1}x_{n-1}$ and from
$x_{n-1}a\in (x_1, \ldots, x_{n-2})$ we get
$a_{n-1}x_{n-1}^2\in (x_1, \ldots, x_{n-2})$ and therefore
  $a_{n-1}x_{n-1}\in (x_1, \ldots, x_{n-2})$, as desired.

(e) is proved in \cite[Theorem 1.2]{Hu81}, converse in \cite{Kuhl}.
\end{proof}

An interesting consequence is:

\begin{Corollary} \label{dseqaresup} \index{d-sequences versus superficial
sequeces} Let $\RR$ be a Noetherian local ring and $M$ a finitely
generated $\RR$-module.
 If $\xx$ is a system of parameters that is a
d-sequence on the module $M$ then $\xx$ is a superficial sequence
with respect to $M$.
\end{Corollary}

\begin{Definition} {\rm Let $\xx= \{x_1, x_2, \ldots, x_n\}$ be a sequence of
elements in $\RR$. $\xx$ is a  proper sequence if
\[ x_{i+1}\H_j(x_1, x_2, \ldots, x_i)=0, \quad
\mbox{for } i=0, 1, \ldots, n-1, \ j>0,\]
where $\H_j(x_1, x_2, \ldots, x_i)$ is the Koszul homology associated to
the subsequence $\{x_1, x_2, \ldots, x_i\}$.\index{proper sequence}
}\end{Definition}

\begin{Proposition}\label{propseqprop} Suppose
 that $\RR$ is a Noetherian ring.
Let   $\xx = \{x_1, x_2, \ldots, x_n\}$ be a sequence in $\RR$.
\begin{enumerate}
\item[{\rm (a)}] If $\xx$ is a d-sequence, then $\xx$ is a proper
sequence.

\item[{\rm (b)}] If $\xx$ is a proper sequence, then \[
x_{k}\H_j(x_1, x_2, \ldots, x_i)=0, \quad \mbox{for } i=0,1, \ldots,
n-1, \ j>0,\ k>i.\]
 \item[{\rm (c)}] Suppose that $\RR$ is a local
ring of dimension $n >0$. If $\xx$ is a proper sequence that is also
a system of parameters in $\RR$, then $\H_j(x_1, x_2, \ldots, x_i)$
is a module of finite length for $i=0,1, \ldots, n$ and $j>0$.
\end{enumerate}
\end{Proposition}

The following provides for a ready source of these sequences.

\begin{Proposition}\label{acidseq} Let $\RR$ be Noetherian  local and
 $\xx=\{x_1,\ldots, x_d \}$ a parameter ideal. If $(x_1, \ldots,
 x_{d-1})$ has grade $d-1$, then
\begin{enumerate}
\item[{\rm (a)}] $\xx$ is a proper sequence;

\item[{\rm (b)}] $\xx$ is a d--sequence if and only if
$ (x_1, \ldots, x_{d-1}): x_d^2=
(x_1, \ldots, x_{d-1}): x_d.$
\end{enumerate}
\end{Proposition}

\subsection{Euler characteristics versus Hilbert coefficients}

In this section we  derive the explicit values for these partial Euler
characteristics in the case of proper sequences. Later in this
section we express the Chern numbers of parameter ideals generated by
d-sequences in terms of the
Euler characteristics.

\medskip

When we apply Proposition~\ref{growthseq} to a proper sequence $\xx=\{x_1, \ldots, x_n\}$,
we obtain (set $\RR_i=\RR/(x_1, \ldots, x_i)$):
\[ 0\rar (0):_{\RR_{i}} x_{i+1} \rar \RR_{i} \stackrel{x_{i+1}}{\rar}
\RR_{i} \rar \RR_{i+1} \rar 0, \]
\[0 \rar \H_1(x_1, \ldots, x_i)\rar
\H_1(x_1, \ldots, x_{i+1}),\]
\[0 \rar \H_j(x_1, \ldots, x_i)\rar
\H_j(x_1, \ldots, x_{i+1})
 \rar \H_{j-1}(x_1, \ldots, x_i)\rar 0, \quad j>1.
\]

\begin{Proposition} \label{chiproper} Let $\RR$ be a Noetherian local ring  of
dimension $d>1$ and
let $\xx=\{x_1, x_2, \ldots, x_n\}$ be a set of generators of an
ideal of codimension $n$. Suppose that for each subset $\xx'=\{x_1,
\ldots, x_i\}$, $1\leq i\leq n$, the Koszul homology modules
$\H_q(\xx')$, $q\geq 1$, have finite length.
  Write  $\xx=\{\xx_0,y\}$ with $\xx_0$ the initial
subsequence up to $x_{n-1}$ and $y=x_n$. Then
\[ \chi_1(\xx) =\lambda_{\RR}((0):_{\RR/(\xx_0)} y).\]
\end{Proposition}

\begin{proof}
 It makes still sense to
define
\[ \chi_1(\xx) = \sum_{i=1}^{n} (-1)^{i-1}\rmh_i(\xx),
\]
since the Koszul homology modules $\H_i(\mathbb{K}(\xx))$, for $i\geq 1$, have finite length.

\medskip

In Proposition~\ref{growthseq},
set ${\mathbf C}=\mathbb{K}(\xx_0)$,
${\mathbb F}=\mathbb{K}(y)$, and
${\mathbb K}=\mathbb{K}(\xx_0,y)$.
Since   the modules
$\H_q(\mathbb{K}(\xx_0))$ have finite length for $q>0$,
the exact sequence
\[ 0\rar \H_1(\mathbb{K}(\H_q(\mathbb{K}(\xx_0);x_{n})))\rar
\H_q(\mathbb{K}(\xx_0))\rar
\H_q(\mathbb{K}(\xx_0))\rar
\H_0(\mathbb{K}(\H_q(\mathbb{K}(\xx_0);x_{n})))
\rar 0
\] implies
that the two modules
\begin{eqnarray*}
A_q &  = & \H_0(\mathbb{K}(\H_q(\mathbb{K}(\xx_0);x_{n})))\\
 B_q &  = & \H_1(\mathbb{K}(\H_q(\mathbb{K}(\xx_0);x_{n})))
\end{eqnarray*}
 have
the same length.
 This gives rise to the
  equalities
\begin{eqnarray*}
\rmh_1(\mathbb{K}) &=& \lambda_\RR(A_1)+ \lambda_{\RR}((0):_{\RR/(\xx_0)} y),\\
\rmh_i(\mathbb{K}) &=& \lambda_\RR(A_i)+ \lambda_\RR(B_{i-1}) , \quad i\geq 2,\\
\rmh_{n}(\mathbb{K}) &=& \lambda_\RR(B_{n-1})  \\
\end{eqnarray*}
which leads to the assertion.
\end{proof}

\begin{Corollary} \label{subsys} Let $\xx=\{x_1, \ldots, x_n\}$ be a set of elements of
one of the following kinds:
\begin{enumerate}
\item[{\rm (a)}] $\xx$ is an initial subsequence of a system of parameters
forming a proper sequence;

\item[{\rm (b)}] $\xx$ is an initial subsequence of a system of parameters
forming a d--sequence;

\item[{\rm (c)}] $\RR$ is a generalized Cohen-Macaulay local ring and $\xx$ is a
subsystem of parameters.
\end{enumerate}
If $\xx$ is any sequence as above and $\xx=\{\xx_0,y\}$, then
\[\chi_1(\xx)=\lambda_{\RR}((0):_{\RR/(\xx_0)}y).\]
Moreover, in case $(2)$,
$\chi_1(\xx)= 0$  if and only if $\xx$ is a regular
sequence.
\end{Corollary}

\begin{proof} The second assertion follows from
Proposition~\ref{dseqprop}(d).
\end{proof}

%\begin{Remark} If $\RR$ is a generalized Cohen-Macaulay local ring, then
%some of the same assertions apply. For any initial subsequence
%$\xx_0$ of a system of parameters,  the modules
%$\H_s(\mathbb{K}(\xx_0))$ have finite length for $s>0$. One derives
%a similar formula for $\chi_1(\xx)$.
%\end{Remark}

\begin{Theorem} \label{chidseq}
 Let $\RR$ be a Noetherian local ring  of
dimension $d>1$ and
let $\xx=\{x_1, x_2, \ldots, x_d\}$ be a system of parameters that is a
 d-sequence in $\RR$.
For $-1\leq i \leq d-1$  set $\xx_{i}=\{x_1, \ldots, x_{d-i-1}\}  $.
Then
\begin{enumerate}
\item[{\rm (a)}]
 $\chi_1(\xx_i)=\lambda_{\RR}((0):_{\RR/(\xx_{i-1})}x_{d-i})$.
\ In particular
$\chi_1(\xx_d)=\lambda_{\RR}((0):_{\RR} x_1)=\lambda_{\RR}(\H_{\m}^0(\RR))$.
\medskip

\item[{\rm (b)}] For $i\geq 1$,
$ (-1)^i \rme_i(\xx)= \chi_1(\xx_{i-1})- \chi_1(\xx_i)\geq 0$.

\end{enumerate}
\end{Theorem}

We first recall the following property of ideals generated by
d-sequences
(see \cite[Theorem 5.6]{HSV3}, \cite[Corollary
12.4]{HSV1}).

\begin{Theorem} \label{dseqhseries}
Let $\xx$ be a system of parameters that is also a d-sequence. The
approximation complex  $\mathcal{M}(\xx)$ associated to $\xx$
\[ 0\rar \H_d(\mathbb{K})\otimes \SS[-d] \rar \cdots \rar
\H_1(\mathbb{K})\otimes \SS[-1] \rar  \H_0(\mathbb{K})\otimes \SS \rar
G\rar 0  \]
($\SS=\RR[\TT_1, \ldots,
\TT_d]$)  is acyclic. For $Q=({\xx})$ and
$G=\gr_Q(\RR)$,
 the
Hilbert series of $G$ is
\[ \H(G, \lambda)= \frac{\sum_{i=0}^d (-1)^i \rmh_i(Q)
\lambda^i}{(1-\lambda)^d}= \frac{\rmh(\lambda)}{(1-\lambda)^d}. \]
\end{Theorem}

\begin{proof} [Proof of {\rm Theorem~\ref{chidseq}}:] We only have to
prove (2). For clarity we consider the cases $i=1$ and $i=2$.
Theorem~\ref{dseqhseries}
 gives for $-\rme_1(\xx)$ the value
\[ -\rme_1(\xx) = -\rmh'(1)= \sum_{i=1}^d (-1)^{i-1}{\cdot} i {\cdot} \rmh_i(Q). \]
In the special case of $d=2$, we have
\begin{eqnarray*}
-\rme_1(\xx) &=& \rmh_1(x_1,x_2)-2{\cdot}\rmh_2(x_1,x_2)\\
 &=& \chi_1(x_1,x_2)-\chi_1(x_1),
\end{eqnarray*}
since $x_2\H_1(x_1)=0$.

For $d>2$, we assemble $\rme_1(\xx)$ from the relations in
Proposition~\ref{chiproper}.

\medskip

 From the Hilbert series, we
have
$ \rme_i(\xx) = \frac{\rmh^{(i)}(1)}{i!}$.
 Set $i=2$. Expanding $\rme_2(\xx)$ and
expressing the $\rmh_i(\xx)$ in terms of the $\rmh_i(\xx_0)$, and these in
terms of the $\rmh_i(\xx_1)$
 we have
\begin{eqnarray*}
\rme_2(\xx) & = & \rmh_2(\xx) -3\rmh_3(\xx) + 6 \rmh_4(\xx) -  \cdots
\\
&=& (\rmh_1(\xx_0)+ \rmh_2(\xx_0)) -3(\rmh_2(\xx_0) + \rmh_3(\xx_0)) +
6(\rmh_3(\xx_0) + \rmh_4(\xx_0)) - \cdots \\
&=&  \rmh_1(\xx_0) -2 \rmh_2(\xx_0) + 3\rmh_3(\xx_0)+ \cdots \\
&=& (\rmh_1(\xx_1) + \lambda(0:_{\RR/(\xx_1)}: x_{d-2}))-2(\rmh_1(\xx_1) +
\rmh_2(\xx_1)) + 3(\rmh_2(\xx_1) + \rmh_3(\xx_1)) - \cdots \\
&=&\lambda_\RR(0:_{\RR/(\xx_1)} x_{d-2})  -\chi_1(\xx_1) \\
&=&\chi_1(\xx_0) -\chi_1(\xx_1).	
\end{eqnarray*}

The proof of the cases $i>2$ is similar.

As for the  assertion of  positivity,
\[ (-1)^i \rme_i(\xx) =
 \lambda(((x_1, \ldots, x_{i-1}):x_i)/(x_1, \ldots,x_{i-1}))-
\lambda((((x_1, \ldots, x_{i-2}):x_{i-1})/(x_1, \ldots,x_{i-2})).
\]
Now observe that since $\xx$ is d-sequence,
$(x_1, \ldots, x_{i-2}):x_{i-1} \subset (x_1, \ldots,
x_{i-2}):x_{i}$, so
there is a natural mapping of
\[(x_1, \ldots, x_{i-2}):x_{i-1})/(x_1, \ldots,x_{i-2}))\lar
(x_1, \ldots, x_{i-1}):x_{i})/(x_1, \ldots,x_{i-1})\]
whose kernel is zero since
$((x_1, \ldots, x_{i-2}):x_{i-1})\cap (x_1, \ldots,x_{i-1})=
(x_1, \ldots,x_{i-2})$. This gives the formula
\[ (-1)^i \rme_i(\xx) =
\lambda(((x_1, \ldots,x_{i-1}):x_i)/((x_1, \ldots,x_{i-1})
+ ((x_1, \ldots,x_{i-2}):x_i))).
\]

\end{proof}

\begin{Corollary} Let $\RR$ be a Noetherian local ring of
dimension $d>1$.
If $\xx$ is a system of parameters that is a d-sequence in $\RR$,
$\chi_1(Q)= -\rme_1(Q)$ if and only if $\operatorname{depth} \RR\geq d-1$, where $Q = (\xx)$.
\end{Corollary}

Equivalent formulas for $\rme_i(Q)$ in terms of local cohomology
modules are given in \cite[Proposition 3.4]{GO2}:
\begin{Proposition}\label{eiform1} Let $Q$ be a parameter ideal generated by a
$d$-sequence $(a_1, \ldots, a_d)$.
 For $0 \le i \le d$ set $Q_i=(a_1, \ldots, a_i)$. Then
\[ (-1)^{i}\rme_i(Q) =  \left\{
\begin{array}{ll}
\lambda_A(A/Q)-\lambda_A \left([Q_{d-1}:a_d]/Q_{d-1}\right) & \quad \mathrm{if} \quad i=0,
\vspace{2mm}\\
\rmh^0(A/Q_{d-i})-\rmh^0(A/Q_{d-i-1}) & \quad \mathrm{if} \quad 1 \le i \le d-1,
\vspace{2mm}\\
\rmh^0(A) & \quad \mathrm{if} \quad i=d.\\
\end{array}
\right.\]
\end{Proposition}

\begin{Question}  This discussion suggests that  just as Buchsbaum
rings are those for which every system of parameters is a d-sequence,
it would be of interest to have a better understanding of the local
rings for which every parameter ideal can be generated by a proper
system.
\end{Question}

Theorem~\ref{chidseq}
settles the question of the comparison between $\chi_1(Q)$ and
$\rme_1(Q)$ for parameter ideals generated by d-sequences.
In particular it asserts that $\chi_1(Q)\geq - \rme_1(Q)$ for those
ideals.
 Let us note an example that the converse does not hold (see also
 \cite[2.4]{chern6}).

\begin{Example}\label{2.3}
Let $\RR$ be a regular local ring with maximal ideal $\n$ and $d =
\operatorname{dim} \RR \ge 3$. Choose an ideal $I$ in $\RR$ so that $\dim
\RR/I = 1$ and $\RR/I$ is a Cohen-Macaulay ring but not a discrete
valuation ring. We look at the idealization $$A = \RR \ltimes \RR/I$$ of
$\RR/I$ over $\RR$. Then for every parameter ideal $Q$ in $A$ we have
\[\chi_1(Q) > 0 = -\rme_1(Q),\] while the parameter ideal $Q = \n
A$ cannot be generated by a $d$--sequence of length $d$, which shows
the converse of Theorem \ref{chidseq} is not true in general.
\end{Example}

\begin{proof} Let $\q = Q\RR$ and put $\cl{\RR} = \RR/I$.
 Then, thanks to the canonical exact sequence
\[0 \to (0) \times \cl{\RR} \to A \overset{p}{\to} \RR \to 0 \] where
$p(r,m) = r$ for each $(r,m) \in A$, we have
\[\lambda_A(A/Q^{n+1}) =
\rme_0(\q)\binom{n+d}{n} + \lambda_\RR(\cl{\RR}/\q^{n+1}\cl{\RR})\] for
all $n \ge 0$, so that \[ (-1)^{i}\rme_i(Q) = \left\{
\begin{array}{ll}
\rme_0(\q) & \quad \mbox{if} \quad i=0,
\vspace{2mm}\\
0 & \quad \mathrm{if} \quad 1 \le i \le d-2,
\vspace{2mm}\\
\rme_0(\q\cl{\RR}) & \quad \mathrm{if} \quad i = d-1,
\vspace{2mm}\\
-\rme_1(\q\cl{\RR}) & \quad \mathrm{if} \quad i=d.
\end{array}
\right.\] Hence $$\chi_1(Q) = \lambda_A(A/Q) - \rme_0(Q) > 0 =
-\rme_1(Q),$$ since $A$ is not a Cohen-Macaulay ring (notice that
$\dim A = d$ and $\depth A = 1$). We have $\rme_1(\q \cl{\RR}) = 0$
if and only if $[\q + I]/I$ is a principal ideal in $\RR/I$, because
$\RR/I$ is a Cohen-Macaulay ring of dimension one. Taking $Q = \n A$,
the maximal ideal $[\n + I]/I$ in $\RR/I$ is not principal, since $\RR/I$
is not a discrete valuation ring. Hence the parameter ideal $Q = \n
A$ in $A$ cannot be generated by any $d$--sequence of length $d$;
otherwise $$\rme_1(\q\RR/I) = (-1)^{d-1}\rme_d(Q) = -\rmh^0(A) = 0$$
by Corollary~\ref{chidseq}, since $\depth A >0$.
\end{proof}

%%%%%%%%%%%%%%%%%%%%%%%%%%%%%%%%%%%%%%%%%%%%%%%%%%%%%%%%%%%%%%%%%%%%%%%%%%%%%%%%

\subsubsection*{A converse}

In this segment we will be dealing mostly with local rings $A$ of
dimension $d\geq 1$ with $\depth A\geq d-1$.
% If the residue field of
%$A$ is infinite and $Q$ is a parameter ideal then we can choose
%$Q=(\xx)=(x_1, \ldots, x_d)$ so that $x_1, \dots, x_{d-1}$ is a regular
%sequence. Note that $\xx$ is a proper sequence
%(Proposition~\ref{acidseq}). The following
For a parameter ideal $Q$,
the key to comparing the values of $\chi_1(Q)$ and $\rme_1(Q)$ is the
following.

 \begin{Proposition}\label{dminus1} Let $A$ be a local ring of dimension
 $d\ge 1$ and infinite residue field. If $Q$ is a parameter ideal
 then it has a generating system $\xx=\{a_1, \ldots, a_d\}$ such that
 $\{a_1, \ldots, x_{d-1}\}$ is a superficial sequence of $A$ relative
 to $Q$. For $Q_{d-1}=(a_1, \ldots, a_{d-1})$
\begin{eqnarray*}
\chi_1(Q) &= &\lambda_A([Q_{d-1}:a_d]/Q_{d-1})\\
 \rme_1(Q) &=&- \lambda_A(\H_{\m}^0(A/Q_{d-1})).
\end{eqnarray*}
Moreover, $\chi_1(Q)=-\rme_1(Q)$ if and only if $\xx$ is a d-sequence.
\end{Proposition}

\begin{proof} Since the residue field of $A$ is infinite, we may
choose a system $a_1, a_2, \ldots, a_d$ of parameters in $A$ so that
$Q = (a_1, a_2, \ldots, a_d)$ and $a_1, a_2, \ldots, a_{d-1}$ forms a
superficial sequence for $A$ with respect to $Q$. Then $a_1, a_2,
\ldots, a_{d-1}$ is an $A$--regular sequence, as $\depth A \ge d-1$.
Thus $\{a_1, \ldots, a_d \}$ is a proper sequence by
Proposition~\ref{acidseq}, and therefore by
Proposition~\ref{chiproper}, $\chi_1(Q) =
\lambda_A([Q_{d-1}:a_d]/Q_{d-1}).$

For the expression for $\rme_1(Q)$, the assertion is clear if $d=1$.
If $d>1$
 let $\cl{A} = A/(a_1)$ and $\cl{Q}= Q/(a_1)$. Then the
system $a_2, a_3, \ldots, a_d$ of parameters for $\cl{A}$ is reducing
also for $\cl{A}$.   $\rme_Q^i(A) =
\rme_{\cl{Q}}^i(\cl{A})$ for all $0 \le i \le d-1$; recall that $a_1$
is $A$--regular and superficial for $A$ with respect to $Q$.
It suffices to recall that
\[\H_{\m}^0(A/Q_{d-1})= \bigcup_{s\geq 1} (0:_{A/Q_{d-1}} {a_d}^s).\]

The last assertion is equivalent to $Q_{d-1}:a_d=\bigcup_{s\geq 1}
(Q_{d-1}:{a_d}^s)$, and thus by Proposition~\ref{acidseq} that $\xx$ is
a $d$-sequence.

\end{proof}

We are now in position to give the following converse to
Theorem~\ref{chidseq}:

\begin{Theorem}\label{chie1dminus1}
Suppose that $d \ge 2$, $\depth A \ge d-1$, and the residue class
field of $A/\m$ of $A$ is infinite. Let $Q$ be a parameter ideal in
$A$. Then the following conditions are equivalent.
\begin{enumerate}
\item[{\rm (a)}] $\chi_1(Q) \ge - \rme_1(Q)$.

\item[{\rm (b)}] $\rme_2(Q) = 0$.

 \item[{\rm (c)}] There exists a system $a_1, a_2,
\ldots, a_d$ of parameters in $A$ which generates $Q$ and forms a
$d$--sequence in $A$.
\end{enumerate}
\end{Theorem}

\subsection{Bounding Euler characteristics}

The relationship between partial Euler characteristics and
superficial elements make for a straightforward comparison with
extended degree functions.

\medskip

 We will  prove that
 Euler characteristics can be uniformly bounded by homological degrees.
The basic tool is the following observation.

\begin{Proposition}\label{basicchi1} Let $\RR$ be a Noetherian local ring
and $M$ a finitely generated $\RR$-module. If $n\geq 2$ and $\xx=\{x_1, \ldots, x_n\}$
is a system of parameters for $M$, set $\xx'=\{x_2, \ldots, x_n\}$.
Then
\[ \chi_1(\xx;M)= \chi_1(\xx';M/x_1M)+ \chi_0(\xx'; 0:_M{x_1}).\]

\end{Proposition}

This is a proved assertion  in the proof of \cite[Theorem 4.6.10(a)]{BH}.

\begin{Theorem}\label{chi1hdeg}
Let $\RR$ be a Noetherian local ring of  infinite residue field
and $M$ a finitely generated $\RR$-module of positive dimension. For every system of
parameters $\xx$ for $M$,
\[ \chi_1(\xx;M)\leq \hdeg_{(\xx)}(M)-\deg_{(\xx)}(M).\]
\end{Theorem}

\begin{proof}
We may replace the system of parameters by another generating the
same ideal but formed by a superficial sequence $\xx$ for
$\hdeg_{(\xx)}$.
\medskip

If $n=1$, $\chi_1(\xx;M)=\lambda(0:_Mx_1)$ while $\hdeg_{(\xx)}(M)=
\lambda(\rmH_{\m}^0(M))+ \deg_{(\xx)}(M)$, and therefore the assertion holds.

\medskip

If $n\geq 2$, $0:_Mx_1$ has finite length and thus $\rme_0(\xx';0
:_M x_1)=0$ so that  by Proposition~\ref{basicchi1},
\[ \chi_1(\xx;M)= \chi_1(\xx';M/x_1M).\]
We iterate this $n-1$ steps so that we have
\[\chi_1(\xx;M)= \chi_1(\xx';M/x_1M).\]
reducing to the dimension $1$ case:
\[ \chi_1(x_n;M/(x_1,\ldots, x_{n-1})M)\leq
\lambda(\rmH_{\m}^0(M/(x_1,\ldots, x_{n-1})M))\leq
\hdeg_{(\xx)}(M)-\deg_{(\xx)}(M), \]
the last inequality by
Theorem~\ref{boundtorsion}.
\end{proof}

\begin{Corollary} Let $\RR$ be a Noetherian local ring. For all systems of
parameters $\xx$ with the same integral closure $Q$, the set
\[ \Xi_1(Q)=\{\chi_1(\xx)\mid \xx\quad \mbox{\rm parameter
ideal with $\bar{(\xx)}=Q$}\}\]
is finite.
\end{Corollary}

\begin{proof} Both $\hdeg_{(\xx)}(\RR)$ and $\deg_{(\xx)}(\RR)$ depend
only on the integral closure of $(\xx)$.
\end{proof}

In analogy with the properties of Chern coefficients, the modules
with bounded Euler characteristics are those which are Cohen-Macaulay
on the punctured spectrum.
The study of these modules was started in \cite{CST} and
comprehensively developed
 in \cite{T}. Our treatment for the following result of \cite{CST}
works for the special case of equidimensional, catenary rings and was intuited from the analogous result for the
function $\rme_1(I,M)$.

\bigskip

We first give an analogue of Proposition~\ref{genhs} for Euler
characteristics.

\begin{Proposition} \label{genhschi}
Let $(\RR, \m)$ be a Noetherian local ring of infinite residue field
and dimension $d\geq 2$. Let $M$ be a finitely generated $\RR$-module
of dimension $d$ and let $\xx=\{x_1, \ldots, x_d\}$ be a system of
parameters that is filter regular.
\begin{enumerate}

\item[{\rm (a)}] $\chi_1(x_{d-1}, x_d; M/(x_1, \ldots,x_{d-2})M) =\chi_1(\xx;M)$.

\item[{\rm (b)}] Let $M_0=\H_{\m}^0(M)$ and set $M'=M/M_0$. Then
$\chi_1(M)\geq \chi_1(M')$.

\end{enumerate}

\end{Proposition}

\begin{proof} (1) follows by iteration of Proposition~\ref{basicchi1}.

\medskip

(2) The Euler characteristics satisfy
\begin{eqnarray*}
\chi_0(\xx;M_0)-\chi_0(\xx;M) + \chi_0(\xx;M') &=& 0, \\
\end{eqnarray*}
and since $\dim M_0< d$, $\chi_0(\xx;M_0)=0$.

Consider the tail of the homology exact sequence
\[ 0 \rar L \lar M_0/(\xx)M_0 \lar M/(\xx)M \lar M'/(\xx)M'
\rar 0.
\] A straightforward calculation shows that
\[ \lambda(L)-\chi_1(\xx;M')+\chi_1(\xx;M) - \chi_1(\xx;M_0)=0.\]
Thus
\begin{eqnarray*}
 \chi_1(\xx;M') -\chi_1(\xx;M) &=& \lambda(L)-\chi_1(\xx;M_0) \\
&\leq & \lambda(M_0/(\xx) M_0)- \chi_1(\xx;M_0)= \chi_0(\xx;M)= 0.
\end{eqnarray*}

\end{proof}

\begin{Theorem}\label{chi1gCM}
Let $(\RR, \m)$ be a Noetherian   local ring
and $M$ an unmixed $\RR$-module.
Then $M$ is a generalized Cohen-Macaulay module if and only if
\[ \Xi_1(M)=\{\chi_1(\xx;M)\mid \xx\quad \mbox{\rm parameter
ideal of $M$}\}
\] is  finite.
\end{Theorem}

\begin{proof} For simplicity of notation we take $M=\RR$. (Or replace
$\RR$ by $\RR/\ann(M)$.)
Assume that $\Xi_1(\RR)$ is a finite set.
Let $\dim \RR=n$. We may assume that $n\geq 2$. Let $\p\neq \m$ be
a prime ideal of positive height and let $y\in \p$ be a parameter,
that is $y$ is not contained in any  minimal prime of
$\RR$.
\medskip
\begin{itemize}
\item[{\rm (a)}] Let $r$ be such that $0:y^r=0:y^{r+1}$, and set $x_1=y^r$.
$x_1$ is a parameter in $\p$ and $0:x_1=0:x_1^2$.

\item[{\rm (b)}] By assumption, there is a bound for $\chi_1(\xx)$, for any system of parameters $\xx=x_1, x_2,
\ldots, x_n$. From Proposition~\ref{basicchi1},
\[ \chi_1(\xx)= \chi_1(\xx';\RR/x_1\RR)+ \chi_0(\xx'; 0:_\RR{x_1}),\]
which implies that $\bar{\RR}=\RR/(x_1)$ has the same property.

\item[{\rm (c)}] Furthermore, if $\dim 0:x_1=n-1$, its multiplicity is bounded
for parameters ideals of the form $x_2^m, \ldots, x_n$ and all $m\geq
1$, which is clearly impossible since $0:x_1$ has dimension $n-1\geq
1$. Therefore, $\dim 0:x_1\leq n-2$, which means that
\[ \height (0:(0:x_1))\geq 2.\]

\item[{\rm (d)}] To show that $\RR_{\p}$ Cohen-Macaulay, consider two cases. If
$\height \p=1$, $(0:(0:x_1))\not \subset \p$, which implies that
$x_1$ is a regular element of $\RR_{\p}$. Thus $\RR_{\p}$ is
Cohen-Macaulay.

\medskip
If $\height \p> 1$, by induction $\RR/(x_1)$ is a generalized
Cohen-Macaulay ring. If  $(0:(0:x_1))\not \subset \p$, as above we
would have that $x_1$ is a regular element of $\RR_{\p}$, and the proof
would be complete.

\medskip

\item[{\rm (e)}]
Assume then
 $(0:(0:x_1)) \subset \p$. Note that $0:x_1\not \subset (x_1)$ (and
 this property is preserved in $\RR_{\p}$)
as $0:x_1=0:x_1^2$, by the choice of $x_1$. This means that the ideal
$(0:(0:x_1))/(x_1)\subset \RR/(x_1) $
has positive height and is annihilated by $(0:x_1)\RR/(x_1)$. But this not
possible if $\RR/(x_1)$ is a generalized Cohen-Macaulay ring: If
$I=(a_1, \ldots, a_n)$ is an ideal of positive codimension in the
generalized Cohen--Macaulay local ring $\AA$, there is an embedding
\[ \AA/0:I \hookrightarrow \AA^n,  \quad 1\mapsto (a_1, \cdots,
a_n),\] and the associated primes of $\AA/0:I$ are minimal primes of
$\AA$ or its maximal ideal, which contradicts $(0:(0:x_1))\subset \p$.

\end{itemize}
We leave the simpler proof of the converse to the reader. An
alternative is to appeal to Theorem~\ref{chi1hdeg}, or make use of
\cite[Example 2.80]{icbook}.

\end{proof}

\subsection{A tale of four invariants} Let $\RR$ be a Noetherian
local ring with an infinite residue field and $M$ a finitely
generated $\RR$-module of dimension $n\geq 2$. Let $\mathcal{P}(M)$ be the
collection of parameter ideals $\xx=\{x_1, \ldots, x_n\}$ for $M$. In
\cite{chern3} and \cite{chern5}, the authors studied numerical
functions
\[ \ff: \mathcal{P}(M) \lar \bbz,\]
on emphasis on the nature of its range \[\XX_\ff(M)= \{\ff(\xx)\mid
\xx\in \mathcal{P}(M)\}.\]

For the two functions $\e_1(\xx,M)$ and $\chi_1(\xx;M)$, more
properly $\ff_1(\xx)=-\e_1(\xx,M) $,
and $\ff_2(\xx)=\chi_1(\xx;M)$, respectively. Identical assertions
about $\XX_\ff(M)$ are expressed in the following table:

\begin{center}
\begin{tabular}{|c|c|c|c|} \hline
\rule{0pt}{3.2ex}  $\XX_\ff(M)\subseteq [0,\infty )$   & $M$ &  \cite{MSV10}
&
  \cite[Appendice II]{Serrebook}     \\ \hline
\rule{0pt}{3.2ex}   $0\in \XX_\ff(M)$   & $M$ Cohen-Macaulay &
\cite{chern3, chern5}$^*$ &
\cite[Appendice II]{Serrebook}
     \\ \hline
\rule{0pt}{3.2ex}   $|\XX_\ff(M)| < \infty$  &  $M$ generalized
Cohen-Macaulay & \cite{chern3, chern5}$^*$ & \cite{CST, chern5}    \\ \hline
\rule{0pt}{3.2ex}   $|\XX_\ff(M)|=1$ & $M$ Buchsbaum & \cite{GO1,
chern5}$^*$ &
\cite{StuckradVogel}    \\ \hline
\rule{0pt}{3.2ex}   $|\{\ff(\xx)\mid
\bar{Q}=\bar{(\xx)}\}| < \infty$  & $Q$   &
\cite{chern5}  & \cite{chern5}  \\ \hline
\end{tabular}
\end{center}
\medskip

{\small TABLE 2. Properties of a finitely generated  module $M$ carried by the values of
either function. An adorned reference [XY]$^*$ requires that the
module $M$ be unmixed.
}

\medskip

A street view of the action of these functions is the following.
First, rewrite the system of parameters $\xx=x_1, \ldots, x_n$ so
that it is a superficial sequence. We may further require that it be
 superficial with respect to a finite set of other  modules and functions such as
$\hdeg_{(\xx)}(M) $ and $\bdeg_{(\xx)}(M)$. Now write
\[ \xx=\{x_1, \ldots, x_{n-2}, x_{n-1}, x_n\}= \{\xx',h, \alpha\}
\quad h=x_{n-1}, \alpha = x_n.\] Then
for $M'=M/(\xx')M$ the fundamental properties of $\e_1$ and $\chi_1$
(Proposition~\ref{genhs} and Proposition~\ref{basicchi1})
lead to the rules of computation
\begin{eqnarray*}
\e_1(\xx,M) & = &\e_1(h,\alpha;M')  \\
\chi_1(\xx;M) & = & \chi_1(h,\alpha;M').\\
\end{eqnarray*}
More precisely, we have
\begin{eqnarray*}
\e_1(\xx,M) & = &-\lambda(\H_{\m}^0(M'/hM'))+ \lambda(0:_{M'}h)  \\
\chi_1(\xx;M) & = &
\lambda(0:_{M'/hM'} \alpha).\\
\end{eqnarray*}

%We are now in position to give a proof of a result of \cite{chern3},
%\cite{chern5}.

%\begin{Theorem}\label{e1gCM}
%Let $(\RR, \m)$ be a Noetherian   local ring
%and $M$ an unmixed $\RR$-module.
%Then $M$ is a generalized Cohen-Macaulay module if and only if
%\[ \Lambda(M)=\{e_1(\xx;M)\mid \xx\quad \mbox{\rm parameter
%ideal of $M$}\}
%\] is  finite.
%\end{Theorem}

%\begin{proof}
%Without changing $\Lambda(M)$ or $\Xi(M)$, we may assume $\dim M=2$,
%$\xx=\alpha,h$.
%We pass to $M'=M/\H_{\m}^0(M)$
%which gives by  Propositions~\ref{genhs} and \ref{genhschi}
%\begin{eqnarray*}
%e_1(\xx;M) &=&\rme_1(\xx;M')=-\lambda(\H_{\m}^0(M'/hM')) \\
%\chi_1(\xx;M) &\leq & \chi_1(\xx;M')=\lambda(0:_{M'/hM'} \alpha).
%\end{eqnarray*}

%Altogether we have
%\[
%\chi_1(\xx;M)\leq \chi_1(\xx;M')\leq -e_1(\xx;M')=-e_1(\xx;M).
%\]
%Thus if $\Lambda(M)$ is bounded, $\Xi(M)$ is bounded as well, and by
%Theorem~\ref{chi1gCM}, $M$ is a generalized Cohen--Macaulay module.
%\end{proof}

It is at this point that cohomological degree functions play a role.
In the particular case of $\hdeg$, one has
\[ \hdeg_{(\xx)}(M') \leq \hdeg_{(\xx)}(M).\]
In dimension $2$, $\hdeg_{(\xx)}(M')$ controls the terms in the
expressions of $\e_1(\xx,M)$
and $\chi_1(\xx;M)$, which led in \cite{chern5} to bounds of the form
\[ \ff(\xx) \leq \hdeg_{(\xx)}(M)- \deg_{(\xx)}(M).\]

There is another explanation for this inequality. For that we recall
two properties of the cohomological degree $\bdeg$. If $\dim M\geq
3$, \[ \bdeg(M)\geq \bdeg(M/(\xx')M),
\] while for $L=\H_{\m}^0(M')$ and $M_0 = M'/L$
\begin{eqnarray*}
 \bdeg(M') &= & \lambda(L) + \bdeg(M_0) \\
&=& \lambda(L) + \deg(M_0) +  \lambda(\H_{\m}^0(M_0/hM_0)).\\
\end{eqnarray*}
A quick calculation will give
\begin{eqnarray*}
\bdeg(M') &=& \deg(M') + \lambda(L)
+ \lambda(\H_{\m}^0(M'/hM'))-\lambda(0:_{M'}h)\\
& = & \deg(M') + \lambda(L) - \e_1(\xx,M').
\end{eqnarray*}

From this we obtain for each of the functions $\ff(\xx)$ above that
\begin{eqnarray*}
\ff(\xx) &\leq & \bdeg_{(\xx)}(M)-\deg_{(\xx)}(M)\\
&\leq & \hdeg_{(\xx)}(M)- \deg_{(\xx)}(M),
\end{eqnarray*}
the last inequality from the minimality of $\bdeg$ among cohomological
degree functions.\index{hdeg and Chern number}\index{hdeg and Euler
number}

\begin{Theorem} Let $\RR$ be a Noetherian local ring and $M$ a
finitely generated $\RR$-module.
For any system of parameters $\xx$ of $M$,
\[ \max\{ \chi_1(\xx,M), -\e_1(\xx,M)\} \leq \hdeg_{(\xx)}(M)-
\deg_{(\xx)}(M).
\]
\end{Theorem}

A major roadblock in the use of
the functions $\e_1(\xx,M)$, $\chi_1(\xx,M)$,
$\hdeg_{(\xx)}(M)$ and $\deg_{(\xx)}(M)$,
  if $\dim M\geq 3$,
is the step requiring the reprocessing of the given system of
parameters into another that is {\em superficial} with respect to
various filtrations.
We will treat classes of parameter ideals that require minimal
reprocessing. In particular
for parameter ideals generated by $d$-sequences,
 the calculation of $\e_1(\xx,M)$ and $\chi_1(\xx,M)$ required no
 preliminary reprocessing (in the case of $\chi_1(\xx,M)$ a broader
 class of parameter ideals can be similarly treated).

\medskip

Another issue is the elucidation of the relationships between the
functions displayed
in the  diagram
\[
\diagram
 & \hdeg_{(\xx)}(M)\dto & \\
 & \bdeg_{(\xx)}(M)\dlto\drto & \\
\e_1(\xx, M)\rrto|<{\rotate\tip}^{?} & & \chi_1(\xx,M)
\enddiagram
\]
We give a comparison between $\e_1(\xx,M)$ and $\chi_1(\xx,M)$ in
cases of interest.

\medskip

\subsubsection*{Quasi--cohomological
degree}\index{quasi--cohomological degree}

\medskip

Let $(\RR,\m)$ be a   Noetherian local ring of dimension $d$ and let
$M$ be a finitely generated $\RR$-module of dimension $d$.
We are going to define a numerical function on systems of parameters
 $\xx=\{x_1, \ldots, x_d\}$ which are akin to homological degrees.

\begin{definition}{\rm A function $\hh(\xx;M):\mathcal{P}(M)\rar \bbz$ is a
{\em quasi-cohomological degree} if the following conditions are met.
\begin{enumerate}
\item[{\rm (a)}] Set $\xx'=\{x_2, \ldots, x_d\}$. If $x_1$ is a superficial
element
and  $\depth M\geq 1$,
\[ \hh(\xx;M)=\hh(\xx';M/x_1M).\]

\item[{\rm (b)}] If $M_0=\H_{\m}^0(M)$ and $M'=M/M_0$, then
\[ \hh(\xx;M)=\hh(\xx;M') + \lambda(M_0).
\]

\end{enumerate}
}\end{definition}

We now define a very special function and restrict it to some systems
of parameters.
 For any finitely generated $M$
$\RR$--module of dimension $d$, set
\[ \hh(\xx;M)= \sum_{i=0}^d (-1)^i \e_i(\xx;M).\]

\begin{Proposition} Over systems of parameters which are $d$-sequences, this
function satisfies the rules of a quasi-cohomological degree.
\end{Proposition}

%\begin{Proposition} Suppose $\xx$ is d--sequence relative to $M$.
%\begin{enumerate}
%\item[{\rm (a)}] Set $\xx'=\{x_2, \ldots, x_d\}$.
%If $\depth M\geq 1$, then
%\[ \hh(\xx;M)=\hh(\xx';M/x_1M).\]

%\item[{\rm (b)}] If $M_0=\H_{\m}^0(M)$ and $M'=M/M_0$, then
%\[ \hh(\xx;M)=\hh(\xx;M') + \lambda(M_0).
%\]

%\end{enumerate}
%\end{Proposition}

\begin{proof}
They are based on three observations about $d$--sequences.
\begin{itemize}
\item[$\bullet$] $\e_d(\xx;M)=(-1)^d \lambda(M_0)$, according to
Proposition~\ref{eiform1}. In particular $\e_d(\xx;M)=0$ if and only
if $\depth M\geq 1$.

\item[$\bullet$] To prove (a), first note that  according to
Proposition~\ref{dseqprop}, $\xx'$ is
a $d$--sequence on $M/x_1 M$ while by the previous item
$\e_d(\xx;M)=0$. Thus to prove the equality
$\hh(\xx;M)=\hh(\xx';M/x_1M)$ it suffices to show that
for all $i\geq 0$,
\[\e_i(\xx;M)=\e_i(\xx';M/x_1M).\]
\item[$\bullet$]
This will follow from Proposition~\ref{genhs} for $i<d$ once we
recognize $x_1$ as a superficial element relative to $M$ (we already
have that $\e_d(\xx;M)=0$).
For that we appeal to a general assertion of \cite{Hu81} (see
also \cite{Kuhl}) that if $\xx$ is a $d$-sequence on $M$ then the
sequence of $1$--forms of $\gr_{\xx}(\RR)$ generated by the $\bar{x_i}$ is
a $d$--sequence on the associated graded module $\gr_{\xx}(M)$.

On the other hand the approximation complex of $M$ relative to $\xx$
is exact and of length $<d$, in particular $\depth \gr_{\xx}(M)\geq 1$.
Thus the $1$--form $\bar{x_1}$ is regular on $\gr_{\xx}(M)$, which
implies that $x_1$ is a superficial  element on $M$.

\end{itemize}

\end{proof}

\begin{Corollary}
If $\xx$ is a $d$--sequence relative to $M$, then
$\hh(\xx;M)=\lambda(M/(\xx)M)$. In particular  $\hh(\xx;M)\geq
\e_0(\xx;M)$ with equality if and only if $M$ is  Cohen--Macaulay.

\end{Corollary}

\begin{proof}
This follows from the telescoping in the formulas for the
coefficients $\e_i(\xx;M)$ in Theorem~\ref{chidseq}.
\end{proof}

\begin{Corollary} If $\xx$ is a $d$--sequence relative to $M$, then
the Betti numbers $\beta_i^{\RR}(M)$ satisfy
\[ \beta_i^{\RR}(M)\leq \lambda(M/(\xx)M)\cdot \beta_i^{\RR}(k).\]
\end{Corollary}

\begin{proof}
It follows from the argument of Theorem~\ref{Degandbetti} where in the
induction part we use the properties of $\hh(\xx;M)$.

\end{proof}

\section{Generalized Hilbert Coefficients}

\subsection{Introduction}

Let $(\RR, \m)$ be a Noetherian local ring, $\AA$ a standard graded $\RR$-algebra 
 and $M$ a finitely generated graded $\AA$-module. 
 The submodule $\H=\H_{\m}^{0}(M)$ of $M$ is annihilated by a power $\m^k$ of $\m$ and therefore is a finitely generated 
 graded module over $\AA/(\m^k \AA)$,
 which led Achilles and Manaresi 
 (\cite{AMa93}) to consider the
   Hilbert function 
$ n \mapsto \lambda(\H_n)$. We focus on the  
 corresponding Hilbert series and Hilbert polynomial 
% will be still written as
%$P(M; \ttt)$ and $H(M; \ttt)$. If $\dim H = r$,
% \[ H(M;n) = \sum_{i = 0}^{r-1} (-1)^i j_i(M) {{n+ r-i-1}\choose{r-i-1}}. \]
% If $r=1$ this polynomial does not provide for $j_1(M)$, so we use instead
on  the sum transform
$n \mapsto \sum_{k\leq n} \lambda(\H_k)$. Its Hilbert polynomial
\[ \sum_{i = 0}^{r} (-1)^i j_i(M) {{n+ r-i}\choose{r-i}} \]
 will be  referred  to as the $j$-polynomial of $M$.
% 
%If
%$\H \neq 0$, for effect of comparisons if need be, we regrade $\H$ so that $\H_n=0$ for $n< 0$ and $\H_0 \neq 0$. 
Its coefficients $j_i(M)$ are integers but unlike the usual case of an Artinian local ring $\RR$ are
very hard to get hold of being  related to $M$ via  a mediating construction. 
We will only consider structures arising as follows. Let $(\RR, \m)$ be a Noetherian local ring and let $M$ be a finitely generated $\RR$-module. For an ideal $I$, the module
\[ \gr_I(M) = \bigoplus_{n\geq 0} I^nM/I^{n+1}M\]
is a finitely generated module over the associated graded ring $\gr_I(\RR)$ of $I$. 
The module $\H=\H_{\m}^0(\gr_I(M))$ is somewhat removed from $M$.
To emphasize this  point it suffices  to ask
for the dimension of $\H$ and its relationship to $M$ and $I$. 
One of the main uses of these coefficients is in the  general study of the integral closure of ideals and 
modules (add several refs).

%This fuzziness will lead us to refine the
%numbering of the  coefficients $j_i(M)$ for a special class of ideals $I$.

%In addition, some general relationships that are known to exist between the standard coefficients
%$\e_0, \e_1, \e_2$, for instance, are not known. 

%It is worth to observe that this polynomial may have degree $<r$ but the notation of $j_i(M)$

%\medskip

The algebras and modules we study here arise as follows. Let $(\RR, \m)$ be a Noetherian local ring
of dimension $d>0$ and let $M$ be a finitely
generated $\RR$-module. Let $\xx= \{x_1, \ldots, x_r\}\subset \m$ be a partial system of parameters (s.o.p. for brief) for $M$, that
 is $\dim M = r + \dim M/(\xx)M$. The associated graded module of $M$ relative to $\xx$ is
 \[ \gr_{\xx}(M) = \bigoplus_{n\geq 0} (\xx)^nM/(\xx)^{n+1}M.\] It has a natural module structure over
 $\RR[\TT_1, \ldots, \TT_r]$. The module $\H = \H_{\m}(\gr_{\xx}(M))$ is annihilated by an $\m$--primary ideal $J$ and therefore $\H$ is a finitely generated graded module over $\RR/J[\TT_1, \ldots, \TT_r]$ of dimension
 $s\leq r$. 
 
 \medskip
 
% {\bf Achtung:}  
%By abuse of notation, we shall denote the coefficients $j_i(\gr_{\xx}(M))$ by $j_i(\xx;M)$. 
%Whatever the dimension of $H$, we will write its Hilbert polynomial as a polynomial of
% pseudo degree $r$, thus $j_0(\xx;M)$ is the leading coefficient even if it vanishes. 
%We are led to this because our main questions relate to the values of $j_1(\xx;M)$.
%\medskip

A general issue is what the values of $j_1(\xx;M)$ say about $M$ itself.
In \cite{chern}, \cite{chern2}, \cite{chern3}, the authors studied the values of a special class of these coefficients.
For a local Noetherian ring $\RR$ and a finitely generated $\RR$--module $M$ we considered the Hilbert coefficients
$\e_i(\xx;M)$ associated to filtrations defined by systems of parameters $\xx$ for $M$, more precisely to the Hilbert
functions
\[ n \mapsto \lambda(M/(\xx)^{n+1}M).\]
The significant distinction between $\gr_{\xx}(M)$ and $\H_{\m}^0(\gr_{\xx}(M))$ is that the latter may not
be homogeneous and therefore the vanishing of some of its Hilbert coefficients does not place them
entirely in the context of \cite{chern}, \cite{chern2} and \cite{chern3}. 
Nevertheless we still want to examine
$j_1(\xx;M)$ in the case of $r=\dim M$
as the means to detect various properties of $M$ (e.g. Cohen-Macaulay, Buchsbaum, finite cohomology etc). 
Here we seek to extend these probes to cases when $r< \dim M$. 
We illustrate one of these issues with a series of  questions. Let $\RR$ be a Noetherian local ring and 
let $I=(x_1, \ldots, x_r)$, $r\leq d = \dim \RR$, be an ideal generated by a partial system of 
parameters, s.o.p.  for short. Let $G$ be the associated graded ring of $I$, $G= \gr_I(\RR)$. The module
$\H= \H_{\m}^0(G)$ has dimension $\leq r$.
We list some questions similar to those  raised in \cite{chern} for a [full] system of parameters:
%. The following %questions are 
%are experimental, later we will make them precise:

\begin{itemize}

\item[{\rm (i)}] What are the possible values of $\dim \H $? Note that $\H=0$ may happen or $\H\neq 0$ but of dimension zero.
%It is worth to observe that this polynomial may have degree $<r$.
% but the notation of $j_1(\xx;M)$ will
%always 
% stand for the term of degree $r-1$.
 
% \medskip
 
 %{\bf For the record should find examples of modules with $\dim H$ any number $\leq r$, for any $r$.}

\item[{\rm (ii)}] What is the signature of $j_1(\xx;\RR)$?
 If $r= \dim \H$, is $j_1(\xx;\RR) \leq 0$? The answer is affirmative if $\H$ is generated in 
 degree $0$ because $j_1(\xx;\RR) = \e_1(\TT, \H)$ and
these coefficients are always non-positive according to \cite{MV}.

\item[{\rm (iii)}] If $\RR$ is unmixed,   $r=\dim \H$ and  $j_1(\xx;\RR)=0$ is $(\xx)$  a complete intersection?
The answer is obviously no. What additional restriction is required?
%Perhaps should also require  of other
%coefficients.

%\item[{\rm (i)}]

\end{itemize}

%}\end{Example}

Let us state one of  our motivating  questions:

\begin{Conjecture}\label{j1is0}{\rm Let $M$ be an $\RR$-module satisfying condition (C). 
 Suppose $\xx = \{x_1, \ldots, x_r\}$ is a partial system of parameters for $M$. If 
 $\dim \H^0_{\m}(\gr_{\xx}(M))=r$ and $j_1(\xx;M)= 0$, then $\xx$ is a regular sequence on $M$. 

\medskip

{ What should (C) be?} At a minimum $M$ should be an ummixed module.

\medskip

 There are two versions
of the conjecture, the {\bf strong} version as stated and a {\bf weak} one. The latter assumes that
$j_1(\xx;M) = 0$ for {\bf all} partial s.o.p.'s of $r$ elements.

}\end{Conjecture}

In outline our approach describe some of the modules $\H_{\m}^0(\xx;M))$ is the following.
 For   a finitely generated $\RR$-module $M$, we call a partial s.o.p.
$\xx=\{x_1, \ldots, x_r\}$ {\em amenable} if for $i\geq 1$ the Koszul homology modules  $\H_i(\xx;M)$ have finite
support. First we establish the ubiquity of such systems that are also $d$-sequences. More precisely,
in Proposition~\ref{ubiquityth} we show that every partial s.o.p. of $r$ elements contains another partial
s.o.p.  of $r$ elements that is amenable for $M$ and generated by a $d$-sequence.
Once provided with such partial s.o.p. $\xx$, $\H_{\m}^0(\gr_{\xx}(M))$ occurs (Theorem~\ref{jdseqcx}) as the homology of the complex obtained by applying the section functor $\H_{\m}^0(\cdot)$ to
  the approximation complex $\mathcal{M}(\xx;M)$ (\cite{HSV3}), 
  \[\H_{\m}^0(\gr_{\xx}(M))= \H_{\m}^0(\mathcal{M}(\xx;M)).\] 

%\end{document}

\subsection{Amenable partial systems of parameters} \index{amenable partial system of parameters}

In this section we introduce a class of partial system of parameters that provides a gateway to the calculation
of numerous $j$-polynomials.

\begin{Definition}{\rm 
  Let $\RR$ be a Noetherian local ring and let
$M$ be a finitely generated $\RR$-module.
The partial s.o.p. $\xx=\{x_1, \dots, x_r\}$ 
 is {\em amenable} to $M$ if  
the Koszul homology modules $\H_i=\H_i(\xx;M)$, $i>0$, have finite support. 
}\end{Definition}

%\begin{Question}{\rm  What is the relationship between these modules and sequentially Cohen-Macaulay ones?
%}\end{Question}

A simple example is that of a local ring $\RR$ of dimension $2$, with $x$ satisfying
$0:x = 0:x^2$ and $x$ not contained in any associated prime $\pp$ of $\RR$ of codimension at most one.
Similarly if $M$ is Cohen-Macaulay on the punctured support of $M$, then any partial s.o.p. $\xx$ is amenable.

The main source of amenable partial s.o.p.'s is found in a property of $d$-sequences. 
Let us recall briefly their definition. We pair it the another notion.
 They are
  extensions of regular sequences, the first the notion of  {\em
  $d$-sequence} \index{d-sequence} invented by Huneke (\cite{Hu1a}),
  and of a {\em proper sequence} (\cite{HSV3II}). They play natural
  roles in the theory of the {\em approximation complexes}
  (see \cite{HSV3II}).
We now collect appropriate properties of these
sequences. (The notion extends {\it ipsis literis} to modules.)
 Let $\RR$ be a commutative ring.

  \begin{Definition} {\rm Let $\xx= \{x_1, x_2, \ldots, x_n\}$ be a sequence of
elements in $\RR$. $\xx$ is a  {\em $d$-sequence} if
\[ (x_1, x_2, \ldots, x_i): x_{i+1}x_k = (x_1, x_2, \ldots, x_i): x_k, \quad
\mbox{for } i=0,\ldots, n-1, \ k\geq i+1.\]
}\end{Definition}

\begin{Definition} {\rm Let $\xx= \{x_1, x_2, \ldots, x_n\}$ be a sequence of
elements in $\RR$. $\xx$ is a  {\em proper sequence} if
\[ x_{i+1}\H_j(x_1, x_2, \ldots, x_i)=0, \quad
\mbox{for } i=0, 1, \ldots, n-1, \ j>0,\]
where $\H_j(x_1, x_2, \ldots, x_i)$ is the Koszul homology associated to
the subsequence $\{x_1, x_2, \ldots, x_i\}$.\index{proper sequence}
}\end{Definition}

\begin{Proposition} Suppose
 that $\RR$ is a Noetherian ring.
Let   $\xx = \{x_1, x_2, \ldots, x_n\}$ be a sequence in $\RR$.
\begin{enumerate}
\item[{\rm (a)}] If $\xx$ is a d-sequence, then $\xx$ is a proper
sequence.

\item[{\rm (b)}] If $\xx$ is a proper sequence, then \[
x_{k}\H_j(x_1, x_2, \ldots, x_i)=0, \quad \mbox{for } i=0,1, \ldots,
n-1, \ j>0,\ k>i.\]
 \item[{\rm (c)}] Suppose that $\RR$ is a local
ring of dimension $n >0$. If $\xx$ is a proper sequence that is also
a system of parameters in $\RR$, then $\H_j(x_1, x_2, \ldots, x_i)$
is a module of finite length for $i=0,1, \ldots, n$ and $j>0$.
\end{enumerate}
\end{Proposition}

We make use of (c) to build amenable partial s.o.p.'s of $r$ elements. We use the proof of the 
following fact to build a system of parameters that is a $d$-sequence and whose first $r$
elements belong to $I$.

\medskip

\subsubsection*{Ubiquity of amenable sequences} The following shows the profusion of amenable partial systems of parameters.

\begin{Proposition} \label{ubiquityth}  Let $\RR$ be a Noetherian local ring, $M$ a finite $\RR$-module and $I\subset \RR$
an ideal. If $\height I\geq n$ and contains a partial s.o.p. of $M$ with $n$ elements, then $I$ contains
another partial s.o.p. of $M$ with $n$ elements that is a $d$-sequence on $M$.

\end{Proposition} 

The assertion is based on \cite[Proposition 2.7]{HSV3II}. That construction applied only to
$M=\RR$ and had a missing sentence. For completeness we write a full proof.

\bigskip

\begin{proof}
Let us recall the definition of $d$-sequences first. Let $\xx=\{x_1,\ldots, x_n\}$ be a sequence of 
elements of $\RR$ and $M$ an $\RR$-sequence. $\xx$ is a $d$-sequence relative to $M$
if \[
(x_1, \ldots, x_i)M:_M x_{i+1}x_k = (x_1, \ldots, x_i)M:_M x_k,  \textrm{for} \; i=0, \ldots, n-1, \; 
\textrm{and} \; k\geq i+1.\] 

The construction uses the following notation. Let $\{P_1, \ldots, P_a\}$ be the minimal primes of $I$.
 If $C$ is an $\RR$-module we denote by
$z(C) $ the union of its associated primes of height $\leq  n-1$.

\medskip

\begin{itemize}

\item We may assume $n>0$. Let $y_1\in I \setminus [z(R)\cup z(M)]$. Pick $m$ large enough so that 
$0:M y_1^m = 0:_M y_1^{m+1}$. We refer to this step as {\em saturation}. Set $x_1 = y_1^m$. If
$n=1$, $\{x_1\}$ will do.

\medskip

\item If $n>1$, consider the ideal $J_1 = 0:_{\RR} (0:_M x_1)$. We claim that $J_1$ is not contained in
any of the primes defining   $z(\RR)\cup z(M)$.  For such a prime $P$ we would  have:
If $x_1$ is $M_P$-regular, $0:_M x_1=0$, which is a contradiction; if $x_1$ is not $M_P$ regular, $x_1$
would be contained in some associated prime of $M_P$ which would be necessarily
contained in $P$, contradicting the choice of $y_1$ and therefore of $x_1$.

\medskip

\item Note that $J_1$ has height at least $n$: $0:_M x_1$ is submodule of $M$ therefore its associated primes (each of
which contains $x_1$), if any,
of height $\leq n-1$ would be in $z(M)$, contradicting the choice of $x_1$. Thus  the annihilator of
$0:_Mx_1$ has height at least $n$.

\medskip

\item Now pick 
\[y_2\in [I\cap J_1]\setminus [z(\RR)\cup z(\RR/(x_1))\cup z(M)\cup z(M/x_1M)]\]
and stabilize it, that is pick $m$ large enough so that 
$0:_My_2^m = 0:_My_2^{m+1}$ and $(x_1M):_M y_2^{m} = (x_1M):_M y_2^{m+1}$. Put $x_2=y_2^m$.
$\{x_1, x_2\}$ is a desired sequence.

\medskip

\item Let us consider the case $n=3$ since it involves one additional construction. Let
\[ J_2 = (x_1)M:_{\RR} (x_1)M:_M x_2,\] an ideal of height at $n$ by the argument above. Pick
\[ y_3 \in [I\cap J_1\cap J_2]\setminus[z(\RR)\cup z(\RR/(x_1))\cup z(\RR/(x_1,x_2))\cup z(M) \cup z(M/x_1M)\cup z(M/(x_1,x_2)M)]
\] and saturate it as follows: pick large $m$ so that
$0:_M y_3^m = 0:_My_3^{m+1}$,
$x_1M:_M y_3^m = x_1M:_My_3^{m+1}$ and
$(x_1,x_2)M:_M y_3^m = (x_1,x_2)M:_My_3^{m+1}$. Set $x_3 = y_3^m$. 

\medskip

\item We claim that $\xx = \{x_1, x_2, x_3\}$ is a $d$-sequence. From the construction both
$\{x_3\}$ and $\{x_1, x_3\}$ are $d$-sequences. 
Therefore the only step to check is the equality
\[ x_1M:_M x_{2} x_3 = x_1M:x_3.\]    
    Let $u$ be an element in the first module; then $x_{3} u \in x_1M:x_{2}$, and 
since $x_3 \in J_{2}$, we have $x_3^2 u\in x_1M$. The assertion now 
follows from the saturation step.

\medskip

\item In general, if $x_1, \ldots, x_s$ have been chosen and $n>s$, one picks $x_{s+1}$ in
the following manner. Again it is easy to see that 
\[ J_s = (x_1, \ldots, x_{s-1})M:_{\RR} (x_1, \ldots, x_{s-1})M:_M x_s\]
is not contained in any prime properly contained in one of the $P_i$, and that $J_s$ has height
at least $n$. Now pick
$y_{s+1}$ in
\[[I\cap J_1\cap \cdots \cap J_s] \setminus 
 [z(\RR)\cup z(\RR/(x_1))\cdots \cup z(\RR/(x_1, \ldots, x_s))
\cup z(M)\cup z(M/x_1M)\cup\cdots \cup z(M/(x_1, \ldots, x_s)M) ]
\]
and saturate it. Now this means picking $m$ large enough so that for each subsequence
$\{x_1, \ldots, x_j\}$, $j\leq s$,
= $(x_1, \ldots, x_j)M:_M y_{s+1}^m = 
(x_1, \ldots, x_j)M:_M y_{s+1}^{m+1}$.
    
\medskip    
    
\item We claim that $\{x_1, \ldots, x_n\}$ is a $d$-sequence relative to $M$. Note that the elements 
`as chosen  generate an ideal of height $n$. Note also that any
subsequence $\{x_1, \ldots, x_s, x_n\}$, $s<n-1$, is  a $d$-sequence  by
induction and the  construction. 
Therefore the only step to check is the equality
\[ (x_1, \ldots, x_{n-2})M: x_{n-1} x_n = (x_1, \ldots, x_{n-2})M:x_n.\]    
    
Let $u$ be an element in the first module; then $x_{n} u \in (x_1, \ldots, x_{n-2})M:x_{n-1}$, and 
since $x_n \in J_{n-1}$, we have $x_n^2 u\in (x_1, \ldots, x_{n-2})M$. The assertion now 
follows from the saturation step. 

%\item Hopefully QED!

%\medskip

%\item This proof has repeated steps: For instance $n=3$. Also whether we need to consider the terms
%$z(\RR/(x_1, \ldots, x_s))$ at all. To check.

\medskip

\end{itemize}
\end{proof}

If $\RR$ is a Noetherian local ring of dimension $d$, any $d$-sequence $\xx$ minimally
generated by $n$ elements satisfies $n\leq d=\dim \RR$. 
An interesting issue is when such sequences  with $n< d$ elements can be extended to a longer $d$-sequence $\yy = \{\xx, x_{n+1}\}$.

An obstruction is the fact that $d$-sequences are {\em proper} sequences and therefore satisfy:
If $\{\yy\} = \{\xx, x_{r+1},\ldots, x_n\}$ is a $d$-sequence then $(\yy)\cdot \H_i(\xx) = 0$ for $i>0$. The 
construction above has the following implication. Suppose $\xx$ is a partial sop that is a $d$-sequence. If
$\H_i(\xx)$,  for $i>0$, have finite support then $\xx$ can be extended to a sop that is a $d$-sequence.

\medskip

\subsubsection*{Amenable sequences and regular sequences}
 
The next elementary observation shows the strength of this condition to decide which partial sop are already regular
sequences.

\begin{Proposition} \label{amenablearereg}
Let $(\RR, \m)$ be a Cohen-Macaulay local ring of dimension $d$ and $M$ a finitely generated $\RR$-module
of projective dimension $s$. Let $\xx=\{x_1, \ldots, x_r\}$, $r\leq d-s$, be a partial sop for $\RR$ and $M$.
If $\xx$ is amenable for $M$ then  it is  a regular sequence on $M$.
\end{Proposition}

\begin{proof}
Let $K$ be the module of syzygies of $M$,
\[ 0 \rar K \lar F \lar M \rar 0.\]
Since $\H_i(\xx;F)=0$, $i\geq 1$, we have
\[ \H_i(\xx;K) = \H_{i+1}(\xx;M), \quad i\geq 1.\]
This means that $\xx$ is amenable for $K$, and by induction on the projective dimensions, $\xx$ is a regular
sequence on $K$.
Since $r < \depth K$, $\H_{\m}^0(K/(\xx)K) = 0$. From the
the exact sequence
\[ 0 \rar \H_1(\xx;M) \lar K/(\xx)K \lar F_0/(\xx)F_0 \lar M/(\xx)M \rar 0\]
 we have that $\H_{\m}^0(\H_1(\xx;M)) = 0$, which by assumption means that $\H_1(\xx;M)=0$.
 \end{proof}

\subsection{Explicit formulas for $j$-coefficients}

In this section 
we derive a general formula for the $j$-multiplicities of a general class of modules.
% inspired by generalized Cohen-Macaulay modules,
%more precisely modules that are Cohen-Macaulay on the punctured spectrum.  Let
%$M$ be a finitely generated $\RR$-module and $\xx=\{x_1, \dots, x_r\}$ a partial sop of $M$ that is %also
% a $d$-sequence relative
%to $M$. 
One of  our main techniques to determine $j$-polynomials is the following:

\begin{Theorem}\label{jdseqcx} Let $\RR$ be a Noetherian local ring and $M$ a finitely generated
$\RR$-module. 
Let $\xx=\{x_1, \ldots, x_r\}$ be a partial sop that is a $d$-sequence relative to $M$,
set $\H_i = \H_{\m}^0(\H_i(\xx;M))$ and 
$\SS=\RR[\TT_1, \ldots, \TT_r]$. 
There exists a complex of $\SS$-modules
\[ 0 \rar \H_r\otimes_{\RR} \SS[-r] \rar \cdots \rar \H_1\otimes_{\RR} \SS[-1] \rar \H_0 \otimes_{\RR}
 \SS \rar \H_{\m}^0(\gr_{(\xx)}(M)) \rar 0.\]
If $\xx$ is amenable for $M$ this complex is exact.
\end{Theorem}

\begin{proof}
We refer to \cite{HSV3II} for details about the approximation complexes $\mathcal{M}(\xx;M)$ used here. The 
$\mathcal{M}$-complex  is an acyclic complex of  $\SS$-modules (unadorned tensor products are over $\RR$)  
\[ 0 \rar \H_r(\xx;M)\otimes \SS[-r] \rar \cdots \rar \H_1(\xx;M)\otimes \SS[-1] \rar \H_0(\xx;M) \otimes
 \SS \rar \gr_{(\xx)}(M) \rar 0.\] 
 
 Our complex arises from applying $\H_{\m}^0$ to $\mathcal{M}(\xx;M)$. 
 Now assume that $\xx$ is amenable on $M$, so that $ \H_i(\xx;M) =
 \H_{\m}^0(\H_i(\xx;M))$, $i\geq 1$. 
We note that the image $L$ of $\H_1\otimes \SS[-1]$ in $\H_0(\xx;M) \otimes \SS$ has support in $\m$, and therefore
$\H_{\m}^1(L)=0$. Since all the $\H_i \otimes_{\RR} \SS[-i]$, $i\geq 1$, are supported in $\m$, we obtain the 
exact complex $\H_{\m}^0(\mathcal{M}(\xx;M))$
\[ 0 \rar \H_r\otimes \SS[-r] \rar \cdots \rar \H_1\otimes \SS[-1] \rar \H_{\m}^0(\H_0(\xx;M)) \otimes
 \SS \rar \H_{\m}^0(\gr_{(\xx)}(M)) \rar 0,\]
as asserted.
\end{proof}

A first application of this complex is the following criterion:

\begin{Corollary} \label{CMreqseq} $\xx$-$\depth M$ $=$ $\TT$-$\depth \H_{\m}^0(\gr_{\xx}(M))$.

\end{Corollary}

\begin{proof} By \cite[Proposition 3.8]{HSV3II},  $\xx$-$\depth M$ $=$ $\TT$-$\depth \gr_{\xx}(M)$. But this value
is also obtainable from the $\TT$-depth $L$,  a module that appears in both complexes. 
\end{proof}

%\begin{Proposition} \label{CMregseq}
%Let $M$ and $\xx$ be as above.  If $\H_{\m}^0(\gr_{\xx}(M))$
%is a maximal Cohen-Macaulay module then $\xx$ is a regular sequence on $M$.

%\end{Proposition}

%\begin{proof} We first remark on the obvious converse.
% If $\xx$ is a regular sequence on $M$, $\H_{\m}^0(\gr_{\xx}(M)) = \H_{\m}^0(\H_0(\xx;M))\otimes \SS%%$, which is either trivial or a maximal
%  Cohen-Macaulay module.
%\medskip

%$H = \H_{\m}(\gr_{\xx}(M))$ is a module over $(\RR/J)[\TT]$ with $\gr_{\TT}(H) = H$ 
%($J$ taken as an $\m$-primary ideal annihilating the coefficient modules
%of the complex). If $H$ is a maximal Cohen-Macaulay module, reducing modulo $\TT$ gives an exact %complex
%%of graded $\RR$-modules
%\[ 0 \rar \H_r[-r] \rar \cdots \rar \H_1[-1] \rar \H_{\m}^0(\H_0(\xx;M)) 
%  \rar \H_{\m}^0(\gr_{(\xx)}(M))_0 \rar 0.\]
%However exactness can only occur if $\H_r = \cdots = \H_1 = 0$. 
%\end{proof}

%{\bf We will see later that it suffices to assume $\depth \H_{m}^0(\gr_{\xx}(M))\geq r-1$.}

%\medskip

There are several other observations arising from this complex, beginning with:

\begin{Corollary}\label{jdseqCor} If $\xx$ and $M$ are as above the following assertions hold:

\begin{itemize}

%\item[{\rm (i)}] If $\xx$-$depth$ $M = s$ then $\depth$ $\H\geq r-s$. [Will delete as it is contained in previous %Cor] 

\item[{\rm (i)}]
If $\depth \H_0(\xx;M) > 0$ then $\H_{\m}^0(\gr_{\xx}(M)) = 0$ and $\xx$ is a $M$-regular sequence.

\item[{\rm (ii)}] {\bf [Hilbert series]} Set
$h_i = \lambda(\H_i(\xx;M))$, $i>0$, and $h_0 = \lambda(\H_{\m}^0(\H_0(\xx;M))$. The   Hilbert series of 
$\H_{\m}^0(\gr_{\xx}(M))$ is
\[ [\![\H]\!]= \frac{\sum_{i=0}^r (-1)^i h_i \ttt^i}{(1-\ttt)^{r+1}}, \]
and its  Hilbert polynomial  is
\[ \sum_{i=0}^r (-1)^i j_i(\xx;M) {{n+r-i}\choose{r-i}}=
 \sum_{i=0}^r (-1)^i h_i{{n+r-i}\choose{r}}.
\]
In particular,
\begin{eqnarray*}
j_0(\xx;M) &=& \sum_{i=0}^r (-1)^i h_i,\\
j_1(\xx;M) & = & \sum_{i=1}^r(-1)^i\cdot i\cdot  h_i,\\
j_r(\xx;M) & = & (-1)^rh_r. 
\end{eqnarray*}

\item[{\rm (iii)}] $\xx$ is a regular sequence on $M$ if and only if $j_1(\xx;M) = \cdots = j_r(\xx;M) =0$.

\item[{\rm (iv)}] $\H_{\m}^1(\gr_{\xx}(M)) = \H_{\m}^1(\H_0(\xx;M))\otimes \SS$, in particular every superficial
element of $\xx$ relative to $M$ induces an injective multiplication on $\H_{\m}^1(\gr_{\xx}(M))$.

\item[{\rm (v)}] If $\dim H =r$, $j_1(\xx;M)\leq 0$. 
\end{itemize}
\end{Corollary}

\begin{proof} (iii): The vanishing of the $j_i(\xx;M)$'s, $i\geq 1$,  is equivalent to the vanishing of the $h_i$'s for $i \geq 1$.

\medskip
(iv): Start with the exact sequence
\[ 0 \rar L \lar \H_0(\xx;M)\otimes \SS \lar \gr_{\xx}(M) \rar 0\]
and take local cohomology. Since $L$ has finite support, $\H_{\m}^i(L)=0$, for $i\geq 1$, which gives
the   isomorphism. The other assertion follows because the action of $x^*$ corresponds to multiplication by a linear
form in the $\TT$'s with some invertible coefficient. 

\medskip
(v): 
$H = \H_{\m}(\gr_{\xx}(M))$ is a module over $(\RR/J)[\TT]$ with $\gr_{\TT}(H) = H$ 
($J$ taken as an $\m$-primary ideal annihilating the coefficient modules
of the complex). If
$\dim H = r$, $j_1(M) = \e_1(\TT;H)$ and therefore if $H$ is unmixed by \cite{chern3} 
its vanishing means that $H$ is Cohen-Macaulay. It will follow that $h_i=0$, $i>0$ and thus
$\xx$ is $M$-regular. {\bf However we do not know that $H$ is unmixed.}
\end{proof}

Note that if $M$ is unmixed,  the depth of $H$ is positive, so
$\dim H>0$ by the Depth Lemma, in particular $H\neq 0$.  

\medskip
 {\bf Query:} How we guarantee that $H$ is unmixed? We need a condition on $M$. A condition
on $\gr_{\xx}(M)$ is that it satisfies $S_1$ (almost sure) and therefore $H$, being a submodule,
inherits the condition. 

\begin{Example}{\rm 
Let 
  $\RR$ is a Cohen-Macaulay local ring, $F$ is a free $\RR$-module and $M$ a submodule such that
$C = F/M$ has finite length. Let $\xx=\{x_1, \ldots, x_r\} $, $r<\dim \RR$,
 a partial sop that is an amenable $d$-sequence
for $M$ such that $(\xx)C=0$.

%}\end{Example}

%\end{document}

The commutative diagram
\[
\diagram
0 \rto  & (\xx)^{n+1}M \rto\dto & (\xx)^{n+1}F \rto \dto & C_{n+1}
\rto\dto  &0 \\
0 \rto          &  (\xx)^nM \rto                 & (\xx)^nF \rto  & C_{n} \rto  &0
\enddiagram
\]
yields the exact sequence
\[ 0 \rar C_{n+1} \lar (\xx)^nM/(\xx)^{n+1}M \lar (\xx)^nF/(\xx)^{n+1}F \lar C_{n} \rar 0,\] 
as the rightmost vertical mapping is trivial.
This yields $\H_{\m}^0(\gr_{\xx}(M))_n = C_{n+1}$. 

Corollary~\ref{jdseqCor}(ii)  gives  a formula for the $j$-polynomial in terms of the
Koszul homology modules $\H_i(\xx;M)$. Let us try to extract it from the exact sequence
$0\rar M \rar F \rar C \rar 0$. We have the isomorphisms
\begin{eqnarray*}
\H_{i-1}(\xx;M) & = & \H_{i}(\xx;C) = \wedge^{i} \RR^r\otimes C, \quad i\geq 2,\\
\end{eqnarray*}
and hence
\begin{eqnarray*}
h_i(\xx;M) & = & \lambda(C) {{r}\choose{i+1}}, \quad i\geq 1
\end{eqnarray*}
along with 
 the exact sequence
  \[ 0 \rar \H_1(\xx;C) = C^r \lar \H_0(\xx;M) \lar \H_0(\xx;F) \rar \H_0(\xx;C) = C \rar 0.\]
  From this we get
 \begin{eqnarray*}
 \H_{\m}^0(\H_0(\xx;M)) &=& C^r.
 \end{eqnarray*}  
A straightforward computation gives the Hilbert series of $\H_{\m}^0(\gr_{\xx}(M))$,
\[ [\![ \H]\!]= {\frac{\sum_{i=0}^{r-1}(-1)^{i}h_i \ttt^i}{(1-\ttt)^{r+1}}} 
= \lambda(C) {\frac{\sum_{i=0}^{r-1} (-1)^{i} {{r}\choose {i+1}}\ttt^i}{(1-\ttt)^{r+1}}},\]
\begin{eqnarray*}
j_0(\xx;M) & = & \lambda(C),
\end{eqnarray*}
and therefore $\dim \H_{\m}^0(\gr_{\xx}(M)) = r$, while
\begin{eqnarray*}
j_1(\xx;M) & = & -\lambda(C).\end{eqnarray*}

Verification: Let $\ff(r) =  \sum_{i=0}^{r-1} (-1)^{i} i
 {{r}\choose {i+1}}$. Then
\begin{eqnarray*}
\ff(r+1) - \ff(r) & = & 
 \sum_{i=0}^{r} (-1)^{i} i\left({{r+1}\choose {i+1}} - {{r}\choose{i+1}}\right) \\
 &=&
\sum_{i=1}^{r} (-1)^{i} i{{r}\choose {i}} = 
 \sum_{i=1}^{r-1} (-1)^{i} r{{r-1}\choose {i}} \\
 &=& r (1-1)^{r-1}=0.
   \end{eqnarray*}

}\end{Example}

%\begin{Corollary} $\xx$-$\depth M$ $=$ $\TT$-$\depth \H_{\m}^0(\gr_{\xx}(M))$.

\begin{Corollary}\label{jrminus1} 
Let $\RR$ be a Cohen-Macaulay local ring and $M$ an $\RR$-module.
Suppose $\xx = \{x_1, \ldots, x_r\}$ is a partial s.o.p. that is an amenable  $d$-sequence relative to $M$. 
If $\dim \H_{\m}^0(\gr_{\xx}(M)) = r$ and $(\xx)$-$\depth M \geq r-1$ then $\xx$ is a regular sequence on
$M$.
\end{Corollary}

\begin{proof} The complex that defines $H = \H_{\m}^0(\gr_{\xx}(M))$ is
\[ 0 \rar \H_1(\xx;M)\otimes \SS[-1] \lar \H_0 \otimes  \SS \lar H \rar 0.\]

\begin{itemize}
\item $r=1$: $H$ is a module of Hilbert polynomial $j_0(\xx;M)(n+1)$, so all its components have length
$j_0(\xx;M)$. Since it  is generated by $H_0$ and $H_n =\TT^nH_0$, it follows that multiplication by
$\TT: H_n \rar H_{n+1}$ is surjective and therefore it is injective as well. This shows that $H$ is a
Cohen-Macaulay module. From the complex we get by reduction modulo $\TT$ that $\H_1(\xx;M)=0$. 

\medskip

\item $r>1$: Suppose $(\xx)$-$\depth M = r-1$,  that is the complex of modules of modules of dimension $r$ reduced
modulo a generic linear combination of $r-1$ forms $\TT'$ in the $\TT_i$, gives a complex
\[ 0\rar \H_1(\xx;M) \otimes 
\RR[\TT'_n][-1] \lar \H_{\m}^0(\H_0(\xx;M))\otimes \RR[\TT_n'] \lar H/(\TT')H \rar 0.\]  
This leaves the Hilbert polynomial of $H/(\TT')H$ as $j_0(\xx;M)(n+1)$. The argument of the case $r=1$ implies
that $H/(\TT')H$ is Cohen-Macaulay and thus $\H_1(\xx;M)=0$, as desired.

\end{itemize}
\end{proof}

\begin{Corollary} Let $\RR$ be a Cohen-Macaulay local ring of dimension $d\geq 3$ and let $M$ be a torsion-free $\RR$-module of projective dimension two. 
Suppose $\xx = \{x_1, \ldots, x_r\}$ is a partial s.o.p. that is an amenable  $d$-sequence relative to $M$. The following holds
\begin{enumerate}
\item If $r\leq d-2$, $\xx$ is a regular sequence on $M$.
\item If $r\geq d-1$, $j_1(\xx;M)\neq 0$.
\end{enumerate}
\end{Corollary} 

\begin{proof} The first assertion follows from Proposition~\ref{amenablearereg}.

 In the other case, if
$r=d$, $j_1(\xx;M) = \e_1(\xx;M)$ which cannot vanish since $M$ is unmixed but not Cohen-Macaulay.
If $r=d-1$ we apply Proposition~\ref{amenablearereg} to the module $K$ of syzygies to conclude that
$\xx$ is a regular sequence on $K$, and therefore $\xx$-depth $M\geq d-2$. Recalling that $\depth M = d-2$,
it follows that 
$\H_{1}(\xx;M)$ is the only nonvanishing Koszul homology module.
 Corollary~\ref{jrminus1} gives then $j_1(\xx;M) = -\lambda(\H_{1}(\xx;M))$.
 \end{proof}

Let us give another application.
 
 \medskip

\begin{Corollary} Let $M$ be a Buchsbaum module and $r$ a positive integer.  Suppose that for every
partial s.o.p. of $s\leq r$ elements $\xx$, 
 $\dim H = s$ and $j_1(\xx;M)=0$. Then $M$ satisfies the condition
$S_r$ of Serre.
\end{Corollary}

\begin{proof}
 We argue by induction on $r$. We may assume $r>1$ by Corollary~\ref{jrminus1}. By induction
on the case $r-1$, we have $\depth \H\geq r-1$. The second observation then says that $\depth \H=r$, that
completes the proof.
\end{proof}

\subsubsection*{Hyperplane sections}

Let us apply Corollary~\ref{jdseqCor}(iv) to one instance of hyperplane section.  We have the usual setup, the exact sequence
\[ 0 \rar K \lar \gr_{\xx}(M)[-1] \stackrel{x^*}{\lar} \gr_{\xx}(M) \lar \gr_{\xx'}(M') \rar 0,\]
where $\xx' $ is the image of $\xx$ in $\RR = \RR/(x)$ and $M' = M/xM$.
Note that we have $K_n=0$ for $n\gg 0$.

In the conditions of the Proposition, we have
\[ 0 \rar \H_{\m}^0(\gr_{\xx}(M))[-1]_n \stackrel{x^*}{\lar} \H_{\m}^0(\gr_{\xx}(M))_n
 \lar \H_{\m}^0(\gr_{\xx'}(M'))_n \rar 0,\]
for $n\gg 0$.

\begin{Theorem}\label{hsjpoly} Let $\RR$ be a Noetherian local ring and $M$ a finitely generated $\RR$-module.
Suppose $\xx = \{x_1, \ldots, x_r\}$ is a partial s.o.p. that is an amenable  $d$-sequence relative to $M$
and  $\xx'$ for the image of $\xx$ in $\RR'=\RR/(x_1)$ and 
$M' = M/x_1M$.
 If  the $j$-polynomial of $M$ relative to $\xx$
  has degree $r\geq 2$ then
\begin{eqnarray*} 
j_i(\xx;M) & = & j_i(\xx';M'), \quad i\leq r-2. 
\end{eqnarray*}
\end{Theorem}

\begin{proof} The only issue is to prove is that the image $x^*$ of $x$ in $\gr_{\xx}(\RR)$ is a superficial element
on $\gr_{\xx}(M)$. This is a technical point of the theory of the approximation complexes (see \cite{HSV3},
 \cite{HSV3II}
 and its references).
 Precisely, according to \cite[Theorem 12.10]{HSV3}, the images $x_1^*, \ldots, x_r^*$ of $\xx$ in $\gr_{\xx}(\RR)$ is a
$d$-sequence relative to $\gr_{\xx}(M)$.  
We must show that for \[K = 0:_{\gr_{\xx}(M)} x_1^*,\]
$K_n=0$, $n\gg 0$. Since $\xx^*$ is a $d$-sequence relative to $\gr_{\xx}(M)$, $K$ is also the annihilator of
the ideal $(\xx^*)$. On the other hand  $K_{n+1} = (\xx^*)K_n$ for large $n$, thus proving the assertion. 
\end{proof}

{\bf Query:} We need to understand the case $r=2$ in order to examine the vanishing of $j_1(\xx;M)$.  
Let us do a brief exploration.
\medskip

\begin{itemize}

\item Let $x,y$ be an amenable $d$-sequence on $M$. $M_0=0:_Mx = 0:_M (x,y)$ has finite support and 
$\xx=\{x,y\}$ is an amenable $d$-sequence on $M' = M/M_0$. The exact sequence
\[ 0 \rar M_0 \lar M \lar M'\rar 0\]
gives rise to the exact sequence
\[ 0 \rar G \lar \gr_{\xx}(M) \lar \gr_{\xx}(M') \rar 0,\]
where $G$ has finite length. This gives that if the $j$-polynomial of $M$ has degree $2$, then
$j_1(\xx;M) = j_1(\xx;M')$.

Since $\depth M'
\geq r-1$, from an argument below [will rearrange later], we have that if $j_1(\xx;M')$ vanishes, $\xx'$ is a regular sequence on, but not on $M$ itself if $M_0\neq 0$.

\medskip

\item Note that  if $M_0\neq 0$, $j_2(\xx;M)= \lambda(G)$. 

%\medskip

%These observations lead to the following: [{\bf Just $r=2$}] 

%\begin{Corollary}\label{hspoly2} Let $\xx$ and $M$ be as above. Then $\xx$ is a regular sequence
%on $M$ if and only if
%\begin{eqnarray*}
%j_i(\xx;M)=0, \quad i\geq 1.
%\end{eqnarray*}
%\end{Corollary}

%{\bf This show the necessity of requiring that when $r\geq 2$, $\depth M\geq 1$.}

\end{itemize}

Before we explore other consequences we consider special cases.
% We will assume that $M$ is unmixed, 
%so that in particular $M$ has $S_1$.

\begin{itemize}

\item {\bf Unmixed modules:}
$r=2$: Here $\H_2(\xx;M)=0$ so the complex is
\[ 0\rar \H_1(\xx;M) \otimes \SS[-1] \lar \H_{\m}^0(\H_0(\xx;M))\otimes \SS \lar \H_{\m}^0(\gr_{\xx}(M))=H \rar 0.\]  
 The Hilbert series of $H$ is 
\[ \frac{h_0 - h_1\ttt}{(1-\ttt)^3}.\] We may assume that $H\neq 0$ by part (i) of the 
corollary. 

Suppose that $j_1(\xx;M)=0$. If $\dim H=1$, $j_1(\xx;M) $ is actually its multiplicity, while $\dim H=0$ is
impossible since by the Depth Lemma it has positive depth. Thus from the complex we have
that $j_1(\xx;M) = h_1(\xx;M)$. Vanishing means that $\xx$ is $M$-regular.

\medskip

\item $r=3$: Here $\H_3(\xx;M)=0$ so the complex is
\[ 0\rar \H_2(\xx;M)\otimes\SS[-2] \lar
\H_1(\xx;M) \otimes \SS[-1] \lar \H_{\m}^0(\H_0(\xx;M))\otimes \SS \lar \H_{\m}^0(\gr_{\xx}(M))=H \rar 0.\]  
 The Hilbert series of $H$ is 
\[ \frac{h_0 - h_1\ttt + h_2\ttt^2}{(1-\ttt)^4}.\] We may assume that $\dim H\geq  1$ since
$\depth $ $H>0$.
   As before, we cannot have $\dim H=2$ and $j_1(\xx;M)=0$, so $\dim H$ is $1$ or $3$. In the latter we would be done if we could prove that $H$ is unmixed.
   
   \medskip

\item We have that $j_1(\xx;M) = h_1 - 2h_2$. 
Let us consider formulas for the $h_i(\xx)$ in the style of \cite{chern7}. Let $\xx' = \{x_1, x_2\}$. Because $\xx$ is a proper sequence, we have
\begin{eqnarray*}
h_1(\xx;M) &=& h_1(\xx';M) + \lambda ((\xx':x_3)/(\xx'))\\
h_2(\xx;M) &= & h_1(\xx';M)= \lambda  ((x_1:x_2)/(x_1))\\
j_1(\xx;M) & = &  \lambda ((\xx':x_3)/(\xx')) - \lambda((x_1:x_2)/(x_1))
\end{eqnarray*}

\end{itemize}

%\begin{itemize}

%\item[{\rm (v)}] {\bf [Explicit negativity]} $j_1(\xx;M) \leq 0$. {\rm (Explained by \cite{MV}).}
%[Not yet...]

%\item[{\rm (vi)}] {\bf [Vanishing]} $j_1(\xx;M) = 0$ if and only if $\xx$ is a regular sequence relative to $M$.
 % [Not yet...] 

%\end{itemize}

%\end{Corollary}

\begin{proof} The proof of (i) is clear.

\medskip

\begin{itemize}
\item The Hilbert series of $\H_{\m}^0(\gr_{\xx}(M))$ is read off the exact homogeneous complex 
\[ \frac{\sum_{i=0}^r (-1)^i h_i\ttt^i}{(1-\ttt)^{r+1}},\]
from which we obtain the coefficients of the $j$ polynomial.
\medskip

\item We assume that $L\neq 0$. Since $L$ is a submodule of the Cohen-Macaulay module $\H_{\m}^0(M/(\xx)M)\otimes \SS$,
$L$ is unmixed and $\dim L = r$.
  Note that $j_0(\xx;M)= h_0 - \deg(L)= h_0 - \e_0(\mathbf{T}; L)$, the equality $\deg(L) = \e_0(\mathbf{T};L)$ because 
the complex is  homogeneous over the ring $\SS$, 
more precisely, if 
\[J = \ann(\H_{\m}^0(M/(\xx)M))
\bigcap_i \ann(\H_i),\quad  i>0,\]  over ring the ring
$\RR/J [\TT_1, \ldots, \TT_r]  = \RR/J[\TT]$. 
\medskip

\item Suppose that $j_1(\xx;M)=0$. From the exact sequence
\[ 0 \rar L \rar \H_{\m}^0(M/(\xx) M) [\TT] \rar \H \rar 0,\]
we have that $\e_1(\TT;L) + j_1(\xx;M) =0$, and thus $\e_1(\TT;L)=0$. $L$ however is generated in degree $1$, so the vanishing of $\e_1(\TT;L)$ does not yet implies  Cohen-Macaulayness.
Instead we observe that 
\[ \e_1(\TT;L) = \deg(L) - \e_1(\TT;L[+1]),\]
whose vanishing--since $\e_1(\TT;L[+1])\leq 0$,  as $L[+1]$ is  homogeneous and generated in degree $0$--implies
that $\deg L=0$ as $L$, if nontrivial, has dimension $r$.

\medskip

Maybe the equation is
\[ \e_1(\TT;L) = -\deg(L) - \e_1(\TT;L[+1]),\]
which means that
\[ \e_0(\TT;L[+1]) = - \e_1(\TT;L[+1]).\]

\medskip

%\item Finally, if $L=0$ we have an exact complex of maximal Cohen-Macaulay $\RR/J[\TT]$-modules which
%\item To go from
% the vanishing to Cohen--Macaulayness  of $L$ requires
%  requires unmixedness of $L$, which is assured because $L$ is a submodule of the Cohen--Macaulay module
%  $\H^0_{\m}(\H_0(\xx;M))\otimes \SS$. With $L$ Cohen--Macaulay (\cite{chern3}),
%    reducing mod $(\TT)$ gives rise to the graded exact complex
%\[ 0 \rar \H_r[-r] \lar \cdots \lar \H_2[-2] \lar \H_1[-1] \lar L/(\TT)L = 0,\] 
%which implies
%\[ \H_r = \cdots = \H_2 = \H_1=0,\] since  the terms are concentrated in different degrees.
% \QED 

\end{itemize}
\end{proof}

\chapter{Multiplicity--Based Complexity of Derived Functors}

\section*{Introduction}

 Let $(\RR, \mathfrak{m})$ be a Noetherian local ring  and let $A$ and $B$
 be finitely generated $\RR$-modules.
Motivated by the occurrence of the derived functors of $\Hom_{\RR}(A,
\cdot)$ and $A\otimes_{\RR}\cdot$ in several constructions based on
$A$ we seek to develop gauges for the sizes for these modules. In the
case of graded modules, a rich {\bf degree} based theory   has been
developed centered on the notion of Castelnuovo regularity. It is
particularly well-suited to handle complexity properties of tensor
products and modules of homomorphisms.
\medskip

Let us recall two general questions regarding the modules
$\CC=\Hom_{\RR}(A,B)$ and $\DD=A\otimes_{\RR}B$.
\begin{itemize}
\item[{$\bullet$}] Can the [minimal]number of generators $\nu(\CC)$ be estimated
in terms of $\nu(A)$ and $\nu(B)$ and other properties of $A$ and $B$? If $\RR=\bbz$, or one of its
localizations, the answer requires information derived from the structure
theory for those modules. Since this is not available for higher
dimensional rings, we will argue that an answer requires knowledge of
the cohomology of the modules.

\item[{$\bullet$}] In contrast the [minimal] number of generators of
$\DD$ is simply $\nu(\DD)=\nu(A)\cdot \nu(B)$. What is hard about
$\DD$ is to gather information about its torsion, more precisely
about its associated primes. For consider its submodule of finite
support
\[ \H_{\m}^0(\DD)= \H_{\m}^0(A\otimes_{\RR}B),
\]
and denote its length by $\rmh^0(\DD)$. Can one estimate $\rmh^0(\DD)$ in
terms of $A$ and $B$?

\item[{$\bullet$}] Contrasting further, the associated primes of $\CC$ are
well-understood, in particular
\[ \H_{\m}^0(\CC) =\Hom_{\RR}(A, \H_{\m}^0(B))
\] whose length satisfies
\[ \rmh^0(\CC) \leq \nu(A)\cdot \rmh^0(B).
\]

\item[{$\bullet$}] Muddling the issues is how to account for the
interaction between $A$ and $B$. One attempt, that of replacing $A$
and $B$ by their direct sum $A\oplus B$,  is only a temporary fix as it
poses the question of what are the ``self--interactions'' of a module?

\item[{$\bullet$}] This brings us full circle: which properties of
$A$, $B$ and of their iteraction can we bring to the table?

\item[{$\bullet$}] We shall refer to these questions as the {\em
HomAB} and {\em TorsionAB} problems. They make sense even as {\em
pure}
questions of Homological Algebra, but we have in view {\em applied}
versions.\index{HomAB problem} \index{TorsionAB problem}

\end{itemize}

Before we discuss specific motivations linking these questions to the
other topics of these Notes, we want to highlight the importance of
deriving {\rm ring-theoretic} properties of the ring $\CC=
\Hom_{\RR}(A, A)$.
Instances of $\CC$ as a non-commutative desinguralization of
$\Spec(\RR)$ are found
in the recent literature [lookup] and $\nu(\CC)$ may play a role as
an embedding dimension. In these cases, $A$ is a Cohen--Macaulay
module but we still lack
estimates. This is obviously a stimulating question.

\medskip

Partly driven by its use in normalization processes, we seek to
develop a {\em length} based theory ground on the notion of extended
multiplicity. Starting with a fixed Deg function--and the choice will
be $\hdeg$--we examine various scenarios of the following
question: Express in terms of $\hdeg(A)$ and $\hdeg(B)$ quantities
such as
\begin{itemize}
\item[{$\bullet$}] $\nu(\Hom_{\RR}(A,B))$;
\item[{$\bullet$}] $\nu(\Hom_{\RR}(\Hom_{\RR}(A,R),R))$;
\item[{$\bullet$}] $\nu(\Ext_{\RR}^i(A,B)), \ i\geq 1$;
\item[{$\bullet$}] $h_0(A\otimes_{\RR}B)=
\lambda(\H_{\mathfrak{m}}^0(A\otimes_{\RR}B))$.
\end{itemize}

The sought-after estimations are polynomial functions on
$\hdeg(A)$ and $\hdeg(B)$,  whose
coefficients are given in terms of invariants of $\RR$. The first of
these questions was treated in \cite{Dalili} and \cite{DV2}, who
refer to it as the HomAB question.

\medskip

The {\em HomAB} question
  asks for uniform estimates for the number of generators of
 $\Hom_{\RR}(A,B)$ in terms of invariants of $\RR$, $A$ and $B$. The
 extended question asks for the estimates of the number of generators
 of $\Ext_{\RR}^{i}(A,B)$ (or of other functors).
%A variation asks the similar question with $A$ and $B$ replaced by
%functors of the category $\mathcal{C}$ of $\RR$-modules, noteworthy
%being $\Hom_{\mathcal{C}}(\Ext_{\RR}^1(A, \cdot), \Ext_{\RR}^1(B, \cdot))
%$.

\medskip

In addition to the appeal of the question in basic homological
algebra, such modules of endomorphisms appear frequently in several
constructions, particularly in the algorithms that seek the integral
closure of algebras (see \cite[Chapter 6]{icbook}).
The algorithms emloyed tend to use rounds of operations of the form:
\begin{itemize}
\item[$\bullet$] $\Hom_{\RR}(E,E)$: ring extension;
\item[$\bullet$] $A \rar \tilde{A} $: $S_2$-ification of $A$;
\item[$\bullet$] $I:J$.
\end{itemize}

 Another natural
application is for primary decomposition of ideals, given the
prevalence of computation of ideal quotients.

\medskip

The setting for these questions are modules over affine integral
domains, more particularly torsion free modules. The following
elementary observations often permits a change to more amenable
rings.

\begin{proposition}\label{changeofrings} Let $\RR$ be an affine
integral domain over the field $k$, with a Noether normalization
$\RR_0= k[z_1, \ldots, z_d]\subset R$.
Suppose $\SS$ is a hypersurface ring $\RR_0\subset \SS=
\RR_0[t]\subset \RR$
over which $\RR$ is a rational extension. If $A$ and $B$ are
torsion free $\RR$-modules, then
\begin{itemize}
\item[{\rm (a)}] $\Hom_{\RR}(A,B) = \Hom_{\SS}(A,B)$;
\item[{\rm (b)}] $\Hom_{\RR_0}(\Hom_{\RR_0}(A, \RR_0),\RR_0))$ is the
$S_2$-ification of $A$.
%\item[{\rm (i)}]
\end{itemize}
\end{proposition}

Clearly the change is not warranted in the presence of
certain conditions, such as when $A$ is a vector bundle, since this
is not passed to the ring $\RR_0$, or $\SS$. Note also that in the case
$A=B$, $\Hom_{\RR}(A,A)$ can be identified to the centralizer of
the homomorphism induced by multiplication by $t$ in the ring
$\Hom_{\RR_0}(A,A)$.

\medskip

K. Dalili  (\cite{Dalili}) was the first to
 consider various cases of the HomAB problem
 and raised a broadly based conjecture (the HomAB Conjecture)
 about the character of the bounds.
The scales he used to measure the invariants were
 the cohomological degree functions introduced in \cite{DGV}. These are
 extensions (denoted by $\Deg$) of the ordinary multiplicity
function $\deg$ but encoding also
 some of the properties of its cohomology (we will discuss them
 briefly  soon).
The bounds treated have the general format
\begin{eqnarray}\label{homabc}
 \nu(\Hom_{\RR}(A,B)) \leq f(\Deg (A), \Deg (B)),
\end{eqnarray}
where $f(x,y)$ is a quadratic polynomials whose coefficients  depend
on various invariants of $\RR$. Although it can be presented in the
form $f(\Deg(A),\Deg(B))= c(\RR)\Deg(A)\Deg(B)$, we will emphasize
polynomials in various other invariants of the modules, such as
number of generators $\nu(A)$ and the classical multiplicity
$\deg(A)$.

Our outlook here has a slightly different perspective. To obviate the
non-interchangeable
roles $A$ and $B$ play in $\Hom_{\RR}(A,B)$, we consider only the case
$E=A=B$. On  one hand, this choice enables the additional structure
of algebra  in
$\Hom_{\RR}(E,E)$. Naturally
the consideration of $E=A\oplus B$ makes for a rough equivalence.
The real distinction from \cite{Dalili} comes in the classes of
modules chosen for examination. We substitute  the modules of low
dimension treated in \cite{Dalili} by a few classes of modules of
arbitrary dimension, but that still affirm the basic conjecture
(\ref{homabc}).

The setup we employ goes back to
\cite[Theorem 1.3]{HiltonRees}:

\begin{Theorem}\label{HiltonRees} \index{Hilton--Rees, Swan theorem}
There exists a natural homomorphism of $\RR$-algebras
 \[ E^* \otimes_{\RR} E \lar \Hom_{\RR}(E,E) \lar
 \underline{\Hom}_{\RR}(\Ext_{\RR}^1(E,\cdot),\Ext_{\RR}^1(E,\cdot))   \rar 0.\]
\end{Theorem}

The analysis turns on the identification of
  the module of natural
 endomorphisms  $\underline{\Hom}_{\RR}(\Ext_{\RR}^1(E,\cdot),\Ext_{\RR}^1(E,\cdot))$.
with another module of endomorphisms. This is usually carried out
 from the analysis of duality of
\cite{AusBr}: There is a homomorphism of $\RR$-algebras
\begin{eqnarray} \label{ausbr00}
 E^* \otimes_{\RR} E \lar \Hom_{\RR}(E,E) \lar \Tor_1^{\RR}(D(E),E) \rar 0
\end{eqnarray}
where $D(E)$ is the Auslander dual of $E$ (see
Definition~\ref{ausdualdef}). Considering that the
number of generators of $E^*$ can be estimated from certain
invariants of $E$ and $\RR$ (see Theorem~\ref{Kia5.3}), the question
turns on the  understanding of $\Tor_1^{\RR}(D(E),E)$.
\medskip

 In several cases,
$\Tor_1^{\RR}(D(E), E)$ is actually identified to another endomorphism ring
$\Hom_{\RR}(C,C)$, with $C$ having much smaller support than $E$ and
explicitly related to $E$.

\medskip
To track estimates of the number of generators of $\Tor_1^{\RR}(D(E),E)$
to properties of $E$, we make use of extended degree multiplicity,
particularly those labelled $\hdeg$ and $\bdeg$ (see
Definitions~\ref{hdegdef} and \ref{bdegdef}).

\section{The HomAB Problem}\label{homabsection}

\subsection{Introduction}
We are now going to describe our results. They will be framed in
 the  different ways the analysis of $\Tor_1^{\RR}(D(E),E)$ takes place by
examining several classes of modules. Since the problem is a local
question, we may assume that $(\RR, \mathfrak{m})$ is a local
Noetherian ring of dimension $d$. We will see that
in good many cases this ends up in the consideration of Cohen-Macaulay
modules. A successful resolution of this case would have a major
impact on all these problems.
We are going to list the classes of modules considered and describe
the corresponding  estimates obtained.

%We shall now describe the contents of the paper.

The following outlines the topics discussed.

\begin{enumerate}
\item[$\bullet$] {\em Dalili's results:} We will give sketchy
descriptions of these results.

\item[$\bullet$] {\em Vector bundles}:  Those are modules free on
certain distinguished open $\XX$ sets of $\Spec \RR$. For instance, if
$\XX$ is the punctured spectrum of $\RR$, $\Tor_1^{\RR}(D(E),E)$ is a module
of finite length. When $E$ is a module of finite projective dimension,
which is free in codimension $\dim \RR-2$, this is treated in
\cite{Dalili}. This remains unresolved in general, with several
undecided examples.

\item[$\bullet$] {\em Isolated singularities}: When $(\RR, \mathfrak{m})$ is a
Gorenstein local ring essentially of finite type over a field, whose
Jacobian ideal $J$ is $\mathfrak{m}$-primary, nearly all known results
that required the modules to have finite projective resolution are
extended. For that we   introduce  a new extended degree
function that makes use of Samuel's multiplicities relative to $J$.
Thus
Theorem~\ref{MCM01} shows that if $E$ is a maximal
Cohen-Macaulay module,
 \[ \nu(\Hom_{\RR}(E,E)) \leq \deg (E) \cdot  \nu(E) + e(J)\cdot  \nu(E)^2.\]
It is then used to show that for any $\RR$-module, free in
codimension $\leq \dim \RR-2$, bounds similar to those of \cite{Dalili}
(that use $\mbox{\rm proj dim }E<\infty$) are proved.

\item[$\bullet$] {\em Modules of syzygies of near-perfect modules}: If
$M$ is a perfect module (or has this property on the punctured
spectrum of $\RR$) and $E$ is one of its modules of syzygies, we will
exhibit a natural identification $\Tor_1^{\RR}(D(E),E) = \Hom_{\RR}(M,M)$. The
difference of dimensions between $E$ and $M$ leads to several very
explicit formulas for $\nu(\Hom_{\RR}(E,E))$, e.g.
(Theorem~\ref{syzofperfect}): If
  $M$ is a
perfect module and   $E$ is its module of
$k$-syzygies ($k\geq 2$),   there exists an exact sequence
\[ 0 \rar E^*\otimes E \lar \Hom_{\RR}(E,E) \lar
\Hom_{\RR}(M,M) \rar 0.\]

\item[$\bullet$] {\em Ideal modules}: An ideal module is a torsion free
$\RR$-module $E$ such that $E^{**}$ is $\RR$-free. This is a property of
ideals of grade at least two. The quotient module $C= E^{**}/E$ has
the property that $\Tor_1^{\RR}(D(E),E) = \Hom_{\RR}(C,C)$.
One application is ({Corollary}~\ref{idealmoddm3}): Let
 $(\RR, \mathfrak{m})$ be a Gorenstein local
ring of dimension $d$ and let $E$ be a torsion free ideal module  free
in codimension $d-3$. Then
\begin{eqnarray*} \nu( \Hom_{\RR}(E,E)) \leq    (\rank(E)+ \nu(E)) (\hdeg
(E)- \deg (E)) + \rank(E)\cdot \nu(E).\end{eqnarray*}
These bounds are
useful in
the estimation of the number of generators of ideal quotients $I:J$ in
rings of polynomials.

\item[$\bullet$] {\em Higher cohomology modules:} Find
estimates for the values of $\nu(\Ext_{\RR}^i(E,F))$, $i\geq 1$. Using an
appropriate multiplicity function it is possible to reduce it to the
basic
 HomAB question. What is not entirely clear is whether there are
 intrinsic approaches.

\item[$\bullet$] {\em Multiplicity-based complexity of tensor
products: } How to measure the torsion in $E\otimes_{\RR}F$? A candidate
is the function
\[h_0(E\otimes_{\RR} F)=\lambda(\H_{\mathfrak{m}}^0(E\otimes_{\RR}F)),\]
a number that is very sensitive to numerous details of the definition
of $E$ and $F$. We give a treatment of graded modules over low
dimensional rings and vector bundles.
\end{enumerate}

\medskip

%\subsubsection*{Test problems}

\subsubsection*{Canonical module} We are going to examine some test
problems.
Let us consider in detail an important maximal Cohen-Macaulay module
associated to a (non-Gorenstein) Cohen-Macaulay ring.
Let $\RR$ be a Cohen-Macaulay local ring, with a canonical module
$\omega$. Let $a_1, \ldots, a_n$ be a  minimal set of  generators of
$\omega$. Consider the complex
\begin{eqnarray} \label{challenge0}
0 \rar \RR \stackrel{\varphi}{\lar} \omega^{\oplus n} \lar E=
\coker(\varphi) \rar 0,
\quad \varphi(1) = (a_1, \ldots, a_n).
\end{eqnarray}

\begin{Proposition} $E$ is a Cohen-Macaulay module.
\end{Proposition}

\begin{proof} Let ${\mathbf z}= z_1, \ldots, z_d$ be  a maximal regular
sequence of $\RR$. It
suffices to note that $\Tor_1^{\RR}(E, \RR/(\mathbf{z}))= 0$ since
$\omega/\mathbf{z}\omega$ is the canonical module of
$\RR/(\mathbf{z})$ and therefore it is a faithful module.
\end{proof}

\begin{Problem} \label{challenge01}{\rm
Find a bound for $\nu(\Hom_{\RR}(E,E))$. We are going to give an answer
for  a
class of rings that includes those in the linkage class of complete
intersections.
}\end{Problem}

To get an estimation for $\nu(\Hom_{\RR}(E,E))$, we begin by
applying  $\Hom(\cdot, \RR)$ to the defining exact sequence, and
obtain the exact complex
\[ 0 \rar \Hom_{\RR}(E,\RR) \lar \Hom(\omega^{\oplus n}, \RR)
\stackrel{\varphi^*}{\lar} \RR \lar
\Ext_{\RR}^1(E,\RR) \lar \Ext_{\RR}^1(\omega^{\oplus n},\RR)\rar 0.\]
 More precisely, we get
\[ 0\rar  \RR/L \lar \Ext_{\RR}^1(E,\RR) \lar  \Ext_{\RR}^1(\omega^{\oplus n},
\RR)\rar 0,
\]
where $L$ is the {\em trace} ideal of $\omega$.

Now we apply
$\Hom_{\RR}(E, \cdot)$ to the sequence defining $E$ and obtain the exact
complex
\[ 0 \rar \Hom_{\RR}(E,\RR) \lar \Hom_{\RR}(E, \omega^{\oplus r} )\lar
\Hom_{\RR}(E,E) \lar \Ext_{\RR}^1(E,\RR) \rar 0,\] since
$ \Ext_{\RR}^1(E, \omega^{\oplus r})=0$. As a consequence we have
 \begin{eqnarray*}
\nu(\Hom_{\RR}(E,E))&\leq &  \nu(\Hom_{\RR}(E,\omega^{\oplus n})) +
\nu(\Ext_{\RR}^1(E,\RR))\\
 & = & n (n-1) \deg(\RR)+1 + n\nu(\Ext_{\RR}^1(\omega,\RR)),
\end{eqnarray*} since $\Hom_{\RR}(E,\omega)$ is a maximal Cohen-Macaulay
module of rank $n-1$.

Finally, we must examine the module $\Ext_{\RR}^1(\omega, \RR) $. We now
make the assumption that $\RR$ is generically Gorenstein, for instance,
that $\RR$ is an integral domain, so that we may identify $\omega$ to
an ideal of $\RR.$ From the exact sequence
\[ 0 \rar \omega \lar \RR \lar \RR/\omega \rar 0\]
we obtain the exact complex
%\[ 0 \rar \Hom_{\RR}(\omega, \omega) \lar \Hom_{\RR}(\omega,\RR) \lar
% \Hom_{\RR}(\omega, \RR/\omega) \lar
\[0 \rar
\Ext_{\RR}^1(\omega, \RR) \lar \Ext^1_{\RR}(\omega, \RR/\omega) \rar 0
\] since $\Ext_{\RR}^i(\omega, \omega)=0$, for $i\geq 1$.

Let $z$ be a regular element of $\omega$, and consider the exact
sequence induced by multiplication by $z$,
\[ 0 \rar \omega \lar \omega \lar \omega/z\omega \rar 0.\] Applying
$\Hom_{\RR}(\cdot, \RR/\omega)$, we obtain the exact sequences
\[ 0 \rar \Hom_{\RR}(\omega/z\omega, \RR/\omega) \lar \Hom_{\RR}(\omega,
\RR/\omega) \rar 0
\] and

\[ 0 \rar \Hom_{\RR}(\omega, \RR/\omega) \lar \Ext_{\RR}^1(\omega/z\omega,
\RR/\omega) \lar \Ext_{\RR}^1(\omega, \RR/\omega)
\rar 0
\]

Not going anywhere ...

\medskip

\subsection{Auslander dual}
The primary setting  of our calculations is the following construction of Auslander
(\cite{AusBr}).

\begin{Definition} \label{ausdualdef} {\rm Let $E$ be a finitely generated $\RR$-module with a
projective presentation
\[ F_1 \stackrel{\varphi}{\lar} F_0 \lar E \rar 0.\]
The {\em Auslander dual} of $E$ is the module $D(E)=
\coker(\varphi^t)$,
\begin{eqnarray} \label{ausdual}
0\rar E^*\lar  F_0^* \stackrel{\varphi^t}{\lar} F_1^* \lar D(E) \rar 0.
\end{eqnarray}
}\end{Definition}

The module
$D(E)$ depends on the chosen presentation but it is unique up to
projective summands. In particular the values of the functors
$\Ext_{\RR}^i(D(E),\cdot)$ and $\Tor_i^{\RR}(D(E), \cdot)$, for $i\geq 1$,
are independent of the presentation. Its use here lies in the
following result (see \cite[Chapter 2]{AusBr}):

\begin{Proposition}\label{Adual}
 Let $\RR$ be a Noetherian ring and
 $E$  a finitely generated $\RR$-module. There are two exact
sequences of functors:
\begin{eqnarray} \label{adual1}
\hspace{.4in} 0 \rar \Ext_{\RR}^1(D(E),\cdot) \lar E\otimes_{\RR}\cdot \lar
\Hom_{\RR}(E^*,\cdot) \lar \Ext_{\RR}^2(D(E),\cdot)\rar 0
\end{eqnarray}
\begin{eqnarray}  \label{adual2}
\hspace{.4in}   0 \rar \Tor_2^{\RR}(D(E),\cdot) \lar E^*\otimes_{\RR}\cdot \lar
\Hom_{\RR}(E,\cdot) \lar \Tor_1^{\RR}(D(E),\cdot)\rar 0.
\end{eqnarray}
\end{Proposition}

%\begin{Question}{\rm Is there a version of Auslander duals that
%uses $\Hom_{\RR}(\cdot, \omega_{\RR})$?}
%\end{Question}

\begin{Corollary}  Let $\RR$ be a Noetherian ring and
 $E$  a finitely generated $\RR$-module and denote by $D(E)$ its
 Auslander dual. Then
\[ \nu(\Hom_{\RR}(E,E)) \leq \nu(E^*)\nu(E) + \nu(\Tor_1^{\RR}(D(E),E)).\]
\end{Corollary}

In one important case, one can zoom in further (\cite[Theorem
5.3]{Dalili}):

\begin{Theorem} \label{Kia5.3}
 Let $\RR$ be a
Gorenstein local ring of dimension $d$ and let $E$ be a finitely
generated $\RR$-module. Then
\[ \bdeg(E^*) \leq  (\deg (E) + d(d-1)/2)\hdeg (E).\] In particular,
\[ \nu(E^*) \leq  (\deg (E) + d(d-1)/2)\hdeg (E).\]
\end{Theorem}

As a
 consequence,  the focus is placed
on $\Tor_1^{\RR}(D(E),E)$, a module with many interesting properties:

\begin{Corollary} The image of $E^{*}\otimes_{\RR}E$ is a two-sided ideal
of $\Hom_{\RR}(E,E)$, so $\Tor_1^{\RR}(D(E),E)$ has a ring structure.
Moreover,
the module $E$ is projective if and only if
$\Tor_1^{\RR}(D(E),E)=0$. In particular, the support of $\Tor_1^{\RR}(D(E),E)$
determines the free locus of $E$.
\end{Corollary}

\begin{proof} The first assertion is routine, considering the actions of
$\Hom_{\RR}(E,E)$ on $E$ and on $E^*$.
The surjection of the natural mapping $E^*\otimes_{\RR} E \rar
\Hom_{\RR}(E,E)$ gives a representation of the identity of $\Hom_{\RR}(E,E)$,
\[ I = \sum_{i} f_i\otimes e_i, \quad f_i\in E^*, e_i\in E,\]
that is
\[ e = \sum_i f_i(e)e_i, \quad \forall e\in E, \]
one of the standard descriptions of finitely generated projective
modules. 
\end{proof}

%\begin{Remark}{\rm A result of \cite[Theorem 1.3]{HiltonRees} identifies
%$\Tor_1^{\RR}(D(E),E)$ with the module of natural
% endomorphisms  $\Hom_{\RR}(\Ext_{\RR}^1(E,\cdot),\Ext_{\RR}^1(E,\cdot))$.
% (We thank C. Weibel
% for bringing this result to our attention.) As a side-effect, for
% several classes of modules, we represent this module as
% $\Hom_{\RR}(C,C)$, where $C$ is derived from $E$ but has a smaller support.
%}\end{Remark}

%functors of the category $\mathcal{C}$ of $\RR$-modules, noteworthy
%being $\Hom_{\mathcal{C}}(\Ext_{\RR}^1(A, \cdot), \Ext_{\RR}^1(B, \cdot))
%$.

We single out from the examination of (\ref{ausdual}) the following
effective description of $\Hom_{\RR}(E,E)$.

\begin{Proposition} \label{ausdual2} Let $E$ be a finitely generated
$\RR$-module with a presentation
\[ \RR^m \stackrel{\varphi}{\lar} \RR^n \lar E \rar 0.\]
Choose bases $\{e_1, \ldots, e_n\}$ and $\{f_1, \ldots, f_m\}$ in
$\RR^n$ and $\RR^m$. Then $\Hom_{\RR}(E,E)$ is isomorphic to the kernel of the
induced mapping
\[\Phi: E\otimes_{\RR} (\RR^n)^t \stackrel{\II\otimes \varphi^t}{\lar}
 E\otimes (\RR^m)^t.\]
\end{Proposition}
\begin{proof} We let $\rme_i^*$ be the corresponding dual basis of $(\RR^n)^t$. An
element $z= \sum x_i \otimes e_i^*$ lies in the kernel of $\Phi$ if
$ \sum x_i \otimes e_i^* \circ \varphi^t = 0$, a condition that means
$[x_1, \ldots, x_n]\cdot \varphi^t=0$. This will imply that if we let
$a_i$ be the image of $\rme_i$ in $E$, the assignment $a_i\rar x_i$ will
define an endomorphism
$\alpha: E \rar E$.

Conversely, given such $\alpha$, $\zeta = \sum_{i} \alpha(a_i)\otimes
e_i^*\in \ker(\Phi)$. The  mappings are clearly inverses of one
another. 
\end{proof}

\begin{remark}{\rm If $\RR$ is an affine ring and $\SS\rar \RR$ is a
polynomial presentation,  we may replace $\RR$ by $\SS$.
}
\end{remark}

\subsection{Isolated singularities and vector bundles}

Let $(\RR, \mathfrak{m})$ be a Gorenstein local ring which is essentially of finite
type over the field $k$. Denote by $J$ be Jacobian ideal of $\RR$. We
will assume that $J$ is an $\mathfrak{m}$-primary ideal.
For any  maximal Cohen-Macaulay (MCM for short) module $E$, we are going to give a bound for
$\nu(\Hom_{\RR}(E,E))$ in terms of the multiplicity of $E$ and the Samuel
multiplicity of $J$.

We are going to gather information on the number of generators of
$\Tor_1^{\RR}(D(E), E)$ when $E$ is a MCM  module.
First, since $\RR$ is Gorenstein, $E^*= \Hom_{\RR}(E, R)$ is also a MCM
 module of the same multiplicity. In particular we have
\begin{eqnarray} \label{aubrb}
\nu(\Hom_{\RR}(E, E)) & \leq & \nu(E)\deg (E) + \nu(\Tor_1^{\RR}(D(E),E)).
\end{eqnarray}

Let
\[ 0 \rar L \lar F_1 \lar F_0 \lar E \rar 0\]
be a minimal projective  presentation of $E$. Since all these modules are maximal
Cohen-Macaulay and $\RR$ is Gorenstein, dualizing we get another exact
sequence of MCM modules
\[ 0 \rar E^* \lar F_0^* \lar F_1^* \lar D(E)  \rar 0.\]

We now quote in full two results \cite{Wang} that we require.

\begin{Theorem}[{\cite[Theorem 5.3]{Wang}}] \label{Wang1} Let $\RR$ be
a Cohen-Macaulay local ring of dimension $d$, essentially of finite
type over a field, and let $J$ be its Jacobian ideal. Then $J\cdot
\Ext_{\RR}^{d+1}(M, \cdot) =0$ for any finitely generated $\RR$-module $M$,
or equivalently, $J\cdot \Ext_{\RR}^{1}(M, \cdot) =0$ for any finitely
generated maximal Cohen-Macaulay $\RR$-module $M$.
\end{Theorem}

\begin{Proposition}[{\cite[Proposition 1.5]{Wang}}] \label{Wang2}
Let $\RR$ be a commutative ring, $M$ an $\RR$-module, and $x\in R$. If
$x\cdot \Ext_{\RR}^1(M, \cdot)=0$ then
$x\cdot \Tor_1^{\RR}(M, \cdot)=0$.
\end{Proposition}

%\begin{corollary}\label{Wang1Tor}
% Let $\RR$ be
%a Cohen-Macaulay local ring of dimension $d$, essentially of finite
%type over a field, and let $J$ be its Jacobian ideal. Then $J\cdot
%\Tor_{d+1}^{\RR}(M, \cdot) =0$ for any finitely generated $\RR$-module $M$.
%\end{corollary}

Since $D(E)$ is MCM and $\RR$ is Gorenstein, it is a $d+1$
syzygy for an appropriate module $L$, and thus $\Ext_{\RR}^{d+1}(L, \cdot) =
\Ext_{\RR}^1(D(E), \cdot )$, a functor that according to
 Theorem~\ref{Wang1} is annihilated by the Jacobian ideal $J$ of
 $\RR$. On the other hand, appealing to Proposition~\ref{Wang2},
 $J\cdot \Tor^{\RR}_1(D(E),\cdot )= 0$.

\medskip

We now integrate these strands. Let $\mathbf{z}= z_1, \ldots, z_r$ be
a regular sequence in $J$. To provide a modicum of generality, if $J$ has
dimension $1$, we take $r=d-1$, while if $J$ is
$\mathfrak{m}$-primary, we pick $\mathbf{z}$ to be a minimal
reduction of $J$.

Consider the exact sequence induced by multiplication by $z_1$
\[ 0 \rar E \stackrel{z_1}{\lar} E \lar E/z_1E \rar 0. \]
Since $J\cdot \Tor^{\RR}_1(D(E),\cdot )= 0$, $\Tor_1^{\RR}(D(E), E)$ embeds in
$\Tor_1^{\RR}(D(E), E/z_1E)$. Iterating we will get
\[
\Tor^{\RR}_1(D(E),E )\hookrightarrow  \Tor^{\RR}_1(D(E), E/\mathbf{z}E ).
 \]
Now we derive a bound for the number of generators of
$\Tor^{\RR}_1(D(E),E)$ in terms of the properties of
$\Tor^{\RR}_1(D(E),E/\mathbf{z}E)$.

Let
\[ \cdots \lar G_2 \lar G_1 \lar G_0 \lar D(E) \rar 0 \]
be a minimal projective resolution of $D(E)$. To tensor with
$E/\mathbf{z}E$, we first do with $\RR/(\mathbf{z})$, which gives a
minimal free resolution of $D(E)/\mathbf{z}D(E)$ over the ring
$ \overline{\RR}=\RR/(\mathbf{z})$.
We now tensor with $\overline{E}=E/(\mathbf{z})E$
 over $\overline{\RR}$. Denote by
$K$ the kernel of
\[  \overline{G_1} \otimes \overline{E} \lar
\overline{G_0} \otimes \overline{E},\]
and by $B$ the corresponding module of boundaries
\[ 0 \rar B \lar K \lar \Tor_1^{\RR}(D(E), \overline{E}) \rar 0.\]

We apply the next lemma to obtain
\[ \nu(\Tor_1^{\RR}(D(E), E)) \leq \Deg (K) \leq \Deg(\overline{G_1}
\otimes \overline{E}),\]
since $\dim \overline{\RR}\leq 1$.

\begin{Lemma} \label{Degdim1} Let $\RR$ be a local Noetherian ring of
dimension $1$ and consider the diagram of finitely  generated
$\RR$-modules
\[
\diagram
 & A \rto^{\varphi} & B \rto & O, \\
 &  & C \uto|<\hole|<<\ahook &
\enddiagram
\]
that is, $C$ is a subquotient of $A$. Then $\nu(C)\leq \Deg (A)$.
\end{Lemma}

\begin{proof} Let $D= \varphi^{-1}(C)$. Since in dimension $1$ $\Deg(D)\leq
\Deg (A)$, and as $\nu(D) \leq \Deg (D)$ always, the assertion follows.
\end{proof}

We note that $G_1 = F_0^*$ and that $\overline{G_1}\otimes
\overline{E}$ is a Cohen-Macaulay that is a homomorphic image of
$\overline{G_1}\otimes  \overline{F_0}$.
 If $\dim \overline{\RR}= 0$, it is an Artin ring of length
$e(J)$, the Samuel multiplicity of the ideal $J$. We have the bound
\[ \lambda(\overline{G_1} \otimes \overline{E}) \leq e(J)\nu(E)^2.\]

We can put together these considerations in the following.

\begin{Theorem} \label{MCM01} Let $(\RR, \mathfrak{m})$ be a Gorenstein
local ring as above, let
$E$ be a MCM module  and let $J$ be the Jacobian ideal of $\RR$.
\begin{itemize}
\item[{\rm (a)}]
 If  $J$  is
$\mathfrak{m}$-primary,
\[ \nu(\Hom_{\RR}(E,E)) \leq \deg (E) \cdot  \nu(E) + e(J)\cdot  \nu(E)^2.\]
\item[{\rm (b)}]
If  $J$ has dimension $1$, then for any regular sequence
$\mathbf{z}\subset J$
of length $\dim \RR-1$,
\[ \nu(\Hom_{\RR}(E,E)) \leq \deg (E) \nu(E) + \deg(\RR/(\mathbf{z})) \nu(E)^2.\]
\end{itemize}
\end{Theorem}

\begin{Example}{\rm Let $\RR= k[x,y,u,v]/(xy-uv)$. Then for any MCM
$\RR$-module $E$ of rank $r$,
\begin{eqnarray*}
\nu(\Hom_{\RR}(E,E)) \leq 4\cdot r^2.
\end{eqnarray*}
}\end{Example}
To prove this directly, let $z$ be an element of
$\mathfrak{m}$ such that $\deg(E/zE)=\deg E$. Consider the exact
sequence
\[ 0 \rar E \stackrel{z}{\lar} E \lar E/zE \rar 0.\]
Consider the cohomology exact sequence
\begin{small}
\[ 0 \rar \Hom_{\RR}(E,E)
 \stackrel{z}{\lar} \Hom_{\RR}(E,E)
 \lar \Hom_{\RR}(E,E/zE)
 \lar \Ext_{\RR}^1(E,E)
 \stackrel{z}{\lar} \Ext_{\RR}^1(E,E).\]
 \end{small} Since $\mathfrak{m}$ is the  Jacobian ideal,
 $z\cdot\Ext_{\RR}^1(E,E)=0$,
 and  we have exact sequence
\[ 0 \rar \Hom_{\RR}(E,E)/z\Hom_{\RR}(E,E) \lar \Hom_{\RR}(E/zE,E/zE) \lar
\Ext_{\RR}^1(E,E) \rar 0.
\]
The mid module is Cohen-Macaulay of multiplicity $2r^2$, while the
module on the right is annihilated by $\mathfrak{m}$, and therefore
its length is at most $\nu(\Hom_{\RR}(E/zE, E/zE))\leq 2r^2$.
Now we apply Proposition~\ref{hs4}(d) to get
\[ \hdeg( \Hom_{\RR}(E,E)/z\Hom_{\RR}(E,E))\leq 4r^2.\]

\begin{Example} \label{syzofk} {\rm Let $(\RR, \mathfrak{m}$) be a hypersurface
singular ring of dimension $d$. Let
\[ 0 \rar E \lar F_{d-1} \lar \cdots \lar F_1\lar F_0 \rar k \rar 0\]
be a minimal free presentation of the residue field $k$. $E$ is a MCM
of known rank and number of generators,  since the Tor algebra of $k$ is known
(according to Assmus),
\[\Tor^{\RR}(k,k) = \wedge k^{d+1}\ [z],\]
where $z$ is a variable of degree $2$. In particular, $\nu(E) =
\beta_{d}(k) = 2^d$.
}\end{Example}

We now determine $\Hom_{\RR}(E,E)$ via the sequence
\[ E^* \otimes E \lar \Hom_{\RR}(E,E) \lar \Tor_1^{\RR}(D(E),E)= \underline{\Hom}_{\RR}(\Ext^1_{\RR}(E,
\cdot), \Ext_{\RR}^1(E, \cdot)\rar 0.
\]
As $\Ext_{\RR}^1(E,\cdot)=\Ext_{\RR}^d(k,\cdot)$, and it is easy to see that
this functor is right exact and $\Ext_{\RR}^d(k,\RR)=k$, it will follow
that $\Tor_1^{\RR}(D(E),E)=k$.

Finally, to determine $E^*$, we consider the complex
\[ 0\rar E \lar F_{d-1} \lar B_{d-2} \rar 0,\]
whose cohomology gives
\[ 0\rar B_{d-2} \lar F_{d-1}^* \lar E^* \lar \Ext_{\RR}^1(B_{d-2}, \RR) =
\Ext_{\RR}^{d-1}(k,\RR)=0.\]
Thus $\nu(E^*) \leq \beta_{d-1}(k)\leq 2^d-1$. It follows that
\[\nu(\Hom_{\RR}(E,E))\leq 2^{d+1}.  \]
{\bf We will check whether $E$ has no free summand.}

\bigskip

We are now going to formulate two results of \cite{Dalili} on modules
of finite projective dimension. They are \cite[Theorem 6.5]{Dalili},
asserting that the conjectural bound (\ref{homabc}) holds for modules
which are free on the punctured spectrum of $\RR$, and \cite[Theorem
7.3]{Dalili}, where the same bound is established for torsion free
modules if $\dim \RR=4$.

\begin{Proposition} \label{TorMCM} Let $(\RR, \mathfrak{m})$ be a Gorenstein local ring
essentially of finite type over the field $k$ and denote by $J$ its
Jacobian ideal. If $J$ is $\mathfrak{m}$-primary and $i$ is a fixed
positive integer, then for every maximal
Cohen-Macaulay module $M$ and every finitely generated $\RR$-module $B$,
\[ \lambda(\Tor_i^{\RR}(M,B)) \leq \beta_i(M) \Deg(B).\]
\end{Proposition}

\begin{proof} We will argue by induction on $\dim B$, the case $\dim B=0$
being clear. Let $\H_{\mathfrak{m}}^0(B)=B_0\neq 0$; the exact
sequence
\[ 0 \rar B_0 \lar B \lar B' \rar 0\]
yields $\Deg(B) = \lambda(B_0) + \Deg(B')$. On the other hand, the
long exact homology sequence
\[ \Tor_i^{\RR}(M,B_0) \lar  \Tor_i^{\RR}(M,B)  \lar  \Tor_i^{\RR}(M,B')
\] gives
\[\lambda ( \Tor_i^{\RR}(M,B))\leq  \lambda ( \Tor_i^{\RR}(M,B_0)) +\lambda (
\Tor_i^{\RR}(M,B')).\] It suffices to consider the case of modules of
positive depth.
Let $z\in J$ be a generic hyperplane section for the module $B$:
\[ 0 \rar B \stackrel{\cdot z}{\lar} B \lar \overline{B} \rar 0.\]
By Proposition~\ref{Wang2}, $z\Tor_i^{\RR}(M,B) =0$, and therefore
\[ \Tor_i^{\RR}(M,B) \hookrightarrow \Tor_i^{\RR}(M, \overline{B}).\] Since
$\dim \overline{B}= \dim B-1$ and
$\Deg(\overline{B})\leq \Deg(B)$, the induction is complete. 
\end{proof}

We now treat the mentioned  results of \cite{Dalili}, with the
assumption of finite projective dimension removed.
Let $E$ be a finitely generated $\RR$-module with a minimal presentation
\[ 0 \rar F_{d} \rar F_{d-1} \lar \cdots \lar F_1 \lar F_0 \lar E
\rar 0,\]
$d= \dim \RR$, where $F_i$ is a free $\RR$-module for $i\leq d-1$.
$F_{d}$ is a MCM module of multiplicity controlled by $\Deg(E)$.
 Applying $\Hom_{\RR}(\cdot, \RR)$ we get a complex which is broken up into
 short exact sequences:
\[ 0 \rar Z_0 \lar F_0^* \lar B_1 \rar 0\]
\[ \vdots \]
\[ 0 \rar Z_{i} \lar F_i^* \lar B_{i+1} \rar 0\]
\[ \vdots \]
\[ 0 \rar B_{d} \lar F_d^* \lar \Ext_{\RR}^d(E,\RR) \rar 0\]
along with the exact sequences, $1\leq i \leq d-1$,
\[ 0 \rar B_i \lar Z_i \lar \Ext_{\RR}^i(E,\RR) \rar 0\] and
\[ 0 \lar \Ext_{\RR}^1(E,\RR) \lar D(E) \lar B_2 \rar 0.\]
It is the last exact sequence that will become the focus of
interest, as we need to estimate $\nu(\Tor_1^{\RR}(D(E),E))$.

\medskip

 Suppose $E$ is free in  dimension $2$. This implies that
the modules $\Ext_{\RR}^i(E,\RR)$, $i\geq 1$, have dimension at most $1$.
Tensoring the last sequence by $E$, we obtain a complex
\begin{eqnarray} \label{dim1tor}
&& \Tor_1^{\RR}(\Ext_{\RR}^1(E,\RR),E)
 \lar \Tor_1^{\RR}(D(E),E) \lar \Tor_1^{\RR}(B_2,E)
\end{eqnarray}
for which we gather information from the homology of the other exact
sequences
\[ \Tor_2^{\RR}(\Ext_{\RR}^2(E,\RR),E) \lar \Tor_1^{\RR}(B_2, E) \lar \Tor_1^{\RR}(Z_2,
E)\]
\[ \Tor_1(Z_2,E) = \Tor_2^{\RR}(B_3,E) \]
\[ \Tor_3^{\RR}(\Ext_{\RR}^3(E,R),E) \lar \Tor_2^{\RR}(B_3, E) \lar
\Tor_2^{\RR}(Z_3,E) \]
\[ \Tor_2(Z_3,E) = \Tor_3^{\RR}(B_4,E) \]
\[ \vdots \]
\[ \Tor_{d}^{\RR}(\Ext^d_{\RR}(E,\RR),E)\lar \Tor_{d-1}^{\RR}(B_d, E)\lar
\Tor_{d-1}^{\RR}(F_d^*,E).
\]

For each $i\geq 1$,
$\Tor_i^{\RR}(\Ext_{\RR}^i(E,\RR),E)$ is a subquotient of $\RR^{\beta_i(E)}\otimes
\Ext_{\RR}^i(E,\RR)$, a module whose $\Deg$ value is
$\beta_i(E)\Deg(\Ext_{\RR}^i(E,\RR))$.
 By
Lemma~\ref{Degdim1}, $\dim \Ext_{\RR}^i(E,\RR)\leq 1$, and thus  every submodule of
$\Tor_i^{\RR}(\Ext_{\RR}^i(E,\RR),E)$ is generated by
 not more than $\beta_i(E)\Deg(\Ext_{\RR}^i(E,\RR))$ elements.

Finally, note that $\Tor_{d-1}^{\RR}(B_d, E)$ must be handled differently.
To estimate a bound for the number of generators for its submodules:
to $\beta_{d}(E)\Deg(\Ext_{\RR}^d(E,\RR))$ we must add
the bound for $\nu(\Tor_{d-1}^{\RR}(F_d^*,E))$, derived in
Proposition~\ref{TorMCM}, since $F_d^*$ is a MCM module.

\begin{Theorem} \label{freed2} Let $(\RR, \mathfrak{m})$ be a Gorenstein local ring
essentially of finite type over the field $k$ and denote by $J$ its
Jacobian ideal.
If $J$ is $\mathfrak{m}$-primary, then for any $\RR$-module that is free
in dimension $d-2$,
\begin{eqnarray*} \nu(\Hom_{\RR}(E,E)) &\leq &\nu(E)(\deg E +
d(d-1)/2)\hdeg (E)+
 \hdeg(E)\sum_{i=1}^d \beta_i(E) \\ &+&
\beta_{d-1}(k)\beta_{d}(E)\deg (E)\hdeg_J(E).
\end{eqnarray*}
\end{Theorem}

\begin{proof} The first summand corresponds to the number of generators for
$E^*\otimes E$, according to Theorem~\ref{Kia5.3}.
 The second arises from the consideration of short
exact sequences of modules of dimension at most $1$,
\[ A \lar B \lar C,\]
where
the bounds $a$ and $c$ for the numbers of generators of $A$ and $C$,
respectively, we get that $a+c$ bounds the number of generators for
any submodule of $B$.

To obtain the last term, we apply Proposition~\ref{TorMCM} to the module
$F_d^*$, a Cohen-Macaulay module with $\deg(F_d)\leq
\beta_{d}(E)\deg(\RR)$.
\end{proof}

\begin{Remark}{\rm
Now we outline quickly the proof of the conjecture for rings of
dimension $4$ and torsion free modules. (More generally, for rings of
arbitrary dimension and torsion free modules  with
the condition $S_{d-3}$.)  These modules  have a presentation
\[ 0 \rar F_3 \lar F_2 \lar F_1 \lar F_0 \lar E\rar 0, \] where $F_i$
are free, for $i\leq 2$, and $F_3$ is a MCM module.
The only significant deviation from the proof above, in the exact
sequence
\[ \Tor_1^{\RR}(\Ext_{\RR}^1(E,\RR),E) \lar \Tor_1^{\RR}(B_2, E) \lar \Tor_1^{\RR}(Z_2,
E),\]
we must deal with the fact that $\Ext_1^{\RR}(E,\RR)$ has dimension at most
two,   the inequality
 $\hdeg (\Ext_{\RR}^1(E,\RR))< \hdeg(E)$, and all the other homology modules having dimension
 at most $1$. This is handled in \cite{Dalili}
to give
$\nu(\Tor_1^{\RR}(\Ext_{\RR}^1(E,\RR),E))\leq \hdeg(\Ext_{\RR}^1(E,R))\nu(E)$.

}\end{Remark}

%\newpage

\subsubsection*{Isolated singularities}
We are going to relax the conditions on $M$: Assume that it has
finite projective dimension, has
codimension $n$ and is Cohen-Macaulay on the punctured spectrum.
Consider its
 minimal free resolution:

\[ \KK: \quad 0 \rar K_m \rar \cdots \rar K_n\rar K_{n-1} \rar K_{n-2} \rar \cdots
\rar K_2\rar K_1\rar K_0 \rar M \rar 0.\]
The conditions on $M$ imply that $\Ext_{\RR}^i(M,\RR)=0$ for $i<n$, and
$\Ext_{\RR}^i(M,\RR)$ has finite support for $i>n$.  The latter means that
the syzygies of $M$ of order $n$, or higher, are vector bundles on the
punctured spectrum.

Let $E$ be the module of $k$-cycles of $\KK$:
\[ 0 \rar K_m \rar \cdots \rar K_{k+1} \stackrel{\varphi}{\rar}
K_{k} \rar E \rar 0,\]
\[ 0 \rar E \rar K_{k-1} \rar \cdots \rar K_1 \rar K_0 \rar M \rar 0.\]
If $k\geq n$, $E$ being a
vector bundle of finite projective dimension, the module $\Hom_{\RR}(E,E)$ was
dealt with in \cite{Dalili}. For this reason we will focus on the
case $k<n$. Denote the dual of $\KK$ by $\LL$. We will set $\LL_i=
\KK_{m-i}^*$. We have the complexes
\[ 0 \rar K_0^* \rar K_1^* \rar \cdots \rar K_{k-1}^* \rar E^* \rar 0,\]
\[ 0 \rar E^* \rar K_{k}^* \stackrel{\varphi^t}{\rar} K_{k+1}^*
\rar \cdots \rar K_m^* \rar 0.
 \]
The Auslander dual $D(E)$ is expressed in the two exact sequences
\[
 0 \rar E^* \rar K_{k}^* \stackrel{\varphi^t}{\rar} K_{k+1}^* \rar
 D(E) \rar 0,\]
\[ 0 \rar  \Ext_{\RR}^{k}(M,\RR) \rar D(E) \rar Z_{m-k-2}(\LL) \rar
\Ext_{\RR}^{k+2}(M,\RR)\rar 0,
\] which we will make use of to calculate $\Tor_1^{\RR}(D(E),E)$.

Let us examine the different cases. If $k+2<n$,
$\Ext_{\RR}^k(M,\RR)=\Ext_{\RR}^{k+2}(M,\RR)=0,$ and from the exact complex
\[ 0 \rar Z_{m-k-2} \rar L_{m-k-2} \rar \cdots \rar L_{m-n} \rar
 B_{m-n} \rar 0\]
we get
\[ \Tor_1^{\RR}(D(E), E)= \Tor_1^{\RR}(Z_{m-k-2}, E)= \Tor_{n-k}^{\RR}(B,E)=
\Tor_n^{\RR}(B,M),
\] the last isomorphism since $E$ is the module of $k$-cycles of the
resolution of $M$.

\medskip

If $M$ is Cohen-Macaulay, $B=\Ext_{\RR}^n(M,\RR)$, and we would get the
formula of Theorem~\ref{syzofperfect}. In the general case, we have a
series of short exact sequences derived from the boundaries $B_i$ and
cycles $Z_i$ of the complex $\LL$: For $i\geq 1$
\[ 0 \rar \Ext_{\RR}^{n}(M,\RR) \rar B \rar B_{m-n-1} \rar 0, \]
\[ 0 \rar B_i \rar Z_i \rar \Ext_{\RR}^{n+i}(M,\RR) \rar 0, \]
\[ 0 \rar Z_{i+1} \rar L_{i+1} \rar B_i \rar 0 .\]

From  these, we would get the exact complexes whose
significant parts we highlight:
\begin{small}
\[\Tor_n^{\RR}(\Ext_{\RR}^n(M,\RR), M) =
\Tor_{n-k}^{\RR}(\Ext_{\RR}^n(M,\RR), E)\rar
\Tor_{n-k}^{\RR}(B,E) \rar
\Tor_{n-k}^{\RR}(B_{m-n-1}, E) ,
\]
\[
 \Tor_{n-k}^{\RR}(\Ext_{\RR}^{n+i}(M,\RR), E)
  \rar   \Tor_{n-k}^{\RR}(B_i E) \rar  \Tor_{n-k}^{\RR}(Z_i,E) \rar
  \Tor_{n-k-1}^{\RR}( \Ext_{\RR}^{n+i}(M,R),E) , \]
\[    \Tor_{n-k}^{\RR}(B_i,E)= \Tor_{n-k-1}^{\RR}(Z_{i+1} , E) .\]
\end{small}
All the terms that involve $\Ext_{\RR}^{n+i}(M,\RR)$ have finite length
that can be assembled in a way bounded by  a quadratic polynomial on
$\Deg (E)$.

\subsection{Ideal modules}

Let $(\RR, \mathfrak{m})$ be a Cohen-Macaulay local ring of dimension
$d$ and let $E$ be a finitely generated torsion free $\RR$-module. $E$
is said  be an {\em ideal module}
 if $E^*$ is a
free $\RR$-module (see \cite{ram1}). Such modules are basically the
first order syzygies of modules
of codimension at least $2$.
This property endows $E$ with a natural embedding into a free
$\RR$-module
\begin{eqnarray}&& 0 \rar E \lar (E^*)^* \lar C \rar 0.
\label{idealmod2} \end{eqnarray}
Another property of these modules that is relevant in the HomAB
problem   is that the   Auslander dual  $D(E)$ has
projective dimension at most two: Dualizing (\ref{ausdual}) gives
\[0\rar E^*\lar  F_0^* \stackrel{\varphi^t}{\lar} F_1^* \lar D(E)
\rar 0.\]

\begin{Proposition} Let $E$ be an ideal module of rank $\ell$ as
above. There exists an exact sequence
\[ 0 \rar R^{\ell} \otimes_{\RR} E \lar \Hom_{\RR}(E,E) \lar \Hom_{\RR}(C,C) \rar
0. \]
\end{Proposition}

\begin{proof} Apply $\Hom_{\RR}(\cdot, E)$ to the exact sequence (\ref{idealmod2})
(note $(E^*)^*\simeq R^{\ell}$) to get the exact sequence
\[ 0 \rar \Hom_{\RR}(\RR^{\ell}, E) \lar \Hom_{\RR}(E,E) \lar \Ext_{\RR}^1(C,E) \lar
\Ext_{\RR}^1(\RR^{\ell},E)=0.\]

Since $E$ is free in codimension $1$, the annihilator $I$ of $C$ has
codimension at least $2$. Let $f$ be a regular element in $I$;
reducing (\ref{idealmod2}) modulo $f$, we get the exact complex
\[ 0 \rar C \lar E/fE \lar R^{\ell}/fR^{\ell} \lar C \rar 0,\]
that identifies $C$ to the submodule of $E/fE$ supported in dimension
at most $d-2$. If we take this fact into
 the isomorphism $\Ext_{\RR}^1(C,E)\simeq  \Hom_{\RR}(C,E/fE)$, we obtain
$\Ext_{\RR}^1(C,E) = \Hom_{\RR}(C,C)$, as desired. 
\end{proof}

We can relate $\Deg (C)$ to $\Deg (E)$, at least in the case when $\Deg=
\hdeg$. Considering that a bound for $\nu(\Hom_{\RR}(C,C))$ is established
in \cite{Dalili} for all modules of dimension at most $2$, we can
apply it to ideal modules that are free in dimension at most $3$. The
formula obtained is similar to the one we are going to derive now by
 treating directly
 the support  $\Tor_1^{\RR}(D(E),E)$.

\begin{Proposition}
 Let $(\RR, \mathfrak{m})$ be a Cohen-Macaulay local
ring of dimension $d$ and let $E$ be a torsion free ideal module that
is free in dimension $2$. Then
\[\hdeg (\Tor_1^{\RR}(D(E),E))\leq
\rank(E)\cdot (\hdeg (E)- \deg (E)) .\]
\end{Proposition}

\begin{proof} We begin by pointing out  that the support of  $\Tor_1^{\RR}(D(E),E)$ has
dimension at most $1$ since $E$ is free in codimension $d-2$. By the
same token,
in  the natural embedding
\[ 0 \rar E \lar E^{**}\lar C \rar 0 \]
the module $C$ has dimension at most $1$. Since $E^{**}$ is free, one
has $\hdeg (C) \leq \hdeg (E) - \deg (E)$.
We also have that $\Tor_2^{\RR}(D(E),C) = \Tor_1^{\RR}(D(E),E)$.
The first of these  modules is the homology of the mapping of modules of
dimension at most $1$,
\[ 0 \rar \Tor_2^{\RR}(D(E),C) \lar  E^{*}\otimes C \lar F_0^* \otimes C,\]
and therefore \[\hdeg (\Tor_2^{\RR}(D(E), C))\leq \rank(E)\cdot  \hdeg(C) =
\rank(E)\cdot (\hdeg (E) - \deg (E)).\]
\end{proof}

\begin{Corollary} Let $(\RR, \mathfrak{m})$ be a Cohen-Macaulay local
ring of dimension $d$ and let $E$ be a torsion free ideal module that
is free in dimension $2$. Then
\[ \nu( \Hom_{\RR}(E,E)) \leq  \rank(E)\cdot (\hdeg (E) + \nu(E)-\deg (E)). \]
\end{Corollary}

One should attempt to cut down  the free locus to codimension $d-3$.
It would be of interest  if we could still have $\dim \Tor_1^{\RR}(D(E),E)\leq
1$, for even though $C$ could have dimension $2$ its component $C_1$
of dimension at most $1$ would satisfy $\hdeg (C_1) \leq \hdeg (E)-\deg
(E)$
(a fact not difficult to verify).

\bigskip

Let us  go through the calculations in the case of  free locus in
codimension $d-3$.  Now the module $C$ in
\[0 \rar E \lar E^{**} \lar C \rar 0\]
has dimension at most $2$. Let $z$ be a generic element for the
following purposes: $\hdeg (\RR)
\geq \hdeg (\RR/(z))$, $\hdeg (E)\geq \hdeg (E/zE)$, $\dim (C/zC) \leq 1$,
\[\dim \Tor_1^{\RR}(D(E), E)/z\Tor_1^{\RR}(D(E), E)\leq 1,\] and
the kernels ${}_zC$ and $L$ of multiplication by $z$ on $C$ and
$D(E)\otimes E$ have
finite length. In addition, since we may assume $\dim R\geq 4$ and
$D(E)$ has projective dimension at most $2$, $z$ can be assumed to be
regular on $D(E)$.

Let us see how reduction modulo $(z)$ affects the degrees data.
Tensoring
\[ 0 \rar E \stackrel{\cdot z}{\lar} E \lar \overline{E} \rar 0\]
by $D(E)$, we have the exact sequence
\[ 0 \rar  \Tor_1^{\RR}(D(E), E)/z\Tor_1^{\RR}(D(E), E)\lar \Tor_1^{\RR}(D(E),
\overline{E}) \lar L \rar 0. \]
This shows that $\Tor_1^{\RR}(D(E),\overline{E})$ has dimension at most
one, and by Nakayama Lemma
\begin{eqnarray} \label{eq1a}
 \nu (\Tor_1^{\RR}(D(E), E)) &\leq& \hdeg (\Tor_1^{\RR}(D(E), \overline{E})).
\end{eqnarray}

To get hold of the last degree, we tensor the defining sequence of $C$ by
$\RR/(z)$ to get
\[ 0 \rar {}_zC \lar \overline{E} \lar \overline{E^{**}} \lar
\overline{C} \rar 0.\]
It shows that ${}_zC$ is the submodule of finite support of $\overline{E}$
and therefore the image $E'$ of $\overline{E}$ satisfies $\hdeg (E')\leq
\hdeg (\overline{E})\leq \hdeg (E)$, and from the exact sequence
\[ 0 \rar E' \lar \overline{E^{**}} \lar \overline{C} \rar 0\]
we have \[\hdeg ( \overline{C})\leq \hdeg (E') - \deg (E')\leq \hdeg (E) - \deg
(E).\]

Now we get hold of $\hdeg (\Tor_1^{\RR}(D(E), \overline{E}))$. From the
exact sequence
\[ 0 \rar {}_zC \lar \overline{E} \lar E' \rar 0\]
we have the exact sequence
\[
\Tor_1^{\RR}(D(E), {}_zC) \lar \Tor_1^{\RR}(D(E), \overline{E}) \lar
\Tor_1^{\RR}(D(E), E').\]
It shows that
\begin{eqnarray}
\hdeg (\Tor_1^{\RR}(D(E)), \overline{E}) &\leq & \hdeg (\Tor_1^{\RR}(D(E), E')) +
\hdeg (\Tor_1^{\RR}(D(E)), {}_zC) \nonumber  \\
&\leq & \hdeg (\Tor_1^{\RR}(D(E), E')) + \rank(F_0)\cdot
\lambda({}_zC)\nonumber \\
& \leq &
\hdeg ( \Tor_1^{\RR}(D(E), E')) + \rank(F_0)\cdot (\hdeg(E)-\deg
(E)).\label{eq2a}
\end{eqnarray} The last equation arises because  $\hdeg (\overline{E})= \hdeg
({}_zC) + \hdeg (E')$ and therefore
\[ \hdeg ({}_zC) \leq \hdeg (\overline{E})-\hdeg (E') \leq \hdeg (E) -
\deg (E). \]
%\Tor_1^{\RR}(D(E), {}_zC)
%\Tor_1^{\RR}(D(E), \overline{E})
%\Tor_1^{\RR}(D(E), E')

We must now deal with $\hdeg (\Tor_1^{\RR}(D(E), E'))$ in the manner of the
earlier case. Starting from
\[ 0 \rar E' \lar  \overline{E^{**}} \lar \overline{C} \rar 0,\]
we have the exact sequence
\[ \Tor_2^{\RR}(D(E), \overline{E^{**}}) \lar \Tor_2^{\RR}(D(E),
\overline{C}) \lar \Tor_1^{\RR}(D(E), E') \lar \Tor_1^{\RR}(D(E),
\overline{E^{**}}).
\]
Since $E^{**}$ is $\RR$-free and $z$ is regular on $D(E)$, the two
modules at the ends vanish and we have $\Tor_1^{\RR}(D(E), E') =
\Tor_2^{\RR}(D(E), \overline{C})$.  As $\overline{C}$ has dimension at
most $1$ and we have bounds for $\hdeg (\overline{C})$, we make use of
the free presentation of $D(E)$ and finally obtain
\begin{eqnarray}
 \hdeg (\Tor_1^{\RR}(D(E), E')) &=& \hdeg (\Tor_2^{\RR}(D(E), \overline{C} ))
 \leq  \rank(E)\cdot \hdeg
(\overline{C}) \nonumber \\
& \leq & \rank(E)\cdot (\hdeg (E)-\deg (E)). \label{eq3a}
\end{eqnarray}

We collect the calculation into:

\begin{Proposition}
 Let $(\RR, \mathfrak{m})$ be a Gorenstein local
ring of dimension $d$ and let $E$ be a torsion free ideal module that
is free in dimension $3$. Then
\[\nu( \Tor_1^{\RR}(D(E),E))\leq
(\rank(E)+ \nu(E)) (\hdeg (E)- \deg (E)).\]
\end{Proposition}

\begin{proof} We have from (\ref{eq1a})
\[
\nu(\Tor_1^{\RR}(D(E),E)) \leq \hdeg (\Tor_1^{\RR}(D(E), \overline{E})),
\] while from (\ref{eq2a}) we have
\[ \hdeg (\Tor_1^{\RR}(D(E), \overline{E}))\leq
 \hdeg (\Tor_1^{\RR}(D(E), E')) + \rank(F_0)\cdot \hdeg(E).\]
Finally, from (\ref{eq3a}) we obtain
\[ \hdeg (\Tor_1^{\RR}(D(E), \overline{E}))\leq (\rank(E)+ \nu(E)) (\hdeg
(E)-\deg
(E)).\]
%Hopefully these substitutions are right!
\end{proof}

\begin{Corollary} \label{idealmoddm3} Let $(\RR, \mathfrak{m})$ be a Gorenstein local
ring of dimension $d$ and let $E$ be a torsion free ideal module  free
in codimension $d-3$. Then
\[ \nu( \Hom_{\RR}(E,E)) \leq    (\rank(E)+ \nu(E)) (\hdeg
(E)- \deg (E)) + \rank(E)\cdot \nu(E).\]
\end{Corollary}

%\newpage

\subsubsection*{Primary decomposition}

Primary decomposition of ideals in a  Noetherian ring
$\RR$, and several other operations as well, involve the computation of
ideal quotients: $I\colon J=\{r\in R\mid rI\subset J\}$.
It may be relevant to have an estimation for the number of generators
of $I\colon J$ given $I$ and $J$.

The typical ring for us will be a polynomial ring over a field, but given the
well-established technique (\cite{Swan}) to convert local bounds for
number of generators into global ones, we focus on  Gorenstein
local rings  and on ideals of codimension at least two.

We may consider $I:J=\Hom_{\RR}(I,J)$ directly or the variant
$\Hom_{\RR}(E,E)$ for $E= J \oplus (I,J)$. Then, since by assumption $\codim
J\geq 2$,
\[ \Hom_{\RR}(E,E) = (J:J)\oplus J:(I,J) \oplus (I,J):(I,J) \oplus
(I,J):J)= \RR^3\oplus I:J,\]
since $(I,J):J= I:J$. Note that $E$ is an ideal module of rank $2$.

\medskip

As a direct consequence of Corollary~\ref{idealmoddm3}, we obtain:

\begin{Theorem}
Let $(\RR, \mathfrak{m})$ be a Gorenstein local
ring of dimension $d$ and let $J$ be an ideal of codimension at least
$d-2$. Then for any ideal $I$,
\[ \nu( I:J) \leq    (2+ \nu(J)+ \nu(I,J)) (\hdeg (J) + \hdeg(I,J)
- 2\deg (\RR)) + 2\cdot (\nu(J)+ \nu(I,J)) - 3.\]
\end{Theorem}

\begin{Remark}{\rm
 We have to calculate $\hdeg (\Ext_{\RR}^i(\RR/I, \RR))$ for
$i=d-2, d-1, d $. If $\codim J = d-2$, $\hdeg (
\Ext_{\RR}^{d-2}(\RR/J, \RR))=\deg
(\RR/J)$, and $\hdeg (\Ext_{\RR}^d(\RR/J, \RR))= \lambda(\H^0_{\mathfrak{m}}(\RR/J))$.

\medskip

Let us consider some cases in detail. If $J$ is
$\mathfrak{m}$-primary, $\Hom_{\RR}(\RR/I, \RR/J)= I:J/J$ is a module whose length is
at most $\lambda(\RR/J)$, and thus $ \nu(I:J)\leq
\lambda(\RR/J)+\nu(J)$. Meanwhile,   given
 that $\hdeg (J)= \deg (\RR) + \lambda(\RR/J)$,
the bound from the formula is
\[ \nu(I:J) \leq (2 + \nu(J) + \nu(I,J))(2+\lambda(\RR/J)+
\lambda(\RR/(I,J))-7,
\]
which is poorer.
}\end{Remark}

The advantage  comes in case $\dim \RR/I=1$, or $2$.
 Let us assume that $I$ and $J$ are Cohen-Macaulay ideals of
dimension $2$. To make the computation of $\hdeg$ simpler, if we take
the ideal module: $E=J \oplus I$.

\begin{Theorem}
Let $(\RR, \mathfrak{m})$ be a Gorenstein local
ring of dimension $d$ and let $J$ and $I$ be Cohen-Macaulay  ideals
of codimension
 $d-2$. Then
\[ \nu( I:J) + \nu(J:I)  \leq    (2+ \nu(J)+ \nu(I)) (\deg  (\RR/J) + \deg
(\RR/I)) + 2\cdot (\nu(J)+ \nu(I)) - 2.\]
\end{Theorem}

%\newpage

\subsubsection*{Vector bundles}

Let $(\RR, \mathfrak{m})$ be a Noetherian local ring and set $\XX=\Spec
\RR$. A case of HomAB question is that of a module $E$ that is locally
free on a subset $\YY\subset \XX$. When $\YY$ is the punctured
spectrum  $ \Spec \RR\setminus \{\mathfrak{m}\}$ of $\RR$, it was treated
in \cite{Dalili}, and a solution given for all modules of finite
projective dimension that are locally free in codimension $\leq
\dim \RR-2$. To remove the finite projective dimension requirement
requires a good understanding of the MCM $\RR$-modules.

Here we will consider this issue but also the case of the complement
of the hyperplane $\YY= \Spec \RR \setminus V(z)$, where $\RR/(z)$ is a
regular local ring.

The emphasis on
Auslander duals puts a great burden on the algebra
  \[\Tor_1^{\RR}(D(E),E)= \Hom_{R}(\Ext_{\RR}^1(E,\cdot),
  \Ext_{\RR}^1(E,\cdot)),\]
 to
determine effective bounds for its number of generators. The
 cases considered in the previous sections all  had a great
control over the support of  $\Tor_1^{\RR}(D(E),E)$.
It might be
worthwhile to consider the following general question.

\begin{Conjecture} \label{smalltor} {\rm  Let $E$ be a torsion free
module. If $\dim \Tor_1^{\RR}(D(E),E)\leq 1$, \[\nu(\Hom_{\RR}(E,E)) \leq
\rank(E)\cdot (\hdeg E + \nu(E)-\deg E).\]
}\end{Conjecture}

 If $\dim \Tor_1^{\RR}(D(E),E)\leq 0$, $E$ is a vector bundle
on the punctured spectrum of $\RR$. In \cite{Dalili}, this is dealt even
more generally if
   $E$ is free in codimension $d-2$, provided
$\mbox{\rm proj dim }E < \infty$, while in Theorem~\ref{freed2}, the
case of isolated singularity is dealt with.

%\newpage

\subsection{Higher cohomology modules} Given a Gorenstein local ring
$\RR$, and two finitely generated $\RR$-modules $A$ and $B$, we look at
the problem of bounding the number of generators of the modules
$\Ext_{\RR}^i(A,B)$, for $i>0$. The approach we use is
 straightforward: Consider a free presentation,
\[ 0 \rar L \lar F \lar A \rar 0,\]
pass to $L$ the given $\hdeg$ information on $A$, and use d\'{e}calage
to compare $\nu(\Ext_{\RR}^i(A,B))$ to
 $\nu(\Ext_{\RR}^{i-1}(L,B))$.
This is allowed since by the  cohomology exact sequences, we have the
short exact sequences
\[ 0 \rar \Hom_{\RR}(A,B) \lar \Hom_{\RR}(F,B) \lar \Hom_{\RR}(L,B) \lar
\Ext_{\RR}^1(A,B) \rar 0,
\] and
\[ \Ext_{\RR}^i(A,B) \simeq \Ext_{\RR}^{i-1}(L,B), \quad i>1. \]

We relate the degrees of $A$ to those of $L$. We shall assume that
the rank of $F$ is $\nu(A)$,
$\dim A=\dim \RR=d$, so that if $L\neq 0$, $\dim L=d$. This gives
\[ \deg(L) = \deg(F)-\deg(A).\]

We have the two expressions for $\hdeg(A)$ and $\hdeg(L)$
(\ref{hdegdef}):

\begin{eqnarray*}
 \hdeg(A) & = &\deg (A) +
 \sum_{i=1}^{d} {{d-1}\choose{i-1}}\cdot
 \hdeg(\mbox{\rm Ext}^i_{\RR}(A,R)),\\
 \hdeg(L) & = &\deg (L) +
 \sum_{i=1}^{d-1} {{d-1}\choose{i-1}}\cdot
 \hdeg(\mbox{\rm Ext}^i_{\RR}(L,\RR)),
\end{eqnarray*} since $\depth L> 0$. If we set $a_i=\hdeg
(\Ext_{\RR}^i(A,\RR))$, these formulas can be rewritten as
\begin{eqnarray*}
 \hdeg(A) & = &\deg (A) + a_1 +
 \sum_{i=2}^{d} {{d-1}\choose{i-1}}\cdot
 a_i,\\
 \hdeg(L) & = &\deg (L) +
 \sum_{i=1}^{d-1} {{d-1}\choose{i-1}}\cdot
 a_{i+1} = \deg(L) + \sum_{i=1}^{d-1} \frac{i}{d-i}{{d-1}\choose{i}}\cdot
 a_{i+1}.
\end{eqnarray*}

This gives
\[
 \frac{1}{d-1} \sum_{i=2}^{d}
{{d-1}\choose{i-1}}\cdot a_i\leq \hdeg(L)-\deg(L)
\leq  (d-1)\cdot \sum_{i=2}^{d}
{{d-1}\choose{i-1}}\cdot a_i.
\]

We now collect these estimations:

\begin{Proposition} \label{hdegsyz} Let $\RR$, $A$ and $L$ be as above.
Then \[
\deg(L) \leq  \nu(A)\deg(\RR)-\deg(A),\]
and if $c=\hdeg(A)-\deg(A)-\hdeg(\Ext^1(A,\RR))$
\[ \frac{1}{d-1} c
\leq \hdeg(L)-\deg(L) \leq
{(d-1)} c.
\]
\end{Proposition}

\begin{corollary}
Let $\RR$ be a Gorenstein local ring. If there are polynomials
$\ff(x,y)$ such that for any two finitely generated $\RR$-modules
$A,B$, $\nu(\Hom_{\RR}(A,B))\leq \ff(\hdeg(A), \hdeg(B))$, there are also
polynomials $\ff_i(x,y)$, $i\geq 1$, of the same degree such that for $i\geq 1$,
\[ \nu(\Ext_{\RR}^i(A,B)) \leq \ff_i(\hdeg(A), \hdeg(B)).\]

\end{corollary}

\subsubsection*{Biduals}
Let $\RR$ be a Noetherian local ring with a canonical module
$\omega_{\RR}$. For a torsion free $\RR$-module by its {\em bidual} \index{bidual
of a module} we shall mean the module \[ E^{**}=
\Hom_{\RR}(\Hom_{\RR}(E,\omega_{\RR}),\omega_{\RR}).\] There is a great deal of
flexibility in the construction of this module. For instance, suppose
$A\hookrightarrow R$ is a Noether normalization, that is, $\RR$ is
finite over the regular local subring $A$. Then
$E^{**}\simeq \Hom_A(\Hom_A(E,A),A)$. This means that it is natural
to assume that $E$ is a module of finite projective dimension over
the base ring.

\medskip

The construction of $E^{**}$ is
of utmost importance in view of its use in the processes of
$S_2$-ification would be estimations for $\nu(E^{**})$ dependent upon
data on $E$. One of the main results of \cite{Dalili} is the
assertion that if $\RR$ is a Gorenstein local ring of dimension $d$,
\[ \bdeg(E^*) \leq  (\deg (E) + d(d-1)/2)\hdeg (E).\]
Regrettably, it cannot be iterated to pass $\hdeg$ information to
$E^*$. Nevertheless some cases are easy to handle.

\medskip

Let $E$ be a module of projective dimension $1$, or more generally of
depth $d-1$, that is suppose
\[ 0 \rar F_1 \lar F_0 \lar E \rar 0,\]
is exact and $F_0,F_1$ are MCM modules.
Dualizing, we get exact sequences
\[ 0 \rar E^* \lar F_0^* \lar K \lar 0,\]
and
\[ 0 \rar K \lar F_1^* \lar \Ext_{\RR}^1(E,\RR) \rar 0.\]
The second sequence allow us to obtain a bound for $\hdeg(K)$ from
$\Ext_{\RR}^1(E,\RR)$, which taken into the other sequence leads to a bound
for $\hdeg(E^*)$. Now we apply Theorem~\ref{Kia5.3}, to bound
$\bdeg(E^{**})$.

\medskip

Another case is of torsion free modules $E$ of finite projective
dimension over   local ring $\RR$ of dimension
$3$. Consider an embedding
\[ 0 \rar F \lar E \lar C \rar 0,\]
where $F$ is a free module of rank $\rank(E)$. The module $C$ has
dimension $2$ and has projective dimension $1$ at every localization
of codimension one. According to Theorem~\ref{Fittdeg}, $\deg(C) =
\deg(\RR/\Fitt(C))$. In this case the divisorial class of $\Fitt(C)$ is
a principal ideal $(f)$.

By taking an embedding  $F\hookrightarrow E$ that
gives a minimum to $\deg(C)$, one may view this number as an
invariant of $E$.

Dualizing we get the exact sequences
\[ 0 \rar E^* \lar F^* \lar K \rar 0,\]
and
\[ 0 \rar K \lar \Ext^1_{\RR}(C,\RR) \lar \Ext_{\RR}^1(E,\RR) \rar 0.\] In the
second, $\Ext_{\RR}^1(E,\RR)$ has dimension at most one since $\dim R=3$
and $E$ is free in codimension one; thus $\deg(K)= \deg(D)$.
In turn, $D=\Ext_{\RR}^1(C,\RR)$ is a
Cohen-Macaulay of dimension two, as one easily shows that
\begin{eqnarray*} D = \Hom_{\RR/(f)}(C, \RR/(f)).\end{eqnarray*}
In particular, $\deg(D) = \deg(C)$.

Since $\Ext_{\RR}^2(K,\RR) = \Ext_{\RR}^3(D,\RR)$ and
$\Ext_{\RR}^3(K,\RR)=0$, we can
bound $\hdeg(K)$ in terms of $\deg(C)$ and $\hdeg(E)$. This allows
bounding $\hdeg(E^*)$ in terms of $\deg(C)$ and $\hdeg(E)$ as well,
and therefore will be in position to apply Theorem~\ref{Kia5.3}.

\section{Length Complexity of Tensor Products}

\subsection{Introduction}
 An analogue of the HomAB
question for tensor products is the following. Let $(\RR, \mathfrak{m})$ be a
Noetherian local ring, and let $A$ and $B$ be finitely generated
$\RR$-modules.

\begin{question} {\rm Can either
$\lambda(\H^0_{\mathfrak{m}}(A\otimes_{\RR}B))$ or
$\nu(\H^0_{\mathfrak{m}}(A\otimes_{\RR}B))$ be estimated in terms of
multiplicity invariants of $A$ and $B$? For instance, we look for the
existence of
polynomials $\ff(x,y)$ with rational coefficients such that
\[ \lambda(\H^0_{\mathfrak{m}}(A\otimes_{\RR}B))\leq \ff(\Deg(A), \Deg(B)).\]
More generally, we look for  bounds for  \[
\lambda(\H^0_{\mathfrak{m}}(\Tor_i^{\RR}(A,B))), i\geq 1.\]
}\end{question}

We label this question the {\em TorAB} Problem.\index{TorAB problem}
We shall consider special cases of these questions but
in low dimensions. Let us begin with
the following observation. We denote $\H^0_{\mathfrak{m}}(E)=E_0$, and
$\lambda(\H^0_{\mathfrak{m}}(A\otimes_{\RR}B))$ by $h_0(A \otimes_{\RR} B)$.

\begin{proposition}\label{AmodA0} Let $(\RR, \mathfrak{m})$ be a Noetherian local
ring and $A, B$ finitely generated $\RR$-modules. Then
\begin{eqnarray*}
\nu(\H_{\mathfrak{m}}^0(A\otimes_{\RR} B)) &\leq & \nu(A_0)\nu(B) +2
\nu(A)\nu(B_0) + \nu(\H_{\mathfrak{m}}^0(A/A_0\otimes_{\RR} B/B_0) \\
 h_0(A\otimes_{\RR} B)&\leq &
h_0(A)\nu(B) + h_0(B)\nu(A) + h_0(A/A_0\otimes_{\RR}
B/B_0).
\end{eqnarray*}
\end{proposition}

\subsubsection{Dimension $1$} Suppose $\RR$ is a local domain of
dimension $1$. We start our discussion with the case of two ideals,
$I,J\subset R$. Consider the commutative diagram
\[
\diagram
0 \rto  & L  \rto\dto & F \rto \dto & IJ
\rto\dto  &0 \\
0 \rto          &  T \rto                 & I\otimes_{\RR}J \rto  & IJ \rto  &0
\enddiagram
\]
where $F$ is a free presentation of $I\otimes_{\RR}J$, and therefore it
has rank $\nu(I)\nu(J)$. $L$ is a torsion free  module with
$\rank(L)=\rank(F)-1$,  and therefore it can be generated by $\deg(L) =
(\rank(F)-1)\deg(\RR)$ elements.

\begin{proposition} If $\RR$ is a   local domain of dimension one,
essentially of finite type over a field   and $I,J$ are
$\RR$-ideals, then
\begin{eqnarray*}
\nu(\H_{\mathfrak{m}}^0(I\otimes_{\RR} J)) &\leq & (\nu(I)\nu(J)-1)\deg (\RR).\\
h_0(I\otimes_{\RR}J)  &\leq & (\nu(I)\nu(J)-1)\deg (\RR) \lambda(\RR/(I,J,K)),
\end{eqnarray*}
where $K$ is the Jacobian ideal of $\RR$.
\end{proposition}

\begin{proof} The second formula arises because $L$ maps onto $T$, the torsion submodule of
$I\otimes_{\RR}J$, which is also annihilated by $I, J$
and $K$. The last assertion
is a consequence of the fact that $T= \Tor_1^{\RR}(I,\RR/J)$, and by
Theorem~\ref{Wang1} and Proposition~\ref{Wang2}, $K$ will annihilate
it. 
\end{proof}

A small enhancement occurs since one can replace $I$ and $J$ by
isomorphic ideals.
In other words, the last factor, $\lambda(\RR/(I,J,K))$, can be replaced
by  $\lambda(\RR/(\tau(I),\tau(J),K))$, where $\tau(I) $ and $ \tau(J)$
are their {\em trace} ideals ($\tau(I) = \mbox{\rm image
$I\otimes_{\RR}\Hom_{\RR}(I,\RR)\rar \RR$}$).

%\bigskip

The version for modules is similar:

\begin{proposition} If $\RR$ is a   local domain of dimension one,
essentially of finite type over a field   and let $A,B$ be finitely
generated  torsion free
$\RR$-modules. Then
\begin{eqnarray*}
\nu(\H_{\mathfrak{m}}^0(A\otimes_{\RR} B)) &\leq &
(\nu(A)\nu(B)-\rank(A)\rank(B))\deg (\RR).\\
h_0(A\otimes_{\RR}B)  &\leq & (\nu(A)\nu(B)-\rank(A)\rank(B))\deg (\RR) \lambda(\RR/K),
\end{eqnarray*}
where $K$ is the Jacobian ideal of $\RR$.
\end{proposition}

To extend this estimation of $h_0(A\otimes_{\RR}B)$ to finitely generated torsion free
$\RR$-modules that takes into account annihilators
we must equip the modules--as is the case of ideals--with a
privileged embedding into free modules.

\begin{lemma} Let $\RR$ be a Noetherian integral domain of dimension
$1$, with finite integral closure. Let $A$ be a torsion free
$\RR$-module of rank $r$ with the embedding $A\rar F=\RR^r$. Let $I$ be
the ideal $\mbox{\rm image }(\wedge^rA \rar \wedge F = \RR)$. Then $I$
annihilates $F/A$.
\end{lemma}

\subsubsection{Vector bundles} Let $(\RR, \mathfrak{m})$ be a regular local
ring of dimension $d\geq 2$. For a vector bundle $A$ (that is, a
finitely generated $\RR$-module that is free on the punctured
spectrum), we consider estimates of the form
\begin{eqnarray} \label{vbtorsion}
h_0(A \otimes_{\RR} B) &\leq & c(\RR) \cdot \hdeg(A) \cdot  \hdeg(B),
\end{eqnarray}
where $c(\RR)$ is a constant depending on $\RR$ and $B$ is a finitely
generated $\RR$-module.

\medskip

From the general observations above, we may assume that $\depth A$
and $\depth B$ are positive.

\medskip

We make some reductions beginning with the following. Since $A$ is
torsion free, consider the natural exact sequence
\[ 0 \rar A \lar A^{**} \lar C \rar 0.
\] Note that
$C$ is a module of finite support, and that $A^{**}$ is a vector
bundle. Furthermore, since $A^{**}$ has depth at least $2$, a direct
calculation yields
\[ \hdeg(A) = \hdeg(A^{**}) + \hdeg(C).\]

Tensoring the exact sequence by $B$, gives the exact complex
\[ \Tor_1^{\RR}(C,B) \lar A\otimes_{\RR}B \lar
A^{**} \otimes_{\RR}B,\]
from which we obtain
\[ h_0(A\otimes_{\RR}B) \leq  h_0(A^{**}\otimes_{\RR}B) +
\lambda(\Tor_1^{\RR}(C,B)).
\]
As $C$ is a module of length $\hdeg(C)$,
\[ \lambda(\Tor_1^{\RR}(C,B))\leq \beta_1(k)\cdot \hdeg(C)\cdot \nu(B).
\] We recall that $\nu(B)\leq \hdeg(B)$.

\bigskip

Let $(\RR, \mathfrak{m})$ be a Gorenstein local ring of dimension $d>0$
and let $A$ be a finitely generated $\RR$-module that is free on the
punctured spectrum, and has finite projective dimension.
 We seek to estimate $h_0(A\otimes B)$ for various $\RR$-modules $B$.

\medskip

Let $B$ be an $\RR$-module with $\depth B>0$
with the minimal free resolution
\[ 0\rar F_{d-1} \lar F_{d-2} \lar \cdots \lar F_1 \lar F_0 \lar B \rar 0, \]
%where up to $F_{d-1}$ is a minimal free presentation of $A$.
Tensoring by $A$ gives a complex
\[ 0\rar F_{d-1}\otimes A \lar F_{d-2}\otimes A  \lar \cdots \lar
F_1\otimes  A\lar F_0\otimes A \lar B\otimes A \rar 0, \]
whose homology $H_i=\Tor_i(B,A)$, $d>i>0$, has finite support. Denoting by $B_i$
and $Z_i$ its modules of boundaries and cycles,
we have several exact sequences
\[ 0 \rar B_0 \lar F_0\otimes A \lar B\otimes A \rar 0,\]
\[ 0 \rar B_i \lar Z_i \lar H_i \rar 0,\]
\[ 0 \rar Z_{i} \lar F_{i}\otimes A \lar B_{i-1} \rar 0.\]
Taking local cohomology, we obtain the following acyclic complexes of
modules of finite length
\[
 \H^0_{\mathfrak{m}}(F_0\otimes A) \lar \H^0_{\mathfrak{m}}(B\otimes A) \lar
\H^1_{\mathfrak{m}}(B_0), \]
%\[   \H^1_{m} (B_1) \lar \H^1_{m}(Z_1)\rar 0,  \]
%\[   \H^i_{m} (B_i) \lar \H^i_{m}(Z_i), \quad i>1 \]

\[  \H^i_{\mathfrak{m}}(F_{i}\otimes A) \lar \H^i_{\mathfrak{m}}(B_{i-1})
\lar \H^{i+1}_{\mathfrak{m}}(Z_{i}), \quad i\geq 1. \]

Now we collect the inequalities of length, starting with
\[ h_0(B\otimes A)\leq \beta_0(B)h_0(A) + h_1(B_0)
\]
and
\[
 h_i(B_{i-1})\leq \beta_i(B)h_i(A) + h_{i+1}(Z_{i})
\leq \beta_i(B)h_i(A) + h_{i+1}(B_{i}),
\]
where we replace the rank of the modules $F_i$ by $\beta_i(B)$.

%\H^_{m}(

\begin{Theorem} Let $B$ be a module of projective dimension $< d$ and let
$A$ be a module that is
free on the punctured spectrum.
Then
\[ h_0(A\otimes B) \leq \sum_{i=0}^{d-1} \beta_i(B) h_i(A).\]
\end{Theorem}

Let us rewrite this inequality in case $\RR$ is a regular local ring.
Since by local duality $h_i(A)= \lambda(\Ext_{\RR}^{d-i}(A,\RR))$ and
$\beta_i(B)\leq \beta_i(k)\hdeg(B)$, we have
\begin{eqnarray*} h_0(A\otimes B) &\leq  & \hdeg(B )\cdot \sum_{i=0}^{d-1}{d\choose
i}\lambda(\Ext_{\RR}^{d-i}(A,\RR))\\
&\leq & d\cdot \hdeg(B)\sum_{i=1}^{d}
 {{d-1}\choose{i-1}}
\lambda(\Ext_{\RR}^{d-i}(A,\RR))\\
& \leq & d\cdot \hdeg(A)\cdot \hdeg(B).
\end{eqnarray*}

\begin{Theorem} \label{h0vb} Let $\RR$ be a regular local ring of
dimension $d$. If $A$ is a f.g. module free on the punctured
spectrum, then for any f.g. $\RR$-module $B$
\begin{eqnarray*}
h_0(A\otimes B) \leq d\cdot \hdeg(A)\cdot \hdeg(B).
\end{eqnarray*}
\end{Theorem}

\begin{proof} It suffices
to add to $\hdeg(B)$ the correction $h_0(B)$ as given in
Proposition~\ref{AmodA0}. 
\end{proof}

\subsubsection{Dimensions $2$ and $3$}
Let $(\RR, \mathfrak{m})$ be a regular local ring of dimension $2$ (or
a polynomial ring $k[x,y]$ over the field $k$). For two $\RR$-modules
$A$ and $B$ we are going to study
\[ h_0(A\otimes_{\RR}B) = \lambda(\H^0_{\mathfrak{m}}(A\otimes_{\RR}B))\]
through a series of reductions on $A$, $B$ and $\RR$. Eventually the
problem will settle on the consideration of a special class of
one-dimensional rings.

\medskip

We already familiar with the stripping away from $A$ and $B$ of their
submodules of finite support, so we may assume that these modules
have depth $\geq 1$. Let $A$ be a module of dimension $2$,
and denote by $A_0$ the torsion submodule of
$A$. Consider the natural exact sequence
\[0 \rar A_0 \lar A \lar A'\rar 0,\]
with $A'$ torsion free. If $A_0\neq 0$, it is a Cohen-Macaulay module
of dimension $1$. We have the exact sequence
\[ 0 \rar \Ext_{\RR}^1(A',\RR) \lar \Ext_{\RR}^1(A,\RR) \lar
\Ext_{\RR}^1(A_\RR)
\rar 0,
\] that yields

 \begin{eqnarray*}
\deg A &=& \deg A' \\
\hdeg(\Ext_{\RR}^1(A,\RR)) &=&
\hdeg(\Ext_{\RR}^1(A',\RR)) +\deg A_0,
\end{eqnarray*}
in particular
\begin{eqnarray*}
\hdeg(A) &=& \deg A_0 + \hdeg(A').
\end{eqnarray*}

As a consequence, $\hdeg(A')$ and $\hdeg(A_0)$ are bounded in
terms of $\hdeg(A)$. Now we tensor the sequence by $B$ to get the
complex
\[ \Tor_1^{\RR}(A',B) \lar A_0\otimes_{\RR} B \lar
 A\otimes_{\RR} B \lar  A\otimes_{\RR} B \rar 0. \]
If we denote by $L$ the image of $A_0\otimes_{\RR} B$ in $A\otimes_{\RR} B$,
since $\Tor_1^{\RR}(A',B)$ is a module of finite length, we
have
\begin{eqnarray*} h_0(A\otimes_{\RR}B) &\leq & h_0(A'\otimes_{\RR}B) + h_0(L),\\
h_0(L) & \leq & h_0(A_0\otimes_{\RR}B).
\end{eqnarray*}

If we apply a similar reduction to $B$ and combine, we get
\begin{eqnarray*} h_0(A\otimes_{\RR}B) &\leq h_0(A'\otimes_{\RR}B') +
h_0(A_0\otimes_{\RR}B') +
h_0(A'\otimes_{\RR}B_0)
+  h_0(A_0\otimes_{\RR}B_0).
\end{eqnarray*}

 \begin{eqnarray*}
\deg A &=& \deg A' \\
\hdeg(\Ext_{\RR}^1(A,\RR)) &=&
\hdeg(\Ext_{\RR}^1(A',\RR)) +\deg A_0,
\end{eqnarray*}
in particular
\begin{eqnarray*}
\hdeg(A) &=& \deg A_0 + \hdeg(A').
\end{eqnarray*}

A term like $h_0(A'\otimes B_0)$ is easy to estimate. Let
\[ 0 \rar F_1 \lar F_0 \lar A' \rar 0\] be a minimal resolution of
$A'$. Tensoring with $B_0$, we get the exact sequence
\[ 0 \rar \Tor_1^{\RR}(A',B_0) \lar F_1\otimes_{\RR}B_0
\lar F_0\otimes_{\RR}B_0
\lar A'\otimes_{\RR}B_0 \rar 0.
\] Since
$\Tor_1^{\RR}(A',B_0)$
has finite support and $F_1\otimes_{\RR}B_0$ has positive depth,
$\Tor_1^{\RR}(A',B_0)=0$. To compute $h_0(A'\otimes_{\RR}B_0)$, consider the
a minimal resolution of $B_0$,
\[ 0 \rar G \lar G \lar B_0 \rar 0,\]
and
 the exact sequence
\[ 0 \rar A'\otimes_{\RR} G \lar A' \otimes_{\RR} G \lar A'\otimes_{\RR}B_0 \rar 0.
\] From the cohomology exact sequence, we have the surjection
\[ \Ext_{\RR}^1(A'\otimes G,R) \lar \Ext_{\RR}^2(A'\otimes_{\RR} B_0),\]
and therefore, since $\Ext_{\RR}^1(A,\RR)$ is a module of finite support,
\[ h_0(A'\otimes_{\RR}B_0) \leq \nu(G)\lambda(\Ext_{\RR}^1(A',\RR)).\]
This shows that
\[h_0(A'\otimes_{\RR}B_0) \leq \nu(G)\cdot (\hdeg(A')-\deg A') <
\hdeg(A')\cdot \hdeg(B_0)
.\]

The reductions thus far lead us to assume that $A$ and $B$ are
$\RR$-modules of positive depth and dimension $1$. Let
\[ 0 \rar F \stackrel{\varphi}{\lar} F \lar A \rar 0\]
be a minimal free resolution of $A$. By Theorem~\ref{Fittdeg},
\[ \deg A = \deg(\RR/\det(\varphi)).\]
Since $\det (\varphi)$ annihilates $A$, we could view $A\otimes_{\RR}B $
as a module of over $\RR/(\det(\varphi \circ \psi))$ where $\psi$ is
the corresponding matrix in the presentation of $B$.

\medskip

To avoid dealing with two matrices, replacing $A$ by $A\oplus B$, we
may consider $h_0(A\otimes_{\RR}A)$, but still denote by $\varphi$ the
presentation matrix (instead of $\varphi\oplus \psi$), and set
$\SS=\RR/(\det(\varphi))$; note that $\deg \SS = \deg A$.

\begin{Example}{\rm
We consider a cautionary family of examples to show that other
numerical readings must be incorporated into the estimates for
$h_0(A\otimes_{\RR}A)$.

Let $A$ be a module generated by two elements, with a free resolution
\[ 0 \rar F \stackrel{\varphi}{\lar} F \lar A \rar 0.\]
Suppose $k$ is a field of characteristic $\neq 2$. To calculate
$h_0(A\otimes_{\RR}A)$, we make use of the decomposition
\[ A \otimes_{\RR}A = S_2(A) \oplus \wedge^2A.\]
Given a matrix representation,
\[ \varphi = \left[ \begin{array}{ll}
a_{11} & a_{12} \\
a_{21} & a_{22} \\
\end{array} \right],
\]
one has
\[\wedge^2 A \simeq R/I_1(\varphi)=R/(a_{11}, a_{12}, a_{21},
a_{22}).\]

The symmetric square of $A$, $S_2(A)$, has a resolution
\[ 0 \rar R \stackrel{\phi}{\lar} F\otimes_{\RR}F \stackrel{\psi}{\lar}
S_2(F),\]
where
\begin{eqnarray*}
\psi(u\otimes v) &=& u \cdot  \varphi(v) + v\cdot  \varphi(u) \\
\phi (u\wedge v) & = & \varphi'(u)\otimes v-\varphi'(v)\otimes u \\
\end{eqnarray*}
where $\varphi'$ is the matrix obtained from $\varphi$ by dividing
out its entries by their gcd $a$, $\varphi=a\cdot \varphi'$.

A straightforward calculation will give
\[ \Ext_{\RR}^2(S_2(A), \RR) = \RR/I_1(\varphi').\]
This shows that
\[ h_0(A\otimes_{\RR}A) = h_0(\RR/I_1(\varphi))+h_0(\RR/I_1(\varphi'))=
h_0(a\RR/aI_1(\varphi'))+h_0(\RR/I_1(\varphi'))=
2\cdot \lambda(\RR/I_1(\varphi')).
\]

Thus  the matrix
\[ \varphi = \left[ \begin{array}{ll}
x & y^n \\
0 & x \\
\end{array} \right],
\]
will define a module $A$, with $\deg(A)=2$, but
$h_0(A\otimes_{\RR}A)=2n$. This means that we must take into account the
degrees of the
entries of $\varphi$ itself.

}\end{Example}

\medskip

For the final reduction, write $\SS=\RR/(\xx)$, where  $\xx$ has the
primary decomposition
$\xx=x_1\cdots x_n$. Setting $z_i= \xx/x_i$, consider the exact
sequence
\[ 0 \rar \RR/(\xx) \lar \RR/(z_1) \oplus \cdots \oplus \RR/(z_n) \lar C
\rar 0,\]
induced by the mapping $1 \mapsto (z_1, \ldots, z_n)$. $C$ is a
module of finite length, and making use of duality and the inequality
(\ref{degxx}),
\[ \lambda(C) \leq \frac{1}{2} ((\deg(\RR/(\xx))^2-n).\]

Tensoring this sequence by $A$, gives
\[ \Tor_1^{\RR}(A,C) \lar A \lar A_1 \oplus \cdots \oplus A_n \lar
A\otimes_{\RR}C \rar 0,\]
and since $\depth A>0$, we have the exact sequence
\[ 0 \rar A\lar A_1 \oplus \cdots \oplus A_n \lar A\otimes_{\RR}C  \rar 0,\]
where $A_i = A/x_iA$ and $\lambda(A\otimes C)\leq \nu(A) \lambda(C)$.
These relations give that
\begin{eqnarray*}
\deg A &=& \sum_{i=1}^n \deg A_i \\
\lambda(A\otimes_{\RR} C) &\geq & \sum_{i=1}^n h_0( A_i).
\end{eqnarray*}
These inequalities show that we are still tracking the $\hdeg(A_i)$ is
terms of $\deg A$.

Tensoring the last  exact sequence by $A$, we obtain the exact
complex
\[ \Tor_1^{\RR}(A, A\otimes_{\RR}C) \lar A\otimes_{\RR}A \lar A_1\otimes_{\RR} A_1
\oplus \cdots \oplus A_n\otimes_{\RR}A_n, \]
from which we have
\begin{eqnarray*} h_0(A\otimes_{\RR}A) &\leq &
\sum_{i=1}^n h_0(A_i\otimes_{\RR}A_i) + \lambda
(\Tor_1^{\RR}(A, A\otimes_{\RR}C)) \\
&\leq & \sum_{i=1}^n h_0(A_i\otimes_{\RR}A_i) + \beta_1(A)\cdot \nu(A)
\cdot \lambda(C). \\
\end{eqnarray*}

Let us sum up these reductions as follows:

\begin{proposition} Let $\RR$ be a two-dimensional regular local ring
and let $A$ be a Cohen-Macaulay $\RR$-module of dimension one. Then
\[ h_0(A\otimes_{\RR}A) \leq 3 \cdot  \hdeg(A)^4,\]
provided
\[ h_0(A\otimes_{\RR}A) \leq 2 \cdot  \hdeg(A)^4\]
 if $\ann A$ is a primary ideal.
\end{proposition}

\begin{proof} Note that $\beta_1(A) \leq \beta_1(k) \cdot \hdeg(A)$,  $\nu(A)\leq
\hdeg(A)$, and $\lambda(C)< \frac{1}{2} \hdeg(A)^2$.
\end{proof}

\subsubsection*{Dimension $3$}
Theorem~\ref{h0vb} can be used to deal with torsionfree modules of
dimension three.

\begin{Theorem} \label{h0dim3}
Let $\RR$ be a regular local ring of dimension $3$,
and let   $A$ and $B$ be torsionfree $\RR$--modules.
Then
\begin{eqnarray}
h_0(A\otimes B)& < & 4\cdot \hdeg(A)\cdot \hdeg(B).
\end{eqnarray}
\end{Theorem}

\begin{proof}
 Consider the
natural exact sequence
\[ 0 \rar A \lar A^{**} \lar C \rar 0.\]
A straightforward calculation will show that
\begin{eqnarray} \label{h0eq0}
\hdeg(A) &=& \hdeg(A^{**}) + \hdeg(C).
\end{eqnarray}
Note that $A^{**}$ is a vector bundle of projective dimension at most
$1$ by the Auslander-Buchsbaum equality \cite[th]{BH}), and $C$ is a
module of dimension at most $1$. Tensoring by
the torsionfree $\RR$-module $B$, we have
the exact sequence
\[ \Tor_1^{\RR}(A^{**}, B) \lar \Tor_1(C,B) \lar A\otimes B \lar
A^{**}\otimes B \lar C\otimes B \rar 0,
\] where $\Tor_1^{\RR}(A^{**}, B)=0$, since $\mbox{\rm proj dim
}A^{**}\leq 1$ and $B$ is torsionfree.

\medskip

From the exact sequence, we have
\begin{eqnarray} \label{h0eq1}
h_0(A\otimes B) &\leq & h_0(A^{**}\otimes B) + h_0(\Tor_1^{\RR}(C,B)).
\end{eqnarray}
Because $A^{**}$ is a vector bundle, by Theorem~\ref{h0vb},
\begin{eqnarray} \label{h0eq2}
h_0(A^{**}\otimes B) &\leq & 3\cdot
\hdeg(A^{**})\cdot \hdeg(B).
\end{eqnarray}

For the module $\Tor_1^{\RR}(C,B)$, from a minimal free presentation of
$B$,
\[ 0 \rar L \lar F \lar B\rar 0,\]
we have an embedding $\Tor_1^{\RR}(C,B) \rar C\otimes L$, and therefore
\begin{eqnarray} \label{h0eq3}
h_0(\Tor_1^{\RR}(C,B)) &\leq & h_0(C\otimes L) \leq 3 \cdot \hdeg(L)\cdot
\hdeg(C),
\end{eqnarray}
because $L$ is a vector bundle. In turn
\begin{eqnarray*}
\deg(L) &=& \beta_1(B)-\deg(B)\leq
\beta_1(\RR/\mathfrak{m})\cdot \hdeg(B)-\deg(B) \\
&=& 3\cdot \hdeg(B)-\deg(B)
\end{eqnarray*}
by Theorem~\ref{Degandbetti},
and since $\Ext_{\RR}^1(L,R)= \Ext_{\RR}^2(B,R)$,
 \begin{eqnarray*}
\hdeg(L) = \deg(L)+ \hdeg(\Ext_{\RR}^1(L,R)) & < & 4\cdot \hdeg(B).
\end{eqnarray*}

Finally we collect (\ref{h0eq2}) and (\ref{h0eq3}) into
(\ref{h0eq1}), along with (\ref{h0eq0}),
\begin{eqnarray*}
h_0(A\otimes B) &< & 3\cdot \hdeg(A^{**})\cdot \hdeg(B) + 4\cdot
\hdeg(B)\cdot \hdeg(C) \\
 &< & 4\cdot \hdeg(A) \cdot \hdeg(B),
\end{eqnarray*}
as asserted.
\end{proof}

\subsubsection*{Graded modules.} To argue for the validity of the
conjecture, we give a rough (high degree) estimate for the case of
graded modules over $\RR=k[x,y]$. We may assume that $A$ is not a
cyclic module. Furthermore, we shall assume that $A$ is
equi-generated.

\medskip

We can apply to
the exact sequence
\[ 0\rar A_0\lar A \lar A'\rar 0\]  Proposition~\ref{addiofreg},
to obtain bounds for $\reg(A_0)$ is terms of $\reg(A)$ and therefore
bound $\reg(A_0)$ in terms of $\hdeg(A)+\alpha(A)$,
 according to
Theorem~\ref{Nagel}.

Let us apply it to the graded $k[x,y]$-module $A$ of depth $> 0$. In
the exact sequence
\[ 0 \rar A_0 \lar A \lar A'\rar 0 \]
we already remarked that $\hdeg(A)=\hdeg(A') + \deg A_0$. It is also
the case that if $A$ is generated by elements of degree $\leq
\alpha(A)$, then by the proposition above and Theorem~\ref{Nagel},
$\reg(A_0)< \hdeg(A) + \alpha(A)$. Actually, since $A_0$ is
Cohen-Macaulay,  a direct calculation will show that $\reg(A_0)\leq
\reg(A)$.

\medskip

We may assume that $A$ is a one-dimensional graded $\RR$-module with a
minimal resolution
\[ 0 \rar F \stackrel{\varphi}{\lar} F \lar A \rar 0.\]
A presentation of $A\otimes_{\RR}A$ is given by
\[ F\otimes_{\RR}F \oplus F\otimes_{\RR}F \stackrel{\psi}{\lar} F\otimes_{\RR}F,  \]
where
$\psi = \varphi\otimes I - I\otimes \varphi$. The kernel of $\psi$ contains the
image of
\[ \phi: F\otimes_{\RR}F \lar F\otimes_{\RR}F \oplus F\otimes_{\RR}F, \quad \phi=
I \otimes \varphi \oplus \varphi \otimes I.\]

Since $\RR=k[x,y]$, $L=\ker(\psi)$ is a free $\RR$-module of rank $r^2$,
$r=\rank(F)$. Because $\phi$ is injective, its image $L_0$ is a free
$\RR$-submodule of $\FF=F\otimes_{\RR} F \oplus F\otimes_{\RR} F$ of the same rank as
$L$, $L_0\subset L$, $\phi': F\otimes F\rar \FF$.
 It follows that the degrees of the entries of
$\phi'$ cannot be higher than those of $\phi$.

To estimate  $h_0(A\otimes_{\RR}A)$, note that
$\Ext_{\RR}^2(A\otimes_{\RR}A,\RR)$ is the cokernel of map of $\phi$. This is a
module generated by $r^2$ elements, annihilated by the maximal minors
of $\phi$.
We already have that $\ff=\det(\varphi)$ annihilates $A$. Now we look for
an element  $\hh$ in the ideal of maximal minors of $\phi'$
 so that $(\ff, \hh)$ has finite colength, and as a
consequence we would have
\[ h_0(A\otimes_{\RR}A) \leq r^2\cdot \lambda(\RR/(\ff,\hh))=
r^2\cdot \deg(\RR/(\ff))\deg(\RR/(\hh)).\]

The entries of $\varphi$ have degree $\leq \deg(A)-1$, so the minors
$\hh$ of $\phi$ have degree
\[\deg \hh \leq r^2\cdot (\deg(A)-1).
\]

\begin{proposition} If $A$ is a graded $\RR$-module of dimension $1$,
then
\[ h_0(A\otimes_{\RR}A) \leq r^4\cdot \deg(A)(\deg(A)-1)< \deg(A)^6.\]
\end{proposition}

\subsection{Graded modules}

Kia Dalili has made use of the work of Chardin, Ha and Hoa
(\cite{CHH9}) to give an affirmative answer to the HomAB question and
the related question on tensor products.

\begin{Theorem}[Dalili] Let $\RR=k[x_1, \ldots, x_n]$ be a standard graded algebra over
the field $k$. There exist polynomials $\ff_i$ and $\g2_i$ such that for
any two finitely generated graded modules $A$ and $B$
\begin{itemize}
\item[{\rm (a)}] $\bdeg(\Ext_{\RR}^i(A, B))\leq \ff_i(\bdeg(A), \bdeg(B), \alpha(A),
\alpha(B))$
\item[{\rm (b)}] $\bdeg(\Tor_i^{\RR}(A, B))\leq \g2_i(\bdeg(A), \bdeg(B), \alpha(A),
\alpha(B))$.
\end{itemize}
\end{Theorem}

\begin{proof}
We will just sketch his argument for $\bdeg(\Hom_{\RR}(A,B))$.
\begin{itemize}
\item[$\bullet$] According to Theorem~\ref{CHH2.5}, there is a polynomial
$\ff_0$ so that
\[ \reg(\Hom_{\RR}(A,B))\leq \ff_0(\reg(A), \reg(B), \mbox{\rm Betti
numbers of $A$ and $B$}).\]
On the other hand, by Theorem~\ref{Nagel},
\[ \reg(P) \leq \Deg(P) + \alpha(P).\]
for any extended degree $\Deg$.  Since
each Betti number of a graded module $P$ satisfies
\[ \beta_i(P)\leq \Deg(P)\binom{n}{i}\]
the inequality leads to another polynomial inequality
\[ \reg(\Hom_{\RR}(A,B)\leq \ff(\Deg(A), \Deg(B),\alpha(A),
\alpha(B)).\]

\item[$\bullet$] We now apply to properties of $\bdeg$. If $F_0\rar A$ is a
minimal
free presentation of $A$, $\Hom_{\RR}(A,B)\subset \Hom_{\RR}(F_0, B)
$,
by Proposition~\ref{Gu2},
\[ \bdeg (\Hom_{\RR}(A,B))\leq \sum_{j=0}^{\reg(\Hom_{\RR}(A,B))}
H_{\Hom_{\RR}(F_0,B)}(j).
\]

\end{itemize}
\end{proof}

\subsection{Open questions}

We will now mention some questions and open problems. We have already
mentioned unsolved problems in context. The questions here represent
boundary cases.

\begin{Problem} \label{challenge2}{\rm Let $\RR$ be a Cohen-Macaulay
local ring with a canonical module $\omega$, and let $E$ be the
maximal Cohen-Macaulay defined in (\ref{challenge0}).
Find a bound for $\nu(\Hom_{\RR}(E,E))$. The problem shows the difficulty
of dealing with non-Gorenstein rings.
}\end{Problem}

\begin{Problem}{\rm
A different challenge arises if $E$ is a self-dual module, that is if
$E\simeq E^*$. (Examples are rank two reflexive modules over
factorial domains.) For these modules $\hdeg (D(E))$ is completely
determined by $\hdeg (E)$.
}\end{Problem}

\begin{Problem}{\rm Let $E$ be a Cohen-Macaulay module of dimension
$3$.  It would be helpful
 to understand this
corner case.
}\end{Problem}

\begin{Problem}{\rm If $\RR$ is a regular local of
dimension $2$(!), what is the relationship
between  $\bdeg$ and $\hdeg$?
}\end{Problem}

\begin{Problem}{\rm
 Let $\RR$ be a Noetherian local domain and $A$ a finitely
generated torsionfree $\RR$-module. Is there an integer $e=e(\RR)$
guaranteeing that if $M$ is not $\RR$-free, then the tensor power
$M^{\otimes {e}}$ has nontrivial torsion? The motivation is
 a result of Auslander
(\cite{Aus61}, see also \cite{Licht66}) that asserts that $e=\dim R$ works
for all regular local rings.
For instance, if $\RR$ is a one-dimensional domain, will $e=2$ work?
A more realistic question is, if $\RR$ is a Cohen-Macaulay local domain
of dimension $d$ and multiplicity $\mu$, will
\[ e= d+\mu-1\]
suffice? Note that if we make no attempt to determine uniform bounds
for $e$, if $\bbq\subset \RR$, then for a module $M$ of rank $r$ and
minimal number of generators $n$, then the embedding
\[ 0\neq \wedge^n M \hookrightarrow M^{\otimes n}\]
shows the existence of a test power.
}
\end{Problem}

\section{Local Modules}\index{local module}
\subsection{Introduction}

Let $(\RR, \m)$ be a Noetherian local ring and $E$ a finitely
generated $\RR$-module. The ring $\CC=\Hom_{\RR}(E,E)$ expresses many
properties of $E$ some of which are well hidden.
Let us list some of these properties:

\begin{itemize}
\item[$\bullet$] The HomAB question: Besides the discussion in
Section~\ref{homabsection}, there are
some noteworthy endomorphism rings
are beginning to appear in the literature. One class of the most
intriguing are those of finite global dimension; see  \cite{DaoHu10},
\cite{VdBergh}.
\index{noncommutative
regular ring}

\medskip

\item[$\bullet$] Degree representation of a module question: For a module $E$, its degree
representation is the smallest integer $r$--when it exists--such that
there is an embedding of $\RR$-algebras $\varphi: \Hom_{\RR}(E,E)\lar
 M_r(\RR)$.  Which modules have such representation?

\medskip

\item[$\bullet$] Cancellation questions: Suppose $\RR$ is a normal local domain.
Given reflexive $\RR$--modules $E_1$ and $E_2$ with $\RR$-isomorphic
endomorphism rings $\End(E_1)\simeq \End(E_2)$, what is the
relationship between them?

\medskip

\item[$\bullet$] Construction of local modules:
We will focus on an extreme case.
The terminology {\em Local Module } stands for a finitely generated
module $E$ over a local Noetherian ring $(\RR, \m)$ such that
the endomorphism ring $\End_{\RR}(E)=\Hom_{\RR}(E,E)$ has a unique two-sided
maximal ideal. These modules are useful in building other modules
whose modules of endomorphisms allow for the determination of their
Jacobson radicals.
 We want to focus on
some of these properties related to its Jacobson radical.

\end{itemize}

\medskip

\begin{example}{\rm Let us indicate some elementary examples.
\begin{enumerate}
\item[{\rm (a)}]  If $E$ is a free $\RR$-module, $E$ is local.

\item[{\rm (b)}]
If $\RR=k[x], x^2=0$, $k$ a field (that is, the ring of dual
numbers), then any f.g. module has the form  $E=\RR^m \oplus k^n$; it
follows that $E$ is local.

\item[{\rm (c)}]
 If $(\RR, (x))$ is a discrete valuation ring and $E=F\oplus T$,
$F$ free and $T$ torsion,
\[ \Hom(E,E)= \left[ \begin{array}{ll}
\Hom(F,F) & \Hom(F,T) \\
\Hom(T,F)=0 & \Hom(T,T)
\end{array} \right]
\]
Given that $T=\bigoplus C_i$, $C_i$ cyclic, and follows easily that
\[ \Hom(E,E)/\mbox{\rm radical}
= \left[ \begin{array}{ll}
\AA & 0 \\
0 & \BB
\end{array} \right]
\]
where  $\AA$ and $\BB$ are matrix rings over the residue field of
$\RR$. Thus $E$ is local if it is a free or a torsion module.

\end{enumerate}
}\end{example}

\begin{question}{\rm  Let $(\RR,\m)$ be a local Noetherian ring of residue
field $k$.
\begin{enumerate}
\item[{\rm (a)}] How to  tell local modules?
If $E$ is a module of finite length, then the length of $\CC$ can be
bound in two ways:
\[ \lambda(\CC)\leq \max\{ \nu(E)\cdot \lambda(E), \ t(E)\cdot
\lambda(E)
\},
\]
where $t(E)$ is the Cohen-Macaulay type of $E$. Which invariants of
$E$ refine this estimate?

\item[{\rm (b)}] How to build local modules? We will consider a contruction
using Auslander duality.

\item[{\rm (c)}] What are the Artinian local modules like?
\end{enumerate}
}
\end{question}

\subsection{Jacobson radical}

We now treat conditions for the algebra $\CC=\Hom_{\RR}(E,E)$
to have a unique two-sided maximal ideal. Observe that
by Nakayama Lemma, $\m \CC$ is contained in the Jacobson radical
$\JJ$. We now identify other subideals of $\JJ$.

\begin{proposition}\label{jrad1} Let $(\RR, \mathfrak{m})$ be a Noetherian local
ring and let $E$ be a finitely generated $\RR$-module.
\begin{enumerate}
\item[{\rm (a)}]
 If $E$ has no free summand, then the image of $E^*\otimes E$ in
$\CC$ is a two-sided ideal contained in the Jacobson radical.

\item[{\rm (b)}]  $\Hom_{\RR}(E, \m E)$ is a two-sided ideal contained in $\JJ$.
\end{enumerate}
\end{proposition}

\begin{proof} The actions of $\CC$ on $E^*\otimes E$ are as follows. For
$\hh\in \CC$, $(f\otimes e)\hh= f\circ \hh \otimes e$ and
$\hh(f\otimes e)= f\otimes \hh(e)$.
In other words, $E$ is a left $\CC$-module and
$E^*$ is a right $\CC$-module.

 Let $\II$ be the identity of $\CC$. To prove that
\[ \hh=\II + \sum_{i=1}^n f_i\otimes e_i\]
is invertible, note that for each $e\in E$,
\[ \hh(e)=e + \sum_{i=1}^n f_i(e) e_i\in e + \mathfrak{m}E,\]
since $f_i(e) \in \mathfrak{m}$ as $E$ has no free summand.
From the Nakayama Lemma, it follows that $\hh$ is a surjective
endomorphism, and therefore must be invertible.

\medskip

The proof of (ii) is similar.
\end{proof}

The ideal $E^*(E)= \tau(E)=(f(e), f\in E^*, e\in E)$ is the {\em trace
} of $E$. It is of interest, according to the next observation, 
 to describe classes of modules of
syzygies with $\tau(E)\subset \m$.

\begin{Remark}{\rm
From the sequence
\[ 0 \rar \m E \lar E \lar E/\m E=\bar{E}= k^n \rar 0,
\]
we have the exact sequence
\[
0 \rar \Hom_{\RR}(E, \m E) \lar \CC \lar M_n(k) \lar \Ext_{\RR}(E, \m E).
\]
This shows that $\dim \CC/\JJ\leq n^2$, since $\Hom_{\RR}(E,\m
E)\subset \JJ$, giving an explicit  control over the number of
maximal two-sided ideals of $\CC$.}
\end{Remark}

\subsubsection{Decomposition}

Let $(\RR,\m)$ be a local ring and $E$ a finitely generated $\RR$-module.
Suppose $E=E_1\oplus E_2$ is a nontrivial decomposition. We seek
to express the Jacobson of $\CC=\Hom(E,E)$ in terms  of those of the
subrings $\CC_i=\Hom(E_i,E_i)$, $i=1,2$ and the relationships
between $E_1$ and $E_2$.
A special   case is that
where $E_1= F=\RR^n$ and $E_2=M$  is a
module without free summands.
% (We look at $\End(E_1\oplus E_2)$ as
%$2\times 2$ matrices with entries in modules of homomorphisms.)

\begin{Proposition} Let $\RR$ be a local ring and $E$ a  module with
a decomposition $E=E_1\oplus E_2$. Set $\CC_i=\Hom_\RR(E_i, E_i)$ and
$\JJ'_i=\Hom_\RR(E_i, \m E_i)$. Assume the transition conditions
\[
\Hom_\RR(E_j, E_i)\cdot \Hom_\RR(E_i,E_j)\subset
\Hom_\RR(E_j, \m E_j), \quad i\neq j.\] Then the Jacobson
radical of $\CC=\End(E)$ is
\[ \JJ=\left[ \begin{array}{cc}
\JJ_1 &\Hom(E_2,E_1) \\
\Hom(E_1, E_2) & \JJ_2
\end{array}
\right],\]
where  $\JJ_i$ is the Jacobson radical of
$\End(E_i)$.
\end{Proposition}

\begin{proof} Note that $\JJ_i'$ is a subideal of $\JJ_i$, and for
$A=\Hom(E_1, E_2)$ and $B=\Hom(E_2, E_1)$ it holds that $A\cdot
B\subset \JJ_2'$ and $B\cdot A\subset \JJ_1'$. From these it follows
that
\[ \LL=\left[ \begin{array}{cc}
 \JJ_1 &  B \\
A &  \JJ_2
\end{array}
\right]\] is a two-sided ideal of $\CC$.

To show $\LL$ is the radical it is enough to show that for $\Phi\in
\LL$, $\II+\Phi$ is invertible, or equivalently it is surjective
endomorphism of $E$. In other words, for $a\in \JJ_1$, $b\in A$,
$c\in B$ and $d\in \JJ_2$ the system of equations
\begin{eqnarray*}
x + ax + by &=& u\\
c x + y + dy &=& v
\end{eqnarray*}
for $u\in E_1$, $v\in E_2$ has always a solution $x\in E_1$ and $y\in
E_2$.
It is enough to observe that in the formal solution
\begin{eqnarray*}
x & = & (\II_1 + a)^{-1}(u-by)\\
y &=& (\II_2 + d -c(\II_1 + a)^{-1}b)^{-1}(v-c(\II_1+a)^{-1}u),
\end{eqnarray*}
$c(\II_1+a)^{-1}b\in \JJ_2$.

\end{proof}

\begin{Corollary}\label{Jacsemi} Let $\RR$ be a local ring and $M$ a  module without
free summands.  If $F$ is a free $\RR$-module of rank $n$ then the Jacobson
radical of $\CC=\End(F\oplus M)$ is
\[ \JJ=\left[ \begin{array}{cc}
\m \cdot \End(F) & \Hom(M,F) \\
\Hom(F,M) & \JJ_0
\end{array}
\right],\]
where  $\JJ_0$ is the Jacobson radical of
$\End(M)$.
\end{Corollary}

This shows that a non-free  module with a free summand is never local.
It permits, to understand the Jacobson radical of $\End(E)$, to peel
away from $E$ a free summand of maximal rank.

\medskip

\subsubsection*{Modules without free summands}

Let $(\RR,\m)$ be a Cohen--Macaulay local ring of dimension $d$ and
$E$ a finitely generated $\RR$--module.

\begin{example}{\rm
We will now describe two classes of modules of syzygies that do not
have free summands.

\medskip

\begin{enumerate}

\item[{\rm (a)}] {\cite[Lemma 1.4]{Herzog78}} Let $\RR$ be a Cohen--Macaulay
local ring and $E$ a module of
syzygies
\[ 0 \rar E \lar F \lar M \rar 0,\]
$E\subset \m F$ and $M$ a maximal Cohen--Macaulay module. Then $E$
has no free summand.

\begin{proof} Suppose $E=\RR\epsilon \oplus E'$, $0\neq \epsilon\in \m F$.
Setting $M'=F/\RR\epsilon$ we have the exact sequence
\[ 0 \rar E' \lar M' \lar M \rar 0 \] showing that $M'$
is a maximal
Cohen--Macaulay module. But $M'$ is a module with a free resolution, so must be free
by Auslander-Buchsbaum equality. This means that $\RR\epsilon$ is a summand of $F$, which is impossible.

\end{proof}

\item[{\rm (b)}] Suppose $\RR$ is a Gorenstein local ring and $E$ is a module of
syzygies
\[ 0 \rar E \stackrel{\varphi}{\lar} F_s \lar \cdots \lar F_1 \lar F_0 \lar M \rar 0,
\] $E\subset \m F_s$. If $1\leq s\leq \height(\ann(M))-2$ then $E$
has no free summand.

\begin{proof}  $E$ is a reflexive being a second syzygy module over
the Gorenstein ring $\RR$. It will be enough to show the dual module
of $\RR\epsilon$ splits off $F_s^*$.

\medskip

By assumption $\Ext_{\RR}^i(M,\RR)=0$ for $i\leq s+1$. Applying
$\Hom(\cdot, \RR)$ to the complex we obtain an exact complex
\[ 0 \rar F_0^* \lar F_1^* \lar \cdots \lar F_s^*
\stackrel{\varphi^*}{\lar E^*} \rar 0,
\]
to prove the assertion since $\varphi=\varphi^{**}$.
\end{proof}

The assumption that $\RR$ is Gorenstein was used to guarantee that
second syzygies modules are reflexive. This could be achieved in many
other ways. For example, with $\RR$ Cohen--Macaulay and $\height
\ann(M)\geq 2$.

\medskip

\item[{\rm (c)}] Taken together these two observations give an overall picture
of the syzygies of $M$ which are without free summands.

\medskip

\item[{\rm (d)}] A special case is that of the syzygies in the Koszul complex of
a regular sequence.

\end{enumerate}
}
\end{example}

\subsection{Construction of local modules} \index{construction of
local modules}

The setup we employ is derived from the analysis of duality of
\cite{AusBr}: 
There is an exact complex of $\RR$-algebras
\begin{eqnarray} \label{ausbr0}
 E^* \otimes_{\RR} E \lar \Hom_{\RR}(E,E)=\End_{\RR}(E) \lar \Tor_1^{\RR}(D(E),E) \rar 0
\end{eqnarray}
where $D(E)$ is the Auslander dual of $E$. We note that $E^*\otimes_{\RR} E$ maps onto the endomorphisms of $E$
that factor through projective modules.
 
%For modules without free summands
%the  question
%turns on the  
We now turn to 
the understanding of $\Tor_1^{\RR}(D(E),E)$ as a more amenable  ring of  endomorphisms. 
 According to
\cite[Theorem 1.3]{HiltonRees} there is a natural surjective (anti-)
homomorphism \[ \Hom_{\RR}(E,E) \mapsto  \Hom_{\RR}(\Ext_{\RR}^1(E,\cdot),\Ext_{\RR}^1(E,\cdot)).\]
Its kernel consists of the maps $f: E \rar E$ that factor through projective modules, a fact  that can be seen
with a brief calculation by considering a presentation $0 \rar M \rar P \rar E \rar 0$, $P$ projective,  and assuming 
for $f: E \rar E$ that $f^*: \Ext_{\RR}^1(E,M) \rar \Ext_{\RR}^1(E,M)$ vanishes. 
The useful point is that 
the  
 endomorphism ring 
 $\Hom_{\RR}(\Ext_{\RR}^1(E,\cdot),\Ext_{\RR}^1(E,\cdot))$ may have a much smaller support than $E$ itself. Furthermore
 if $E$ has no free summands, $E^*\otimes_{\RR} E$ maps into the Jacobson radical of $\Hom_{\RR}(E,E)$ and therefore
$E$ is a local module if 
$\Hom_{\RR}(\Ext_{\RR}^1(E,\cdot),\Ext_{\RR}^1(E,\cdot))$ has a unique maximal two-sided ideal.

\medskip

Let us cast these observations in a format convenient for reference:

\begin{Theorem}\label{AusEndo} Let $\RR$ be a commutative Noetherian ring and $E$ a finitely generated 
$\RR$-module. 
There is an exact complex of $\RR$-algebras
\begin{eqnarray*} \label{ausbr01}
 E^* \otimes_{\RR} E \lar \Hom_{\RR}(E,E)\lar \Hom_{\RR}(
 \Ext_{\RR}^1(E, \cdot), \Ext_{\RR}^1(E, \cdot))
  \rar 0,
\end{eqnarray*}
where  $E^*\otimes_{\RR} E$ maps onto the endomorphisms of $E$
that factor through projective modules.
\end{Theorem}

In each application of this setup we look for which of
 $\Hom_{\RR}(\Ext_{\RR}^1(E,\cdot),\Ext_{\RR}^1(E,\cdot))$ or
$\Tor_1^\RR(D(E),E)$ is more amenable. Let us consider two examples.
%\bigskip

%\bigskip

\begin{example} 
{\rm Let $\RR$ be a Cohen--Macaulay local ring of dimension
$d\geq 3$,
and  let $E$
 be a module defined by one relation
\begin{eqnarray}\label{onerelmod}
 0\rar \RR \stackrel{\varphi}{\lar} \RR^n \lar E \rar 0. 
\end{eqnarray}
We assume that the entries of $\varphi$ define an ideal $I$ of height
$\geq 3$, minimally generated by $n$ elements. This makes $E$ a reflexive module without free summands.
Let us determine the
number of generators of $\CC=\Hom_{\RR}(E,E)$, and some of its other
properties.
\medskip

The functor $\Ext_{\RR}^1(E,X)= \Ext_{\RR}^1(E, \RR)\otimes X= \RR/I\otimes
X$, so that
  \[\Hom_{\RR}(\Ext_{\RR}^1(E,\cdot),\Ext_{\RR}^1(E,\cdot))=\Hom_{\RR}(\RR/I,\RR/I)=\RR/I.\]
Tensoring (\ref{onerelmod}) by $E^*$, we get the complex
\[ 0 \rar E^* \lar E^*\otimes \RR^n \lar E^*\otimes E \rar 0,\]
which is exact by acyclicity lemma.  Since $E^*$ is reflexive and $\RR$ is Cohen-Macaulay, $E^*$ has the
condition $S_2$ of Serre and therefore $E^*\otimes E$ has the condition $S_1$ of Serre, in particular it
is torsionfree.
Since $E^*$ has no
free summand,  by Proposition~\ref{jrad1} it follows that $E$ is a
local module. As for the number of generators of $\CC$,
\[\beta_0(I)\beta_1(I)-\beta_0(I)+1 \leq
 \nu(\CC)\leq \nu(E)\nu(E^*)+1= \beta_0(I)\beta_1(I)+1,\]
 where $\beta_i(\cdot)$ denotes Betti numbers.}
\end{example}

\subsubsection*{Syzygies of perfect modules}
\index{perfect module}

 Let $\RR$ be a Gorenstein  local ring of
dimension $d$.
Let us consider some modules with
a very rich structure--the modules of syzygies of Cohen-Macaulay
modules, or of mild generalizations thereof.

\medskip

Let us begin this discussion with an example,
 the modules of cycles of a Koszul complex
$\KK(\xx)$ associated to a regular sequence $\xx=\{x_1, \ldots,
x_n\}$, $n\geq 5$:
%\begin{small}
\[ \KK(\xx): \quad 0 \rar K_n \rar K_{n-1} \rar K_{n-2} \rar \cdots
\rar K_2\rar K_1\rar K_0 \rar 0.\]
%\end{small}

%First order syzygies have a more general treatment in the next
%section.
For simplicity we take for module $E$ the $1$-cycles $Z_1$ of $\KK$.
There is a pairing in the subalgebra  $\ZZ$ of cycles leading to
\[ Z_1 \times Z_{n-2}\rar Z_{n-1} = \RR,\]
that identifies $Z_{n-2}$ with the dual $E^*$ of $E$.

We are now ready to put this data into the framework of the
Auslander dual. Dualizing the projective presentation of $E$,
\[ 0 \rar K_n \lar \cdots \lar K_{3}\lar K_{2} \lar E=Z_1 \rar 0,\]
gives us the exact  complex
\[0 \rar E^* \lar K_{2}^* \lar K_{3}^* \lar D(E) \rar 0.\]
In other words,  to the identification  $D(E)=Z_{n-4}$:
\[ 0 \rar Z_{n-2} \lar K_{n-2} \lar K_{n-3} \lar Z_{n-4}  \rar 0.\]

Now for the computation of $\Tor_1^{\RR}(D(E),E)$:

\begin{eqnarray*} \Tor_1^{\RR}(D(E),E)&=& \Tor_1^{\RR}(Z_{n-4},E) =
\Tor_2^{\RR}(Z_{n-5}, E) =
\cdots \\
&=& \Tor_{n-3}^{\RR}(Z_0,E)= \Tor_{n-2}^{\RR}(\RR/I,E) = K_{n}\otimes
\RR/I=\RR/I.
\end{eqnarray*}

\begin{Remark}{\rm The number of generators of $\Hom_{\RR}(E,E)$ is bounded
by
\[ \nu(E)\nu(E^*) + 1 .\]

For the purpose of a comparison, let us evaluate $\hdeg (E)$. The
multiplicity of $E$ is $(n-1)\deg(\RR)$. Applying $\Hom_{\RR}(\cdot, R)$ to the
projective resolution of $E$, we get
 $\Ext_{\RR}^{n-2}(E,\RR) =\Ext_{\RR}^n(\RR/I,\RR) = \RR/I $ is Cohen-Macaulay,
and its contribution in the formula for $\hdeg (E)$ becomes
\[ \hdeg (E) = \deg (E) + {{d-1}\choose{n-3}}
\hdeg (\Ext_{\RR}^{n-2}(E,\RR))= (n-1) +
{{d-1}\choose{n-3}}
\deg (\RR/I).
\]

}\end{Remark}

It is clear that in appealing to \cite[Theorem 5.2]{Dalili}, to get
information about $\nu(E^{*})$, a similar calculation can be carried
out for any module of cycles of a projective resolution of broad
classes
of Cohen-Macaulay modules.

\medskip

Let $(\RR, \mathfrak{m})$ be a Gorenstein local ring of dimension $d$.
Let $M$ be a perfect $\RR$-module with a minimal free resolution

\[ \FF: \quad 0 \rar F_n \rar F_{n-1} \rar F_{n-2} \rar \cdots
\rar F_2\rar F_1\rar F_0 \rar M \rar 0.\]
We observe that dualizing $\FF$ gives a minimal projective resolution
$\LL$ of $\Ext_{\RR}^{n}(M,\RR)$.
Let $E$ be the module $Z_{k-1}=Z_{k-1}(\FF)$ of $k-1$-cycles of $\FF$,
\[ F_{k+1} \lar F_{k} \lar E \rar 0.\]
Dualizing, to define the Auslander dual $D(E)$, gives the complex
\[ 0 \rar E^{*} \lar L_{n-k} \lar L_{n-k-1} \lar D(E)\rar 0.\]
It identifies $E^*$ with the $(n-k)$--cycles of $\LL$, and $D(E)$ with
its $(n-k-2)$--cycles. In particular this gives $\nu(E)=
\beta_{k}(M)$ and $\nu(E^*)=\beta_{k-1}(M)$.

\medskip

Let us now determine  $\hdeg (E)$. From the complex
\[ 0 \rar E \lar F_{k-1} \lar \cdots \lar F_1 \lar F_0 \lar M \rar 0,\]
we obtain that $\mbox{\rm proj. dim. } E = n-k$ and that the
Cohen-Macaulay property of $M$ gives that
$\Ext_{\RR}^n(M,\RR)= \Ext_{\RR}^{n-k}(E,\RR) $ is also Cohen-Macaulay. This
 is all that is required:
\begin{eqnarray*} \hdeg (E)& = &\deg (E) + {{d-1}\choose{n-k-1}} \deg
(\Ext_{\RR}^n(M,\RR))\\
&=& \rank(E) \deg (\RR) +  {{d-1}\choose{n-k-1}} \deg (M).
\end{eqnarray*}
Note that $\rank(E)= \beta_0(M)-\beta_1(M) + \cdots
(-1)^{k-1}\beta_{k-1}(M)$.

\medskip

Now
to make use of the Auslander dual setup, we seek some
control over $\Tor_1^{\RR}(D(E),E)$. We make use first of the complex
\[ 0 \rar D(E) \lar L_{n-k-2} \lar \cdots \lar L_0 \lar
\Ext_{\RR}^n(M,\RR)
\rar 0,\]
to get
\[ \Tor_1^{\RR}(D(E),E) \simeq \Tor_{n-k}^{\RR}(\Ext^n(M,\RR),E), \]
and then of the minimal resolution of $E$  to obtain
\[ \Tor_1^{\RR}(D(E),E) \simeq \Tor_{n}^{\RR}(\Ext^n(M,R),M). \]
In particular, $\Tor_1^{\RR}(D(E),E)$ is independent of which module of
syzygies was taken. Furthermore, the calculation  shows that
 $\Tor_2^{\RR}(D(E),E)=0$.

Placing these elements together, we have the exact sequence
\[ 0 \rar E^*\otimes E \lar \Hom_{\RR}(E,E) \lar
\Tor_n^{\RR}(\Ext_{\RR}^n(M,\RR),M)
\rar 0.\] Note that $E^*\otimes E$ is a torsion free $\RR$-module.

Finally, it follows from Proposition~\ref{ausdual2} that
$\Tor_n^{\RR}(\Ext_{\RR}^n(M,\RR),M)$
can be identified to $\Hom_{\RR}(M,M)$.

\medskip

Let us sum up these observations in the following:

\begin{Theorem} \label{syzofperfect} Let $\RR$ be a Gorenstein local ring and  let $M$ be a
perfect module with a minimal resolution $\KK$. For the module $E$ of
$k$-syzygies of $M$, $1<k<n$,   there exists an exact sequence
\[ 0 \rar E^*\otimes E \lar \Hom_{\RR}(E,E) \lar
\Hom_{\RR}(M,M) \rar 0.\]
\end{Theorem}

This gives the bound
\[ \nu(\Hom_{\RR}(E,E)) \leq \beta_{k}(M)\beta_{k-1}(M)+\nu(\Hom_{\RR}(M,M)).\]
If $M$ is cyclic, or $\dim M = 2 $, the estimation is easy. The
formula also shows up the case of MCM modules to be a corner case for
the general HomAB problem.

\bigskip

\subsubsection*{Syzygies of cyclic modules}
A more general class of local modules arises  as follows.
We recall Theorem~\ref{syzofperfect}:

\begin{Theorem} \label{syzofperfect2} Let $\RR$ be a Gorenstein local ring and  let $M$ be a
perfect module with a minimal resolution. For the module $E$ of
$k$-syzygies of $M$, $k< \mbox{\rm proj dim $M$}$,   there exists an exact sequence
\[ 0 \rar E^*\otimes E \lar \Hom_{\RR}(E,E) \lar
\Hom_{\RR}(M,M) \rar 0.\]
In particular, if $M=\RR/I$, $E$ is local if it has no free summand.
\end{Theorem}

\subsubsection*{Generic modules}\index{generic module}
Let $\varphi$ be a generic $n\times m$ matrix, $m>n$ and consider the
module $E$
\[ 0\rar \RR^n \stackrel{\varphi}{\lar} \RR^m \lar E \rar 0.
\]

We study the ring $\CC=\Hom_\RR(E,E)$. Let $I=I_n(\varphi)$ be the
ideal of maximal minors of $\varphi$. $\RR/I$ is a normal
Cohen--Macaulay domain and many of its properties
are deduced from the Eagon--Northcott complexes (\cite{EN}). We are
going to use them to obtain corresponding properties of $E$:

\begin{itemize}
\item[$\bullet$] If $n<m\geq 3$ (which we assume throughout), $E$ is a reflexive module.
\item By dualizing we obtain the
Auslander dual $D(E)$
 of $E$:
\[ 0\rar E^* \lar \RR^m \stackrel{\varphi^*}{\lar} \RR^n \lar
D(E)=\Ext_\RR^1(E, \RR)\rar 0.
\]

\item[$\bullet$] $\Ext_\RR^1(E,\RR)$ is isomorphic to an ideal of $\RR/I$:
$\Ext_{\RR}^1(E,\RR)$ is a perfect $\RR$--module, annihilated $I$.
Thse assertions follow directly from the perfection of the complexes.
Finally, localizing at $I$ the complex is quasi-isomorphic to Koszul
complex of a regular sequence which yield
\[ \Ext_\RR(E,\RR)_I \simeq (\RR/I)_I.\]

\end{itemize}

\begin{Theorem} $\CC=\Hom_{\RR}(E,E)$ is a local ring generated by
$\nu(E^*)\nu(E)+1=n\times {m\choose n+1}+1$ elements.
\end{Theorem}

\begin{proof}
As in the discussion of example above,
The functor $\Ext_{\RR}^1(E,X)= \Ext_{\RR}^1(E, \RR)\otimes X= \RR/I\otimes
X$, so that
  \[\Hom_{\RR}(\Ext_{\RR}^1(E,\cdot),\Ext_{\RR}^1(E,\cdot))
=\Hom_{\RR/I}(\Ext_\RR^1(E,\RR),\Ext_\RR^1(E,\RR))
=\RR/I,\]
as $\Ext_\RR^1(E,\RR)$ is isomorphic to an ideal of the normal domain
$\RR/I$.

Furthermore, since $E$ is reflexive and $E^*$ is a module of syzygies
in the Eagon--Northcott complex and clearly has no free summand, $E$
has the same property.
The final assertion comes from the complex again.
\end{proof}

%\bigskip

%\subsubsection*{Modules of dfferentials}

\subsection{Homological properties of local modules}

Let $(\RR,\m)$ be a Noetherian local ring and $E$
a finitely generated $\RR$-module. Set $\CC=\End(E)$ and  consider
the natural  action of $\CC$  making $E$ a $\CC$--module. In this section $\mbox{\rm mod}(\CC)$ denotes the category of finitely generated
 left $\CC$-modules.

\medskip
\subsubsection*{Projective generation}
%We want to determine the modules with the following property.
Let us begin with an observation.

\begin{Proposition} \label{projgen}
Let $E$ be a local module and $\CC$ its ring of endomorphisms.
If $M$ is a nonzero  projective $\CC$-module, then $M$ is
a generator of $\mbox{\rm mod}(\CC)$.
\end{Proposition}

\begin{proof}
 Let $\JJ$ be the Jacobson
radical of $\CC$. Up to isomorphism there is only one indecomposable
$\CC/\JJ$-module, which we denote by $L$.

Let $F$ be a finitely generated left $\CC$--module. Then $M/\JJ M$
and $F/\JJ F$ are isomorphic to direct sums $L^n$ and $L^m$ of $L$.
It follows that for some integer $r>0$ there is a surjection
\[ \varphi: M^r \rar F/\JJ F. \]
Since $M$ is projective, $\varphi$ can lifted to a mapping
$\Phi: M^r\rar F$ such that $F = \Phi(M^r) + \JJ F$. By Nakayama
Lemma $\Phi$ is surjective, as desired.
\end{proof}

Let us make a related observation. Let $(\RR, \m)$ be a Noetherian local ring and $\CC$ be a finite $\RR$-algebra
that has a unique $2$-sided maximal ideal $\JJ$. 

\begin{Proposition} \label{ABeq}
Suppose $\m$--depth of $\CC$ is zero. If $F$ is a finitely generated left $\CC$-module of finite
projective dimension, then $F$ is $\CC$-projective.
\end{Proposition}

\begin{proof} Let us argue by contradiction. We may assume that $F$ has projective dimension $1$. In the proof of 
Proposition~\ref{projgen} applied to $\CC$ itself, from $\CC/\JJ=L^n$ and $F/\JJ F=L^m$, we have an 
isomorphism $(\CC/\JJ)^m \simeq	 F^n/\JJ F^n$ which can be lifted into a surjection
\[ \CC^m \rar F^n \rar 0.\]
The kernel, by assumption, is a projective $\CC$-module $P $ that is contained 
in $ \JJ \CC^m$.
We note that since $P$ is a generator its left annihilator is trivial.

Let $p$ the least integer such that $\JJ^p $ has a left annihilator. We claim that such integer exists.
Since $\CC/\m \CC$ is Artinian, there exists $r$ with $\JJ^r\subset \m \CC$. But $\CC$ has depth zero by
hypothesis so there is a nonzero $w\in \CC$ with $\m w=0$. This element will annihilate $\JJ^r$ and thus the
integer $p$ exists.
From $\JJ^{p-1} P \subset \JJ^p \CC^m$ it will follow that $\JJ^{p-1} P$ has a nonzero left annihilator. 
Since $P$ is a generator, there is a surjection $P^s\rar \CC\rar 0$ and thus another 
surjection $\JJ^{p-1} P^s \rar \JJ^{p-1}\CC \rar 0$, which will contradict the choice of $p$.
\end{proof}

This opens the way of the extension of
the Auslander-Buchsbaum equality 
(\cite[Theorem 1.3.3] {BH}) to 
  $\CC$-modules, that is if $A$ is a finitely generated left $\CC$-module of
finite projective dimension  then
   \[ \mbox{\rm proj dim}_{\CC} A + \mbox{\rm depth } A = \mbox{\rm depth } \CC.\]

The proof follows exactly the same steps of the commutative case, with Proposition~\ref{ABeq} providing the base case in the induction process:

\medskip

\begin{itemize}

\item[{$\bullet$}] We may assume $\depth \CC > 0$. Let $A$ be a finitely generated $\CC$-module of
projective dimension $r$, $1\leq r < \infty$. Consider an exact sequence $0 \rar P \rar F \rar A \rar 0$,
with $F$ a projective $\CC$-module.

\item[{$\bullet$}] $P$ has projective dimension $r-1$. Let us argue that $\depth P = \depth A + 1$. If
$\depth A =0$, let $x\in \m$ be a regular element on $\CC$,
and therefore a regular element on $P$ as well.
 Reduction mod $x$ gives rise the an exact sequence 
 \[ 0 \rar 0:_A x \rar P/xP \rar \CC/x\CC \rar A/xA \rar 0.\]
Since $\depth A =0$, $0:_A x$ contains elements annihilated by $\m$. This shows that $\depth P/xP = 0$,
and therefore $\depth P =1$. Similarly, in general we have $\depth P = \depth A + 1$. 

\item[{$\bullet$}] By induction the Auslander-Buchsbaum equality holds for the $\CC/x\CC$-module $P/xP$. It
remains to show that 
$\mbox{\rm proj dim}_{\CC/x\CC} P/xP = 
\mbox{\rm proj dim}_{\CC} P $. To deny equality it suffices to work with the penultimate module of 
$M$ of syzygies of $P$, a module of projective dimension $1$ and assume that $M/xM$ is $\CC/x\CC$-projective.
But then $\mbox{\rm prom dim}_{\CC} M/xM \leq 1$ which is readily contradicted by
the exactness of 
\[\Ext_{\CC}^1(M, X) \stackrel{\cdot x}{\lar} \Ext_{\CC}(M, X) \lar \Ext_{\CC}^2(M/xM,X) \rar 0,\]  
 the non-vanishing of $\Ext_{\CC}^1(M,X)$ and Nakayama Lemma.

\end{itemize}

\begin{Theorem} \label{homollocal} Let $(\RR, \m)$ be a normal Noetherian local ring
and $E$ a local
module.
 If $E$ is $\CC$--projective and
it is $\RR$-torsionfree whose rank is invertible in $\RR$
   then $E$ is $\RR$-free.
\end{Theorem}

\begin{proof}

According to Corollary~\ref{Jacsemi},
 $E$ has no free summand.
By Proposition~\ref{projgen}, $E$ is a projective generator and in
particular we have an isomorphism of left $\CC$-modules
\begin{eqnarray}\label{stablegen}
 E^s \simeq \CC\oplus H.
\end{eqnarray}
%Since $\RR$ is a local ring this gives rise to a decomposition of
%$\RR$-modules,
%\[ E^r \simeq \CC,\]
%where $r$ is necessarily
%the rank of $E$ as an $\RR$-modules.
% If $\JJ$ is the Jacobson radical of $\CC$, we have
%\[ E^s/\JJ E^s \simeq \CC/\JJ\oplus H/\JJ H.\] By Krull-Schmidt, this leads
%to an isomorphism of $\CC$--modules
%\[ E^r/\JJ E^r  \simeq \CC/\JJ,\] which lifts to an isomorphism of $\CC$--modules
%\[ E^r \simeq \CC\].

\medskip

To prove the assertion, we recall how the {\em trace} of the elements of $\CC$
may be defined. Let $\SS$ be the field of fractions of $\RR$.
Consider that canonical ring homomorphism
\[ \varphi: \End_\RR(E)\lar \End_\SS(\SS\otimes E).
 \] For $\ff\in \End_\RR(E)$, define
\[ \mathrm{tr}(\ff)= \mathrm{tr}_\SS(\varphi(\ff)),\]
where $\mathrm{tr}_\SS$ is the usual trace of a matrix
representation. It is independent of the chosen $\SS$-basis of
$\SS\otimes E$.

To show that $\mathrm{tr}(\ff)\in \RR$ uses a standard argument.
 For each prime ideal
$\p$ of $\RR$ of height $1$, $E_{\p}$ is a free $\RR_{\p}$-module of rank
$r=\mbox{\rm rank }(E)$,  and picking one of its basis gives also a $\SS$-basis for
$\SS\otimes E$ while showing that $\mathrm{tr}(\ff)\in \RR_{\p}$. Since
$\RR$ is normal, $\bigcap_{\p} \RR_{\p}=\RR$, and thus $\mathrm{tr}(\ff)\in
\RR$.

Note that this defines an element of $\Hom_\RR(\CC, \RR)$.
Now looking at the isomorphism (\ref{stablegen}), yields that the
trace of $\CC$ as an $\RR$-module is contained in the trace
ideal of $E$, which is a proper ideal of $\RR$ since $E$ has no $\RR$-free summand.
But this is impossible since $\mathrm{tr}(\II)=r$, which is a unit of
$\RR$ by assumption.

\end{proof}

\begin{Corollary}\label{mainlocalcor}
Let $\RR$ be a Noetherian local ring and  $E$  a local  $\RR$-module such that $\CC$ has finite global dimension.
\begin{enumerate}
\item[{\rm (a)}] \mbox{\rm (\cite[Theorem 2.5]{BH84}, \cite[Theorem 3.1]{RegAlg})}
  $\CC$ is a maximal Cohen-Macaulay
$\RR$-module.

\item[{\rm (b)}]  If $E$  is a maximal Cohen-Macaulay $\RR$--module then it is a projective $\CC$--module.

\item[{\rm (c)}] Moreover, if
$\RR$ is integrally closed  and
  the rank $\mathrm{r}$ of $E$ is invertible in
 $\RR$  then $\RR$ is a regular local ring.
\end{enumerate}
\end{Corollary}

\begin{proof} (a)
Since $\CC$ is a local  ring, that is has a unique maximal two-sided ideal,   by  \cite[Theorem
2.5]{BH84} or \cite[Theorem 3.1]{RegAlg}, $\CC$
 is a Cohen-Macaulay $\RR$-module. 
 Let us argue this point. Let $\mbox{\rm gl dim } \CC=n$. Then by the Auslander-Buchsbaum inequality,
 $\depth \CC= n$.
 Let $\xx=\{x_1, \ldots, x_n\}\subset \m$ be a regular sequence on $\CC$
 and set $\BB = \CC/(\xx)\CC$.
 For any $\CC$-module, $\Ext_{\CC}^{n+1}(M, \CC)=0$ since $\mbox{\rm proj dim}_{\CC}(M)\leq n$.
 For any $\BB$-module $M$, by Rees's theorem (\cite{Rees56}),
 \[ \Ext_{\BB}^1(M, \BB) = \Ext_{\CC}^{n+1}(M,\CC) =0,\]
 which implies that $\BB$ is a left injective $\BB$-module.  
 
 We claim that $\BB$ is an Artinian module. Let $\BB_0= \H_{\m}^0(\BB)$ and set $\BB' = \BB/\BB_0$. 
  If  $\BB'\neq 0$ there is an exact sequence
 of $\BB$-modules  
  \[ 0 \rar \BB' \stackrel{\cdot x}{\lar} \BB' \lar \BB'/x \BB'\rar 0.\] 
  Applying $\Hom_{\BB}(\cdot, \BB)$ we have a surjection
  \[ \Hom_{\BB}(\BB', \BB) \stackrel{\cdot x}{\lar} \Hom_{\BB}(\BB', \BB)\rar 0,\] 
  and thus $\Hom_{\BB}(\BB', \BB)=0$ by Nakayama lemma. But this module cannot vanish since $\BB$
  contains a simple submodule $L$ and there are surjections
  $\BB' \rar \BB'/\JJ \BB' \rar L$, and $\Hom_{\BB}(\BB',L) \subset \Hom_{\BB}(\BB', \BB)$.
  
  \medskip
  
  (b) Let
  $\xx=x_1, \ldots, x_n$ be a maximal
$\RR$-sequence for $E$ and $\CC$, $n= \mbox{\rm gl dim }\CC$,
 $E/(\xx)E$ is a $\CC$-module with nonzero socle. It follows from
Auslander-Buchsbaum formula 
that
\[ \mbox{\rm proj dim}_\CC E/(\xx)E=
 \mbox{\rm proj dim}_\CC E + n =\mbox{\rm gl dim }\CC=n.\]

(c)
By Theorem~\ref{homollocal} $E$ is a free $\RR$-module. This
means that
$\CC$ is a matrix ring over $\RR$ of finite global dimension
 and therefore $\RR$ is a regular
local ring.
\end{proof}

\begin{Conjecture} Let $E$ be a local $\RR$--module. If $E$ is $\CC$--cyclic of finite projective dimension over $\CC$, then $E$ is $\CC$--projective.
\end{Conjecture}

\begin{Question}{\rm Let $\AA$ be an affine algebra over a field $K$ of characteristic  zero and $E$ be either
(i) the module of $K$-differentials $\Omega_K(\AA)$, or (ii)  the module of $K$-derivations 
$\mbox{\rm Der}_K(\AA)$. If $E_\pp$ is  a local $\AA_\pp$ module for each prime $\pp$, is $E$ a projective $\AA$-module?
}\end{Question}

\chapter{The Equations of Ideals}\index{equation of an ideal}

%\section*{Introduction}

%(\section{The Equations of Filtrations of Ideals} \label{eqnf}

\section{Introduction}

\noindent
Let $\RR$ be a Noetherian ring  and let $I$ be an
 ideal. By the {\em equations} of $I$ is meant a free presentation of
 the Rees algebra $\RR[It]$ of $I$,
\[ 0 \rar \LL \lar \SS = \RR[\TT_1, \ldots, \TT_m] \stackrel{\psi}{\lar}
\RR[It] \rar 0,  \quad \TT_i \mapsto f_it .\]
More precisely, we refer to the ideal $\LL$ as the {\em ideal of
equations} of $I$.\index{equations of an ideal}
The ideal $\LL$ depends on the chosen set of generators of $I$, but all
of its significant cohomological properties, such as the integers that bound the
degrees of minimal generating sets of $L$, are independent of the presentation
$\psi$.
The examination of $\LL$ is one pathway to the
unveiling of
the properties of $\RR[It]$.
It codes the syzygies of all powers of $I$, and therefore
 is a carrier of not just algebraic properties of $I$, but of  analytic
ones
as well.
This is a notion that can be extended to Rees algebras $\Rees= \sum_{n\geq
0}I_nt^n$ of more general
filtrations. Among these  we will be interested in the integral
closure of $\RR[It]$.
Our approach to $\LL$ involves the examination of some its graded components, particularly $L_1$
 and $L_2$, and some corresponding quotients, such as $T=\LL/(L_1)$, the so-called
 {\em module of nonlinear equations}\index{module of nonlinear equations} of $I$.

\medskip

We are going to make a hit list of properties that would be desirable
to have about $\LL$, and discuss its significance.
Foremost among them  is to estimate the degrees of the generators of
$L$,
the so-called {\em relation type} of the ideal $I$: $\reltype(I)$.
If $\LL$ is generated by forms of degree $1$, $I$ is said to be of linear type
 (this is independent of the set of generators). The Rees algebra
 $\RR[It]$ is then the symmetric algebra $\Symi=\Sym(I)$ of $I$.
Such is the case when
 the $f_i$ form a regular sequence, $\LL$ is then generated by the
 Koszul forms $f_i\TT_j-f_j\TT_i$, $i<j$.
There are other kinds of kinds of sequences with this property (see
\cite{HSV1} for a discussion).

\medskip

The general theme we are going to pursue is the following: How to
bound the relation type of an ideal in terms of its numerical
invariants known {\it a priori}? Given the prevalence of various
multiplicities among these invariants most of our questions depend on
them.

\begin{enumerate}
\item[{\rm (a)}] Let $(\RR, \mathfrak{m})$ be a Noetherian local ring of dimension
$d$ and let $I$
be an $\mathfrak{m}$-primary ideal. Is there a low polynomial $f(e_0,e_1,
\ldots, e_d)$ on the Hilbert coefficients of $I$ such that
\[ \reltype(I)\leq f(e_0,\rme_1, \ldots, e_d)?\]

\item[{\rm (b)}] Among the equations in $L$, an important segment arises from
the {\em reductions} of $I$. That is, subideals $J\subset I$ with the
property that $I^{r+1}= JI^r$, for some $r$, the so-called {\em
reduction number} of $I$ relative to $J$; in notation, $\red_J(I)$.
 Minimal reductions are
important choices here and there are various bounds for the
corresponding values of $\red_J(I)$ ($\red(I)$ will be the minimum
value of $\red_J(I)$).  Among these formulas, for Cohen-Macaulay
local rings $(\RR, \mathfrak{m})$ of dimension $d$, and
$\mathfrak{m}$-primary ideal $I$, one has
\[ \red(I) \leq \frac{d}{o(I)} e_0(I)-2d+1,
\] where $o(I)$ is the largest integer $s$ such that $I\subset
\mathfrak{m}^s$. These formulas are treated in \cite[Section
2.2]{icbook}. There are also versions of these formulas for arbitrary Noetherian local rings
(\cite{chern5}).

\item[{\rm (c)}] If the ideal $I$ is not $\mathfrak{m}$-primary, which {\em
multiplicities} play the roles of the Hilbert coefficients?
 We plan to carry this out using {\em extended multiplicities}.
It will require a broadening of our understanding of related issues:
Sally modules of reductions.

\item[{\rm (d)}] A related question, originally raised by C. Huneke, attracted
considerable attention: whether for ideals $I$ generated by systems
of parameters, there is a uniform bound for $\reltype(I)$.
Although the answer was ultimately negative,
significant classes of rings do have this property. Wang obtained
several beautiful bounds for rings of low dimensions
(\cite{Wang97},\cite{Wang97a}), while Linh and Trung
(\cite{LinhTrung}) derived uniform bounds for all generalized
Cohen-Macaulay rings.

\item[{\rm (e)}]
Let $\mathcal{F}$ be the special fiber of $\Rees$,
that is $\mathcal{F}= \Rees\otimes \RR/\mathfrak{m}$. Since
\[\reltype(\mathcal{F})\leq \reltype(\Rees),\] it may be useful to
 constrain
 $\reltype(\mathcal{F})$ in terms of the Hilbert coefficients of the
 affine algebra $\mathcal{F}$. A rich vein to examine is that of
 ideals which are near complete intersections, more precisely in case
 of ideals of deviation at most two.
Among other questions, one should examine the problem of {\em
variation} of the Hilbert polynomials of ideals $I\subset I' \subset
\overline{I}$, where the latter is the integral closure of $I$. For
ideals of finite co-length, one of the questions (given the equality
$\rme_0(I)=e_0(I')$), is to study the change $\rme_1(I')-e_1(I)$ in terms
of the degrees of $\mathcal{F}$ and $\mathcal{F}'$.

\item[{\rm (f)}] Almost all the previous questions can be asked for Rees
algebras of modules or even of filtrations whose components are
modules. The requisite multiplicities have already been built.

\item[{\rm (g)}]
We will treat mainly ideals of small deviations--particularly
  ideals of codimension $r$ generated
by $r+1$ or $r+2$ elements--because of their ubiquity in the
construction
of birational maps.
Often, but not
  exclusively, $I$ will be an ideal of finite
co-length in a local ring, or in a ring of polynomials over a field.
Our focus on $\Rees$ is shaped by the following fact. The class of
ideals $I$ to be considered will have the property that both  its
symmetric algebra $\Sym(I)$ and the normalization $\Rees'$ of $\Rees$
have amenable properties, for example, one of them (when not both) is
Cohen-Macaulay. In such case,  the diagram
\[ \Symi \twoheadrightarrow \Rees \subset \Rees'\]
gives a convenient platform from which to examine $\Rees$.

%\item[{\rm (h) }]

\end{enumerate}

%There are several motivations for looking at [and for] these
%equations. We are going to list some:

%\begin{itemize}

%\item[{$\bullet$}] Ideals which are almost complete intersections
%occur in some of the more notable birational maps and in geometric
%modeling.

%\item[{$\bullet$}] Their Rees algebras and special fibers interact
%more closely than for more general classes of ideals [Hilbert
%functions, etc].

%\item[{$\bullet$}] Some questions raised by David Cox at a recent
%talk in Luminy. They are addressed in \cite{syl1}.
%They are part of a general program of finding algorithms that produce
%all the equations of an ideal, or at least some distinguished
%polynomial (e.g. the `elimination equation' in it).

%\item[{$\bullet$}] Will attempt to give a local version of
%birationality.

%\end{itemize}

\subsection{Syzygies and Hilbert functions}

Let $\RR$ be a Noetherian ring  and let $I$ be an
 ideal. By the {\em equations} of $I$ it is meant a free presentation of
 the Rees algebra $\RR[It]$ of $I$,
\begin{eqnarray}
 0 \rar \LL \lar \SS = \RR[\TT_1, \ldots, \TT_m] \stackrel{\psi}{\lar}
\RR[It] \rar 0,  \quad \TT_i \mapsto f_it .
\end{eqnarray}
More precisely,  $\LL$ is the defining ideal of the Rees
algebra of $I$ but we refer to it simply as the {\em ideal of
equations} of $I$.
The ideal $L$ depends on the chosen set of generators of $I$, but all
of its significant cohomological properties, such as the integers that bound the
degrees of minimal generating sets of $L$, are independent of the presentation
$\psi$.
The examination of $\LL$ is one pathway to the
unveiling of
the properties of $\RR[It]$.
It codes the syzygies of all powers of $I$, and therefore
 is a carrier of not just algebraic properties of $I$, but of  analytic
ones
as well. It is also a vehicle to understanding geometric properties of
several constructions built out of $I$, particularly of rational maps.

The search for these equations and their use has attracted
considerable interest by a diverse group of researchers. We just
mention some that have  directly influenced this work. One main
source lies in the work of L. Bus\'e, M. Chardin, D. Cox, J. Hong,
 J. P.  Jouanolou  and A. Simis who
have charted, by themselves or with co-workers, numerous roles of
resultants and other elimination techniques to obtain these
equations (\cite{BuCh}, \cite{BuCoDa}, \cite{BuJou}, \cite{Dandrea},
and references therein.) Another important development was given by
A. Kustin, C. Polini and B. Ulrich, who provided a comprehensive
analysis of the equations of ideals (in the binary case), not
just necessarily of almost complete intersections, but still of ideals
whose syzygies   are
almost all linear (\cite{KPU}). Last,
has been the important work of D. Cox not only for its theoretical value to the
understanding of these equations, but for the  role
it has played in bridging the fields of commutative algebra and
of geometric modelling (\cite{CHW}, \cite{Cox08}).

The material here is based on \cite{syl1} and \cite{syl2}.
  It
deals
with  questions in higher dimensions that
 were triggered  in \cite{syl1}, but is mainly concerned with the
more general issues of the structure of the Rees algebras of almost
complete intersections.
Our underlying metaphor  here is to focus on distinguished sets of
equations by examining 4 Hilbert functions associated to the ideal
$I$ and to the coefficients of its syzygies. It brings considerable
effectivity to the methods by developing explicit formulas for some
of the equations.

 There are  natural and technical  reasons to focus on almost intersections. A good
deal of elimination theory is intertwined with birationality
questions. Now, a regular sequence of forms of fixed degree $\geq 2$
never defines a birational map. Thus, the first relevant case is the
next one, namely, that of an almost complete intersection.  Say,
$I\subset R=k[x_1,\ldots,x_d]$ is minimally generated by forms
$a_1,\ldots, a_d, a_{d+1}$ of fixed degree, where $a_1,\ldots, a_d$
form a regular sequence.  If these generators define a birational map
of ${\mathbb P}^{d-1}$ onto its image in ${\mathbb P}^d$ then
any set of forms of this degree containing
$a_1,\ldots, a_d, a_{d+1}$ still defines a birational map onto its
image.  Thus, almost complete intersections give us in some sense the
hard case.

 Note that these almost complete intersections have maximal
codimension, i.e., the ideal $I$ as above is $\m$-primary, where
$\m=(x_1,\ldots,x_d)$. The corresponding rational map with such base
ideal is a regular map with image a hypersurface of ${\mathbb
P}^{d}$. However, one can stretch the theory to one more case,
namely, that of an almost complete intersection of
$I=(a_1,\ldots,a_{d-1}, a_d)$ of submaximal height $d-1$.  Here the
corresponding rational map $\Psi_I:{\mathbb P}^{d-1}\dasharrow {\mathbb
P}^{d-1}$ is only defined off the support $V(I)$ (whose geometric
dimension is $0$) and one can ask when this map is birational---thus
corresponding to the notion of a Cremona transformation of ${\mathbb
P}^{d-1}$.

\begin{Notation}\rm To describe the problems treated and the solutions given, we give a
modicum of notation and terminology.
In order to avoid some degeneracies in the theory, we will assume throughout that $k$
has characteristic zero.
\begin{itemize}
\item[$\bullet$] $\RR=k[x_1, \ldots, x_d]_{\m}$ with $d \geq 2$ and $\m=(x_1,
\ldots, x_d)$.  

\item[$\bullet$] $I=(a_1, \ldots, a_d, a_{d+1})$ with
$\deg(a_i)=n$ for all $i$, $\height(I)=d$ and $I$ minimally generated
by these forms.  

\item[$\bullet$] Assume that $J=(a_1, \ldots, a_d)$ is a
minimal reduction of $I=(J, a)$, that is $I^{r+1}=JI^r$ for some
natural number $r$.

  \item[$\bullet$] $ \bigoplus_{j} \RR(-n_j)
\stackrel{\varphi}{\longrightarrow} \RR^{d+1}(-n) \longrightarrow I
\longrightarrow 0$ is a free minimal presentation of $I$.  

\item[$\bullet$] $\SS=\RR[\TT_1, \ldots, \TT_{d+1}]$ and $L=\ker(\SS\surjects \RR[It])$ via
$T_j\mapsto a_jt$, where $\RR[It]$ is the Rees algebra of $I$; note that
$L$ is a homogeneous ideal in the standard grading of $\SS$ with
$\SS_0=\RR$. 

 \item[$\bullet$] $L_i$ : $\RR$--module generated by forms of
$L$ of
degree $i$ in $\TT_j$'s. For example, the degree $1$ component of $L$
is the ideal of entries of the matrix product
\begin{eqnarray*}
L_1=I_1([\TT_1 \quad \cdots \quad
\TT_{d+1}]\cdot \varphi).\end{eqnarray*}

\item[$\bullet$] $\nu(L_i)$ denotes the minimal number of {\em fresh} generators
of $L_i$. Thus $\nu(L_2)$ is the minimal number of generators of the
$\RR$-module $L_2/S_1L_1$.

\item[$\bullet$]
The {\em elimination degree} of $I$ is
\begin{eqnarray} \label{edeg}
\edeg(I)&=&\inf
\{ i \;|\; L_i \nsubseteq \m \SS  \}.
\end{eqnarray}
One knows quite generally that $L= (L_1):
\mathfrak{m}^{\infty}$.
 A {\em secondary elimination
degree} is an integer $r$ such that $L= (L_1):
\mathfrak{m}^r$, i.e., an integer at least as large as the stabilizing exponent
of the saturation.

\item[$\bullet$]  The {\em special fiber} of $I$
is the ring $\mathcal{F}(I)= \RR[It]\otimes \RR/\mathfrak{m}$. This is
a hypersurface ring
\begin{eqnarray} \label{elimequ}
 \mathcal{F}(I) = k[T_1, \ldots, T_{d+1}]/(\ff(\TT)),\end{eqnarray}
 where
$\ff(\TT)$ is an irreducible polynomial of degree $\edeg(I)$.
$\ff(\TT)$ is called the {\em elimination equation}\index{elimination
equation} of $I$, and may
be taken as  an element of $L$.
\end{itemize}
\end{Notation}

Besides the syzygies of $I$, $\ff(\TT)$ may be considered the most
significant of the {\em equations} of $I$. Determining it, or at
least its degree, is one of the main goals of the subject.
The enablers, in our treatment, are   four Hilbert functions associated
to $I$, those of $\RR/I$, $\RR/J:a$, $\RR/I_1(\varphi)$ and the
Hilbert-Samuel function defined by $I$. Each encodes singly or in
conjunction  different
aspects of $L$.

\medskip

 In order to
outline  the other relationships between the invariants of the ideal $I$ and
its equations, we make use of the following diagram:

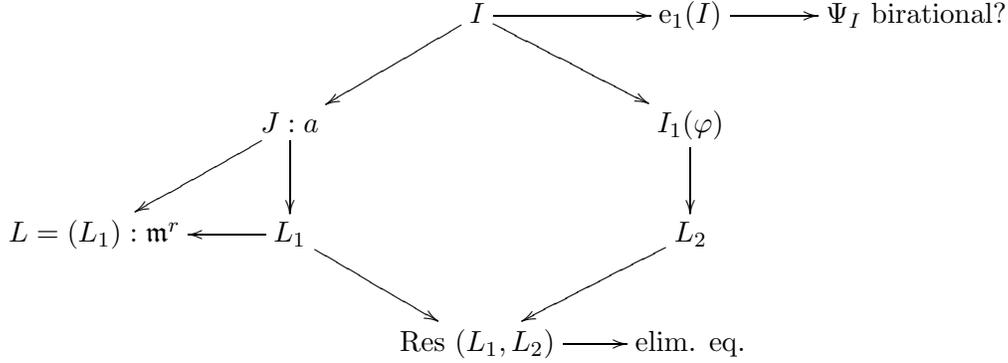
\begin{eqnarray*}
\diagram
& & I \dlto \rto \drto &\rme_1(I)\rto & \mbox{\rm $\Psi_I$ birational?} \\
& J:a \dlto \dto & & I_1(\varphi)\dto & \\
L=(L_1): \mathfrak{m}^r &  L_1\lto \drto & & L_2 \dlto  & \\
& & \mbox{\rm Res }(L_1,L_2)\rto & \mbox{\rm elim. eq.} &
\enddiagram
\end{eqnarray*}

At the outset and throughout there is the role
 played by
 the Hilbert function of $J:a$, which besides that of directing all the
 syzygies of $I$,  is the encoding  of an integer $r$
such that $L= (L_1):\mathfrak{m}^r$. According to
Theorem~\ref{secondedeg}, $r$ can be taken as $r=\epsilon +1$,
where $\epsilon$ is  the {\em socle degree} of
$\RR/J:a $. Since $L$ is expressed as the quotient of two
Cohen-Macaulay ideals, it has shown in practice to be an effective
tool for saturation.

A persistent  question is that of how to obtain higher degree
generators from the syzygies of $I$.
We will provide an iterative approach to generate the successive
components of $L$:
\[ L_1 \leadsto L_2 \leadsto L_3 \leadsto \cdots .\]
This is not an effective process, except for the step
$ L_1\leadsto L_2$.
 The more interesting development  seems to be the use
of the syzygetic lemmas to
obtain $\delta(I)=L_2/\SS_1L_1$ in the case of almost complete
intersections. This is a formulation of the {\em method of moving
quadrics} of several authors. The novelty here is the use of the
Hilbert function of the
ideal $I_1(\varphi)$ to understand and conveniently express
$\delta(I)$. Such  level of detail was not present even in earlier
treatments of $\delta(I)$.
The {\em syzygetic lemmas}   are observations based on the Hilbert
functions of $I$ and of $I_1(\varphi)$ to allow a description of
$L_2$. It converts the expression (\ref{deltaI})
\begin{eqnarray} \delta(I) = \Hom_{R/I}(\RR/I_1(\varphi),
H_1(I)),\end{eqnarray}
where $H_1(I)$ is the canonical module of $\RR/I$ (given by the
syzygies of $I$), into a set of generators of $L_2/\SS_1L_1$.
It is fairly effective in the case of binary, ternary and some
quaternary forms, as we
shall see. In these cases, out of $L_1$ and $L_2$ we will be able to
write the elimination equation in the form of a resultant
\begin{eqnarray}
\textrm{\rm Res }(L_1,L_2),\end{eqnarray} or as one of its factors.
 We  prove the
non-vanishing of this determinant under three different situations.
In the case of binary ideals, whose syzygies are of arbitrary
degrees, we give a far-reaching generalization of \cite{syl1}. For
an ideal $I$ of $k[x,y]$, generated by $3$ forms of degree $n$,
we build out of
$L_1$ and $L_2$, in a straightforward manner,  a nonzero
polynomial of $k[\TT_1, \TT_2, \TT_3]$ in $L$, of degree $n$ (Theorem~\ref{detB2plus}).

The other results, in higher dimension,
  require that the content of syzygies ideal $I_1(\varphi)$ be
a power of the maximal ideal $\mathfrak{m}$. Thus, in the case of
ternary forms of degree $n$ whose content ideal $I_1(\varphi)=
\mathfrak{m}^{n-1}$, our main result (Theorem~\ref{detB3s}) proves the
non-vanishing of $\mbox{\rm Res }(L_1,L_2)$ without appealing to
conditions of genericity (but with degree constraints). A final result
is  very special to quaternary forms (Theorem~\ref{detB42}).

\subsubsection{Preliminaries}

Before engaging in the above questions proper, we will outline the
basic homological and arithmetical data involved in these
ideals.

\subsubsection*{The resolution}

The general format of the resolution of $I$ goes as follows. First,
note that $J:a$ is a Gorenstein ideal, and that $I=J: (J:a)$. Since
Gorenstein ideals, at least in low dimensions, have an amenable
structure, it may be desirable to look at $J:a$
as a building block to $I$  and its equations. In this arrangement
the syzygies of $I$ will be organized in terms of those of $J$ and of
$J:a$.

\medskip

Let us recall how the resolution of $I$ arises as a mapping cone of
the Koszul complex $\mathbb{K}(J)$ and a minimal resolution of
$J:a$ (this was first given in \cite{pszpiro}; see also
\cite[Theorem A.139]{compu}).

\begin{Theorem} \label{Dubreil}
 Let $\RR$ be a Gorenstein local ring, let ${\mathfrak a}$ be a
perfect ideal of height $g$ and let $\mathbb{F}$ be a minimal free
resolution of $\RR/{\mathfrak a}$. Let ${\bf z}= z_1, \ldots, z_g \subset {\mathfrak a}$ be a
regular sequence, let
$\mathbb{K}=\mathbb{K}({\bf z}; R)$ be the corresponding Koszul complex,
 and let $u\colon  \mathbb{K} \rar \mathbb{F}$ be the
comparison mapping induced by the inclusion $({\bf z}) \subset L$.
Then the dual $\mathbb{C}(u^*)[-g]$ of the mapping cone of $u$, modulo the
subcomplex $u_0\colon R \rar R$, is a free resolution of length $g$ of
$\RR/({\bf z})\colon {\mathfrak a}$.
Moreover, the canonical module of
$\RR/({\bf z})\colon {\mathfrak a}$ (modulo shift in the graded case)
is ${\mathfrak a}/({\bf z})$.
\end{Theorem}

\begin{Remark} \label{Dubrmk} {\rm For later reference,
we point out
three observations when  the ideal
is the above almost complete intersection $I$.
\begin{enumerate}

\item[$\bullet$] If $J:a = (b_1, \ldots, b_d)$, writing
\[ [a_1, \ldots, a_d] = [b_1, \ldots, b_d]\cdot \phi, \]
gives
\[ I = (J, \det(\phi)),\]
so that $I$ is a Northcott ideal.

\item[$\bullet$]If $I$ is a generated by forms of degree $n$, a choice for $J$
is simply a set of forms $a_1, \ldots, a_d$ of degree $n$ generating
a regular sequence. Many of the features of
$I$--such as the ideal $I_1(\varphi)$--can be read off $J:a$.
Namely, $I_1(\varphi)$ is the sum of $J$ and  the coefficients arising
in the expressions of $a(J:a)\subset J$.

\item[$\bullet$] Since $J:a$ is a  Gorenstein ideal, its rich structure in
dimension $\leq 3$ is fundamental to the study in these cases.
%Discuss the Pfaffians (degrees, etc.)
\end{enumerate}

}\end{Remark}

\subsubsection*{Hilbert functions}

There are four Hilbert functions related to the ideal $I$ that are
relevant here, the first three of the Artinian modules
$\RR/I$, $\RR/I_1(\varphi)$ and $\RR/J:a$, and are  therefore of easy
manipulation. Their interactions will be a mainstay of the paper.

\medskip

 The first elementary observation, whose proof is well-known as to be omitted,
 will be useful when we need the Hilbert function of the
canonical module of $\RR/I$.

\begin{Proposition} \label{Dubrmk2}
If $I$
\ is generated by forms of degree $n$ of $k[x_1,
\ldots, x_d]$, the Hilbert function of $\RR/I$ satisfies
\begin{eqnarray*}  H_{R/I}(t) & = &  H_{R/J}(t) - H_{I/J}(t) \\
&=& H_{R/J}(t) - H_{R/J:a}(t-n).
\end{eqnarray*}
\end{Proposition}

The fourth Hilbert function is that of the associated graded ring
\[ \gr_I(\RR) = \bigoplus_{m\geq 0} I^m/I^{m+1}.\]
It affords the Hilbert--Samuel polynomial ($m\gg 0$)

\[ \lambda(\RR/I^{m+1}) = e_0(I){{d+m}\choose{d}} -\rme_1(I)
{{d+m-1}\choose{d-1}} + \textrm{lower degree terms of $m$},\]
where $\rme_0(I)$ is the
multiplicity of the ideal $I$. A related Hilbert polynomial is that
associated to the integral closure filtration $\{\overline{I^m}\}$:
 \[ \lambda(\RR/\overline{I^{m+1}}) =
 \overline{e}_0(I){{d+m}\choose{d}} - \overline{e}_1(I)
{{d+m-1}\choose{d-1}} + \textrm{lower degree terms of $m$}.\]
For an $\mathfrak{m}$-primary ideal $I$ generated by forms of degree $n$, $\overline{I^m}=
\mathfrak{m}^{nm}$, so the latter coefficients are really invariants of the ideal
$\mathfrak{m}^{n}$ and one has
\begin{eqnarray*}
e_0(I) &=& \overline{e}_0(I) = n^d \\
e_1(I) &\leq & \overline{e}_1(I) =   {\frac{d-1}{2}}(n^d-n^{d-1}).\\
\end{eqnarray*}

The case of equality $\rme_1(I)=\overline{e}_1(I)$ has a straightforward
(and general) interpretation in terms of the corresponding Rees
algebras.

\begin{Proposition} \label{r1} Let $(\RR, \mathfrak{m})$ be an
analytically unramified normal local domain of dimension $d$, and
 let $I$ be an
$\mathfrak{m}$-primary ideal. Then $\rme_1(I)=\overline{e}_1(I)$ if and
only if $\Rees(I)$ satisfies the condition $(R_1)$ of Serre.
\end{Proposition}

\begin{proof} Let $\AA= \Rees(I)$, set $\BB$ for its integral closure. $\BB$
is a finitely generated $\AA$-module. Consider the exact sequence
\[ 0 \rar \AA \lar \BB \lar C = \BB/\AA \rar 0. \]
Since $\rme_0(I)= \overline{e}_0(I)$, $C$ is a graded $\AA$-module of
dimension $\leq d$. If $\dim C=d$, its multiplicity $\rme_0(C)$ is derived from
the Hilbert polynomials above as $\rme_0(C) = \overline{e}_1(I)-e_1(I)$.
This sets up the assertion since $\AA$ and $\BB$ are equal in
codimension one if and only if $\dim C< d$. 
\end{proof}

This permits stating \cite[Proposition 3.3]{syl1} as follows (see also
\cite{ehrhart} for degrees formulas).

\begin{Proposition}
Let $\RR=k[x_1, \ldots, x_d]$ and let $I = (f_1,
\ldots, f_{d+1})$ be an ideal of finite colength, generated by forms
of degree $n$. Denote by $\mathcal{F}$ and $\mathcal{F}'$ the special
fibers of $\Rees(I)$ and $\Rees(\mathfrak{m}^n)$, respectively. The following
conditions are equivalent:
\begin{enumerate}
\item[{\rm (a)}] $[\mathcal{F}':\mathcal{F}]=1$, that is, the rational mapping
\[ \Psi_I: [f_1: f_2: \cdots : f_{d+1}] : \mathbb{P}^{d-1}\dasharrow
\mathbb{P}^d\]
is birational onto its image$\,${\rm ;}

\item[{\rm (b)}] $\deg \mathcal{F}= n^{d-1}\,${\rm ;}

\item[{\rm (c)}] $\rme_1(I) =  {\frac{d-1}{2}}(n^d-n^{d-1})\,${\rm ;}

\item[{\rm (d)}] $\Rees(I)$ is non-singular in codimension one.
\end{enumerate}
\end{Proposition}

A great deal of this investigation is to determine $\deg
\mathcal{F}$, which as we referred to earlier is the {\em elimination degree} of
$I$ (in notation, edeg$(I)$)\index{elimination degree}. Very often this turns into
explicit formulas for the elimination equation.

\subsection{Relation type}
We assume that $I$ is homogeneous, so that the ideal $LL$ is
bi-graded,
\[ \LL: \quad \ff= \sum_{\alpha}
{c_{\alpha}}{\TT^{\alpha}},
 \]
and we seek bounds on ${\deg(c_{\alpha})}$ and
${\deg(\TT^{\alpha})}$.
 We develop
 general bounds for each of  the two bi-degrees:
\[ \mbox{\rm reltype}(	I)\leq (r_{\xx}, r_{\TT}).\]
 $r_{\TT}$
 is still too large, being not specific.
We will discuss separately the two bounds $r_{\xx}$ and $r_{\TT}$.

\subsubsection*{Secondary elimination degrees}

A solution to some of questions raised above resides in the
understanding of the primary decomposition of $(L_1)$, the defining
ideal of the symmetric algebra ${\rm Sym}(I)$ of $I$ over $\SS=\RR[{\bf T}]$.
As pointed out earlier, we know that $\LL=(L_1):\mathfrak{m}^{\infty}$
as quite generally a power of $I$ kills $\LL/(L_1)$.
The $\LL$-primary component is therefore $\LL$ itself and $(L_1)$ has
only two primary components, the other being $\mathfrak{m}S$-primary.
On the other hand, according to
\cite[Proposition 2.2]{syl1}, $(L_1)$ is Cohen-Macaulay so its
primary components are of the same dimension.
Write
\[ (L_1)= \LL\cap Q,\]
where $Q$ is $\mathfrak{m}S$-primary. It allows for the expression of
$L$ as a saturation of $(L_1)$ in many ways. For example,
 for nonzero $\alpha \in I$ or $\alpha \in
I_1(\varphi)$, we have $\LL=(L_1) :\alpha^{\infty}$.

\medskip

We give now an explicit saturation by exhibiting integers $r$
 such that $\LL=
(L_1):\mathfrak{m}^r$, which as we referred to earlier are {\em secondary
elimination degrees}\index{secondary elimination degree}.
Finding its least value is one of our goals in
individual cases.
In the actual practice we have found the
computation effective, perhaps because $\LL$ is given as the quotient
of two Cohen-Macaulay ideals, the second generated by monomials.

\begin{Theorem}\label{secondedeg}
Let $\RR=k[x_1, \ldots, x_d]$ and let $I = (f_1,
\ldots, f_{d+1})$ be an ideal of finite colength, generated by forms
of degree $n$. Suppose $J=(f_1, \ldots, f_{d})$ is a minimal
reduction, and set $a=f_{d+1}$. Let $\epsilon$ be the socle degree of
$\RR/J:a$, that is the largest integer $m$ such that $(\RR/J:a)_m\neq 0$.
If $r=\epsilon+1$, then
\[ \LL= (L_1): \mathfrak{m}^r.\]
Moreover,  if some
form $\ff$ of $L_i$ has
coefficients in $\mathfrak{m}^r$, then $\ff\in (L_1)$. In particular,
the torsion elements of $\Sym(I)$ have $\RR$-degree $<r$.

\end{Theorem}

\begin{proof} The assumption on $r$ means that $(J:a)_m=({\mathfrak{m}^m})_m$ for
$m\geq r$, that is if $J:a$ has initial degree $d'$ then
\[ \mathfrak{m}^{r} = \sum_{i\geq d'}(J:a)_i \mathfrak{m}^{r-i}.\]

Now any element $\pp\in \LL$ can be written as
\[ \pp = \hh_p \TT_{d+1}^p + \hh_{p-1}\TT_{d+1}^{p-1} + \cdots + \hh_0, \]
where the $\hh_i$ are polynomials in $\TT_1, \ldots, \TT_{d}$.
If $\hh_1,\ldots, \hh_{p}$ all happen to vanish, then $\hh_0\in (L_1)$ since $f_1, \ldots,
f_{d}$ is a regular sequence.
The proof will consist in reducing to this situation.

To wit, let $u$ be a form of degree $r$ in $ \mathfrak{m}^r$. To prove that
$u\pp\in (L_1)$, we may assume that $u=v\alpha$, with $v\in
\mathfrak{m}^{r-i}$ and $\alpha$ a minimal generator in $(J:a)_i$.

  We are going to replace $v\alpha \pp$ by an
equivalent element of $\LL$, but of lower degree in $\TT_{d+1}$.
Since $\alpha\in (J:a)$, there is an induced form $\g2\in L_1$
\[ \g2 = \alpha \TT_{d+1} + \hh, \quad \mbox{\rm $\hh$ linear form in $\TT_1, \ldots, \TT_{d}$}.\]
Upon substituting $\alpha \TT_{d+1}$ by $\g2-\hh$, we get an equivalent
form
\[ \qq = \g2_{p-1} \TT_{d+1}^{p-1} + \g2_{p-2}\TT_{d+1}^{p-2} + \cdots + \g2_0, \]
where the coefficients of the $\g2_i$ all lie in $\mathfrak{m}^r$.
Further reduction of the individual terms of $\qq$ will eventually
lead to a form only in the $\TT_1, \ldots, \TT_{d}$.

The last assertion just reflects the nature of the proof.
\end{proof}

\begin{Remark}{\rm An a priori bound for the smallest secondary elimination
degree arises as follows. Let $\RR=k[x_1, \ldots, x_d]$,  and $I=
(J,a)$ as above. Since $J$ is generated by $d$ forms $f_1, \ldots, f_d$ of degree $n$,
the socle of $\RR/J$ is determined by the Jacobian of the $f_i$ which
has degree $d(n-1)$. Now from the exact sequence
\[0 \rar (J:a)/J \lar \RR/J \lar \RR/J:a\rar 0,\]
the socle degree of $\RR/J:a$ is smaller than $d(n-1)$, and therefore
this value gives the bound.  In fact, all the examples we examined
had  $\epsilon+1 = \min\{\; i \mid \LL=(L_1):\m^{i} \}$, where
$\epsilon$  is the socle degree of $\RR/(J:a)$.
}\end{Remark}

\begin{Example}{\rm Let  $\RR=k[x_1,x_2,x_3, x_4]$ and
 $\mathfrak{m}=(x_1, x_2, x_3, x_4)$. Let
\begin{enumerate}
%\item[$\bullet$] $\RR=k[x_1,x_2,x_3, x_4]$
%\item[$\bullet$] $\mathfrak{m}=(x_1, x_2, x_3, x_4)$.
\item[$\bullet$] $I=(x_1^3,\;\; x_2^3,\;\; x_3^3,\;\; x_4^3,\;\; x_1^2x_2+x_3^2x_4)$ with $a=x_1^2x_2+x_3^2x_4$.
\item[$\bullet$] $J=(x_1^3,\;\; x_2^3,\;\; x_3^3,\;\; x_4^3)$.
\item[$\bullet$] $J:a=(J,\; x_1x_3,\;\; x_1x_4^2,\;\; x_2^2x_3,\;\; x_1^2x_2-x_3^2x_4,\;\; x_2^2x_4^2)$.
\item[$\bullet$] The presentation of $I$ is of the form  $\RR(-5)\oplus R^{13}(-6)\oplus R(-7) \stackrel{\varphi}{\lar} R^5(-3) \stackrel{\tau}{\lar} I \lar 0$
and $I_1(\varphi)=(x_1x_3,\; x_1x_4,\; x_1^3,\; x_1x_2^2,\; x_1^2x_2,\; x_2x_3,\; x_2x_4^2,\; x_2^2x_4,\; x_2^3,\; x_3x_4^2,\; x_3^2x_4,\; x_3^3,\; x_4^3)$.
\item[$\bullet$] Hilbert series:
             \begin{enumerate}
             \item[$\diamond$] Hilbert series of $I$ \, : \, $1+ 4t+10t^2+15t^3+15t^4+7t^5+t^6$.
             \item[$\diamond$] Hilbert series of $(J:a)$ \, : \, $1+ 4t+9t^2+9t^3+4t^4+t^5$.
             \item[$\diamond$] Hilbert series of $I_1(\varphi)$ \, : \, $1+ 4t+7t^2$.
             \end{enumerate}
\end{enumerate}
A run with {\it Macaulay2} showed:
\begin{enumerate}
\item[{\rm (a)}] $L=(L_1):\mathfrak{m}^6$, exactly as predicted from the Hilbert
function of $J:a$;

\item[{\rm (b)}] The calculation yielded $\edeg(I)=9$; in particular, $\Psi_I$ is not
birational.
\end{enumerate}
}
\end{Example}

\subsection{The syzygetic lemmas}

The following material complements and refines some  known facts
(see  \cite{HSV1}, \cite{SV1}, \cite[Chapter 2]{alt}).
Its main purpose is an application to almost complete
intersections.
Since its contents deal with arbitrary ideals, we will momentarily change notation.
Let $I\subset \RR$ be an ideal generated by $n$ elements $a_1,\ldots, a_n$.
Consider a free presentation of $I$
\begin{equation}\label{presentation}
\RR^m\stackrel{\varphi}{\lar} \RR^n\lar I\rar 0
\end{equation}
and let $(L_1)\subset \LL=\bigoplus_{d\geq 0} L_d \subset \BB=\RR[\TT]=\RR[\TT_1, \ldots,\TT_n]$ denote as
before the presentation ideals of the symmetric algebra
and of the Rees algebra of $I$, respectively, corresponding to the chosen presentation.

\medskip

A starting point is the following
  observation. Suppose $\ff(\TT)= \ff(\TT_1, \ldots,
\TT_n)\in L_d$; write it as
\[ \ff(\TT_1, \ldots, \TT_n) = \ff_1(\TT) \TT_1 + \cdots +\ff_n(\TT)
\TT_n,\]
where $\ff_i(\TT)$ is a form of degree $d-1$.

\medskip

Evaluating $\TT$ at the vector $\aa= (a_1, \ldots, a_n)$ gives a
syzygy of $\aa$ \[ z=(\ff_1(\aa), \ldots, \ff_n(\aa)) \in Z_1,\] the
module of syzygies of $I$, \[ z\in Z_1 \cap I^{d-1}R^n,\] that is,
$z$ is a syzygy with coefficients in $I^{d-1}$.  Note that \[
\widehat{\ff}(\TT)=a_1\ff_1(\TT) + \cdots + a_n \ff_n(\TT) \in
L_{d-1}\cap I\cdot \BB_{d-1}.\] Conversely, any form ${\hh}(\TT)$ in
$L_{d-1}\cap I\cdot \BB_{d-1}$ can be lifted to a form $\HH(\TT)$ in
$L_d $ with $\widehat{\HH}(\TT)=\hh(\TT)$.

\medskip

 Such  maps are referred to \index{downgrading and upgrading maps}
as {\em downgrading} and {\em
upgrading}, although they are not always well-defined. However, in some case
it opens the opportunity to calculate some of the higher $L_d$.
\medskip

Here is an useful observation.

\begin{Lemma} \label{syzlemma1} Let
 $\ff(\TT)\in L_d$ and write
\[  \ff(\TT) = \ff_1(\TT) \TT_1 + \cdots + \ff_n(\TT) \TT_n.\]
If
 \[ \ff_1(\TT) a_1 + \cdots +\ff_n(\TT) a_n=0,\]
then $\ff(\TT) \in \BB_{d-1}L_1$.
\end{Lemma}

\begin{proof} The assumption is that $\vv=(\ff_1(\TT), \ldots, \ff_n(\TT))$
is a syzygy of $a_1, \ldots, a_n$ over the ring $\RR[\TT]$. By
flatness, \[ \vv = \sum_{j} \hh_j(\TT) \zz_j,\] where $\hh_j(\TT)$
are forms of degree $d-1$ and $\zz_j\in Z_1$. Setting $\zz_j =
(z_{ij}, \ldots, z_{nj})$, \[ \ff_i(\TT) = \sum_{j} \hh_j(\TT)
z_{ij},\] and therefore
\begin{eqnarray*}
  \ff(\TT) &=& \ff_1(\TT) \TT_1 + \cdots +\ff_n(\TT) \TT_n\\
&=& \sum_i (\sum{_j} \hh_j(\TT) z_{ij})\TT_i \\
&=& \sum_{j}\hh_j(\TT) (\sum_{i} z_{ij}\TT_i) \in \BB_{d-1}L_1.\\
\end{eqnarray*}
\end{proof}

\begin{Corollary} Let
$\hh_j(\TT)$, $1\leq j\leq m$, be a  set of generators
of  $L_{d-1}\cap
I\BB_{d-1}$. For each $j$, choose a form $\FF_j(\TT)\in L_d$ such that
under one of the operations above $\widehat{\FF}(\TT)=\hh_j(\TT)$. Then
\[ L_d = (\FF_1(\TT), \ldots, \FF_{m}(\TT), L_1\BB_{d-1}).\]
\end{Corollary}

\begin{proof}  For $\ff(\TT)\in L_d$, $\ff(\TT)=
\sum_{i}T_i\ff_i(\TT)$, write
\[  \sum_{i}a_i\ff_i(\TT) = \sum_j c_{j}\hh_j(\TT).\]
Applying Lemma~\ref{syzlemma1} to the   polynomial
\[ \ff(\TT) - \sum_{j}c_j\FF_j(\TT)\]
will give the desired assertion.
\end{proof}

This leads to the iterative procedure to find the equations  $\LL=(L_1, L_2,
\ldots)$ of the ideal $I$.

\bigskip

Let $Z_1=\ker(\varphi)\subset \RR^n$, where $\varphi$ is as in (\ref{presentation})
and let $B_1$ denote the submodule of $Z_1$ whose elements come from the Koszul relations
of the given set of generators of $I$.
The $\RR$-module
\[ \delta(I)=  Z_1\cap I\RR^n/B_1 \]
has been introduced in \cite{Si} in order to understand the Koszul homology with
coefficients in $I$. It is independent of the free presentation of $I$
and as such it has been dubbed the {\em syzygetic module} of $I$.
 The following basic result has been proved in \cite{SV1}.

\begin{Lemma}[The main syzygetic lemma] \label{syzlemma} Let $I$ be an ideal with
presentation as above.  Then
\begin{equation}\label{main_syzygetic}
\delta(I)\, \stackrel{\phi}{\simeq} \, {L}_2/ {L}_1\BB_1.\\
\end{equation}
\end{Lemma}
The mapping $\phi$ is given as follows: For $z= \sum \alpha_i\TT_i$,
$(\alpha_1, \ldots, \alpha_n)\in Z_1\cap IR^n$,
$\alpha_i=\sum_{j=1}^n c_{ij}a_j$,
\[ \phi([z]) = \sum_{ij}c_{ij}\TT_i\TT_j \in
 {L}_2/\BB_1 {L}_1.\]
 In particular, ${L}_2/\BB_1L_1$ is also independent of the free presentation of $I$ and
\[ \nu( {L}_2/\BB_1 {L}_1)= \nu(\delta(I)).\]

We refer to the process of writing the $\alpha_i$ as linear
combination of the $a_j$ as the {\em extraction} of $I$. The
knowledge of the degrees of the $c_{ij}$ is controlled by the
degrees of $\delta(I)$. Note that, quite generally, the kernel of the natural surjection
$\Symi(I) \lar \Rees(I)$
in degree $d$ is $L_d/L_1\BB_{d-1}$.
However, a more iterative form of bookkeeping of $\LL$ is through the modules $L_d/\BB_1L_{d-1}$
that represent the {\em fresh} generators in degree
$d$. Unfortunately, except for the case $d=2$, one knows no explicit expressions for
these modules, hence the urge for different methods to approach the problem.

\bigskip

\subsubsection*{Almost complete intersections}
We now go back to the particular setup of almost complete intersections.
As before, $(\RR, \mathfrak{m})$ denotes the standard graded polynomial ring
$k[x_1,\ldots, x_d]$ and its irrelevant ideal and $I=(a_1, \ldots,
a_d, a_{d+1})$ is an $\mathfrak{m}$-primary ideal minimally generated
by $d+1$ forms of the same degree. We assume that $J=(a_1, \ldots, a_d)$ is a minimal
reduction of $I$; set $a=a_{d+1}$.

\medskip

Considerable  numerical information about $ {L}_2$ is readily available
in this setup.

\begin{Proposition} \label{syzlemma0} Let $I$ be as above  and
let $\varphi$ be a minimal presentation map as in {\rm (\ref{presentation})}. Then
\[ \nu( {L}_2/\BB_1 {L}_1)= \nu (I_1(\varphi):
\mathfrak{m})/I_1(\varphi)).\] In particular, if
$I_1(\varphi)=\mathfrak{m}^{s}$, $s\geq 1$, then
\[ \nu( {L}_2/\BB_1 {L}_1)= {{d+s-2}\choose{d-1}}.\]
Moreover, $\delta(I)$ is generated by the last $s$ graded components of the
first Koszul homology module $H_1(I)$.
\end{Proposition}

\begin{proof} Consider the so-called syzygetic sequence of $I$
\[ 0 \rar \delta(I)\lar \H_1(I) \lar (\RR/I)^{d+1} \lar I/I^2 \rar 0. \]
Note that  $\H_1(I)$ is isomorphic to the canonical module of $\RR/I$. Dualizing with
respect to $\H_1(I)$ gives the exact complex
\[ \H_1(I)^{d+1} \lar \Hom_{\RR/I}(\H_1(I),\H_1(I))\simeq  \RR/I \lar
\Hom_{\RR/I}(\delta(I), \H_1(I)\rar 0.\] The image in $\RR/I$ is the ideal
generated by $I_1(\varphi)$, and since $I\subset I_1(\varphi)$,  one has
\begin{eqnarray} \label{deltaI}
 \delta(I) \simeq \Hom_{\RR/I}(\RR/I_1(\varphi), \H_1(I)).
\end{eqnarray}
It follows that $\delta(I)$ is isomorphic to the canonical module of
$\RR/I_1(\varphi)$,
and therefore $\nu(\delta(I))$ is the
Cohen-Macaulay type
of $\RR/I_1(\varphi)$.

In case $I_1(\varphi)=\mathfrak{m}^{s}$,
$\Hom_{R/I}(\RR/\mathfrak{m}^s, \H_1(I))$ cannot have a nonzero element
$u$ in  higher degree as otherwise $\mathfrak{m}^{s-1}u$ would lie in
the socle of $\H_1(I)$, a contradiction.
 It follows that $\delta(I)$  is generated
by ${d+s-2}\choose{d-1}$ elements.
\end{proof}

To help identify the generators of $\delta(I)$
requires information about the Hilbert function of $\H_1(I)$.
For reference  we use
 the socle degree of
$\RR/I_1(\varphi)$, which we denote by $p$.
 We
recall that if $(1, d, a_2, \ldots, a_r)$ is the Hilbert function of
$\RR/I$, that of $\H_1(I)$ is $(a_r, \ldots, d, 1)$, together with an
appropriate
shift. Since $\delta(I)$ is a graded submodule of $\H_1(I)$, it is convenient to
organize a table as follows:
\begin{eqnarray*}
\H_1(I) & : & (a_r, \ldots, a_s, \ldots, a_2, d, 1) \\
\delta(I) &: & (b_p, \ldots, b_1,1),
\end{eqnarray*}
where $b_i\leq a_i$. Note that the degrees are increasing. For
example, the degree of the $r$th component of $\H_1(I)$ is the initial
degree $t$ of the ideal $J:a$, while the degree of the $p$th component of
$\delta(I)$ is $t+r-p$.

\bigskip

We can also formulate the previous result in the following form.

\begin{Proposition} \label{acinotsyzygetic} Let $\RR$ be a Gorenstein local ring 
and $I$ an almost complete intersection of dimension $0$. Then  $\delta(I)\neq 0$, that is $I$ is 
not syzygetic. 
\end{Proposition}

\begin{proof}
This comes from the formula $\delta(I) = \Hom_{\RR/I)}(\RR/I_1(\varphi), \H_1(\RR/I))$. 
\end{proof}

%\bigskip

\subsubsection*{Balanced ideals}

Let us introduce the following concept for ease of  reference:

\begin{Definition}\rm
An $\mathfrak{m}$-primary ideal $I\subset \RR$ minimally generated
by $d+1$ elements of the same degree is $s$-{\em balanced} if
 $I_1(\varphi)=\m^s$, where $\varphi$ is the matrix of syzygies of
 $I$.
\end{Definition}

Note  that, due to the Koszul relations,  $s$ is at most the common
degrees of the generators of $I$
The basic result about these ideals goes as follows.

\begin{Theorem} \label{socdegree} Let $\RR=k[x_1, \ldots, x_d]$ and
 $I$ an almost complete  intersection of  finite colength
 generated by forms of degree $n$. If $I$ is $s$--balanced,
  then:
 \begin{enumerate}
 \item[{\rm (a)}]  The socle degree of $\H_1(I)$ is $d(n-1)$;
 
 \item[{\rm (b)}]  $\mathfrak{m}^{d(n-1)-s+1} = I
 \mathfrak{m}^{(d-1)(n-1)-s}$;
 
 \item[{\rm (c)}] The coefficients of $L_2$ have degree
 $r(I)=(d-1)(n-1)-s$;
 
\item[{\rm (d)}] $\RR_{n+r(I)}= I_{n+r(I)}$.

\end{enumerate}
\end{Theorem}
\begin{proof}
(i)  Let $J$ be the minimal reduction of $I$ defined earlier.
The socle degree of $H_1(I)$ is determined from the natural
embedding
\begin{eqnarray*}
\H_1(I) \simeq J:a/J \hookrightarrow \RR/J,
\end{eqnarray*} where the socle of $\RR/J$, which is also the socle of
any of its nonzero submodules,
 is defined by the Jacobian
determinant
of $d$ forms of degree $n$.

\medskip

(ii) We write $H_1(I)$ and $\delta(I)$ as graded modules (set
$\epsilon=d(n-1)$)
\begin{eqnarray*}
\H_1(I) & = &  h_s \oplus h_{s+1} \oplus \cdots \oplus h_{\epsilon} \\
\delta(I) &= & f_{\epsilon -s+1} \oplus \cdots \oplus f_{\epsilon},
\end{eqnarray*} dictated by the fact that the two modules share the
same socle, $h_{\epsilon}=f_{\epsilon}$. One has
$\mathfrak{m}^{\epsilon-s+1} \H_1(I) = 0$, hence
$ \mathfrak{m}^{\epsilon-s+1} \RR/I = 0$,
or equivalently, \[\mathfrak{m}^{\epsilon-s+1}=
I\mathfrak{m}^{\epsilon-n-s+1}.\]

(iii) The degree $r(I)$ of the coefficients of $L_2$ is obtained from the
elements of $f_{\epsilon -s+1}$, and writing them as syzygies with
coefficients in I, that is
\[ r(I)= \epsilon-s+1-n = (d-1)(n-1)-s.\]

(iv) The last assertion follows from (ii) and (iii).
\end{proof}

Let us give some consequences of this analysis which will be used later.

\begin{Corollary} \label{wellbal} Let $\RR=k[x_1, \ldots, x_d]$ be a ring of
polynomials, and $I$ an almost complete intersection of finite
colength generated in degree $n$. Suppose that $I$ is $s$--balanced.
\begin{enumerate}
\item[{\rm (a)}] Let $d=3$.
If  $s=n-1$, $L_2$ is generated by ${{s+1}\choose{2}}$ forms with
coefficients of degree $s$ and
there are precisely $n$  minimal syzygies of
degree $n-1$, and $n\leq 7\,${\rm ;}

\item[{\rm (b)}] If $d=4$ and $s=n=2$, then there are precisely $15$ minimal
syzygies of degree $2$.
\end{enumerate}
\end{Corollary}

\begin{proof}  We begin by observing the values of $r(I)$. In
 case (a), $r(I)= (3-1)(n-1)-s= (3-2)(n-1)=s$, while in (b)
$r(I)= (4-1)(2-1)-2=1$.

\medskip

The first assertion of (a) comes from Proposition~\ref{syzlemma0} and
the value $r(I)=s$. As for the number of syzygies,
 according to Theorem~\ref{socdegree}(d),
$\RR_{n+s}= I_{n+s}$ in case (a), and $\RR_{n+1}=I_{n+1}$ in case (b),
which
will permit the determination of the
dimension of the syzygies of degree $s$, or higher in case (a), and
for all degrees in case (b).

\medskip

Let us focus on the case $r(I)=s$.
  Consider the exact sequence corresponding to the generators of
$I$,
\[ \RR^{d+1} \stackrel{\pi}{\lar} \RR \lar \RR/I\rar 0,\]
and write $\psi_s$ for the
vector space map
induced by $\pi$ on the homogeneous component
of degree $s$ of $\RR^{d+1}$. We have an exact sequence of $k$-vector spaces and $k$-linear
maps
\[ \RR_s^{d+1} \stackrel{\psi_s}{\lar} \RR_{n+s} \lar \RR_{n+s}/I_{n+s}\rar 0\]
and $\ker (\psi_s)$ is the $k$-span of the syzygies of $I$ of $\RR$-degree $s$.

One easily has
\begin{equation}\label{dimker} \dim_k (\ker (\psi_s))=
(d+1){{s+d-1}\choose{d-1}} - \dim_k(I_{n+s}).
\end{equation}

In this case, one gets
\begin{equation}\label{dimker2} \dim_k (\ker (\psi_s))=
(d+1){{s+d-1}\choose{d-1}} - {{s+n+d-1}\choose{d-1}}.
\end{equation}
If $I$ is moreover $s$-balanced for some $s\geq 1$ then it must be the case that
\begin{equation}\label{balanced_inequality} (d+1)\dim_k (\ker (\psi_s))\geq
\nu (\mathfrak{m}^s)=\dim_k(\RR_s)=
{{s+d-1}\choose{d-1}}.
\end{equation}

(1) For $d=3$ and $I$ $s$-balanced with $s=n-1$ the equality (\ref{dimker2}) gives
$ \dim_k (\ker (\psi_{n-1}))=n$ while the inequality (\ref{balanced_inequality})
easily yields $n\leq 7$.

\medskip

(2) This follows immediately from (\ref{dimker2}).
\end{proof}

The numerical data alone give a bird eye view of the generators of the
graded pieces $L_1$ and $L_2$ of the ideal of equations of $L$.
This corollary is suitable in other cases, even when $I_1(\varphi)$ is a
less well packaged ideal.

\subsection{The number of equations}
Let $\RR=k[x_1, \ldots, x_{d}]$ be a ring of polynomials, and
 $I$ be a homogeneous  almost complete intersection $I=(a_1,
\ldots, a_d, a_{d+1})$ of finite colength.
We shall now calculate the role of the secondary elimination degree
of the ideal $I$ in estimations of the number of generators of the
ideal $L$.

\medskip

%\begin{Theorem} For every natural number $r$,
%\[ \nu((L_1):\mathfrak{m}^r) \leq \beta_1(I) + 2 {{r+d-1}\choose {d}} ,\]
%where $\beta_1(I)$ is the first Betti number of the ideal $I$.
% \end{Theorem}

We start our discussion
by first giving an application of the $\mathcal{Z}$--complex to the
symmetric algebra $\Sym(I)$.

\begin{proposition} $\Sym(I)$ is a Cohen-Macaulay ring of type at
most $\dim \RR -1$.
\end{proposition}

\begin{proof} The first assertion is a direct consequence of
{\cite[Theorem 10.1]{HSV1}}:

\begin{Corollary}\label{Zaci3} Let $\RR$ be a Cohen-Macaulay local and let $I$ be an
almost complete intersection. If $\depth \RR/I\geq \dim \RR/I-1$, then
the complex $\mathcal{Z}$ is acyclic and $\Sym(I)$ is a
Cohen-Macaulay algebra.
\end{Corollary}

The determination of the Cohen-Macaulay type comes about as follows.
Set $\BB=\RR[\TT_1, \ldots, \TT_{d+1}] $, and let $\KK$ be the Koszul
complex of $\aa=\{a_1, \ldots, a_{d+1}\}$.
Denoting by $Z_i$ the $i$--cycles of $\KK$, the
$\mathcal{Z}$--complex
is
\[ 0 \rar Z_d \otimes \BB[-d] \rar \cdots \rar Z_1\otimes \BB[-1] \rar
(L_1)
\rar 0,\]
\[ 0 \rar (L_1) \rar \BB \rar \Sym(I) \rar 0.\]
The canonical module of $\Sym(I)$ is (up to a shift) the module
\[ \Ext_{\BB}^{d}(\Sym(I), \BB)\simeq \Ext_{\BB}^{d-1}((L_1), \BB).\]

After observing that $Z_d=\RR$,
we break the long exact sequence into the short exact sequences
\[ 0\rar \SS \rar Z_{d-1}\otimes \BB[-(d-1)]\rar  \H_{d-2}\rar 0,
\]
\[ 0 \rar \H_{d-2}\rar Z_{d-2}\otimes \BB[-(d-2)]\rar  \H_{d-3}\rar 0,
\]
\[ \vdots \]
\[ 0 \rar H_{2}\rar Z_{2}\otimes \BB[-2]\rar  H_1\rar 0,
\]
\[ 0 \rar H_{1}\rar Z_{1}\otimes \BB[-1]\rar  (L_1)\rar 0.
\]

Noting that for each $i\geq 2$, the strands of the Koszul complex
\[ 0\rar K_{d+1} \rar K_d \rar \cdots \rar K_{i+1}\rar Z_i \rar 0
\]
are acyclic, while
\[ 0\rar K_{d+1} \rar K_d \rar \cdots \rar K_{2}\rar B_1 \rar 0
\] is a resolution of the module of $1$-boundaries, and
\[0 \rar B_1 \rar Z_1 \rar \H_1(I) \rar 0,\] where $H_1(I)$ is the
canonical module of $\RR/I$,
we are now in position to examine $\Ext_{\BB}^{d-1}((L_1),
\BB)$.

\bigskip

\begin{lemma} $\depth \H_i \geq i+1$ for $1\leq d-2$.
\end{lemma}

\begin{proof} The assertion is true for $i=d-2$ since $Z_d\simeq \RR$ and
$\depth Z_{i}=i$. By induction and the depth lemma applied the
sequence
\[ 0 \rar \H_{i+1}\rar Z_{i+1}\otimes \BB[-(i+1)]\rar  \H_i\rar 0,\]
we prove the assertion.
\end{proof}

This gives rise to the following exact sequences:

\[
\Ext_\BB^{d-2}(\H_1, \BB)\rar \Ext_\BB^{d-1}((L_1), \BB ) \rar
\Ext_\BB^{d-1}(Z_1\otimes \BB[-1], \BB)=0,
\]
\[
\Ext_\BB^{d-3}(\H_2, \BB)\rar \Ext_\BB^{d-2}(\H_1, \SS) \rar
\Ext_\BB^{d-2}(Z_2\otimes \BB, \BB)= \RR/I\otimes \BB[+2] \rar \Ext_\BB^{d-2}(\H_2,
\BB)=0\]
\[ \vdots \]
\[ \Hom_\BB(\BB,\BB)\rar \Ext_\BB^{1}(\H_{d-2}, \BB) \rar
\Ext_\BB^{1}(Z_{d-1}\otimes \BB[-(d-1)],
\BB)=\RR/I\otimes \BB[d-1] \rar 0.
\]
Stitching together the maps
\[\Ext_\BB^{1}(\H_{d-2}, \BB)\rar
\cdots \rar
 \Ext_\BB^{d-3}(\H_2, \BB)\rar
 \Ext_\BB^{d-2}(\H_1, \BB) \rar
 \Ext_\BB^{d-1}((L_1), \BB),
\]
it gives that $\Ext_\BB^{d-1}((L_1), \BB)$ is generated by $d-1$
elements. \end{proof}

\subsubsection*{Binary ideals}

Let $\RR=k[x,y]$ (instead of the
general notation $\RR=k[x_1,x_2]$).
Let $I\subset \RR=k[x, y]$ be an $(x,y)$-primary ideal generated by three forms of
degree $n$.  $I$ has a minimal free resolution \[ 0 \rar \RR^{2}
\stackrel{\varphi}{\lar} \RR^3 \lar I \rar 0.\]
We will assume throughout the section that the first   column of
$\varphi$ has  degree $r$, the other degree $s\geq r$. We note
$n=r+s$.

\medskip

Here we give a general format of the elimination equation of $I$ up to a power,
thus answering several questions raised in \cite{syl1}.

\subsubsection*{Elimination equation and degree}

Set $\BB=\RR[\TT_1, \TT_2, \TT_3]$ as before.
Notice that $\BB$ is standard bigraded over $k$.
We denote by $f$ and $g$ the defining forms of the symmetric algebra
of $I$, i.e., the generators of the ideal $(L_1)\subset \BB$ in the earlier
notation. We write this in the form
\[ [f,g] = [\TT_1, \TT_2,\TT_3]\cdot  \varphi.\]
In the standard bigrading, by assumption, $f$ has bidegree
$(\RR,1)$, $g$ bidegree $(s,1)$.
According to Lemma~\ref{syzlemma}, the component $L_2$ could be determined from
$(L_1):I_1(\varphi)$. In dimension two it is more convenient to get
hold of a smaller quotient, $N=(L_1):(x,y)^{r}$. We apply basic
linkage theory to develop some properties of this ideal.

\begin{enumerate}
\item[$\bullet$]
 $N$, being a
direct link of the Cohen-Macaulay  ideal $(x,y)^{r}$, is a perfect
Cohen-Macaulay ideal of codimension two.
The canonical module of $\BB/N$ is generated by $(x,y)^{r}\BB/(f,g)$, so
that its Cohen-Macaulay type is $r+1$, according to
Theorem~\ref{Dubreil}.

\item[$\bullet$] Therefore, by
the Hilbert-Burch theorem, $N$ is the ideal of maximal minors of an
$(r+2)\times (r+1)$ matrix $\zeta$ of homogeneous forms.

\item[$\bullet$] Thus,  $N=(f,g) : (x,y)^{r}$ has a
presentation $ 0 \rar \BB^{r+1} \stackrel{\zeta}{\lar} \BB^{r+2} \lar N
\rar 0$, where $\zeta$ can be written in the form \[ \zeta = \left[
\begin{array}{c} \sigma \\ ------- \\ \tau \end{array}\right], \]
with $\sigma$ is a $2 \times (r+1)$ submatrix with rows  whose entries are
biforms of bidegree $(0,1)$ and $(s-r, 1)$; and $\tau$
is an $r \times (r+1)$ submatrix
whose entries are  biforms of bidegree $(1,0)$.
\item[$\bullet$] Since $N\subset L_2$, this shows that in $L_2$ there are
$r$ forms  $\hh_i$ of degree
$2$ in the $\TT_i$ whose $\RR$-coefficients are forms in $(x,y)^{s-1}$.

\item[$\bullet$]
If $s=r$, write
\begin{equation}\label{biforms}
[\hh]= [\hh_1 \;\; \cdots \;\;
\hh_{r}] = [x^{r-1}\;\; x^{r-2}y \;\; \cdots \;\; xy^{r-2} \;\; y^{r-1}]
\cdot \CC,
\end{equation}
 where $\CC$ is an $r \times r$ matrix whose entries belong
to $k[\TT_1, \TT_2, \TT_3].$

\item[$\bullet$]
If $s> r$, collecting the $s-r$ forms of the product
$f(x,y)^{s-r-1}$, write
\begin{equation}\label{biforms2}
[\ff;\hh]= [f(x,y)^{s-r-1}; \;\;\hh_1 \;\; \cdots \;\;
\hh_{r}] = [x^{s-1}\;\; x^{s-2}y \;\; \cdots \;\; xy^{s-2} \;\; y^{s-1}]
\cdot \CC,
\end{equation}
 where $\CC$ is an $s \times s$ matrix whose entries belong
to $k[\TT_1, \TT_2, \TT_3]$.

\end{enumerate}

A first consequence of this analysis is one of our main results:

\begin{Theorem} \label{detC} In both cases,
 $\det \CC $ is a polynomial of degree $n$.
\end{Theorem}

\begin{proof}
Case $s=r$:
 Suppose that $\det \CC=0$. Then there exists a nonzero
vector $\left[\begin{array}{c} {\bf a}_1 \\ \vdots \\ {\bf a}_{r}
\end{array}\right]$  whose entries are in $k[\TT_1, \TT_2, \TT_3]$ such
that  $\CC \cdot \left[\begin{array}{c} {\bf a}_1 \\ \vdots \\ {\bf
a}_{r}
\end{array}\right]=0$.  Hence $[\hh_1 \;\; \cdots \;\; \hh_{r}]
\left[\begin{array}{c} {\bf a}_1 \\ \vdots \\ {\bf a}_{r}  \end{array}\right]=0$.
Since the relations of the $0,0, \hh_1,\ldots, \hh_{r}$ are $\BB$-linear combinations of
the columns of $\zeta$, we get a
contradiction.
\medskip

The assertion on the degree follows since the degree of $\det \CC$ is
$2r=r+s=n$.

\medskip

Case $s>r$:
 Suppose that $\det \CC=0$. Then there exists a nonzero
vector $\left[\begin{array}{c} {\bf a}_1 \\ \vdots \\ {\bf a}_s
\end{array}\right]$  whose entries are in $k[\TT_1, \TT_2, \TT_3]$ such
that  $\CC \cdot \left[\begin{array}{c} {\bf a}_1 \\ \vdots \\ {\bf
a}_s
\end{array}\right]=0$.  Hence $[\ff; \; \; \hh_1 \;\; \cdots \;\;
\hh_{r}]
\left[\begin{array}{c} {\bf a}_1 \\ \vdots \\ {\bf a}_s
\end{array}\right]=0$.
We write this relation as follows
\[
\sum_{i=1}^{s-r} {\bf a_i} f x^{s-r-i} y^{i-1} + \sum_{j=1}^{r} {\bf
a}_{s-r+j}\hh_j=
{\bf a}f
 + \sum_{j=1}^{r} {\bf
a}_{s-r+j}\hh_j=0,  \]
where \[ {\bf a}=
\sum_{i=1}^{s-r} {\bf a_i}  x^{s-r-i} y^{i-1}.\]

Since the relations of the $f,0, \hh_1,\ldots, \hh_{r}$ are $\BB$-linear combinations of
the columns of $\zeta$,
\begin{itemize}
\item[$\bullet$] ${\bf a}_{s-r+j}\in (x,y)\BB$, for $1\leq j\leq r $; and
\item[$\bullet$] ${\bf a}\in (x,y)^{s-r}\BB$,
\end{itemize}
and therefore ${\bf a}_i\in (x,y)\BB$, for all $i$. This gives a
contradiction.
\medskip

The assertion on the degree follows since the degree of $\det \CC$ is
$(s-r) + 2r =r+s=n$.
\end{proof}

\begin{Example}{\rm Let $\RR=k[x, y]$ and $I$ the ideal defined by
 {$\varphi=\left[ \begin{array}{ll}
x^2 & y^4 \\ xy & x^3y+x^4 \\ y^2 & xy^3
\end{array} \right]$.}
\begin{enumerate}

\item[$\bullet$] $N= (f,g) : \m^2 = (f,g, h_1, h_2)$, where
\[\begin{array}{lll}
f &= & x^2\TT_1+xy\TT_2+y^2\TT_3 \\
g &= & y^4\TT_1+(x^3y+x^4)\TT_2 + xy^3\TT_3 \\
\hh_1 &= & y^3\TT_1^2-x^3\TT_2^2-x^2y\TT_2^2+xy^2\TT_1\TT_3-x^2y\TT_2\TT_3-xy^2\TT_2\TT_3 \\
\hh_2 &= & xy^2\TT_1^2+y^3\TT_1\TT_2-x^3\TT_2\TT_3-x^2y\TT_2\TT_3-y^3\TT_3^2
\end{array}
\]

\item[$\bullet$] $[xf\;\; yf \;\; h_1 \;\; h_2 ] = \m^3 \CC$, where
\[
\CC=\left[
\begin{array}{llll}
\TT_1 & 0 & -\TT_2^2 & -\TT_2\TT_3 \\
\TT_2 & \TT_1 & -\TT_2^2-\TT_2\TT_3 & -\TT_2\TT_3 \\
\TT_3 & \TT_2 & \TT_1\TT_3 -\TT_2\TT_3 & \TT_1^2 \\
0 & \TT_3 & \TT_1^2 & \TT_1\TT_2-\TT_3^2
\end{array}
\right]
\]

\item[$\bullet$]  $\det\CC$ is the elimination equation.
\end{enumerate}
}\end{Example}

\begin{Remark}{\rm The polynomials $\hh_1, \ldots, \hh_{r}$ were also obtained
in \cite{syl1} by a direct process involving Sylvester elimination, in
the cases $s=r$ or $s=r+1$. In
\cite{syl1} though they did not arrive with the elements of
structure--that is with their relations--provided in the
Hilbert-Burch matrix. It is this fact that opens the way in the
binary case to a greater generality to the ideals treated and a more
detailed understanding of the ideal $\LL$.
}\end{Remark}

The next result provides a secondary elimination degree for these
ideals.

\begin{Corollary}\label{balanced_exponent} $\LL= (L_1): (x,y)^{n-1}$.
\end{Corollary}

\begin{proof}
With the previous notation, let $N=\LL\cap Q$ be the primary decomposition of $N$,
where $Q$ is $(x,y)\BB$-primary. Writing $\beta=\det \CC$, we then have $N:\beta=Q$.
On the other hand, $(L_1,{\bf h})\subset N$ by construction and
$(x,y)^{s-1}S\subset (L_1,{\bf h}):\beta$
since by (\ref{biforms} and \ref{biforms2}) the biforms $\hh$, or
$\ff, \hh$,  must effectively involve all monomials of degree $s-1$ in $x,y$.
It follows that $(x,y)^{s-1}\BB\subset Q$, hence
$$\LL(x,y)^{s-1}\subset \LL Q\subset \LL\cap Q=N=(L_1):(x,y)^{r},$$
thus implying that $\LL(x,y)^{s+r-1}\subset (L_1)$.
This shows the assertion.
\end{proof}

%\subsubsection*{Elimination equation up to a power}

\begin{Theorem} \label{detB2plus} Let $I$ be as above and $\beta=\det \CC$. Then
$\beta$ is a power of the elimination equation of $I$.
\end{Theorem}

\begin{proof}
Let $\pp$ denote the elimination equation of $I$.
Since $\pp$ is irreducible it suffices to show that $\beta$ divides
a power of $\pp$.

The associated primes of $N=(L_1): (x,y)^{r}$ are
the defining ideal $L$ of the Rees algebra and
 $\mathfrak{m}\BB= (x,y)\BB$.
 We have a primary decomposition
\[ N= \LL \cap Q,\]
where $Q$ is $\mathfrak{m}\BB$-primary.
  From the proof of
  Theorem~\ref{detC}, localizing at $\mathfrak{m}\BB$ gives
$ (x,y)^{s-1}\BB\subset Q$. (Equality will hold
  when $r=s$.)

 The equality $N= \LL\cap Q= (L_1, \hh) $ implies that
$(x,y)^{s-1}\pp \subset (L_1,\hh)$.
Since $f,g$ are of bidegrees $(r,1)$ and $(s,1)$, it must be
the case that each polynomial
 $x^{i}y^{s-1-i}\pp $ lies in the span of the $(f(x,y)^{s-r-1}, \hh)$
alone.   This gives a representation
\[ \pp[(x,y)^{s-1}]= [\ff;\hh] \cdot \AA,\]
(or simply $\pp[(x,y)^{r-1}]= [\hh] \cdot \AA$, if $s=r$)
where $\AA$ is an $s\times s$ matrix with entries in $k[\TT_1, \TT_2, \TT_3]$.
Replacing $[\ff;\hh]$ by $[(x,y)^{s-1}] \cdot \CC$, gives the matrix
equation
\[ [(x,y)^{s-1}] \big( \CC\cdot \AA-\pp \mathbf{I}\big)=0, \]
where $\mathbf{I}$ is the $s\times s$ identity matrix.

Since the minimal syzygies of $(x,y)^{s-1}$ have coefficients in
$(x,y)$, we must have
 \[\CC\cdot \AA=\pp \mathbf{I}, \]
so that $\det \CC\cdot \det \AA= \pp^{r}$, as desired. 
\end{proof}

\subsubsection*{Ternary ideals}\label{ternary}

We outline a conjectural scenario that we expect many such
ideals  to conform to. Suppose $I$ is an ideal of $\RR=k[x_1,x_2,x_3]$
generated by forms
$a_1,a_2, a_3, a$ of degree $n\geq 2$, with $J=(a_1,a_2,a_3)$ being a
minimal reduction. This approach is required because the linkage
theory method lacks the predictability  of the  binary ideal case.

\bigskip

\subsubsection*{Balanced ternary ideals}

Suppose that $I$ is $(n-1)$--balanced, where $n$ is the degree
 of the generators of $I$.

By Corollary~\ref{wellbal}(a)
 there are $n$ linear forms $\ff_i$, $1\leq
i\leq n$,
\[ \ff_i= \sum_{j=1}^4 c_{ij}\TT_j\in L_1,\]
arising  from the syzygies of $I$ of degree $n-1$. These syzygies,
according to Theorem~\ref{Dubreil}, come from the syzygies of $J:a$,
which by the structure theorem of codimension three Gorenstein
ideals, is given by the Pfaffians of a skew-symmetric matrix $\Phi$,
of size  at most $2n-1$.

According to
Proposition~\ref{syzlemma0}, there are
$n\choose 2$ quadratic forms $\hh_k$  ($1\leq k\leq {{n}\choose
{2}}$):
\[ \hh_k= \sum_{1\leq i\leq j\leq 4} c_{ijk}\TT_i\TT_j\in L_2,\]
with $\RR$-coefficients of degree $n-1$.

Picking a basis for $\mathfrak{m}^{n-1}$ (simply denoted by
$\mathfrak{m}^{n-1}$),
 and writing the $\ff_i$ and $\hh_k$ in matrix format, we
\begin{eqnarray} \label{acs3sedeg}
 [\ff_1, \ldots, \ff_n, \hh_1, \ldots, \hh_{{n}\choose{2}}]&=&
\mathfrak{m}^{n-1}\cdot
\BB,
\end{eqnarray}
where $\BB$ is the corresponding content matrix (see \cite[???]{dual}).
Observe that $\det\BB$ is either zero, or a polynomial of degree
\[ n+ 2{{n}\choose{2}}= n^2.\]
It is therefore a likely candidate for the {\em elimination
equation}. Verification
consists in checking that $\det\BB$ is irreducible for an ideal in
any
given {\em generic} class. We will make this more precise
 on a quick analysis of the lower degree cases.

\begin{Theorem} \label{detB3s} If $I\subset R=k[x_1,x_2,x_3]$ is a  $(n-1)$-balanced almost complete
intersection ideal  generated by forms of degree $n$ {\rm (}$n\leq 7${\rm )}, then
\begin{eqnarray*}
 \det \BB\neq 0.
\end{eqnarray*}
\end{Theorem}

\begin{proof}  Write each of the
quadrics ${\hh}_j$ above in the form
\[ {\hh}_j = c_j\TT_4^2 + \TT_4\ff(\TT_1,\TT_2, \TT_3)+
\g2(\TT_1,\TT_2, \TT_3),
 \]
where $c_j$ is a form of $\RR$ of degree $n-1$, and similarly,
\[ \ff_i= c_{i4} \TT_4+ \sum_{k=1}^3 c_{ik}\TT_k \in (L_1).\]

\medskip

Write $\mathfrak{c}=(c_j, c_{i4})$ for
the ideal of $\RR$ generated by the leading coefficients of
$\TT_4$ in the $\ff_i$'s and of $\TT_4^2$ in the ${\hh}_j$'s. It is apparent that if $\mathfrak{c}=\mathfrak{m}^{n-1}$, there will be a non-cancelling
term $\TT_4^{n^2}$ in the expansion of $\det \BB$.

To argue that indeed $\mathfrak{c}=\mathfrak{m}^{n-1}$ is the case, assume otherwise.
Since the $\ff_i$ are minimal generators
that contribute to $(J:a)$, we may assume that the $c_{i4}$ are
linearly independent. This implies that  we may replace one of the ${\hh}_j$ by a
form
\begin{eqnarray*}{\hh} &=&  \TT_4\ff(\TT_1,\TT_2, \TT_3)+
\g2(\TT_1, \TT_2, \TT_3)\\
&=&
  \TT_4(\RR_1\TT_1+  r_2\TT_2+r_3 \TT_3)+
\TT_1\g2_1+  \TT_2\g2_2+ \TT_3\g2_3,
\end{eqnarray*}
where the $r_i$ are $(n-1)$-forms in $\RR$ and the $\g2$'s are $\TT$-linear
involving  only $\TT_1,\TT_2,\TT_3$.

Evaluate now $\TT_i$ at the corresponding generator of $I$ to get \[
(ar_1+\g2_1(a_1,a_2,a_3) )\, a_1 + (ar_2+\g2_2(a_1,a_2,a_3) )\, a_3+
(ar_3+\g2_3(a_1,a_2,a_3) )\, a_3=0,\] a syzygy of the ideal $J=(a_1,a_2,a_3)$.
Since $J$ is a complete intersection, $ar_i+\g2_i(a_1,a_2,a_3)={\bf
u}_i(a_1,a_2,a_3)\in J$ for $i=1,2,3$, with ${\bf u}_i$ a linear form
in $\TT_1,\TT_2,\TT_3$ with coefficients in $\RR$.  These are syzygies
of the generators of $I$, so lifting back to $1$-forms in $\TT$ and
substituting yields $\hh=\hh'+{\bf k}$, where
\begin{eqnarray*} \hh'&=&(\RR_1\TT_4+\g2_1(\TT_2,\TT_2,\TT_3)-{\bf u}_1(\TT_1,\TT_2,\TT_3))\,\TT_1\\
&+& (\RR_2\TT_4+\g2_2(\TT_1,\TT_2,\TT_3)-{\bf u}_2(\TT_1,\TT_2,\TT_3))\,\TT_2\\
&+& (\RR_3\TT_4+\g2_3(\TT_1,\TT_2,\TT_3)-{\bf u}_3(\TT_1,\TT_2,\TT_3))\,\TT_3
\end{eqnarray*}
is an element of $(L_1)$, and because $\hh$ is a relation then so is the term ${\bf k}$.
But the latter only involves $\TT_1,\TT_2,\TT_3$, hence it belongs to the defining ideal of the symmetric
algebra of the complete intersection $J$. But the latter is certainly contained in $(L_1)$.
Summing up we have found that $\hh\in (L_1)$, which is a contradiction
since $\hh$ is a minimal
generator of $L_2/\SS_1L_1$. 
\end{proof}

\begin{Example} {\rm Let  $I=(J,a)$, where
$J=(x_1^3,\,x_2^3,\,x_3^3)$ and $a=x_1x_2x_3$.

\begin{enumerate}
\item[$\bullet$] $\nu(J:a)=3$.

\item[$\bullet$] The Hilbert series of $\RR/(J:a)$ is $1+3t+3t^2+t^3$.

\item[$\bullet$] $(L_1) : \m^4 = (L_1) : \m^5$ while $(L_1) : \m^3 \neq (L_1) : \m^4$

\item[$\bullet$] Equations of $I$:
\[\left\{\begin{array}{ll} L=(L_1,L_2,L_3)& \\
 \nu(L_1)=6; \nu(L_2)=3; L_3=kF, &
\end{array}
\right.
\]
where $F=-\TT_1\TT_2\TT_3+\TT_4^3$ is the elimination equation.

\item[$\bullet$]The corresponding map $\Psi_I$ is not birational onto its image:
the elimination equation has
degree $3$, as one easily computes by hand: $F=-\TT_1\TT_2\TT_3+\TT_4^3$.

\item[$\bullet$] $I_1(\varphi)=\m^2$, i.e., $I$ is $2$--balanced.

\noindent This shows that there
are $(n-1)$--balanced ideals $I\subset k[x_1,x_2,x_3]$ generated in degree $n$
such that the corresponding map $\Psi_I$ is not birational onto its image.

\item[$\bullet$] Let $\ff_1, \ff_2, \ff_3$ be a part of generators of
$L_1$ with coefficients in $\m^2$ and $L_2=(h_1, h_2, h_3)$. Then we
write $[\ff_1, \ff_2, \ff_3, \hh_1, \hh_2, \hh_3]=\m^2 \BB$ (Equation
(\ref{acs3sedeg}), where \[ \BB = \left[\begin{array}{rrrrrr} 0& 0&
-\TT_4 & \TT_3 & 0 & 0 \\ 0 & -\TT_4 & 0 & 0 & \TT_2 & 0 \\ -\TT_4 & 0& 0 & 0 &
0& \TT_1 \\ \TT_2 \TT_3 & 0 & 0& 0& 0 &-\TT_4^2 \\ 0 & \TT_1\TT_3
&0&0 &-\TT_4^2 &0
\\ 0 & 0 & \TT_1\TT_2 &-\TT_4^2 &0 & 0 \end{array} \right] \]

\item[$\bullet$] $\det \BB = F^3$.
\end{enumerate}
}\end{Example}

\begin{Remark}{\rm To strengthen Theorem~\ref{detB3s} to the
assertion that $\det \BB$ is a power of the elimination equation, one
needs more understanding of the ideal $(L_1):\mathfrak{m}^{n-1}$.
Here is one such instance.
}\end{Remark}

\begin{Proposition}\label{being_a_power}
Let $I\subset R=k[x_1,x_2,x_3]$ be an $(n-1)$-balanced almost complete
intersection ideal  generated by forms of degree $n$.
Keeping the notation introduced at the
beginning of this section, we obtain the following :
\begin{enumerate}
\item[{\rm (a)}] $2n-2$ is a secondary elimination degree of $I$.

\item[{\rm (b)}] If $(L_1):\mathfrak{m}^{n-1} = ({\bf f}, {\bf h})=(\ff_1, \ldots, \ff_n, \hh_1, \ldots, \hh_{{n}\choose{2}})$, then  $\det \BB$ is a power of the elimination equation.
\end{enumerate}
\end{Proposition}

\begin{proof}
Let $(L_1):\mathfrak{m}^{n-1}=L\cap Q$, $Q$ an
$\mathfrak{m}S$-primary ideal. As in  (\ref{acs3sedeg}) one has
$ [\ff, \hh]= \m^{n-1}\cdot \BB$. Notice that  $(\ff, \hh) \subset ( (L_1) : \m^{n-1}) $.
Write $\beta=\det \BB$. Since $\det \BB \neq 0$,
it follows that
\[ \m^{n-1}S \subset (\ff, \hh) : \beta \subset ( (L_1) : \m^{n-1}) : \beta =Q.\]
This implies that
\[ L \cdot \mathfrak{m}^{n-1} \subset  LQ
\subset (L_1):\mathfrak{m}^{n-1}.\]
Hence
$L= (L_1):\mathfrak{m}^{2n-2}$, which proves (i).

Now suppose that $(L_1):\mathfrak{m}^{n-1} = ({\bf f}, {\bf h})$. By (i), we have
 \[\pp \in (L_1):\mathfrak{m}^{2n-2} = ( (L_1):\m^{n-1} ) : \m^{n-1}. \] Therefore
$\pp \mathfrak{m}^{n-1}\subset (\ff, \hh)$, which gives a representation
\[ \pp [\mathfrak{m}^{n-1}]= [\ff, \hh] \cdot \AA,\]
where $\AA$ is a square matrix with entries in $\SS$.
Replacing $[\ff,\hh]$ by $[\mathfrak{m}^{n-1}] \cdot \BB$, gives the matrix
equation
\[ [\mathfrak{m}^{n-1}] \big( \BB\cdot \AA-\pp \mathbf{I}\big)=0, \]
where $\mathbf{I}$ is the  identity matrix. Since the minimal syzygies of $\mathfrak{m}^{r-1}$ have coefficients in
$\mathfrak{m}$, we must have
 \[\BB\cdot \AA=\pp \mathbf{I}, \]
so that $\det \BB\cdot \det \AA= \pp^m$, for certain $m$; this proves (ii).
 \end{proof}

\subsubsection*{Ternary quadrics}

We apply the preceding discussion to the situation where the ideal
$I$ is generated by $4$ quadrics of the
polynomial ring $\RR=k[x_1,x_2,x_3]$. The socle of
$\RR/J$ is generated by the Jacobian determinant of $a_1,a_2,a_3$,
which implies that $\lambda(I/J)\geq 2$. Together we obtain that
$\lambda(\RR/I)=6$. The Hilbert function of $\RR/I$ is $(1,3, 2)$. Since
we cannot have $u\mathfrak{m}\subset I$ for some $1$-form $u$, the
type of $I$ is $2$ and its socle is generated in degree three.

The canonical module of $\RR/I$ satisfies  $\lambda((J:a)/J)=6$, hence
$\lambda(\mathfrak{m}/(J:a))=1$, that is to say
$J:a=(v_1,v_2,v_3^2)$, where the $v_i$ are linearly independent
$1$-forms.
$\RR/I$ has a free presentation
\[ 0 \rar R^2\lar R^5\stackrel{\varphi}{\lar} R^4\lar R/I \rar 0.\]

The ideal $I_1(\varphi)$ is either
$\mathfrak{m}$ or  $(v_1,v_2,v_3^2)$.
 In the first case, $\delta(I)$ is the socle of
$H_1(I)$, an element of degree $3$,
 therefore its image in $ {L}_2$ is a $2$-form $h_1$
with linear coefficients. In particular reduction number cannot be
two.
Putting it together with the two
linear syzygies of $I$, we
\[ [\ff_1, \ff_2, \hh_1] = [x_1,x_2,x_3]\cdot \BB, \]
where $\BB$ is a $3\times 3$ matrix with entries in $k[\TT_1,
\TT_2,\TT_3,\TT_4]$, of column degrees $(1,1,2)$. The quartic $\det \BB$,
 is the elimination equation of $I$.

In the other case, $\delta(I)$ has degree two, so its image in
$ {L}_2$ is a form with coefficients in the field. The
corresponding mapping  $\Psi_I$ is not birational.

\medskip

We sum up the findings in this case:

\begin{Theorem} \label{aci32} Let $\RR=k[x_1, x_2, x_3]$ and $I$ a
$\mathfrak{m}$-primary
almost complete intersection generated by quadrics. Then
\begin{enumerate}
\item[{\rm (a)}]  If $I_1(\varphi)=\mathfrak{m}$ the corresponding
mapping $\Psi_I$ is
birational onto its image.

\item[{\rm (b)}]  If  $I_1(\varphi)\neq\mathfrak{m}$ then $I_1(\varphi)= (v_1, v_2, v_3^2)$, where $v_1,v_2,
v_3$ are linearly independent $1$-forms,  and the mapping $\Psi_I$ is
not birational onto its image.
\end{enumerate}
\end{Theorem}

Here is a sufficiently general example fitting the first case in the above theorem
(same behavior as $4$ random quadrics):
\begin{equation}\label{quadrics}
J=(x_1^2, x_2^2, x_3^2), \quad a=x_1x_2+x_1x_3+x_2x_3.
\end{equation}
A calculation shows that $\det \BB$  is irreducible.

\bigskip

An example which is degenerate (non-birational) as in the second case above is

\begin{equation}\label{quadrics_degenerate}
J=(x_1^2, x_2^2, x_3^2), \quad a=x_1x_2.
\end{equation}
In this case $\det \BB $ is a square, so we still recover the
elimination equation from $\BB$.

\subsubsection*{Ternary cubics and quartics}

Let $I$ be an ideal generated by $4$ cubics and suppose that $I$ is
$2$-balanced (i.e.,
$I_1(\varphi)=\mathfrak{m}^2$).
Using this (see the beginning of Section~\ref{ternary}) and the fact that
 $J:a$ is a codimension $3$ Gorenstein ideal, it follows that $J:a$ is minimally generated by the Pfaffians of
a skew-symmetric matrix of sizes $3$ or $5$.

In the first case, $J:a$
is generated by $3$ quadrics and $I$ is a Northcott ideal. In the
second  case, $J:a$
cannot be generated by $5$ quadrics, as its Hilbert function would be
$(1,3,1)$ and therefore the Hilbert function of $\RR/I$ would have to be
\[ (1,3,6,7,6,3,1)- (0,0,0, 1,3,1)= (1,3,6,6, 3,2,1),\]
giving that the canonical module of $\RR/I$ had a generator in degree
$0$. Thus $J:a$ must be generated by
 $3$ quadrics and $2$ cubics.

\medskip

In both cases the ideal $I$ is well balanced.

\begin{Theorem} \label{aci33and4} Let $\RR=k[x_1, x_2, x_3]$ and $I$ a
generic well balanced
$\mathfrak{m}$-primary
almost complete intersection. Then
\begin{enumerate}
\item[{\rm (a)}] If $I$
 generated by cubics and  $I_1(\varphi)=
\mathfrak{m}^2$,
 the polynomial
$\det \BB$,
defined by the equation {\rm (\ref{acs3sedeg})}, is the elimination
equation of $I$.

\item[{\rm (b)}] If $I$
 generated by quartics and  $I_1(\varphi)=
\mathfrak{m}^3$,
 the polynomial
$\det \BB$,
defined by the equation {\rm (\ref{acs3sedeg})}, is the elimination
equation of $I$.
\end{enumerate}
\end{Theorem}

\begin{proof} In the following table we give typical cubics. In the third
example, $\det \BB $ is divisible by $ (\TT_1\TT_2\TT_3-\TT_4^3)$.

\bigskip

\begin{small}
\label{tab:1}
\begin{tabular}{|l|c|c|l|}
\hline\noalign{\smallskip}
$I$ & $\nu(J:a)$ & $\edeg(I)$ & $L$      \\
\hline\noalign{\smallskip}
$(x_1^3+x_2^2x_3,\,x_2^3+x_1x_3^2,\, x_3^3+x_1^2x_2, \,  a=x_1x_2x_3)$
 &  $3$  & $9$ & $(L_1,L_2,L_5, L_9=\det \BB)$ \\
\hline
$(x_1^3,\,x_2^3,\,x_3^3,  a=x_1^2x_2+x_2^2x_3+x_1x_3^2)$ & $5$ & $9$ &
$(L_1, L_2, L_4, L_9=\det \BB)$
\\
\hline
 $ (x_1^3,\,x_2^3,\,x_3^3,  a=x_1x_2x_3)$ & $3$ & $3$ & $(L_1, L_2,
 L_3)$ \\
\hline
\end{tabular}
\end{small}
\end{proof}

\begin{Example}
{\rm Let  $I=(J, a)$, where
$J=(x_1^3+x_2^2x_3,\,x_2^3+x_1x_3^2,\, x_3^3+x_1^2x_2)$ and $a=x_1x_2x_3$.

\begin{enumerate}
\item[$\bullet$] $\nu(J:a)=3$ (all generators are quadrics)

\item[$\bullet$] $I_1(\varphi)=\m^2$

\item[$\bullet$] The Hilbert series of $\RR/(J:a)$ is $1+3t+3t^2+t^3$

\item[$\bullet$] $(L_1) : \m^4 = (L_1) : \m^5$ while $(L_1) : \m^3 \neq (L_1) : \m^4$

\item[$\bullet$] Equations of $I$:
\[\left\{\begin{array}{ll} L=(L_1,L_2,L_5,L_9)& \\
 \nu(L_1)=6; \nu(L_2)=3; \nu(L_5)=15; L_9=k\det \BB, &
\end{array}
\right.
\]
where $\det \BB$ is the elimination equation.

\item[$\bullet$] $\edeg(I)=9=3^2$ (the mapping $\Psi_I$ is birational onto its image)

\end{enumerate}
}
\end{Example}

\begin{Example}
{\rm Let $I=(J, a)$, where
$J=(x_1^3,\,x_2^3,\,x_3^3)$ and $a=x_1^2x_2+x_2^2x_3+x_1x_3^2$.

\begin{enumerate}
\item[$\bullet$] $\nu(J:a)=5$ ($3$ quadratic and $2$ cubic Pfaffians))

\item[$\bullet$] $I_1(\varphi)=\m^2$

\item[$\bullet$] Hilbert series of $\RR/(J:a)$ is $1+3t+3t^2+t^3$.

\item[$\bullet$] $(L_1) : \m^4 = (L_1) : \m^5$ while $(L_1) : \m^3 \neq (L_1) : \m^4$.

\item[$\bullet$] Equations of $I$:
$$\left\{\begin{array}{ll} L=(L_1,L_2,L_4,L_9)& \\
 \nu(L_1)=7; \nu(L_2)=3; \nu(L_4)=6; L_9=k\det \BB, &
\end{array}
\right.
$$
where $\det \BB$ is the elimination equation.

\item[$\bullet$] $\edeg(I)=9=3^2$ (the mapping $\Psi_I$ is birational onto its image)

\end{enumerate}

}
\end{Example}

\begin{Example}{\rm
In the case of $3$-balanced quartics,
 $\det \BB$ is a polynomial of degree $16$,
\begin{eqnarray} \label{acs34edeg}
 [\ff_1, \ldots, \ff_4, \hh_1, \ldots, \hh_6]&=&
\mathfrak{m}^{3}\cdot
\BB.
\end{eqnarray}

A typical example is
$I= (x^4,y^4,z^4, x^3y+y^3z+z^3x) $,
$I_1(\varphi)=\mathfrak{m}^3$ and $\det\BB$ is irreducible.

}
\end{Example}

\subsubsection*{Quaternary quadrics}

Let $\RR=k[x_1, x_2,x_3, x_4]$ with $\m=(x_1, x_2,x_3, x_4)$.
Let $I$ be generated by $5$ quadrics $a=a_1, a_2, a_3, a_4, a_5 $ of
$\RR$, $J=(a_2, a_3, a_4, a_5)$.
The analysis of this case is
less extensive than the case of ternary ideals. Let us assume that
$I_1(\varphi)=\mathfrak{m}^2$ -- that is,  the balancedness exponent equals the degree
of the generators, a case that occurs generically in this degree.

We claim that the Hilbert function of $\RR/I$ is $(1,4,5)$.
Since
we cannot have $u\mathfrak{m}\subset I$ for some $1$-form $u$,
 its socle is generated in degree two or higher. First, we argue that
 $J:a\neq \mathfrak{m}^2$; in fact,  otherwise $\RR/I$  would be of length $11$
 and type $6$ as $\mathfrak{m}^2/J$
 is its canonical module. But then the Hilbert function of  $\RR/I$ would be $(1,4,5,1)$, and the
 last two graded components would be in the  socle, which is impossible.

Thus it must be the case that $\lambda(\RR/I)=10$ and the Hilbert function of $H_1(I)$ is
$(5,4,1)$. Note that $\nu(\delta(I))=4$.

The last two graded components are of degrees $3$ and $4$. In degree
$3$ it leads to $4$ forms $\qq_1, \qq_2, \qq_3, \qq_4$ in $L_2$, with linear
coefficients:
\[ [\qq_1, \qq_2, \qq_3, \qq_4]= [x_1, x_2, x_3, x_4]\cdot \BB.\]

\begin{Theorem} \label{detB42}
$\det \BB\neq 0$.
\end{Theorem}

\begin{proof} We follow the pattern of the proof of Theorem~\ref{detB3s}.  Write each of the
quadrics $\qq_i$ as
\[ \qq_i = c_i\TT_5^2 + \TT_5\ff(\TT_1 \ldots, \TT_4)+
\g2(\TT_1, \ldots, \TT_4),
 \]
where $c_i$ is a linear form of $\RR$. If $(c_1,c_2,c_3, c_4)= \mathfrak{m}$
we  can take the $c_i$ for indeterminates in order to obtain the
corresponding $\det \BB$. In this case the occurrence of a non-cancelling term
$\TT_5^8$ in $\det\BB$ would be clear.

By contradiction, assume that the forms $c_i$ are not linearly
independent. In this case, we could replace one of the $\qq_i$ by a
form \[ \qq=  \TT_5\ff(\TT_1, \ldots, \TT_4)+
\g2(\TT_1, \ldots, \TT_4). \]
Keeping in mind that $\qq$ is a minimal generator we are going to argue
that $\qq \in (L_1)$. For that end, write
the form  $\qq$ as
\[ \qq=
  \TT_5(r_1\TT_1+  \cdots+r_4 \TT_4)+
\TT_1\g2_1+  \cdots + \TT_4\g2_4,  \]
where the $r_i$ are $1$-forms in $\RR$ and $\g2_i$'s are $\TT$--linear involving only $\TT_1, \ldots, \TT_4$.

Evaluate now the leading $\TT_i$ at the ideal to get the syzygy
\[ \hh= a(r_1\TT_1+  \cdots+r_4 \TT_4)+
a_1\g2_1+  \cdots + a_4\g2_4,  \]
of $I$, but actually of the ideal $J$. Since $J$ is a complete
intersection, all the coefficients of $\hh$ lie in $J$.
This implies that $r_ia\in J$ for all $r_i$. But this is impossible
since $J:a\in I_1(\varphi)=\mathfrak{m}^2$, unless all $r_i=0$.
This would imply that $\qq\in (L_1)$, as asserted. 
\end{proof}

\begin{Theorem} \label{aci42} Let $\RR=k[x_1, x_2, x_3, x_4]$ and $I$ a
generic
$\mathfrak{m}$-primary
almost complete intersection.
 If $I$ is
 generated by quadrics and  $I_1(\varphi)=
\mathfrak{m}^2$,
 the polynomial
$\det \BB$,
defined by the equation $({\rm \ref{acs3sedeg}})$ is divisible by the
elimination
equation of $I$.
\end{Theorem}

\subsubsection*{Primary decomposition}

We now derive a value for the secondary elimination degree via a
primary decomposition.

\begin{Proposition} Let $I=(a,b,c,d,e)$ be an ideal of
$\RR=k[x_1,x_2,x_3,x_4]$ generated by quadrics such that $I_1(\varphi)=
\mathfrak{m}^2$.
Let $L$ be the defining ideal of the Rees algebra of $I$.
Then \[ (L_1): \mathfrak{m}^2= L \cap \mathfrak{m}S.\]
In particular,
\[L= (L_1): \mathfrak{m}^3. \]
\end{Proposition}

By Theorem~\ref{secondedeg}, it will suffice to show:

\begin{Lemma} $\mathfrak{m}^3= (J:a)\mathfrak{m}$.
\end{Lemma}

\begin{proof} We make use of the diagram

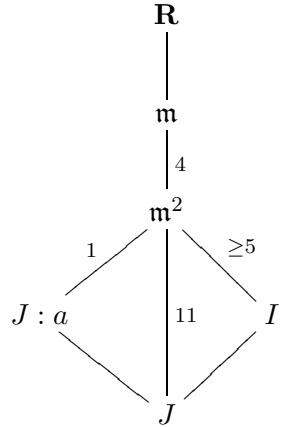
\begin{eqnarray*}
\diagram
& \RR \dline
  & \\
& \mathfrak{m}\dline^{4} & \\
& \mathfrak{m}^2\dlline_{1} \drline^{\geq 5} \ddline^{11} & \\
 J:a \drline & & I \dlline \\
 & J & \\
\enddiagram
\end{eqnarray*}

Let  $J=(b,c,d,e)$ be a
minimal reduction of $I$.
Since $J:a$ is a Gorenstein ideal, it cannot be $\mathfrak{m^2}$.
On the other hand,
$\lambda(J:a/J)= \lambda(\RR/I)\geq 10$. Thus
$\lambda(\mathfrak{m}^2/(J:a))=1$ and
$(J:a):\mathfrak{m}=\mathfrak{m}^2$ defines the socle of $\RR/(J:a)$.
The Hilbert function of $\RR/J:a$ is then $(1,4,1)$ which implies that
$\mathfrak{m}^3= (J:a)\mathfrak{m}$. 
\end{proof}

\subsubsection*{Examples of quaternary quadrics}

We give a glimpse of the various cases.

\begin{Example}{\rm Let $I=(J, a)$, and $J=
(x_1^2,\,x_2^2,\,x_3^2,\, x_4^2)$ and $a=x_1x_2+x_2x_3+x_3x_4$

\begin{enumerate}
\item[$\bullet$] $\nu(J:a)=9$ (all generators are quadrics)

\item[$\bullet$] $I_1(\varphi)=\m^2$  (i.e., $2$-balanced)

\item[$\bullet$] $\edeg(I)=8=2^3$ (i.e., birational)

\item[$\bullet$] As explained above, there are
four Rees equations $q_j$ of bidegree $(1,2)$ coming from the syzygetic principle.
One has:
\[\left\{\begin{array}{ll}
L=(L_1,L_2,L_9),& \\   %\mbox{where   $\nu(L)=20$}\\
\nu(L_1)=15; \nu(L_2)=4;  L_9=k\pp, & \mbox{where $\pp=\det\BB$}
\end{array}
\right.
\]
\end{enumerate}
}
\end{Example}

The following example is similar to the above example, including the syzygetic principle, except that
$F$ is now the square root of $\det\BB$.

\begin{Example}{\rm Let $I=(J,a)$, where
$J=(x_1^2,\,x_2^2,\,x_3^2,\, x_4^2)$ and  $a=x_1x_2+x_3x_4.$

\begin{enumerate}
\item[$\bullet$] $\nu(J:a)=9$ (all generators are quadrics)

\item[$\bullet$] $I_1(\varphi)=\m^2$

\item[$\bullet$] The Hilbert series of $\RR/(J:a)$ is $1+4t+t^2$

\item[$\bullet$] $(L_1) : \m^3 = (L_1) : \m^4$ while $(L_1) : \m^2 \neq (L_1) : \m^3$

\item[$\bullet$] Equations of $I$:
\[\left\{\begin{array}{ll} L=(L_1,L_2,L_4)& \\
 \nu(L_1)=15; \nu(L_2)=4; L_4=k\pp. &
\end{array}
\right.
\]

\item[$\bullet$] $\edeg(I)=4$ (the mapping $\Psi_I$ is non-birational onto its image)

\item[$\bullet$] Let $\qq_1, \qq_2, \qq_3, \qq_4$  be $4$ generating forms in $L_2$, with linear
coefficients. Write
\[ [\qq_1, \qq_2, \qq_3, \qq_4]= [x_1, x_2, x_3, x_4]\cdot \BB,\]
Then $\det \BB =\pp^2$.
\end{enumerate}
}\end{Example}

\begin{Example}{\rm Let $I=(J,a)$, where
$J=(x_1^2,\,x_2^2,\,x_3^2,\, x_4^2)$ and $a=x_1x_2+x_2x_3+x_3x_4+x_1x_4.$

\begin{enumerate}
\item[$\bullet$] $\nu(J:a)=4$ (complete intersection of $2$ linear equations and $2$ quadrics)

\item[$\bullet$] $I_1(\varphi)=\m$ (i.e., $1$-balanced, not $2$-balanced)

\item[$\bullet$] $\nu(Z_1)=9$ ($2$ linear syzygies and $7$ quadratic ones)

\item[$\bullet$] Hilbert series of $\RR/(J:a)$ : $1+2t+t^2$.

\item[$\bullet$] $L_1 : \m^3 = L_1 : \m^4$ but  $L_1 : \m^2 \neq L_1 : \m^3$

\item[$\bullet$] $\edeg(I)=8=2^3$ (i.e., birational)

\item[$\bullet$] The syzygetic principle does not work to get the
elimination equation:
One has:
\[\left\{\begin{array}{ll}
L=(L_1, L_2, L_3, L_8),& \\   %\mbox{where   $\nu(L)=20$}\\
\nu(L_1)=9; \nu(L_2)=1;  \nu(L_3) =2, L_8=k\pp, &
\end{array}
\right.\]
Nevertheless there are relationships among these numbers that are
understood from the syzygetic discussion. Thus
Theorem~\ref{socdegree}(c) implies that $L_2/\SS_1L_1$ is generated
by one form with linear coefficients.

\item[$\bullet$] Let $f_1$ and $f_2$ be in $L_1$ with linear coefficients, i.e.,
\begin{enumerate}
\item[$\diamondsuit$] $f_1=(x_1+x_3)(\TT_2-\TT_4)+(x_2-x_4)(-\TT_5)$

\item[$\diamondsuit$] $f_2=(x_1-x_3)(-\TT_5)+(x_2+x_4)(\TT_1-\TT_3)$
\end{enumerate}

\item[$\bullet$] Let $h+L_1\SS_1$ be a generator of $L_2/L_1\SS_1$. We
observed that, as a coset,  $h$ can be written as $h=h_1+\SS_1L_1$ and
$h=h_2+\SS_1L_1$, with $h_1$ and $h_2$ forms of bidegree $(2,2)$, with
$\RR$-content contained in the contents of $f_1$ and $f_2$, respectfully.

\item[$\bullet$] Let
$h_1=(-2x_1\TT_4\TT_5-2x_3\TT_4\TT_5+x_4\TT_5^2)(x_1+x_3)+(x_3\TT_1\TT_2
-x_3\TT_2\TT_3-x_3\TT_1\TT_4+x_3\TT_3\TT_4-
x_4\TT_1\TT_5+2x_2\TT_3\TT_5-x_4\TT_3\TT_5-x_3\TT_5^2)(x_2-x_4)$. Then $h_1 \in
L_2$ and $2h+h_1 \in L_1\SS_1$.  Hence we may choose $h_1+L_1\SS_1$
to be a generator of $L_2/L_1\SS_1$.

\item[$\bullet$] Let
$h_2=(x_4\TT_1\TT_2-x_4\TT_2\TT_3-x_4\TT_1\TT_4+x_4\TT_3\TT_4
-x_3\TT_2\TT_5+2x_1\TT_4\TT_5-x_3\TT_4\TT_5-x_4\TT_5^2)(x_1-x_3)+
(-2x_2\TT_3\TT_5-2x_4\TT_3\TT_5+x_3\TT_5^2)(x_2+x_4)$. Then $h_2 \in L_2$ and
$2h+h_2 \in L_1\SS_1$.  Hence we may choose $h_2+L_1\SS_1$ to be a
generator of $L_2/L_1\SS_1$.

\item[$\bullet$] Write
\[
[ f_1 \;\; h_1]  =  [x_1+x_3 \;\; x_2-x_4]\BB_1 \quad \mbox{\rm and} \quad
[ f_2 \;\; h_2]  = [x_1-x_3 \;\; x_2+x_4] \BB_2
\]
Then $\det \BB_1$ and $\det \BB_2$ form a minimal generating set of
$L_3/\SS_1L_2$.

\item[$\bullet$] Write $[f_1\;\; f_2\;\; \det \BB_1 \;\; \det \BB_2]
= \  [x_1 \;\; x_2 \;\; x_3 \;\; x_4]\BB$. Then
$\det \BB = \pp$.
\end{enumerate}

}\end{Example}

\subsection{Almost complete intersections of dimension $1$}

Let $(\RR,\mathfrak{m})$ be a Gorenstein local ring of dimension
$d\geq 1$, and let $I$ be a Cohen-Macaulay almost complete
intersection of dimension $1$.
Consider the syzygetic sequence of one of these ideals,
\begin{small}
\[ 0 \rar \delta(I) \lar \H_1(I) \lar (\RR/I)^d \lar I/I^2 \rar 0.\]
\end{small}
 $\H_1(I)$ is the canonical module of $\RR/I$. We can express
 $\delta(I)$ in the same manner that was done in the case of almost
 complete intersections of finite colength.

\begin{Proposition} Let $I$ be a Cohen-Macaulay ideal of dimension
one, generated
 by
$d$ elements. Then
\[ \delta(I) = \Hom_{R/I}(\RR/I_1(\varphi), \H_1(I)),\]
where $I_1(\varphi)$ is the ideal generated by the entries of the
presentation matrix of $I$. Moreover, if $I_1(\varphi)$ is mixed, one
can replace it in the formula for $\delta(I)$ by its component $I_1$ of
dimension one. In particular
$\delta(I)$ is generated by $r$ elements, where $r$ is the
Cohen-Macaulay type of $I_1$.
\end{Proposition}

\begin{proof} Since $\dim \RR/I=1$ and $I$ is Cohen-Macaulay, the image $L$ of
$\H_1(I) \rar (\RR/I)^d$ is a Cohen-Macaulay module. Dualizing with
$\H_1(I)$, gives the exact sequence
\[ 0\rar \Hom_{\RR/I}(L, \H_1(I)) \lar \Hom_{\RR/I}(\H_1(I), \H_1(I))
=\RR/I
\lar \Hom_{\RR/I}(\delta(I), \H_1(I)) \rar 0, \]
%\[ 0\rar \Hom_{R/I}(I/I^2, H_1(I)) \lar \Hom_{R/I}((\RR/I)^d, H_1(I))
%\lar \Hom_{R/I}(L, H_1(I)) \lar \Ext_{R/I}(I/I^2,H_1(I)\rar 0. \]
so $\Hom_{\RR/I}(\delta(I), \H_1(I)) = \RR/J$. Because this module is
Cohen-Macaulay and of dimension $1$, to determine $J$ we can localize
at the minimal prime ideals $\mathfrak{p}$ of $I$. Therefore, we have reduced the
question to the case of ideals of finite colength and obtain that
$J_{\mathfrak{p}}=I_1(\varphi)_{\mathfrak{p}}$.
By duality,
\[ \delta(I) = \Hom_{\RR/I}(\RR/J, \H_1(I)),\]
the canonical module of $\RR/J$.
The last assertion follows from the representation of $\delta(I)$ as
$\Hom_{\RR/I}(\RR/I_1(\varphi), \H_1(I)) $. 
\end{proof}

\section{Almost Cohen--Macaulay Algebras}\index{aCM: almost Cohen--Macaulay algebra}

\subsection{Introduction}
Our goal in this section is to treat the following class of rings:

\begin{Definition}{\rm Let $\RR$ be a local Noetherian ring (or a graded algebra over a local ring). If $\m$ is
its maximal ideal, we say that $\RR$ is an {\em almost} Cohen--Macaulay algebra \index{almost Cohen--Macaulay algebra}
if $\depth \RR \geq \dim \RR -1$.
}\end{Definition}

For simplicity we refer to them as aCM algebras. There is no apparent  for a separate designation for this
subclass of Noetherian rings. We will argue that in the context of certain Rees algebras this property
appears with frequency. Furthermore, the methods one uses in studying them and the numerical invariants they carry
are very similar to Cohen--Macaulay algebras themselves.
The properties we will focus on will be the calculation
of multiplicities and the Castelnuovo--Mumford regularity.
The results in this section can be found in \cite{acm}.

\subsection{The canonical presentation}

Let $\RR$ be a Cohen--Macaulay local ring of dimension $d$, or a polynomial ring
 $\RR=k[x_1, \ldots, x_d]$ for
$k$ a field.
By an {\em almost complete intersection} we mean an ideal
  $I=(a_1, \ldots, a_g, a_{g+1})$ of codimension $g$ where the subideal $J=(a_1, \ldots, a_g)$ is a
  complete intersection    and  $a_{g+1}\notin J$.
    By the {\em equations} of $I$ it is meant a free presentation of
 the Rees algebra $\RR[It]$ of $I$,
\begin{eqnarray}\label{presRees1}
 0 \rar \LL \lar \BB = \RR[\TT_1, \ldots, \TT_{g+1}] \stackrel{\psi}{\lar}
\RR[It] \rar 0,  \quad \TT_i \mapsto f_it .
\end{eqnarray}
More precisely,  $\LL$ is the defining ideal of the Rees
algebra of $I$ but we refer to it simply as the {\em ideal of
equations} of $I$. We are particularly interested in establishing the properties of $\LL$
when $\RR[It]$ is Cohen--Macaulay or has almost maximal depth.
This broader view requires a change of focus from $\LL$ to one of its quotients.
 We are going to study some classes of ideals whose
Rees algebras have these properties. They tend to occur in classes where the reduction number $\red_J(I)$ attains an extremal value.

\medskip

We first set up the framework to deal with properties of $\LL$ by a standard decomposition. We keep the notation of above, $I = (J,a)$. The presentation ideal $\LL$ of $\RR[It]$
 is a graded ideal
$\LL = L_1 + L_2 + \cdots$, where $L_1$ are linear forms in the $\TT_i$ defined by a	matrix $\phi$ of the syzygies of $I$,
$ L_1 = [\TT_1, \ldots, \TT_{g+1}] \cdot \phi $.
Our basic prism is given by the exact sequence
\[ 0 \rar \LL/(L_1) \lar \BB/(L_1) \lar \RR[It] \rar 0.\]
Here
$\BB/(L_1)$ is a presentation of the symmetric algebra of $I$ and
$\SS = \Symi(I)$ is a Cohen--Macaulay ring under very broad conditions,
including when $I$ is an ideal of finite colength. The emphasis here will be entirely on $T = \LL/(L_1)$, which we call the
{\em module of nonlinear} relations of $I$.\index{module of nonlinear relations of an ideal}
Let $\RR$ be a Cohen--Macaulay local domain of dimension $d\geq 1$ and let $I$ be an almost complete intersection as in Corollary~\ref{Zaci}. The ideal of equations $\LL$ can be studied  in two stages: $(L_1)$ and $\LL/(L_1)= T$:
\begin{eqnarray} \label{canoseque}
 0 \rar T \lar \SS = \BB/(L_1)= \Sym(I) \lar \RR[It] \rar 0.\end{eqnarray}
We will argue that this exact sequence is very useful. Note that $\Sym(I)$ and $\RR[It]$ have dimension $d+1$, and that $T$ is the $\RR$--torsion submodule of $\SS$.
Let us give some of its properties.

\begin{Proposition} \label{canoseq} Let $I$ be an ideal as above.
\begin{enumerate}

\item[{\rm (a)}]  $(L_1)$ is a Cohen-Macaulay ideal of $\BB$.

\item[{\rm (b)}] $T$ is a Cohen--Macaulay $\SS$--module if and only if $\depth \RR[It]\geq d$.

\item[{\rm (c)}] If $I$ is $\m$--primary then $\mathcal{N}=T\cap \m\SS$ is the nil radical of
$\SS$ and $\mathcal{N}^s = 0$ if and only if $\m^sT =0$. This is equivalent to saying that $\mbox{\rm sdeg}(I)$ is the index of nilpotency of $\Symi(I)$.

\item[{\rm (d)}] $T=\mathcal{N} + \mathcal{F}$, where $\mathcal{F}$ is a lift in $\SS$ of the relations in $\SS/\m \SS$ of the special
fiber ring $\mathcal{F}(I)=\RR[It] \otimes \RR/\m$.  In particular if $\mathcal{F}(I)$ is a hyersurface ring,
$T = (f, \mathcal{N})$.

\end{enumerate}

\end{Proposition}

\begin{proof} (i) comes from Corollary~\ref{Zaci}.
\medskip

\noindent (ii)  In the defining sequence of $T$,
\[ 0 \rar (L_1) \lar \LL \lar T \rar 0,\]
since $(L_1)$ is a Cohen--Macaulay ideal of codimension $g$, as an $\BB$--module, we have $\depth (L_1) = d+2$, while
$\depth \LL = 1 + \depth \RR[It]$. It follows that $\depth T =\min \{d+1,  1 +\depth \RR[It]\}$.
Since $T$ is a module of Krull dimension $d+1$ it is a Cohen--Macaulay module if and only if
$\depth \RR[It] \geq d$.

\medskip

\noindent (iii) $\m\SS$ and $T$ are both minimal primes and for large $n$, $\m^nT=0$. Thus $T$ and $\m\SS$ are the only minimal
primes of $\SS$, $\mathcal{N} = \m \SS \cap T$.
To argue the equality of the two indices of nilpotency, let $n$ be such that $\m^nT=0$. The ideal $\m^n\SS + T$
has positive codimension, so contains regular elements since $\SS$ is Cohen--Macaulay. Therefore to show
 \[\m^s T=0 \Longleftrightarrow \mathcal{N}^s =0\]
it is enough to multiply both expressions by $\m^n\SS + T$. The verification is immediate.

\medskip

\noindent (iv) Tensoring the sequence (\ref{canoseque}) by $\RR/\m$ gives the exact sequence
\[ 0 \rar \m\SS\cap T/\m T = \mathcal{N}/\m T \lar T/\m T \lar \SS/\m \SS \lar \mathcal{F}(I) \rar 0. \]
By Nakayama Lemma we may ignore $\m T$ and recover $T$ as asserted. 
\end{proof}

The main intuition derived from these basic observations is (ii): Whenever methods are developed to study the equations of
$\RR[It]$ when this algebra is Cohen--Macaulay, should apply [hopefully] in case they are almost Cohen--Macaulay.

\begin{Remark}{\rm If $I$ is not $\m$--primary
but still satisfies one of the other conditions of Corollary~\ref{Zaci},
the nilradical $\mathcal{N}$ of $\SS$ is given by $T\cap N_0\SS$, where $N_0$ is the intersection of the minimal
primes $\p$ for which $I_{\p}$ is not of linear type.

}
\end{Remark}

For almost complete intersections $I$ of finite colength there is another symmetric algebra that is an
ancestor to $\RR[It]$. Let $L_1^*$ be the submodule of $L_1$ defined by the Koszul relations,
\[ 0 \rar T^* = \LL/(L_1^*) \lar \SS^*= \BB/(L_1^*) \lar \RR[It]\rar  0.\]

This definition gives rise to the exact sequence
 \[ 0 \rar \H_1(I)\otimes \BB \lar T^* \lar T \rar 0,\]
where $\H_1(I)$ is the one-dimensional Koszul homology of $I$.

\begin{Proposition} $T$ is a Cohen--Macaulay ideal of $\SS$ if and only if $T^*$ is a Cohen--Macaulay ideal
of $\SS^*$.

\end{Proposition}

\begin{proof} This follows from Corollary~\ref{Zaci} but applied instead to the $\mathcal{B}$ complex instead. Both conditions
are equivalent to $\RR[It]$ being almost Cohen--Macaulay. 
\end{proof}

\subsubsection*{Reduced symmetric algebras}
\index{reduced symmetric algebra}

\begin{Proposition}\label{linkofmax} Let $\RR$ be a Gorenstein local domain of dimension $d$ and let $J$ be a parameter
ideal. If $J$ contains two minimal generators in $\m^2$,
then the Rees algebra of $I=J\colon \m$ is Cohen--Macaulay and
$\LL = ( L_1, \ff)$ for some quadratic form $\ff$.
\end{Proposition}

\begin{proof}
The equality $I^2 = JI$  comes from \cite{CPV1}. The Cohen--Macaulayess of $\RR[It]$ is a general argument
(in \cite{CPV1} and probably elsewhere). Let $\ff$ denote the quadratic form
\[ \ff = \TT_{d+1} ^2 + \mbox{\rm lower terms}.\]

Let us show that $\m T=0$.
Reduction modulo $\ff$ can be used to present any element in $\LL$ as
\[ F =  \TT_{d+1} \cdot A + B\in \LL, \]
where $A$ and $B$ are forms in $\TT_1,\ldots, \TT_d$. Since $I = J:\m$,
 any element in $\m \TT_{d+1}$ is equivalent, modulo
$L_1$, to a linear form in the other variables.
Consequently
\[\m F \subset (L_1, \RR[\TT_1, \ldots, \TT_d]) \cap \LL \subset (L_1),\]
as desired.

By Proposition~\ref{canoseq}, we have the exact sequence
\[ 0\rar T \lar \SS/\m \SS \lar \mathcal{F}(I) \rar 0.\]
But $T$ is a maximal Cohen--Macaulay $\SS$--module, and so it is also a maximal Cohen--Macaulay $\SS/\m \SS$--module as
well. It follows that $T$ is generated by a monic polynomial that divides the image of $\ff$ in $\SS/\m\SS$.
It is now clear that $T=(\ff)\SS$. 
\end{proof}

%\begin{Remark}{\rm Don't recall the extra fact that $\mathcal{F}(I)$ is Cohen--Macaulay was known. }
%\end{Remark}

\begin{Corollary} For the ideals above $\Sym(I)$ is reduced.
\end{Corollary}

%The condition on the generators of $J$ is only required if $\RR$ is a regular local ring.

We now discuss a generalization, but since we are still developing the examples, we are somewhat informal.

\begin{Corollary}\label{sdegs}
Suppose the syzygies of $I$ are contained in $\m^s\BB$ and that $\m^s\LL \subset  (L_1)$.  We have the exact sequence
\begin{eqnarray} \label{canoseque2}
 0 \rar T \lar \SS/\m^s\SS \lar \RR[It]\otimes \RR/\m^s = \mathcal{F}_s(I) \rar 0.
\end{eqnarray}
If $\RR[It]$ is almost Cohen--Macaulay, $T$ is a Cohen--Macaulay module that is an ideal of the polynomial
ring $\CC = \RR/(\m^s)[\TT_1, \ldots, \TT_{d+1}]$, a ring of multiplicity
${s+d-1}\choose{d}$.
Therefore we have that $\nu(T) \leq {{s+d-1}\choose {d}}$.
\end{Corollary}

Note that also here $\mathcal{F}_s(I)$ is Cohen--Macaulay. We wonder whether $\mathcal{F}(I)$ is Cohen--Macaulay.

\bigskip

Let $(\RR, \m)$ be a Cohen--Macaulay local ring and $I$ an almost complete intersection as in Corollary~\ref{Zaci}.
We examine the following surprising fact.

\begin{Theorem} \label{reducedsymi}    Suppose  that $\RR$ is a Cohen--Macaulay  local ring and  $I$ is
an $\m$--primary
 almost  complete intersection such that  $\SS=\Sym(I)$ is reduced.  Then
 $\RR[It]$ is an almost Cohen-Macaulay algebra.
 \end{Theorem}

 \begin{proof}
  Since
$0 =\mathcal{N} = T\cap \m \SS$,   on one hand from (\ref{canoseque}) we have that $T$ satisfies the condition $S_2$ of Serre, that is
\[ \depth T_{P} \geq \inf\{2, \dim T_P\}\]
for every prime ideal $P$ of $\SS$. On the other hand,
  from (\ref{canoseque2}) $T$ is an ideal of the polynomial ring $\SS/\m \SS$. It follows that
$T= (\ff )\SS$, and consequently $\depth \
 \RR[It] \geq d$.
\end{proof}

\begin{Example}{\rm  If $\RR = \mathbb{Q}[x,y]/(y^4-x^3)$, $J = (x)$ and $I = J: (x,y)= (x, y^3)$, $\depth \RR[It] = 1$.
}\end{Example}

There are a number of immediate observations.

\begin{Corollary} If $I$ is an ideal as in {\rm Theorem~\ref{reducedsymi}}, then the special fiber ring $\mathcal{F}(I)$ is Cohen-Macaulay.
\end{Corollary}

\begin{Remark}{\rm  If $I$ is an almost complete intersection as in (\ref{Zaci}) and its radical is a regular prime
ideal $P$, that is  $\RR/P$ is regular local ring, the same assertions will apply if $\Symi(I)$ is reduced.
}\end{Remark}

\subsubsection{Direct links of Gorenstein ideals}

We briefly outline a broad class of extremal Rees algebras.
Let $(\RR, \m)$ be a
Gorenstein local ring of dimension $d\geq 1$. A natural source of almost complete intersections in $\RR$ direct links of Gorenstein ideals. That
is, let $K$ be a Gorenstein ideal of $\RR$ of codimension $s$, that is $\RR/K$ is a Gorenstein ring of dimension $d-s$.
If $J = (a_1, \ldots, a_s)\subset K$ is a complete intersection of codimension $s$, $J\neq K$, $I = J:K$ is an almost complete intersection, $I = (J, a)$. Depending on $K$, sometimes these ideals come endowed with very good
properties. Let us recall one of them.

\begin{Proposition}\label{sourceofacis} Let $(\RR, \m)$ be a Noetherian local ring of dimension $d$.
\begin{enumerate}
\item[{\rm (a)}] {\rm (\cite[Theorem 2.1]{CPV1})} Suppose $\RR$ is a Cohen--Macaulay local ring and  let $\p$ be a prime ideal of codimension $s$ such that $\RR_{\p}$ is a Gorenstein ring and let $J$ be a complete intersection of
codimension $s$ contained in $\p$. Then for $I= J:\p$ we have $I^2= JI$ in the following two cases: {\rm (a)} $\RR_{\p}$ is not a regular local ring; {\rm (b)} if $\RR_{\p}$ is a regular local ring two of the elements $a_i$ belong to $\p^{(2)}$.

\item[{\rm (b)}] {\rm (\cite[Theorem 3.7]{CHV})}  	Suppose $J$ is an irreducible $\m$--primary ideal. Then {\rm (a)} either there exists a minimal set of generators $\{x_1, \ldots, x_d\}$ of $\m$ such that $J = (x_1,\ldots, x_{d-1}, {x_d}^r)$, or {\rm (b)}
$I^2 = JI$ for $I = J:\m$.

\end{enumerate}
\end{Proposition}

\medskip

The following criterion is a global version of Corollary~\ref{F2Cor}

\begin{Proposition}\label{rednumone} Let $\RR$ be a Gorenstein local ring and $I=(J,a)$ an almost complete intersection {\rm[}when we write $I=(J,a)$ we always mean that $J$ is a reduction{\rm]}. If $I$ is an unmixed ideal {\rm[}height unmixed{\rm]} then
$\red_J(I) \leq 1$ if and only if $J:a = I_1(\phi)$.
\end{Proposition}

\begin{proof} Since the ideal $JI$ is also unmixed, to check the equality $J:a= I_1(\phi)$ we only need to check at the minimal primes of $I$ (or, of $J$, as  they are the same). Now Corollary~\ref{F2Cor} applies. 
\end{proof}

If in Proposition~\ref{sourceofacis} $\RR$ is a Gorenstein local ring and $I$ is a Cohen--Macaulay ideal, their
associated graded rings are Cohen--Macaulay, while the
Rees algebras are also Cohen--Macaulay if $\dim \RR\geq 2$.

\begin{Theorem} \label{reesoflink} Let $\RR$ be a Gorenstein local ring and $I$ a Cohen--Macaulay ideal that is an  almost complete intersection. If $\red_J(I)\leq 1$ then in the canonical representation
\[ 0 \rar T \lar \SS \lar \RR[It] \rar 0,\]
\begin{itemize}
\item[{\rm (a)}] If $\dim \RR \geq 2$ $\RR[It]$ is Cohen--Macaulay.

\item[{\rm (b)}] $T$ is a Cohen--Macaulay module over $\SS/(I_1(\phi))\SS$, in particular
\[ \nu(T) \leq \deg \RR/I_1(\phi).\]
\end{itemize}
\end{Theorem}

\begin{Example}{\rm
Let $\RR=k[x_1, \ldots, x_d]$, $k$ an algebraically closed field, and let $\p$ be a homogeneous prime ideal of codimension $d-1$. Suppose $J= (a_1, \ldots, a_{d-1})$
is a complete intersection of codimension $d-1$ with at least two generators in $\p^2$. Since $\RR/\p$ is regular,
$I=J:\p$ is an almost complete intersection and $I^2 = JI$. Since $\p$ is a complete intersection, say
$\p = (x_1-c_1x_{d} , x_2-c_2x_d, \ldots, x_{d-1} -c_{d-1}x_d)$, $c_i\in k$,  we write  the matrix equation $J= \AA \cdot \p$, where
$\AA$ is a square matrix of size $d-1$. This is the setting where the Northcott ideals occur, and therefore
$I = (J, \det \AA)$.

\medskip

By Theorem~\ref{reesoflink}(b), $\nu(T) \leq \deg (\RR/\p) = 1$. Thus $\LL$ is generated by the syzygies of $I$ (which are well-understood) plus a quadratic equation.
}\end{Example}

\subsection{Metrics of aCM algebras}\index{metrics of aCM algebras}

Let $(\RR, \m)$ be a Cohen--Macaulay local ring of dimension $d$ and let $I$ be an almost complete intersection of
finite colength. We assume that $I=(J,a)$, where $J$ is a minimal reduction of $I$. These assumptions will hold
for the remainder of the section. We emphasize that they apply to the case when $\RR$ is a polynomial ring over
a field and $I$ is a homogeneous ideal.

\medskip

In the next statement we highlight the information about the equations of $I$ that is a direct consequences of the aCM hypothesis. In the next segments we begin to obtain the required data in an explicit form. As for
 notation, $\BB = \RR[\TT_1, \ldots, \TT_{d+1}]$ and
  for a graded
 $\BB$-module $A$, $\deg(A)$ denotes the  multiplicity relative to the maximal homogeneous ideal $\mathcal{M}$ of $\BB$,
 $\deg(A)= \deg( \gr_{\mathcal{M}}(A))$. In actual computations $\mathcal{M}$ can be replaced by a reduction. For instance, if $E$ is a graded $\RR$--module and $A=E\otimes_{\RR} \BB$, picking a reduction  $J$ for $\m$ gives the reduction
 $\mathcal{N} = (J, \TT_1, \ldots, \TT_{d+1})$ of $\mathcal{M}$. It will follow that $\deg(A) = \deg(E)$.

\begin{Theorem} \label{aCM1} If the algebra $\RR[It]$ is almost Cohen--Macaulay, in the canonical sequence
\[ 0 \rar T \lar \SS \lar \RR[It] \rar 0\]
\begin{itemize}
\item[{\rm (a)}] $\reg(\RR[It]) = \red_J(I) + 1$.

\item[{\rm (b)}] $\nu(T) \leq \deg(\SS) - \deg(\RR[It])$.

\end{itemize}
\end{Theorem}

\begin{proof} (a) follows from Corollary~\ref{Sallyrel}. As for (b), since $T$ is a Cohen--Macaulay module, $\nu(T) \leq \deg(T)$.
\end{proof}

The goal is
 to find $\deg(T)$, $\deg(\RR[It])$ and $\deg(\SS)$ in terms of more direct metrics of $I$. This will be answered in Theorem~\ref{degSymi}.

\bigskip

\subsubsection*{Cohen--Macaulayness of the Sally module}

Fortunately there is a simple criterion to test whether $\RR[It]$ is an aCM algebra: It is so
if and only if  it satisfies the Huckaba Test:
\[ \rme_1(I) = \sum_{j\geq 1}\length(I^j/JI^{j-1}). \]
Needless to say, this is exceedingly effective if you already know
$e_1(I)$,  in particular there is no need to determine the equations of $\RR[It]$ for the
purpose.

\medskip

Let $\RR$ be a Noetherian ring, $I$ an ideal and $J$ a reduction of $I$. The Sally module of
$I$ relative to $J$, $S_J(I)$, is defined by the exact sequence of $\RR[Jt]$--modules
\[ 0 \rar I \RR[Jt] \lar I \RR[It] \lar S_J(I) = \bigoplus_{j\geq 2} I^j/IJ^{j-1} \rar 0.\]
  The definition applies more broadly to other filtrations. We refer the reader to \cite[p. 101]{icbook} for
  a discussion. Of course this module depends on the chosen reduction
 $J$, but its Hilbert function and its depth are independent of $J$.
There are extensions of this construction to more general reductions--and we employ one below.

 \medskip

If $\RR$ is a Cohen--Macaulay local ring and $I$ is $\m$--primary with a minimal reduction, $S_J(I)$ plays
a role in mediating among properties of $\RR[It]$.

 \begin{Proposition} \label{Sallyelem} Suppose $\RR$ is a Cohen--Macaulay local ring of dimension $d$. Then
 \begin{enumerate}
 \item[{\rm (a)}] If $S_J(I) = 0$ then $\gr_I(\RR)$ is Cohen-Macaulay.
 
 \item[{\rm (b)}] If $S_J(I)\neq 0$ then $\dim S_J(I) = d$.
 \end{enumerate}
 \end{Proposition}

Some of the key properties of the Sally module  are in
display in the next result (\cite[Theorem 3.1]{Huc96}). It converts the property of $\RR[It]$ being almost Cohen--Macaulay into the property of $S_J(I)$ being Cohen--Macaulay.

\begin{Theorem}[Huckaba Theorem] \label{Huckaba}\index{Huckaba test  of Cohen--Macaulayness}
 Let $(\RR,\m)$ be a Cohen--Macaulay local ring of
dimension $d \geq 1$ and $J$ a
parameter ideal. Let $\mathcal{A}=\{I_n, n\geq 0\}$ be an filtration
of $\m$-primary ideals such that $J\subset I_1$ and $\BB=\RR[I_nt^n,
n\geq 0]$ is $\AA=\RR[Jt]$-finite.
Define the Sally module $S_{\BB/\AA}$
of $\BB$ relative to $\AA$ by the exact sequence
\[ 0 \rar I_1 \AA \lar I_1\BB \lar S_{\BB/\AA}\rar 0.\]
Suppose $S_{\BB/\AA}\neq 0$.
Then
\begin{enumerate}
\item[{\rm (a)}]  $e_0(S_{\BB/\AA})=e_1(\BB)-\lambda(I_1/J) \leq \sum_{j\geq 2}\lambda(I_j/JI_{j-1})$.
\item[{\rm (b)}]  The following conditions are equivalent:
\begin{itemize}
\item[{\rm (i)}] $S_{\BB/\AA}$ is Cohen-Macaulay;

\item[{\rm (ii)}] $\depth \gr_{\mathcal{A}}(\RR)\geq d-1$;

\item[{\rm (iii)}] $e_1(\BB)= \sum_{j \geq 1}\lambda(I_j/JI_{j-1})$;

\item[{\rm (iv)}] If $I_j = (I_1)^j$, $\BB$ is almost Cohen--Macaulay.

\end{itemize}

\end{enumerate}

\end{Theorem}

\begin{proof}  If $J=(\xx)=(x_1,\ldots, x_d)$,
$S_{\BB/\AA}$ is a finite module over the ring $\RR[\TT_1, \ldots,
\TT_d]$, $\TT_i \rar x_it$. Note that
\[\lambda( S_{\BB/\AA}/\xx S_{\BB/\AA})=
\sum_{j\geq 2}\lambda(I_j/JI_{j-1}),
\]
which shows the first assertion.

\medskip

For the equivalencies, first note that equality means that the first
Euler characteristic $\chi_1(\xx;S_{\BB/\AA})$ vanishes, which
by Serre's theorem (\cite[Theorem 4.7.10]{BH}) says that $S_{\BB/\AA}$ is Cohen--Macaulay.
The final assertion comes from the formula for the multiplicity of
$S_{\BB/\AA}$ in terms of $e_1(\BB)$ (\cite[Theorem 2.5]{icbook}).
\end{proof}

%%---------------------
% \begin{Question}{\rm Denote by $f_i$ the terms of the summation. We noticed examples were $e_1(I) +1 = \sum_{j\geq 1} f_i$: Does this imply that $\depth \gr_I(\RR) \geq d-2$?
%% }\end{Question}
%%-----------------------

\subsubsection*{Castelnuovo regularity}

The Sally module encodes also  information about the Castelnuovo regularity $\mbox{\rm reg}(\RR[It])$ of
the Rees algebra. The following Proposition and its Corollary are extracted from the literature (\cite{Huck87}, \cite{Trung98}), or proved directly by adding the exact sequence that defines $S_J(I)$ (note that $I\RR[Jt]$ is a maximal Cohen--Macaulay module) to the canonical  sequences relating $\RR[It])$ to $\gr_I(\RR)$ and $\RR$ via $I\RR[It]$ (see
\cite[Section 3]{Trung98}).

\begin{Proposition} \label{Sallyreg} Let $\RR$ be a Cohen--Macaulay local ring, $I$ an $\m$--primary ideal and $J$ a minimal reduction.
Then
\[ \mbox{\rm reg}(\RR[It]) = \mbox{\rm reg}(S_J(I)). \]
In particular
\[ \mbox{\rm reltype}(I) \leq  \mbox{\rm reg}(S_J(I)).\]
\end{Proposition}

\begin{Corollary} \label{Sallyrel}
If $I$ is an almost complete intersection and $\RR[It]$ is almost Cohen--Macaulay, then
\[ \mbox{\rm reltype}(I) = \red_J(I) + 1.\]
\end{Corollary}

\subsubsection*{The Sally fiber of an ideal}\index{Sally fiber of an ideal}

To help analyze the problem, we single out an extra structure.
Let $(\RR, \m)$ be a Cohen--Macaulay local ring of dimension $d>0$, $I$ an $\m$--primary ideal and $J$
one of its minimal reductions.

\begin{Definition}[Sally fiber] The Sally fiber of $I$ is the graded module\[ F(I) = \bigoplus_{j\geq 1} I^j/JI^{j-1}. \]
\end{Definition}

$F(I)$ is an Artinian $\RR[Jt]$--module whose last non-vanishing component is $I^r/JI^r$, $r=\red_J(I)$. The
equality $e_1(I) = \length(F(I))$ is the condition for the almost Cohen--Macaulayness of $\RR[It]$.
We note that $F(I)$ is the fiber of $S_J(I)$ extended by the term $I/J$.
To
obtain additional control over $F(I)$ we are going to endow it with additional structures in
cases of interest.

\bigskip

Suppose $\RR$ is a Gorenstein local ring, $I=(J,a)$. The modules $F_j = I^j/JI^{j-1}$ are cyclic modules
over the Artinian Gorenstein ring $ \AA = \RR/J:a$. We turn $F(I)$ into a graded module over the  polynomial ring $\AA[s]$ by defining
\[ a^j \in F_j \mapsto  s\cdot a^j = a^{j+1} \in F_{j+1}.\]
This is clearly well-defined and has $s^r\cdot F(I)=0$. Several of the properties of the $F_n$'s arise from this
representation, for instance the length of $F_j$ are non-increasing.
Thus $F(I)$ is a graded module over the Artinian Gorenstein
ring $\BB= \AA[s, s^r=0]$.

\medskip

\begin{Remark}\label{Fvasneweq}{\rm
The variation of the values of  the $F_j$ is connected to the degrees of the generators of $\LL$.  For convenience we set $I=(J,a)$ and $\BB = \RR[u, \TT_1, \ldots, \TT_d]$, with $u$ corresponding to $a$. For example:

\medskip

\begin{itemize}

\item[{\rm (a)}] Suppose that for some $s$, $\ff_s = \lambda(F_s) = 1$. This means that we have $d$ equations of the form
\[\hh_i= x_i u^s + \g2_i \in \LL_s\]
where  $\g2_i\in (\TT_1, \ldots, \TT_{d})\BB_{s-1}$. Eliminating the $x_i$, we derive a nonvanishing monic equation in $\LL$ of degree $d\cdot s$. Thus $\red_J(I) \leq ds -1$.

\medskip

\item[{\rm (b)}] A more delicate observation,  is that whenever $\ff_s > \ff_{s+1}$ then
there are {\bf fresh} equations in  $\LL_{s+1}$. Let us explain why this happens:
$\ff_s = \lambda(JI^{s-1}: I^s)$, that is the ideal $L_s$ contains elements of the form
\[ c\cdot u^s + \mathbf{g}, \quad c\in JI^{s-1}:I^s,  \quad  \mathbf{g} \in (\TT_1, \ldots, \TT_d)\BB_{s-1}.\]
Since $\ff_{s+1} < \ff_s$, $JI^s: I^{s+1}$ contains properly $JI^{s-1}:I^s$, which means that we must have elements in
$L_{s+1}$
\[ d\cdot u^{s+1} + \mathbf{g}, \quad\]
with $d \notin JI^{s-1}: I^s$ and $\mathbf{g} \in (\TT_1, \ldots, \TT_d)\BB_{s}$. Such elements cannot belong to
$L_s\cdot \BB_1$, so they are fresh generators.

The converse also holds.

\end{itemize}

}\end{Remark}

\subsubsection*{A toolbox}

We first give a simplified version of \cite[Proposition 2.2]{CPV1}.
Suppose $\RR$ is a Gorenstein local ring of dimension $d$. Consider the two
exact sequences.

\[0 \rar J/JI= (\RR/I)^d \lar \RR/JI \lar \RR/J\rar 0\]
 and the [old] syzygetic sequence
\[ 0 \rar \delta(I) \lar H_1(I) \lar (\RR/I)^{d+1} \lar I/I^2 \rar 0.
\]
The first
gives
\[\length(\RR/JI)= d\cdot \length(\RR/I) +\length(\RR/J),\]
 the other
\[\length(\RR/I^2) =(d+2)\length(\RR/I)-\length(\H_1(I)) + \length
(\delta(I)).\]
Thus
\[\length(I^2/JI)= \length(I/J) -\length(\delta(I))\]
since $\H_1(I)$ is
the canonical module of $\RR/I$. Taking into account the syzygetic
formula in \cite{syl2} we finally have:

\begin{Proposition} \label{F2} Let $(\RR, \m)$ be a
Gorenstein local ring of dimension $d>0$, $J= (a_1, \ldots, a_d)$ a
parameter ideal and $I=(J,a)$ and $a\in \m$. Then
\begin{eqnarray*} \length(I^2/JI) & = &\length(I/J)
-\length(\RR/I_1(\phi))\\ &=& \length(\RR/J:a)  -\length(\RR/I_1(\phi))\\
& = & \length(\RR/J:a) - \length(\Hom(\RR/I_1(\phi), \RR/J:a) \\
& = & \length(\RR/J:a) - \length((J:a): I_1(\phi))/J:a) \\
& = & \length (\RR/(J:a):I_1(\phi)).
\end{eqnarray*}

\end{Proposition}

Note that in dualizing $\RR/I_1(\phi)$ we made use of the fact that $\RR/J:a$ is a Gorenstein ring.

\begin{Corollary} \label{F2Cor} $I^2= JI$ if and only if $J:a= I_1(\phi)$. In this
case, if $d>1$ the algebra $\RR[It]$ is Cohen--Macaulay.
\end{Corollary}

\begin{Corollary}\label{F3Cor} If $\RR[It]$ is an aCM algebra and $\red_J(I) = 2$, then
$e_1(I) = 2\cdot \lambda(I/J) - \lambda(\RR/I_1(\phi))$.
\end{Corollary}

\begin{Remark}{\rm We could enhance these observations considerably if formulas for $\lambda(JI^2:I^3)$ were to be developed. More precisely, how do the syzygies of $I$ affect $JI^2:I^3$?
}
\end{Remark}

\subsection{Multiplicities and number of relations}

To benefit from Theorem~\ref{aCM1}, we need to have effective formulas for $\deg(\SS)$ and $\deg(\RR[It])$.
We are going to develop them now.

\begin{Proposition}\label{multirees} Let $\RR=k[x_1, \ldots, x_d]$ and $I$ an almost complete intersection as above,
$I=(f_1, \ldots, f_d, f_{d+1})=(J, f_{d+1})$ generated by
forms
 of degree $n$. Then
$ \deg(\RR[It])  = \sum_{j=0}^{d-1} n^j.$
\end{Proposition}

\begin{proof}  After an elementary observation, we make use of one of the beautiful multiplicity formulas of \cite{HTU}.
Set $A=\RR[It]$, $A_0 = \RR[Jt]$, $\mathcal{M} = (\m, It)A$ and $\mathcal{M}_0 = (\m, Jt)A_0$.
Then
\[ \deg(\gr_{\mathcal{M}_0}(A_0))= \deg(\gr_{\mathcal{M}_0}(A))= \deg(\gr_{\mathcal{M}}(A)), \]
the first equality because $A_0 \rar A$ is a finite rational extension, the second is because $(\m, Jt)A$ is a reduction of $(\m, It)A$. Now we use \cite[Corollary 1.5]{HTU} that gives $\deg(A_0)$.
\end{proof}

\subsubsection*{The multiplicity of the symmetric algebra}\index{multiplicity of symmetric algebras}

We shall now prove one of our main results, a formula for $\deg S(I)$ for ideals generated by forms of the same degree.
Let $\RR=k[x_1, \ldots, x_d]$, $I = (\ff) = (f_1, \ldots, f_d, f_{d+1})$ an almost complete intersection generated by forms of degree $n$. At some point we assume, harmlessdly,  that $J = (f_1, \ldots,  f_d)$ is a complete intersection.
There will be a slight change of notation in the rest of this section. We set $\BB = \RR[\TT_1, \ldots, \TT_{d+1}]$ and
$\SS = \Symi(I)$.

\begin{Theorem}[{\bf Degree Formula}] \label{degSymi} $\deg \SS = \sum_{j=0}^{d} n^j - \lambda(\RR/I)$.
\end{Theorem}

\begin{proof}
Let $\mathbb{K}(\ff) = \bigwedge^{d+1} \RR^{d+1}(-n)$ be the Koszul complex associated to $\ff$,
\[ 0 \rar {K}_{d+1} \rar K_{d} \rar \cdots \rar K_2 \rar K_1 \rar K_0 \rar 0,\]
and consider the associated $\mathcal{Z}$--complex
\[ 0 \rar Z_d\otimes \BB(-d) \rar Z_{d-1} \otimes \BB(-d+1)
\rar \cdots \rar Z_2 \otimes \BB(-2) \rar Z_1 \otimes \BB(-1) \stackrel{\psi}{\rar} \BB \rar 0.\]

This complex is acyclic with $\H_0(\mathcal{Z}) = \SS = \Symi(I)$.
Now we introduce another complex obtained by replacing $Z_1 \otimes \BB(-1) $ by $B_1 \otimes \BB(-1)$, where
$B_1$ is the module of $1$--boundaries of $\mathbb{K}(\ff)$,
 followed by
the restriction of $\psi$ to $B_1 \otimes \BB(-1)$.

 \medskip

 This defines another acyclic complex, $\mathcal{Z}^*$, actually the $\mathcal{B}$--complex of $\ff$,  and we set $\H_0(\mathcal{Z}^*) = \SS^*$. The relationship
 between $\SS$ and $\SS^*$ is given in the following observation:
 
 \end{proof}

 \begin{Lemma} $\deg \SS^* = \deg \SS + \lambda(\RR/I)$.
 \end{Lemma}

 \begin{proof} Consider the natural mapping between $\mathcal{Z} $ and $\mathcal{Z}^*$:

 \[
 \diagram
 0 \rto  & Z_{d}\otimes \BB(-d) \rto \dto_{\phi_d} & \cdots \rto  & Z_2\otimes \BB(-2) \rto \dto_{\phi_2} & B_1 \otimes \BB(-1) \rto \dto  & \BB \rto \dto & \SS^* \rto \dto & 0 \\
 0 \rto & Z_{d}\otimes \BB(-d) \rto & \cdots \rto & Z_2\otimes \BB(-2) \rto & Z_1 \otimes \BB(-1) \rto & \BB \rto &
 \SS \rto & 0.
 \enddiagram
   \]
 The maps $\phi_2, \ldots, \phi_2 $ are isomorphisms while the other maps are defined above. They induce the short exact sequence of modules of dimension $d+1$,
 \[ 0 \rar (Z_1/B_1)\otimes \BB(-1) \lar \SS^* \lar \SS \rar 0.\]

Note that $Z_1/B_1= \H_1(\mathbb{K}(\ff))$ is the canonical module of $\RR/I$, and therefore has the same length
as
$\RR/I$. Finally, by the additivity formula for the multiplicities (\cite[Lemma 13.2]{Eisenbudbook}),
\[ \deg \SS^* = \deg \SS + \lambda(Z_1/B_1),\]
as desired. 
\end{proof}

We are now give our main calculation of multiplicities.

\begin{Lemma} $\deg \SS^* = \sum_{j=0}^{d} n^j.$

\end{Lemma}

\begin{proof} We note that the $\mathcal{Z}^*$--complex is homogeneous for the total degree [as required for the computation of multiplicities] provided the $Z_i$'s and $B_1$ have the same degree.
We can conveniently write $B_i$ for $Z_i$, $i\geq 2$.
This is clearly the case since they are
 generated in degree $n$. This is not the case for $Z_1$. However when  $\ff$ is a regular sequence, all the $Z_i$ are equigenerated, an observation we shall make  use of below.

\medskip

Since the modules of $\mathcal{Z}^*$ are homogeneous
we have that the Hilbert series of $\SS^*$ is given as
\[ H_{\SS^*}(\ttt) = {\frac{\sum_{i=0}^d (-1)^{i} h_{B_i}(\ttt) \ttt^i}{(1-\ttt)^{2d+1}}} = {\frac{h(\ttt)}{(1-\ttt)^{2d+1}}},
\]
where $h_{B_i}(\ttt)$ are the $h$--polynomials of the $B_i$. More precisely, each of the terms of $\mathcal{Z}^*$ is a $\BB$--module of the form $A\otimes \BB(-r)$ where $A$ is generated in a same degree. Such modules are isomorphic to their associated  graded modules.

\medskip

The multiplicity of $\SS^*$ is given by the standard formula
\[\deg \SS^* = (-1)^d   {\frac{h^{(d)}(1)}{ d!}}.\]
We now indicate how the $h_{B_i}(\ttt)$ are  assembled. Let us illustrate
the case when $d =4$ and $i=1$. $B_1$ has  a free resolution of the strand of the Koszul complex
\[ 0 \rar \RR(-3n) \rar \RR^{5}(-2n) \lar \RR^{10}(-n) \lar \RR^{10} \lar B_1 \rar 0,\]
so that
\[ h_{B_1}(\ttt) = 10  - 10\ttt^{n} + 5\ttt^{2n}-\ttt^{3n},\]
and similarly for all $B_i$.

\medskip

We are now ready to make our key observation. Consider a complete intersection $P$
generated by $d+1$ forms of degree $n$ in a polynomial ring of dimension $d+1$
and set $\SS^{**} = \Symi(P)$. The corresponding approximation complex now has $B_1=Z_1$.
 The approach above would for the new $Z_i$ give the
same $h$--polynomials of the $B_i$ in the case of an almost complete intersection (but in dimension $d$).
This means that the Hilbert series of $\SS^{**}$ is given by
\[ H_{\SS^{**}}(\ttt) = {\frac{h(\ttt)}{(1-\ttt)^{2d+2}}}.\]
 It follows
 that
$\deg \SS^*$ can be computed as the degree of the symmetric algebra generated by a regular sequence of $d+1$ forms
of degree $n$, a  result that is given in \cite{HTU}. Thus,
\[ \deg \SS^* = \deg \SS^{**} = \sum_{j=0}^d n^j,\]
and the calculation of $\deg \SS$ is complete. 
\end{proof}

We will now write Theorem~\ref{degSymi} in a more convenient formulation for applications.

\begin{Theorem}\label{degSymibis} Let $\RR =k[x_1, \ldots, x_d]$ and  $I = (f_1, \ldots, f_d, f_{d+1})$
is an ideal  of forms of degree $n$. If $J = (f_1, \ldots, f_d)$ is a complete intersection, then
\[ \deg \SS = \sum_{j=0}^{d-1} n^j + \lambda(\RR/J:I).\]
\end{Theorem}

\begin{proof} The degree formula gives
\begin{eqnarray*}
 \deg \SS & = & \sum_{j=0}^{d-1}n^j + [n^d - \lambda(R/I)] = \sum_{j=0}^{d-1}n^j + [\lambda(\RR/J) - \lambda(\RR/I)]\\
 &=& \sum_{j=0}^{d-1} n^j +\lambda(I/J) = \sum_{j=0}^{d-1}n^j + \lambda(\RR/J:I).
 \end{eqnarray*}
 \end{proof}

\begin{Corollary} \label{degT} Let $I=(J,a)$ be an ideal of finite colength as above. Then the module of
linear relations satisfies
 $\deg T = \length (I/J)$. In particular if $\RR[It]$ is almost Cohen--Macaulay, $T$ can be generated by $\lambda(I/J)$ elements.
\end{Corollary}

\begin{proof} From the sequence of modules of the same dimension
\[ 0 \rar T \lar \SS \lar \RR[It]\rar 0\]
we have
\[ \deg T = \deg \SS - \deg \RR[It] = \lambda(I/J).\]
\end{proof}

The last assertion of this Corollary   can also  be obtained
from \cite[Theorem 4.1]{MPV12}.

\bigskip

\subsubsection*{The Cohen--Macaulay type of the module of nonlinear relations}

We recall the terminology of Cohen--Macaulay type of a module.
 Set $\BB=\RR[\TT_1, \ldots, \TT_{d+1}]$. If $E$ is a
  finitely generated $\BB$--module of codimension $r$, we say that $\Ext_{\BB}^r(E, \BB)$ is its canonical module. It is the first non vanishing $\Ext_{\BB}^i(E, \BB)$ module denoted by $\omega_E$. The minimal number of the generators of $\omega_{E}$ is called   the {\em Cohen--Macaulay type} of $E$ and is denoted by  $\mbox{\rm type}(E)$.
When $E$ is graded and Cohen--Macaulay, it gives the last Betti number of a projective resolution of $E$. It can be expressed in different ways, for example for
the module of linear relations $\omega_T=\Ext_{\BB}^d(T, \BB) = \Hom_{\SS}(T,\omega_{\SS})$.

\begin{Proposition}\label{typeofT} Let $\RR$ be a Gorenstein local  ring of dimension $d \geq 2$
 and $I=(J,a)$ an ideal of finite colength as above. If $\RR[It]$ is an aCM algebra and $\omega_{R[It]}$ is Cohen--Macaulay, then the type of the  module $T$ of nonlinear relations satisfies
\[ \mbox{\rm type}(T) \leq \mbox{\rm type }(S_J(I)) + d-1,\]
where $S_J(I)$ is the Sally module.
\end{Proposition}

\begin{proof} We set $\Rees = \RR[It]$ and $\Rees_0 = \RR[Jt]$.
First apply $\Hom_{\BB}(\cdot, \BB) $ to the basic presentation
\[ 0 \rar T \lar \SS \lar \Rees \rar 0,\]
to obtain the cohomology sequence
\begin{eqnarray}\label{type1}
 0 \rar \omega_{\Rees} \lar \omega_{\SS} \lar \omega_T \lar \Ext_{\BB}^{d+1}(\Rees, \BB) \rar 0.
 \end{eqnarray}

Now apply the same functor to the exact sequence of $\BB$--modules
\[ 0 \rar I\cdot \Rees[-1]  \lar \Rees \lar \RR \rar 0\]
to obtain the exact sequence
\[ 0 \rar \omega_{\Rees} \stackrel{\theta}{\lar} \omega_{I\Rees[-1]} \lar \Ext_{\BB}^{d+1}(\RR, \BB) = \RR
  \lar \Ext_{\BB}^{d+1}(\Rees, \BB) \lar
\Ext_{\BB}^{d+1}(I\Rees[-1], \BB) \rar 0.\]
%By Corollary~\ref{canofRees}(b),
Since $\omega_{\Rees}$ is Cohen--Macaulay and $\dim \RR \geq 2$,
 the cokernel of $\theta$ is either $\RR$ or an $\m$-primary ideal that satisfies the condition $S_2$ of Serre. The only choice is $\coker(\theta)=\RR$. Therefore
 \[\Ext_{\BB}^{d+1}(\Rees, \BB)\simeq \Ext_{\BB}^{d+1}(I\Rees[-1], \BB).\]

 Now we approach the module $ \Ext_{\BB}^{d+1}(I\Rees, \BB)$ from a different direction. We note that $\Rees$---but not $\SS$ and $T$---is also a finitely generated $\BB_0=\RR[\TT_1, \ldots, \TT_d]$--module as it is annihilated by a monic polynomial $\ff$ in $\TT_{d+1}$ with coefficients in $\BB_0$. By Rees Theorem
  we have  that for all $i$, $\Ext_{\BB}^i(\Rees, \BB) =\Ext_{\BB/(\ff)}^{i-1}(\Rees, \BB/(\ff))$, and a similar observation applies to $I\cdot \Rees$.

  \medskip

  Next consider the finite, flat morphism $\BB_0 \rar \BB/(\ff)$. For any $\BB/(\ff)$--module $E$ with a projective
  resolution $\mathbb{P}$, we have that $\mathbb{P}$ is a projective $\BB_0$--resolution of $E$. This
  means that the isomorphism of complexes
  \[ \Hom_{\BB_0}(\mathbb{P}, \BB_0) \simeq \Hom_{\BB/(\ff)}(\mathbb{P}, \Hom_{\BB_0}(\BB/(\ff), \BB_0))=
 \Hom_{\BB/(\ff)}(\mathbb{P}, \BB/(\ff))
  \]
  gives  isomorphisms for all $i$
  \[ \Ext_{\BB_0}^i(E, \BB_0) \simeq  \Ext_{\BB/(\ff)}^i(E, \BB/(\ff)). \]
 Thus
 \[ \Ext_{\BB}^i(\Rees, \BB) \simeq \Ext_{\BB_0}^{i-1}(\Rees, \BB_0).\]
In particular, $\omega_{\Rees} = \Ext_{\BB_0}^{d-1}(\Rees, \BB_0)$.

%%----
%   In turn this last module can also be expressed using $\omega_{\SS_0}$, $\SS_0=\Symi(J)$ and
% particular  $\Ext_{\SS}^1(\Rees, \omega_{\SS}) = \Ext_{\SS_0}^1(\Rees, \omega_{\SS_0})$ that uses the fact that
% $\BB_0 \rar \BB/(\ff)$ is a finite morphism of Cohen-Macaulay algebras.
%%----

\medskip

Finally apply $\Hom_{\BB_0}(\cdot, \BB_0)$ to the exact sequence of $\BB_0$--modules
 and examine its   cohomology sequence.

\[ 0 \rar I\cdot \Rees_0 \lar I\cdot \Rees \lar S_J(I) \rar 0\]
is then
\[ 0 \rar \omega_{I\Rees} \lar \omega_{I\Rees_0} \lar \omega_{S_J(I)}\lar \Ext_{\BB_0}^d(\Rees, \BB_0)= \Ext_{\BB}^{d+1}(\Rees, \BB) \rar 0.\]
%%-------------
% From our discussion above, we have
% \[ \Ext_{\BB_0}(I\cdot \Rees [-1], \BB_0)\simeq \Ext_{\BB_0}^d(\Rees, \BB_0) \simeq
% \Ext_{\BB/(\ff)}^d(\Rees, \BB/(\ff)) \simeq \Ext_{\BB}^{d+1}(\Rees, \BB). \]
%%-----
 Taking this into
(\ref{type1}) and the  that $\mbox{\rm type}(\SS)=d-1$ gives the desired estimate.
 \end{proof}

\begin{Remark}{\rm
A class of ideals with $\omega_{\Rees}$ Cohen--Macaulay is discussed in Corollary~\ref{canofRees}(b).
}\end{Remark}

%\subsection{Distinguished aCM algebras}

%This section treats several classes of Rees algebras which are almost Cohen--Macaulay.

\subsection{Equi-homogeneous acis}\index{equi-homogeneous ACIs}

We shall now treat an important class of extremal Rees algebras.
Let $\RR = k[x_1, \ldots, x_d]$ and let $I=(a_1, \ldots, a_d, a_{d+1})$ be an ideal
of finite colength, that is,  $\m$--primary. We further assume that the first
 $d$ generators form a regular sequence and $a_{d+1} \notin  (a_1, \ldots, a_d)$.
    If $\deg a_i=n$, the integral closure
 of $J=(a_1, \ldots, a_d)$ is the ideal $\m^n$, in particular $J$ is a minimal reduction of $I$. The integer
 $\edeg(I) = \red_J(I) +1$ is called the {\em elimination degree} of $I$.
 The study of the equations of $I$, that is, of $\RR[It]$, depends on a comparison between
 the metrics  of $\RR[It]$ to
 those of $\RR[\m^n t]$, which are well known.

\begin{Proposition} {\rm (\cite{syl2})} \label{birideal}
The following conditions are equivalent:
\begin{itemize}
\item[{\rm (a)}] $\Phi$ is a birational mapping, that is the natural embedding $\mathcal{F}(I) \hookrightarrow \mathcal{F}(\m^n)$ is an isomorphism of quotient fields;

\item[{\rm (b)}] $\red_J(I) = n^{d-1}-1$;

\item[{\rm (c)}] $e_1(I) = {\frac{d-1}{2}}(n^d - n^{d-1})$;

\item[{\rm (d)}] $\RR[It]$ satisfies the condition $R_1$ of Serre.
\end{itemize}
\end{Proposition}

For lack of a standardized terminology, we say that these ideals are {\em birational}\label{birational ideal}.

\begin{Corollary}\label{canofRees} For an ideal $I$ as above, the following hold:
\begin{itemize}

\item[{\rm (a)}] The algebra $\RR[It]$ is not Cohen--Macaulay except when $I = (x_1, x_2)^2$.

\item[{\rm (b)}] The canonical module of $\RR[It]$ is Cohen--Macaulay.

\end{itemize}

\end{Corollary}

\begin{proof} (i) follows from the condition of Goto--Shimoda (\cite{GS82}) that the reduction number of a Cohen--Macaulay Rees algebra $\RR[It]$ must satisfy $\red_J(I) \leq \dim \RR -1$, which in the case
$n^{d-1}- 1 \leq d-1$
is only met if $d =n =2$

\medskip

\noindent (ii) The embedding $\RR[It] \hookrightarrow \RR[\m^n t]$ being an isomorphism in codimension one, their canonical modules are isomorphic. The canonical module of a Veronese subring such as
$\RR[\m^n t]$ is well-known (see \cite[p. 187]{HV85}, \cite{HSV87}; see also \cite[Proposition 2.2]{BR}).
\end{proof}

\subsubsection*{Binary ideals}

These are the ideals of $\RR=k[x,y]$ generated by $3$ forms of degree $n$. Many of their Rees algebras are almost Cohen--Macaulay. We will showcase the technology of the Sally module in treating a much studied class of ideals. First
we discuss a simple case (see also \cite{syl1})).

\begin{Proposition}\label{22} Let $\phi$ be a $3\times 2$ matrix of quadratic forms in $\RR$ and $I$ the ideal
given by its $2\times 2$ minors. Then $\RR[It]$ is almost Cohen--Macaulay.
\end{Proposition}

\begin{proof} These ideals have reduction number $1$ or $3$. In the first case all of its Sally modules vanish and
$\RR[It]$ is Cohen--Macaulay.

\medskip

In the other case $I$ is a birational ideal and $e_1(I) = {4\choose 2} = 6$. A simple calculation shows that
$\lambda(\RR/I) = 12$, so that $\lambda(I/J) = 16-12 = 4$. To apply Theorem~\ref{Huckaba}, we need to verify
the equation
\begin{eqnarray}\label{Sally22}
 f_1 + f_2 + f_3 = 6.
 \end{eqnarray}
We already have $f_1=4$. To calculate $f_2$ we need to take $\lambda(R/I_1(\phi))$ in Corollary~\ref{F2}. $I_1(\phi)$ is an ideal generated by $2$ generators or $I_1(\phi) = (x,y)^2$. But in the first case the Sylvester resultant
of
 the linear equations of $\RR[It]$ would be a quadratic polynomial, that is $I$ would have reduction number $1$, which
 would contradict the assumption. Thus by Corollary~\ref{F2}, $f_2 = 4-\lambda(\RR/I_1(\phi)) = 1$. Since
 $f_2\geq f_3>0$ we have $f_3=1$ and the equation (\ref{Sally22}) is satisfied. 
 \end{proof}

We have examined higher degrees examples of birational ideals of this type which are/are not almost Cohen--Macaulay.
Quite a lot is known about the following ideals.
$\RR = k[x,y]$ and $I$ is a codimension $2$ ideal given by that $2\times 2$ minors of a $3\times 2$ matrix
with homogeneous entries of degrees $1$ and $n-1$.

\begin{Theorem} \label{2birideal}
If $I_1(\phi) = (x, y)$ then:
\begin{enumerate}

\item[{\rm (a)}] $\deg \mathcal{F}(I) = n$, that is $I$ is birational.

\medskip

\item[{\rm (b)}]  $\RR[It]$ is almost Cohen--Macaulay.

\medskip

\item[{\rm (c)}]  The equations of $\LL$ are given by a straightforward algorithm.

\end{enumerate}
\end{Theorem}

\begin{proof} The proof of (i) is in \cite{CHW},  and in other sources (\cite{CdA}, \cite{KPU}; see also
\cite[Theorem 2.2]{DHS} for a broader statement in any characteristic and \cite[Theorem 4.1]{Simis04} in
characteristic zero),
 and of  (ii)  in \cite[Theorem 4.4]{KPU}, while (iii)  was conjecturally  given in  \cite[Conjecture 4.8]{syl1} and proved in \cite{CHW}. We give a combinatorial proof of (ii) below (Proposition~\ref{aCMofbin}).
\end{proof}

We note that
$\deg(\SS) = 2n$, since $\SS$ is a complete intersection defined by two forms of [total] degrees $2$ and $n$, while $\RR[Jt]$ is defined by one equation of degree $n+1$. Thus $\nu(T) \leq 2n-(n+1)= n-1$, which is the number of generators given in the algorithm.

\bigskip

We point out a property of the module $T$. We recall that an $\AA$--module is an {\em Ulrich} module if it is a maximal Cohen--Macaulay module with $\deg M=\nu(M)$ (\cite{HK}).

\begin{Corollary} $T$ is an Ulrich $\SS$--module.
\end{Corollary}

%Let $\RR= k[x,y]$ and $\phi$ a $3\times 2$ matrix with  columns of forms of degree $p$ and $q$, whose
%$2\times 2$ minors define an ideal of codimension $2$. We may assume that the minors sharing the
%top row define a reduction $J$ of $I= (J,a)$. If $\alpha$ and $\beta$ are the entries of the top row,
%$\AA = \RR/(\alpha, \beta)$, a ring whose representations can be very complex.
%In one case however

Considerable numerical information in the Theorem~\ref{2birideal} will follow from:

\begin{Proposition}\label{aCMofbin} If $\deg \alpha =1$ and $\deg \beta = n-1$, then $\length(F_j) = n -j$. In particular, $\RR[It]$ is almost Cohen--Macaulay.
\end{Proposition}

\begin{proof}
%Before we start on the proof proper, we note some properties of the structure of $F$(I).
%We may, by changing notation, assume that $\alpha = x$ and $\beta = y^{n-1}$. The ring $\AA= k[y]/(y^{n-1})$. This is
%an uniserial and
%\[ 0 \subset (y^{n-2}) \subset \cdots \subset (y) \subset \AA\]
%is its unique composition series. The corresponding series of quotient modules
%can be represented by the same series in reverse order
%\[ \AA \supset  (y) \supset \cdots \supset (y^{n-2}) \supset 0.  \]
%Each of the components of $F(I)$ can be identified to one of these
%\[ F_j =(y^{d_j}) s^j \AA.\]
%According to the definition of the action of $s$ on $F_j$ above, mapping the canonical generator of $F_j$ into the canonical generator of $F_{j+1}$.
Note that  the ideal is birational, $F_{n-1}\neq 0$. On the other hand, $\LL$ contains fresh generators in all degrees $j\leq n$. This means that for $f_j = \lambda(F_j)$,
\[ f_j> f_{j+1}>0, \quad j<n.\]
Since $f_1 = n-1$, the decreasing sequence of integers
\[ n-1 = f_1 > f_2 > \cdots > f_{n-2} > f_{n-1}> 0\]
implies that $f_j = n-j$.
Finally, applying Theorem~\ref{Huckaba} we have that $\RR[It]$ is an aCM algebra since $\sum_{j}f_j = e_1(I)= {n\choose 2}$. 
\end{proof}

\subsubsection*{Quadrics}\index{quaternary quadrics}

Here we explore sporadic classes of aCM algebras defined by   quadrics in
$k[x_1, x_2, x_3, x_4]$.

\medskip

  First we
 use Proposition~\ref{F2} to look at  other cases of quadrics. For $d=3$, $n=2$, $\edeg(I)=2$ or $4$. In the first case
$J:a= I_1(\phi)$. In addition $J:a\neq \m$ since the socle degree of
$\RR/J$ is $3$. Then $\RR[It]$ is Cohen--Macaulay. If $\edeg(I)=4$ we must have [and conversely!]
 $\length(\RR/J:a)=2$ and $I_1(\phi)=\m$. Then
$\RR[It]$ is almost Cohen--Macaulay.

\medskip

Next we treat almost complete intersections of finite colength generated by quadrics of $\RR= k[x_1, x_2, x_3, x_4]$. Sometimes we denote the variables by $x,y,\ldots$, or use these symbols to denote  [independent] linear forms. For notation we
use $J= (a_1, a_2, a_3, a_4)$, and $I = (J, a)$.

\bigskip

Our goal  is to address the following:
\begin{Question}{\rm  Let $I$ be an almost complete intersection generated by $5$ quadrics of $x_1, x_2, x_3, x_4$. If
  $I$ is a birational ideal, in which cases   is $\RR[It]$ is an almost Cohen--Macaulay algebra? In this case, what are the generators of the its module of nonlinear relations?
}\end{Question}

In order to make use of Theorem~\ref{Huckaba}, our main tools are Corollary~\ref{F2} and \cite[Theorem 2.2]{syl2}.
They make extensive use of the syzygies of $I$. The question forks into three cases, but our analysis is
complete in only one of them.

\bigskip

\subsubsection*{The Hilbert functions of quaternary quadrics}

We make a quick classification of the Hilbert functions of the ideals $I=(J,a)$. Since $I/J \simeq \RR/J:a$ and $J$ is a complete intersection, the problem is equivalent to determine the Hilbert functions of $\RR/J:a$, with $J:a$ a Gorenstein ideal.
The Hilbert function $H(\RR/I)$ of $\RR/I$ is $H(\RR/J)-H(\RR/J:a)$. We will need the Hilbert function of the corresponding canonical module in order to make use of \cite[Proposition 3.7]{syl2} giving information about $L_2/\BB_1L_1$.

\medskip

We shall refers to the  sequence   $(f_1, f_2, f_3, \ldots, )$, $f_i = \lambda(I^i/JI^{i-1})$, as the $\ff$--sequence of     $(I,J)$. We recall that these sequences are monotonic and that
 if $I$ is birational, $\sum_{i\geq  1} f_i=e_1(I)=12$.

\begin{Proposition}
Let  $\RR= k[x_1, x_2, x_3, x_4]$ and $I=(J,a)$ an almost complete intersection generated by $5$ quadrics, where $J$ is a complete intersection,.  Then
$L=\lambda(\RR/J:a)\leq 6$ and the possible Hilbert functions of $\RR/J:a)$ are:
\begin{eqnarray*}
L = 6 & : & (1,4,1), \quad (1,2,2,1)^{*}, \quad (1,1,1,1,1,1)^{*} \\
L = 5 & : & (1,3,1), \quad (1,1,1,1,1)^{*} \\
L = 4 & : & (1,2,1), \	\quad (1,1,1,1)^{*} \\
L = 3 & : & (1,1,1)^{*} \\
L = 2 & : & (1,1)^{**} \\
L = 1 & : & (1)^{**}
\end{eqnarray*}
If $I$ is a birational ideal, the corresponding Hilbert function is one of the unmarked sequences above.
\end{Proposition}

\begin{proof} Since $\lambda(\m^2/I)\geq 5$ and $\lambda(\m^2/J)=11$, $L=\lambda(\RR/J:a)\leq 6$. Because the Hilbert function of $R/(J:a)$ is symmetric and $L \leq 6$, the list includes all the viable Hilbert functions.

\medskip

Let us first rule out those marked with ${\mbox{\rm a}}^{*}$, while those marked with ${\mbox{\rm a}}^{**}$ cannot be birational.
In each of these $J:a$ contains at least $2$ linearly independent linear forms, which we denote by $x,y$, so that $J:a/(x,y)$ is a Gorenstein ideal of the regular ring $\RR/(x,y)$. It follow that $(J:a)/(x,y)$ is a complete intersection. In the case of $(1,2,2,1)$, $J:a = (x,y, \alpha, \beta)$, where $\alpha$ is a form of degree $2$ and $\beta$ a form of degree $3$, since $\lambda(\RR/J:a)=6$. Since $J\subset J:a$, all the generators of $J$ must be contained in $(x,y,\alpha)$, which is impossible by Krull theorem.
Those  strings with at least  $3$ $1$'s are also excluded since $J:a$ would
 have the form $(x,y,z,w^s)$, $s\geq 3$, and the argument above applies. The case $(1,1)$, $J:a = (x,y,z,w^2)$. This means that $I_1(\phi)=J:a$, or $J:a = \m$. In the first case, by Corollary~\ref{F2Cor}, $I^2 = JI$. In the second case, $I_1(\phi) = \m$. This will imply that
$ \lambda(I^2/JI) = 2 - 1 = 1$, and therefore $I$ will not be birational (need the summation to add to $12$).
\end{proof}

\subsubsection*{Hilbert function $(1,4,1)$}

If $R/J:a$ has Hilbert function $(1,4,1)$, $J:a \subset \m^2$ but we cannot have equality since $\m^2$ is not a Gorenstein ideal. We also have
 $I_1(\phi) \subset \m^2$. If they are not equal, $I_1(\phi) = J:a$, which by Corollary~\ref{F2Cor} would mean that $\red_J(I) = 1$.

\begin{Theorem} \label{141} Suppose $I$  that is birational and $I_1(\phi)\subset \m^2$.
Then
$\RR[It]$ is almost Cohen--Macaulay.
\end{Theorem}

\begin{proof}
The assumption $I_1(\phi)\subset \m^2$ means that the Hilbert function of $J:a$ is $(1,4,1)$ and vice-versa.
 Note also that by assumption  $\length(I^7/JI^6)\neq 0$. Since $\lambda(I/J) = \lambda(\RR/J:a) = 6$, it suffices to show that $\lambda(I^2/JI)=1$.
 From $\lambda(\m^2/I) = 5$, the module $\m^2/I$ is of length $5$ minimally generated by $5$ elements. Therefore $\m^3\subset I$, actually $\m^3=\m I$.

\medskip

There is an isomorphism  $\hh: \RR/J:a \simeq I/J$, $r\mapsto ra$. It moves the socle of
$\RR/J:a$ into the socle of $I/J$. If $a\notin J:a$, then $\m^2 = (J:a, a)$ and $a$ gives the
socle of $\RR/J:a$, thus it is mapped by $\hh$ into the socle of
$I/J$, that is $\m\cdot a^2\in J$. Thus,
$ \m\cdot a^2\in \m^3 J \subset JI$. On the other hand, if  $a\in J:a$, then, since $a^2\in J$, we have $a^2\in \m^2 J$ and
$\m\cdot a^2\in \m^3 J\subset JI$.

\medskip

An example is $J= (x^2, y^2, z^2, w^2)$, $a =  xy + xz + xw + yz$.
 \end{proof}

\subsubsection*{Hilbert function $(1,3,1)$}

Our discussion about this case is very spare.

\begin{itemize}

\item[$\bullet$]
For these Hilbert functions,
$ J:a = (x, P)$,
where $P$ is a Gorenstein ideal in a regular local ring of dimension $3$--and therefore
is given by the Pfaffians of a
skew-symmetric matrix, necessarily $5\times 5$.
 Since $J \subset J:a$, $L$ must contain forms of degree $2$. In addition,   $P$ is given by $5$  $2$--forms (and $(x,\m^2)/(x,P)$ is the socle of $\RR/J:a$).

\item[$\bullet$] If $I$ is birational,  $\RR[It]$ is almost Cohen--Macaulay if and only if $\lambda(I^2/JI) = 2$ and $\lambda(I^3/JI^2)=1$.   The first equality, by Proposition~\ref{F2}, requires $\lambda(\RR/I_1(\phi)) = 3$ which gives that $I_1(\phi)$ contains the socle of $J:a$ and another independent linear form. In all it  means that $I_1(\phi) = (x,y, (z,w)^2)$. On the other hand
 $\lambda(I^2/JI)=2$ means that $JI:I^2 = (x,y,z,w^2)$ (after more label changes).

\item[$\bullet$]  An   example is  $J=(x^2, y^2, z^2, w^2)$ with $a=  xy + yz + zw + wx + yw$. The ideal $I= (J,a)$ is birational.

 \end{itemize}

\bigskip

  \subsubsection*{Hilbert function $(1,2,1)$}

 We do not have the full analysis of this case either.

 \begin{itemize}

 \item[$\bullet$]  An example is  $J=(x^2, y^2, z^2, w^2)$ and $a = xy+yz+xw+zw$. The ideal $I=(J,a)$  is birational.
The expected $\ff$--sequence of such  ideals is $(4,3,1,1,1,1,1)$.

\item[$\bullet$]
If $I$ is birational then $I_1(\phi)=\m$. We know that $I_1(\phi)\neq J:a$, so $I_1(\phi) \supset \m^2$, that is
$I_1(\phi)= (x,y,\m^2)$,  $(x,y,z,\m^2)$, or $\m$. Let us exclude the first two cases.

$(x,y,\m^2)$: This leads to two equations
\begin{eqnarray*}
xa &=& xb + yc\\
ya &=& xd + ye,
\end{eqnarray*}
with $b,c,d,e\in J$. But this gives the equation $(a-b)(a-e) -dc=0$, and $\red_J(I)\leq 1$.

\medskip

  $(x,y,z,\m^2)$: Then the Hilbert function of $\RR/I_1(\phi)$ is $(1,1)$. According to \cite[Proposition 3.7]{syl2}, $L_2$ has a form of bidegree $(1,2)$, with coefficients in $I_1(\phi)$, that is, in
  $(x,y,z)$. This  gives $3$ forms with coefficients in this ideal, two in degree $1$, so by elimination we get a monic
 equation of degree $4$.

\end{itemize}

We summarize the main points of these observations into a normal form assertion.

\begin{Proposition}\label{121} Let $I$ be a birational ideal and the  Hilbert function of $\RR/J:a$ is
$(1,2,1)$. Then up to a change of variables
to $\{ x,y,z,w\}$, $I$ is a Northcott ideal, that is there is a $4\times 4$ matrix $\AA$,
\[ \AA = \left[ \begin{array}{ccc}
& \BB &  \\
 \hline
 & \mathbf{C} & \\
 \end{array} \right] . \]
 where $\BB$ is a $2\times 4$ matrix whose entries are linear forms and $\mathbf{C}$ is a matrix with scalar entries and
 $\mathbf{V} = [x,y, \alpha, \beta]$, where $\alpha, \beta$ are quadratic forms in $z,w$ such that
 \[ I = (\mathbf{V}\cdot \AA, \det \AA).\]
 \end{Proposition}

 \begin{proof} There are two independent linear forms in $J:a$ which we denote by $x,y$.
 We observe that $(J:a)/(x,y)$ is a Gorenstein ideal in a polynomial ring of
  dimension two, so it is a complete intersection: $J:a = (x,y, \alpha,\beta)$, with $\alpha$ and $\beta$ forms of degree $2$ (as $\lambda(\RR/J:a)=4$), from which we remove the terms in $x,y$, that is we may assume $\alpha, \beta\in (z,w)^2$.

  \medskip

  Since $J \subset J:a$, we have a matrix $\AA$,
 \[ J = [x,y,\alpha, \beta] \cdot \AA =\mathbf{V} \cdot \AA.\]
 By duality $I=J:(J:a)$, which by Northcott theorem (\cite{Northcott}) gives
 \[ I= (J, \det \AA).\]
 Note that $a$ gets, possibly, replaced by $\det \AA$.
  The statement about the degrees of the entries of $\AA$ is clear. 
  \end{proof}

 \begin{Example}{\rm
 Let
 \[ \AA = \left[ \begin{array}{rrrr}
 x+y & z + w & x-w & z \\
 z & y+w & x - z & y \\
 1 & 0 & 2 & 3 \\
 0 & 1 & 1 & 2 \\
 \end{array} \right], \quad {\mathbf v} = \left[ x, y,  z^2 + zw + w^2, z^2-w^2\right] .\]
This  ideal is birational but $\RR[It]$ is not aCM. This is unfortunate but opens the question of when such ideals are birational. The $\ff$--sequence here is $(4,3,3,1,1,1,1)$.

 }\end{Example}

 \bigskip

\subsubsection*{The degrees of $\LL$}

We examine how the Hilbert function of $\RR/J:a$ organizes
 the generators of $\LL$. We denote the presentations variables by
$u, \TT_1, \TT_2, \TT_3, \TT_4$, with $u$ corresponding to $a$.

\begin{itemize}

\item[$\bullet$] $(1,4,1)$: We know (Theorem~\ref{141}) that $JI:I^2= \m$. This means that we have forms
\begin{eqnarray*}
\hh_1  & = & xu^2 + \cdots \\
\hh_2  & = & yu^2 + \cdots \\
\hh_3 & = & zu^2 + \cdots \\
\hh_4 & = & w u^2 + \cdots
\end{eqnarray*}
with the $(\cdots)$ in $(\TT_1, \TT_2, \TT_3, \TT_4)\BB_1$. The corresponding resultant, of degree $8$, is nonzero.

\item[$\bullet$] $(1,2,1)$: There are two forms of degree $1$ in $\LL$,
\begin{eqnarray*}
\ff_1 & = & x u + \cdots\\
\ff_2 & = & y u + \cdots
\end{eqnarray*}
The forms in $L_2/\BB_1L_1$ have coefficients in $\m^2$. This will follow from $I_1(\phi)= \m$. We need a way to generate
two forms of degree $3$. Since we expect $JI:I^2 = \m$, this would mean the presence of two forms in $L_3$,
\begin{eqnarray*}
\hh_1^{*} & = & z u^3 + \cdots\\
\hh_2^{*} & = & w u^3 + \cdots ,
\end{eqnarray*}
which together with $\ff_1$ and $\ff_2$ would give the nonzero degree $8$ resultant.%

\item[$\bullet$] $(1,3,1)$: There is a form $\ff_1 = xu + \cdots \in L_1$ and two forms in $L_2$
\begin{eqnarray*}
\hh_1 & = & y u^2 + \cdots \\
\hh_2 & = & zu^2+ \cdots
\end{eqnarray*}
predicted by \cite[Proposition 3.7]{syl2}  if $I_1(\phi) = (x,y,z,w^2)$. (There are indications that this is always the case.)
 We need a cubic equation $\hh_3^{*} = wu^3+ \cdots$ to launch the nonzero resultant of degree $8$.

\end{itemize}

For all quaternary quadrics with $\RR[It]$ almost Cohen--Macaulay, Corollary~\ref{degT} says that
$\nu(T) \leq \lambda(\RR/J:a)$. Let us compare to the actual number of generators in the examples discussed above:

\[  \left[
\begin{array}{ccc}
\nu(T) &  & \lambda(I/J)\\
5 & (1,4,1) & 6\\
4 & (1,3,1) & 5\\
4 & (1,2,1)& 4
\end{array} \right]
\]
We note that in the last case, $T$ is an  Ulrich module.

\subsection{Monomial ideals}\index{aCM algebras defined by binomial ideals}

Monomial ideals of finite colength which are almost complete intersections have a very simple description. We examine a narrow class of them.
Let $\RR=k[x,y,z]$ be  a polynomial ring over and infinite field and let $J$ and $I$ be $\RR$--ideals such that
\[{\ds J=(x^{a},\; y^{b},\; z^{c}) \subset (J,\; x^{\alpha} y^{\beta} z^{\gamma})=I. }\]
This is the general form of almost complete intersections of $\RR$ generated by monomials. Perhaps the most interesting cases are those
where
${\ds \frac{\alpha}{a} + \frac{\beta}{b} + \frac{\gamma}{c} <1}$. This inequality
ensures that $J$ is not a reduction of $I$.
Let
\[{\ds Q =(x^{a}-z^{c},\; y^{b}-z^{c},\;x^{\alpha} y^{\beta} z^{\gamma} ) }\]
 and suppose that ${\ds a > 3 \alpha,\; b > 3 \beta,\;
c > 3 \gamma }$. Note that $I=(Q,\; z^c)$.  Then $Q$ is a minimal reduction of $I$ and the reduction number $\mbox{\rm red}_{Q}(I) \leq
2$. We label these ideals $I(a,b,c,\alpha, \beta, \gamma)$.

\medskip

We will examine in detail the case $a=b=c=n\geq 3$ and
$\alpha=\beta = \gamma =1$. We want to argue that
$\RR[It]$ is almost Cohen--Macaulay. To benefit from the monomial generators in using {\em Macaulay2} we set
$I = (xyz, x^n, y^n,z^n)$. Setting $\BB=\RR[u, \TT_1, \TT_2, \TT_3]$, we claim that
\[
\LL= (z^{n-1}u - xy\TT_3, y^{n-1}u - xz\TT_2, x^{n-1}u - yz\TT_1, z^n\TT_2-y^n\TT_3,  z^n\TT_1-x^n\TT_3,
y^n\TT_1-x^n\TT_2,
\]
\[ y^{n-2}z^{n-2}u^2 - x^2\TT_2\TT_3, x^{n-2} z^{n-2} u^2 -y^2 \TT_1\TT_3,   x^{n-2}y^{n-2} u^2 - z^2\TT_1\TT_2,
x^{n-3}y^{n-3}z^{n-3}u^3 - \TT_1\TT_2\TT_3).
  \]
 We also want to show that these ideals define an almost Cohen--Macaulay Rees algebra.

\bigskip

There is a natural specialization argument. Let $X$, $Y$ and $Z$ be new indeterminates and let
$\BB_0 = \BB[X,Y,Z]$. In this ring define the ideal $\LL_0$ obtained by replacing in the list above of generators
of $\LL$, $x^{n-3}$ by $X$ and accordingly $x^{n-2}$ by $xX$, and so on; carry  out similar actions on the other variables:

\[
\LL_0= (z^2 Zu - xy\TT_3, y^{2}Yu - xz\TT_2, x^{2}Xu - yz\TT_1, z^3Z\TT_2-y^3Y\TT_3,  z^3Z\TT_1-x^3X\TT_3,
y^3 Y\TT_1-x^3X\TT_2,
\]
\[ yzYZu^2 - x^2\TT_2\TT_3, x z XZ u^2 -y^2 \TT_1\TT_3,   xyXY u^2 - z^2\TT_1\TT_2,
XYZu^3 - \TT_1\TT_2\TT_3).
  \]

Invoking {\em Macaulay2} gives a (non-minimal) projective resolution
\[ 0 \rar \BB_0^4 \stackrel{\phi_4}{\lar}
 \BB_0^{17} \stackrel{\phi_3}{\lar}
\BB_0^{22} \stackrel{\phi_2}{\lar}
 \BB_0^{10} \stackrel{\phi_1}{\lar}
  \BB_0 \lar \BB_0/\LL_0 \rar 0.
 \]

We claim that the specialization $X \rar x^{n-3}$, $Y \rar y^{n-3}$, $Z \rar z^{n-3}$ gives a projective resolution
of $\LL$.

\begin{itemize}

\item[$\bullet$] Call $\LL'$ the result of the specialization in $\BB$. We argue that $\LL' = \LL$.

\medskip

\item[$\bullet$] Inspection of the Fitting ideal $F$ of  $\phi_4$ shows that it contains $(x^3, y^3,z^3, u^3, \TT_1\TT_2\TT_3)$.
From standard theory, the radicals of the Fitting ideals of $\phi_2$ and $\phi_2$  contain $\LL_0$, and therefore
the radicals of the Fitting ideals of these mappings after specialization will contain the ideal $(L_1)$ of $\BB$, as
$L_1 \subset \LL'$.

\medskip

\item[$\bullet$] Because $(L_1)$ has codimension $3$,
 by the acyclicity
 theorem (\cite[1.4.13]{BH})
  the complex gives a projective resolution of $\LL'$. Furthermore, as $\mbox{\rm proj. dim }\BB/\LL' \leq 4$,
 $\LL'$ has no associated primes of codimension $\geq 5$. Meanwhile
the Fitting ideal of $\phi_4$ having codimension $\geq 5$,  forbids the existence of associated primes
of codimension $4$. Thus $\LL'$ is
unmixed.

\medskip

\item[$\bullet$] Finally, in $(L_1) \subset \LL' $, as $\LL'$ is unmixed its associated primes are minimal primes of
$(L_1)$, but by Proposition~\ref{canoseq}(c), there are just two such, $\m\BB$ and $\LL$. Since
$\LL' \not \subset \m\BB$, $\LL$ is its unique associated prime. Localizing at $\LL$ gives
the equality of $\LL' $ and $\LL$ since $\LL$ is a primary component of $(L_1)$.

\end{itemize}

Let us sum up this discussion:

\begin{Proposition} \label{nnn111} The Rees algebra of $I(n, n, n, 1, 1, 1)$, $n\geq 3$, is almost Cohen--Macaulay.
\end{Proposition}

\begin{Corollary} $e_1(I(n,n,n,1,1,1)) = 3(n+1)$.
\end{Corollary}

\begin{proof} Follows easily since $e_0(I) = 3n^2$, the colengths of the monomial ideals $I$ and $I_1(\phi)$ directly
calculated and $\red_J(I) = 2$ so that
\[ e_1(I) =  \lambda(I/J) + \lambda(I^2/JI) = \lambda(I/J) + [\lambda(I/J) - \lambda(\RR/I_1(\phi))]=
(3n-1) + 4.\]
\end{proof}

\begin{Remark}{\rm
We have also experimented with other cases beyond those with $xyz$ and in higher dimension as well.
\begin{itemize}

\item[{\rm (a)}] In $\dim \RR=4$, the ideal $I= I(n,n,n,n,1,1,1,1)= (x_1^n, x_2^n, x_3^n, x_4^n, x_1x_2x_3x_4)$, $n\geq 4$, has
a Rees algebra $\RR[It]$ which is almost Cohen--Macaulay.

\item[{\rm (b)}] The argument used was a copy of the previous case, but we needed to make an adjustment in the last step
to estimate the codimension of the last Fitting ideal $F$ of the corresponding mapping $\phi_5$. This is a large matrix,
so it would not be possible to find the codimension of $F$ by looking at all its maximal minors. Instead, one argues as follows.
Because $I$ is $\m$--primary, $\m\subset \sqrt{F}$, so we can drop the entries in $\phi_5$ in $\m$ Inspection will give
$u^{16} \in F$, so dropping all $u$'s gives an additional minors in $\TT_1, \ldots, \TT_4$, for $\height (F)\geq 6$.
This suffices to show that $\LL = \LL'$.

%\item It seems tempting to conjecture that in all dimensions these ideals have Rees algebras which are aCM.

\end{itemize}

}\end{Remark}

\begin{Conjecture} \label{monoaciRacm}
{\rm Let $I$ be a monomial ideal of $k[x_1, \ldots, x_n]$. If $I$ is an almost complete intersection
of finite colength its Rees algebra $\RR[It]$ is almost Cohen-Macaulay.

\begin{Question} \label{GorRat}
{\rm Our discussion of Rees algebras of almost complete intersections show the important role
that Gorenstein ideals play in their analyses. Can these roles be reversed, that is, what are
the properties of (bi-) rational morphisms defined by Gorenstein ideals? For example, let
$\RR = k[x,y,z] $ and $I$ be the ideal of codimension three generated by $5$ quadrics; what are
the properties of $\RR[It]$?

}\end{Question}

}\end{Conjecture}

\section{Special Fiber of Filtrations}

\subsection{Introduction}
Let $\RR$ be a Noetherian ring and let $M$ be a finitely generated
$\RR$-module and let
$ \mathcal{A}=\{A_n, n\geq 0\}$ be a decreasing filtration of ideals
of $\RR$ and
$ \mathcal{M}=\{M_n, n\geq 0\}$  a decreasing
filtration of submodules of $M$ such that $A_m\ M_n\subset M_{m+n}$.
 Denote by
\[ \gr_{\mathcal{A}}(\RR) = \bigoplus_{n} A_n/A_{n+1}, \quad
 \gr_{\mathcal{M}}(M) = \bigoplus_{n} M_n/M_{n+1}\]
the corresponding  associated graded modules.

\begin{definition} If $(\RR, \mathfrak{m})$ is a local ring, the {\em
special fiber} \index{special fiber of a module} of $M$ is the
 graded module
\begin{eqnarray*}
\mathcal{F}({M}) &=& \gr_{\mathcal{M}}(M) \otimes
\RR/\mathfrak{m}
\end{eqnarray*}
over the graded ring $\mathcal{F}(\RR)$.
\end{definition}

The main class of filtrations we will deal with will be $I$-good
filtrations associated to an ${\mathfrak{m}}$-primary ideal $I$ of a
local ring $\RR$:  $A_{n+1}= IA_n$ and $M_{n+1}= IM_n$ for $n\gg
0$.
\medskip

One of our motivations is the following. Let $\ff= \{f_1, \ldots,
f_m\}\subset \RR$ be a set of elements of the local ring $(\RR,
\mathfrak{m})$. There are
several algebras defined by $\ff$, noteworthy the Rees algebra
$\AA=\RR[It]$, $I=(\ff)$, and its special fiber $\mathcal{F}(I)=
\RR[It]\otimes_{\RR}(\RR/\mathfrak{m})$.
Let $\BB=\bar{\AA}$. If $I$ is $\mathfrak{m}$-primary, the associated
graded rings of $\AA$ and $\BB$ have the same multiplicity $\rme_0(I)$.
However the multiplicities of their special fibers--denoted by
$f_0(I)$ and $\bar{f}_0(I)$--will often differ.
Another change may occur in the values of the Hilbert coefficients
$\rme_1(I)$ and $\bar{e}_1(I)$.
% We are interested in comparing the
%changes: $f_0(I)\rar \bar{f}_0(I)$ and $\rme_1(I)\rar \bar{e}_1(I)$.

\subsection{The reduction number formula}
\index{reduction number formula}

%\subsection{Extended degrees and reduction numbers}\index{extended
%degree!reduction number}

One can use the philosophy of $\Deg$ to obtain generalizations of the
bounds for reduction numbers of ideals in Cohen-Macaulay local rings to
arbitrary Noetherian local rings.

In order to make use of Theorem~\ref{e1hs}, we need information about
the reduction number of $L$ in terms related to multiplicity. Let us
recall \cite[Theorem 2.45]{icbook}:

\begin{Theorem}\label{mainresredCM}\index{reduction number formula:
Cohen--Macaulay case}
 Let $(\RR,\m)$ be a Cohen-Macaulay local ring of dimension $d$ and
infinite residue field. For an $\m$-primary ideal $I$, \[ \mbox{\rm
red}(I) \leq \frac{d\cdot \l(\RR/J)}{o(I)} - 2d+1,\] where $J$ is a minimal
reduction of $I$ and $o(I)$ is the $\m$-adic order of $I$.
\end{Theorem}

To establish such a result for arbitrary Noetherian rings (\cite{chern4}), we proceed
differently.
The version of the following lemma for Cohen-Macaulay rings can be
found in \cite[Chapter 3, Theorem 1.1]{Sal78}.

\begin{Lemma} \label{lengthdim1}
Let $(\RR, \m)$ be a Noetherian local ring of dimension one.
\begin{enumerate}
\item[{\rm (a)}]
 For any
ideal $I$ and for any parameter $x$, $\nu(I) \leq \l(\RR/(x))$. More
generally, for any submodule $U$ of a module $E$ generated by $n$
elements, $\nu(U)\leq \l(\RR/(x))\cdot n$.

\item[{\rm (b)}]
If $\RR$ is Cohen--Macaulay and $x$ belongs to $\m^s$, then ${\ds \nu(U) \leq \frac{\l(\RR/(x))}{s}\cdot n}$.
\end{enumerate}
\end{Lemma}

\begin{proof} (a) {\bf Case 1:} Let $I$ be an $\RR$--ideal and $x$ a parameter of
$\RR$.
Set $L =\H_{\m}^0(\RR)$ and consider the two exact sequences
\[ 0 \rar L \lar \RR \lar \RR' \rar 0, \hspace{1 in} 0 \rar L_0 =I\cap L \lar I \lar I' \rar 0.  \]
Note that $\RR'$ has positive depth.
Tensoring by $\RR/(x)$, we get the exact sequences
\[  0 \rar L/xL \lar \RR/x\RR \lar \RR'/x\RR' \rar 0,  \hspace{0.5 in}   0 \rar L_0/xL_0 \lar I/xI \lar I'/xI' \rar 0.\]
The second of the preceding sequences gives
\[ \nu(I) =\l(I/ \m I )  \leq \l(I/xI) = \l(L_0/xL_0) + \l(I'/xI').\]
Note that ${\ds \l(I'/xI') = \l(\RR'/x\RR') }$. Combining these inequalities we obtain :
\[\begin{array}{lll}
\nu(I) & \leq &  \l(L_0/xL_0) + \l(I'/xI') \\ && \\
       &=& \l(0:_{L_0}x) +  \l(I'/xI') \\ && \\
       &\leq &  \l(0:_{L}x) +  \l(I'/xI') \\ && \\
       &=& \l(L/xL)  +  \l(I'/xI') \\ && \\
       &=& \l(\RR/x\RR) - \l(\RR'/x\RR') +  \l(I'/xI') \\ && \\
       &=& \l(\RR/x\RR).
\end{array}
\]

\medskip

\noindent {\bf Case 2 :} Let $U$ be a submodule of a free
$\RR$--module $F=\RR e \oplus \RR^{n-1}$. We use induction on $n$ to
show that ${\ds \nu(U) \leq \l(\RR/x\RR)\cdot n }$. If $n=1$, it is done
by Case I. Suppose that $n \geq 2$. Consider the following short
exact sequence \[ 0 \lar \RR e \lar F=\RR e \oplus \RR^{n-1}
\stackrel{\pi}{\lar} \RR^{n-1} \lar 0, \] which is split. Let $U_0=U
\cap \RR e \subset \RR e$ and $U'=\pi(U) \subset \RR^{n-1} $. By
induction hypothesis, we get \[ \nu(U) \leq \nu(U_{0}) + \nu(U') \leq
\l(\RR/x\RR) + \l(\RR/x\RR) \cdot (n-1) = \l(\RR/x\RR)\cdot n.  \]

\medskip

\noindent {\bf Case 3 :} Let $U$ be a submodule of $E$ with $\nu(E)=n$. Then there is a natural surjection :
\[  \pi : \RR^n \surj E \lar 0.   \] Then $\pi^{-1}(U)$ is a
submodule of $\RR^n$ so that $\nu(\pi^{-1}(U)) \leq \l(\RR/x\RR)\cdot n$  by Case II.
Therefore we obtain
\[ \nu(U) \leq \nu(\pi^{-1}(U)) \leq \l(\RR/x\RR) \cdot n. \]

\medskip

\noindent (b) Assume that $\RR$ is Cohen--Macaulay and let $y\RR$ be a
minimal reduction of $\m$.
%%%%%%%%%%%%%%%%%%%%%%%%%%%%%%%%%%%%%%%%%%%
% For any integer $t$,
%\[  \l(R/x^t) \geq \l(\m^{st}). \]
%%%%%%%%%%%%%%%%%%%%%%%%%%%%%%%%%%%%%%%%%
For sufficiently large $t$, we have \[t \l(\RR/x\RR) =\l(\RR/x^t\RR) \geq
\l(\m^{st}) = e_{0}(\m) st - e_{1}(\m) = \l(\RR/y\RR) st - - e_{1}(\m).
\] Comparing the leading coefficients of polynomials in $t$, we get
\[ \l(\RR/x\RR) \geq \l(\RR/y\RR) s. \] Applying part(a), we obtain \[ \nu(U)
\leq \l(\RR/y\RR) \cdot n \leq \frac{\l(\RR/x\RR)}{s} \cdot n.\] 
\end{proof}

%%%%%%%%%%%%%%%%%%%%%%%%%%%%%%%%%%%%%%%%%%%%%%%%%%%%%%%%%%%%%%%%%%%%%%%%%%%%%%%%%%%%%%%%%%%%%%%%%%%%%%%%%%%%%%%
%%% {\bf Comment:} To be useful we must have an idea of the reduction number of $L$ in terms of the data. See the next frame.
%%%%%%%%%%%%%%%%%%%%%%%%%%%%%%%%%%%%%%%%%%%%%%%%%%%%%%%%%%%%%%%%%%%%%%%%%%%%%%%%%%%%%%%%%%%%%%%%%%%%%%%%%%%%%%%

\begin{Theorem}\label{mainresrednCM}
 Let $(\RR,\m)$ be a Noetherian local ring of dimension $d$ and
infinite residue field. For an $\m$-primary ideal $I$ there is  a minimal reduction $J$ of $I$, \[ \mbox{\rm
red}_{J}(I) \leq d\cdot \l(\RR/J) - 2d+1.\]
%%%%%%%%%%%%%%%%%%%%%%
% unless $\RR$ is a regular local ring and $I=\m$.
%%%%%%%%%%%%%%%%%%%%%%%%%%
\end{Theorem}

\begin{proof} Let $J=(x_1, \ldots, x_{d})$ be a minimal reduction of $I$ and
let $J_0= (x_1, \ldots, x_{d-1})$. Then \[ \nu(I^n) \leq \nu(J_{0}^n)
+ \nu(I^n/J_{0}^n).  \] We need to estimate ${\ds \nu(I^n/J_{0}^n)}$,
where $I^n/J_0^n$ is an ideal of the $1$-dimensional local ring
$\RR/J_0^n$.

\medskip

%%%%%%%%%%%%%%%%%%%%%%%%%%%%%%%%%%%%%%%%%%%%%%%%%%%%%%%%%%%%%%%%%%%%%%%%%%%%%%%%%%%%%%%%%%%
% \[\begin{CD}
%  0 @>>> {\ds  \frac{J_{0}^{n-1} \cap I^{n}}{J_{0}^{n}}} @>>> {\ds \frac{I^n}{J_{0}^n } }  @>>> {\ds
% \frac{I^{n}}{J_{0}^{n-1} \cap I^{n}} = \frac{I^n + J_{0}^{n-1}}{J_{0}^{n-1}}}   @>>> 0 \\ &&\\
% 0 @>>> {\ds   \frac{(I^n + J_{0}^{n-1}) \cap J_{0}^{n-2}}{J_{0}^{n-1}}  }  @>>> {\ds \frac{I^n +
% J_{0}^{n-1}}{J_{0}^{n-1}}}  @>>> {\ds \frac{I^n + J_{0}^{n-2}}{J_{0}^{n-2}}} @>>> 0\\ && \\
% && \vdots \\ && \\
% 0 @>>> {\ds  \frac{(I^n + J_{0}^{n-i}) \cap J_{0}^{n-i-1}}{J_{0}^{n-i}}  }  @>>> {\ds \frac{I^n +
% J_{0}^{n-i}}{J_{0}^{n-i}}}  @>>> {\ds \frac{I^n + J_{0}^{n-i-1}}{J_{0}^{n-i-1}}} @>>> 0\\ && \\
% && \vdots \\ && \\
% 0 @>>> {\ds  \frac{(I^n + J_{0}^{2}) \cap J_{0}}{J_{0}^{2}}  }  @>>> {\ds \frac{I^n + J_{0}^{2}}{J_{0}^2}}  @>>> {\ds % \frac{I^n + J_{0}}{J_{0}}} @>>> 0.
% \end{CD}
%\]
%%%%%%%%%%%%%%%%%%%%%%%%%%%%%%%%%%%%%%%%%%%%%%%%%%%%%%%%%%%%%%%%%%%%%%%%%%%%%%%%%%%%%%%%%%%

\noindent Let ${\ds M_{i} = \frac{(I^n + J_{0}^{n-i+1}) \cap J_{0}^{n-i}}{J_{0}^{n-i+1}} }$ and ${\ds N_{i}=\frac{I^n + J_{0}^{n-i+1}}{J_{0}^{n-i+1}}}$. Then  we obtain the following series of exact sequences :
\[ 0 \lar M_{i} \lar N_{i} \lar N_{i+1} \lar 0,  \]
where $i=1, \ldots, n-1$. Note that for each $i$, $M_{i}$ is a submodule of ${\ds J_{0}^{n-i}/J_{0}^{n-i+1}}$ as an $\RR/J_{0}$--module.
Hence by Lemma~\ref{lengthdim1}, for each $i=1, \ldots, n-1$,
\[ \nu(M_{i}) \leq \l(\RR/J) \nu(J_{0}^{n-i}/J_{0}^{n-i+1}) = \l(\RR/J) {\ds {{d+n-2-i}\choose{d-2}}},    \]
Since ${\ds N_{n} = \frac{I^n + J_{0}}{J_{0}} }$ is a submodule of $\RR/J_{0}$, by Lemma~\ref{lengthdim1}, we get
\[ \nu(N_{n}) \leq \l(\RR/J).    \]
Therefore we obtain
\[ \begin{array}{lll}
\nu(I^n) &\leq & \nu(J_{0}^n) + \nu(I^n/J_{0}^n) \\ && \\
         % &\leq & \nu(J_{0}^n) + \nu(M_{1}) + \nu(N_{2}) \\ && \\
         %&\leq & \nu(J_{0}^n) + \nu(M_{1}) + \nu(M_{2}) +\nu(N_{2}) \\ && \\
         % &     & \vdots \\ && \\
         &\leq & \nu(J_{0}^n) + \nu(M_{1}) + \nu(M_{2}) + \cdots + \nu(M_{n-1})+ \nu(N_{n}) \\ && \\
         &\leq & {\ds {{d+n-2}\choose{d-2}} + \l(\RR/J) \sum_{i=1}^{n}  {{d+n-2-i}\choose{d-2}}} \\ && \\
         &=& {\ds {{d+n-2}\choose{d-2}} + \l(\RR/J) {{d+n-2}\choose{d-1}.}   }
\end{array}
\]
Recall that if
\[ {\ds  \nu(I^n) < {{n+d}\choose{d}},  } \]
then the minimal reduction number $r=r(I)$ of $I$ satisfies $r \leq n-1$  (\cite{ES}, \cite[Theorem 2.36]{icbook}).
At this point we may need to replace the reduction $J$ by another one.
Hence by solving the inequality
\[{\ds   {{d+n-2}\choose{d-2}} + \l(\RR/J) {{d+n-2}\choose{d-1}}   <  {{n+d}\choose{d}},  }  \]
we obtain the desired relation.
\end{proof}

%%%%%%%%%%%%%%%%%%%%%%%%%%%%%%%%%%%%%%%%%%%%%%%%%%%%%%%%%%%%%%%%%%%%%%%%%%%%%%%%%%%%%%%%%%%%%%%%%%%%%%
% \noindent Proof of the last equality by induction on $n$. We want to show that
% \[{\ds \sum_{i=1}^{n}  {{d+n-2-i}\choose{d-2}} =  {{d+n-2}\choose{d-1}}.     }    \]
% It is clear when $n=1$. Suppose that $n \geq 2$. Then
% \[\begin{array}{lll}
% {\ds \sum_{i=1}^{n+1}  {{d+n-1-i}\choose{d-2}}} &= & {\ds {{d+n-2}\choose{d-2}} + {{d+n-3}\choose{d-2}} + \cdots + {{d-2}\choose{d-2}}} \\ && \\ &=& {\ds
% {{d+n-2}\choose{d-2}} +  \sum_{i=1}^{n}  {{d+n-2-i}\choose{d-2}}  } \\ && \\
% &=& {\ds {{d+n-2}\choose{d-2}} +  {{d+n-2}\choose{d-1}}  } \\ && \\
% &=& {\ds {{d+n-1}\choose{d-1}}  }
% \end{array}
% \]
%%%%%%%%%%%%%%%%%%%%%%%%%%%%%%%%%%%%%%%%%%%%%%%%%%%%%%%%%%%%%%%%%%%%%%%%%%%%%%%%%%%%%%%%%%%%%%%%%

In the formula of Theorem~\ref{mainresrednCM},
\[ e_0(I)\leq \l(\RR/J) \leq \operatorname{hdeg}_I(\RR),\]
where $\operatorname{hdeg}_I(\RR)$ is the extended multiplicity of
$\RR$ in the theory of cohomological multiplicities
(see \cite{icbook}). Let us derive an inequality for
$\operatorname{red}(I)$
that shows its dependence on $a=e_0(I)$ and
$b=I(\RR)=\operatorname{hdeg}_I(\RR)-e_0(I)$, the {\em Cohen-Macaulay
 defficency} of $\RR$ relative to this multiplicity.

\medskip

The argument is similar to the one above but uses the following
majorization  for $d\geq 2$ (\cite[Theorem 2.1]{RTV}):
\[ \l(I^n/I^{n+1})\leq a {{n+d-1}\choose {d-1}} + b {{n+d-2}\choose
{d-2}}.
\]

As  above, it suffices to solve for $n$ in
\[ \nu(I^n)\leq \l(I^n/I^{n+1})\leq a {{n+d-1}\choose {d-1}} + b {{n+d-2}\choose
{d-2}}\leq {n+d\choose d}.
\]
In other words, any integer larger than the larger root of the
polynomial
\[ n^2 -(da -2d+1)n - d(d-1)(a+b-1)\]
will do. For instance, if $b$ is $\neq 0$ and $a+b\geq d+1$, then
$n=d(a+b-2)$ will work. Note that this is slightly better than
$d(a+b)-2d+1$, the value predicted in
 Theorem~\ref{mainresrednCM}.

\begin{Corollary}
Let $(\RR, \m)$ be a Noetherian local ring of dimension $d>0$ and assume $J = (a_1, \ldots, a_d)$ is a minimal reduction of $\m$.
Let $\m=(J, a_{d+1}, \ldots, a_g)$ be a minimal set of generators of $\m$.
Then
\[ e_{1}(\m) \leq e_{1}(\m) - e_{1}(J) \leq \l(\RR/J:\m) \left[{\ds {{g+\l(\RR/J)d -3d+1}\choose{g-d}}} -1 \right] f_{0}(J). \]
\end{Corollary}

%%%%%%%%%%%%%%%%%%%%%%%%%%%%%%%%%%%%%%%%%%%%%%%%%%%%%%%%%%%%%%%%%%%%%%%%%%%%%%%%%%%
% {\bf Comment:} The two results can be combined to give an expression for $\rme_1(L)$ in terms of various data on $I$ and % the number $m$. It is not cute.
%%%%%%%%%%%%%%%%%%%%%%%%%%%%%%%%%%%%%%%%%%%%%%%%%%%%%%%%%%%%%%%%%%%%%%%%%%%%%%%%%%%%%%%%%%%%

\subsection{Special fibers and birationality}

 Let $\RR$ be local ring and $I \subset L$ be ideals
with $L$ integral over $I$. This implies, of course, that $\BB= \RR[Lt]$ is
rational over $\AA=\RR[It]$. What should  it mean to say that special
fiber $\mathcal{F}(L)$ is birational over $\mathcal{F}(I)$?
Perhaps simply that the cokernel of $\mathcal{F}(I) \rar
\mathcal{F}(L)$ has dimension $< \dim \mathcal{F}(I)$.

\bigskip

One can develop a brief theory here, some of which can be extended to
non-primary ideals. 
We are going to carry out a comparison between $\deg \mathcal{F}(I)$
and
$\deg \mathcal{F}(L)$. Consider the exact sequence
\[ 0 \rar \Rees= \RR[It] \lar \Rees' = \RR[Lt] \rar D \rar 0.\]
Reducing mod $\mathfrak{m}$ we have the exact sequence
\[ \Tor_1^{\RR}(D, k) \lar \mathcal{F}= \bar{\Rees}\lar  \mathcal{F}'=
\bar{\Rees'} \lar \bar{D} \rar 0.
\]

We know that $\dim D=\dim \bar{D}$, and it is easy to see [just look
at the associated primes of $D$] that
$\dim \Tor_1^{\RR}(D,k)\leq \dim D$.

\medskip

If $\dim D< d$ [a condition that is equivalent to $\rme_1(I)=e_1(L)$],
then
$\deg \mathcal{F}= \deg \mathcal{F}'$.

\medskip

Let us assume $\dim D=d$. In this case, $\deg D =\rme_1(L)-e_1(I)$.
 From the homology sequence above, we have
\begin{eqnarray} \label{degbar1}
 \deg \mathcal{F} + \deg \bar{D} &\leq &
\deg \mathcal{F}' + \deg \Tor_1^{\RR}(D,k).
\end{eqnarray}
Let us estimate  $\deg \Tor_1^{\RR}(D,k)$. Let
\[ \cdots \lar F_2 \lar F_1 \lar \RR=F_0 \lar k \rar 0 \] be a minimal
resolution of $k$. Tensoring with $D$ we have the complex
\[ \cdots \lar D\otimes F_2 \lar  D\otimes  F_1 \lar D= D\otimes  F_0
\lar \bar{D}  \rar 0. \] In the subcomplex
\[ D\otimes F_2 \stackrel{\varphi}{\lar} D\otimes F_1
\stackrel{\psi}{\lar}
 \mathfrak{m}D \rar 0\]
$\Tor_1^{\RR}(D,k)= Z/L$, where $Z = \ker \psi $ and $L= \mbox{\rm image
}\varphi$, and therefore $\deg \Tor_1^{\RR}(D,k)\leq \deg Z$ (note that
$\dim Z\leq \dim D).$
On the other hand, from the exact sequence
\[ 0 \rar Z \lar D\otimes F_1 \lar \mathfrak{m}D \rar 0,\]
we have
\begin{eqnarray}
\label{degbar2}
\deg Z & = & \deg D\otimes F_1 - \deg \mathfrak{m}D.
\end{eqnarray}

\medskip

We can now assemble these calculations:

\begin{Theorem} Let $I\subset L$ be $\mathfrak{m}$-primary ideals
with the same integral closure. Then
\[ \deg \mathcal{F} -\deg \mathcal{F}' \leq
(\beta_1(k)-1)(\rme_1(L)-\rme_1(I)).
\]
\end{Theorem}

\begin{proof} From (\ref{degbar1}) and (\ref{degbar2}), we have
\begin{eqnarray*} \label{degbar3}
 \deg \mathcal{F} -
\deg \mathcal{F}'
& \leq &  \deg \Tor_1^{\RR}(D,k)- \deg \bar{D}\\
& \leq & \deg Z -\deg \bar{D} \\
&=& \deg D\otimes F_1 - \deg \mathfrak{m}D - \deg \bar{D} \\
&=& \deg D \otimes F_1 -\deg D = (\rank(F_1)-1)\deg D \\
&=& (\beta_1(k)-1)(e_1(L)-\rme_1(I)).
\end{eqnarray*}
\end{proof}

\begin{Remark}{\rm For future examination, we make the following
observations.

\begin{enumerate}

\item[{\rm (a)}] We are aware that if $\mathcal{F}$ is a domain,
$\mathcal{F}\subset \mathcal{F}' $, so the formula above would add
nothing to it. We need an example with
 $\deg \mathcal{F}> \deg \mathcal{F}' $.

\item[{\rm (b)}] It clear that the same formulas hold for general integral
filtrations [subfiltrations of the integral closure filtrations, not
just $I$-adic filtrations]. In this respect, a straightforward fact is
that if $\AA$ is such an algebra and $\BB$ is its $S_2$-ification, then
\[ \deg \mathcal{F}(\AA) =  \deg \mathcal{F}(\BB).
\]

\item[{\rm (c)}] Looking over the proof above, it is clear that several pieces
of homology were disregarded. Maybe this weakened the final
inequality but we could not see how. This may warrant another visit
to the proof. [The original hand-written has some unclear comments.]

\item[{\rm (d)}] If this type of comparison really works, it might be worthwhile
to compare them with the estimates for $\deg \mathcal{F}$ to be found
in \cite[Chapter 2]{icbook}.

\end{enumerate}
}\end{Remark}

\subsection{Variation of Hilbert coefficients}
\index{variation of Hilbert coefficients}
We examine a list of questions about the
changes that $\rme_0(I)$ and $\rme_1(I)$ undergo when we enlarge $I$. A
must case is: \[ \rme_0(I),\rme_1(I) \lar e_0(L),\rme_1(L), \quad L=(I,x).
\] Clearly the optimal baseline is that of an ideal $I$ generated by
a system of parameters, but we will consider very general cases.
 As will be seen, some relationships involve
the Hilbert coefficients $f_0(I)$ and $f_1(I)$ of the special fiber.
We will follow te discussion of \cite{chern4}.

\medskip

In our calculations we make repeated used of
 the following elementary observation.

\begin{Lemma}\label{tensor}
If $(\RR, \m)$ is a Noetherian local ring and $M$ is a module of finite length, then
\[ \l(M\otimes N)\leq \l(M)\cdot \nu(N).\]
\end{Lemma}

\begin{proof} Induct on $\l(M)=n$. If $n=1$, then $M\simeq \RR/\m$ so that the
assertion is clear. Suppose that $n \geq 2$ and let \[{\ds M=M_{0}
\supsetneq M_{1} \supsetneq \cdots \supsetneq M_{n-1} \supsetneq
M_{n}=0 }\] be the composition series of $M$. By tensoring $0\rar
M_{n-1} \rar M \rar M/M_{n-1} \rar 0$ with $N$, we get the right
exact sequence \[ M_{n-1} \otimes N \rar M\otimes N \rar
M/M_{n-1}\otimes N \rar 0.\] Since $\l(M_{n-1})=1$ and
$\l(M/M_{n-1})=n-1$, by induction hypothesis we have \[\l(M_{n-1}
\otimes N) \leq \nu(N),\;\; \mbox{\rm and}\;\; \l(M/M_{n-1}\otimes N)
\leq (n-1)\nu(N). \] Therefore, \[ \l(M \otimes N) \leq \l(M_{n-1}
\otimes N) + \l(M/M_{n-1}\otimes N) \leq (1+n-1)\nu(N)= \l(M)
\nu(N).\] 
\end{proof}

\begin{Theorem}[Variation of the multiplicity]  \index{variation of the multiplicity} Let $(\RR, \m)$ be a Noetherian local ring of dimension
$d$ and infinite residue field, and let $I\subset L=(I,h)$ be
$\m$--primary ideals. Then \[ e_0(I) - e_0(L) \leq \l(\RR/I:L)\cdot
f_0(I).\]
\end{Theorem}

\begin{proof} For $n \in \bbn$, consider the following filtration :
\[\begin{array}{lll}
{\ds I^{n}=M_{0} \subset M_{1}=(M_{0},\; I^{n-1}h)} \subset \cdots &\subset& {\ds M_{r-1}=(M_{r-2},\; I^{n-r+1}h^{r-1}) }\\ && \\ &\subset&   {\ds M_{r}=(M_{r-1},\; I^{n-r}h^{r})} \\ && \\
 &\subset&  {\ds  \cdots \subset M_{n}=(M_{n-1}, h^n)=L^n.}
\end{array}\]
Then we obtain
\begin{eqnarray*}
\l(\RR/I^n)-\l(\RR/L^n) &=& \l((I,h)^n/I^n)= \l((I^n,hI^{n-1}, \ldots, h^{n-1}I, h^n)/I^n) \\
&= & \sum_{r=1}^{n}  \l(M_r/M_{r-1}).
\end{eqnarray*}
For each $r$, $M_r/M_{r-1}$ is generated by $h^rI^{n-r}+M^{r-1}$.]
Consider the natcityural surjection
\[ \zeta: \RR/(I:h)\otimes I^{n-r} \surj M_r/M_{r-1} =<h^rI^{n-r}+M^{r-1} >, \]where $\zeta(\cl{r} \otimes x)=rh^rx+M^{r-1}$. Using Lemma~\ref{tensor}, we have
\[{\ds \l(M_r/M_{r-1}) \leq \l( \RR/(I:h)\otimes I^{n-r}) = \l(\RR/(I:h)) \nu(I^{n-r}). }\]
It follows that
\[ \l(\RR/I^n)-\l(\RR/(I,h)^n) \leq \lambda(\RR/I:h) \sum_{r=0}^{n-1} \nu(I^r).\]
 The iterated Hilbert function ${\ds \sum_{r=0}^{n-1} \nu(I^r)}$ is
 of polynomial type of degree $d$ with leading (binomial)
 coefficient $f_{0}(I)$. Also , for $n\gg 0$, ${\ds
 \l(\RR/I)^n-\l(\RR/(I,h)^n)}$  is the difference of two polynomials of
 degree $d$  and leading (binomial) coefficients $\rme_0(I)$ and $\rme_0(L)$. Hence
\[ e_0(I) - e_0(L)  \leq \l(\RR/I:L)\cdot f_0(I).\]
\end{proof}

\begin{Theorem}[Variation of the Chern number] \index{variation of
the Chern number}
\label{e1h}
 Let $(\RR, \m )$ be a Noetherian local ring and let $I\subset
L=(I,h)$ be $\m$--primary ideals. If $h$ is integral over $I$, then
\[\rme_1(L) -\rme_1(I) \leq \mbox{\rm red}_I(L)\cdot \l(\RR/I:L)\cdot
f_0(I).\]
\end{Theorem}

\begin{proof} Let $s$ be the reduction number of $L$ with respect to $I$.
Then \[{\ds h^{s+1}\in I(I,h)^s. }\] Then we obtain the following
filtration : \[\begin{array}{lll} I^{n}=M_{0} \subset M_{1}=(M_{0},\;
I^{n-1}h) \subset \cdots &\subset & M_{r}=(M_{r-1},\; I^{n-r}h^{r})
\\ && \\ &\subset & \cdots \subset M_{s}=(M_{s-1},\; I^{n-s}h^s)=L^n .
\end{array}
\] Therefore \[{\ds \l(\RR/I^n)-\l(\RR/(I,h)^n) = \l(L^{n}/I^{n}) =
\sum_{r=1}^{s} \l(M_{r}/M_{r-1}) \leq \l(\RR/I:h) \sum_{r=1}^s
\nu(I^{n-r}). }\]

Now for $n\gg 0$, $\l(\RR/I^n)-\l(\RR/(I,h)^n)$ is the difference of two
polynomials of degree $d$ and with same leading (binomial)
coefficients $\rme_0(I)$ and $\rme_0(L)$, therefore it is at most a
polynomial of degree $d-1$ and leading coefficient $\rme_1(L)-e_1(I)$.
On the other hand, for $n\gg 0$, we have \[{\ds \l(\RR/I:h)
\sum_{r=1}^s \nu(I^{n-r}) \leq \l(\RR/I:h) \cdot s \cdot
\sum_{i=0}^{d-1} (-1)^{i} f_{i}(I){{n+d-i-1}\choose{d-i-1}} }\] This proves
that \[\rme_1(L)-e_1(I) \leq \mbox{\rm red}_I(L)\cdot \l(\RR/I:L)\cdot
f_0(I). \] 
\end{proof}

\begin{Corollary}
 Let $(\RR, \m )$ be a Noetherian local ring. Let $I\subset L=(I,h)$ be $\m$--primary ideals such that $I$ is a minimal reduction of $L$. Then
\[\rme_1(L)  \leq \mbox{\rm red}_I(L) \cdot \l(\RR/I:L) .\]
\end{Corollary}

\begin{Example}\label{e0Ih}{\rm (\cite[Example 7.36]{icbook})
Let $\RR=k[x,y,z]_{(x,y,z)}$ and let $I$ and $L$ be $\RR$--ideals such
that \[{\ds I=(x^{a},\; y^{b},\; z^{c}) \subset (I,\; x^{\alpha}
y^{\beta} z^{\gamma})=L, }\] where ${\ds \frac{\alpha}{a} +
\frac{\beta}{b} + \frac{\gamma}{c} <1}$. This inequality condition
ensures that $h=x^{\alpha} y^{\beta} z^{\gamma} \notin \cl{I}$. Then
we have \[e_{0}(I) - e_{0}(L) = abc - (ab \gamma +bc \alpha + ac
\beta) = {\ds abc \left(1- \frac{\alpha}{a} - \frac{\beta}{b} -
\frac{\gamma}{c} \right).} \] Since ${\ds (I : L) = (I: x^{\alpha}
y^{\beta} z^{\gamma}) =(x^{a-\alpha},\; y^{b-\beta},\; z^{c-\gamma})
}$, we obtain \[\begin{array}{lll} \l(\RR/(I:L))f_{0}(I) &= &
(a-\alpha)(b-\beta)(c-\gamma) \\ && \\ &= & abc -bc \alpha -ac \beta
- ab \gamma + a \beta \gamma + b \alpha \gamma + c \alpha \beta -
\alpha \beta \gamma \\ && \\ &=& e_{0}(I) - e_{0}(L) + \alpha \beta
\gamma {\ds \left( \frac{a}{\alpha}+ \frac{b}{\beta} +
\frac{c}{\gamma} -1 \right)} \\ && \\ & > & e_{0}(I) - e_{0}(L)
\end{array}
 \] Let ${\ds J=(x^{a}-z^{c},\; y^{b}-z^{c},\;x^{\alpha} y^{\beta}
z^{\gamma} ) }$ and suppose that ${\ds a > 3 \alpha,\; b > 3 \beta,\;
c > 3 \gamma }$. Note that $L=(J,\; z^c)$.  Then $J$ is a minimal
reduction of $L$ and the reduction number $\mbox{\rm red}_{J}(L) \leq
2$.  We can estimate $\rme_{1}(L)$ : \[ e_{1}(L)=e_1(L)-e_1(J) \leq 2
\l(\RR/J:L) . \] }\end{Example}

We resume the discussion of variation of Hilbert coefficients.

\begin{Proposition}\label{f0h}
Let $L=(I,h)$ and suppose $h$ is integral over $I$. Then
\[ f_0(L) \leq (1+ \mbox{\rm red}_I(L) \cdot \l(\RR/I:h))f_0(I) .\]
\end{Proposition}

\begin{proof} Let $s$ be the reduction number of $L$ with respect to $I$.  By
tensoring the following exact sequence with $\RR/\m$ \[ 0 \lar I^n \lar
L^n \lar L^n/I^n \lar 0,\] we obtain
 \[I^n/\m I^n \lar L^n/\m L^n
\lar (L^n/I^n) \otimes \RR/\m \rar 0. \] Therefore
\begin{eqnarray*}
  \l(L^n/\m L^n)
-\l(I^n/\m I^n) &\leq &\l((L^n/I^n) \otimes \RR/\m) \\
&\leq & \l(L^n/I^n) =\l(\RR/I^n)-\l(\RR/L^n). \end{eqnarray*}
 Hence we have the
inequalities of the leading coefficients (in degree $\dim(\RR)-1$)
\[f_0(L) - f_{0}(I) \leq\rme_1(L) -\rme_1(I). \] Using
Theorem~\ref{e1h}, we obtain \[ f_0(L) - f_{0}(I) \leq\rme_1(L) -
e_1(I) \leq s \cdot \l(\RR/I:h) \cdot f_{0}(I), \] which completes the
proof.  
\end{proof}

\begin{Remark}{\rm
Note that the formulas for the variations of $\rme_1(I)$ and $f_0(I)$
require that the new ideal has the same integral closure as $I$.
}\end{Remark}

Write $L= (I, H)$, where $\nu(H) = \nu(L/I)$, and note that $I:L =
I:H$.  Write the difference of Hilbert functions using the obvious
filtration
\begin{eqnarray*}
\l(\RR/I^n)-\l(\RR/L^n) &=&\l((I,H)^n/I^n)= \l((I^n,HI^{n-1}, \ldots, H^{n-1}I,H^n)/I^n) \\
&= & \sum_{r=1}^{n}  \l(M_r/M_{r-1}),
\end{eqnarray*}
where $M_r=(I^n,HI^{n-1}, \ldots, H^{r-1}I^{n-r+1}, H^{r}I^{n-r})$.
Note that $M_r/M_{r-1}$ is generated by the images of $H^rI^{n-r}$.
More precisely, if $L= (h_1, \ldots, h_s, I)$, then $M_r/M_{r-1}$ is
generated by batches of elements, difficult to control (not entirely
the case for the estimation of $\rme_1$ when $L$ is integral over $I$).
This filtration has been used by several authors when $I$ is a system
of parameters. Our formulation wraps it differently to accommodate our
data.

\begin{Theorem}\label{e1hs} Let $(\RR, \m)$ be a Noetherian local ring of
dimension $d$, let $I$ be an $\m$-primary ideal and let $L=(h_1,
\ldots, h_m, I)$ be integral over $I$ of reduction number $s$. Then
\[\rme_1(L)-e_1(I) \leq \l(\RR/I:L)\cdot \left[{{m+s}\choose
{s}}-1\right]\cdot f_0(I).\]
\end{Theorem}

\begin{proof} We have already given parts of the proof. The remaining part is
to estimate the growth of the length of ${M_r/M_{r-1}} = (h_1,
\ldots,h_m)^r I^{n-r}+ M_{r-1}/M_{r-1} $.  We note that this module
is annihilated by $I:L$ and is generated by the `monomials' in the
$h_i$ of degree $r$, with coefficients in $I^{n-r}$.  There is a
natural surjection \[ \Phi: \RR/(I:L) \otimes \RR^{b_{r}} \otimes I^{n-r}
\lar M_r/M_{r-1},\] where $\Phi(\cl{r} \otimes e_{i} \otimes x)= r
\alpha_{i}x + M_{r-1}$ and ${\ds b_r = {{m+r-1}\choose{r}} }$.
Therefore for $n>>0$, \[\begin{array}{lll} \l(\RR/I^n)-\l(\RR/L^n) &=& {\ds
\sum_{r=1}^{s} \l(M_r/M_{r-1}) } \\ && \\ &\leq & {\ds \sum_{r=1}^{s}
\l(\RR/I:L) \nu(I^n) {{m+r-1}\choose{r}} } \\ && \\ &=& {\ds \l(\RR/I:L)
\nu(I^n) \left( {{m+s}\choose{s}} -1 \right), }
\end{array}
\] which completes the proof.
\end{proof}

\begin{Corollary}
Let $(\RR, \m)$ be a Noetherian local ring of dimension $d>0$. Let
$L=(I,h)$ be an $\RR$-ideal.
\begin{enumerate}
\item[{\rm (a)}] If $h$ is integral over $I$ and $\red_I(L)=s$, then
\[\rme_1(L)-e_1(I) \leq \l(\RR/I:h)\cdot
s\cdot f_0(I).\]

\item[{\rm (b)}] Furthermore, if $\RR$ is a
Gorenstein ring and $I$ is a parameter ideal, then
\[\rme_1(L)\leq \l(\RR/I:h)\cdot s= [e_0(L)-\l(\RR/L)]\cdot f_0(L).
\]
\end{enumerate}

\end{Corollary}

\begin{proof} (b) We have only to observe that \[
\l(\RR/I:h)=\l(\RR/I)-\l(I:h/I)= e_0(I)-\l(I:h/I)=e_0(L)-\l(\RR/L),\]
since
$h$ is integral over $I$ and
 $I:h/I$ is the canonical module of $\RR/L$.
\end{proof}

The values of $\rme_1$  are also related to the multiplicity
of certain Sally modules even in the non Cohen--Macaulay case, according to \cite{C09}.
Let $(\RR, \m)$ be a Noetherian local ring of dimension $d \geq 1$ with
infinite residue field. Let $I$ be an $\m$--primary ideal and $J$ a
minimal reduction of $I$. Suppose that $\dim(S_{J}(I))=d$ and that
$\H^{0}_{\m}(\RR) \subset I$, then the multiplicity $s_0$ of the Sally
module $S_{I}(I)$ is \[s_0(J,I) =\rme_1(I)-e_0(I)-e_1(J)+ \l(\RR/I).  \]

\begin{Corollary}\label{Sally}
Let $(\RR, \m)$ be a Noetherian local ring of dimension $d \geq 1$ with
infinite residue field. Let $I$ be an $\m$--primary ideal and $J$ a
minimal reduction of $I$. Suppose that $\dim(S_{J}(I))=d$ and that
$\H^{0}_{\m}(\RR) \subset I$, then the multiplicity $s_0$ of the Sally
module $S_{I}(I)$ satisfies \[s_0(J,I) \leq -e_{0}(I) + \l(\RR/I) +
\l(\RR/J:I)\left[{{\nu(I) -d +r}\choose {r}}-1\right], \] where $r$ is
the reduction number.
\end{Corollary}

\begin{Example}{\rm  Let $\RR=k[x,y]_{\m}$, where $\m=(x,y)$, and let $L=(a_1, \ldots, a_{n}, a_{n+1})$ be an $\m$--primary $R$--ideal with a resolution:
\[ 0 \lar \RR^{n} \stackrel{\varphi}{\lar} \RR^{n+1} \lar L \lar 0, \]
where every entry of $\varphi$ is linear and $n \geq 2$. Note that
$L=\m^n=\cl{L}$.
%%%%%%%%%%%%%%%%%%%%%%%%%%%
%% $L = \m^{n}$
%%%%%%%%%%%%%%%%%%%%%%%%%
We may assume that $J=(a_1, a_2)$ is a minimal reduction of $L$.
Let $I=(a_1, \ldots, a_{n})$. In particular, $L=\cl{L}=\cl{J} \subset \cl{I} \subset \cl{L}=L$ so that $L$ is integral over $I$.
Then since each $a_i$ is a $n$--form, the presentation for $I$ is of the form:
\[ 0 \lar \RR^{n-1} \stackrel{\psi}{\lar} \RR^n \lar I \lar 0,   \]
where $\psi = \left[ \psi_1 \;\; \psi_2 \;\; \cdots \;\; \psi_{n-1}   \right]$ with the first column $\psi_1$ consisting of $2$-forms and all the other columns $\psi_2, \ldots, \psi_{n-1}$ consisting of linear forms.
Using ${\ds e_{1}(L)=e_{1}(\m^n) = \frac{1}{2}n(n-1)}$, and
\[{\ds e_{1}(L) - e_{1}(I) \leq  \mbox{\rm red}_{I}(L) \l(\RR/I:a_{n+1}) f_{0}(I),   }\]
we obtain
\[{ \ds\rme_1(I) \geq  \frac{1}{2}n(n-1) -  \mbox{\rm red}_{I}(L) \l(\RR/I:a_{n+1}) f_{0}(I) }    \]
}\end{Example}

One situation that may be amenable to further analysis is when $L= I: \m$, or
more generally $L=I:m^s$ for some values of $s$. We refer to $L$ as a
{\em socle extension} of $I$.

\begin{Example}{\rm
Let $(\RR, \m)$ be a Cohen-Macaulay local ring of dimension $d$.
Suppose $\RR$ is not a regular local ring. Let $I=(x_1, \ldots, x_d)$ be a parameter
ideal and $L = I:\m= (I, h_1, \ldots, h_m)$. According to \cite{CPV1}
(see also \cite[Theorem 1.106]{icbook}), $L^2=IL$. In particular
$\rme_0(I)=e_0(L)$. Using Theorem~\ref{e1hs}, we estimate $\rme_1(L)$:
%%%%%%%%%%%%%%%%%%%%%%%%%%%%%%%%%%%%%%%%%%%%%%%%%%%%%%%%%%%%%%%%%%%%%%%%%%%%%%%%%%%%%%%%%%%%%%%%%%
%%% consider \[ L^n/I^n = (I^n + HI^{n-1})/I^n\] from which we get
%%%%%%%%%%%%%%%%%%%%%%%%%%%%%%%%%%%%%%%%%%%%%%%%%%%%%%%%%%%%%%%%%%%%%%%%%%%%%%%%%%%%%%%%%%%%%%%%%%%
\[ e_{1}(L)=\rme_1(L)-e_1(I) \leq \l(\RR/I:L) \cdot \left[{{m+1}\choose
{1}}-1\right] \cdot f_0(I)= m. \] As for $f_0(L)$, we have ${\ds
f_0(L) \leq 1 + m}$.  }\end{Example}

These elementary inequalities should be accompanied by an analysis of the
possible equality. Let us argue that in the Cohen-Macaulay case this
is very rare.

\subsection{Special fiber and multiplicity}

We seek relationships between the Hilbert coefficients of $\gr(\BB)$
of their special fibers. We first recall the
  formula of C. Lech (\cite{Lech})
that gives rise to {\it a priori} bounds for $e(I)$.

\begin{Theorem} \label{Lechformula} If $(\RR, \mathfrak{m})$ is a
Noetherian local ring of dimension $d$ and  $I$ is an
$\mathfrak{m}$-primary ideal, then
\begin{eqnarray} \label{Lech}
%\psboxit{box 1.0 setgray fill}{\fbox{$
%\begin{array}{c}
%\ \\
\  e_0(I) \leq d!\cdot \deg(\RR)\cdot \lambda(\RR/\overline{I})%  \\
%\ \\
%\end{array}  $}}
\end{eqnarray}
where $\bar{I}$ is the integral closure of $I$.
\end{Theorem}

\begin{Proposition}
Let $(\RR,\mathfrak{m})$ be a normal local domain and let $I$ be an
$\mathfrak{m}$-primary ideal. Suppose that $\CC= \overline{\RR[It]}$ is
finite over $\RR[It]$. We denote by
$f_0(I)$ the multiplicity of the special fiber of $I$, and set
$\overline{f}_0(I)
= \deg(\CC/\mathfrak{m}\CC)$.
One has
\[e_0(I) \leq \min\{f_0(I)\cdot \lambda(\RR/I), \overline{f}_0(I)\cdot
\lambda(\RR/\overline{I})\}.
\]
\end{Proposition}

\begin{proof} We first observe that $C_{n+1}= IC_n= \overline{I}C_n$, for $n\gg
0$. In particular, in that range, $C_{n+1}\subset \mathfrak{m}C_n$.
Consider now the corresponding exact sequence
\[ 0
\rar \mathfrak{m}C_n/C_{n+1} \lar C_n/C_{n+1} \lar C_n/\mathfrak{m}
C_n \rar 0.
\]
Counting multiplicities, we have
\[\deg(\gr(\CC))= e_0(I) \leq \deg(\mathfrak{m}\CC/\overline{I}\CC) +
\deg(\CC/\mathfrak{m}\CC)\leq\overline{f}_0(I)(
\lambda(\mathfrak{m}/\overline{I})+1)
, \]
as  desired.

The other inequality, $\rme_0(I)\leq f_0(I)\cdot \lambda(\RR/I)$, has a
similar proof. 
\end{proof}

\subsection{Ideals generated by quadrics}

\chapter{Appendix}

\section*{Introduction}

\section{Approximation Complexes}

\subsection{Introduction} These are complexes derived from Koszul
complexes, and arise as follows (for details, see \cite{HSV3},
\cite{HSV3II}, \cite{HSV1}, \cite[Chapter 4]{alt}).

\bigskip

\subsubsection*{The $\mathcal{Z}$ and $\mathcal{M}$
complexes}\index{approximation complex}

Let $\RR$ be a commutative ring, $F$ a free $\RR$-module of rank $n$,
with a basis  $\{e_1, \ldots, e_n\}$,
 and
$\varphi: F\rar \RR$ a homomorphism. The exterior algebra $\bigwedge F$ of
$F$ can be endowed with a differential
\[ \partial_{\varphi} : \bigwedge^{r} F
{\lar} \bigwedge^{r-1}F,
\]
\[ \partial_{\varphi}(v_1\wedge v_2\wedge \cdots \wedge v_{r}) =
\sum_{i=1}^{r} (-1)^i\varphi(v_i) (v_1\wedge \cdots \wedge
\widehat{v_i} \wedge  \cdots \wedge v_{r}).  \]

The complex $\bbk(\varphi)=\{ \bigwedge F, \partial_{\varphi}\}$ is called the {\em
Koszul complex} of $\varphi$. \index{Koszul complex} Another notation
for it is: Let $\xx=\{\varphi(e_1), \ldots, \varphi(e_n)\}$, denote
the Koszul complex by $\bbk(\xx)$.

\medskip

Let $\SS = S(F)= \Sym(F) = \RR[\TT_1, \ldots, \TT_n]$, and consider the
exterior algebra of $F\otimes_{\RR} S(F)$. It can be viewed as a Koszul
complex obtained from $\{\bigwedge F, \partial_{\varphi}\}$ by change
of scalars $\RR \rar \SS$, and another complex defined by the
$\SS$-homomorphism
\[ \psi: F \otimes_{\RR}S(F) \lar S(F), \quad \psi(e_i)= \TT_i.
\]
The two differentials $\partial_{\varphi}$ and $\partial_{\psi}$
satisfy
\[ \partial_{\varphi}\partial_{\psi}+
\partial_{\psi}\partial_{\varphi}=0,\]
which leads directly to the construction of several
complexes.

\subsection{Construction}

\begin{Definition}{\rm  Let $\ZZ$, $\BB$ and $\HH$ be the modules of
cycles, boundaries  and
 the homology of $\varphi$. 
\begin{itemize}
\item[{\rm (a)}] The $\mathcal{B}$-complex of $\varphi$ is
$\mathcal{B}= \{\BB\otimes_{\RR} \SS, \partial\}$;  \index{approximation complex; 
$\mathcal{B}$-complex}

\item[{\rm (b)}] The $\mathcal{Z}$-complex of $\varphi$ is
$ \mathcal{Z}=\{ \ZZ\otimes_{\RR} \SS, \partial\}$; \index{approximation
complex; $\mathcal{Z}$--complex}

\item[{\rm (c)}] The $\mathcal{M}$-complex of $\varphi$ is
$ \mathcal{M}=\{ \HH\otimes_{\RR} \SS, \partial\}$,
\end{itemize}
where $\partial$ is the differential induced by $\partial_{\psi}$.
}\end{Definition}

These are complexes of graded modules over the polynomial ring
$\SS=\RR[\TT_1, \ldots, \TT_n]$.

\begin{Proposition} Let $I = \varphi(F)$. Then
\begin{enumerate}
\item[{\rm (a)}]  The homology of
$\mathcal{Z}$  and of
$\mathcal{M}$
depend only on $I$;

\item[{\rm (b)}] $ H_0(\mathcal{Z})= \Sym(I)$;

\item[{\rm (c)}] $  H_0(\mathcal{M})= \Sym(I/I^2)$. 
\end{enumerate}
\end{Proposition}

The two complexes are closely related:

\begin{Proposition}[{\cite[Proposition 4.4]{HSV1}}] For each positive integer $m$ there exists an
exact sequence of $\RR$-modules
\[ \cdots \lar H_{r}(\mathcal{M}_{m+1}) \lar H_{r}(\mathcal{Z}_{m})
\lar H_{r}(\mathcal{M}_{m}) \lar H_{r-1}(\mathcal{Z}_{m+1}).\]
\end{Proposition}

\begin{Corollary} The following conditions are equivalent for the
complexes $\mathcal{M}$ and $\mathcal{Z}$:
\begin{itemize}
\item[{\rm (a)}] $\mathcal{M}$ is acyclic;
\item[{\rm (b)}] $\mathcal{Z}$ is acyclic and $\Sym(I)=\RR[It]$.
\end{itemize}
\end{Corollary}

If the natural surjection $\Sym(I)\rar \RR[It]$ is an isomorphism, $I$
is said to be of {\em linear type}\index{ideal of linear type}. The
following assertion on the local number of generators  is a requirement:

\begin{Proposition} If the $\RR$-ideal $I$ is of linear type, then
\[ \nu(I_{\mathfrak{p}})\leq \height \mathfrak{p}, \quad  \mbox{for
every \ }
\mathfrak{p}\supset I.
\]
\end{Proposition}

\subsubsection*{Acyclicity}
The homology of the Koszul complex $\bbk(\varphi)$ is not fully
independent of $I$, for instance, it matters the number of
generators. More distinctive are the nature of the ideals for
which
  $\mathcal{Z}$ and  $\mathcal{M}$ are acyclic.

Part of the significance of $d$-sequences lies in the following
result (\cite{Hu1a}, \cite{Valla1}):

\begin{Theorem} Let $I$ be an ideal generated by a $d$-sequence. Then
the symmetric and  Rees algebras of $I$ are isomorphic,
\[  \Sym(I)\simeq  \RR[It].\]
\end{Theorem}

\begin{Theorem}[Acyclicity Theorem]  Let $(\RR, \mathfrak{m})$ be a Noetherian local ring of
infinite residue field, and let $I$ be an $\RR$-ideal.
\begin{enumerate}
\item[{\rm (a)}] $I$ is generated by a $d$-sequence
 if and only if the approximation complex $\mathcal{M}$ is acyclic.

\item[{\rm (b)}] $I$ is generated by a  proper sequence
 if and only if the approximation complex $\mathcal{Z}$ is acyclic.
\end{enumerate}
\end{Theorem}

\subsubsection*{Koszul homology} Let $I$ be an ideal generated by the
sequence $\xx=\{x_1, \ldots, x_n\}$. The corresponding approximation
complexes are defined in terms of the modules of the cycles and the
homology of $\bbk(\xx)$:
\begin{eqnarray*}
\mathcal{Z}: \quad
0 \rar Z_n\otimes_\RR\SS[-n] \rar Z_{n-1}\otimes_\RR\SS[-n+1] \rar
\cdots \rar Z_1\otimes_\RR\SS[-1] \rar \SS \rar 0
\end{eqnarray*}
\begin{eqnarray*}
\mathcal{M}: \quad
0 \rar \H_n\otimes_\RR\SS[-n] \rar \H_{n-1}\otimes_\RR\SS[-n+1] \rar
\cdots \rar \H_1\otimes_\RR\SS[-1] \rar \H_0\otimes_\RR \SS \rar 0.
\end{eqnarray*}
There is some resemblance to projective resolutions of
$\H_0(\mathcal{Z})= \Sym(I)$ as a $\SS$-module, and of
$\H_0(\mathcal{M})= \Sym(I/I^2)$ as a $\SS/I\SS$-module. This is
reinforced when the complexes are acyclic and the components
$Z_i\otimes_\RR\SS$ and $\H_i\otimes_\RR\SS$ have high depths as $\SS$-modules.

\medskip

A source of large classes of ideals where the $H_i$ (and consequently
the modules $Z_i $ of cycles) have large depths is given by the
following (\cite{Hu5}):

\begin{Theorem} Let $\RR$ be a Gorenstein local ring. If the ideal $I$
lies in the linkage class of a complete intersection then the Koszul
homology modules $\H_i(I)$ are Cohen-Macaulay of dimension $\dim \RR/I$.
\end{Theorem}

\subsubsection*{Cohen-Macaulay algebras}
Some of the applications of the complexes $\mathcal{Z}$ and
$\mathcal{M}$ make use of information on the Koszul homology of the
ideal.

\begin{Theorem}[{\cite[Theorem 10.1]{HSV1}}] Let $\RR$ be a Cohen-Macaulay local ring and let $I$ be an
ideal of positive height. Assume:
\begin{itemize}
\item[$\bullet$] $\nu(I_{\mathfrak{p}})\leq \height \mathfrak{p}+1$ for
every \ $\mathfrak{p} \supset I$;

\item[$\bullet$] $\depth (\H_i)_{\mathfrak{p}} \geq \height
\mathfrak{p}-\nu(I_{\mathfrak{p}})+1 $ for every \ $\mathfrak{p}\supset
I$ and every $0\leq i \leq \nu(I_{\mathfrak{p}})-\height
I_{\mathfrak{p}} $.
\end{itemize}
Then
\begin{itemize}
\item[{\rm (a)}] The complex $\mathcal{Z}$ is acyclic.
\item[{\rm (b)}] $\Sym(I)$ is a  Cohen-Macaulay ring.
\end{itemize}
\end{Theorem}

\begin{Corollary}\label{Zaci} Let $\RR$ be a Cohen-Macaulay local and let $I$ be an
almost complete intersection.
The complex $\mathcal{Z}$ is acyclic and $\Sym(I)$ is a
Cohen-Macaulay algebra in the following cases:
\begin{itemize}
\item[{\rm (a)}] $I$ is $\m$-primary.
\item[{\rm (b)}]  $\height I= d-1$.
Furthermore if $I$ is
generically a complete intersection then $I$ is of linear type.
\item[{\rm (c)}]  $\height I= d-2$ and
 $\depth \RR/I\geq 1$. Furthermore if $\nu(I_{\mathfrak{p}})\leq
  \height \mathfrak{p}$ for $I\subset \mathfrak{p}$ then $I$ is of
  linear type.
\end{itemize}
\end{Corollary}

\begin{Theorem}[{\cite[Theorem 9.1]{HSV1}}] \label{HSV19.1} Let $\RR$ be a Cohen-Macaulay local ring and let $I$ be an
ideal of positive height. Assume:
\begin{itemize}
\item[$\bullet$] $\nu(I_{\mathfrak{p}})\leq \height \mathfrak{p}$ for
every \ $\mathfrak{p} \supset I$;

\item[$\bullet$] $\depth (\H_i)_{\mathfrak{p}} \geq \height
\mathfrak{p}-\nu(I_{\mathfrak{p}})+1 $ for every \ $\mathfrak{p}\supset
I$ and every $0\leq i \leq \nu(I_{\mathfrak{p}})-\height
I_{\mathfrak{p}} $.
\end{itemize}
Then
\begin{itemize}
\item[{\rm (a)}] The complex $\mathcal{M}$ is acyclic.
\item[{\rm (b)}] $\Sym(I)\simeq  \RR[It]$ and $\Sym(I/I^2) \simeq
\gr_I(\RR) $.
\item[{\rm (c)}] $\Sym(I)$ and $\Sym(I/I^2)$ are Cohen-Macaulay rings.
\item[{\rm (d)}] If it is further the case that $\RR$ is a Gorenstein ring
and the $\H_i$ are Cohen-Macaulay, then $\Sym(I/I^2)$ is a Gorenstein
ring.
\end{itemize}
\end{Theorem}

\begin{Corollary} Let $\RR$ be a Gorenstein local ring of dimension $d$
and $I$ an almost complete intersection of codimension $g$  such that
$\depth \H_1(I)\geq d-g-1$. If $I$ is syzygetic then $I$ is of linear
type and $\gr_I(\RR)$ is Cohen-Macaulay.
\end{Corollary}

\begin{proof} By Corollary~\ref{Zaci}, $\mathcal{Z}$ is acyclic.  To verify
that $\mathcal{M}$ is acyclic, by Theorem~\ref{HSV19.1} it suffices
to check that for every minimal prime ideal $\mathfrak{p}$ of $I$,
$I_{\mathfrak{p}}$ is a complete intersection. This follows from the
exact sequence
\[ 0 \rar \delta(I)=0 \rar \H_1(I_{\mathfrak{p}})\lar
(\RR/I)_{\mathfrak{p}}^{g+1} \lar (I/I^{2})_{\mathfrak{p}} \rar 0.\]
Since $\H_1(I_{\mathfrak{p}})$ is an injective
$(\RR/I)_{\mathfrak{p}}$-module, the sequence splits and
$I_{\mathfrak{p}}$ is generated by $g$ elements. 
\end{proof}

\subsection{Approximation complexes with coefficients} \index{approximation complexes 
with coefficients} 

There is another family of approximation complexes built to uncover properties of one
other module introduced by David Rees. Let $\RR$ be a commutative ring, $I$ an ideal and $M$ an $\RR$-module. The Rees module of $I$ relative to $M$ is
\[ \mathcal{R}(I;M) = \bigoplus_{n\geq 0} I^nM.\]  
It is a construction that plays a strikingly clarifying role in the current standard proof of the Artin-Rees Lemma. $\mathcal{R}(I;M)$ is also the ancestor of the associated graded
module of $I$ relative to $M$
\[ \mathcal{R}(I;M)\otimes \RR/I = \gr_I(M) = \bigoplus_{n\geq 0} I^nM/I^{n+1}M.\] 
It is important to derive homological properties of these modules particularly for
the role $\gr_I(M)$ has in various theories of multiplicity.

\bigskip

The construction of the complexes follows the script above. Let $\xx=\{x_1, \ldots, x_n\}$ be a sequence generating the 
ideal $I$ and denote by 
\[ \{\mathbb{K}(\xx;M), \partial\}  =\{\mathbb{K}(\xx;\RR)\otimes_{\RR} M, \partial\}
\]
 the corresponding  Koszul complex. Set $B(\xx;M)$, $Z(\xx;M)$ and 
$\H(\xx;M)$ for their modules of boundaries, cycles and homologies.  

\medskip

Attaching coefficients  we obtain the complex
\[ \{\mathbb{K}(\xx;M)\otimes \SS, \partial\otimes \SS\}= \mathbb{K}(\xx; M\otimes \SS), \partial\otimes \SS .\]
The module $\mathbb{K}(\xx;M\otimes \SS)$ supports another differential $\partial'$
defined by the sequence $\TT = \{\TT_1, \ldots, \TT_n\}$. As previously $\partial$ and $\partial'$
(anti-) commute which permits the definitions of the complexes ($\H_i = \H_i(\xx;M)$)

\begin{eqnarray*}
\mathcal{Z}(I;M): \quad
0 \rar Z_n\otimes\SS[-n] \rar Z_{n-1}\otimes\SS[-n+1] \rar
\cdots \rar Z_1\otimes\SS[-1] \rar Z_0\otimes \SS \rar 0
\end{eqnarray*}
\begin{eqnarray*}
\mathcal{M}(I;M): \quad
0 \rar \H_n\otimes\SS[-n] \rar \H_{n-1}\otimes\SS[-n+1] \rar
\cdots \rar \H_1\otimes\SS[-1] \rar \H_0\otimes \SS \rar 0.
\end{eqnarray*}

When these complexes are acyclic, $\H_0(\mathcal{Z}(I;M)) \simeq \mathcal{R}(I;M)$ and
$\H_0(\mathcal{M}(I;M)) \simeq \gr_I(M)$, in which case a great deal of information passes between
the Koszul homology modules $\H_i(\xx;M)$ and the Rees modules.

\begin{Proposition}[{\cite[Proposition 3.8]{HSV3II}}] If the complex $\mathcal{M}(I;M)$ is acyclic then 
\[ \mbox{\rm $I$-depth $M = \SS_{+}$-depth $\gr_I(M)$}.\]
\end{Proposition}

The sequential criterion of acyclicity is similar to the case $M =\RR$:

\begin{Theorem}[{\cite[Theorem 4.1]{HSV3II}}] \label{HSV3II4.1} Let $(\RR, \m)$
 be a Noetherian local ring with infinite residue field. Let
  $I$ be an
ideal and let $M$ be a finitely generated $\RR$-module. The following conditions are
equivalent: 
\begin{itemize}
\item[{\rm (a)}] The complex $\mathcal{M}(I;M)$ is acyclic;
 
\item[{\rm (b)}] $I$ is generated by a $d$-sequence relative to $M$.
\end{itemize}
\end{Theorem}

\section{Elimination}

\subsection{Introduction}
Let $\RR$ be a Noetherian ring, usually $\RR=k[x_1, \ldots, x_n]$ for
$k$ a field,
 and let $I$ be an
 ideal, and denote by $\LL
$ its ideal of equations
\begin{eqnarray*}
 0 \rar \LL \lar \SS = \RR[\TT_1, \ldots, \TT_m] \stackrel{\psi}{\lar}
\RR[It] \rar 0,  \quad \TT_i \mapsto f_it .
\end{eqnarray*}

\subsubsection*{Sylvester forms}\index{Sylvester form}

 Let
$\mathbf{f}= \{f_1, \ldots, f_m\} $ be a set of polynomials in
 $\LL\subset \SS=\RR[\TT_1, \ldots , \TT_n]$ and let $\mathbf{a}= \{a_1, \ldots,
a_m\}\subset \RR $. If  $f_i \in (\mathbf{a}) \SS$ for all $i$, we can
write
\[ \mathbf{f}= [f_1 \cdots f_m] = [a_1 \cdots a_m] \cdot\AA =
\mathbf{a}\cdot \AA,\] where $\AA $ is an $m\times m$ matrix with entries
in $\SS$. We call $(\aa)$ a $\RR$-content of $\mathbf{f}$.
 Since  $\aa\not \subset L$, then $\det(\AA)\in \LL$.

  By an abuse of terminology\index{Sylvester form}, we refer to $\det(\AA)$ as a
{\em
Sylvester form}, or the {\em Jacobian
dual}\index{Jacobian dual},  of $\mathbf{f}$ relative to $\mathbf{a}$, in
notation
\[ \det(\mathbf{f})_{(\mathbf{a})}=\det(\AA).\]

\begin{Proposition} \label{Sylvester} Let $\RR=k[x_1, \ldots, x_n]$, $\m=(x_1, \ldots,
x_d)$,
and $\ff_1, \ldots, \ff_s$ be forms in $\SS=\RR[\TT_1, \ldots,
\TT_m]$. Suppose  the $\RR$-content of the $\ff_i$
(the ideal
generated by the coefficients in $\RR$) is generated by forms $\aa_1,
\ldots, \aa_q$ of the same degree. Let
\[ \ff = [\ff_1, \ldots, \ff_s]= [\aa_1, \ldots, \aa_q]\cdot \AA=
\aa\cdot \AA
\] be  the corresponding Sylvester decomposition.
If $s\leq q$ and the {first-order syzygies of $\ff$ have coefficients in
$\m$}, then $I_s(\AA)\neq 0$.
\end{Proposition}

 Note that $I_s(\AA)$ is generated by forms of $k[\TT_1, \ldots,
\TT_m]$.
 When we apply it to forms in $(L_1): \mathfrak(c)$,
$\mathfrak{c}\subset \RR$, $I_s(\AA)\subset \LL$.
\bigskip

\begin{proof} The condition
 $I_p(\AA)=0$ means the columns of $\AA$ are linearly dependent
over $k[\TT_1, \ldots, \TT_m]$, thus for some nonzero column vector
$\cc\in k[\TT_1, \ldots,\TT_m]^s $
\[ \AA\cdot \cc = 0.\]
Therefore \[ \ff\cdot \cc= \aa\cdot \AA \cdot \cc=
0\]
is a syzygy of the $\ff_i$ whose content in not in $\m$, against
the assumption.
\end{proof}

We refer to the condition on the content of syzygies as a
{\em syzygetic certificate}\index{syzygetic certificate}
 for $I_s(\AA)\neq 0$. (We use the word as
a quick to check hypothesis.)
 This particular
certificate will be called the {\em Sylvester certificate}.
The following is just a convenient formulation for the nonvanishing
of the special Sylvester determinant we will make use of.

\begin{Proposition} \label{ppower}
 Let $I$ be a homogeneous almost complete
intersection  as above.
Let $(\bb)\subset I_1(\varphi)$ and
suppose there are  forms of degree $q$, $\aa=\{\aa_1, \ldots,
\aa_s\}\subset \RR$, such that in bidegree $(q,*)$ we have elements
\[(\hh_1, \ldots, \hh_s)
\subset (L_1): (\bb)\subset (L_1,L_2).\]
 If $(\aa)\cdot (\bb)\subset \m^r$
for some secondary elimination degree $r$,
  then   $\beta =\det \AA $ is
 a power of the elimination equation of $I$.
\end{Proposition}

\begin{proof}
Let $\pp$ denote the elimination equation of $I$.
Since $\pp$ is irreducible it suffices to show that $\beta$ divides
a power of $\pp$.
From $\pp(\aa)\subset (L_1): \bb$,
this gives a representation
\[ \pp[\aa]= [\hh] \AA= [\aa]\BB\cdot \AA,\]
where $\AA$ is $s\times s$ matrix with entries in $k[\TT_1, \ldots,
\TT_m]$, leading to the
equation
\[ [\aa] \big( \BB\cdot \AA-\pp \mathbf{I}\big)=0, \]
where $\mathbf{I}$ is the $s\times s$ identity matrix.

Since the minimal syzygies of $\aa$ have coefficients in
$\m$, we must have
 \[\BB\cdot \AA=\pp \mathbf{I}, \]
so that $\det \BB\cdot \det \AA= \pp^s$, as desired. 
\end{proof}

%\begin{Remark}\rm We note that $\beta=\det\BB$ equals the elimination equation if and only
%if the associated rational map to $I$ is birational onto the image. By the theorem, this is
%the case if and only if $\beta$ is an irreducible polynomial over $k$.
%\end{Remark}

In this construction it is desirable that the ideal $(\mathbf{a})$ have few
generators. Thus, we would like to suggest the use of irreducible
decompositions since the ideals that arise as components are often
complete intersections.
To see how this occurs, we note the following.\index{irreducible
decomposition}

\begin{Theorem} \label{irrdec} Let $(\RR,\mathfrak{m})$ be a Gorenstein local ring and let
$I$ be an $\mathfrak{m}$--primary ideal. Let $J\subset I$ be a
subideal generated by a system of parameters and let $E = (J:I)/J$ be
the canonical module of $\RR/I$. If $E = (e_1, \ldots, e_s$, $\rme_i\neq
0$, and $I_i= \ann(e_i)$, then $I_i$ is an irreducible ideal and
\[ I = \bigcap_{i=1}^{s} I_i.\]
\end{Theorem}

The statement and its proof will apply to ideals of rings of
polynomials over a field.

\bigskip

\begin{proof} The module $E$ is the injective envelope of $\RR/I$, and therefore it is a
faithful $\RR/I$--module (see \cite[Section 3.2]{BH} for relevant
notions). For each $\rme_i$, $\RR e_i$ is a nonzero submodule
of $E $ whose socle is contained in the socle of $E$ (which is
isomorphic to $\RR/\mathfrak{m}$) and therefore
its annihilator $I_i$ (as an $\RR$-ideal) is irreducible. Since the
intersection of the $I_i$ is the annihilator of $E$, the asserted
equality follows.
\end{proof}

A variation is the following:

\begin{Theorem} \label{irrdec2} Let $(\RR,\mathfrak{m})$ be a
Gorenstein local  ring and let
$I$ be an irreducible  Cohen-Macaulay ideal of height $g$. Let $J\subset I$ be a
subideal generated by a regular sequence of $g$ elements and let $E = (J:I)/J$ be
the canonical module of $\RR/I$. Then $E\simeq \RR/I$, that is, $\RR/I$ is
a Gorenstein ring.
\end{Theorem}

\begin{proof} The module $E$ is a
faithful $\RR/I$--module (see again \cite[Section 3.2]{BH} for relevant
notions). Suppose $E$ is minimally generated by $\{e_1, \ldots,
e_{\RR}\}$. Let $\mathfrak{p}$ be  the radical of $I$; since $I$ is
irreducible, it is a primary ideal.
 Let $I_i$ be the  annihilator of $\RR e_i$, so that
$ I = \bigcap_{i=1}^{s} I_i$. The embedding $\RR/I_i
 \hookrightarrow E$ shows that $I_i$ is $\mathfrak{p}$-primary;
 further, localizing $E$ at $\mathfrak{p}$ gives that
 all $I_i$ are irreducible (with the same radical). Thus we
 have that $I$ equals one of the $I_i$, say $I=I_1$ and $I\subset
 I_i$.

For each $\rme_i$, $\RR e_i$ is a nonzero submodule
of $E $ whose socle is contained in the socle of $E$ (which is
isomorphic to $\RR/\mathfrak{m}$) and therefore $I_i$ is irreducible.
\end{proof}

\begin{example}{\rm \index{irreducible decomposition: explicit}
In the case where $I$ is a codimension $2$ ideal with a
free resolution
\[ 0 \rar \RR^{n-1} \stackrel{\varphi}{\lar} \RR^n \lar I \rar 0,\]
\[ \varphi = \left[ \begin{array}{lcr}
& & \\
& \varphi' & \\
& & \\
\hline
a_{n-1,1} & \cdots & a_{n-1,n-1} \\
a_{n,1} & \cdots & a_{n,n-1} \\
\end{array} \right],
\]
where the last two maximal minors $\Delta_{n-1}, \Delta_n$ of $\varphi$
 form a regular sequence,  then
\[  (e_1, \ldots, e_{n-1})=(\Delta_{n-1}, \Delta_n): I = I_{n-2}(\xi') \]
and each ideal $I_i=(\Delta_{n-1}, \Delta_n):e_i$ is a complete
intersection of codimension $2$.

 The fact that the irreducible ideal $I_i$ is a complete
intersection is an observation of Serre (see \cite[Corollary
21.20]{Eisenbudbook}). The explicit decomposition above is that of
 \cite[Proposition 21.24]{Eisenbudbook}.

\medskip

If $\RR=k[s,t]$ and $I$ is an irreducible ideal, it is still true that
$I$ is a complete intersection (see \cite{Fe}). In this case however
there is no quick way to find its two generators (obviously if $I$ is
not homogeneous).

\medskip

Another class of ideals to which these observations apply are
monomial ideals $I$ of polynomial rings over fields, $\RR= k[x_1,
\ldots, x_n]$. Here  in an irreducible decomposition
$ I = \bigcap_{i=1}^{s} I_i$, the $I_i$ can be chosen to be
monomial. It is easy to show that  irreducible monomial ideals
are complete intersections.

}\end{example}

\subsection{Dimension two: an extended
example}\label{dimtwo} \index{elimination: dimension two}

%Here, as a long exercise, we discuss results of \cite{syl}.

Given an ideal $I\subset \RR= k[s,t]$,
generated by  3
 forms of the same degree, from its syzygies
\[0\rar  \RR^2 \stackrel{\varphi}{\lar} \RR^3 \lar I \rar 0,\]
 we extract selected
 collections of
Sylvester forms.
Starting out from  $2$ forms,
  the defining equations
of $\Symi$,
\begin{eqnarray*} \label{genf}
f &=& \alpha_1\TT_1+ \beta_1\TT_2+ \gamma_1\TT_3 \\ \label{geng}
g &=& \alpha_2\TT_1+ \beta_2\TT_2+ \gamma_2\TT_3,
\end{eqnarray*}
following \cite{Cox}, we obtain by elimination
higher degrees forms in the defining ideal of
$\Rees(I)$.

\medskip

To examine the homology of these collections of Sylvester forms,
we  make use of a
 computer-assisted methodology to show that these algorithmically
 specified sets
 generate the  ideal of definition $M$ of $\Rees(I)$ in several cases of
 interest--in particular answering some questions raised
 \cite{Cox}. More precisely, the so-called ideal of moving forms $M$
 is given when $I$ is generated by forms of degree at most $5$. In
 arbitrary degree, the algorithm  will provide the elimination equation in
 significant cases.

\medskip

\subsubsection*{Basic Sylvester forms in dimension $2$}
The initial operation we make use of is the following.
Let $\RR=k[s,t]$ and let $f,g\in \BB= \RR[s,t,\TT_1,\ldots, \TT_{\RR}]$. If $f,g\in
(u,v)\BB$, for some ideal $(u,v)\subset \RR $,
the form derived from
\[
\left[ \begin{array}{r}
f \\g \\
\end{array} \right]=
\left[ \begin{array}{rr}
a & b \\c & d \\
\end{array} \right]
\left[ \begin{array}{r}
u \\v \\
\end{array} \right],\]
\[ h=ad-bc = \det(f,g)_{(u,v)}, \]
will be called a basic Sylvester form. An inspection of the
degrees  gives that
if $\deg_{\RR}(f)< \min\{\deg_{\RR}(g),
\deg_{\RR}(u)+\deg_{\RR}(v)\}$, then
\begin{eqnarray*} \deg_{\RR}(h)&<& \deg_{\RR}(g),\\
\deg_{\TT}(h)&=& \deg_{\TT}(f)+\deg_{\TT}(g).
\end{eqnarray*}

\begin{Remark}{\rm
In our applications, $I= C(f,g)$,
 the content ideal of $f,g$.
In some of these cases, $C(f,g) = (s,t)^n$, for some $n$,
an ideal which admits the irreducible
decomposition
\[ (s,t)^n = \bigcap_{i=1}^n (s^i, t^{n+1-i}).\]
One can then
 process $f,g$ through all the pairs $\{s^i, t^{n-i+1}\}$,
and collect  the
determinants for the next round of elimination.  As in the
classical Sylvester forms, the inclusion $C(f,g)\subset (s,t)^n$ may
be used anyway to start the process,  although without the measure of control of degrees
afforded  by the  equality of ideals.
}\end{Remark}

\subsubsection*{Cohen-Macaulay algebras}\index{Cohen--Macaulay algebras:
dimension two}

The following is  a well-known test of Cohen-Macaulayness (see
\cite{syl1}):

\begin{Theorem} \label{cmtest} If $\dim R\geq 2$ and
  the reduction number of $I$ is
$\leq 1$, that is $I^2=JI$, then $\RR[It]$ is Cohen-Macaulay. The
converse holds if $\dim \RR=2$.
\end{Theorem}

We pointed out in Theorem~\ref{cmtest} that the
 basic control of Cohen-Macaulayness of a Rees algebra
of an ideal $I\subset  k[s,t]$ is that its reduction number be at
most $1$. We next give a mean of checking this property directly off
a free presentation of $I$.

\begin{Theorem}\label{macaulay_case} Let $I\subset \RR$ be an ideal  of codimension $2$, minimally
generated by $3$ forms of the same degree. Let
\[ \varphi= \left[ \begin{array}{rr}
\alpha_1 & \alpha_2 \\
\beta_1 & \beta_2 \\
\gamma_1 & \gamma_2\\
\end{array} \right] \]
be the Hilbert-Burch presentation matrix of $I$. Then $\Rees$ is
Cohen-Macaulay if and only if the equalities of ideals of $\RR$ hold
\[ (\alpha_1, \beta_1, \gamma_1)=(\alpha_2, \beta_2, \gamma_2)=(u,v),
\] where $u,v$ are forms.
\end{Theorem}

\begin{proof} Consider  the presentation
\[ 0 \rar\ {\LL} \lar \Symi(I)=\RR[\TT_1, \TT_2, \TT_3]/(f,g) \lar \Rees \rar 0,\]
where $f,g$ are the $1$-forms
\[
\left[
\begin{array}{r}
f \\ g
\end{array} \right]
=
\left[
\begin{array}{rrr}
\TT_1 & \TT_2 & \TT_3 \\
\end{array} \right] \cdot \varphi.
\]

If $\Rees$ is Cohen-Macaulay, the reduction number of $I$ is $1$ by
Theorem~\ref{cmtest}, so there must be a nonzero quadratic form $h$
with coefficients in $k$ in  the presentation ideal $M$ of $\Rees$.
In addition to $h$, this ideal contains $f,g$,  hence in order to
produce such terms its Hilbert-Burch matrix must be of the form
\[
\left[
\begin{array}{ll}
u & v \\
 p_1 & p_2 \\
q_1 & q_2 \\
\end{array} \right]
\] where $u,v$ are forms of $k[s,t]$, and the other entries are
$1$-forms of $k[\TT_1,\TT_2,\TT_3]$. Since $p_1,p_2$ are $q_1,q_2$ are
pairs of linearly independent $1$-forms, the assertion about the
ideals defined by the columns of $\varphi$ follow.
\end{proof}

\subsubsection*{Base ideals generated in degree $4$}

This is the case treated by D. Cox in his Luminy lecture (\cite{Cox}).
We accordingly change the notation to $\RR=k[s,t]$, $I= (f_1, f_2, f_3)$,
forms of degree $4$. The field $k$ is infinite, and we further assume that $f_1, f_2$ form a regular
sequence so that $J=(f_1, f_2)$ is a reduction of $I$ and of
$(s,t)^4$. Let
\begin{eqnarray}\label{es4}
0 \rar \RR(-4-\mu)\oplus \RR(-8+\mu)  \stackrel{\varphi}{\lar}
\RR^3(-4) \lar \RR\lar
\RR/I\rar 0,
\end{eqnarray}
\begin{eqnarray*}
 \varphi= \left[ \begin{array}{rr}
\alpha_1 & \alpha_2 \\
\beta_1 & \beta_2 \\
\gamma_1 & \gamma_2\\
\end{array} \right]
\end{eqnarray*}
 be the Hilbert-Burch presentation of $I$. We
obtain the equations of $f_1, f_2, f_3$ from this matrix.

\medskip

Note that $\mu$ is the degree of the first column of $\varphi$,
$4-\mu$ the other degree. Let us first consider (as in \cite{Cox}) the
case $\mu=2$.

\medskip

\subsubsection*{Balanced case}

We shall now give a computer-assisted treatment of the {\em
balanced} case, that is when the resolution (\ref{es4}) of the ideal
$I$ has $\mu=2$  and the content ideal of the syzygies is $(s,t)^2$.
Since $k$ is infinite, it is easy to show that there is a change of
variables, $\TT_1, \TT_2, \TT_3 \rar x,y,z$, so that $(s^2, st, t^2)$ is a
syzygy of $I$. The forms $f,g$ that define the symmetric algebra of
$I$ can then be written
\[ [f \quad g] =
 [s^2 \quad st \quad t^2] \left[\begin{array}{rr} x & u \\
y & v \\
z & w \\
\end{array} \right],
\]
where $u,v,w$ are linear forms in $x,y,z$.
Finally, we will assume that
 the ideal
 $I_2\left(\left[ \begin{array}{rrr}
x & y & z \\
u & v & w \\
\end{array} \right]\right)$ has codimension two.
Note that this is a generic condition.

%We will say that $I$ is {\em generic} if on writing

\medskip

We introduce now the {\em equations} of $I$.

\bigskip

 $\bullet$ Linear equations $f$ and $g$:

\begin{eqnarray*}
[f \quad g ] &=& [x \; y \; z]\; \varphi = [x \; y \; z] \left[ \begin{array}{rr}
\alpha_1 & \alpha_2 \\
\beta_1 & \beta_2 \\
\gamma_1 & \gamma_2\\
\end{array} \right] \\
&=& [s^2 \quad st \quad t^2]  \left[\begin{array}{rr} x & u \\
y & v \\
z & w \\
\end{array} \right],
\end{eqnarray*}
where $u,v,w$ are linear forms in $x,y,z$.

\bigskip

{ $\bullet$ Biforms $h_1$ and $h_2$:}

\bigskip

Write $\Gamma_1 $ and $\Gamma_2$ such that
\[
[f \quad g ]=[x \; y \; z]\; \varphi =[\;s \quad t^2\;]\;
 \Gamma_1 = [\;s^2 \quad t\;]\; \Gamma_2.
\]

Then $h_1=\det \Gamma_1$ and $h_2=\det \Gamma_2$.

\bigskip
{ $\bullet$ Implicit equation $F=\det \Theta$,
 where $[ h_1 \quad h_2] = [s \quad t]\; \Theta$.}

\bigskip

Using generic entries for $\varphi$, in place of the true $k$-linear
forms in old variables $x,y,z$, we consider the ideal of
$\AA=k[s,t,x,y,z,u,v,w]$ defined by

\begin{eqnarray*}  f &=& s^2x+sty+t^2z \\
 g &=&  s^2u+stv+t^2w\\
h_1&=& - syu - tzu + sxv + txw\\
h_2 &=& - szu - tzv + sxw + tyw\\
F &=& -z^2u^2 + yzuv-xzv^2-y^2uw+2xzuw+xyvw-x^2w^2
\end{eqnarray*}

\begin{Proposition} \label{case24}
If $I_2\left(\left[ \begin{array}{rrr}
x & y & z \\
u & v & w \\
\end{array} \right]\right)
$ specializes to a codimension two ideal of $k[x,y,z]$, then
 $ L = (f,g,h_1,h_2,F)\subset \AA = \RR[x,y,z,u,v,w]$
 specializes
to the defining ideal of $\Rees$.
\end{Proposition}
%{\bf \% \% The specialization assumption in the statement at least
%implies that the original column entries of original $\varphi$ do
%not generate the same $2$-generated ideal of $R$ (because of
%Theorem~\ref{macaulay_case}); it also implies that the
%parametrization is proper. The question is whether it is equivalent
%to these. In any case, something has to be said or clarified here.
%Same for degree $5$ and ff.

%Aron: It seems better to just assert the fact. What do you say/}

\medskip

\begin{proof}
{\em Macaulay2} (\cite{Macaulay2}) gives a resolution
\[ 0 \rar \AA \stackrel{d_2}{\lar} \AA^5 \lar \AA^5 \lar L \rar 0\]
where
\[ d_2 = \left[
\begin{array}{c}
zv-yw \\
zu-xw \\
-yu+xv \\
\kern-2pt-t \\
s\\
\end{array} \right].
\]
The assumption on
 $I_2\left(\left[ \begin{array}{rrr}
x & y & z \\
u & v & w \\
\end{array} \right]\right)$ says that the entries of $d_2$ generate
an ideal of codimension four and thus
implies that the specialization $LS$ has projective
dimension two and that it is unmixed. Since $LS \not \subset (s,t)S$,
 there is an element  $q\in (s,t)\RR$
 that is regular modulo $S/LS$. If
\[ LS = Q_1 \cap \cdots \cap Q_{r}\]
is the primary decomposition of $LS$, the
 localization
 $LS_q$ has the corresponding decomposition since $q$ is not contained
 in any of the $\sqrt{Q_i}$. But now $\Sym_q=\Rees_q$, so
 $LS_q=(f,g)_u$, as $I_q=\RR_q$.
\end{proof}

\subsubsection*{Non-balanced case}

%\subsubsection*{Non-generic case}

We shall now give a similar computer-assisted treatment of the
non-balanced
case, that is when the resolution (\ref{es4}) of the ideal $I$ has
$\mu=3$. This implies hat the content ideal of the syzygies is
$(s,t)$.
Let us first indicate how the proposed algorithm would
behave.

\begin{itemize}

\item[{$\bullet$}] Write the forms $f,g$ as
\begin{eqnarray*}
f & = & a s + b t \\
g & = & c s + d t, \\
\end{eqnarray*}
where
\[
 \left[ \begin{array}{r}
c \\[5pt]
d \\
\end{array}\right]
=
 \left[ \begin{array}{rrr}
x & y & z \\[5pt]
u & v & w \\
\end{array}\right]
 \left[ \begin{array}{l}
s^2 \\[5pt]
st \\[5pt]
t^2
\end{array}\right]
\]

\item[{$\bullet$}] The next form is the Jacobian of $f,g$ with
respect to $(s,t)$
\[ h_1 = \det(f,g)_{(s,t)}= ad - bc =
 - bxs^2 - byst - bzt^2 + aus^2  + avst + awt^2.
\]
\item[{$\bullet$}] The next two generators
\[ h_2 = \det(f,h_1)_{(s,t)} =
 b^2x s + b^2y t - abzt - abus - abvt + a^2w t
\]
and the elimination equation
\[ h_3 = \det(f,h_2)_{(s,t)} =
 - b^3 x + ab^2 y - a^2 bz + ab^2 u - a^2 bv + a^3 w.
\]
%\item[{$\bullet$}]

%\item[{$\bullet$}]

\end{itemize}

\begin{Proposition} \label{case41}
$ L = (f,g,h_1,h_2,h_3)\subset \AA = k[s,t, x,y,z,u,v,w]$
  specializes
to the defining ideal of $\Rees$.
\end{Proposition}

\begin{proof}
{\em Macaulay2} (\cite{Macaulay2}) gives the following resolution
of $L$
\[ 0 \rar \AA^2 \stackrel{\varphi}{\lar} \AA^6 \stackrel{\psi}{\lar}
\AA^5
\lar L \rar 0,\]
\begin{scriptsize}
\[
\varphi =
\left[
\begin{array}{rr}
s & 0 \\
t & 0 \\
-b & s \\
a & t \\
0 & -b \\
0 & a \\
\end{array}
\right], \]
\end{scriptsize}
\begin{tiny}
\[ \psi =
\left[
\begin{array}{cccccc}
-b^2x+abu & -b^2y+abz+abv-a^2w & -bsx-bty+asu+atv & -btz+atw &
   -s^2x-sty-t^2z & -s^2u-stv-t^2w \\[5pt]
t         &       -s           &        0         &    0     &
      0           &         0      \\[5pt]
a         &        b           &        t         &    -s    &
      0           &         0      \\[5pt]
0         &        0           &        a         &     b    &
      t           &        -s      \\[5pt]
0         &        0           &        0         &     0    &
      a           &         b      \\
\end{array} \right]
\]
\end{tiny}

The ideal of $2\times 2$ minors of $\varphi$ has codimension $4$, even
after we specialize from $A$ to $S$ in the natural manner.
Since $LS$ has projective dimension two, it will be
 unmixed. As $LS\not \subset (s,t)$, there is an element  $u\in (s,t)R$
 that is regular modulo $S/LS$. If
\[ LS = Q_1 \cap \cdots \cap Q_{\RR}\]
is the primary decomposition of $LS$, the
 localization
 $LS_u$ has the corresponding decomposition since $u$ is not contained
 in any of the $\sqrt{Q_i}$. But now $\Sym_u=\Rees_u$, so
 $LS_u=(f,g)_u$, as $I_u=R_u$.
 \end{proof}

\subsubsection*{Degree $5$ and above}

It may be worthwhile to extend this to arbitrary degree, that is
assume that $I$ is defined by $3$ forms of degree $n+1$ (for
convenience in the notation to follow). We first consider the case
$\mu=1$.
 Using the procedure above, we
would obtain the  sequence of polynomials in $\AA= \RR[a,b,x_1, \ldots,
x_n, y_1,\ldots, y_n]$

\begin{itemize}

\item[{$\bullet$}] Write the forms $f,g$ as
\begin{eqnarray*}
f & = & a s + b t \\
g & = & c s + d t, \\
\end{eqnarray*}
where
\[
 \left[ \begin{array}{r}
c \\
d \\
\end{array}\right]
=
 \left[ \begin{array}{rrr}
x_1 & \cdots  & x_n \\
y_1 & \cdots  & y_n \\
\end{array}\right]
 \left[ \begin{array}{c}
s^{n-1} \\
s^{n-2}t \\
\vdots \\
st^{n-2}\\
t^{n-1}
\end{array}\right]
\]

\item[{$\bullet$}] The next form is the Jacobian of $f,g$ with
respect to $(s,t)$
\[ h_1 = \det(f,g)_{(s,t)}= ad - bc
\]
\item[{$\bullet$}] Successively we would set
\[ h_{i+1} = \det(f,h_i)_{(s,t)}, \quad 1 < n.
\]
\item[{$\bullet$}]  The polynomial
\[h_{n} = \det(f,h_{n-1})_{(s,t)}\]
is   the elimination equation.
%\item[{$\bullet$}]

%\item[{$\bullet$}]

\end{itemize}

\begin{Proposition} \label{case15} $ L = (f,g,h_1,\ldots, h_5)\subset
\AA $
 specializes
to the defining ideal of $\Rees$.
\end{Proposition}

In {\em Macaulay2}, we checked the
degrees $5$ and $6$ cases. In both cases, the ideal $L$ (which has
one more generator in degree $6$) has a projective resolution of
length $2$ and the ideal of maximal minors of the last map has
codimension four.

\begin{Conjecture} \label{mu1dimtwo}{\rm For arbitrary $n$,
$ L = (f,g,h_1,\ldots, h_n).$}
\end{Conjecture}

In degree $5$, the interesting case is when the Hilbert-Burch matrix
$\phi$ has degrees $2$ and $3$. Let us describe the proposed
generators. For simplicity, by a change of coordinates, we assume
that the coordinates of the degree $2$ column of $\varphi$ are
$s^2,st,t^2\,$

\begin{eqnarray*}
f & = & s^2x+sty+t^2z\\
g & = &
(s^3w_1+s^2tw_2+st^2w_3+t^3w_4)x+(s^3w_5+s^2tw_6+st^2w_7+t^3w_8)y\\
& + & (s^3w_9+s^2tw_{10}+st^2w_{11}+t^3w_{12})z
\end{eqnarray*}

Let
\[\begin{array}{lllll}
\left[ \begin{array}{c} f \\[5pt]
 g \end{array} \right] &=&
\left[ \begin{array}{lll} \kern3pt x &\quad\quad\quad y &\quad\;  z
\\[5pt]
sA &\quad  sB+tC &\quad tD
\end{array}
\right]
\left[ \begin{array}{c} s^2 \\
 st \\
  t^2 \end{array}
  \right] &=&
  \phi\left[ \begin{array}{c} s^2 \\
   st \\
    t^2 \end{array} \right] \\
    &&&&\\
&=& \left[ \begin{array}{ll}\kern15pt x &\kern17pt ys+zt  \\[5pt]
 sA+tB &\quad stC+t^2D
 \end{array} \right]
 \left[ \begin{array}{c} s^2 \\[5pt]
  t \end{array} \right]&=& B_1
  \left[ \begin{array}{c} s^2 \\[5pt]
   t \end{array} \right] \\&&&&\\
&=& \left[ \begin{array}{ll}\kern8pt xs+yt &\kern25pt z  \\[5pt]
 s^2A+stB &\quad
sC+tD \end{array} \right] \left[ \begin{array}{c} s \\[5pt]
 t^2
\end{array} \right] &=& B_2\left[ \begin{array}{c} s \\[5pt]
 t^2
\end{array} \right]
\end{array},
\]
where $A,B,C,D$ are $k$-linear forms in $x,y,z$.

\bigskip

\[\begin{array}{lll}
h_1 &= & \det (B_1) \\
&=& s^2(-yA)+st(xC-yB-zA)+t^2(xD-zB) \\
&=& s^2(-yA)+ t (xCs-yBs-zAs+xDt-zBt) \\
&=& s(-yAs+xCt-yBt-zAt )+t^2(xD-zB ), \\ &&\\
h_2 &=& \det (B_2) \\
&=& s^2(xC-zA)+st(xD+yC-zB)+t^2(yD) \\
&=& s^2(xC-zA)+t(xDs+yCs-zBs+yDt ) \\
&=& s(xCs-zAs+xDt+yCt-zBt )+t^2(yD ).
\end{array}
\]

\bigskip

\[\begin{array}{lllll}
\left[ \begin{array}{c} f \\[5pt]
 h_1 \end{array} \right] &=&
 \left[\begin{array}{ll} \kern12pt x &\kern62pt ys+zt \\
  -yA &\quad xCs-yBs-zAs+xDt-zBt  \end{array}\right]
  \left[ \begin{array}{c} s^2 \\[5pt]
   t \end{array} \right] &=& C_1
  \left[ \begin{array}{c} s^2 \\[5pt]
   t \end{array} \right] \\
    &&&&\\
&=&\left[\begin{array}{ll}\kern45pt xs+yt &\kern30pt z\\[5pt]
 -yAs+xCt-yBt-zAt &\quad
xD-zB   \end{array}\right]
\left[ \begin{array}{c} s \\[5pt]
 t^2 \end{array} \right]&=&C_2
 \left[ \begin{array}{c} s \\[5pt]
  t^2 \end{array} \right]
\end{array}
\]

\bigskip

\[\begin{array}{lllll}
\left[ \begin{array}{c} f \\[5pt]
 h_2 \end{array} \right] &=&
 \left[\begin{array}{ll} \kern18pt x &\kern62pt ys+zt \\[5pt]
  xC-zA  &\quad xDs+yCs-zBs+yDt \end{array}\right]
  \left[ \begin{array}{c} s^2 \\[5pt]
   t \end{array} \right] &=& C_3
   \left[ \begin{array}{c} s^2 \\[5pt]
    t \end{array} \right] \\ &&&&\\
&=&\left[\begin{array}{ll}\kern45pt xs+yt &\kern15pt z  \\[5pt]
 xCs-zAs+xDt+yCt-zBt &\quad yD   \end{array}\right]
 \left[ \begin{array}{c} s \\[5pt]
  t^2 \end{array} \right]&=&C_4
  \left[ \begin{array}{c} s \\[5pt]
   t^2 \end{array} \right]
\end{array}
\]

 \bigskip

 \[\begin{array}{lllll}
 c_1 &=& \det(C_1) &=& x^2(Cs+Dt)+xy(-Bs)+xz(-As-Bt)+yz(At)+y^2(As)\\ &&&&\\
 c_2 &=& \det(C_2) &=& x^2(Ds) + xy(Dt) + xz(-Bs-Ct)+yz(As)+z^2(At)\\&&&&\\
 c_3 &=& \det(C_3) &=& x^2(Ds) + xy(Dt) +xz(-Bs-Ct) + yz(As)+ z^2(At)\\ &&&&\\
 c_4 &=& \det(C_4) &=&  xy(Ds) + xz(-Cs-Dt) + yz(-Ct)+ z^2(As+Bt) + y^2(D)
 \end{array}
 \]

 \bigskip

  \[
 \left[ \begin{array}{c} f \\[5pt]
  h_1 \\
   h_2 \end{array} \right]=
 \left[ \begin{array}{lll}\kern15pt x &\kern38pt y &\kern30pt z \\[5pt]
 \kern5pt -yA &\quad xC-yB-zA &\quad xD-zB \\[5pt]
   xC-zA &\quad xD+yC-zB &\kern25pt yD  \end{array}\right]
   \left[ \begin{array}{c} s^2 \\[5pt]
    st \\[5pt]
     t^2 \end{array} \right]
 \]

 \bigskip

 Then $F= -x^3D^2+x^2yCD+xy^2(-BD) + x^2z(2BD-C^2)+xz^2(2AC-B^2)
 +xyz(BC-3AD)+y^2z(-AC)+yz^2(AB)+y^3(AD)+z^3(-A^2)$, an equation of degree $5$.
 In particular, the parametrization is birational.

\begin{Proposition} \label{case25}
 $L=(f,\, g,\, h_1,\, h_2,\, c_1,\, c_2,\, c_4,\, F)$
  specializes
to the defining ideal of $\Rees$.
\end{Proposition}

\begin{proof} Using {\it Macaulay2},
the ideal $L$ has a resolution:
\[
0 \longrightarrow S^1 \stackrel{d_3}{\longrightarrow} S^6
\stackrel{d_2}{\longrightarrow}  S^{12}
\stackrel{d_1}{\longrightarrow}   S^8 \longrightarrow L \longrightarrow 0.
\]

\[
d_3= [-z \;\; y \;\; x \;\; -t \;\; s \;\; 0   ]^t
\]
\begin{scriptsize}
\[
d_2= \left[ \begin{array}{llllll} y & z&  0&  0 & 0&  \kern130pt 0 \\
  x & 0 & z&  0&  0 &  \kern130pt 0   \\
\kern-5pt -v & 0&  0&  z&  0&  \kern75pt x^2w_4-xzw_7+xyw_8+xzw_{12} \\
u&  0 & 0&  0 & z & -xzw_3+xyw_4+z^2w_6-yzw_7+y^2w_8-xzw_8-z^2w_{11}+yzw_{12} \\
 0 &  x& \kern-5pt -y& 0&  0& \kern130pt 0\\
 0 & \kern-5pt -v& 0& \kern-5pt -y& 0&  \kern60pt xzw_1-x^2w_3+yzw_5+z^2w_9-xzw_{11}\\
 0 & u & 0 & 0 &\kern-5pt -y & \kern70pt xzw_2-x^2w_4+z^2w_{10}-xzw_{12}  \\
 0 & 0&  u & 0& \kern-5pt -x & \kern60pt xzw_1+yzw_5-xzw_6+x^2w_8+z^2w_9\\
 0 & 0&  0&  u&  v& \kern130pt 0\\
 0 & 0&  v&  x&  0 & \kern15pt -xyw_1+x^2w_2-y^2w_5+xyw_6-x^2w_7-yzw_9+xzw_{10} \\
 0 & 0&  0&  0&  0& \kern125pt -t\\
 0 & 0&  0&  0&  0& \kern130pt s \end{array} \right]
\]
\end{scriptsize}

The ideals of maximal minors give $\codim I_1(d_3)=5$ and $\codim
I_{5}(d_2)= 4$ after specialization. As we have been arguing, this
suffices to show that the specialization is a prime ideal of
codimension two. 
\end{proof}

\subsubsection*{Elimination forms in higher degree}

In degrees greater than $5$, the methods above are not very suitable.
However, in several cases they are still supple enough to produce the
elimination equation. We have already seen this when one of the
syzygies is of degree $1$. Let us describe two other cases.

\bigskip

$\bullet$ Degree $n=2p$, $f$ and $g$ both of degree $p$. We use the
decomposition
\[ (s,t)^p = \bigcap_{i=1}^p (s^i, t^{p+1-i}).\]
For each $1\leq i \leq p$, let
\[ h_i = \det(f,g)_{(s^i, t^{p+1-i})}.\]
These are quadratic polynomials with coefficients in $(s,t)^{p-1}$.
We set
\[ [h_1, \cdots, h_p] = [s^{p-1}, \cdots, t^{p-1}]\cdot \AA,\]
where $\AA$ is a $p\times p$ matrix whose entries are $2$-forms in
$k[\TT_1,\TT_2,\TT_3]$. The Sylvester form of degree $n$, $F = \det(\AA)$, is the
required elimination equation.

\bigskip

$\bullet$ Degree $n=2p+1$, $f$  of degree $p$. We use the
decomposition
\[ (s,t)^p = \bigcap_{i=1}^p (s^i, t^{p+1-i}).\]
For each $1\leq i \leq p$, let
\[ h_i = \det(f,g)_{(s^i, t^{p+1-i})}.\]
These are quadratic polynomials with coefficients in $(s,t)^{p}$.
We set
\[ [f, h_1, \cdots, h_p] = [s^{p}, \cdots, t^{p}]\cdot \AA,\]
where $\AA$ is a $(p+1)\times (p+1)$ matrix with one column whose
entries are linear forms and the remaining columns with entries
$2$-forms in $k[\TT_1,\TT_2,\TT_3]$. The Sylvester form $F = \det(\AA)$ is the
required elimination equation.

\printindex
\end{document}